\documentclass[10pt,reqno,a4paper]{amsbook}
\pdfoutput=1
\usepackage{etex}
\usepackage[utf8]{inputenc}
\usepackage[english]{babel}
\usepackage{etoolbox,chngcntr}
\usepackage{amsmath,amsthm,amssymb,array,stmaryrd,color,graphicx,mathtools,multirow,setspace}
\usepackage{soul}\setul{0.3ex}{}
\usepackage{bussproofs}
\usepackage{xspace}
\usepackage{longtable}
\usepackage{booktabs}
\usepackage[protrusion=true,expansion=true]{microtype}
\usepackage[bookmarksdepth=2,pdfencoding=auto]{hyperref}

\expandafter\def\csname opt@stmaryrd.sty\endcsname
{only,shortleftarrow,shortrightarrow}
\usepackage{extpfeil}
\newextarrow{\xbigtoto}{{20}{20}{20}{20}}
   {\bigRelbar\bigRelbar{\bigtwoarrowsleft\rightarrow\rightarrow}}

\usepackage[all]{xy}
\usepackage{tikz}
\usetikzlibrary{calc,shapes.callouts,shapes.arrows,matrix,patterns}
\newcommand{\hcancel}[5]{%
    \tikz[baseline=(tocancel.base)]{
        \node[inner sep=0pt,outer sep=0pt] (tocancel) {#1};
        \draw[red, line width=0.3mm] ($(tocancel.south west)+(#2,#3)$) -- ($(tocancel.north east)+(#4,#5)$);
    }%
}

\usepackage[natbib=true,style=numeric,maxnames=10]{biblatex}
\usepackage[babel]{csquotes}
\bibliography{bibliography}

\theoremstyle{definition}
\newtheorem{defn}{Definition}[section]
\newtheorem{ex}[defn]{Example}

\theoremstyle{plain}
\newtheorem{prop}[defn]{Proposition}
\newtheorem{cor}[defn]{Corollary}
\newtheorem{lemma}[defn]{Lemma}
\newtheorem{thm}[defn]{Theorem}
\newtheorem{scholium}[defn]{Scholium}

\theoremstyle{remark}
\newtheorem{rem}[defn]{Remark}

\newtheorem{speculation}[defn]{Speculation}
\newtheorem{caveat}[defn]{Caveat}
\newtheorem{conjecture}[defn]{Conjecture}

\newcommand{\ZZ}{\mathbb{Z}}
\newcommand{\FF}{\mathbb{F}}
\renewcommand{\AA}{\mathbb{A}}
\newcommand{\A}{\mathcal{A}}
\newcommand{\B}{\mathcal{B}}
\renewcommand{\C}{\mathcal{C}}

\newcommand{\E}{\mathcal{E}}
\newcommand{\F}{\mathcal{F}}
\renewcommand{\G}{\mathcal{G}}
\let\acuteH\H
\newcommand{\konig}{K\acuteH onig}
\renewcommand{\H}{\mathcal{H}}
\renewcommand{\O}{\mathcal{O}}
\newcommand{\K}{\mathcal{K}}
\newcommand{\N}{\mathcal{N}}
\newcommand{\M}{\mathcal{M}}
\renewcommand{\L}{\mathcal{L}}
\renewcommand{\P}{\mathcal{P}}
\newcommand{\R}{\mathcal{R}}
\newcommand{\T}{\mathcal{T}}
\newcommand{\I}{\mathcal{I}}
\newcommand{\J}{\mathcal{J}}
\renewcommand{\S}{\mathcal{S}}
\renewcommand{\U}{\mathcal{U}}
\newcommand{\V}{\mathcal{V}}
\newcommand{\NN}{\mathbb{N}}
\newcommand{\PP}{\mathbb{P}}
\newcommand{\RR}{\mathbb{R}}
\newcommand{\CC}{\mathbb{C}}
\newcommand{\QQ}{\mathbb{Q}}
\newcommand{\GG}{\mathbb{G}}
\newcommand{\TT}{\mathbb{T}}
\newcommand{\aaa}{\mathfrak{a}}
\newcommand{\bbb}{\mathfrak{b}}
\newcommand{\ccc}{\mathfrak{c}}
\newcommand{\ppp}{\mathfrak{p}}
\newcommand{\qqq}{\mathfrak{q}}
\newcommand{\mmm}{\mathfrak{m}}
\newcommand{\nnn}{\mathfrak{n}}
\newcommand{\Hom}{\mathrm{Hom}}
\newcommand{\HOM}{\mathcal{H}\mathrm{om}}
\newcommand{\EXT}{\mathcal{E}\mathrm{xt}}
\newcommand{\TOR}{\mathcal{T}\mathrm{or}}
\newcommand{\Ind}{\mathrm{Ind}}
\newcommand{\id}{\mathrm{id}}

\newcommand{\placeholder}{\underline{\quad}}
\let\oldul\ul
\renewcommand{\ul}[1]{\text{\oldul{$#1$}}}
\newcommand{\Set}{\mathrm{Set}}

\newcommand{\Sh}{\mathrm{Sh}}
\newcommand{\PSh}{\mathrm{PSh}}
\newcommand{\Zar}{\mathrm{Zar}}
\newcommand{\Et}{\mathrm{\acute{E}t}}
\newcommand{\fppf}{\mathrm{fppf}}
\newcommand{\ph}{\mathrm{ph}}
\newcommand{\surj}{\mathrm{surj}}
\newcommand{\Sch}{\mathrm{Sch}}
\newcommand{\Aff}{\mathrm{Aff}}
\newcommand{\Mod}{\mathrm{Mod}}
\newcommand{\Alg}{\mathrm{Alg}}
\newcommand{\Ring}{\mathrm{Ring}}
\newcommand{\LocRing}{\mathrm{LocRing}}
\newcommand{\RL}{\mathrm{RL}}
\newcommand{\LRL}{\mathrm{LRL}}
\newcommand{\LRS}{\mathrm{LRS}}
\newcommand{\LRT}{\mathrm{LRT}}
\newcommand{\RT}{\mathrm{RT}}
\newcommand{\pt}{\mathrm{pt}}
\newcommand{\tors}{\mathrm{tors}}
\newcommand{\lfp}{\mathrm{lfp}}
\newcommand{\fp}{\mathrm{fp}}
\DeclareMathOperator{\Spec}{Spec}
\DeclareMathOperator{\Proj}{Proj}

\newcommand{\RelSpec}{\operatorname{\ul{\mathrm{Spec}}}}
\newcommand{\RelProj}{\operatorname{\ul{\mathrm{Proj}}}}
\newcommand{\op}{\mathrm{op}}
\DeclareMathOperator*{\colim}{colim}
\DeclareMathOperator{\rank}{rank}
\DeclareMathOperator{\Ann}{Ann}
\DeclareMathOperator{\Int}{int}
\DeclareMathOperator{\Clos}{cl}
\DeclareMathOperator{\Kernel}{ker}
\DeclareMathOperator{\cok}{cok}
\DeclareMathOperator{\im}{im}
\DeclareMathOperator{\supp}{supp}
\newcommand{\Ass}{\mathrm{Ass}}
\newcommand{\Sym}{\mathrm{Sym}}
\newcommand{\Gr}{\mathrm{Gr}}
\newcommand{\Open}{\T}
\newcommand{\?}{\,{:}\,}
\newcommand{\hg}{\mathbin{:}}  
\renewcommand{\_}{\mathpunct{.}\,}
\newcommand{\speak}[1]{\ulcorner\text{\textnormal{#1}}\urcorner}
\newcommand{\Ll}{\vcentcolon\Longleftrightarrow}
\newcommand{\notat}[1]{{!#1}}
\newcommand{\lra}{\longrightarrow}
\newcommand{\lhra}{\ensuremath{\lhook\joinrel\relbar\joinrel\rightarrow}}

\newcommand{\brak}[1]{{\llbracket{#1}\rrbracket}}
\newcommand{\sdense}{{\widehat\Box}}
\newcommand{\sdenseother}{\Box}
\newcommand{\ie}{i.\,e.\@\xspace}

\newcommand{\vs}{vs.\@\xspace}
\newcommand{\resp}{resp.\@\xspace}
\newcommand{\inv}{inv.\@}
\newcommand{\notnot}{\emph{not~not}\xspace}
\newcommand{\affl}{\ensuremath{{\ul{\AA}^1_S}}\xspace}
\newcommand{\afflx}{\ensuremath{{\ul{\AA}^1_X}}\xspace}
\newcommand{\afflt}{\ensuremath{{\ul{\AA}^1_T}}\xspace}
\newcommand{\afflvi}{\ensuremath{{\ul{\AA}^1_{V_i}}}\xspace}
\newcommand{\affla}{\ensuremath{{\ul{\AA}^1_{\Spec A}}}\xspace}
\newcommand{\afflz}{\ensuremath{{\ul{\AA}^1_{\Spec Z}}}\xspace}
\newcommand{\xra}{\xrightarrow}
\newcommand{\stacksproject}[1]{\cite[{\href{https://stacks.math.columbia.edu/tag/#1}{Tag~#1}}]{stacks-project}}
\newcommand{\apart}{\mathrel{\#}}
\newcommand{\fieldext}{\mathrel{|}}
\newcommand{\effective}{ef{}fective\xspace}

\newenvironment{indentblock}{%
  \list{}{\leftmargin\leftmargin}%
  \item\relax
}{%
  \endlist
}

\newcommand{\XXXh}[1]{}

\newcommand{\defeq}{\vcentcolon=}
\newcommand{\defequiv}{\vcentcolon\equiv}
\newcommand{\seq}[1]{\mathrel{\vdash\!\!\!_{#1}}}

\definecolor{gray}{rgb}{0.7,0.7,0.7}

\title{Using the internal language of toposes in algebraic geometry}
\author{Ingo Blechschmidt}
\email{iblech@web.de}

\makeatletter
\counterwithout{section}{chapter}
\counterwithout{footnote}{chapter}
\counterwithout{table}{chapter}
\counterwithout{figure}{chapter}
\patchcmd{\@thm}{\let\thm@indent\indent}{\let\thm@indent\noindent}{}{}
\patchcmd{\@thm}{\thm@headfont{\scshape}}{\thm@headfont{\bfseries}}{}{}
\patchcmd{\@makechapterhead}{\chaptername}{Part}{}{}
\patchcmd{\@chapter}{\chaptername}{Part}{}{}
\patchcmd{\@schapter}{\chaptername}{Part}{}{}
\addto\captionsenglish{\renewcommand\chaptername{Part}}
\def\l@section{\@tocline{1}{0pt}{1pc}{}{}} 
\def\l@chapter{\@tocline{-1}{12pt}{0pt}{}{\bfseries}}

\newcommand{\nocontentsline}[3]{}
\newcommand{\tocless}[1]{\let\addcontentsline=\nocontentsline}
\parindent=\normalparindent
\renewenvironment{proof}[1][\proofname]{\par
  \pushQED{\qed}%
  \normalfont \topsep6\p@\@plus6\p@\relax
  \trivlist
  \item[\hskip\labelsep
        \itshape
    #1\@addpunct{.}]\ignorespaces
}{%
  \popQED\endtrivlist\@endpefalse
}
\let\@afterindenttrue\@afterindentfalse
\def\subsection{\@startsection{subsection}{2}%
  {0pt}{.5\linespacing\@plus.7\linespacing}{-.5em}%
  {\normalfont\bfseries}}
\makeatother

\begin{document}

\newcommand{\HRule}{\rule{\linewidth}{.6pt}}

\begin{center}
\thispagestyle{empty}

\
\bigskip\bigskip

\HRule \\[0.4cm]
{\huge \bfseries Using the internal language of toposes \\ in algebraic geometry\par}\vspace{0.4cm}
\HRule

\bigskip\bigskip\bigskip\bigskip

\Large

Dissertation
\smallskip

zur Erlangung des akademischen Grades
\bigskip\bigskip

Dr.~rer.~nat.
\bigskip\bigskip

eingereicht an der
\bigskip\bigskip

Mathematisch-Naturwissenschaftlich-Technischen Fakultät
\smallskip

der Universität Augsburg
\bigskip\bigskip

von
\bigskip\bigskip

{\large Ingo Blechschmidt}
\bigskip\bigskip

\bigskip\bigskip
\bigskip\bigskip
\bigskip\bigskip
\bigskip\bigskip
\bigskip\bigskip
\bigskip\bigskip

\includegraphics[width=0.4\textwidth]{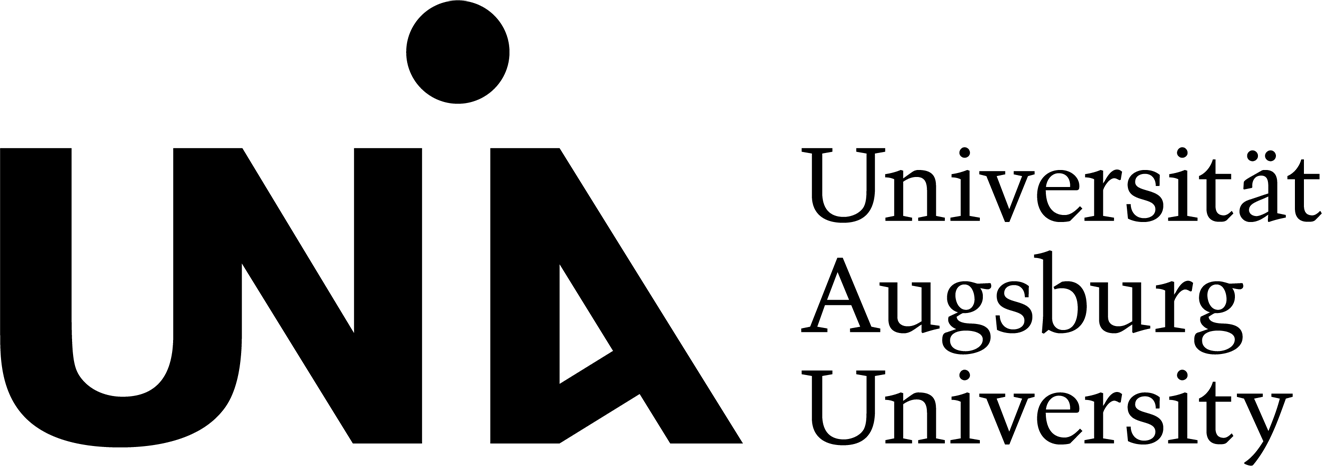}

\bigskip\bigskip
\bigskip\bigskip

Juni 2017

\end{center}

\newpage

\vspace*{54em}

\begin{tabbing}
  External reviewers:\quad \= \kill
  Supervisor: \> Marc Nieper-Wißkirchen, University of Augsburg \\
  \\
  External reviewers: \> Thierry Coquand, University of Gothenburg \\
  \> Frank Herrlich, Karlsruhe Institute of Technology \\
  \\
  Oral examination: \> October 16th, 2017
\end{tabbing}

{\tocless

\chapter*{Abstract}

Any scheme has its associated little and big Zariski toposes. These toposes
support an internal mathematical language which closely resembles the usual
formal language of mathematics, but is ``local on the base scheme'': For
example, from the internal perspective, the structure sheaf looks like an
ordinary local ring (instead of a sheaf of rings with local stalks) and vector
bundles look like ordinary free modules (instead of sheaves of modules
satisfying a local triviality condition). The translation of internal statements and
proofs is facilitated by an easy mechanical procedure.

We investigate how the internal language of the little Zariski topos
can be exploited to give simpler definitions and more conceptual
proofs of the basic notions and observations in algebraic geometry.
To this end, we build a dictionary relating internal and external notions and
demonstrate its utility by giving a simple proof of Grothendieck's generic
freeness lemma in full generality. We also employ this framework to state a
general transfer principle which relates modules with their induced quasicoherent
sheaves, to study the phenomenon that some
properties spread from points to open neighborhoods, and to compare general
notions of spectra.

We employ the big Zariski topos to set up the foundations of a synthetic account
of scheme theory. This account is similar to the synthetic account of
differential geometry, but has a distinct algebraic flavor. Central to the
theory is the notion of synthetic quasicoherence, which has no analogue in
synthetic differential geometry. We also discuss how various common subtoposes
of the big Zariski topos can be described from the internal point of view and
derive explicit descriptions of the geometric theories which are
classified by the fppf and by the surjective topology.


\chapter*{Acknowledgments}

I'm immensely grateful to Marc Nieper-Wißkirchen for his excellent supervision and
careful guidance, for trusting me with a lot of freedom, and for shaping me
mathematically. In this last regard I also acknowledge Jost-Hinrich Eschenburg,
Hansjörg Kielhöfer, and the contributors of the nLab.

These notes owe a special debt to Mike Shulman's work on the stack semantics of
toposes~\cite{shulman:stack}, without which this project wouldn't have been
started, and they greatly profited from insightful discussions with several people,
including Peter Arndt, Filip Bar, Andrej Bauer, Martin Brandenburg, Olivia
Caramello, Joost van Dijk, Adam Epstein, Felix Geißler, Simon Henry, Pol van Hoften, Matthias Hutzler, Simon
Kapfer, Georg Lehner, Guilhermo Frederico Lima de Carvalho e Silva, Tadeusz Litak, Zhen Lin
Low, Giovanni Morando, and Alexander Oldenziel. I'm grateful for their many valuable suggestions. I
particularly thank Adam Epstein for his great hospitality during my stay in
Warwick and Jiří Adámek, Jürgen Koslowski, and Ravi Vakil for their
much-appreciated encouragement.

I thank Tim Baumann, Joost van Dijk, Kathrin Gimmi, Matthias Hutzler, Sori Lee, Alexander Oldenziel, Lukas Stoll, and Marcus Zibrowius for their
careful readings of earlier drafts.

Dear to my heart is doing mathematics with children and school students. It was
a pleasure to work together with Kathrin Helmsauer and Sven Prüfer over the
last few years to establish and run Augsburg's Matheschülerzirkel, and I'm
grateful to the many people who invested their time into this project and shared our
vision.

In particular, I thank my close family consisting of my sisters Dorothea and
Veronika Pörtge and my mother Gotlind Blechschmidt, who have accompanied me for
all my life and always encouraged and inspired me, Meru Alagalingam, Martin Frieb,
Michael Hartmann, Kathrin Helmsauer, Matthias Hutzler,
Christian Hübschmann, Jorid Kretzschmar, Jörg Lehmann, Alex\-an\-der Mai, Ninwe
Ninos, Sven Prüfer, Peter Quast, Lisa Reisch\-mann, Anna Rubeck, Ma\-xi\-mi\-li\-an
Schlögel, Peter Uebele, Philipp Wacker, my office mates Felix Geißler, Hedwig Heizinger,
Stephanie Zapf, and the wonderful participants of our courses in Papenburg, who
each supported me in a unique way. I'm indebted to Audrey Tang, who had a
strong and lasting influence on me many moons ago.

Here in Augsburg, we have a strong bond between undergraduates, PhD students,
and members of faculty. I'll try to carry our mathematical pizza seminars and
our non-mathematical and slightly nonstandard leisure time activities forward.

}

\setcounter{tocdepth}{1}
\tableofcontents

\chapter{Basics}

\section{Introduction}

{\tocless

\subsection*{Internal language of toposes}

A \emph{topos} is a category which shares certain categorical properties with
the category of sets; the archetypical example is the category of sets, and
the most important example for the purposes of this thesis is the category of
set-valued sheaves on a topological space.

Any topos~$\E$ supports an \emph{internal language}. This is a device which
allows one to \emph{pretend} that the objects of~$\E$ are plain sets and that
the morphisms are plain maps between sets, even if in fact they are not. For
instance, consider a morphism~$\alpha : X \to Y$ in~$\E$. From the \emph{internal
point of view}, this looks like a map between sets, and we can formulate the
condition that this map is surjective; we write this as
\[ \E \models \forall y\?Y\_ \exists x\?X\_ \alpha(x) = y. \]
The appearance of the colons instead of the usual element signs reminds us that
this expression is not to be taken literally --~$X$ and~$Y$ are objects of~$\E$
and thus not necessarily sets. The definition of the internal language is made
in such a way so that the meaning of this internal statement is that~$\alpha$
is an epimorphism. Similarly, the translation of the internal statement
that~$\alpha$ is injective is that~$\alpha$ is a monomorphism.

Furthermore, we can \emph{reason} with the internal language. There is a
metatheorem to the effect that if some statement~$\varphi$ holds from the
internal point of view of a topos~$\E$ and if~$\varphi$ logically implies some
further statement~$\psi$, then~$\psi$ holds in~$\E$ as well. As a simple
example, consider the elementary fact that the composition of surjective maps
is surjective. Interpreting this statement in the internal language of~$\E$, we
obtain the more abstract result that the composition of epimorphisms in~$\E$ is
epic.

There is, however, a slight caveat to this metatheorem. Namely, the internal
language of a topos is in general only \emph{intuitionistic}, not
\emph{classical}. This means that internally, one cannot use the law of
excluded middle~($\varphi \vee \neg\varphi$), the law of double negation
elimination~($\neg\neg\varphi \Rightarrow \varphi$), or the axiom of choice.
For instance, one rendition of the axiom of choice is that any vector space is
free. But it need not be the case that a vector space internal to a topos
is free as seen from the internal perspective: By the technique explained in
this thesis, this would imply the absurd statement that any sheaf of modules on
a reduced scheme is locally free.

The restriction to intuitionistic reasoning is not as confining as it might first
appear, in particular because there is a widely applicable metatheorem ensuring
that statements of a certain form are provable classically if and only if they
are provable intuitionistically. We will discuss practical consequences below (on
page~\pageref{sect:appreciating-intuitionistic-logic}).

\subsection*{Algebraic geometry}
We apply the internal language of toposes to algebraic geometry in two
different ways, corresponding to the two different toposes associated to a
scheme~$X$: the \emph{little Zariski topos} which is just the topos~$\Sh(X)$ of
set-valued sheaves on~$X$, and the \emph{big Zariski topos} which we introduce
below.

The internal language of the little Zariski topos can be used as follows.
The structure sheaf~$\O_X$ of a scheme~$X$ is a sheaf of rings in that its sets of
local sections carry ring structures and these ring structures are compatible
with restriction. From the internal point of view of~$\Sh(X)$,
the structure
sheaf~$\O_X$ looks much simpler: It looks just like a plain ring (and
not a sheaf of rings). Similarly, a sheaf of~$\O_X$-modules looks just like a
plain module over that ring.

This allows to import notions and facts from basic linear and commutative
algebra into the sheaf setting. For instance, it turns out that a sheaf
of~$\O_X$-modules is of finite type if and only if, from the internal
perspective, it is finitely generated as an~$\O_X$-module. Now consider the
following fact of linear algebra: If in a short exact sequence of modules the two
outer ones are finitely generated, then the middle one is too. The usual proof of
this fact is intuitionistically valid and can thus be interpreted in the
internal language. It then \emph{automatically} yields the following more advanced
proposition: If in a short exact sequence of sheaves of~$\O_X$-modules the
two outer ones are of finite type, then the middle one is too.

This example was not in any way special: \emph{Any (intuitionistically valid) theorem
about modules yields a corresponding theorem about sheaves of modules.}

The internal language machinery thus allows us to understand the basic notions
and statements of scheme theory as notions and statements of linear and
commutative algebra, interpreted in a suitable sheaf topos. This brings
conceptual clarity and reduces technical overhead.

In Section~\ref{sect:internal-language}, we explain how the internal language
machinery works, and then develop in Part~\ref{part:little-zariski} a
\emph{dictionary} relating common notions of scheme theory and corresponding
notions of algebra. Once built, this dictionary can be used arbitrarily often.
We stress that no in-depth knowledge of topos theory or categorical logic is
necessary to apply this apparatus.

In simple cases, the internal language can be regarded as a tool for ensuring
that certain kinds of ``fast and loose reasoning'' in algebraic geometry can be
rigorously justified. For instance, when trying to quickly gauge whether some
plausible-looking statement holds for schemes and sheaves, we might content
ourselves to check that the statement holds for rings and modules and then trust
that it also holds in the general case. Or when trying to construct a certain
sheaf of modules, we might content ourselves to construct it over affine open
subsets and then appeal to some gluing lemma, without meticulously checking the
details.

The internal language apparatus ensures that this kind of reasoning will never
result in wrong conclusions, provided that one can formulate the statements and
constructions in the internal language and that the correctness proof in the
affine setting is intuitionistically valid.

We believe that already this application of the internal language is useful to
working algebraic geometers. However, more advanced applications are also
possible. They result from considering internal statements whose logical form
is more complex, in particular from statements which quantify over subsets or which
contain implication and negation signs.

For instance, if~$X$ is a reduced scheme,
the internal universe of~$\Sh(X)$ has the peculiar feature that~$\O_X$ is
Noetherian and a field, even if~$X$ is not locally Noetherian and (as will
almost always be the case) the local rings~$\O_{X,x}$ are not fields. This fact
has no simple external counterpart; it's rather an intricate statement about
the interplay between the rings~$\Gamma(U, \O_X)$ for varying open subsets~$U
\subseteq X$.

Thanks to this particular feature, linear and commutative algebra over~$\O_X$
are particularly simple from the
internal point of view. For instance, Grothendieck's generic freeness lemma,
which is usually proved using a somewhat involved series of reduction steps,
admits a short, easy, and conceptual proof with this technique.

To briefly indicate a part of this, let~$\F$ be a sheaf of~$\O_X$-modules of finite
type. A basic version of Grothendieck's generic freeness lemma then states
that~$\F$ is locally free on some dense open subset of~$X$; this fact
is stated in Vakil's lecture notes as an ``important hard
exercise''~\cite[Exercise~13.7.K]{vakil:foag}. In fact, this proposition is just the
interpretation of the following basic statement of intuitionistic linear algebra in
the sheaf topos: Any finitely generated vector space is \notnot free.
The proof of this statement is entirely straightforward.\footnote{Intuitionistically,
the statement that any finitely generated vector space is \emph{free} is stronger than
the doubly negated version and cannot be shown. It would imply that any sheaf
of finite type is not only locally free on some dense open subset, but locally
free on the entire space. We discuss this example in more detail in
Section~\ref{sect:upper-semicontinuous-functions} and in particular in
Lemma~\ref{lemma:locally-free-dense}. A proof of Grothendieck's generic
freeness lemma in its full form is given in
Section~\ref{sect:generic-freeness}.

For concreteness, here is the standard intuitionistic proof that any finitely
generated vector space~$V$ is \notnot free. Let~$(x_1,\ldots,x_n)$ be a
generating family. If~$n = 0$, we are done. Else it's \notnot the case that
either some~$x_i$ can be expressed as a linear combination of the other
vectors, or not. The former implies
that~$(x_1,\ldots,x_{i-1},x_{i+1},\ldots,x_n)$ is a generating family, whereby
we can appeal to induction to obtain that~$V$ is \notnot free. The latter implies that~$(x_1,\ldots,x_n)$ is
linearly independent and therefore a basis. In both cases it follows that~$V$
is \notnot free, therefore~$V$ is indeed \notnot free.

In this argument, we used the intuitionistically valid proof
scheme~$\neg\neg\varphi \wedge (\varphi \Rightarrow \neg\neg\psi) \Longrightarrow
\neg\neg\psi$. We expand on this in
Section~\ref{sect:appreciating-intuitionistic-logic}.}

It is in this way that the internal language unlocks new approaches: by
making concepts accessible which would otherwise be too unwieldy to manage and
by allowing to import a huge corpus of prior work, namely the entire literature
on constructive algebra.

The internal language also sheds light on the phenomenon that
sometimes, truth of a property at a point~$x$ spreads to some open
neighborhood of~$x$; and in particular that sometimes, truth of a property at
the generic point spreads to some dense open subset. For instance, if the stalk
of a sheaf of finite type is zero at some point, the sheaf is even zero on some
open neighborhood; but this spreading does not occur for general sheaves which
may fail to be of finite type.

We formalize this by introducing a \emph{modal operator}~$\Box$ into the
internal language, such that the internal statement~$\Box\varphi$ means
that~$\varphi$ holds on some open neighborhood of~$x$. Furthermore, we
introduce a simple operation on formulas, the~\emph{$\Box$-translation}
$\varphi \mapsto \varphi^\Box$, such that~$\varphi^\Box$ means that~$\varphi$
holds at the point~$x$. This translation is defined on a purely syntactical
level. The question whether truth at~$x$ spreads to truth on a
neighborhood can then be formulated in the following way: Does~$\varphi^\Box$
intuitionistically imply~$\Box\varphi$?

This allows to deal with the question in a simpler, logical way, with the
technicalities of sheaves blinded out. We also give a metatheorem which
covers a wide range of cases. Namely, spreading occurs for all those properties
which can be formulated in the internal language without
using~``$\Rightarrow$'',~``$\forall$'', and~``$\neg$''.

To take up the example above, consider the property of a module~$\F$ being
the zero module. In the internal language, it can be formulated as~$(\forall x\?\F\_ x = 0)$.
Because of the appearance of~``$\forall$'', the metatheorem is not
applicable to this statement. But if~$\F$ is of finite type, there are
generators~$x_1,\ldots,x_n\?\F$ from the internal point of view, and the
condition can be reformulated as~$x_1 = 0 \wedge \cdots \wedge x_n = 0$; the
metatheorem is applicable to this statement.

\subsection*{Synthetic algebraic geometry}
All of the applications mentioned above employ the little Zariski topos of the base
scheme~$X$, the topos of sheaves on the underlying topological space of~$X$.
Its internal language simplifies the treatment of sheaves of rings and modules
over~$X$, but the treatment of~\emph{schemes} over~$X$ is simplified only a
little bit: From the internal point of view of~$\Sh(X)$, a morphism~$T \to X$
of schemes looks like a morphism~$T \to \pt$. Therefore relative scheme theory
is turned into absolute scheme theory (over the ring~$\O_X$), but it still
requires the machinery of locally ringed spaces.

The internal language of the \emph{big} Zariski topos of~$X$ allows for a more far-reaching
change of perspective. It incorporates Grothendieck's functor-of-points
philosophy in order to cast modern algebraic geometry, relative to the arbitrary
base scheme~$X$, in a naive \emph{synthetic} language reminiscent of the classical
Italian school.

The synthetic approach is best explained by contrasting it with the usual
approach to scheme theory, which is to layer it upon some standard form of set theory:
to give a scheme means to firstly give a set of points; then
to describe a topology on this set; and finally to equip the resulting space
with a local sheaf of rings. Basic objects of study in algebraic geometry, such
as closed subschemes of projective spaces, are in this way encoded using a large
amount of machinery.

There is also a somewhat lesser used, but philosophically rewarding and more
``economical'' approach within set theory: Grothendieck's functorial approach.
In this account of scheme theory, to give a scheme means to give a functor
from the category of commutative rings to the category of sets. For instance,
the Fermat scheme is given by the functor
\[ A \longmapsto \{ (x,y,z) \in A^3 \,|\, x^n + y^n - z^n = 0 \}, \]
that is, by a \emph{scheme} in the colloquial sense for prescribing a set of
solutions for any ring.

This approach requires fewer preparations and involves
only objects of intrinsic interest to algebraic geometry: $A$-valued points,
where~$A$ ranges over all rings. These tend to be better behaved, for instance in
that the set of~$A$-valued points of a product of schemes is isomorphic to the
product of the sets of~$A$-valued points, and are more fundamental from a
geometric point of view. In contrast, the set-theoretical points of a scheme in
the approach using locally ringed spaces actually parameterize irreducible
closed subsets, not points in an intuitive sense.

Canonical references for the functorial approach are lectures notes by
Gro\-then\-dieck~\cite{grothendieck:functorial-ag} and the book by Demazure and
Gabriel~\cite{demazure:gabriel}. A summary in English, including a proof of the
equivalence with the approach using locally ringed spaces, is contained in the
first chapter of~\cite{vezzani:fun}. At the Secret Blogging Seminar, there was
an insightful long-running discussion on the merits of the functorial
approach~\cite{secret-blogging-seminar:fpov}, and further philosophical background
is contained in~\cite{mclarty:ontology}. The thesis of Zhen Lin Low~\cite{low:local-models}
contains recent developments on an abstract theory of gluing local models.

The description of basic objects can still be somewhat involved in the
functorial approach. For instance, while the functor associated to
projective~$n$-space is given on fields by the simple expression
\[ \begin{array}{r@{}c@{}l}
  K &{}\longmapsto{}& \text{the set of lines through the origin in~$K^{n+1}$} \\
  && \qquad \cong \{ [x_0 \hg \cdots \hg x_n] \,|\, \text{$x_i \neq 0$ for some~$i$} \},
\end{array} \]
on general rings it is given by
\[ \begin{array}{r@{}c@{}l}
  A &{}\longmapsto{}& \text{the set of quotients~$A^{n+1} \twoheadrightarrow P$,
  where~$P$ is projective of rank~1,} \\
  && \qquad \text{modulo isomorphism}.
\end{array} \]

On the one hand, typically only field-valued points admit a simple description.
On the other hand, the $A$-valued points for more general rings~$A$ are crucial
in order to impart a meaningful sense of cohesion on the field-valued
points. They therefore can't simply be dropped.\footnote{For instance,
let~$\ul{\AA}^1 : A \mapsto A$ be the functor associated to the affine line.
The Yoneda lemma guarantees that the set of morphisms~$\ul{\AA}^1 \to \ul{\AA}^1$
in the functor category~$[\Ring,\Set]$ is in canonical bijection with the
set~$\ZZ[U]$, as one would expect: Algebraic functions~$\AA^1 \to \AA^1$ should
be given by polynomials. (The discussion could also be relativized so that the
answer is the polynomial ring~$k[U]$, where~$k$ is some base field.) However, if we compute the
set of morphisms in~$[\mathrm{Field},\Set]$ we obtain~$\int_{K \in
\mathrm{Field}} \Hom(K,K)$, a set which contains pathological functions such as
some which permute the elements of the prime fields in arbitrary ways.}

We can resolve the tension by incorporating an automatic management of the
\emph{stage of definition}, the rings~$A$ such that we're
considering~$A$-valued points, into our language. Such a language is provided
by the internal language of the big Zariski topos. It allows for the Fermat
scheme to be given by the naive expression
\[ \{ (x,y,z) \? (\ul{\AA}^1)^3 \,|\, x^n + y^n - z^n = 0 \} \]
and for projective~$n$-space to be given by either of the expressions
\begin{multline*}
  \qquad\text{the set of lines through the origin in~$(\ul{\AA}^1)^{n+1}$}
  \quad\text{or} \\
  \{ [x_0 \hg \cdots \hg x_n] \,|\, \text{$x_i \neq 0$ for some~$i$} \}.\qquad
\end{multline*}
This is not a specialized trick to give short descriptions of some schemes:
Like with the internal universe of any topos, the full power of intuitionistic
logic is available to reason about the objects constructed in this way.

We can thus add an approach to the list of ways of giving a rigorous foundation
to algebraic geometry, the synthetic approach which layers scheme theory not
upon a classical set theory, but rather directly encodes schemes as sets and
morphisms of schemes as maps of sets in the nonclassical universe provided by
the big Zariski topos of a base scheme. We can therefore use a simple,
element-based language to talk about schemes.

This is similar to synthetic approaches to other fields of mathematics, such as
differential geometry~\cite{kock:sdg}, domain
theory~\cite{hyland:synthetic-domain-theory}, computability
theory~\cite{bauer:synthetic-computability-theory}, and more recently and very
successfully homotopy theory~\cite{hott} and related
subjects~\cite{schreiber:cohesive,schreiber-shulman:qgft,riehl-shulman:synthetic-infinity-categories}.
The synthetic approaches allow in each case to encode the
objects of study directly as (nonclassical) sets, with geometric,
domain-theoretic, computability-theoretic, or homotopic structure being
automatically provided for.

The implicit algebro-geometric structure has visible consequences on the
internal universe of the big Zariski topos and endows it with a distinctive
algebraic flavor. For instance, the statement
``\emph{any} map~$\ul{\AA}^1 \to \ul{\AA}^1$ is a polynomial function''
holds from the internal point of view. This is also a property which sets the
internal universe of the big Zariski topos apart from the toposes studied in
synthetic differential geometry.

If one is content with building upon classical scheme theory, the big Zariski
topos~$\Zar(X)$ of a base scheme~$X$ can be constructed as the topos of
sheaves on the Grothendieck site~$\Sch/X$ of~$X$-schemes.\footnote{Some care is
needed in order to avoid set-theoretical issues of size. We discuss this fine
point in Section~\ref{sect:proper-choice-of-site}. If one is interested in
foundational questions and doesn't merely want to use the big Zariski topos in
order to employ its convenient internal language, one can rest assured that
there's a way to construct it without resorting to classical scheme theory. We
sketch this in Section~\ref{sect:big-zariski-without-classical-scheme-theory}.}
Explicitly, an object of~$\Zar(X)$ is a functor~$F : (\Sch/X)^\op \to \Set$
satisfying the gluing condition with respect to Zariski coverings:
If~$T = \bigcup_i U_i$ is a cover of an~$X$-scheme~$T$ by open subsets, the
diagram
\[ F(T) \longrightarrow \prod_i F(U_i) \xbigtoto{} \prod_{j,k} F(U_j \cap U_k) \]
should be an equalizer diagram. A premier example of an object of~$\Zar(X)$ is
the functor~$\ul{Y}$ of points associated to an~$X$-scheme~$Y$, mapping
an~$X$-scheme~$T$ to~$\Hom_X(T, Y)$. It satisfies the gluing condition since
one can glue morphisms of schemes in the Zariski topology.

The object~$\ul{\AA}^1$ which already appeared is the functor of points
of the affine line over~$X$, the~$X$-scheme~$\AA^1_X \defeq X \times_{\Spec\ZZ}
\ZZ[U]$. Its value on an~$X$-scheme~$T$ is
\[ \afflx(T) = \Hom_X(T, \AA^1_X) \cong \Hom_{\Spec\ZZ}(T, \Spec\ZZ[U]) \cong
\Gamma(T,\O_T). \]
This object has a canonical structure as a ring object in~$\Zar(X)$. In fact,
from the internal point of view of~$\Zar(X)$, it is a local ring and even a
field in the sense that nonzero elements are invertible. In the case~$X =
\Spec\ZZ$, this fact was first observed by Kock~\cite{kock:univ-proj-geometry}. At
the same time, it is not a reduced ring -- a feat possible only in an
intuitionistic context. This curious interplay is quite important, since the
sets
\[ \{ x\?\afflx \,|\, x = 0 \} \quad\text{and}\quad
  \{ x\?\afflx \,|\, x^2 = 0 \} \]
should and do describe two different~$X$-schemes: the first is isomorphic
to~$X$ while the second is an infinitesimal thickening of~$X$, the vanishing scheme
of~$U^2$ in~$\AA^1_X$. In contrast, the sets~$\{ x\?\afflx \,|\, x \neq 0 \}$
and~$\{ x\?\afflx \,|\, x^2 \neq 0 \}$ should and do coincide. By the field
property, both conditions are equivalent to~$x$ being invertible.

The synthetic spectrum of an~$\afflx$-algebra~$A$ can be defined as
\[ \Spec(A) \defeq
  \text{the set of~$\afflx$-algebra homomorphisms~$A \to \afflx$}. \]
On first sight, this definition seems to overlook potential non-maximal prime
ideals of~$A$, since it only gives the~$\afflx$-valued points. But in fact,
this description correctly reflects the relative spectrum construction.
It yields a simple correspondence between synthetic affine schemes and solution
sets of polynomial equations. For instance, it's easy to show that there is a
canonical isomorphism
\begin{multline*}
  \qquad\Spec(\afflx[U_1,\ldots,U_n]/(f_1,\ldots,f_m)) \cong \\
  \{ (u_1,\ldots,u_n) \? (\afflx)^n \,|\,
    f_1(u_1,\ldots,u_n) = \cdots = f_m(u_1,\ldots,u_n) = 0 \}.\qquad
\end{multline*}
We give internal descriptions of further constructions of relative scheme
theory in Section~\ref{sect:basic-constructions}.

In order to be able to \emph{reason} internally (in contrast to only using the
internal language to describe~$X$-schemes and more general spaces in a simple
language), it's crucial to have strong and meaningful axioms available. One
such axiom posits that~$\afflx$ is a local and \emph{synthetically
quasicoherent} ring and implies all known ring-theoretic properties
of~$\afflx$. Synthetic quasicoherence is the internal analogue of the usual
condition on a sheaf of modules to be quasicoherent. This notion doesn't have a
counterpart in synthetic differential geometry and is central to our account of
synthetic algebraic geometry, since we derive all of its basic concepts such as
open and closed immersions and synthetic schemes from it.

Modal operators are useful in the big topos setting as well. For instance,
there is a modal operator~$\Box_\text{ét}$ in the big Zariski topos such that
the internal statement~$\Box_\text{ét} \varphi$ roughly means that~$\varphi$
holds on an étale covering and such that the translated
formula~$\varphi^{\Box_\text{ét}}$ means that~$\varphi$ holds in the \emph{big étale
topos} familiar from étale cohomology. In this way, we can access the internal
universe of the big étale topos from within the big Zariski topos. The
ring~$\afflx$ enjoys additional properties when studied in the étale topos,
where it is separably closed, in the fppf topos, where it is \emph{fppf-local},
and in the ph~topos, where it satisfies a strong form of algebraic
closure.

\subsection*{Limitations} The internal language is \emph{local}, in the sense
that if~$X = \bigcup_i U_i$ is an open covering and an internal statement
holds in the sheaf toposes~$\Sh(U_i)$, it holds in~$\Sh(X)$ as well. On the one
hand, this property is very useful. But on the other hand, it causes an inherent
limitation of the internal language:
Global properties of sheaves of modules like ``being generated by global
sections'', ``being ample'', or ``having vanishing sheaf cohomology'' and global properties of schemes like ``being
quasicompact'' can \emph{not} be
expressed in the internal language.

Thus for global considerations, the internal language of~$\Sh(X)$ is only
useful in that local subparts can be simplified. Also, some global features
reflect themselves in certain metaproperties of the internal language. For
instance, a scheme is quasicompact if and only if the internal language
has a weak version of the so-called disjunction property of mathematical
logic (Section~\ref{sect:compactness}).

The locality limitation only refers to locality with respect to the base
scheme. For instance, the little and big Zariski toposes of~$X$ \emph{can}
distinguish between affine and projective~$n$-space over~$X$, even though these
are locally isomorphic.

The internal languages of both toposes can be used on a case-by-case
basis, employing them as part of longer arguments in the context of ordinary scheme
theory where it's useful to do so. However, if one wants to stay solely in one
of the provided internal universes and not use ordinary scheme theory at all,
then one will of course run into the further limitation that internal scheme
theory, as put forward in this thesis, is only developed to a small extent.

\subsection*{Introductory literature} This text is intended
to be self-contained, requiring only basic knowledge of scheme theory. In
particular, we assume no prior familiarity with topos theory or formal logic.

Nevertheless, a gentle
introduction to topos theory is an article by
Leinster~\cite{leinster:introduction}. Standard references for the internal
language of toposes include the book of Mac~Lane and
Moerdijk~\cite[Chapter~VI]{moerdijk-maclane:sheaves-logic}, the book of
Goldblatt~\cite[Chapter~14]{goldblatt:topoi}, Caramello's and Streicher's lecture
notes~\cite{caramello:preliminaries,streicher:ctcl}, the book of
Borceux~\cite[Chapter~6]{borceux:handbook3}, and Part~D of
Johnstone's Elephant~\cite{johnstone:elephant}. Motivation and background on
the internal language can also be found in Chapter~0 of Shulman's lecture
notes~\cite{shulman:categorical-logic}.

In the 1970s, there was a
flurry of activity on applications of the internal language. An article by
Mulvey~\cite{mulvey:repr} of this time gives a very accessible
introduction to the topic, culminating in an internal proof of the Serre--Swan
theorem (with just one external ingredient needed).

\subsection*{Related work} The internal language of toposes was applied to algebraic geometry before. For
instance, Wraith used it to construct (and verify the universal property
of) the little étale topos of a scheme by internally developing the theory of
strict henselization~\cite{wraith:generic-galois-theory}. However, to the best
of our knowledge, systematically building a dictionary relating external and
internal notions has not been attempted before, and the use of modal operators
to study the spreading of properties from points to neighborhoods seems to be
new as well.

In particular, Tierney remarked in 1976 that a certain property of the internal
universe of the little Zariski topos ``is surely important, though its precise
significance is still somewhat obscure''~\cite[p.~209]{tierney:spectrum}. This
property can now be recognized as a small shadow of the internal
characterization of quasicoherence. We expand on this in
Section~\ref{sect:field-properties}.

In some regards, this thesis is an extended answer to a MathOverflow question
by Gubkin~\cite{mo:gubkin}, who inquired about uses of the internal language of
toposes in algebraic geometry.

Brandenburg put forward a related program of internalization in his PhD
thesis~\cite{brandenburg:tensor-foundations}. However, he internalizes
constructions of algebraic geometry not in toposes, but in tensor categories.
There is some overlap in working out precise universal properties, particularly
when dealing with the big Zariski topos.

In other branches of mathematics, the internal language of toposes is used as well. For
instance, there is an ongoing effort in mathematical physics to understand
quantum mechanical systems from an internal point of view: To any quantum
mechanical system, one can associate a so-called Bohr topos containing an
internal mirror image of the system. This mirror image looks like a
system of classical mechanics from the internal perspective, and therefore
tools like Gelfand duality can be used to construct an internal
phase space for the system~\cite{bohr1,bohr2,bohr3}.

In stochastics, the usefulness of an internal language was recently stressed by
Tao~\cite{tao:analysis-rel-measure-space}. Such a language makes the
common notational practice of dropping the explicit dependence of the
value~$X(\omega)$ of a random variable on the sample~$\omega$ completely
rigorous and simplifies the basic theory. Tao also highlighted how a suitable
language can be used to simplify ``$\varepsilon$/$\delta$ management'' in
analysis~\cite{tao:cheap-nsa}. Furthermore, there is a topos-theoretic approach to
measure theory, in which the sheaf of measurable real functions on
a~$\sigma$-algebra looks like the ordinary set of real numbers from an internal point
of view~\cite{jackson:sheaf-theoretic-measure-theory}; this has applications in
noncommutative geometry~\cite{henry:measure-theory-boolean-toposes}.

Intuitionistic methods have found many applications in computer science.
Recently, the internal language of a topos of trees and a suitable modal
operator was used to study guarded recursion, encompassing, for instance, an
internal Banach fixed-point theorem~\cite{birkedal:al:sgdt}.

In constructive mathematics, the internal language of toposes is routinely used
to obtain models of intuitionistic theories fulfilling certain anti-classical
axioms. For instance, there are toposes in which the axiom ``any map~$\RR \to
\RR$ is continuous'' (appropriately formulated) holds~\cite{kock:sdg,moerdijk:reyes:models}
and toposes in which the Church--Turing thesis ``any map~$\NN \to \NN$ is
computable'' holds (certain realizability toposes).
The internal language can also be used to extract computational content
out of classical constructions. To cite just one recent example, Mannaa and
Coquand used it to implement algorithms for working with the algebraic closure
of an arbitrary field of characteristic zero~\cite{mannaa:coquand:alg-closure}.

One way this thesis contributes to the program of constructive
mathematics is that intuitionistic mathematics gains new areas of application.
For instance, the constructive account of the theory of Krull dimension was
originally developed to remove Noetherian hypotheses, extract computational meaning, and
simplify proofs~\cite{dyn:krull-integral,dyn:char-krull}. It can now also be used to
reason about the dimension of schemes, since the topological dimension of a
scheme~$X$ coincides with the Krull dimension of the structure sheaf~$\O_X$
regarded as an ordinary ring from the internal perspective of~$\Sh(X)$
(Section~\ref{sect:krull-dimension}).

We obtained a second contribution to constructive mathematics as a byproduct of
deducing transfer principles which relate a module over a ring~$A$ with its
induced quasicoherent sheaf on~$\Spec A$: Using the internal language of the
little Zariski topos we can algorithmically turn certain non-constructive
arguments concerning prime ideals into constructive ones. We discuss this in
Section~\ref{sect:eliminating-prime-ideals}; it is related to the
\emph{dynamical methods in algebra} explored by Coquand, Coste, Lombardi, Roy,
and others~\cite{clr:dynamicalmethod,cl:logical}.

Caramello uses topos theory to build bridges between different mathematical
subjects, in a certain precise sense~\cite{caramello:1,caramello:2}. She
exploits that toposes can admit presentations by sites of different character. Our
contribution is certainly related to her grand research program in spirit, but since we
focus only on specific presentations of a few specific toposes associated to
schemes, there are as yet only few direct technical connections.

\subsection*{Notational conventions} To stress that a discussion takes place in
an intuitionistic context, we occasionally write~``$\forall x\?X$''
or~``$\exists x\?X$'' instead of~``$\forall x \in X$'' and~``$\exists x \in
X$'' not only in internal statements, where it's proper to do so, but also when
not reasoning internally.

If~$X$ and~$Y$ are sets, we mean by~``$[X,Y]$'' the set of all maps from~$X$
to~$Y$. This expression will often occur in internal formulas; its external
meaning will then be the Hom sheaf. We write pairs~$\langle a, b \rangle$ using
angle brackets. The preimage of a set~$M$ under a map~$f$ is
written~``$f^{-1}[M]$''. Similarly, the image is written~``$f[N]$''.

The constant sheaf with stalks~$M$ is written~``$\ul{M}$''.

}

\section{The internal language of a sheaf topos}\label{sect:internal-language}

At its heart, the internal language of a topos provides a coherent way of
translating any mentions of set-theoretical elements to
\emph{generalized elements}, carefully keeping track of and adapting
the stage of definition. We want to illustrate this with a simple example
before giving the formal definition.

A map~$f : X \to Y$ of sets is injective if and only if
\begin{equation}\label{inj}
  \forall x,x' \in X\_ f(x) = f(x') \Longrightarrow x = x'.
\end{equation}
This condition can not only be interpreted in~$\Set$, but in any category~$\C$ whose
objects are structured sets and whose morphisms are maps between the underlying
sets. If we want to go beyond such kind of categories, we have to restate the
condition in purely category-theoretic language:
\begin{equation}\label{injcat}
  \forall (1 \xra{x} X), (1 \xra{x'} X)\_ f \circ x = f \circ x'
  \Longrightarrow x = x'.
\end{equation}
This condition makes sense in all categories which contain a terminal
object~$1$, and is equivalent to condition~\eqref{inj} in the case~$\C = \Set$.
This has a deeper reason: The one-element set~$1 = \{ \star \}$ is a
\emph{separator} of~$\Set$, that is objects of~$\Set$ are uniquely determined
by their \emph{global elements}, morphisms from the terminal object.

However, in categories in which the terminal object is not a separator,
condition~\eqref{injcat} is not very meaningful. This is for instance the case
if~$\C$ is the category of vector spaces over a field or if~$\C$ is the category~$\Sh(X)$ of set-valued sheaves on a topological
space~$X$. Global elements of a sheaf~$\F$ are in natural one-to-one
correspondence with global sections~$s \in \F(X)$ (hence the name), whereby
condition~\eqref{injcat} only states that~$f$ is \emph{injective on global
sections}. Since many interesting sheaves admit no or only few global sections,
this statement is typically not very substantial.

A basic tenet of category theory is therefore to not only refer to global
elements~$1 \to X$, but also to \emph{generalized elements}~$A \to X$,
where~$A$ ranges over all objects. The domain~$A$ is called the \emph{stage of
definition} in this context. Bearing this principle in mind, a better
translation of the injectivity condition is the statement
\begin{equation}\label{injgen}
  \forall \text{objects $A$ in~$\C$}\_ \forall (A \xra{x} X), (A \xra{x'\!} X) \text{ in~$\C$}\_\
  f \circ x = f \circ x' \ \Longrightarrow\ x = x'.
\end{equation}
This statement expresses that~$f$ is a monomorphism and therefore
correctly captures the structural essence of injectivity.

Unlike this manual translation guided by trial and error and categorical
philosophy, the internal language provides a purely mechanical translation
scheme.  It is fully formal, can be analyzed rigorously, works smoothly with
arbitrarily convoluted statements, and most importantly can be trusted to
support \emph{reasoning}: If a statement formulated in a naive element-based
language intuitionistically implies a further such statement, then the
translation of the former implies the translation of the latter.

The power of the internal language doesn't unfold in basic situations like with
the example above, where one can easily translate statements and even proofs by
hand. It unfolds when considering more complex statements. For instance, the
short proof of Grothendieck's generic freeness lemma promised in the
introduction rests on the internal statement~``any ideal
of~$\O_{\Spec(R)}[U_1,\ldots,U_n]$ is \notnot finitely generated'',
where~$R$ is a reduced ring. For the proof of Grothendieck's generic freeness
lemma it's not necessary to actually perform the translation of this statement into external
language, but for definiteness we display the translation here nevertheless:
\begin{indentblock}\label{page:convoluted-statement}
For any element~$f \in R$ and any (not necessarily quasicoherent) sheaf of
ideals~$\J \hookrightarrow \O_{\Spec(R)}[U_1,\ldots,U_n]|_{D(f)}$: If
\begin{indentblock}
for any element~$g \in R$ the condition that
\begin{indentblock}
the sheaf~$\J$ is of finite type on~$D(g)$
\end{indentblock}
implies that~$g = 0$,
\end{indentblock}
then~$f = 0$.
\end{indentblock}
This statement is obviously quite convoluted, and its proof is even more so;
therefore it probably wouldn't occur to one to base a proof of Grothendieck's
generic freeness lemma on this statement. The internal language is thus of real
use here. We'll expand on this example in Section~\ref{sect:noetherian} and in
Section~\ref{sect:generic-freeness}.\footnote{The statement can be proven by
hand, but it's much simpler to only verify the case~$n = 0$ (and even reduce
this case to simple other properties which~$\O_{\Spec(R)}$ enjoys from the
internal point of view) and then to apply Hilbert's basis theorem. Hilbert's
basis theorem is famous for admitting only a nonconstructive proof, and
nonconstructive proofs can't be translated by the internal language machinery;
but this is only true for the conclusion ``any ideal is finitely generated''.
The intuitionistically weaker conclusion ``any ideal is \notnot finitely
generated'' does admit a constructive proof, and is all what's needed here.}

\subsection{Internal statements}
Let~$X$ be a topological space. Later, $X$ will be the underlying space of a
scheme. The meaning of internal statements is given by a set of rules, the
\emph{Kripke--Joyal semantics} of the topos of sheaves on~$X$.

\begin{defn}\label{defn:kripke-joyal}The meaning of
\[ U \models \varphi \quad\text{(``$\varphi$ holds on $U$'')} \]
for open subsets~$U \subseteq X$ and formulas~$\varphi$ over~$U$ is given by
the rules listed in Table~\ref{table:kripke-joyal}, recursively in the structure of~$\varphi$.
In a \emph{formula over~$U$} there may appear sheaves defined on~$U$ as domains
of quantifications,~$U$-sections of sheaves as terms, and morphisms of sheaves
on~$U$ as function symbols. If~$V \subseteq U$ is an open subset, then formulas
over~$U$ can be pulled back to formulas over~$V$. The symbols~``$\top$'' and~``$\bot$'' denote truth
and falsehood, respectively. The universal and existential quantifiers come in
two flavors: for bounded and unbounded quantification.
The translation of~$U \models \neg\varphi$ does not have to be separately defined, since
negation can be expressed using other symbols: $\neg\varphi \defequiv (\varphi
\Rightarrow \bot)$. If we want to emphasize the particular topos, we write
\[ \Sh(X) \models \varphi \quad\Ll\quad X \models \varphi. \]
\end{defn}

\begin{table}
  \centering
  \[ \renewcommand{\arraystretch}{1.3}\begin{array}{@{}lcl@{}}
    U \models s = t \? \F &\Ll& s|_U = t|_U \in \Gamma(U, \F) \\
    U \models s \in \G &\Ll& s|_U \in \Gamma(U,\G) \quad\quad\text{($\G$ a
    subsheaf of~$\F$, $s$ a section of~$\F$)} \\
    U \models \top &\Ll& U = U \text{ (always fulfilled)} \\
    U \models \bot &\Ll& U = \emptyset \\
    U \models \varphi \wedge \psi &\Ll&
      \text{$U \models \varphi$ and $U \models \psi$} \\
    U \models \bigwedge_{j \in J} \varphi_j &\Ll&
      \text{for all~$j \in J$: $U \models \varphi_j$} \quad\quad\text{($J$ an
      index set)} \\
    U \models \varphi \vee \psi &\Ll&
      \hcancel{\text{$U \models \varphi$ or $U \models \psi$}}{0pt}{3pt}{0pt}{-2pt} \\
    && \text{there exists a covering $U = \bigcup_i U_i$ such that for all~$i$:} \\
    && \quad\quad \text{$U_i \models \varphi$ or $U_i \models \psi$} \\
    U \models \bigvee_{j \in J} \varphi_j &\Ll&
      \hcancel{\text{$U \models \varphi_j$ for some~$j \in J$}}{0pt}{3pt}{0pt}{-2pt}
      \quad\quad\text{($J$ an index set)} \\
    && \text{there exists a covering $U = \bigcup_i U_i$ such that for all~$i$:} \\
    && \quad\quad \text{$U_i \models \varphi_j$ for some~$j \in J$} \\
    U \models \varphi \Rightarrow \psi &\Ll&
      \hcancel{\text{$U \models \varphi$ implies $U \models \varphi$}}{0pt}{3pt}{0pt}{-2pt} \\
    && \text{for all open~$V \subseteq U$:
      $V \models \varphi$ implies $V \models \psi$} \\
    U \models \forall s \? \F\_ \varphi(s) &\Ll&
      \text{for all sections~$s \in \Gamma(V, \F)$ on open $V \subseteq U$: $V \models
      \varphi(s)$} \\
    U \models \exists s \? \F\_ \varphi(s) &\Ll&
      \hcancel{\text{there exists a section~$s \in \Gamma(U,\F)$ such that $U
      \models \varphi(s)$}}{0pt}{3pt}{0pt}{-2pt} \\
    &&
      \text{there exists an open covering $U = \bigcup_i U_i$ such that for all~$i$:} \\
    && \quad\quad \text{there exists~$s_i \in \Gamma(U_i, \F)$ such that
    $U_i \models \varphi(s_i)$} \\
    U \models \forall \F\_ \varphi(\F) &\Ll&
      \text{for all sheaves $\F$ on open $V \subseteq U$: $V \models \varphi(\F)$} \\
    U \models \exists \F\_ \varphi(\F) &\Ll&
      \text{there exists an open covering $U = \bigcup_i U_i$ such that for all~$i$:} \\
    && \quad\quad \text{there exists a sheaf~$\F_i$ on~$U_i$ such that
    $U_i \models \varphi(\F_i)$}
  \end{array} \]
  \caption{\label{table:kripke-joyal}The Kripke--Joyal semantics of a sheaf
  topos.}
\end{table}

\begin{rem}The last two rules in Table~\ref{table:kripke-joyal}, concerning
\emph{unbounded quantification}, are not part of the classical Kripke--Joyal
semantics. They are part of Mike Shulman's stack semantics~\cite{shulman:stack},
a slight but important extension. They are needed so that we can formulate universal
properties in the internal language. (Prior work in the same direction include
the topos models explored by Pitts~\cite[Section~3]{pitts:polymorphism} and, in
the context of set theory, work by Awodey, Butz, Simpson, and
Streicher~\cite{awodey-butz-simpson-streicher:bist}, which was carried out
independently and published after Shulman's paper.)
\end{rem}

\begin{ex}\label{ex:injective-surjective}
Let~$\alpha : \F \to \G$ be a morphism of sheaves on~$X$. Then
$\alpha$ is a monomorphism of sheaves if and only if, from the internal
perspective,~$\alpha$ is simply an injective map:
\allowdisplaybreaks
\begin{align*}
  & X \models \speak{$\alpha$ is injective} \\[0.5em]
  \Longleftrightarrow\
  & X \models \forall s\?\F\_ \forall t\?\F\_ \alpha(s) = \alpha(t) \Rightarrow s = t \\[0.5em]
  \Longleftrightarrow\ &
    \text{for all open~$U \subseteq X$, sections $s \in \Gamma(U, \F)$:} \\
  &\qquad\qquad \text{for all open~$V \subseteq U$, sections $t \in \Gamma(V, \F)$:} \\
  &\qquad\qquad\qquad\qquad
      V \models \alpha(s) = \alpha(t) \Rightarrow s = t \\[0.5em]
  \Longleftrightarrow\ &
    \text{for all open~$U \subseteq X$, sections $s \in \Gamma(U, \F)$:} \\
  &\qquad\qquad \text{for all open~$V \subseteq U$, sections $t \in \Gamma(V, \F)$:} \\
  &\qquad\qquad\qquad\qquad
      \text{for all open~$W \subseteq V$:} \\
  &\qquad\qquad\qquad\qquad\qquad\qquad
        \text{$\alpha_W(s|_W) = \alpha_W(t|_W)$ implies $s|_W = t|_W$} \\[0.5em]
  \Longleftrightarrow\ &
    \text{for all open~$U \subseteq X$, sections $s, t \in \Gamma(U, \F)$:} \\
  &\qquad\qquad
        \text{$\alpha_U(s|_U) = \alpha_U(t|_U)$ implies $s|_U = t|_U$} \\[0.5em]
  \Longleftrightarrow\ &
    \text{$\alpha$ is a monomorphism of sheaves}
\end{align*}
The corner quotes ``$\speak{\ldots}$'' indicate that translation into formal
language is left to the reader. Similarly, the morphism~$\alpha$ is an epimorphism of
sheaves if and only if, from the internal perspective,~$\alpha$ is a
surjective map. Notice that injectivity and surjectivity are
notions of a simple element-based language. The Kripke--Joyal semantics
takes care to properly handle \emph{all} sections, not only global ones.
\end{ex}

The rules are not all arbitrary. They are finely concerted to make the
following two propositions true, which are crucial for a proper appreciation of the
internal language.

\begin{prop}[Locality of the internal language]
\label{prop:locality-of-the-internal-language}
Let~$U = \bigcup_i U_i$ be covered by open subsets. Let~$\varphi$
be a formula over~$U$. Then
\[ U \models \varphi \qquad\text{iff}\qquad
  \text{$U_i \models \varphi$ for each $i$}. \]
\end{prop}
\begin{proof}Induction on the structure of~$\varphi$. The canceled
rules in Table~\ref{table:kripke-joyal} would make this proposition false.\end{proof}

As a corollary, one may restrict the open coverings and universal
quantifications in the definition of the Kripke--Joyal semantics
(Table~\ref{table:kripke-joyal}) to open subsets of some basis of the topology.
For instance, if~$X$ is a scheme, one may restrict to affine open subsets.

Furthermore, Proposition~\ref{prop:locality-of-the-internal-language} shows that the internal language is monotone in
the following sense: If~$U \models \varphi$, and~$V$ is an open subset of~$U$,
then~$V \models \varphi$. (This follows by applying the proposition to the
trivial covering~$U = V \cup U$.)

\begin{prop}[Soundness of the internal language]
\label{prop:soundness-of-the-internal-language}
If a formula~$\varphi$ implies a further formula~$\psi$ in intuitionistic logic, then
$U \models \varphi$ implies $U \models \psi$.
\end{prop}
\begin{proof}
Proof by induction on the structure of formal intuitionistic proofs; we are to
show that any inference rule of intuitionistic logic is satisfied by the
Kripke--Joyal semantics. For instance, there is the following rule governing
disjunction:
\begin{quote}
If~$\varphi \vee \psi$ holds, and both $\varphi$ and $\psi$ imply a further
formula~$\chi$, then~$\chi$ holds.
\end{quote}
So we are to prove that if~$U \models \varphi \vee \psi$, $U \models (\varphi
\Rightarrow \chi)$, and $U \models (\psi \Rightarrow \chi)$, then $U \models \chi$.
This is done as follows: By assumption, there exists a covering~$U = \bigcup_i
U_i$ such that for each index~$i$, $U_i \models \varphi$ or $U_i \models \psi$.
Again by assumption, we may conclude that~$U_i \models \chi$ for each~$i$. The statement
follows because of the locality of the internal language.

A complete list of which rules are to prove is
in Appendix~\ref{appendix:inference-rules}.
\end{proof}

In particular, if a formula~$\psi$ has an unconditional intuitionistic proof,
then~$U \models \psi$.

The restriction to intuitionistic logic is really necessary at this point. We
will encounter many examples of classically equivalent internal statements whose
translations using the Kripke--Joyal semantics are wildly different. To
anticipate just one example, the statement
\[ X \models \speak{$\F$ is finite free}, \]
referring to a sheaf~$\F$ of~$\O_X$-modules, means that~$\F$ is finite locally
free. The statement
\[ X \models \neg\neg(\speak{$\F$ is finite free}) \]
instead means that~$\F$ is finite locally free on a dense open subset of~$X$.

In
particular, our treatment of modal operators to understand spreading of
properties from points to neighborhoods depends on having the ability to make
finer distinctions -- distinctions which are not visible in classical logic.
In Section~\ref{sect:appreciating-intuitionistic-logic} there is a discussion of what the restriction to
intuitionistic logic amounts to in practice.

Because of the multitude of quantifiers, literal translations of internal statements
can sometimes get slightly unwieldy. There are simplification rules for certain
often-occurring special cases:
\begin{prop}\label{prop:simplification}
    \[ \renewcommand{\arraystretch}{1.3}\begin{array}{@{}lcl@{}}
      U \models \forall s\?\F\_ \forall t\?\G\_ \varphi(s,t)
      &\Longleftrightarrow&
      \text{for all open~$V \subseteq U$,} \\
      && \text{sections~$s \in \Gamma(V,\F)$, $t \in \Gamma(V,\G)$:
      $V \models \varphi(s,t)$} \\[0.3em]
      U \models \forall s\?\F\_ \varphi(s) \Rightarrow \psi(s)
      &\Longleftrightarrow&
      \text{for all open~$V \subseteq U$, sections~$s \in \Gamma(V,\F)$:} \\
      &&\qquad\qquad \text{$V \models \varphi(s)$ implies $V \models \psi(s)$}
      \\[0.3em]
      U \models \exists!s\?\F\_ \varphi(s)
      &\Longleftrightarrow&
      \text{for all open~$V \subseteq U$,} \\
      &&
      \text{there is exactly one section~$s \in \Gamma(V,\F)$ with:} \\
      &&\qquad\qquad V \models \varphi(s)
    \end{array} \]
\end{prop}
\begin{proof}Straightforward. By way of example, we prove the existence claim
in the ``only if'' direction of the last rule. (This rule formalizes
the saying ``unique existence implies global existence''.) By definition of~$\exists!$, it
holds that
\[ U \models \exists s\?\F\_ \varphi(s)
  \qquad\text{and}\qquad
  U \models \forall s,t\?\F\_ \varphi(s) \wedge \varphi(t) \Rightarrow s = t. \]
Let~$V \subseteq U$ be an arbitrary open subset. Then there exist local
sections~$s_i \in \Gamma(V_i,\F)$ such that~$V_i \models \varphi(s_i)$, where~$V
= \bigcup_i V_i$ is an open covering. By the locality of the internal language,
on intersections it holds that~$V_i \cap V_j \models \varphi(s_i)$, so by the
uniqueness assumption, it follows that the local sections agree on intersections.
They therefore glue to a section~$s \in \Gamma(V,\F)$. Since~$V_i \models
\varphi(s)$ for all~$i$, the locality of the internal language allows us to
conclude that~$V \models \varphi(s)$.
\end{proof}

\begin{rem}Note that~$\Sh(X) \models \neg\varphi$ is in general a much stronger
statement than merely saying that~$\Sh(X) \models \varphi$ does not hold:
The former always implies the latter (unless~$X = \emptyset$, in which case
\emph{any} internal statement is true), but the converse does not hold: The
former statement means that~$U = \emptyset$ is the \emph{only} open subset on
which~$\varphi$ holds, that is that~$\varphi$ holds \emph{nowhere}. In
contrast, the statement~$\Sh(X) \not\models \varphi$ only means that~$\varphi$
does \emph{not hold everywhere}.

For instance, let~$X$ be a scheme and let~$f \in
\Gamma(X,\O_X)$ be a global function on~$X$. We will see in
Section~\ref{sect:sheaves-of-rings} that ``$\Sh(X) \not\models \speak{$f$ is
invertible}$'' means that~$f$ is not an invertible element in the
ring~$\Gamma(X,\O_X)$. In contrast,~``$\Sh(X) \models \neg(\speak{$f$ is
invertible})$'' means that~$f$ is \emph{nowhere invertible}, not even on
smaller non\-empty open subsets. This implies that~$f$ is nilpotent.
\end{rem}

It's instructive to consider the special case that~$X$ is the one-point space.
In this case~$\Sh(X) \simeq \Set$, by taking constant sheaves.
Let~$\varphi$ be a formula in which arbitrary sets and elements may occur as
parameters. Then we can regard~$\varphi$ as a formula over the open set~$X$ by
substituting any occuring set by the induced constant sheaf and any occuring
element by the induced global section, whereby it's meaningful to
write~``$\Set \models \varphi$''. One can check that
\[ \Set \models \varphi \qquad\text{iff}\qquad \text{$\varphi$ in the usual
mathematical sense}. \]
(A proof is presented in Lemma~\ref{lemma:properties-of-constant-sheaves}.)
Hence the internal language of the topos~$\Set$ is just the ordinary
mathematical language.\footnote{Readers familiar with Tarski's theorem on
undefinability of truth will recognize that we're sweeping a subtle issue under
the rug, as it's not actually possible, within the confines of a chosen formal
system as foundation of mathematics, to formally define what it means for a
formula to ``hold in the usual mathematical sense''. This issue can be solved
in a number of ways, for instance by performing the Kripke--Joyal translation
at the meta level (at the price of not having arbitrary set-indexed
conjunctions and disjunctions available in the internal language, at least if
the meta level is taken to be something like~$\mathrm{PRA}$) or by choosing a
universe~$U$ and replacing~$\Set$ by the full subcategory of those sets
contained in~$U$. The latter option can even be carried out without an increase
in consistency strength, by employing a system
like Feferman's~$\mathrm{ZFC}/\mathrm{S}$~\cite{shulman:set-theory,feferman:set-foundations}.}

\subsection{Internal constructions}
\label{sect:internal-constructions}
The Kripke--Joyal semantics defines the
interpretation of internal \emph{statements}. The interpretation of internal
\emph{constructions} is given by the following definition.

\begin{defn}\label{defn:interpretation-internal-constructions}
The interpretation of an internal construction~$T$
is denoted by~$\brak{T} \in \Sh(X)$ and is given by the following rules.
\begin{itemize}\item If~$\F$ and~$\G$ are sheaves, $\brak{\F \times \G}$ is the
categorical product of~$\F$ and~$\G$ (\ie their product as presheaves).
\item If~$\F$ and~$\G$ are sheaves, $\brak{\F \amalg \G}$ is the categorical
coproduct of~$\F$ and~$\G$, \ie the sheafification of the presheaf
$U \mapsto \Gamma(U,\F) \amalg \Gamma(U,\G)$.
\item If~$\F$ is a sheaf, the interpretation~$\brak{\P(\F)}$ of the power set
construction is the sheaf given by
\[ \text{$U \subseteq X$ open} \quad\longmapsto\quad \{ \G \hookrightarrow \F|_U \}, \]
\ie sections on an open set~$U$ are subsheaves of~$\F|_U$ (either literally
or isomorphism classes of arbitrary monomorphisms into~$\F|_U$).
\item If~$\F$ is a sheaf and~$\varphi(s)$ is a formula containing a free
variable~$s\?\F$, the interpretation~$\brak{\{s\?\F\,|\,\varphi(s)\}}$ is given
by the subpresheaf of~$\F$ defined by
\[ \text{$U \subseteq X$ open} \quad\longmapsto\quad \{ s \in \Gamma(U,\F) \ |\
  U \models \varphi(s) \}. \]
By the locality of the internal language, this presheaf is in fact a
sheaf.
\end{itemize}
\end{defn}

The definition is made in such a way that, from the internal perspective, the
constructions enjoy their expected properties. For instance, it holds that
\[ \Sh(X) \models
  \bigl(\forall x\?\brak{\{s\?\F \,|\, \varphi(s)\}}\_ \psi(x)\bigr)
  \Longleftrightarrow
  \bigl(\forall x\?\F\_ \varphi(x) \Rightarrow \psi(x)\bigr). \]
We gloss over several details here. See~\cite[Section~D4.1]{johnstone:elephant} for
a proper treatment.

Morphisms can internally be constructed by appealing to the \emph{principle of
unique choice}: Let~$\varphi(s,t)$ be a formula with free variables of
type~$s\?\F$, $t\?\G$. Assume
\[ \Sh(X) \models \forall s\?\F\_ \exists!t\?\G\_ \varphi(s,t). \]
Then there is one and only one morphism~$\alpha : \F \to \G$ of sheaves such
that for any local section~$s \in \Gamma(U,\F)$, $\Sh(X) \models
\varphi(s,\alpha(s))$. This follows from the meaning of unique existence with
the Kripke--Joyal semantics (Proposition~\ref{prop:simplification}).

An important application is showing that two sheaves~$\F$ and~$\G$ are
isomorphic (usually as objects with more structure, for instance sheaves of
modules). To this end, it suffices to give a formula~$\varphi(s,t)$ satisfying,
in addition to the condition above, the condition
$\Sh(X) \models \forall t\?\G\_ \exists! s\?\F\_ \varphi(s,t)$,
expressing that the induced morphism~$\alpha$ is a bijective map from the
internal perspective. This implies the statement
\[ \Sh(X) \models \exists \alpha\?\HOM(\F,\G)\_ \speak{$\alpha$ is bijective},
\]
but this statement is strictly weaker: Its interpretation with the
Kripke--Joyal semantics is that the sheaves~$\F$ and~$\G$ are \emph{locally}
isomorphic.

\subsection{Geometric formulas and
constructions}\label{sect:geometric-formulas-and-constructions}
In formal and categorical logic so-called geometric formulas play a
special role. They are named that way because, in a sense which can be made
precise, their meaning is preserved under pullback along geometric morphisms.
\begin{defn}\label{defn:geometric-formulas}
A formula is \emph{geometric} if and only if it consists only of
\[ {=} \quad {\in} \quad {\top} \quad {\bot} \quad {\wedge} \quad {\vee} \quad
{\bigvee} \quad {\exists}, \]
but not~``$\bigwedge$'' nor ``$\Rightarrow$'' nor~``$\forall$'' (and thus
not~``$\neg$'' either, since negation is defined using~``$\Rightarrow$'').
A \emph{geometric implication} is a formula of the form
\[ \forall \cdots \forall\_ (\cdots) \Rightarrow (\cdots) \]
with the bracketed subformulas being geometric.
\end{defn}
The \emph{parameters} of a formula~$\varphi$ are the sheaves
being quantified over, sections of sheaves appearing as terms, and morphisms of
sheaves appearing as function symbols in~$\varphi$.
We say that a formula~$\varphi$ holds \emph{at a point~$x \in X$} if and only
if the formula obtained by substituting all parameters in~$\varphi$ with their
stalks at~$x$ holds in the usual mathematical sense.

\begin{lemma}\label{lemma:geometric-stalk-neighborhood}
Let~$x \in X$ be a point. Let~$\varphi$ be a geometric formula (over some open
neighborhood~V of~$x$).
Then~$\varphi$ holds at~$x$ if and only if there exists an open neighborhood~$U
\subseteq X$ of~$x$ (contained in~V) such that~$\varphi$ holds on~$U$.
\end{lemma}
\begin{proof}This is a very general instance of the phenomenon that sometimes,
truth at a point spreads to truth on a neighborhood. It can be proven by
induction on the structure of~$\varphi$, but we will give a more conceptual
proof later (Corollary~\ref{cor:geometric-spreading}).
\end{proof}

This lemma is a very useful metatheorem. We will properly discuss its
significance in Section~\ref{sect:spreading}. For now, we just use it to prove a
simple criterion for the internal truth of a geometric implication; we will
apply this criterion many times.

\begin{cor}\label{cor:geometric-implication}
A geometric implication holds on~$X$ if and only if it holds at
every point of~$X$.\end{cor}
\begin{proof}For notational simplicity, we consider a geometric implication of
the form
\[ \forall s\?\F\_ \varphi(s) \Rightarrow \psi(s). \]
For the ``only if'' direction, assume that this formula holds on~$X$ and let~$x
\in X$ be an arbitrary point. Let~$s_x \in \F_x$ be the germ of an arbitrary
local section~$s$ of~$\F$ and assume that~$\varphi(s)$ holds at~$x$. By
Lemma~\ref{lemma:geometric-stalk-neighborhood}, it follows that~$\varphi(s)$ holds on some open neighborhood of~$x$. By
assumption,~$\psi(s)$ holds on this neighborhood as well. Again by the
lemma,~$\psi(s)$ holds at~$x$.

For the ``if'' direction, assume that the geometric implication holds at every
point. Let~$U \subseteq X$ be an arbitrary open subset and let~$s \in
\Gamma(U,\F)$ be a local section such that~$\varphi(s)$ holds on~$U$. By the
lemma and the locality of the internal language, to show that~$\psi(s)$ holds
on~$U$, it suffices to show that~$\psi(s)$
holds at every point of~$U$. This is clear, since again by the
lemma,~$\varphi(s)$ holds at every point of~$U$.
\end{proof}

\begin{ex}Injectivity and surjectivity are geometric implications (surjectivity
can be spelled~$\forall y\?\G\_ (\top \Rightarrow \exists x\?\F\_ \alpha(x) =
y)$). Thus Corollary~\ref{cor:geometric-implication} gives a deeper reason for the well-known fact that a
morphism of sheaves is a monomorphism \resp an epimorphism if and only if it is
stalkwise injective \resp surjective.\end{ex}

A construction is \emph{geometric}\label{page:geometric-constructions} if and only if it commutes with pullback
along arbitrary geometric morphisms. We do not want to discuss the notion of
geometric morphisms here; suffice it to say that calculating the stalk at a
point~$x \in X$ is an instance of such a pullback. Among others, the following
constructions are geometric:
\begin{itemize}
\item finite product: $(\F \times \G)_x \cong \F_x \times \G_x$
\item finite coproduct: $(\F \amalg \G)_x \cong \F_x \amalg \G_x$
\item arbitrary coproduct: $(\coprod_i \F_i)_x \cong \coprod_i (\F_i)_x$
\item set comprehension with respect to a \emph{geometric} formula~$\varphi$:
\[ \brak{\{ s\?\F \,|\, \varphi(s) \}}_x \cong \{ [s]\in\F_x \,|\,
\text{$\varphi(s)$ holds at $x$} \} \]
\item free module: $(\R\langle \F \rangle)_x \cong \R_x\langle \F_x
\rangle$ ($\R$ a sheaf of rings, $\F$ a sheaf of sets)
\item localization of a module: $\F[\S^{-1}]_x \cong \F_x[\S_x^{-1}]$
\end{itemize}
Compatibility with taking stalks is not sufficient for geometricity.
It is just the most easily visualized requirement.
The following constructions are not in general geometric:
\begin{itemize}
\item arbitrary product
\item set comprehension with respect to a non-geometric formula
\item powerset
\item internal Hom: $\HOM(\F,\G)_x \not\cong \Hom(\F_x,\G_x)$
\end{itemize}

\subsection{Appreciating intuitionistic logic}
\label{sect:appreciating-intuitionistic-logic}
The principal (and only) difference between classical and intuitionistic logic
is that in classical logic, the axioms schemes of \emph{excluded middle} and
\emph{double negation elimination} are added.
\[ \varphi \vee \neg\varphi \qquad\qquad \neg\neg\varphi \Rightarrow \varphi \]
A classically trained mathematician might legitimately wonder why one should
drop these axioms: Are they not obviously true? The pragmatic answer to this
question is that the translations of these axioms with the Kripke--Joyal
semantics are, except for uninteresting special cases of the base space~$X$,
plainly false -- irrespective of one's philosophical convictions. Therefore the
internal language is in general only sound with respect to intuitionistic logic and
not with respect to classical logic. Concretely, there is the following
proposition.
\begin{prop}The internal language of
a~$\mathrm{T}_1$-space~$X$ is \emph{Boolean}, \ie it verifies the classical
axiom schemes displayed above, if and only if~$X$ is discrete.
The internal language of an irreducible or locally Noetherian scheme~$X$ is Boolean if and only if~$X$ has
dimension~$\leq 0$.
\end{prop}
\begin{proof}\label{prop:lang-boolean}
The internal language of~$\Sh(X)$ is Boolean if and only if for
any open subset~$U \subseteq X$ it holds that~$U$ is the only dense open subset
of~$U$. This can be checked manually, by using the definition of the
Kripke--Joyal semantics, but we'll be able to give a more conceptual proof
later (Lemma~\ref{lemma:boolean-dense}). The first claim is then an exercise in
point-set topology, while the second is more difficult
(Corollary~\ref{cor:boolean-dim0}).
\end{proof}

However, there is also a more satisfying answer, which furthermore
illuminates how to intuitively picture intuitionistic mathematics.
Namely, when doing intuitionistic mathematics, we use the same formal symbols as classically, but with
\emph{a different intended meaning}. For instance, the classical reading of an
existential statement like~$\exists x\?A\_ \varphi(x)$ is that there exists
some element~$x \? A$ with the property~$\varphi(x)$. In contrast, its
intuitionistic reading is that such an element can actually be
\emph{constructed}, \ie explicitly given in some form. This is a much stronger
statement. Classically, a proof that it is \emph{not} the case that such an
element does \emph{not} exist -- formally $\neg\neg \exists x\?A\_ \varphi(x)$
(or, equivalently even in intuitionistic mathematics, $\neg\forall x\? A\_
\neg\varphi(x)$) -- suffices to demonstrate the existential statement; this is not
so in intuitionistic mathematics.

Similarly, the intuitionistic meaning of a disjunction~$\varphi \vee \psi$ is
not only that one of the disjuncts is true, but that one can explicitly state
which case holds. It is in general not enough to show that it is impossible
that both~$\varphi$ and~$\psi$ fail.

In this picture, it is obvious that one should not adopt the law of excluded
middle or the principle of double negation elimination as axioms. But we
don't \emph{reject} those axioms in the sense of postulating their
converses, we simply don't use them. Therefore any intuitionistically
true result is also true classically. In fact, for some special instances,
these two classical axioms do hold intuitionistically. For example, any natural
number is zero or is not zero -- this is not a triviality, but can be proven by
induction.\footnote{The analogous statement about real numbers cannot be
shown. Intuitively, for a number given by a decimal expansion starting
with~$0.0000\ldots$ one cannot decide whether the string of zeros will continue
indefinitely or whether eventually a non-zero digit will occur. This argument
can be made rigorous. The analogous statement about algebraic numbers
\emph{can} be proven; the information contained in a witness of algebraicity (a
monic polynomial which the given number is a zero of) suffices to make the
case distinction~\cite[Chapter~VI.1, p.~140]{mines-richman-ruitenburg:constructive-algebra}.}

A consequence of not adopting these axioms is that proofs by contradiction are
not generally justified; they are intuitionistically valid only for those
statements which can be proven to be true or false. A proof of a
\emph{negated formula} is not the same as a proof by contradiction. For
instance, the usual proof that~$\sqrt{2}$ is not rational is
intuitionistically perfectly fine: From the assumption that~$\sqrt{2}$ is
rational one deduces a contradiction~($\bot$). This is exactly the definition
of~$\neg(\speak{$\sqrt{2}$ is rational})$.

A more positive consequence of not adopting the law of excluded middle and the
principle of double negation elimination is that intuitionistically, we can
make \emph{finer distinctions}. For instance, for a formula~$\varphi$, the doubly
negated formula~$\neg\neg\varphi$ (``\notnot~$\varphi$'') is a certain kind of weakening of~$\varphi$:
If~$\varphi$ holds, then~$\neg\neg\varphi$ does as well, while the converse can
not be shown in general.\footnote{A detailed proof of the correct implication
goes as follows: Assume~$\varphi$. We are to show~$\neg\neg\varphi$, \ie
$(\neg\varphi \Rightarrow \bot)$. So assume~$\neg\varphi$, we are to
show~$\bot$. Since~$\varphi$ and~$\varphi \Rightarrow \bot$,~$\bot$ indeed
follows.} An example from everyday life runs as follows: If in the morning you
can't find the key for your apartment, but you know that it must hide
somewhere since you used it to open the door in the evening before, you
intuitionistically know~$(\neg\neg\exists x\_ \speak{the key is at position
$x$})$, but you cannot claim the unnegated proposition. One cannot model this
distinction with pure classical logic.

Double negation also has a concrete geometric meaning with the
Kripke--Joyal semantics. Namely,~$X \models \neg\neg\varphi$ holds if and
only if there is a dense open subset~$U$ of~$X$ such that~$U \models \varphi$.
This is of course a weaker statement than~$X \models \varphi$.
In Section~\ref{sect:modalities}, we will discuss this fact and other
\emph{modal operators} in more detail. For instance, there is a similarly defined modal
operator~$\Box$ such that~$X \models \Box\varphi$ if and only if there is an
open neighborhood~$U$ of a given point~$x$ such that~$U \models \varphi$. Also
there is a different operator~$\Box$ such that~$X \models \Box\varphi$ if and only
if~$\varphi$ holds on a scheme-theoretically dense open subset.

For future reference, we remark that if~$\varphi \Rightarrow \psi$,
then also~$\neg\neg\varphi \Rightarrow \neg\neg\psi$; that weakening
twice has no further effect, \ie~$\neg\neg\neg\neg\varphi \Leftrightarrow
\neg\neg\varphi$;\footnote{In fact, negating thrice is the same as negating
once: Assume~$\neg\neg\neg\varphi$. We are to show~$\neg\varphi$. So
assume~$\varphi$, we are to show~$\bot$. Since~$\varphi$,~$\neg\neg\varphi$.
By~$\neg\neg\neg\varphi$,~$\bot$ follows.} and that the double negation of the
law of excluded middle,~$\neg\neg(\varphi \vee \neg\varphi)$, holds.

A classical mathematician might then ask which classical results are valid
intuitionistically. Firstly, in linear and commutative algebra, most
of the basic theorems stay valid, provided one exercises some caution in
formulating them (for instance, one should not arbitrarily weaken assumptions
by introducing double negations). This is because the proofs of these
statements are usually direct; if intuitionistically unacceptable case
distinctions do occur, they can often be eliminated by streamlining the proof.

Consider as a simple example the proposition that the kernel of a linear map is
a linear subspace. The case distinction ``either the kernel consists just of the
zero vector, in which case the claim is trivial, or otherwise \ldots'' is not
intuitionistically acceptable, but can be entirely dispensed with: The proof
for the general case works in the special case just as well.

Secondly, there is \emph{Barr's theorem}. This metatheorem states that, if a
geometric implication has a proof using classical logic and the axiom of choice
(from certain axioms which too can be formulated as geometric implications), then
it also has a proof using intuitionistic logic. For instance, the
Nullstellensatz (any finite system of polynomial equations over an
algebraically closed geometric field either possesses a solution or a
certificate that there can't be any solutions) is typically proven using
maximal ideals and therefore using some forms of the axiom of choice. By Barr's
theorem, there is also an intuitionistically valid proof.

Barr's theorem only requires restrictions on the form of the statement and the
axioms; the given classical proof can be of any form whatsoever and can freely
use statements which are not geometric implications. The proof of Barr's
theorem itself is not intuitionistically valid, which puts us in the curious
situation that we nonconstructively know that there exists a constructive
proof. Because I strive, as far as possible, that these notes can be
interpreted in a constructive metatheory, and because I personally prefer
direct intuitionistic proofs, we won't use Barr's theorem in what follows.
However, it's useful to quickly assess whether a classically known proposition
also has an intuitionistic proof.

Finally, we should clarify the status of the axiom of choice. This axiom, which
is strictly speaking not part of classical logic, but of a classical set
theory, is not accepted in an intuitionistic context: By \emph{Diaconescu's
theorem}, it implies the law of excluded middle in presence of the other axioms
of set theory.

Standard references for intuitionistic algebra are a textbook by Mines,
Richman and
Ruitenburg~\cite{mines-richman-ruitenburg:constructive-algebra} and a textbook
by Lombardi~\cite{lombardi:quitte:constructive-algebra}. The standard
reference for intuitionistic analysis is a book by Bishop and
Bridges~\cite{bishop-bridges:constructive-analysis}. Further explanations and
pointers to relevant literature can be found in an expository article and a
recorded lecture by Bauer~\cite{bauer:int-mathematics,bauer:video}. A
recent survey of intuitionistic logic from a historical and logical point of
view is~\cite{melikhov:intuitionistic-logic}.

\begin{rem}For much of this text, we work in a classical metatheory. This
means that we allow ourselves to occasionally use the law of excluded middle
and the axiom of choice when reasoning \emph{about} the internal language.
In particular, we have the theory of schemes as commonly presented at our
disposal. This decision has two reasons.

Firstly, we want to connect the internal world with the usual external
framework of algebraic geometry, in order to be directly useful to working
algebraic geometers who work in a classical metatheory. We want to prove
statements like ``a scheme~$X$ as classically defined has this-or-that property
if and only if, from the internal point of view of~$\Sh(X)$, this-or-that
holds''.

Secondly, as of yet, there is no full constructive account of the theory of
schemes with which we could establish a link with the internal language. We
sketch how such an account could be developed, and also why one might want to
do that, in Section~\ref{sect:constructive-scheme-theory}.\end{rem}

\subsection{A fine point on internal natural numbers}

The internal world of~$\Sh(X)$ contains an object which behaves like the set of
natural numbers. For instance, it's possible to prove statements about internal
natural numbers by induction and to construct functions on the internal natural
numbers by recursion. Externally, this object is the constant sheaf~$\ul{\NN}$
of locally constant~$\NN$-valued functions on~$X$.

On the other hand, we can utilize that the internal language supports infinite
conjunctions and disjunctions; therefore we may include expressions
like~``$\bigvee_{n \in \NN}$'', where~$n$ ranges over all external natural
numbers, in internal formulas.

The two approaches are related as follows. Any external natural number~$n \in
\NN$ gives rise to a global section~$\ul{n}$ of the sheaf~$\ul{\NN}$, thus to
an internal natural number. If~$\varphi(n)$ is a formula of the
internal language depending on a parameter~$n \? \ul{\NN}$, then
\[ \Sh(X) \models \exists n\?\ul{\NN}\_ \varphi(n)
  \quad\text{if and only if}\quad
  \Sh(X) \models \bigvee_{n \in \NN} \varphi(\ul{n}), \]
and similarly for~``$\forall$'' and~``$\bigwedge$''. In practice we can
therefore often ignore the subtle difference. A proof of the equivalence rests
on the observation
\[ \Sh(X) \models \forall n\?\ul{\NN}\_
  \bigvee_{m \in \NN} n = \ul{m}, \]
which can be checked by translating this statement using the Kripke--Joyal
semantics.

A similar relation holds for internal polynomials. If~$\R$ is a sheaf of rings
over~$X$, then we can construct, internally in~$\Sh(X)$, the ring of
polynomials over~$\R$. This will yield a certain sheaf~$\brak{\R[T]}$.
If~$\varphi(f)$ is a formula of the internal language containing a free
variable~$f \? \R[T]$, then
\begin{multline*}
  \quad \Sh(X) \models \exists f\?\R[T]\_ \varphi(f)
  \quad\text{if and only if} \\
  \Sh(X) \models \bigvee_{n \in \NN}
    \exists a_0\?\R\_ \ldots \exists a_n\?\R\_
    \varphi(a_nT^n + \cdots + a_1T + a_0). \quad
\end{multline*}

\chapter{The little Zariski topos}\label{part:little-zariski}

\section{Sheaves of rings}
\label{sect:sheaves-of-rings}

Recall that a \emph{sheaf of rings} can be categorically described as a
sheaf of sets~$\R$ together with maps of sheaves $+, \cdot : \R \times \R \to
\R$, $- : \R \to \R$, and global elements~$0, 1$ such that certain axioms hold.
For instance, the axiom on the commutativity of addition is rendered in
diagrammatic form as follows:
\[ \xymatrix{
  \R \times \R \ar[rr]^{\mathrm{swap}} \ar[rd]_{+} && \R \times \R \ar[ld]^{+} \\
  & \R
} \]

From the internal perspective, a sheaf of rings looks just like a plain ring.
This is the content of the following proposition.

\begin{prop}\label{prop:rings-internally}
Let~$X$ be a topological space. Let~$\R$ be a sheaf of sets on~$X$.
Let~$+, \cdot : \R \times \R \to \R$ and $- : \R \to \R$ be maps of sheaves and let~$0, 1$ be
global elements of~$\R$. Then these data define a sheaf of rings if and only
if, from the internal perspective, these data fulfill the usual equational ring
axioms.\end{prop}
\begin{proof}We only discuss the commutativity axiom. The internal statement
\[ \Sh(X) \models \forall x,y\?\R\_ x + y = y + x \]
means that for any open subset~$U \subseteq X$ and any local sections~$x,y \in
\Gamma(U,\R)$, it holds that~$x + y = y + x \in \Gamma(U,\R)$. This is
precisely the external commutativity condition.
\end{proof}

\begin{lemma}\label{lemma:internal-invertibility}
Let~$X$ be a topological space. Let~$\R$ be a sheaf of rings
on~$X$. Let~$f$ be a global section of~$\R$. Then the following statements are
equivalent:
\begin{enumerate}
\item $f$ is invertible from the internal point of view, \ie $\Sh(X) \models
\exists g\?\R\_ fg = 1$.
\item $f$ is invertible in all stalks~$\R_x$.
\item $f$ is invertible in~$\Gamma(X,\R)$.
\end{enumerate}
\end{lemma}
\begin{proof}Since invertibility is a geometric implication, the equivalence of
the first two statements is clear. Also, it is obvious that the third statement
implies the other two. For the remaining direction, note that the
uniqueness of inverses in rings can be proven intuitionistically. Therefore, if~$f$ is invertible
from the internal point of view, it actually holds that
\[ \Sh(X) \models \exists! g\?\R\_ fg = 1. \]
Since unique internal existence implies global existence
(Proposition~\ref{prop:simplification}), this shows that the first statement
implies the third.
\end{proof}

\subsection{Reducedness}\label{sect:reducedness} Recall that a scheme~$X$ is \emph{reduced} if and only
if all stalks~$\O_{X,x}$ are reduced rings. Since the condition on a ring~$R$
to be reduced is a geometric implication,
\[ \forall s\?R\_ \Bigl(\bigvee_{n \geq 0} s^n = 0\Bigr) \Longrightarrow s = 0, \]
we immediately obtain the following characterization of reducedness in the
internal language:
\begin{prop}\label{prop:reduced-ring}
A scheme~$X$ is reduced iff, from the internal point of view, the
ring~$\O_X$ is reduced.\end{prop}

\subsection{Locality} Recall the usual definition of a local ring: a ring
possessing exactly one maximal ideal. This is a so-called \emph{higher-order
condition} since it involves quantification over subsets. It is also not of a
geometric form. Therefore, for our purposes, it is better to
adopt the following elementary definition of a local ring.
\begin{defn}\label{defn:local-ring}A \emph{local ring} is a ring~$R$ such that~$1 \neq 0$ in~$R$ and
for all~$x,y \? R$
\[ \text{$x+y$ invertible} \quad\Longrightarrow\quad
  \text{$x$ invertible}\ \vee\ \text{$y$ invertible}. \]
\end{defn}
In classical logic, it is an easy exercise to show that this definition is
equivalent to the usual one. In intuitionistic logic, we would need to be
more precise in order to even state the question of equivalence, since
intuitionistically, the notion of a maximal ideal bifurcates into several
non-equivalent notions.\footnote{For instance, should a maximal ideal~$\mmm$ be
such that if~$\nnn$ is any ideal with~$\mmm \subseteq \nnn \subsetneq (1)$,
then~$\mmm = \nnn$? Or should the condition be that if~$\nnn$ is any ideal
with~$\mmm \subseteq \nnn$, then~$\mmm = \nnn$ or~$\nnn = (1)$?
Intuitionistically, the latter condition is stronger than the former.}
This is a common phenomenon in intuitionistic
mathematics: Classically equivalent notions may bifurcate into related but
inequivalent notions intuitionistically, each having a unique character and
yielding slightly different theories.

\begin{prop}\label{prop:local-ring}
In the internal language of a scheme~$X$ (or a locally ringed
space), the ring~$\O_X$ is a local ring.\end{prop}
\begin{proof}The condition in Definition~\ref{defn:local-ring} is a conjunction of two geometric
implications (the first one being~$1 = 0 \Rightarrow \bot$, the second being
the displayed one) and holds on each stalk.\end{proof}

\begin{rem}When first exposed to locally ringed spaces, one might ask why the
requirement is that the \emph{stalks}~$\O_{X,x}$ are local rings, instead of the
easier-to-define sets of sections~$\O_X(U)$. This question has of course a good
geometric answer. Using the internal language, it also has a purely formal
answer: The requirement that the stalks are local rings is precisely the
requirement that the ring~$\O_X$ is a local ring from the perspective of the
internal language of~$X$.
\end{rem}

\subsection{Field properties}\label{sect:field-properties}
From the internal point of view, the structure
sheaf~$\O_X$ of a scheme~$X$ is \emph{almost} a field, in the sense that any
element which is not invertible is nilpotent. This is a genuine property of
schemes, not shared with arbitrary locally ringed spaces. It is also a specific
feature of the internal universe: Neither the local rings~$\O_{X,x}$ nor the
rings of local sections~$\Gamma(U,\O_X)$ have this property in general.

\begin{prop}\label{prop:neginvnilpotent}Let~$X$ be a scheme. Then
\[ \Sh(X) \models \forall s\?\O_X\_ \neg(\speak{$s$ invertible}) \Rightarrow
\speak{$s$ nilpotent}. \]
\end{prop}
\begin{proof}By the locality of the internal language and since~$X$ can be
covered by open affine subsets, it is enough to show that for any affine
scheme~$X = \Spec A$ and any global function~$s \in \Gamma(X,\O_X) = A$ it holds
that
\[ X \models \neg(\speak{$s$ invertible}) \quad\text{implies}\quad
  X \models \speak{$s$ nilpotent}. \]
The meaning of the antecedent is that any open subset on which~$s$ is
invertible is empty. This implies in particular that the standard open subset~$D(s)$ is
empty. This means that~$s$ is an element of any prime ideal of~$A$, thus
nilpotent, and therefore implies the a priori weaker statement~$X \models \speak{$s$
nilpotent}$ (which would allow~$s$ to have different indices of nilpotency on
an open covering).
\end{proof}

\begin{rem}In classical logic, the statement ``not invertible implies
nilpotent'' is equivalent to ``any element is invertible or nilpotent''.
However, in intuitionistic logic, the latter is strictly stronger than the
former. We will see in the next section
(Corollary~\ref{cor:scheme-dimension-zero}) that the structure sheaf of a
scheme fulfills the latter condition if and only if the scheme is
zero-dimensional (or empty). An overview of the basic properties of the
intuitionistically different field conditions is contained
in~\cite{johnstone:rings-fields-and-spectra}.
\end{rem}

\begin{cor}\label{cor:field-reduced}
Let~$X$ be a scheme. If~$X$ is reduced, the ring~$\O_X$ is a field
from the internal point of view, in the sense that
\[ \Sh(X) \models \forall s\?\O_X\_ \neg(\speak{$s$ invertible}) \Rightarrow
s=0. \]
Conversely, if~$\O_X$ is a field in this internal sense, then~$X$ is reduced.\end{cor}
\begin{proof}We can prove this purely in the internal language: It suffices to
give an intuitionistic proof of the fact that a local ring which satisfies the
condition of the previous proposition fulfills the stated field condition if
and only if it is reduced. This is straightforward.
\end{proof}

This field property is very useful. We will put it to good use when giving a
simple proof of the fact that~$\O_X$-modules of finite type on a reduced scheme
are locally free on a dense open subset (Lemma~\ref{lemma:locally-free-dense}).
The field property only holds in the precise form as stated;
the classically equivalent condition that any element is invertible or zero is
intuitionistically stronger. This is an instance of the already remarked upon
phenomenon of intuitionistic bifurcation of notions.

The observation that the structure sheaf is (almost) a field is attributed by
Tierney to Mulvey~\cite[p.~209]{tierney:spectrum}.
Tierney also states that ``[it] is surely important, though its precise significance is still somewhat
obscure'' (ibid). We think that it's significant as a special case of the
following more general proposition,
which states that we can deduce a certain unconditional
statement from the premise that, under the assumption that some element~$f\?\O_X$ is invertible, an element~$s\?\O_X$ is zero. This is
interesting on its own, but will be of particular importance in understanding
quasicoherence from the internal point of view (Section~\ref{sect:qcoh}) and
interpreting the relative spectrum as an internal spectrum
(Section~\ref{sect:relative-spectrum}).

\begin{prop}\label{prop:cond-zero}
Let~$X$ be a scheme. Then
\[ \Sh(X) \models
  \forall f\?\O_X\_
  \forall s\?\O_X\_
  (\speak{$f$ \inv} \Rightarrow s = 0) \Longrightarrow
  \textstyle
  \bigvee_{n \geq 0} f^n s = 0. \]
\end{prop}
\begin{proof}It is enough to show that for any affine scheme~$X = \Spec A$ and
any global functions~$f, s \in A$ such that
\[ X \models (\speak{$f$ \inv} \Rightarrow s = 0), \]
it holds that $X \models \textstyle \bigvee_{n \geq 0} f^n s = 0$. This
indeed follows, since by assumption such a function~$s$ is zero on~$D(f)$, \ie $s$
is zero as an element of~$A[f^{-1}]$.
\end{proof}

Proposition~\ref{prop:neginvnilpotent} follows from this proposition by
setting~$s \defeq 1$.

\subsection{Krull dimension}\label{sect:krull-dimension}
Recall that the \emph{Krull dimension} of a
ring is usually defined as the supremum of the lengths of strictly
ascending chains of prime ideals. As with the classical definition of a local ring,
this definition does not lead to a well-behaved notion in an intuitionistic
context. Furthermore, it is a higher-order condition, so interpreting it
with the Kripke--Joyal semantics is a bit unwieldy.

Luckily, there is an elementary definition of the Krull dimension which works
intuitionistically and which is classically equivalent to the usual notion. It
was found by Coquand and Lombardi, building upon work by
Joyal and Español~\cite{dyn:krull-integral,dyn:char-krull}, and can be
used to give a short proof that~$\dim k[X_1,\ldots,X_n] = n$, where~$k$ is a
field~\cite{dyn:krull-dim-polynomial-ring}.

\begin{defn}Let~$R$ be a ring. A \emph{complementary sequence} for a
sequence~$(a_0,\ldots,a_n)$ of elements of~$R$ is a sequence~$(b_0,\ldots,b_n)$
such that the following inclusions of radical ideals hold:
\[ \renewcommand{\arraystretch}{1.3}
\left\{\begin{array}{@{}rcl@{}}
  \sqrt{(1)} &\subseteq& \sqrt{(a_0,b_0)} \\
  \sqrt{(a_0 b_0)} &\subseteq& \sqrt{(a_1,b_1)} \\
  \sqrt{(a_1 b_1)} &\subseteq& \sqrt{(a_2,b_2)} \\
  &\vdots \\
  \sqrt{(a_{n-1} b_{n-1})} &\subseteq& \sqrt{(a_n,b_n)} \\
  \sqrt{(a_n b_n)} &\subseteq& \sqrt{(0)}
\end{array}\right. \]
The ring~$R$ is \emph{of Krull dimension~$\leq n$} if
and only if for any sequence~$(a_0,\ldots,a_n)$ there exists a
complementary sequence. (The ring~$R$ is trivial if and only if it is
of Krull dimension~$\leq -1$.)
\end{defn}
Unlike the usual definition, this definition posits only a condition
on elements and not on ideals. It is thus of a simpler logical form.
(The radical ideals appear only for convenience. We will eliminate them in the
proof of Proposition~\ref{prop:dimension-scheme-ox}.)
Also note that we do not define the Krull dimension of a ring as some natural
number (this is intuitionistically not possible for general rings). Instead, we
only define what it means for the Krull dimension to be less than or equal to
a given natural number.

For the following, no intuition about the definition is needed; however, we
feel that some motivation might be of use. Recall that we can picture inclusions of
radical ideals geometrically by considering standard open subsets~$D(f) = \{
\ppp \in \Spec R \,|\, f \not\in \ppp \}$: The inclusion~$\sqrt{(f)} \subseteq
\sqrt{(g,h)}$ holds if and only if~$D(f) \subseteq D(g) \cup D(h)$, and
intersections are calculated by products, \ie~$D(f) \cap D(g) = D(fg)$.

The condition that~$(b_0,\ldots,b_n)$ is complementary to~$(a_0,\ldots,a_n)$
thus means that~$D(a_0)$ and~$D(b_0)$ cover all of~$\Spec R$; that their
intersection is covered by~$D(a_1)$ and~$D(b_1)$; that in turn their
intersection is covered by~$D(a_2)$ and~$D(b_2)$; \ldots; and that finally, the
intersection of~$D(a_n)$ and~$D(b_n)$ is empty.

For the special case~$n = 0$, the condition that~$R$ is of Krull
dimension~$\leq 0$ means that for any element~$a_0$ there exists an
element~$b_0$ such that~$D(a_0)$ and~$D(b_0)$ cover~$\Spec R$ and are disjoint.

The definition of the Krull dimension can be written in such a way as to mimic the
definition of the inductive Menger--Urysohn dimension of topological
spaces~\cite[Section~1]{dyn:krull-integral}. For an internal characterization
of the dimension of smooth manifolds, we refer the reader to a result of
Fourman~\cite[Section~3]{fourman:reals}.

\begin{thm}Let~$R$ be a ring.
\begin{enumerate}
\item In classical logic, the ring~$R$ is
of Krull dimension~$\leq n$ if and only if its Krull dimension
as usually defined using chains of prime ideals is less than or equal to~$n$.
\item If the ring~$R$ is
of Krull dimension~$\leq n$, the radical of any finitely generated ideal is
equal to the radical of some ideal which can be generated by~$n+1$ elements.
This holds intuitionistically, and there is an explicit algorithm for computing
the reduced set of generators from the given ones. (Kronecker's theorem)
\end{enumerate}
\end{thm}
\begin{proof}See~\cite[Theorem~1.2]{dyn:krull-integral} for the first
statement. The proof relies on the observation that~$\dim R \leq n$ if and only
if~$\dim R[S_x^{-1}] \leq n-1$ for all~$x \in R$, where~$S_x = x^\NN (1+xR)
\subseteq R$. We put the second statement only to demonstrate that the
definition of the Krull dimension is constructively sensible. It follows from
the identity~$\sqrt{(x,a_0,\ldots,a_n)} =
\sqrt{(a_0-xb_0,\ldots,a_n-xb_n)}$, where~$(b_0,\ldots,b_n)$ is a complementary
sequence for~$(a_0,\ldots,a_n)$.
\end{proof}

We can apply the constructive theory of Krull dimension to the structure
sheaf~$\O_X$ of a scheme~$X$ as follows. The condition that a
scheme~$X$ has dimension exactly~$n$ (in the usual sense using ascending chains
of closed irreducible subsets) is not local -- the dimension may vary on
an open cover; hence it isn't possible to characterize this condition in
the internal language. However, the condition that the dimension of~$X$ is less
than or equal to~$n$ \emph{is} local, thus there is hope that it can be
internalized. And indeed, this is the case.

\begin{prop}\label{prop:dimension-scheme-ox}
Let~$X$ be a scheme. Then:
\[ \dim X \leq n \quad\Longleftrightarrow\quad
  \Sh(X) \models \speak{$\O_X$ is of Krull dimension~$\leq n$}
  \]
\end{prop}
\begin{proof}

A condition of the form~``$\sqrt{(f)} \subseteq \sqrt{(g,h)}$''
like in the constructive definition of the Krull dimension is not a geometric
formula when taken on face value. However, it is equivalent to a geometric
condition, namely to
\[ \exists a,b\?\O_X\_ \bigvee_{m \geq 0} f^m = ag + bh. \]
Therefore the condition~$\speak{$\O_X$ is of Krull
dimension~$\leq n$}$ is (equivalent to) a geometric implication and thus holds
internally if and only if it holds at every point~$x \in X$. This in turn means that the
Krull dimension of any stalk~$\O_{X,x}$ is less than or equal to~$n$. This is
equivalent to the (Krull) dimension of~$X$ being less than or equal to~$n$.
\end{proof}

We will state and prove a generalization of this lemma about the dimension of closed
subschemes later, as Lemma~\ref{lemma:dim-closed-subscheme}.

If~$X$ is a reduced scheme, we have seen in Corollary~\ref{cor:field-reduced}
that~$\O_X$ is a field from the internal perspective, in the sense that
non-invertible elements are zero. But fields are well-known to be of Krull
dimension zero. Why is this not a contradiction to the proposition just proven?
Intuitionistically, the notion of a field bifurcates into several
non-equivalent notions:
\begin{enumerate}
\item ``Any element which is not invertible is zero.''
\item ``Any element which is not zero is invertible.''
\item ``Any element is either zero or invertible.''
\end{enumerate}
Only fields in the sense~(3) are automatically of Krull dimension zero.
Fields in the weaker senses can have higher Krull dimension, as exhibited by
the structure sheaf of reduced schemes with positive dimension.

With our conventions, a scheme~$X$ is of dimension~$\leq 0$ if and only if it
is empty or if it's nonempty and of dimension zero.
\begin{cor}\label{cor:scheme-dimension-zero}
Let~$X$ be a scheme. Then:
\[ \dim X \leq 0 \quad\Longleftrightarrow\quad
  \Sh(X) \models \forall s\?\O_X\_ \speak{$s$ \inv} \vee \speak{$s$ nilpotent}.
  \]
If furthermore~$X$ is reduced, this is further equivalent to~$\O_X$ being a
field in the strong sense that any element of~$\O_X$ is invertible or zero.
\end{cor}
\begin{proof}By Proposition~\ref{prop:dimension-scheme-ox} and the fact that~$\O_X$ is a local ring from
the internal perspective, this is an immediate consequence of
interpreting the following standard fact of ring theory in the internal
language of~$\Sh(X)$: A local ring~$R$ is of Krull
dimension~$\leq 0$ if and only if any element of~$R$ is invertible or
nilpotent.

It is well-known that this holds classically; to make sure that it
holds intuitionistically as well (so that it can be used in the internal
universe), we give a proof of the ``only if'' direction. Let~$a \? R$ be
arbitrary. By assumption on the Krull dimension, there exists an element~$b \?
R$ such that~$\sqrt{(1)} \subseteq \sqrt{(a,b)}$ and~$\sqrt{(ab)} =
\sqrt{(0)}$. The latter means that~$ab$ is nilpotent. Since~$R$ is local, the
former implies that~$a$ is invertible or that~$b$ is invertible. In the first
case, we are done. In the second case, it follows that~$a$ is nilpotent, so we
are done as well.
\end{proof}

As a further corollary we note the curious fact that the classicality of the
internal language of~$\Sh(X)$, where~$X$ is a scheme, is tightly coupled with
the properties of the ring~$\O_X$: Internally, the law of excluded middle and
the principle of double negation elimination are ``almost equivalent'' to the
Krull dimension of~$\O_X$ being~$\leq 0$.
\begin{cor}\label{cor:boolean-dim0}
Let~$X$ be a scheme. If the internal language of~$\Sh(X)$ is Boolean, then
$\dim X \leq 0$. The converse holds if~$X$ is irreducible or locally Noetherian.
\end{cor}
\begin{proof}
We show that any element of~$\O_X$ is invertible or nilpotent, therefore
verifying the hypothesis of the previous corollary.
Let~$s\?\O_X$ be given. By assumption, either~$s$ is invertible or~$s$ is not
invertible. In the latter case~$s$ is nilpotent by
Proposition~\ref{prop:neginvnilpotent}.

We defer the converse direction to
Proposition~\ref{prop:boolean-dim0-continued} since we don't want to interrupt
the exposition here with a certain necessary technical condition.
\end{proof}

\subsection{Integrality}\label{sect:integrality}
In intuitionistic logic, the notion of an integral
domain bifurcates into several inequivalent notions. The following two are
important for our purposes:
\begin{defn}\label{defn:integral-domain}
A ring~$R$ is an \emph{integral domain in the weak sense} if and
only if~$1 \neq 0$ in~$R$ and
\[ \forall x,y\?R\_ xy = 0 \Longrightarrow (x = 0) \vee (y = 0). \]
A ring~$R$ is an \emph{integral domain in the strong sense} if and only if~$1
\neq 0$ in~$R$ and
\[ \forall x\?R\_ x = 0 \vee \speak{$x$ is regular}, \]
where~$\speak{$x$ is regular}$ means that~$xy = 0$ implies~$y = 0$ for any~$y \?
R$.\end{defn}

For the following result, recall that a scheme~$X$ (or a ringed space) is
\emph{integral at a point~$x \in X$} if and only if~$\O_{X,x}$ is an integral
domain (in either sense, since we have adopted a classical metatheory).

\begin{prop}\label{prop:internal-integrality}
Let~$X$ be a ringed space. Then:
\begin{enumerate}
\item $X$ is integral at all points if and only if, internally,~$\O_X$ is an
integral domain in the weak sense.
\item If~$X$ is even a locally Noetherian scheme, then~$\O_X$ is an integral
domain in the weak sense iff it is an integral domain in the strong sense from
the internal point of view.
\end{enumerate}
\end{prop}
\begin{proof}The condition on a ring to be an integral domain in the weak sense
is a conjunction of two geometric implications,~``$1 = 0 \Rightarrow \bot$''
and the implication displayed in the definition. Therefore the first statement
is obvious.

For the second statement, we observe the condition on a function~$f \in
\Gamma(U,\O_X)$ to be regular from the internal perspective is open: It holds
at a point~$x \in U$ if and only if it holds on some open neighborhood of~$x$.
We will give a proof of this specific feature of locally Noetherian schemes
later on, when we have developed appropriate machinery to do so easily
(Proposition~\ref{prop:regularity-spreading}). In any case, this openness
property was the essential ingredient for the equivalence between ``holding
internally'' and ``holding at every point''
(Corollary~\ref{cor:geometric-implication}). Therefore~$\O_X$ is an integral
domain in the strong sense from the internal point of view if and only if all
local rings~$\O_{X,x}$ are integral domains. By the first statement, this is
equivalent to~$\O_X$ being an integral domain in the weak sense from the
internal point of view.
\end{proof}

We record the following lemma for later use. The proof presented here is
already simple, but a more conceptual proof is also possible (see
Section~\ref{sect:common-lemmas-transfer-principles}).
\begin{lemma}\label{lemma:regular-affine}
Let~$X = \Spec A$ be an affine scheme. Let~$f \in A$. Then~$f$ is
a regular element of~$A$ if and only if~$f$ is a regular element of~$\O_X$ from
the internal perspective.\end{lemma}
\begin{proof}The Kripke--Joyal translation of internal regularity is:
\begin{quote}For any (without loss of generality: standard) open subset~$U \subseteq X$ and any function~$g \in
\Gamma(U,\O_X)$, $fg = 0$ in~$\Gamma(U,\O_X)$ implies~$g = 0$
in~$\Gamma(U,\O_X)$.\end{quote}
So the ``if'' direction is clear (use~$U \defeq X$). For the ``only if'' direction,
we use that~$\Gamma(U,\O_X)$ is a localization of~$A$ and that regular elements
remain regular in localizations.
\end{proof}

\subsection{Bézout property} Recall that a \emph{Bézout ring} is a ring in
which any finitely generated ideal is a principal ideal. In intuitionistic
mathematics, this is a better notion than that of a principal ideal ring: The
requirement that \emph{any} ideal is a principal ideal is far too strong.
Intuitively, this is because without any given generators to begin with, one
cannot hope to explicitly pinpoint a principal generator.
One can (provably) not even verify this property for the ring~$\ZZ$.\footnote{\label{fn:z-principal-ideal-domain}Assume
that any ideal of~$\ZZ$ is finitely generated. Let~$\varphi$ be an arbitrary
statement; we want to intuitionistically deduce~$\varphi \vee \neg\varphi$.
Consider the ideal~$\aaa \defeq \{ x \in \ZZ \,|\, (x = 0) \vee \varphi \}
\subseteq \ZZ$. The definition is such that~$\varphi$ holds if and only
if~$\aaa$ contains an element other than zero; and that~$\neg\varphi$ holds if
and only if zero is the only element of~$\aaa$.
By assumption,~$\aaa$ is finitely generated. Since~$\ZZ$ is a
Bézout ring, it is therefore even principal:~$\aaa = (x_0)$ for some~$x_0 \in
\ZZ$. Even intuitionistically we have~$(x_0 = 0) \vee (x_0 \neq 0)$ (for the
natural numbers, this can be proven by induction). In the first case, it
follows that~$\aaa$ contains only zero; in the second case, it follows
that~$\aaa$ contains an element other than zero. Thus~$\neg\varphi \vee
\varphi$.

This kind of reasoning is called \emph{exhibiting a Brouwerian
counterexample}. The definition of~$\aaa$ may look slightly dubious,
considering that~$\varphi$ does not depend on~$x$; but we will see that such
definitions actually have a clear geometric meaning -- they can be used to
define extensions of sheaves by zero in the internal language
(Lemma~\ref{lemma:extension-by-zero}).}

\begin{prop}Let~$X$ be a scheme (or a ringed space).
\begin{enumerate}
\item $\O_X$ is a Bézout ring from the internal perspective if and only if all
rings~$\O_{X,x}$ are Bézout rings.
\item $\O_X$ is such that, from the internal perspective, of any two elements,
one divides the other, if and only if all rings~$\O_{X,x}$ are such.
\end{enumerate}
\end{prop}
\begin{proof}Both properties can be formulated as geometric implications:
\begin{multline*}
  \text{(1)}\quad
  \forall f,g\?\O_X\_
  \top \Rightarrow
  \exists d\?\O_X\_
  (\exists a,b\?\O_X\_ d = af + bg) \wedge {} \\
  (\exists u\?\O_X\_ f = ud) \wedge
  (\exists v\?\O_X\_ g = vd)
\end{multline*}
\[
  \text{(2)}\quad
  \forall f,g\?\O_X\_
  \top \Rightarrow
  (\exists u\?\O_X\_ f = ug) \,\vee\,
  (\exists u\?\O_X\_ g = uf) \qquad\qquad\qquad\quad \qedhere
\]
\end{proof}

\begin{cor}\label{cor:dedekind-smith}
Let~$X$ be a Dedekind scheme, \ie a locally Noetherian normal scheme
of dimension~$\leq 1$. Then, from the internal perspective, any matrix
over~$\O_X$ can be put into Smith canonical form, \ie is equivalent to a
(rectangular) diagonal matrix with diagonal entries~$a_1|a_2|\cdots|a_n$
successively dividing each other.
\end{cor}
\begin{proof}It is well-known that such a scheme has principal ideal domains as
local rings~$\O_{X,x}$. For local domains, the Bézout condition is equivalent to the
property that of any two elements, one divides the other. Therefore all local
rings have this property, and by the previous proposition, the internal
ring~$\O_X$ has it as well. The statement thus follows from interpreting the
following fact of linear algebra in the internal universe: Let~$R$ be a ring
such that of any two elements, one divides the other. Then any matrix over~$R$
can be put into Smith canonical form.

The usual proof of this fact is indeed intuitionistically valid: Let a
matrix over~$R$ be given. By induction, one can show that for any finite family
of ring elements, one divides all the others. Hency some matrix entry is a factor
of all the others. We move this entry to the upper left by row and column
transformations and then kill the other entries of the first row and the first
column. After these operations, it is still the case that the entry in the
first row and column is a factor of all other entries. Continuing in this
fashion, we obtain a diagonal matrix. Its diagonal entries already fulfill
the divisibility condition and thus do not have to be sorted.
\end{proof}

Phrases such as ``if by chance the entry in the upper left divides
all the others, we can directly proceed with the next step; otherwise, some
other entry must be a factor of all entries, so \ldots'' may not be included in
a proof which is intended to be intuitionistically valid.
Those phrases assume that one may make the case distinction that for
any two ring elements~$x,y$, either~$x$ divides~$y$ or not. Fortunately, those
case distinctions are in fact superfluous.

A consequence of the corollary is that internally to the sheaf topos of a
Dedekind scheme, the usual structure theorem on finitely
presented~$\O_X$-modules is available. We will exploit this in
Lemma~\ref{lemma:torsion-stuff}, where we give an internal proof of the
fact that on Dedekind schemes, torsion-free~$\O_X$-modules are locally free.

\subsection{Normality} We will discuss the property of a ring to be
\emph{normal}, \ie to be integrally closed in its total field of
fractions, in Section~\ref{sect:normality}, after giving an internal
characterization of the sheaf of rational functions.

\subsection{Special properties of constant sheaves of rings} Let~$R$ be an
ordinary ring and~$\ul{R}$ the associated sheaf of locally constant~$R$-valued
functions on a topological space. If~$R$ is reduced, local, or a field,
then~$\ul{R}$ is so as well, from the internal point of view.

We will prove this in greater generality: Appropriately formulated, a constant
sheaf~$\ul{R}$ has some property~$\varphi$ from the internal point of view if
and only if~$R$ has the property~$\varphi$ externally
(Lemma~\ref{lemma:properties-of-constant-sheaves}).

\subsection{Noetherian conditions}
\label{sect:noetherian}

Recall the usual notion of a Noetherian ring: Any sequence~$\aaa_0 \subseteq
\aaa_1 \subseteq \cdots$ of ideals should stabilize, \ie there should exist a
natural number~$n$ such that~$\aaa_n = \aaa_{n+1} = \cdots$.

Intuitionistically, this definition has two problems. Firstly, without the
axiom of dependent choice, it is often not possible to construct a
\emph{sequence} of ideals: Often, it is only possible to show that there
\emph{exists} a suitable ideal~$\aaa_{n+1}$ depending on~$\aaa_n$. But since in
general there is no canonical choice for this successor ideal, the axiom of dependent choice
would be required to collect those into a sequence, \ie a function from~$\NN$
to the set of ideals.

Secondly, the conclusion that the sequence stabilizes is too strong.
Intuitionistically, one cannot even show that a weakly descending sequence of
natural numbers stabilizes in this sense; the statement that one could is
equivalent to the \emph{limited principle of omniscience for~$\NN$}.
Intuitionistically, it is only true that a weakly descending sequence~$a_0 \geq
a_1 \geq \cdots$ of natural numbers eventually \emph{stalls} in the sense that
there exists an index~$n$ such that~$a_n = a_{n+1}$ (but~$a_{n+1} > a_{n+2}$ is
allowed).\footnote{Classically, the following three statements about a ring are
equivalent: (1)~Every ascending chain of ideals stabilizes. (2)~Every ascending
chain of finitely generated ideals stabilizes. (3)~Every ascending chain of
finitely generated ideals stalls.}

We give two constructively inequivalent notions of Noetherian rings. The first
one is of independent constructive interest and enjoys the property that the structure sheaf
of a scheme~$X$ satisfies the Noetherian condition from the internal point of
view of~$\Sh(X)$ if and only if all stalks~$\O_{X,x}$ are Noetherian.

The second one is quite weak from a constructive point of view, but still
interesting from a geometric point of view and useful enough to derive
nontrivial consequences. It is satisfied by the structure sheaf of any (not
necessarily locally Noetherian) reduced scheme.

There are several proposals for a constructively sensible definition of
Noetherian rings in the literature on constructive algebra, each with unique
advantages and
disadvantages~\cite{richman:noetherian,mines-richman-ruitenburg:constructive-algebra,perdry:noetherian,perdry:lazy,perdry-schuster:noetherian,tennenbaum:hilbert}.
Insightful comments on why this is so can be found in the introduction and more
specifically on page~27 of the textbook by Lombardi and
Quitté~\cite{lombardi:quitte:constructive-algebra}.

{\tocless

\subsection*{Processly Noetherian rings}

\begin{defn}Let~$M$ be a partially ordered set. An \emph{ascending process
with values in~$M$} consists of an initial value~$x_0 \in M$ and a function~$f
: M \to \P(M)$ such that for any~$x \in M$ and any~$y \in f(x)$, $x \preceq y$,
and such that:
\begin{itemize}
\item The set $f(x_0)$ is inhabited.
\item For any~$x_1 \in f(x_0)$, the set $f(x_1)$ is inhabited.
\item For any~$x_1 \in f(x_0)$ and any~$x_2 \in f(x_1)$, the set $f(x_2)$ is inhabited.
\item And so on.
\end{itemize}
Such a process \emph{stalls} if and only if there
exists a step~$n$ and elements~$x_1, \ldots, x_n$ such that~$x_{i+1}
\in f(x_i)$ for~$i = 0,\ldots,n-1$ and such that~$x_n \in f(x_n)$.
The set~$M$ satisfies the \emph{ascending process condition} if and only if every
ascending process with values in~$M$ stalls.
\end{defn}

Intuitively, we picture~$f(x)$ as the set of all possible results of running
the process for a single step, starting with the value~$x$. This set could
be a singleton, in case that the process deterministically produces a single
value, but it may also contain more than one element, for instance if
the process cannot provide the next value in a canonical way. Instead of
arbitrarily choosing a definitive value for its result, the process may instead
collect all the possible values in the set~$f(x)$.

\begin{rem}The usual term for what we call ``to stall'' is
``to halt''~\cite{richman:noetherian,perdry-schuster:noetherian}. However, this
choice of wording is slightly unfortunate, since the phrase ``the process
halts'' intuitively suggests that the process stops and won't produce further
results in the future, even though this is not what is mathematically meant. We
are grateful to Matthias Hutzler for proposing ``to stall'', which seems quite
appropriate.
\end{rem}

\begin{defn}A ring~$A$ is \emph{processly Noetherian} if and only if the
set of finitely generated ideals in~$A$ satisfies the ascending process
condition.\end{defn}

An ascending chain of elements~$a_0 \preceq a_1 \preceq \cdots$ in a partially
ordered set gives rise to an ascending process by setting~$x_0 \defeq a_0$
and~$f(x) \defeq \{ y \,|\, \exists n\_ x = a_n \wedge y = a_{n+1} \}$.
(This process stalls iff there is an index~$n$ such that~$a_n = a_{n+1}$.)
Conversely, the axiom of dependent choice would allow to construct an
ascending chain from an ascending process. In
a classical context, a ring is therefore
processly Noetherian if and only if it is Noetherian in the usual sense.

The notion of a processly Noetherian ring works well in an
intuitionistic context: Important rings such as~$\ZZ$ and more generally~$\O_K$
for any algebraic number field~$K$ are processly Noetherian, and matrices
over Bézout rings which are integral domains in the weak sense and
processly Noetherian can be put into Smith canonical form.\footnote{In the
algorithm to put a matrix into Smith canonical form, one has to repeatedly
\emph{choose} generators for principal ideals and associated Bézout
representations (see for instance~\cite[Section~4]{richman:noetherian}). Since
these choices are not unique, the algorithm doesn't produce a \emph{sequence}
of intermediate ideals, but only a \emph{process}. This example was our main
motivation for the notion of processly Noetherian rings.}

Richman also studied Noetherian rings in a constructive context without
dependent choice~\cite{richman:noetherian}. His notion of \emph{ascending tree
condition} is equivalent to our ascending process condition. His condition
emphasizes the branching nature of a non-deterministic computation, while ours
emphasizes the step-for-step picture of computation.

There are three reasons why we did not define a ring to be processly Noetherian
if and only if the set of all (not necessarily finitely generated) ideals satisfies
the ascending process condition. Firstly, this stricter condition excludes
rings as~$\ZZ$.\footnote{The main ingredient in the proof that~$\ZZ$ is
Noetherian is that any ideal of~$\ZZ$ is a principal ideal, since (looking at
the prime factor decomposition) one can give explicit bounds on the length of
strictly ascending chains of principal ideals. However, as detailed in the
footnote on page~\pageref{fn:z-principal-ideal-domain}, constructively one
cannot show that every ideal of~$\ZZ$ is a principal ideal; one can only verify
that finitely generated ideals are principal. Geometrically, ideals which are
not finitely generated correspond to sheaves of ideals which may fail to
be quasicoherent.} Secondly, restricting to finitely generated
ideals in this context is a well-established procedure in constructive
mathematics~\cite{mines-richman-ruitenburg:constructive-algebra,richman:noetherian}
and suffices for the applications of the Noetherian condition one typically expects.
Thirdly, our definition provides a link to the external condition on a scheme
to be locally Noetherian, as shown by the following proposition.

\begin{prop}\label{prop:internal-noetherianity}
Let~$X$ be a scheme. The following statements are equivalent:
\begin{enumerate}
\item All stalks~$\O_{X,x}$ are Noetherian.
\item From the internal point of view of~$\Sh(X)$, the ring~$\O_X$ is processly
Noetherian.
\end{enumerate}
\end{prop}

\begin{proof}Statement~(1) can be reformulated in a way so it doesn't refer to stalks:
For any open affine subscheme~$U \subseteq X$ and any ascending chain~$\aaa_0
\subseteq \aaa_1 \subseteq \cdots$ of finitely generated ideals in~$\Gamma(U,
\O_X)$ there is a partition of unity~$1 = \sum_i f_i \in \Gamma(U, \O_X)$ such
that for each~$i$ there exists an index~$j$ such that~$\aaa_j = \aaa_{j+1}$ as
ideals of~$\Gamma(U, \O_X)[f_i^{-1}]$.

%

We'll verify the equivalence using this formulation. For proving the
direction ``(1)~$\Rightarrow$~(2)'', we may assume that~$X =
\Spec A$ is affine and that internally, we are given an
ascending process on the set of finitely generated ideals of~$\O_X$.
Externally, this is a finite type sheaf of ideals~$\I$ together with a morphism~$\M
\to \P(\M)$ where~$\M$ is the sheaf whose~$U$-sections are finite type ideal
sheaves of~$\O_X|_U$.

Since~$X \models \speak{$f(\I)$ is inhabited}$, there exists an open covering~$X
= \bigcup_i U_i$ and finite type sheaves of ideals~$\I_i \hookrightarrow
\O_X|_{U_i}$ such that~$U_i \models \I_i \in f(\I)$. Without loss of generality,
we may assume that the open sets~$U_i$ are standard open sets~$D(f_i)$ and that
the covering is finite. Since the sheaves~$\I_i$ are quasicoherent (being of
finite type, they are images of morphisms of the form~$\O_X|_{U_i}^n \to
\O_X|_{U_i}$), they correspond to ideals~$J_i \subseteq A[f_i^{-1}]$. We note for future reference
that for~$D(g) \subseteq D(f_i)$, the restricted sheaf of ideals~$\I_i|_{D(g)}$
corresponds to the extension of~$J_i$ in the further
localized ring~$A[g^{-1}]$.

For each~$i$,~$D(f_i) \models \speak{$f(\I_i)$ is inhabited}$.
Hence there exists an open covering~$D(f_i) = \bigcup_j D(f_{ij})$ and finite
type sheaves of ideals~$\I_{ij} \hookrightarrow \O_X|_{D(f_{ij})}$; these
correspond to ideals~$J_{ij} \subseteq A[f_{ij}^{-1}]$ such that~$J_i
\subseteq J_{ij}$ (where we have suppressed the localization
morphism~$A[f_i^{-1}] \to A[f_{ij}^{-1}]$ in the notation). Equivalently,
writing~$J_i' \defeq A \cap J_i$ and~$J_{ij}' \defeq A \cap J_{ij}$ for the
contractions, we have the inclusions~$J_i' \subseteq J_{ij}'$ of ideals of~$A$.

Continuing in this fashion, we obtain an infinite tree of ideals~$J_{i_1 \cdots i_n}'$.
We now prune this tree in the following fashion: If the node at position~$(i_1,
\ldots, i_n)$ has the property that~$D(f_{i_1 \cdots i_n}) \models
\bigvee_{m=0}^{n-1} (\I_{i_1 \cdots i_m} = \I_{i_1 \cdots i_{m+1}})$, then we cut
of all childs of this node.

The resulting tree doesn't contain an infinite path, since any
sequence~$J_{i_1}' \subseteq J_{i_1 i_2}' \subseteq \cdots$ locally stalls by
assumption on~$A$. Because only finitely many subtrees branch off at each node,
the tree is finite (this is an application of the graph-theoretical
\emph{\konig's lemma}).

The outermost nodes then yield an open covering of~$X$ such that, on each
member of the covering, the internal statement~$\speak{$f$ stalls}$ holds. By
the locality of the internal language, this statement holds on~$X$.

For the converse direction, let an affine open subset~$U \subseteq X$ and an
ascending sequence~$\aaa_0 \subseteq \aaa_1 \subseteq \cdots$ of finitely
generated ideals in~$\Gamma(U, \O_X)$ be given. Internally, we construct the
process
\[ f : \M \lra \P(\M),\ \I \longmapsto \{ \J \? \M \,|\,
  \bigvee_{n \geq 0} (\I = \aaa_n^\sim \wedge \J = \aaa_{n+1}^\sim) \} \]
with initial value~$\aaa_0^\sim$. The assumption that~$f$ stalls yields an open
covering~$U = \bigcup_i D(f_i)$ such that for each~$i$, there is an index~$n$
such that~$\aaa_n^\sim = \aaa_{n+1}^\sim$ on~$D(f_i)$, that is~$\aaa_n =
\aaa_{n+1}$ as ideals of~$\Gamma(U, \O_X)[f_i^{-1}]$.
\end{proof}

\begin{rem}The proof shows that, if the base scheme fulfills the stronger
condition that it is locally Noetherian, then internally speaking even the set
of all quasicoherent ideals (instead of merely the finitely generated
ones) fulfills the ascending process condition. We have not taken this property
as the definition of a processly Noetherian ring since it is a notion not
usually studied in constructive mathematics (compare
Remark~\ref{rem:qcoh-in-constructive-mathematics}).
\end{rem}

Proposition~\ref{prop:internal-noetherianity} looks like it could be a shadow
of a deeper result, since it states that, for the property of being a processly
Noetherian ring, the same relation between truth at every point and internal
truth holds as for geometric implications -- even though being a processly
Noetherian ring is a emphatically a higher-order condition, hence not a
geometric implication. It's conceivable that a metatheorem explaining this
phenomenon is hidden in the background.

There is also an internal characterization of the property that~$X$ is locally
Noetherian (in contrast to the property that all stalks are Noetherian).
However, as described above, the corresponding internal notion is of limited
usefulness.

\begin{prop}Let~$X$ be a scheme. The following statements are equivalent:
\begin{enumerate}
\item The scheme~$X$ is locally Noetherian.
\item From the internal point of view of~$\O_X$, any ascending chain of
finitely generated ideals stabilizes.
\end{enumerate}
\end{prop}

\begin{proof}Similar to, but easier than, the proof of
Proposition~\ref{prop:internal-noetherianity}.
\end{proof}


\subsection*{Anonymously Noetherian rings}

Classically, there is a characterization of Noetherian rings which doesn't
involve ascending sequences: A ring is Noetherian if and only if any of its
ideals is finitely generated. We mentioned in the footnote on
page~\pageref{fn:z-principal-ideal-domain} that this condition is far too
strong from a constructive point of view; not even the ring~$\ZZ$ verifies it.
However, it can be weakened to yield an interesting notion:

\begin{defn}A ring~$A$ is \emph{anonymously Noetherian} if and only if any ideal
of~$A$ is \notnot finitely generated. A module~$M$ is \emph{anonymously Noetherian}
if and only if any submodule of~$M$ is \notnot finitely generated.\end{defn}

\begin{ex}There is an intuitionistic proof that the ring~$\ZZ$ is anonymously
Noetherian: Let~$\aaa \subseteq \ZZ$ be any ideal. Under the assumption that
either there exists a nonzero element in~$\aaa$ or not, the ideal~$\aaa$ is
\notnot finitely generated, even \notnot principal: For in the first case, a
minimal element~$d$ of~$\aaa \cap \NN^+$ (which \notnot exists) witnesses~$\aaa
= (d)$. In the second case the ideal~$\aaa$ is the zero ideal. Since the
assumption is \notnot satisfied, the ideal~$\aaa$ is \notnot \notnot finitely
generated, so \notnot finitely generated. (We remark on this proof scheme on
page~\pageref{proof-scheme-boxed-statements}.) \end{ex}

It appears that this Noetherian condition has not been studied in the
literature on constructive algebra. Indeed, from the philosophical and the
computational point of view on constructive mathematics, the notion of an anonymously
Noetherian ring is not very useful: From an intuitionistic proof that a given
ideal is finitely generated one can mechanically extract explicit generators,
thereby satisfying computational or philosophical demands (``the generators are
really there''). In contrast, an intuitionistic proof that a given ideal is
\notnot finitely generated doesn't contain computational content in general.

Along the same lines, one could dismiss the fact (proved below) that Hilbert's
basis theorem, stating that~$A[X]$ is Noetherian if~$A$ is, holds for the
notion of anonymously Noetherian rings. In fact, one could feel mocked by this
version of Hilbert's basis theorem: It promises that any ideal of~$A[X]$
``has'' finitely many generators in some Platonic sense (without providing any
clue on how one might go on to find the generators -- they remain \emph{anonymous}), provided
that any ideal of~$A$ ``has'' finitely many generators in the same
sense.\footnote{We borrowed the term ``anonymous'' from type theory, where it
is used with a similar meaning (see for
instance~\cite{kraus-escardo-coquand-altenkirch:anonymous}). However, there is
a subtle difference: Unique and anonymous-in-the-sense-of-type-theory existence
implies existence, but \notnot existence does not.}

However, applications in the internal universe of toposes provide a further
motivation for constructive reasoning, related but distinct from computational
or philosophical considerations. A first indication that the notion of anonymously
Noetherian rings is useful is that the structure sheaf~$\O_X$ of any reduced
scheme is anonymously Noetherian from the internal point of view of~$\Sh(X)$. We
exploit this observation in Section~\ref{sect:generic-freeness} to give a short
proof of Grothendieck's generic freeness lemma.

Secondly, internal universes of toposes may satisfy certain classicality
principles which are not generally satisfied in constructive mathematics. If
a topos is set up in an intuitionistically sensible manner, one might then even
be able to extract constructive results. For instance, the structure sheaf of a
reduced scheme satisfies the principle
\[ \Sh(X) \models \forall s\?\O_X\_
  \neg\neg(s = 0) \Longrightarrow s = 0. \]
This principle can be put to use as follows. Let's consider the situation that
we have an intuitionistic proof that some ring element~$s$ is zero under the
assumption that some ideal~$\aaa$ is finitely generated. Hence we also have an
intuitionistic proof that~$s$ is \notnot zero under the assumption that~$\aaa$
is \notnot finitely generated. This assumption could be validated by the anonymously
Noetherian property, yielding an unconditional proof that~$s$ is \notnot zero.
Usually in constructive mathematics, we would be stuck at this point; but
internally in~$\Sh(X)$, we may continue and deduce that~$s$ is actually zero.

This observation is the basis for a fully constructive proof of Grothendieck's
generic freeness lemma and puts our work in line of the general research
program of extracting constructive content from classical
proofs~\cite{coquand:classical,feferman:kreisel,kohlenbach:applprooftheory,kiselyov:lem}.

\begin{thm}\label{thm:hilbert}
Let~$A$ be an anonymously Noetherian ring. Then the polynomial
algebra~$A[X]$ is anonymously Noetherian as well, intuitionistically.
\end{thm}
\begin{proof}Classically, this is precisely the statement of Hilbert's basis
theorem, whose usual accounts do not care about the sensibilities of
constructive mathematics. However, a careful reading of for instance the proof
given in~\cite[Theorem~7.5]{atiyah:macdonald:commutative-algebra} shows that
the theorem holds intuitionistically as stated.
\end{proof}

\begin{lemma}Let~$0 \to M' \to M \to M'' \to 0$ be a short exact sequence of
modules. Intuitionistically, the module~$M$ is anonymously Noetherian if and only
if~$M'$ and~$M''$ are.
\end{lemma}

\begin{proof}The usual proof applies.\end{proof}

\begin{prop}\label{prop:ox-anonymously-noetherian}
Let~$X$ be an arbitrary reduced scheme (not necessarily locally
Noetherian). Then~$\O_X$ is anonymously Noetherian from the internal point of view
of~$\Sh(X)$.\end{prop}
\begin{proof}By Corollary~\ref{cor:field-reduced}, the ring~$\O_X$ fulfills a
suitable field condition from the internal point of view. Therefore it suffices
to give an intuitionistic proof of the following statement: Let~$k$ be a ring such that
any element of~$k$ which is not invertible is zero. Then any ideal of~$k$ is
\notnot finitely generated.

Let~$\aaa \subseteq k$ be an arbitrary ideal. We have~$\neg\neg(1 \in \aaa \vee
1 \not\in \aaa)$. Therefore~$\neg\neg(\aaa = (1) \vee \aaa = (0))$. Thus~$\aaa$
is \notnot finitely generated (even \notnot principal).
\end{proof}

The external translation of the statement
that~$\O_X[U_1,\ldots,U_n]$ is anonymously Noetherian was displayed on
page~\pageref{page:convoluted-statement}, as an example of a convoluted
statement which profits from the simpler internal account.

}

\section{Sheaves of modules}
\label{sect:sheaves-of-modules}

From the internal perspective, a sheaf of~$\R$-modules, where~$\R$ is a sheaf
of rings, looks just like a plain module over the plain ring~$\R$. This is
proven just as the correspondence between sheaf of rings and internal rings
(Proposition~\ref{prop:rings-internally}).

\subsection{Finite local freeness}

Recall that an~$\O_X$-module~$\F$ is \emph{finite locally free} if and only
if there exists a covering of~$X$ by open subsets~$U$ such that on each
such~$U$, the restricted module~$\F|_U$ is isomorphic as an~$\O_X|_U$-module
to~$(\O_X|_U)^n$ for some natural number~$n$ (which may depend on~$U$).

\begin{prop}\label{prop:locally-free}
Let~$X$ be a scheme (or a ringed space). Let~$\F$ be
an~$\O_X$-module. Then~$\F$ is finite locally free if and only if, from the
internal perspective,~$\F$ is a finite free module, \ie
\[ \Sh(X) \models \bigvee_{n \geq 0} \speak{$\F \cong (\O_X)^n$}, \]
or more elementarily
\[ \Sh(X) \models \bigvee_{n \geq 0}
  \exists x_1,\ldots,x_n\?\F\_
  \forall x\?\F\_
  \exists! a_1,\ldots,a_n\?\O_X\_
  x = \textstyle\sum\limits_i a_i x_i. \]
\end{prop}
\begin{proof}By the expression~``$(\O_X)^n$'' in the internal language we mean
the internally constructed object~$\O_X \times \cdots \times \O_X$ with its
componentwise~$\O_X$-module structure. This coincides with the sheaf~$(\O_X)^n$ as
usually understood.

It is clear that the two stated internal conditions are equivalent, since the
corresponding proof in linear algebra is intuitionistically valid. The equivalence with
the external notion of finite local freeness follows because the
interpretation of the first condition with the Kripke--Joyal semantics is the
following: There exists a covering of~$X$ by open subsets~$U$ such that for
each such~$U$, there exists a natural number~$n$ and a morphism of
sheaves~$\varphi : \F|_U \to (\O_X|_U)^n$ such that
\[ U \models \speak{$\varphi$ is~$\O_X$-linear} \quad\text{and}\quad
  U \models \speak{$\varphi$ is bijective}. \]
The first subcondition means that~$\varphi$ is a morphism of sheaves
of~$\O_X|_U$-modules and the second one means that~$\varphi$ is an isomorphism of
sheaves.
\end{proof}

\begin{rem}There are intuitionistic proofs of the following facts:
An~$R$-module is a dualizable object in the monoidal category of all
$R$-modules if and only if it is finitely generated and projective. If~$R$ is
local, then an~$R$-module is finitely generated and projective if and only if
it is finite free. Therefore an~$\O_X$-module is internally dualizable if and
only if is finite locally free.
\end{rem}

\subsection{Finite type, finite presentation, coherence}
Recall the conditions of an~$\O_X$-module~$\F$ on a scheme~$X$ (or a ringed
space) to be of finite type, of finite presentation, and to be coherent:
\begin{enumerate}
\item $\F$ is \emph{of finite type} if and only if there exists a covering of~$X$ by
open subsets~$U$ such that for each such~$U$, there exists an exact sequence
\[ (\O_X|_U)^n \longrightarrow \F|_U \longrightarrow 0 \]
of~$\O_X|_U$-modules.
\item $\F$ is \emph{of finite presentation} if and only if there exists a covering of~$X$ by
open subsets~$U$ such that for each such~$U$, there exists an exact sequence
\[ (\O_X|_U)^m \longrightarrow (\O_X|_U)^n \longrightarrow \F|_U \longrightarrow 0. \]
\item $\F$ is \emph{coherent} if and only if~$\F$ is of finite type and the
kernel of any~$\O_X|_U$-linear morphism~$(\O_X|_U)^n \to \F|_U$, where~$U \subseteq
X$ is any open subset, is of finite type.
\end{enumerate}

The following proposition gives translations of these definitions into the
internal language.
\begin{prop}\label{prop:finite-type-and-co}
Let~$X$ be a scheme (or a ringed space). Let~$\F$ be
an~$\O_X$-module. Then:
\begin{enumerate}
\item $\F$ is of finite type if and only if~$\F$, considered as an ordinary
module from the internal perspective, is finitely generated, \ie if
\[ {\qquad\qquad} \Sh(X) \models
  \bigvee_{n \geq 0}
  \exists x_1,\ldots,x_n\?\F\_
  \forall x\?\F\_
  \exists a_1,\ldots,a_n\?\F\_
  x = \textstyle\sum\limits_i a_i x_i. \]
\item $\F$ is of finite presentation if and only if~$\F$ is a finitely
presented module from the internal perspective, \ie if
\[ {\qquad\qquad} \Sh(X) \models \bigvee_{n,m \geq 0}
  \speak{there is a short exact sequence $\O_X^m \to \O_X^n \to \F \to 0$}.
  \]
\item $\F$ is coherent if and only if~$\F$ is a coherent module from the
internal perspective, \ie if
\begin{multline*}
{\qquad\qquad\qquad}
  \Sh(X) \models \speak{$\F$ is finitely generated} \mathop{\wedge} \\
  \bigwedge_{n \geq 0} \forall \varphi \? \HOM_{\O_X}(\O_X^n,\F)\_
  \speak{$\Kernel \varphi$ is finitely generated}.
\end{multline*}
\end{enumerate}
\end{prop}
\begin{proof}Straightforward -- the translations of the internal statements using
the Kripke--Joyal semantics are precisely the corresponding external
statements.
\end{proof}

\begin{rem}We believe that Proposition~\ref{prop:finite-type-and-co} settles a
question Lawvere raised on the category theory mailing
list~\cite{lawvere:finiteness-question}: ``What concept of finiteness is
appropriate for those important mathematical applications in topology for which
[Kuratowski-finiteness] doesn't seem right? [\ldots\!] Especially, a suitably
`finite' module should be a vector bundle or a [coherent sheaf] in the sense of
Serre so that our simplified topos theory could apply more directly to those
things it should.''
\end{rem}

Recall that an~$\O_X$-module~$\F$ is \emph{generated by global sections} if and
only if there exist global sections~$s_i \in \Gamma(X,\F)$ such that for any~$x
\in X$, the stalk~$\F_x$ is generated by the germs of the~$s_i$.
This condition is of course not local on the base. Therefore there cannot
exist a formula~$\varphi$ such that for any space~$X$ and
any~$\O_X$-module~$\F$ it holds that~$\F$ is generated by global sections if
and only if~$\Sh(X) \models \varphi(\F)$. But still, global generation can be
characterized by a mixed internal/external statement:

\begin{prop}Let~$X$ be a scheme (or a ringed space). Let~$\F$ be
an~$\O_X$-module. Then~$\F$ is generated by global sections if and only if
there exist global sections~$s_i \in \Gamma(X,\F)$, $i \in I$ such that
\[ \Sh(X) \models \forall x\?\F\_ \bigvee\nolimits_{\textnormal{$J=\{i_1,\ldots,i_n\} \subseteq I$ finite}}
  \exists a_1,\ldots,a_n\?\O_X\_
  x = \sum_j a_j s_{i_j}. \]
Furthermore,~$\F$ is generated by finitely many global sections if and only if
there exist global sections~$s_1,\ldots,s_n \in \Gamma(X,\F)$ such that
\[ \Sh(X) \models \forall x\?\F\_ \exists a_1,\ldots,a_n\?\O_X\_ x = \sum_j a_j
s_j. \]
\end{prop}
\begin{proof}The given internal statements are geometric implications, their
validity can thus be checked stalkwise.\end{proof}

\begin{rem}The analogue of Proposition~\ref{prop:finite-type-and-co} for
sheaves of algebras instead of sheaves of modules holds. More precisely,
let~$\A$ be a sheaf of~$\O_X$-algebras on a scheme~$X$ (or a ringed space).
Then:
\begin{enumerate}
\item $\A$ is of finite type if and only if~$\A$, considered as an ordinary
algebra from the internal perspective, is finitely generated, \ie if
\[ {\qquad\qquad} \Sh(X) \models
  \bigvee_{n \geq 0}
  \exists x_1,\ldots,x_n\?\A\_
  \forall x\?\F\_
  \exists p\?\O_X[X_1,\ldots,X_n]\_
  x = p(x_1,\ldots,x_n). \]
\item $\A$ is of finite presentation if and only if~$\A$ is a finitely
presented algebra from the internal perspective, \ie if
\begin{multline*}
  {\qquad\qquad} \Sh(X) \models \bigvee_{n,m \geq 0}
  \exists f_1,\ldots,f_m \? \O_X[X_1,\ldots,X_n]\_ \\
  \speak{$\A \cong \O_X[X_1,\ldots,X_n]/(f_1,\ldots,f_m)$}. \qquad
\end{multline*}
\end{enumerate}
\end{rem}

\subsection{Tensor product and flatness} The tensor product
of~$\O_X$-modules~$\F$ and~$\G$ on a scheme~$X$ (or a ringed space) is usually
constructed as the sheafification of the presheaf
\[ \text{$U \subseteq X$ open} \quad\longmapsto\quad \Gamma(U,\F) \otimes_{\Gamma(U,\O_X)}
\Gamma(U,\G). \]
From the internal point of view,~$\F$ and~$\G$ look like ordinary modules, so
that we can consider their tensor product as usually constructed in
commutative algebra, as a certain quotient of the free module on the elements
of~$\F \times \G$:
\[ \O_X\langle x \otimes y \,|\, x\?\F, y\?\G \rangle / R, \]
where~$R$ is the submodule generated by
\begin{gather*}
  (x+x') \otimes y - x \otimes y - x' \otimes y, \\
  x \otimes (y+y') - x \otimes y - x \otimes y', \\
  (sx) \otimes y - s(x \otimes y), \\
  x \otimes (sy) - s(x \otimes y)
\end{gather*}
with~$x,x'\?\F$, $y,y'\?\G$, $s\?\O_X$.
This internal construction gives rise to the same sheaf
of modules as the externally defined tensor product:

\begin{prop}\label{prop:internal-tensor-product}
Let~$X$ be a scheme (or a ringed space). Let~$\F$ and~$\G$
be~$\O_X$-modules. Then the internally constructed tensor product~$\F
\otimes_{\O_X} \G$ coincides with the external one.
\end{prop}
\begin{proof}
Since the proof of the corresponding fact of commutative algebra is
intuitionistically valid, the internally defined tensor product~$\F \otimes_{\O_X} \G$
has the following universal property: For any~$\O_X$-module~$H$,
any~$\O_X$-bilinear map~$\F \times \G \to H$ uniquely factors over the
canonical map~$\F \times \G \to \F \otimes_{\O_X} \G$.

Interpreting this property with the Kripke--Joyal semantics, we see that the
internally constructed tensor product has the following external property:
For any open subset~$U \subseteq X$ and any~$\O_X|_U$-module~$\H$ on~$U$,
any~$\O_X|_U$-bilinear morphism~$\F|_U \times \G|_U \to \H$ uniquely factors over the
canonical morphism~$\F|_U \times \G|_U \to (\F \otimes_{\O_X} \G)|_U$.

In particular, for~$U = X$, this property is well-known to be the universal
property of the externally constructed tensor product. Therefore the
claim follows.
\end{proof}

A description of the stalks of the tensor product
follows purely by considering the logical form of the construction:
\begin{cor}Let~$X$ be a scheme (or a ringed space). Let~$\F$ and~$\G$
be~$\O_X$-modules. Then the stalks of the tensor product coincide with the
tensor products of the stalks: $(\F \otimes_{\O_X} \G)_x \cong \F_x
\otimes_{\O_{X,x}} \G_x$.\end{cor}
\begin{proof}
We constructed the tensor product using the following operations: product of
two sets, free module on a set, quotient module with respect to a submodule;
submodule generated by a set of elements given by a geometric formula.
All of these operations are geometric, so the tensor product construction is
geometric as well (see Section~\ref{sect:geometric-formulas-and-constructions}). Hence taking stalks commutes with performing the
construction.
\end{proof}

Recall that an~$\O_X$-module~$\F$ is \emph{flat} if and only if all
stalks~$\F_x$ are flat~$\O_{X,x}$-modules. We can characterize flatness in the
internal language.
\begin{prop}\label{prop:flatness}
Let~$X$ be a scheme (or a ringed space). Let~$\F$ be
an~$\O_X$-module. Then~$\F$ is flat if and only if, from the internal
perspective,~$\F$ is a flat~$\O_X$-module.
\end{prop}
\begin{proof}
Recall that flatness of an~$A$-module~$M$ can be characterized without
reference to tensor products by the following condition (using
suggestive vector notation): For any natural number~$p$,
any $p$-tuple~$m \? M^p$ of elements of~$M$ and
any~$p$-tuple $a \? A^p$ of elements of~$A$, it should hold that
\[
  a^T m = 0 \ \Longrightarrow\
  \bigvee\limits_{q \geq 0} \exists n\?M^q, B\?A^{p \times q}\_
  Bn = m \wedge a^T B = 0. \]
The equivalence of this condition with tensoring being exact holds
intuitionistically as
well~\cite[Theorem~III.5.3]{mines-richman-ruitenburg:constructive-algebra}.
This formulation of flatness has the advantage that it is the conjunction of
geometric implications (one for each~$p \geq 0$); therefore it holds internally
if and only if it holds at any point.
\end{proof}

\subsection{Support} Recall that the \emph{support} of an~$\O_X$-module~$\F$ is
the subset~$\supp\F \defeq \{ x \in X \,|\, \F_x \neq 0 \} \subseteq X$. If~$\F$ is
of finite type, this set is closed, since its complement is then open by a
standard lemma. (We will give an internal proof of this fact in
Lemma~\ref{lemma:module-zero-point-neighbourhood}.)

\begin{prop}\label{prop:characterization-support}
Let~$X$ be a scheme (or a ringed space). Let~$\F$ be
an~$\O_X$-module. Then the interior of the complement of the support of~$\F$
can be characterized as the largest open subset of~$X$ on which the internal
statement~$\F = 0$ holds.
\end{prop}
\begin{proof}
For any open subset~$U \subseteq X$, it holds that:
\begin{align*}
  &\ U \subseteq \Int(X \setminus \supp \F) \\
  \Longleftrightarrow&\ U \subseteq X \setminus \supp \F \\
  \Longleftrightarrow&\ U \subseteq \{ x \in X \,|\, \forall s \in \F_x\_ s = 0 \} \\
  \Longleftrightarrow&\ U \models \forall s\?\F\_ s = 0 \\
  \Longleftrightarrow&\ U \models \speak{$\F = 0$}
\end{align*}
The second to last equivalence is because~``$\forall s\?\F\_ s = 0$'' is a
geometric implication and can thus be checked stalkwise.
\end{proof}

\begin{rem}\label{rem:support-sheaf-of-sets}
The support of a sheaf of \emph{sets}~$\F$ is defined as the subset~$\{ x \in X \,|\,
\text{$\F_x$ is not a singleton} \}$. A similar proof shows that the interior
of its complement can be characterized as the largest open subset of~$X$ where
the internal statement~$\speak{$\F$ is a singleton}$ holds.\end{rem}

\subsection{Torsion} Let~$R$ be a ring. Recall that the
\emph{torsion submodule}~$M_\tors$ of an~$R$-module~$M$ is defined as
\[ M_\tors \defeq \{ x\?M \,|\, \exists a\?R\_ \speak{$a$ regular} \wedge ax = 0 \}
\subseteq M. \]
This definition is meaningful even if~$R$ is not an integral domain.
An~$R$-module~$M$ is \emph{torsion-free} if and only if~$M_\tors$ is
the zero submodule; an~$R$-module~$M$ is a \emph{torsion module} if and only
if~$M_\tors = M$.

Recall also that if~$\F$ is a sheaf of~$\O_X$-modules on an integral
scheme~$X$, there is a unique subsheaf~$\F_\tors \subseteq \F$ with the
property that~$\Gamma(U,\F_\tors) = \Gamma(U,\F)_\tors$ for all affine open
subsets~$U \subseteq X$. The content of the following proposition is that
internally constructing the torsion submodule of~$\F$, regarded as a plain
module from the internal perspective, gives exactly this subsheaf. There is
therefore no harm in using the same notation~``$\F_\tors$'' for the result of
the internal construction.

\begin{prop}\label{prop:torsion-int-ext}Let~$X$ be an integral scheme. Let~$\F$ be an~$\O_X$-module. Let~$U
= \Spec A \subseteq X$ be an affine open subset. Let~$s \in \Gamma(U,\F)$ be a local
section. Then
\[ s \in \Gamma(U,\F)_\tors \quad\text{if and only if}\quad
  U \models s \in \F_\tors. \]
\end{prop}
\begin{proof}
The ``only if'' direction is trivial in view of
Lemma~\ref{lemma:regular-affine}: If~$s$ is a torsion element
of~$\Gamma(U,\F)$, there exists a regular element~$a \in \Gamma(U,\O_X)$ such
that~$as = 0$. By the lemma, this element is regular from the internal
perspective as well, so~$U \models \speak{$a$ regular} \wedge as = 0$.

For the ``if'' direction, we may assume that there exists an open covering~$X =
\bigcup_i U_i$ by standard open subsets~$U_i = D(f_i)$ such that there are
sections~$a_i \in \Gamma(U_i,\O_X) = A[f_i^{-1}]$ with~$U_i \models \speak{$a_i$ regular}
\wedge a_i s = 0$. Without loss of generality, we may assume that the
denominators of the~$a_i$'s are ones, that the $f_i$ are
finite in number, and that the~$f_i$ are regular (\ie nonzero, since~$A$ is an
integral domain). By Lemma~\ref{lemma:regular-affine}, the~$a_i$ are
regular in~$A[f_i^{-1}]$ and by regularity of the~$f_i$ also regular in~$A$.
Therefore their product~$\prod_i a_i \in A$ is regular in~$A$ as well and
annihilates~$s$.
\end{proof}
%
%
%
%
%
%

\begin{prop}\label{prop:torsion-submodule-stalks}
Let~$X$ be a locally Noetherian scheme. Let~$\F$ be
an~$\O_X$-module. Let~$x \in X$ be a point. Then~$(\F_\tors)_x =
(\F_x)_\tors$.\end{prop}
\begin{proof}This would be obvious if the condition on an element~$s\?\F$ to
belong to~$\F_\tors$ were a geometric formula. Because of the universal
quantifier, it is not:
\[ s \in \F_\tors \quad\Longleftrightarrow\quad
  \exists a\?\O_X\_ (\forall b\?\O_X\_ ab = 0 \Rightarrow b = 0) \wedge as = 0. \]
But since~$X$ is assumed to be locally Noetherian, regularity is an open
property nonetheless (see Proposition~\ref{prop:regularity-spreading} for an
internal proof of this fact). Thus the claim still follows, just like in the
proof of Proposition~\ref{prop:internal-integrality}.
\end{proof}

\subsection{Kähler differentials}\label{sect:kaehler-differentials}

Let~$A \to B$ be a homomorphism of rings. The~$B$-module~$\Omega^1_{B|A}$ of
Kähler differentials can be constructed as the free~$B$-module on the
basis~$(db)_{b \in B}$ consisting of formal symbols modulo appropriate
relations ensuring that the map~$b \mapsto db$ is~$A$-linear and satisfies the
Leibniz rule; it verifies the universal property that the map~$B \to
\Omega^1_{B|A}$ is the initial~$A$-linear derivation of~$B$.

For constructing the sheaf of Kähler differentials for a morphism of schemes,
one often resorts to the alternative construction as~$I/I^2$, where~$I
\subseteq B \otimes_A B$ is the kernel of the multiplication map~$B \otimes_A B
\to B$. The verification of the universal property is slightly harder for this
construction. Gathmann comments on this situation as follows~\cite[p.~134]{gathmann:ag}:
\begin{quote}
Of course, if~$f : X \to Y$ is a morphism of general (not necessarily affine)
schemes, we want to consider the relative differentials of every restriction
of~$f$ to affine opens of~$X$ and~$Y$, and glue them together to get a quasi-coherent
sheaf~$\Omega_{X|Y}$. To do this, we have to give a different description of the relative
differentials, as the construction [via the free module] does not glue very well.
\end{quote}
After having constructed a global sheaf of Kähler differentials using the
alternative description (only in the case that~$f$ is separated, though with a
little bit of more work this assumption can be dispensed with), he goes on as
follows~\cite[Remark~7.4.8]{gathmann:ag}:
\begin{quote}
It should be stressed that [the definition using the alternative construction]
is essentially useless for practical computations. Its only use is to show that
a global object~$\Omega_{X|Y}$ exists that restricts to the old definition on
affine open subsets. For applications, we will always use [the definition using
the free module and the calculation of the Kähler differentials of a morphism
of the form~$k \to k[x_1,\ldots,x_n]/(f_1,\ldots,f_m)$] on open subsets.
\end{quote}
Vakil chooses a similar route~\cite[Section~21.2]{vakil:foag}.

Regarding the sheaf of Kähler differentials as the conormal sheaf of the
diagonal embedding (and hence using the alternative construction) is of course
essential for further developments and hence very useful. However, if the goal
is just to construct a global sheaf of Kähler differentials, we can employ the
internal language to construct it using only the formal construction via the
free module.

Specifically, if~$f : X \to S$ is a morphism of schemes, we can construct, in
the internal universe of~$\Sh(X)$, the module of Kähler differentials of the
morphism~$f^\sharp : f^{-1}\O_S \to \O_X$. A number of basic properties then
follow purely formally:

\begin{prop}\label{prop:kaehler}Let~$f : X \to S$ be a morphism of schemes (or of locally ringed
spaces). Let~$\Omega^1_{X|S}$ be defined as the interpretation of the internal
construction~$\Omega^1_{\O_X|f^{-1}\O_S}$.
\begin{enumerate}
\item The sheaf of Kähler differentials has the following universal property:
For any open subset~$U \subseteq X$, any~$(f^{-1}\O_S)|_U$-linear
derivation~$\O_X|_U \to \E$ over~$U$ uniquely factors over~$\O_X|_U \to
\Omega^1_{X|S}|_U$.
\item The stalks~$(\Omega^1_{X|S})_x$ are canonically isomorphic
to~$\Omega^1_{\O_{X,x}|\O_{S,f(x)}}$.
\item The sheaf~$\Omega^1_{X|S}$ is quasicoherent.
\end{enumerate}
\end{prop}

\begin{proof}The first claim is just the interpretation of the internal
universal property using the Kripke--Joyal semantics and using the
simplification rule given in Proposition~\ref{prop:simplification}.

The second claim is immediate because the construction of Kähler differentials
via the free module is geometric (page
\pageref{page:geometric-constructions}). Therefore the operations~``taking
Kähler differentials'' and~``taking the stalk at~$x$'' commute.

For verifying the third claim, it suffices to verify that~$\Omega^1_{X|S}$ is
quasicoherent in the case that~$X = \Spec(B)$ and~$S = \Spec(A)$ are affine.
Let~$\varphi : A \to B$ be the homomorphism of rings given by~$f$. We employ the
technique and notation of Section~\ref{sect:generic-filter}:
\begin{align*}
  (\Omega^1_{B|A})^\sim &\cong
  \underline{\Omega^1_{B|A}}[\F^{-1}] \cong
  \Omega^1_{\ul{B}|\ul{A}}[\F^{-1}] \\
  &\cong
  \Omega^1_{\ul{B}[\F^{-1}]|\ul{A}[(\varphi^{-1}[\F])^{-1}]} \cong
  \Omega^1_{\O_{\Spec(B)}|f^{-1}\O_{\Spec(A)}}.
\end{align*}
The first isomorphism is by
Proposition~\ref{prop:tilde-construction-internally}, the second because the
geometric construction ``taking Kähler differentials'' commutes with ``taking
constant sheaves'' (since it can be expressed as pullback along a geometric
morphism), the third by~\stacksproject{00RT}, and the fourth
by~Proposition~\ref{prop:tilde-construction-internally}.
\end{proof}

Incidentally, the Stacks Project too uses the construction via the free
module~\stacksproject{08RL}, however because the Stacks Project doesn't employ
the internal language they have to manually sheafify and keep track of open
subsets.

\subsection{Internal proofs of common lemmas}

\begin{lemma}Let~$X$ be a scheme (or a ringed space). Let
\[ 0 \lra \F \lra \G \lra \H \lra 0 \]
be a short exact sequence of~$\O_X$-modules. If~$\F$ and~$\H$ are of finite
type, so is~$\G$; similarly, if~$\F$ and~$\H$ are finite locally free, so
is~$\G$.
\end{lemma}
\begin{proof}From the internal perspective, we are given a short exact sequence
of modules with the outer ones being finitely generated (\resp finite free)
and we have to show that the middle one is finitely generated (\resp finite
free) as well. It is well-known that this follows; and since the usual proof of
this fact is intuitionistically valid, we are done.
\end{proof}

The proof works very generally, in the context of arbitrary ringed
spaces, and is still very simple. This is common to proofs using the internal
language. Particular features of schemes enter only at clearly recognizable
points, for instance when an internal property specific to the structure sheaf
of schemes is used (such as in Proposition~\ref{prop:neginvnilpotent}).

\begin{lemma}\label{lemma:coherent-stuff}
Let~$X$ be a scheme (or a ringed space).
\begin{enumerate}
\item Let~$0 \to \F \to \G \to \H \to 0$ be an exact sequence
of~$\O_X$-modules. If two of the three modules are coherent, so is the third.
\item Let~$\F \to \G$ be a morphism of~$\O_X$-modules such that~$\F$ is
of finite type and~$\G$ is coherent. Then its kernel is of finite type as well.
\item If~$\F$ is a finitely presented~$\O_X$-module and~$\G$ is a
coherent~$\O_X$-module, the~$\O_X$-modules~$\HOM_{\O_X}(\F,\G)$ and~$\F \otimes_{\O_X} \G$
are coherent as well.
\end{enumerate}
\end{lemma}
\begin{proof}These statements follow directly from interpreting the
corresponding standard proofs of commutative algebra in the internal language.
For those standard proofs, see for instance the lecture notes of Ravi
Vakil~\cite[Section~13.8]{vakil:foag}, where they are given as a series of
exercises.
\end{proof}

\begin{lemma}\label{lemma:kernel-of-epi-fingen}
Let~$X$ be a scheme (or a locally ringed space). Let~$\alpha : \G
\to \H$ be an epimorphism of finite locally free~$\O_X$-modules. Then the
kernel of~$\alpha$ is finite locally free as well.\end{lemma}
\begin{proof}It suffices to give an intuitionistic proof of the following
statement: The kernel of a matrix over a local ring, which as a linear map is
surjective, is finite free.

Let~$M \? R^{n \times m}$ be such a matrix. Since by the surjectivity
assumption some linear combination of the columns is~$e_1$ (the first canonical
basis vector), some linear combination of the entries of the first row of~$M$
is~$1$. By locality of~$R$, at least one entry of the first row is invertible.
By applying appropriate column and row transformations, we may therefore assume that~$M$
is of the form
\[ \left(
  \begin{array}{c|ccc}
    1 & 0 & \cdots & 0 \\ \hline
    0 & \multicolumn{3}{c}{\multirow{3}{*}{\raisebox{-5mm}{\scalebox{1.2}{$\widetilde M$}}}} \\
    \raisebox{2pt}{\vdots} & & & \\
    0 & & &
  \end{array}
\right) \]
with the submatrix~$\widetilde M$ fulfilling the same condition as~$M$.
Continuing in this way, it follows that~$m \geq n$ and that we may assume
that~$M$ is of the form
\[ \left(
  \begin{array}{ccc|c}
    1 & & & \multirow{3}{*}{\raisebox{-3mm}{\scalebox{1.2}{\ 0}}} \\
    & \ddots && \\
    && 1
  \end{array}
\right)\!. \]
The kernel of such a matrix is obviously freely generated by the canonical
basis vectors corresponding to the zero columns. In particular, the rank of the
kernel is~$m-n$.
\end{proof}

\begin{rem}The internal language machinery gives no reason to believe that the
dual statement is true, \ie that the cokernel of a monomorphism of finite
locally free~$\O_X$-modules is finite locally free. This would follow from
an intuitionistic proof of the statement that the cokernel of an injective map
between finite free modules over a local ring is finite free. But this
statement is of course false (consider the exact sequence
$0 \lra \ZZ_{(2)} \stackrel{\cdot 2}{\lra} \ZZ_{(2)} \lra \FF_2 \lra 0$
of~$\ZZ_{(2)}$-modules).
\end{rem}

\begin{lemma}\label{lemma:epi-iso}Let~$X$ be a scheme (or a locally ringed space). Let~$\alpha : \G
\to \H$ be an epimorphism of finite locally free~$\O_X$-modules of the same
rank. Then~$\alpha$ is an isomorphism.\end{lemma}
\begin{proof}It suffices to give an intuitionistic proof of the following
statement: A square matrix over a local ring, which as a linear map is
surjective, is invertible.

This follows from the proof of the previous lemma, since it shows that the
kernel of such a matrix is finite free of rank zero.
\end{proof}

\begin{rem}The conclusion of Lemma~\ref{lemma:epi-iso} also holds if~$X$ is
only assumed to be a ringed space. To show this, it suffices to give an
intuitionistic proof of the following statement: A square matrix over a (not
necessarily local) ring, which as a linear map is surjective, is invertible.
Such a matrix~$A$ possesses a right inverse. Therefore~$\det A$ is invertible.
Thus~$A$ is invertible with inverse~$(\det A)^{-1} \cdot \operatorname{ad} A$.
\end{rem}

\begin{lemma}\label{lemma:mono-iso}Let~$X$ be a scheme (or a locally ringed space).
Let~$\alpha : \G \to \H$ be a morphism of finite locally
free~$\O_X$-modules. If~$\alpha$ is a monomorphism, then the rank of~$\G$ is
less than or equal to the rank of~$\H$.\end{lemma}

\begin{proof}It suffices to give an intuitionistic proof of the following
statement: A matrix over a nontrivial ring, which as a linear map is injective,
has a smaller or equal number of columns than rows. Such a proof can for instance be found
in~\cite[p.~1013]{richman:trivial-rings}.
\end{proof}

\begin{lemma}Let~$X$ be a scheme (or a ringed space). Let
$0 \to \F \to \G \to \H \to 0$ be a short exact sequence of~$\O_X$-modules.
Then for the closures of the supports there holds the equation
$\Clos \supp \G = \Clos \supp \F \cup \Clos \supp \H$.
\end{lemma}
\begin{proof}Switching to complements, we have to prove that
\[ \Int(X \setminus \supp\G) = \Int(X \setminus \supp\F) \cap \Int(X \setminus
\supp\H). \]
By Proposition~\ref{prop:characterization-support}, it suffices to prove
\[ \Sh(X) \models (\G = 0\ \Longleftrightarrow\ \F = 0 \wedge \H = 0); \]
this is a basic observation in linear algebra, valid intuitionistically.
\end{proof}
Of course, a stronger version of this lemma -- about the supports themselves
instead of their closures -- is easily proven without using the internal
language. We included this example only for illustrative purposes.

\begin{lemma}Let~$X$ be a scheme (or a locally ringed space). Let~$\L$ be a line
bundle on~$X$, \ie an~$\O_X$-module locally free of rank~1.
Let~$s_1,\ldots,s_n \in \Gamma(X,\L)$ be global sections. Then these sections
globally generate~$\F$ if and only if
\[ \Sh(X) \models \bigvee_i \speak{$\alpha(s_i)$ is invertible for some
isomorphism~$\alpha : \L \to \O_X$}. \]
\end{lemma}
\begin{proof}It suffices to give an intuitionistic proof of the following fact:
Let~$R$ be a local ring. Let~$L$ be a free~$R$-module of rank~1.
Let~$s_1,\ldots,s_n\?L$ be given elements. Then~$L$ is generated as
an~$R$-module by these elements if and only if for some~$i$, the image of~$s_i$
under some isomorphism~$L \to R$ is invertible.

The choice of such an isomorphism does not matter, since any two such
isomorphisms~$\alpha, \beta : L \to R$ differ by a unit of~$R$: $\alpha(x) =
\alpha(\beta^{-1}(1)) \cdot \beta(x)$ for any~$x\?L$,
and~$\alpha(\beta^{-1}(1)) \cdot \beta(\alpha^{-1}(1)) = 1$ in~$R$.

For the ``if'' direction, we have that some~$\alpha(s_i)$ is a generator
of~$R$. Since~$\alpha$ is an isomorphism, it follows that~$s_i$ generates~$L$,
and thus that in particular, the family~$s_1,\ldots,s_n$ generates~$L$.

For the ``only if'' direction, we have that the unit of~$R$ can be expressed as
a linear combination of the~$\alpha(s_i)$, where~$\alpha : L \to R$ is some
isomorphism (whose existence is assured by the assumption on the rank of~$L$).
Since~$R$ is a local ring, it follows that one of the summands and thus one of
the~$\alpha(s_i)$ is invertible.
\end{proof}

\begin{rem}The canonical ring homomorphism~$\O_{X,x}
\twoheadrightarrow k(x)$ is local. Therefore a germ in~$\O_{X,x}$ is invertible
if and only if its image in~$k(x)$ is not zero. From this one can conclude that
global sections~$s_1,\ldots,s_n \in \Gamma(X,\F)$ generate~$\F$ if and only if,
for any point~$x \in X$, the images~$s_i \in \F|_x$ in the fibers do not vanish
simultaneously.
\end{rem}

\begin{lemma}Let~$X$ be a scheme (or a ringed space). Let~$\L$ be a locally
free~$\O_X$-module of rank~1. Then~$\L^\vee \otimes_{\O_X} \L \cong \O_X$.\end{lemma}
\begin{proof}Recall that the dual is defined by~$\L^\vee \defeq
\HOM_{\O_X}(\L,\O_X)$. Since~``$\HOM$'' looks like~``$\Hom$'' from the internal
point of view, the dual sheaf~$\L^\vee$ looks just like the ordinary dual
module. However, to prove the claim, it does \emph{not} suffice to give an
intuitionistic proof of the following fact of linear algebra: ``Let~$L$ be a
free~$R$-module of rank~1. Then there exists an isomorphism~$L^\vee \otimes_R L
\to R$.'' Since the interpretation of ``$\exists$'' using the Kripke--Joyal
semantics is local existence, this would only show that~$\L^\vee \otimes_{\O_X}
\L$ is \emph{locally} isomorphic to~$\O_X$.

Instead, we have to actually \emph{write down} (\ie explicitly give) a
linear map in the internal language -- not using the assumption that~$L$ is
free of rank~1, as this would introduce an existential quantifier again (see
Section~\ref{sect:internal-constructions}).
So we have to prove the following fact: Let~$L$ be an~$R$-module. Then there
explicitly exists a linear map~$L^\vee \otimes_R L \to R$ such that this map is
an isomorphism if~$L$ is free of rank~1.

This is done as usual: Define~$\alpha : L^\vee \otimes_R L \to R$ by~$\lambda
\otimes x \mapsto \lambda(x)$. Since~$L$ is free of rank~1, there is an
isomorphism~$L \cong R$. Precomposing~$\alpha$
with the induced isomorphism~$R^\vee \otimes_R R \to L^\vee \otimes_R L$,
we obtain the linear map~$R^\vee \otimes_R R \to R$ given by the same
term: $\lambda \otimes x \mapsto \lambda(x)$. One can check that an inverse is given
by~$x \mapsto \id_R \otimes x$.
\end{proof}

\begin{lemma}\label{lemma:torsion-stuff}
Let~$X$ be a scheme (or a ringed space). Let~$\F$ be
an~$\O_X$-module.
\begin{enumerate}
\item Assume~$X$ to be a locally Noetherian scheme. Then $\F$ is torsion-free
(meaning~$\F_\tors = 0$) if and only if all stalks~$\F_x$ are torsion-free.
\item The quotient sheaf~$\F/\F_\tors$ is torsion-free and the torsion
submodule~$\F_\tors$ is a torsion module.
\item The dual sheaf $\F^\vee$ is torsion-free.
\item If~$\F$ is reflexive (meaning that the canonical morphism~$\F \to
\F^{\vee\vee}$ is an isomorphism), then~$\F$ is torsion-free.
\item If~$\F$ is finite locally free, then~$\F$ is reflexive.
\item Assume~$X$ to be a Dedekind scheme and~$\F$ to be of finite presentation.
If~$\F$ is torsion-free, then~$\F$ is finite locally free.
\end{enumerate}
\end{lemma}
\begin{proof}The first statement follows from the observation that~$(\F_\tors)_x
= (\F_x)_\tors$ (Proposition~\ref{prop:torsion-submodule-stalks}). All the
others follow simply by interpreting the corresponding facts of linear algebra
in the internal universe. For concreteness, we give intuitionistic proofs of
the last three statements.

So let~$M$ be a reflexive~$R$-module. We have to show that~$M$ is
torsion-free. To this end, let an element~$x \? M$ and a regular element~$a \?
R$ such that~$ax = 0$ be given. For any~$\vartheta \? M^\vee$, it follows
that~$\vartheta(x) = 0$, since~$a \vartheta(x) = \vartheta(ax) = \vartheta(0) =
0$ and~$a$ is regular. Thus the image of~$x$ under the canonical map~$M \to
M^{\vee\vee}$ is zero. By reflexivity, this implies that~$x$ is zero.

For statement~(5), we have to prove that~$R$-modules of the form~$R^n$ are
reflexive. This is obvious, the required inverse map is~$(R^n)^{\vee\vee} \to
R^n,\ \lambda \mapsto \sum_i \lambda(\vartheta_i)$ where~$\vartheta_i : R^n \to
R,\ (x_j)_j \mapsto x_i$.

In view of Corollary~\ref{cor:dedekind-smith} we can put matrices over~$\O_X$
into Smith canonical form if~$X$ is a Dedekind scheme. Therefore it suffices
to give an intuitionistic proof of the following fact: Let~$R$ be an integral
domain in the strong sense such that matrices over~$R$ can be put into
Smith canonical form. Let~$M$ be a finitely presented torsion-free~$R$-module.
Then~$M$ is finite free.

Such a proof can proceed as follows: Since~$M$ is finitely presented, it is the cokernel of
some matrix. Without loss of generality, we may assume that the presentation matrix is diagonal,
so~$M$ is isomorphic to some finite direct sum~$\bigoplus_i R/(a_i)$.
Since~$M$ is torsion-free, all the summands~$R/(a_i)$ are torsion-free as well.
Since~$R$ is an integral domain in the strong sense,
the~$a_i$ are either zero or invertible. Thus~$R/(a_i)$ is isomorphic to~$R$ or
to the zero module. In any case,~$R/(a_i)$ is finite free and therefore~$M$
is finite free as well.
\end{proof}

\section{Upper semicontinuous functions}
\label{sect:upper-semicontinuous-functions}

\subsection{Interlude on natural numbers}
In classical logic, the natural numbers are complete in the sense that any
inhabited set of natural numbers possesses a minimal element. This statement
cannot be proven intuitionistically -- intuitively, this is because one cannot
explicitly pinpoint the (classically existing) minimal element of an arbitrary
inhabited set;\footnote{Let~$\varphi$ be an arbitrary formula. Assuming that
any inhabited subset of the natural numbers possesses a minimal element, we
want to show that~$\varphi \vee \neg\varphi$. Define the subset $U \defeq \{ n \in
\NN \,|\, (n = 1) \vee \varphi \} \subseteq \NN$, which surely is inhabited by~$1
\in U$. So by assumption, there exists a number~$z \in \NN$ which is the
minimum of~$U$. We have $z = 0$ or $z > 0$. If~$z = 0$, we have~$0 \in U$,
so~$(0 = 1) \vee \varphi$, so~$\varphi$ holds.  If~$z > 0$, then~$\neg\varphi$
holds: If~$\varphi$ were true, zero would be an element of~$U$, contradicting
the minimality of~$z$.} see below for a sheaf-theoretic interpretation.

In intuitionistic logic, the completeness principle can be salvaged in two
essentially different ways: either by strengthening the premise, or by
weakening the conclusion.

\begin{lemma}\label{lemma:minimum-subset-naturals}
Let~$U \subseteq \NN$ be an inhabited subset of the natural
numbers.
\begin{enumerate}
\item Assume~$U$ to be \emph{detachable}, \ie assume that for any natural
number~$n$, either~$n \in U$ or~$n \not\in U$. Then~$U$ possesses a minimal
element.
\item In any case,~$U$ does \notnot possess a minimal element.
\end{enumerate}
\end{lemma}
\begin{proof}
The first statement can be proven by induction on the witness of inhabitation,
\ie the given number~$n$ such that~$n \in U$. We omit further details, since we will
not need this statement in our applications.

For the second statement, we give a careful proof since logical subtleties matter. To simplify the
exposition, we assume that~$U$ is upward-closed, \ie that any number
larger than some element of~$U$ lies in~$U$ as well. Any subset can be closed
in this way (by considering~$\{ n \in \NN \,|\, \exists m \in U\_ n \geq m \}$)
and a minimal element of the closure will be a minimal element for~$U$ as well.

We induct on the number~$n \in U$ given by the assumption that~$U$ is
inhabited. In the case~$n = 0$ we are done since~$0$ is a minimal element
of~$U$. For the induction step~$n \to n+1$, the intuitionistically valid double
negation of the law of excluded middle
gives
\[ \neg\neg(n \in U \vee n \not\in U). \]
Because of the tautologies $(\varphi \Rightarrow \psi) \Rightarrow
(\neg\neg\varphi \Rightarrow \neg\neg\psi)$ and~$\neg\neg\neg\neg\varphi \Rightarrow
\neg\neg\varphi$ (see Section~\ref{sect:appreciating-intuitionistic-logic}), it
suffices to show that~$n \in U \vee n \not\in U$ implies the conclusion.
So assume~$n \in U \vee n \not\in U$.
If~$n \in U$, then~$U$ does \notnot possess a minimal element by the induction
hypothesis. If~$n \not\in U$, then~$n+1$ is a minimal element (and so, in
particular,~$U$ does \notnot possess a minimal element): If~$m$ is
any element of~$U$, we have~$m \geq n+1$ or~$m \leq n$. In the first case,
we're done. In the second case, it follows that~$n \in U$ because~$U$ is
upward-closed and so we obtain a contradiction. From this contradiction we
can trivially deduce~$m \geq n+1$ as well. \qedhere
\end{proof}

If we want to work with a complete partially ordered set (poset) of natural numbers in intuitionistic
logic, we have to construct a suitable completion. The idea of the following
definition is to encode numbers as the (not necessarily existing) minimum of
inhabited upward-closed subsets.
\begin{defn}The \emph{completed poset of natural numbers} is
the set~$\widehat{\NN}$ of all inhabited upward-closed subsets of~$\NN$, ordered by
reverse inclusion. The elements of~$\widehat{\NN}$ are called \emph{generalized natural numbers}.\end{defn}
\begin{lemma}The completed poset of natural numbers is the least poset
containing~$\NN$ and possessing minima
of arbitrary inhabited subsets.\end{lemma}
\begin{proof}
The embedding $\NN \hookrightarrow \widehat\NN$ is given by
\[ n \in \NN \quad\longmapsto\quad {\uparrow}(n) \defeq \{ m \in \NN \,|\, m \geq n \}. \]
If~$M \subseteq \widehat\NN$ is an inhabited subset, its minimum is
\[ \min M = \bigcup M \in \widehat\NN. \]
We omit the proof of the universal property.
\end{proof}

\begin{rem}\label{rem:surjectivity-embedding}
In classical logic, the map~$\widehat\NN \to \NN,\ U \mapsto \min U$
is a well-defined isomorphism of partially ordered sets. In fact, it is the
inverse of the canonical embedding~$\NN \hookrightarrow \widehat\NN$. In
intuitionistic logic, this embedding is still injective, but it cannot be
shown to be surjective: It is only the case that any element of~$\widehat\NN$
does \notnot possess a preimage (by Lemma~\ref{lemma:minimum-subset-naturals}).
\end{rem}

\subsection{A geometric interpretation}
We are interested in the completed natural numbers for the following reason: A
generalized natural number in the topos of sheaves on a topological space~$X$ is
the same as an upper semicontinuous function~$X \to \NN$.

\begin{lemma}\label{lemma:upper-semicontinuous-functions}
Let~$X$ be a topological space. The sheaf~$\widehat\NN$ of
generalized natural numbers on~$X$ is canonically isomorphic to the sheaf of upper
semicontinuous~$\NN$-valued functions on~$X$.\end{lemma}
\begin{proof}
When referring to the natural numbers in the internal language, we actually
refer to the constant sheaf~$\ul{\NN}$ on~$X$. (This is because the
sheaf~$\ul{\NN}$ fulfills the axioms for a natural numbers object,
cf.\@~\cite[Section~VI.1]{moerdijk-maclane:sheaves-logic}.) Recall that its sections on an
open subset~$U \subseteq X$ are continuous functions~$U \to \NN$, where~$\NN$
is equipped with the discrete topology.

Therefore, a section of~$\widehat\NN$ on an open subset~$U \subseteq X$ is
given by a subsheaf~$\A \hookrightarrow \ul{\NN}|_U$ such that
\[ U \models \exists n\?\ul{\NN}\_ n \in \A
  \quad\text{and}\quad
  U \models \forall n,m\?\ul{\NN}\_ n \geq m \wedge n \in \A \Rightarrow m \in
  \A. \]
Since these conditions are geometric implications, they are satisfied if and only if any
stalk~$\A_x$ is an inhabited upward-closed subset of~$\ul{\NN}_x \cong \NN$.
The association
\[ x \in X \quad\longmapsto\quad \min\{ n \in \NN \,|\, n \in \A_x \} \]
thus defines a map~$X \to \NN$. This map is indeed upper semicontinuous, since
if~$n \in \A_x$, there exists an open neighborhood~$V$ of~$x$ such that the constant
function with value~$n$ is an element of~$\Gamma(V,\A)$ and therefore~$n \in
\A_y$ for all~$y \in V$.

Conversely, let~$\alpha : U \to \NN$ be an upper semicontinuous function. Then
\[ \text{$V \subseteq U$ open} \quad\longmapsto\quad \{ f : V \to \NN \,|\, \text{$f$
continuous,\ $f \geq \alpha$ on~$V$} \} \]
is a subobject of~$\ul{\NN}|_U$ which internally is inhabited and upward-closed.
Further details are left to the reader.
\end{proof}

Under the correspondence given by Lemma~\ref{lemma:upper-semicontinuous-functions}, locally \emph{constant}
functions map precisely to the (image of the) \emph{ordinary} internal natural numbers
(in the completed natural numbers).
In a similar vein, the sheaf given by the internal construction of
the set of \emph{all} upward-closed subsets of the natural numbers (not
only the inhabited ones) is canonically isomorphic to the sheaf of
upper semicontinuous functions with values in~$\NN \cup \{ +\infty
\}$.

The correspondence can be used to understand classical facts about
upper semicontinuous functions as features of intuitionistic number theory. For
instance, it is well-known that any upper semicontinuous~$\mathbb{N}$-valued
function on an arbitrary topological space is locally constant on a dense open subset.
This can be explained as follows: The generalized natural number associated to such a
function is \notnot an ordinary natural number from the internal point of view.
Since ``not not'' translates to ``holding on a dense open subset''
(Proposition~\ref{prop:modops-kripke}), it follows that there is a dense open
subset on which the function corresponds to an ordinary internal natural
number, \ie is locally constant.

\subsection{The upper semicontinuous rank function}
Recall that the rank of an~$\O_X$-module~$\F$ on a scheme~$X$ (or
locally ringed space) at a point~$x \in X$ is defined as the~$k(x)$-dimension
of the vector space~$\F_x \otimes_{\O_{X,x}} k(x)$. If we assume that~$\F$ is
of finite type around~$x$, this dimension is finite and equals the minimal
number of elements needed to generate~$\F_x$ as an~$\O_{X,x}$-module by
Nakayama's lemma.

In the internal language, we can define an element of~$\widehat\NN$ by
\begin{multline*}
  \rank\F \defeq \min\{ n \in \NN \,| \\
  \speak{there is a gen.\@ family for~$\F$ consisting of~$n$ elements} \} \in \widehat\NN.
\end{multline*}
If the module~$\F$ is finite locally free, it will be a finite free module from the
internal point of view and the rank defined in this way will be an
actual natural number (see below); but in general, the rank is really an element of the
completion.

\begin{prop}\label{prop:rank-function-internally}
Let~$\F$ be an~$\O_X$-module of finite type on a scheme~$X$ (or a locally ringed
space). Under the correspondence given by Lemma~\ref{lemma:upper-semicontinuous-functions}, the internally
defined rank maps to the rank function of~$\F$.
\end{prop}
\begin{proof}
We have to show that for any point~$x \in X$ and natural number~$n$, there
exists a generating family for~$\F_x$ consisting of~$n$
elements if and only if there exists an open neighborhood~$U$ of~$x$ such that
\[ U \models \speak{there exists a generating family
for~$\F$ consisting of~$n$ elements}. \]
The ``if'' direction is obvious. For the ``only if'' direction, consider
(liftings to local sections of a)
generating family~$s_1,\ldots,s_n$ of~$\F_x$. Since~$\F$ is of finite type,
there also exist sections~$t_1,\ldots,t_m$ on some neighborhood~$V$ of~$x$ which
generate any stalk~$\F_y$, $y \in V$. Since the~$t_i$ can be expressed as a
linear combination of the~$s_j$ in~$\F_x$, the same is true on some open
neighborhood~$U \subseteq V$ of~$x$. On this neighborhood, the~$s_j$ generate
any stalk~$\F_y$, $y \in U$, so by geometricity we have
\[ U \models \speak{$s_1,\ldots,s_n$ generate~$\F$}. \qedhere \]
\end{proof}
Once we understand when properties spread from points to
neighborhoods in logical terms, we will be able to give a simpler proof of Proposition~\ref{prop:rank-function-internally} (see
Lemma~\ref{lemma:gen-family-n}).

\begin{lemma}\label{lemma:rank-functor-locally-constant}
Let~$X$ be a scheme (or a locally ringed space). Let~$\F$ be an~$\O_X$-module of
finite type. If~$\F$ is finite locally free, its rank function is locally
constant. The converse holds if~$X$ is a reduced scheme.
\end{lemma}
\begin{proof}The rank function is locally constant if and only if internally,
the rank of~$\F$ is an actual natural number. Also recall that the structure
sheaf fulfills a certain field condition if~$X$ is a reduced
scheme (Corollary~\ref{cor:field-reduced}). Therefore it suffices to give a
proof of the following fact of intuitionistic linear algebra: Let~$R$ be a
local ring. Let~$M$ be a finitely generated~$R$-module. If~$M$ is finite
free, its rank is an actual natural number. The converse holds if~$R$ fulfills
the field condition that any element which is not invertible is zero.

So assume that such a module~$M$ is finite free. Then it is isomorphic
to~$R^n$ for some actual natural number~$n$; by the internal proof in
Lemma~\ref{lemma:kernel-of-epi-fingen}, the rank of~$M$ is therefore this
number~$n$ (for any surjection~$R^m \twoheadrightarrow R^n$ it holds that~$m
\geq n$).

Conversely, assume that the rank of~$M$ is an actual natural number. Then
there exists a minimal generating family~$x_1,\ldots,x_n\?M$. We can verify that this family is
indeed linearly independent (and thus a basis, demonstrating that~$M$ is finite
free): Let~$\sum_i a_i x_i = 0$ with~$a_i\?R$. If any~$a_i$ were
invertible, the family~$x_1,\ldots,x_{i-1},x_{i+1},\ldots,x_n$ would too
generate~$M$, contradicting the minimality. So each~$a_i$ is not invertible.
By the field property of~$R$, each~$a_i$ is zero.
\end{proof}

\begin{lemma}\label{lemma:locally-free-dense}
Let~$X$ be a reduced scheme. Let~$\F$ be an~$\O_X$-module of
finite type. Then~$\F$ is finite locally free on a dense open subset.\end{lemma}
\begin{proof}Since ``dense open'' translates to ``not not'' in the internal
language (Proposition~\ref{prop:modops-kripke}), it suffices to give an
intuitionistic proof of the following fact: Let~$R$ be a local ring which fulfills an
appropriate field condition. Let~$M$ be a finitely generated~$R$-module.
Then~$M$ is \notnot finite free.

By Remark~\ref{rem:surjectivity-embedding}, the rank of such a module~$M$ is
\notnot an actual natural number. By the last part of the
previous proof, it thus follows that~$M$ is \notnot finite free.
\end{proof}

\begin{rem}Besides basics on natural numbers in an intuitionistic
setting and some dictionary terms (``reduced'', ``finite locally free'',
``finite type``, ``dense open''), this proof does not depend on any further
tools. In particular, it doesn't depend on Nakayama's lemma or on facts about
semicontinuous functions. For the (more complex) standard proof of this fact, see for
instance~\cite{vakil:foag}, where the claim is dubbed an ``important hard
exercise'' (Exercise~13.7.K).\end{rem}

\subsection{The upper semicontinuous dimension function} Recall that the
dimension of a topological space~$X$ at a point~$x \in X$ is defined as the
infimum
\[ \dim_x X \defeq \inf\{ \dim U \,|\, \text{$U$ open neighborhood of~$x$} \}. \]

The map~$X \to \NN \cup \{+\infty\},\ x \mapsto \dim_x X$ is upper
semicontinuous and thus corresponds to an internal generalized (possibly
unbounded) natural number. The following proposition shows that this number has
an explicit description.

\begin{prop}Let~$X$ be a scheme. Then the upper semicontinuous function
associated to the internal number ``Krull dimension of~$\O_X\!$'' is the
dimension function~$x \mapsto \dim_x X$.\end{prop}
\begin{proof}Internally, we define the Krull dimension of~$\O_X$ as the infimum
over all natural numbers~$n$ such that~$\O_X$ is of Krull
dimension~$\leq n$. This infimum need not exist in the natural numbers, of
course; so we really mean the upward-closed set~$\A$ of all those numbers. (It
is inhabited if and only if, from the external perspective, the dimension
of~$X$ is locally finite. In this case, it defines a generalized natural number.)

We thus have to show for any point~$x \in X$:
\[ \inf\{ n \in \NN \cup \{+\infty\} \,|\, n \in \A_x \} =
  \dim_x X. \]
The condition on~$n$ can be expressed as follows, where we write~``$\ul{n}$''
to denote the constant function with value~$n$:
\begin{align*}
  &\ n \in \A_x \\
  \Longleftrightarrow &\
  \text{for some open neighborhood~$U$ of~$x$, $\ul{n} \in \Gamma(U,\A)$} \\
  \Longleftrightarrow &\
  \text{for some open neighborhood~$U$ of~$x$,} \\
  & \qquad\qquad U \models \speak{$\O_X$ is of Krull dimension~$\leq n$} \\
  \Longleftrightarrow &\
  \text{for some open neighborhood~$U$ of~$x$,} \\
  & \qquad\qquad \dim U \leq n
\end{align*}
The second equivalence follows from the external description of internally-defined
subsheaves given in Section~\ref{sect:internal-constructions}.
We thus have:
\[
  \inf\{ n \,|\, n \in \A_x \} =
    \inf\{ \dim U \,|\, \text{$U$ open neighborhood of~$x$} \}
  = \dim_x X. \qedhere
\]
\end{proof}

\section{Modalities}
\label{sect:modalities}

Philosophers and logicians do not only study what is \emph{true}, but also what
is \emph{known}, what is \emph{believed}, what is \emph{possible}, and so on.
Such \emph{modalities} are absent from the usual mathematical practice.
However, it turns out that a specific kind of such modalities plays a role in
understanding when properties spread from points to neighborhoods.

Briefly, this is because for any point~$x$ of a topological space~$X$, there
exists a modal operator~$\Box$ such that for any formula~$\varphi$ of the
internal language of the sheaf topos~$\Sh(X)$, the internal
statement~$\Box\varphi$ means that~$\varphi$ holds on some open neighborhood
of the given point~$x$. In this way, we can reduce sheaf-theoretic questions to
questions of modal intuitionistic (non-sheafy) mathematics.

The techniques developed in this section also enable us to use the internal
language of~$\Sh(X)$ to talk about sheaves on \emph{subspaces} of~$X$ (and more
general \emph{sublocales} of~$X$).

Topological interpretations of modal logic were studied before, for instance by
Awodey and Kishida~\cite{awodey-kishida:modal}. However, they study a
different kind of modal operators, not corresponding to the Lawvere--Tierney
topologies of topos theory, and pursue different goals.

\subsection{Basics on truth values and modal operators}
\label{sect:basics-on-truth-values}

\begin{defn}The \emph{set of truth values~$\Omega$} is the powerset of the
singleton set~$1 \defeq \{\star\}$, where~$\star$ is a formal symbol.\end{defn}

In classical logic, any subset of~$\{\star\}$ is either empty or inhabited, so
that~$\Omega$ contains exactly two elements, the empty set (``false'')
and~$\{\star\}$ (``true''). But
in intuitionistic logic, this cannot be shown; indeed, if we interpret the
definition in the topos of sheaves on a space~$X$, we obtain a (rather large) sheaf~$\Omega$
with
\[ \text{$U \subseteq X$ open} \quad\longmapsto\quad \Gamma(U,\Omega) = \{ V \subseteq U \,|\, \text{$V$
open} \}. \]
(This is because by definition of~$\Omega$ as the power object of the terminal
sheaf~$1$, sections of~$\Omega$ on an open subset~$U$ correspond to
subsheaves~$\F \hookrightarrow 1|_U$, and those are given by the greatest open
subset~$V \subseteq U$ such that~$\Gamma(V,\F)$ is inhabited.)
Obviously, in general, this sheaf has many sections, in particular more than
the binary coproduct~$1 \amalg 1$ (unless any open subset of~$X$ is also
closed).

The \emph{truth value} of a formula~$\varphi$ is by definition the subset
$\{ x \? 1 \,|\, \varphi \} \in \Omega$, where~``$x$'' is a fresh variable not
appearing in~$\varphi$. This subset is inhabited if and only
if~$\varphi$ holds and is empty if and only if~$\neg\varphi$ holds.
Conversely, we can associate to a subset~$F \subseteq 1$ the
proposition~$\speak{$F$ is inhabited}$.

By the above description of~$\Omega$ in
a sheaf topos~$\Sh(X)$, the interpretation of the truth value
of a formula~$\varphi$ in the internal language of~$\Sh(X)$ is a certain open
subset of~$X$. Tracing the definitions, we see that this open subset is
precisely the largest open subset on which~$\varphi$ holds, \ie the union of
all open subsets~$U \subseteq X$ such that~$U \models \varphi$.

Under the correspondence of formulas with truth values, logical operations
like~$\wedge$ and~$\vee$ map to set-theoretic operations like~$\cap$ and~$\cup$
-- for instance, we have
\[ \{ x \? 1 \,|\, \varphi \} \cap \{ x \? 1 \,|\, \psi \} =
  \{ x \? 1 \,|\, \varphi \wedge \psi \}. \]
This justifies a certain abuse of notation: We will sometimes treat elements
of~$\Omega$ as propositions and use logical instead of set-theoretic
connectives. In particular, if~$\varphi$ and~$\psi$ are elements of~$\Omega$,
we will write~``$\varphi \Rightarrow \psi$'' to mean~$\varphi \subseteq \psi$;
``$\bot$'' to mean~$\emptyset$; and~``$\top$'' to mean~$1$.

\begin{defn}\label{defn:modal-operator}
A \emph{modal operator} (or \emph{Lawvere--Tierney topology}) is a map~$\Box : \Omega \to \Omega$ such
that for all~$\varphi, \psi \in \Omega$,
\begin{enumerate}
\item $\varphi \Longrightarrow \Box\varphi$,
\item $\Box\Box\varphi \Longrightarrow \Box\varphi$,
\item $\Box(\varphi \wedge \psi) \Longleftrightarrow \Box\varphi \wedge \Box\psi$.
\end{enumerate}
Syntactically, the symbol~``$\Box$'' binds stronger than any other logical
connective. For instance, axiom~(2) is shorthand for~``$(\Box(\Box(\varphi)))
\Rightarrow (\Box(\varphi))$'' and axiom~(3) is shorthand for~``$(\Box(\varphi
\wedge \psi)) \Leftrightarrow ((\Box(\varphi)) \wedge (\Box(\psi)))$''.
\end{defn}

The intuition is that~$\Box\varphi$ is a certain weakening of~$\varphi$, where
the precise meaning of ``weaker'' depends on the modal operator. By the second
axiom, weakening twice is the same as weakening once.

In classical logic, where~$\Omega = \{ \bot, \top \}$, there are only two modal
operators: the identity map and the constant map with value~$\top$.
Both of these are not very interesting: The identity operator does not weaken
propositions at all, while the constant operator weakens every proposition to
the trivial statement~$\top$.

In intuitionistic logic, there can potentially exist further modal operators.
For applications to algebraic geometry, the following four operators will have
a clear geometric meaning and be of particular importance:
\begin{enumerate}
\item $\Box\varphi \defequiv (\alpha \Rightarrow \varphi)$, where~$\alpha$ is a
fixed proposition.
\item $\Box\varphi \defequiv (\varphi \vee \alpha)$, where~$\alpha$ is a
fixed proposition.
\item $\Box\varphi \defequiv \neg\neg\varphi$ (the \emph{double negation
modality}).
\item $\Box\varphi \defequiv ((\varphi \Rightarrow \alpha) \Rightarrow \alpha)$,
where~$\alpha$ is a fixed proposition.
\end{enumerate}

\begin{lemma}Any modal operator~$\Box$ is monotonic, \ie if~$\varphi
\Rightarrow \psi$, then~$\Box\varphi \Rightarrow \Box\psi$. Furthermore,
a modus ponens rule holds: If~$\Box\varphi$ holds, and if~$\varphi$
implies~$\Box\psi$, then~$\Box\psi$ holds as well.\end{lemma}
\begin{proof}Assume~$\varphi \Rightarrow \psi$. This is equivalent to
supposing~$\varphi \wedge \psi \Leftrightarrow \varphi$. We are to show
that~$\Box\varphi \Rightarrow \Box\psi$, \ie that~$\Box\varphi \wedge
\Box\psi \Leftrightarrow \Box\varphi$. This follows since by the third
axiom on a modal operator, we have~$\Box\varphi \wedge \Box\psi \Leftrightarrow
\Box(\varphi \wedge \psi)$, and~$\Box$ respects equivalence of propositions.

For the second statement, consider that if~$\varphi \Rightarrow \Box\psi$, by
monotonicity and the second axiom on a modal operator it follows
that~$\Box\varphi \Rightarrow \Box\Box\psi \Rightarrow \Box\psi$.
\end{proof}

\label{proof-scheme-boxed-statements}The modus ponens rule justifies the
following proof scheme: When trying to show, given that some boxed
statement~$\Box\varphi$ holds, that some further boxed statement~$\Box\psi$
holds, we may give a proof of~$\Box\psi$ under the stronger
assumption~$\varphi$. Symbolically:
\[ (\Box\varphi \Rightarrow \Box\psi) \Longleftrightarrow
  (\varphi \Rightarrow \Box\psi). \]

\begin{rem}There is some contention on what symbol one should use for modal
operators in the sense of Definition~\ref{defn:modal-operator}. This is
because, in the modal logic community, the symbol~``$\Box$'' usually refers to
the modal operator~``it is necessary that''. For this modal operator, one often
imposes the \emph{reflexivity axiom}~$\Box\varphi \Rightarrow \varphi$ which we
don't impose (and which would trivialize the theory). Conversely, our
axiom~$\varphi \Rightarrow \Box\varphi$ isn't meaningful in the necessity
interpretation. This axiom is meaningful for the modal operator~``it is
possible that'', commonly denoted~``$\Diamond$''; but for this modal operator, the
axiom~$\Diamond(\varphi \wedge \psi) \Leftrightarrow \Diamond\varphi \wedge
\Diamond\psi$ isn't meaningful.

A classical modal operator which matches our axioms is ``I believe that''
under the proviso of ultimate knowledge (``I believe every true
statement'').

Goldblatt chooses the symbol~``$\nabla$'' for the modal operators in our
sense~\cite[Section~14.5]{goldblatt:topoi},~\cite{goldblatt:modality}. The
symbol~``$\varbigcirc$'' is also common, particularly in the hardware
verification community. A discussion of the relationship between these three
kinds of operators is contained in~\cite{pfenning:davies:reconstruction}. We
are grateful to Tadeusz Litak for valuable comments and references pertaining
to this topic.
\end{rem}

\subsection{Geometric meaning}\label{sect:modalities-geometric-meaning}
Let~$X$ be a topological space. As discussed
above, an open subset~$U \subseteq X$ defines an internal truth value (a global
section of the sheaf~$\Omega$). We also denote it by~``$U$'', such that
\[ V \models U \quad\Longleftrightarrow\quad V \subseteq U \]
for any open subset~$V \subseteq X$. (Shortcutting the various intermediate
steps, this can also be taken as a definition of~``$V \models U$''.)
If~$A \subseteq X$ is a closed subset, there is thus an internal truth
value~$A^c$ corresponding to the open subset~$A^c = X \setminus A$. If~$x \in
X$ is a point, we define~``$\notat{x}$'' to denote the truth value
corresponding to~$\Int(X \setminus \{x\})$, such that
\[ V \models \notat{x} \quad\Longleftrightarrow\quad V \subseteq \Int(X
\setminus \{ x \}) \quad\Longleftrightarrow\quad x \not\in V. \]

\begin{prop}\label{prop:modops-kripke}
Let~$U \subseteq X$ be a fixed open and~$A \subseteq X$ be a fixed
closed subset. Let~$x \in X$. Then, for any open subset~$V \subseteq X$, it
holds that:
\[ \renewcommand{\arraystretch}{1.3}\begin{array}{@{}lcl@{}}
  V \models (U \Rightarrow \varphi) &\Longleftrightarrow&
    V \cap U \models \varphi. \\[0.3em]
  V \models (\varphi \vee A^c) &\Longleftrightarrow&
    \textnormal{there is an open subset~$W \subseteq V$} \\
  && \quad\quad \textnormal{containing~$A \cap V$ such that $W \models \varphi$.} \\[0.3em]
  V \models \neg\neg\varphi &\Longleftrightarrow&
    \textnormal{there is a dense open subset~$W \subseteq V$ s.\,th.\@ $W \models
    \varphi$.} \\[0.3em]
  V \models ((\varphi \Rightarrow \notat{x}) \Rightarrow \notat{x}) &\Longleftrightarrow&
    \textnormal{$x \not\in V$ or there is an open neighborhood~$W \subseteq V$} \\
  && \quad\quad \textnormal{of~$x$ such that $W \models \varphi$.}
\end{array} \]
\end{prop}
\begin{proof}
\begin{enumerate}
\item Omitted.

\item Let~$V \models \varphi \vee A^c$. Then there exists an open covering~$V =
\bigcup_i V_i$ such that for each~$i$, $V_i \models \varphi$ or $V_i \subseteq
A^c$. Let~$W \subseteq V$ be the union of those~$V_i$ such that~$V_i \models \varphi$.
Then~$W \models \varphi$ by the locality of the internal language and~$A \cap V
\subseteq W$.

Conversely, let~$W \subseteq V$ be an open subset containing~$A \cap V$ such
that~$W \models \varphi$. Then~$V = W \cup (V \cap A^c)$ is an open covering
attesting~$V \models \varphi \vee A^c$.

\item For the ``only if'' direction, let~$W \subseteq V$ be the largest
open subset on which~$\varphi$ holds, \ie the union of all open subsets
of~$V$ on which~$\varphi$ holds. For the ``if'' direction, we may assume that
the given set~$W$ is also the largest open subset on which~$\varphi$ holds (by
enlarging~$W$ if necessary). The claim then follows by the following chain of
equivalences:
\begin{align*}
  &\ V \models \neg\neg\varphi \\
  \Longleftrightarrow&\ \forall \text{$Y \subseteq V$ open}\_
    \Bigl(\forall \text{$Z \subseteq Y$ open}\_ (Z \models \varphi) \Rightarrow Z
    = \emptyset\Bigr) \Longrightarrow Y = \emptyset \\
  \Longleftrightarrow&\ \forall \text{$Y \subseteq V$ open}\_
    \Bigl(\forall \text{$Z \subseteq Y$ open}\_ Z \subseteq W \Rightarrow Z
    = \emptyset\Bigr) \Longrightarrow Y = \emptyset \\
  \Longleftrightarrow&\ \forall \text{$Y \subseteq V$ open}\_
    Y \cap W = \emptyset \Longrightarrow Y = \emptyset \\
  \Longleftrightarrow&\ \text{$W$ is dense in~$V$.}
\end{align*}

\item Straightforward, since the interpretation of the internal statement with
the Kripke--Joyal semantics is
\[ \forall \text{$Y \subseteq V$ open}\_
  \Bigl(\forall \text{$Z \subseteq Y$ open}\_
    Z \models \varphi \Rightarrow x \not\in Z\Bigr) \Longrightarrow x \not\in
    Y. \qedhere \]
\end{enumerate}
\end{proof}

\subsection{The subspace associated to a modal operator}
\label{sect:subspace-to-modal-operator}
Any modal operator~$\Box : \Omega \to \Omega$ in the sheaf topos of~$X$ induces
on global sections a map
\[ j : \Open(X) \to \Open(X), \]
where~$\Open(X) = \Gamma(X,\Omega)$ is the set of open subsets of~$X$.
Explicitly, it is given by
\begin{align*}
  j(U) &= \text{largest open subset of~$X$ on which~$\Box U$ holds} \\
  &= \bigcup\ \{ V \subseteq X \ |\ \text{$V$ open},\ V \models \Box U \}.
\end{align*}
By the axioms for a modal operator, the map~$j$ fulfills similar such axioms: For any open
subsets~$U, V \subseteq X$,
\begin{enumerate}
\item $U \subseteq j(U)$,
\item $j(j(U)) \subseteq j(U)$,
\item $j(U \cap V) = j(U) \cap j(V)$.
\end{enumerate}
Such a map is called a \emph{nucleus} on~$\Open(X)$. Table~\ref{table:nuclei}
lists the nuclei associated to the four modal operators
of Proposition~\ref{prop:modops-kripke}.

\begin{table}
  \centering
  \renewcommand{\arraystretch}{1.3}
  \begin{tabular}{llll}
    \toprule
    Modal operator & associated nucleus &
      $j(V) = X$ iff \ldots &
      subspace \\\midrule
    $\Box\varphi \equiv (U \Rightarrow \varphi)$ &
      $j(V) = \Int(U^c \cup V)$ & $U \subseteq V$ & $U$ \\
    $\Box\varphi \equiv (\varphi \vee A^c)$ &
      $j(V) = V \cup A^c$ & $A \subseteq V$ & $A$ \\
    $\Box\varphi \equiv \neg\neg\varphi$ &
      $j(V) = \Int(\Clos(V))$ & $V$ is dense in $X$ &
      \multicolumn{1}{p{1cm}}{smallest dense sublocale of~$X$} \\
    $\Box\varphi \equiv ((\varphi \Rightarrow \notat{x}) \Rightarrow \notat{x})$ &
      $\begin{array}{@{}ll@{}}
        j(V) = X \setminus \Clos\{x\}, & \text{if $x \not\in V$} \\
        j(V) = X, & \text{if $x \in V$}
      \end{array}$ &
      $x \in V$ & $\{x\}$ \\
    \bottomrule
  \end{tabular}
  \vspace{0.5em}

  \caption{\label{table:nuclei}List of important modal operators and their
  associated nuclei (notation as in Proposition~\ref{prop:modops-kripke}).}
\end{table}

Any nucleus~$j$ defines a subspace~$X_j$ of~$X$, to be described below, with a small caveat: In
general, the subspace~$X_j$ cannot be realized as a topological subspace, but
only as a so-called \emph{sublocale}; the notion of a locale is a slight
generalization of the notion of a topological space, in which an underlying set
of points is not part of the definition. Instead, a locale is simply given by a
frame (a partially ordered set with certain properties) of arbitrary \emph{opens} satisfying some axioms -- these opens may, but do not necessarily have to,
be sets of points. Sheaf theory carries over to locales essentially unchanged,
since the notions of presheaves and sheaves only refer to open sets and coverings,
but not points.

Accessible introductions to the theory of locales include two notes by
Johnstone~\cite{johnstone:art,johnstone:point}. Locales are also well-known for
a curious application in the theory of
randomness~\cite{simpson:random1,simpson:random2}.

\begin{defn}\label{defn:subspace-by-nucleus}Let~$j$ be a nucleus on~$\Open(X)$.
Then the sublocale~$X_j$ of~$X$ is given by the frame of opens
$\Open(X_j) \defeq \{ U \in \Open(X) \,|\, j(U) = U \}$.
\end{defn}
If~$j$ is induced by a modal operator~$\Box$, we also write~``$X_\Box$''
for~$X_j$. In three of the four cases listed in Table~\ref{table:nuclei}, the
sublocale~$X_\Box$ can indeed be realized as a topological subspace. The only
exception is the sublocale~$X_{\neg\neg}$ associated to the double negation
modality. It can also be described as the \emph{smallest dense sublocale}
of~$X$; this is obviously a genuine locale-theoretic notion, since there
is (in general) no smallest dense topological subspace
(consider~$\RR$ and its dense subsets~$\QQ$ and~$\RR \setminus \QQ$).

The inclusion~$i : X_j \hookrightarrow X$ cannot in general be described on the
level of points, since~$X_j$ might not be realizable as a topological subspace.
But for sheaf-theoretic purposes, it suffices to describe~$i$ on the level of
opens. This is done as follows:
\[ i^{-1} : \Open(X) \lra \Open(X_j),\ \quad U \longmapsto j(U). \]
Thus we can relate the toposes of sheaves on~$X_j$ and~$X$ by the usual
pullback and pushforward functors.
\begin{align*}
  i^{-1} \F &= \text{sheafification of $(U \mapsto \colim_{U \preceq i^{-1}V} \Gamma(V,\F))$} \\
  i_* \G &= (U \mapsto \Gamma(i^{-1}U, \G)) = (U \mapsto \Gamma(j(U), \G))
\end{align*}
As familiar from honest topological subspace inclusions, the pushforward
functor~$i_* : \Sh(X_j) \to \Sh(X)$ is fully faithful and the composition~$i^{-1}
\circ i_* : \Sh(X_j) \to \Sh(X_j)$ is (canonically isomorphic to) the identity.

\subsection{Internal sheaves and sheafification}\label{sect:internal-sheaves}
It turns out that the image of
the pushforward functor~$i_* : \Sh(X_\Box) \to \Sh(X)$, where~$\Box$ is a modal
operator in~$\Sh(X)$, can be explicitly described. Namely, it consists exactly
of those sheaves which from the internal point of view
are so-called~\emph{$\Box$-sheaves}, a notion explained below.

Furthermore, if we identify~$\Sh(X_\Box)$ with its image in~$\Sh(X)$, the
pullback functor is given by an internal sheafification process with respect to
the modality~$\Box$. Thus the external situation of pushforward/pullback
translates to forget/sheafify. This broadens the scope of the internal
language of~$\Sh(X)$: It cannot only be used to talk about sheaves on~$X$ in a simple,
element-based language, but also to talk about sheaves on arbitrary subspaces
of~$X$.

To describe the notion of~$\Box$-sheaves and related ones, we switch to the internal
perspective and thus forget that we're working over a base space~$X$; we are simply given a modal operator~$\Box :
\Omega \to \Omega$ and have to take care that our proofs are intuitionistically valid. A
reference for the material in this subsection is a preprint by de
Vries~\cite{vries:sheafification}.\footnote{On page~5 of that
preprint there is a slight typing error: Fact~2.1(i) gives the
characterization of~$j$-closedness, not~$j$-denseness. The correct
characterization of~$j$-denseness in that context is~$\forall b \in B\_ j(b \in
A)$.}

Recall that a set~$S$ is a \emph{subsingleton} if and only if~$\forall x,y\?S\_
x = y$, and that a set~$S$ is a \emph{singleton} if and only if it is a subsingleton and
it is inhabited (\ie~$\exists x\?S\_ \top$); this amounts to~$\exists!x\?S\_ \top$.

\begin{defn}\label{defn:box-sheaves}
A set~$F$ is \emph{$\Box$-separated} if and only if
\[ \forall x,y\?F\_ \Box(x = y) \Longrightarrow x = y. \]
A set~$F$ is a \emph{$\Box$-sheaf} if and only if it is~$\Box$-separated and
\[ \forall S \subseteq F\_
  \Box(\speak{$S$ is a singleton}) \Longrightarrow
  \exists x\?F\_ \Box(x \in S). \]
\end{defn}

The two conditions can be combined: A set~$F$ is a~$\Box$-sheaf if and only if
\[ \forall S \subseteq F\_
  \Box(\speak{$S$ is a singleton}) \Longrightarrow
  \exists! x\?F\_ \Box(x \in S). \]

Intuitively, reading~``$\Box\varphi$'' as ``locally~$\varphi$'', a set
is~$\Box$-separated if elements which are locally equal are in fact equal. A
set is a~$\Box$-sheaf if furthermore for any set~$S$ of elements which locally
contains just a single element there is an element which is locally contained
in~$S$.

This phrasing is reminiscent of the usual gluing condition, which demands that
any family of sections which locally is just a single section (in that the
sections of the family agree on their common domain of definition) gives rise
to a global section which coincides with the given sections on their respective
domain. Remark~\ref{rem:modal-operator-in-presheaf-topos} below sketches how to
make this relation precise.

\begin{defn}\label{defn:plus-construction}
The \emph{plus construction} of a set~$F$ with respect to~$\Box$ is the set
\[ F^+ \defeq \{ S \subseteq F \,|\, \Box(\speak{$S$ is a singleton}) \}/{\sim},
\]
where the equivalence relation is defined by~$S \sim T \vcentcolon\Leftrightarrow
\Box(S = T)$. There is a canonical map~$F \to F^+$ given by~$x \mapsto
[\{x\}]$. The \emph{$\Box$-sheafi\-fi\-ca\-tion} of a set~$F$ is the
set~$F^{++}$.
\end{defn}

If~$F$ is~$\Box$-separated, then for any subset~$S \subseteq F$ it holds
that
\[ \Box(\speak{$S$ is a singleton}) \quad\Longleftrightarrow\quad
  \speak{$S$ is a subsingleton} \wedge \Box(\speak{$S$ is inhabited}). \]

\begin{rem}\label{rem:modal-operator-in-presheaf-topos}
The topos of \emph{pre}sheaves on a topological space~$X$ admits an
internal language as well~\cite[Section~VI.7, discussion after
Theorem~1]{moerdijk-maclane:sheaves-logic}. In it, there
exists a modal operator~$\Box$ reflecting the topology of~$X$. A presheaf on~$X$ is separated
in the usual sense if, from the internal perspective of~$\PSh(X)$, it
is~$\Box$-separated; and it is a sheaf if, from the internal perspective, it
is a~$\Box$-sheaf. Furthermore, the~$\Box$-sheafification of a presheaf
(considered as a set from the internal perspective) coincides with the usual
sheafification.\end{rem}

\begin{ex}\label{ex:special-sets-sheaves}
Any singleton set is a~$\Box$-sheaf. The empty set is
always~$\Box$-separated (trivially) and is a~$\Box$-sheaf if and only
if~$\Box\bot \Rightarrow \bot$.\end{ex}

We will see geometric examples of~$\Box$-sheaves in further sections.
For instance, on an integral or locally Noetherian scheme~$X$, the structure sheaf~$\O_X$
is~$\neg\neg$-separated and its~$\neg\neg$-sheafification is the sheaf~$\K_X$
of rational functions (Proposition~\ref{prop:kx-is-negneg-sheafification}).

\begin{lemma}For any set~$F$, it holds that: \begin{enumerate}
\item $F^+$ is~$\Box$-separated.
\item The canonical map~$F \to F^+$ is injective if and only if~$F$
is~$\Box$-separated.
\item If~$F$ is~$\Box$-separated, then $F^+$ is a~$\Box$-sheaf.
\item If~$F$ is a~$\Box$-sheaf, then the canonical map~$F \to F^+$ is bijective.
\end{enumerate}
Let $\Sh_\Box(\Set)$ be the full subcategory of~$\Set$ consisting of
the~$\Box$-sheaves. Then it holds that:
\begin{enumerate}
\addtocounter{enumi}{4}
\item The functor~$(\placeholder)^+ : \Set \to \Set$ is left exact.
\item The functor~$(\placeholder)^{++} : \Set \to \Sh_\Box(\Set)$ is left exact and left
adjoint to the forgetful functor~$\Sh_\Box(\Set) \to \Set,\ F \mapsto F$.
\end{enumerate}\end{lemma}
\begin{proof}These are all straightforward, and in fact simpler than their
classical counterparts, since there are no colimit formulas which would have to
be dealt with.
\end{proof}

\begin{rem}\label{rem:epi-in-box-sheaves}
As is to be expected from the familiar inclusion of sheaves in
presheaves on topological spaces, the forgetful functor~$\Sh_\Box(\Set) \to \Set$
does not in general preserve colimits. It is instructive to see why
epimorphisms in~$\Sh_\Box(\Set)$ need not be epimorphisms in~$\Set$: A map~$f:A
\to B$ between~$\Box$-sheaves is an epimorphism in~$\Sh_\Box(\Set)$ if and only
if
\[ \forall y\?B\_ \Box(\exists x\?X\_ f(x) = y), \]
that is preimages do not need to exist, it suffices for them to~``$\Box$-exist''.
(Using results about the~$\Box$-translation, to be introduced below, this
characterization will be obvious.) This condition is intuitionistically weaker
than the condition that~$f$ is an epimorphism in~$\Set$, \ie that~$f$ is
surjective. This should be compared to the failure of the forgetful functor~$\Sh(X)
\to \PSh(X)$ to preserve epimorphisms: A morphism of sheaves does not need to
have preimages for any local section in order to be an epimorphism. Instead, it
suffices for any local section to \emph{locally} have preimages.\end{rem}

\begin{prop}Let~$X$ be a topological space. Let~$\Box$ be a modal operator
in~$\Sh(X)$. Let~$i : X_\Box \hookrightarrow X$ be the inclusion of the
associated sublocale. Corestricting the pushforward functor~$i_* : \Sh(X_\Box) \to
\Sh(X)$ to its essential image, it induces an equivalence~$\Sh(X_\Box) \simeq
\Sh_\Box(\Sh(X))$ between the category of sheaves on~$X_\Box$ and the category
of~$\Box$-sheaves in~$\Sh(X)$.
\end{prop}
\begin{proof}For the further development of the theory, we need the statement
of this proposition, but not the proof, which really is routine in dealing with
subtoposes and modal operators. Nevertheless, a proof can proceed by
combining Example~A4.6.2(a) and Theorem~C1.4.7
of~\cite{johnstone:elephant}, observing that for a topos of sheaves on a
locale~$Y$ it holds that~$\Open(Y) = \Gamma(Y, \Omega_{\Sh(Y)})$, and that the
subobject classifier of~$\Sh_\Box(\Sh(X))$ is~$\{ \varphi : \Omega_{\Sh(X)} \,|\,
\Box \varphi \Leftrightarrow \varphi \}$.
\end{proof}

\begin{rem}It's possible to rewrite the sheaf condition in the following form.
A set~$F$ is~$\Box$-separated if and only if, for any truth value~$\varphi \?
\Omega$ such that~$\Box\varphi$, the canonical map
\[ F \lra F^\varphi, \]
which maps an element~$x\?F$ to the constant map~$\varphi \to X$ with value~$x$
(where~$\varphi$ is considered as a subset of the terminal set~$1$), is
injective. The set~$F$ is a~$\Box$-sheaf if and only if furthermore this map is
surjective for all such truth values.
\end{rem}

\subsection{Sheaves for the double negation modality}
\label{sect:negneg-sheaves}

Recall that if~$\Box$ is the modal operator associated to a sub\emph{space}~$Y$
of a topological space~$X$, then the sheaves on~$X$ which are~$\Box$-sheaves
are easy to describe: These are precisely the sheaves in the essential image of
the pushforward functor~$\Sh(Y) \to \Sh(X)$. For the double negation modality,
the same is true, only that~$Y$ is then the perhaps unfamiliar \emph{smallest
dense sublocale} of~$X$.

The following proposition gives a characterization of~$\neg\neg$-separated
presheaves and~$\neg\neg$-sheaves in explicit terms.

\begin{prop}\label{prop:negneg-sheaves}
Let~$X$ be a topological space. Let~$\F$ be a sheaf on~$X$. Then:
\begin{enumerate}
\item $\F$ is~$\neg\neg$-separated if and only if any two local sections
of~$\F$, which are defined on a common domain and which agree on a dense open
subset of their domain, are already equal.
\item $\F$ is a~$\neg\neg$-sheaf if and only if it is~$\neg\neg$-separated and
for any open~$U \subseteq X$ and any open~$V \subseteq U$ dense
in~$U$, any~$V$-section of~$\F$ extends to an~$U$-section of~$\F$.
\item If~$\F$ is~$\neg\neg$-separated, the sections of $\F^+$ on an open
subset~$U \subseteq X$ can be described by pairs~$\langle V,s \rangle$, where~$V$ is a dense
open subset of~$U$ and~$s$ is a section of~$\F$ on~$V$. Two such pairs~$\langle V,s \rangle,
\langle V',s' \rangle$ determine the same element in~$\Gamma(U,\F^+)$ if and only if~$s$ and~$s'$
agree on~$V \cap V'$.
\end{enumerate}
\end{prop}
\begin{proof}
The first statement is obvious from the definition of~$\neg\neg$-separatedness
(Definition~\ref{defn:box-sheaves} for~$\Box = \neg\neg$) and the geometric
interpretation of double negation (Proposition~\ref{prop:modops-kripke}).

For the second statement, we need to show that, assuming that~$\F$
is~$\neg\neg$-separated, the sheaf~$\F$ has the extension property if and only if
\begin{multline*}
  \Sh(X) \models \forall \S \? \P(\F)\_
  \speak{$\S$ is a subsingleton} \wedge
  \neg\neg(\speak{$\S$ is inhabited}) \Longrightarrow \\
  \exists x\?\F\_ \neg\neg(x \in \S).
\end{multline*}
A section~$\S \in \Gamma(U,\P(\F))$ which internally is a
subsingleton and \notnot inhabited is precisely a subsheaf~$\S \hookrightarrow
\F|_U$ such that all stalks~$\S_x$, $x \in U$ are subsingletons and such that for
some dense open subset~$V \subseteq U$, the stalks~$\S_x$, $x \in V$ are
inhabited. This is precisely the datum of a section of~$\F$ defined on some
dense open subset of~$U$, considering the gluing of the unique germs in~$\S_x$ for
those points~$x$ such that~$\S_x$ is inhabited. (Conversely, a section~$s \in
\Gamma(V,\F)$ defines a subsheaf~$\S$ by setting~$\Gamma(W,\S) \defeq \{ s|_W \,|\,
W \subseteq V \}$.)

In view of this explicit description and the observation that the asserted
existence~(``$\exists x\?\F\_ \neg\neg(x \in \S)$'') is actually a question of
unique existence, the second statement follows.

For the third statement, one can check that the presheaf on~$X$ defined by
\[ \text{$U \subseteq X$ open} \quad\longmapsto\quad
  \{ \langle V,s \rangle \,|\, \text{$V \subseteq U$ dense open},\ s \in \Gamma(V,\F)
  \}/{\sim} \]
is in fact a sheaf (with respect to the topology of~$X$), internally a $\neg\neg$-sheaf,
and that it has the universal property of the~$\neg\neg$-sheafification
of~$\F$.
\end{proof}

The conditions~(1) and~(2) of Proposition~\ref{prop:negneg-sheaves} can be
summarized as follows: A sheaf~$\F$ on a topological space is
a~$\neg\neg$-sheaf if and only if, for any open subset~$U \subseteq X$, the
restriction map~$\Gamma(\Int\Clos U, \F) \to \Gamma(U,\F)$ is
bijective~\cite[Lemma~36]{jackson:sheaf-theoretic-measure-theory}.

In the case that~$X$ contains a \emph{generic point}, that is a point~$\xi \in X$
such that~$\Clos\{\xi\} = X$, we can describe the sublocale~$X_{\neg\neg}$ in
simple terms: In this case, it coincides with the subspace~$\{\xi\}$.
For instance, such a generic point exists and is unique if~$X$ is an irreducible scheme.

\begin{lemma}\label{lemma:negneg-generic-point}
Let~$X$ be a topological space and~$\xi \in X$ be a point such
that~$\Clos\{\xi\} = X$. Then the modal operator~$\Box \defequiv ((\placeholder
\Rightarrow \notat{\xi}) \Rightarrow \notat{\xi})$ coincides with the double
negation modality and~$X_{\neg\neg} = \{\xi\}$ as sublocales of~$X$.\end{lemma}
\begin{proof}The semantics of the formula~$\notat{\xi}$ was defined by the
equivalence
\[ U \models \notat{\xi} \quad\Longleftrightarrow\quad
  \xi \not\in U. \]
By the assumption on~$\xi$, this is equivalent to requiring~$U = \emptyset$.
Thus for any open subset~$U$ the formulas~$\notat{\xi}$ and~$\bot$ have the
same meaning; they are therefore logically equivalent from the internal point of
view. The given modal operator thus simplifies:
\[ \Box\varphi \quad\equiv\quad ((\varphi \Rightarrow \notat{\xi}) \Rightarrow \notat{\xi})
  \quad\Leftrightarrow\quad ((\varphi \Rightarrow \bot) \Rightarrow \bot)
  \quad\Leftrightarrow\quad \neg\neg\varphi. \]
The second claim follows from Table~\ref{table:nuclei}.
\end{proof}

\begin{cor}\label{cor:negneg-generic-point-pushpull}
Let~$X$ be a topological space and let~$\xi \in X$ be a point such
that~$\Clos\{\xi\} = X$. Since~$X_{\neg\neg} = \{\xi\}$, the
category of~$\neg\neg$-sheaves in~$\Sh(X)$ coincides with the category of
sheaves on~$\{\xi\}$ and can therefore be identified with the category of sets.
Under this identification,
\begin{enumerate}
\item sheafifying an object~$\F \in \Sh(X)$ with respect
to the double negation modality (\ie pulling back to~$X_{\neg\neg}$) is the
same as calculating its generic stalk~$\F_\xi$ and
\item pushing forward a set~$M$ along~$X_{\neg\neg} \hookrightarrow X$ is the
same as calculating the constant sheaf associated to~$M$.
\end{enumerate}
\end{cor}
\begin{proof}The first statement follows because pulling back to~$X_{\neg\neg}$
is the same as pulling back to~$\{\xi\}$. The pushforward of a set~$M$,
considered as a sheaf on~$X_{\neg\neg}$, to~$X$ is explicitly given by
\[ U \quad\longmapsto\quad \begin{cases}
  M, & \text{if $U \neq \emptyset$,} \\
  \{\star\}, & \text{else.}
\end{cases} \]
We omit the routine verification that this sheaf coincides with the constant
sheaf~$\underline{M}$ associated to~$M$.
\end{proof}

The following technical lemma will occasionally be handy. It is an internal
reflection of the fact that an open subset of an affine scheme can always be
written as the union of standard open subsets. We will generalize it
to schemes which are not necessarily integral in
Section~\ref{sect:rational-functions} (see
Lemma~\ref{lemma:dense-standard-reflection-generalized}).

\begin{lemma}\label{lemma:dense-standard-reflection}
Let~$X$ be an integral scheme. Let~$\varphi$ be any formula
over~$X$. Then
\[ \Sh(X) \models \neg\neg\varphi \Longrightarrow \exists f\?\O_X\_
  \neg\neg(\speak{$f$ \inv}) \wedge (\speak{$f$ \inv} \Rightarrow \varphi). \]
\end{lemma}
\begin{proof}We may assume that~$X$ is the spectrum of an integral domain~$A$
and that there is a dense open subset~$U \subseteq X$ on which~$\varphi$ holds.
The open set~$U$ may be covered by standard open subsets~$D(f_i)$; since~$X$ is
irreducible, at least one of these is itself
dense. We may take this~$f_i$ as the sought~$f$.
\end{proof}

We can now also follow up on a promise made in Section~\ref{sect:appreciating-intuitionistic-logic} and prove the following
somewhat tangential lemma.
\begin{lemma}\label{lemma:boolean-dense}
Let~$X$ be a topological space. The internal language of~$\Sh(X)$ is Boolean if
and only if for any open subset~$U \subseteq X$ it holds that~$U$ is the only
dense open subset of~$U$.
\end{lemma}
\begin{proof}That the internal language of~$\Sh(X)$ is Boolean amounts to
\[ \Sh(X) \models \forall \varphi\?\Omega\_ \neg\neg\varphi \Rightarrow
\varphi. \]
This is equivalent to the external statement that for any open subset~$U
\subseteq X$ and for any open subset~$V \subseteq U$ it holds that: If~$V$ is
dense in~$U$, then~$V$ is equal to~$U$.
\end{proof}

\subsection{\texorpdfstring{The~$\Box$-translation}{The □-translation}}
\label{sect:box-translation}

In logic, there is certain well-known transformation~$\varphi
\mapsto \varphi^{\neg\neg}$ on formulas, the \emph{double negation
translation}, with the following curious property: A formula~$\varphi$ is
derivable in classical logic if and only if its
translation~$\varphi^{\neg\neg}$ is derivable in intuitionistic logic. The
translation~$\varphi^{\neg\neg}$ is obtained from~$\varphi$ by putting
``$\neg\neg$'' before any subformula, \ie before any~``$\exists$''
and~``$\forall$'', around any logical connective, and around any atomic
statement (``$x=y$'', ``$x \in A$''). For instance, the double negation
translation of ``$f$ is surjective'' is
\[ \neg\neg \forall y\?Y\_ \neg\neg \exists x\?X\_ \neg\neg f(x) = y. \]

We will describe a slight generalization of the double negation translation,
the~$\Box$-translation for any modal operator~$\Box$. It will be pivotal
for using the internal language of a space~$X$ to express internal statements
about sheaves defined on subspaces of~$X$. The~$\Box$-translation has been studied
in other contexts
before~\cite{aczel:russell-prawitz,escardo:oliva:peirce-shift}. To the best of
our knowledge, this application -- expressing the internal language of
subtoposes in the internal language of the ambient topos -- is new.

\begin{defn}The~\emph{$\Box$-translation} is recursively defined as follows.
\newcommand{\optBox}{\textcolor{gray}{\Box}}
\begin{align*}
  (f = g)^\Box &\defequiv \Box(f = g) \\
  (x \in A)^\Box &\defequiv \Box(x \in A) \\
  \top^\Box &\defequiv \Box\top \quad \text{($\Leftrightarrow \top$)} \\
  \bot^\Box &\defequiv \Box\bot \\
  (\varphi \wedge \psi)^\Box &\defequiv \optBox(\varphi^\Box \wedge \psi^\Box) &
  \textstyle (\bigwedge_i \varphi_i)^\Box &\defequiv \textstyle \optBox(\bigwedge_i \varphi_i^\Box) \\
  (\varphi \vee \psi)^\Box &\defequiv \Box(\varphi^\Box \vee \psi^\Box) &
  \textstyle (\bigvee_i \varphi_i)^\Box &\defequiv \textstyle \Box(\bigvee_i \varphi_i^\Box) \\
  (\varphi \Rightarrow \psi)^\Box &\defequiv \optBox(\varphi^\Box \Rightarrow \psi^\Box) \\
  (\forall x\?X\_ \varphi)^\Box &\defequiv \optBox(\forall x\?X\_ \varphi^\Box) &
  (\forall X\_ \varphi)^\Box &\defequiv \optBox(\forall X\_ \varphi^\Box) \\
  (\exists x\?X\_ \varphi)^\Box &\defequiv \Box(\exists x\?X\_ \varphi^\Box) &
  (\exists X\_ \varphi)^\Box &\defequiv \Box(\exists X\_ \varphi^\Box)
\end{align*}
\end{defn}

\begin{defn}A formula~$\varphi$ is \emph{$\Box$-stable} if and only
if~$\Box\varphi$ implies~$\varphi$.\end{defn}

\begin{lemma}\begin{enumerate}
\item Formulas in the image of the $\Box$-translation are~$\Box$-stable,
\ie for any formula~$\varphi$ it holds that
$\Box(\varphi^\Box) \Longrightarrow \varphi^\Box$.
\item In the definition of the~$\Box$-translation, one may omit the boxes
printed in gray.
\end{enumerate}\end{lemma}
\begin{proof}The first statement is obvious, since one of the axioms for a modal
operator demands that~$\Box\Box\varphi \Rightarrow \Box\varphi$ for any
formula~$\varphi$. The second statement follows by an induction on the
formula structure. By way of example, we prove the case for~``$\Rightarrow$'':
\newcommand{\withgray}{\text{$\Box$ with the gray parts}}
\newcommand{\withoutgray}{\text{$\Box$ without the gray parts}}
\begin{align*}
  &\ (\varphi \Rightarrow \psi)^\withgray \\
  \Longleftrightarrow &\ \Box(\varphi^\withgray \Rightarrow \psi^\withgray) \\
  \Longleftrightarrow &\ (\varphi^\withgray \Rightarrow \psi^\withgray) \\
  \Longleftrightarrow &\ (\varphi^\withoutgray \Rightarrow \psi^\withoutgray) \\
  \Longleftrightarrow &\ (\varphi \Rightarrow \psi)^\withoutgray
\end{align*}
The first step is by definition; the second by~$\Box$-stability
of~$\psi^\withgray$ and the intuitionistic tautology~$\Box(\alpha \Rightarrow
\beta) \Leftrightarrow (\alpha \Rightarrow \beta)$ for~$\Box$-stable
formulas~$\beta$; the third by the induction hypothesis; and the fourth by
definition.
\end{proof}

\begin{lemma}\label{lemma:box-translation-sound}
The~$\Box$-translation is sound with respect to intuitionistic logic:
Assume that there exists an intuitionistic proof of an
implication~$\varphi \Rightarrow \psi$. Then there is also an intuitionistic
proof of the translated implication~$\varphi^\Box \Rightarrow \psi^\Box$.
\end{lemma}
\begin{proof}This follows by an induction on the structure of intuitionistic
proofs. We have to verify that we can mirror any inference rule of
intuitionistic logic in the translation. For instance, one of the disjunction
rules justifies the following proof scheme: In order to prove~$\varphi \vee
\psi \Rightarrow \chi$, it suffices to give proofs of~$\varphi \Rightarrow
\chi$ and~$\psi \Rightarrow \chi$. We have to justify the translated proof
scheme: In order to prove~$(\varphi \vee \psi)^\Box \Rightarrow \chi^\Box$, it
suffices to give proofs of~$\varphi^\Box \Rightarrow \chi^\Box$ and~$\psi^\Box
\Rightarrow \chi^\Box$.

So assume that proofs of the two implications are given. Further
assume~$(\varphi \vee \psi)^\Box$, \ie~$\Box(\varphi^\Box \vee \psi^\Box)$.
We want to show~$\chi^\Box$. Since this is a~$\Box$-stable statement, we may
assume that in fact~$\varphi^\Box \vee \psi^\Box$ holds. Then the claim is
obvious by the two given proofs.

The cases for the other rules (see Appendix~\ref{appendix:inference-rules} for
a list) are similar and left to the reader.\end{proof}

\begin{rem}The reader well-versed in formal logic will have noticed that we are
mixing syntax and semantics here. The proper way to state Lemma~\ref{lemma:box-translation-sound} would be
to formally adjoin a box operator to the language of intuitionistic logic,
governed by three inference rules which are modeled on the three axioms for a
modal operator. This formal box operator could then be instantiated by any
concrete modal operator~$\Box : \Omega \to \Omega$.\end{rem}

Soundness of the~$\Box$-translation is important for the following reason.
If~$\varphi$ and~$\varphi'$ are equivalent formulas, we are
accustomed to be able to freely substitute~$\varphi$ by~$\varphi'$ anywhere we
want. Since a modal operator~$\Box$ is semantically defined as a map~$\Omega
\to \Omega$, it is trivially justified that~$\Box\varphi$ and~$\Box\varphi'$
are equivalent: The formulas~$\varphi$ and~$\varphi'$ give rise to the
\emph{same} element~$\{x \? 1 \,|\, \varphi\} = \{x \? 1 \,|\, \varphi'\}$
of~$\Omega$, and therefore their images under~$\Box$ are equal as well.

However, it is \emph{not} clear and in fact wrong in general that the translated formulas~$\varphi^\Box$
and~$(\varphi')^\Box$ are equivalent. This follows only if the soundness
lemma can be applied (two times, once for each direction). We should stress that to apply this
lemma, it is not enough to merely \emph{know} that~$\varphi$ and~$\varphi'$ are
equivalent; instead, there has to be an intuitionistic proof of this
equivalence. This is really a stronger requirement, since an
equivalence~$\varphi \Leftrightarrow \varphi'$ might
hold in a particular model, \ie in the internal language of some particular
topos, without possessing an intuitionistic proof, \ie holding in any topos. We
give an explicit example of this situation below
(Example~\ref{ex:translation-equivalence}).

\begin{lemma}\label{lemma:open-stalk}
Let~$\varphi$ be a formula such that for any subformulas~$\psi$
appearing as antecedents of implications, it holds that~$\psi^\Box \Rightarrow
\Box\psi$. (In particular, this condition is satisfied if there are
no~``$\Rightarrow$'' signs in~$\varphi$ or if~$\varphi$ is a geometric implication.) Then $\Box\varphi \Rightarrow
\varphi^\Box$.\end{lemma}
\begin{proof}We prove this by an induction on the formula structure. All cases
except for~``$\Rightarrow$'' are obvious. For this case, assume~$\Box(\psi
\Rightarrow \chi)$; we are to show that~$(\psi^\Box \Rightarrow \chi^\Box)$.
Since this is a~$\Box$-stable statement, we can in fact assume that~$(\psi
\Rightarrow \chi)$. We then have
\[ \psi^\Box \Longrightarrow \Box\psi \Longrightarrow \Box\chi
\Longrightarrow \chi^\Box, \]
with the first step being by the requirement on antecedents, the second by the
monotonicity of~$\Box$, and the third by the induction hypothesis.
\end{proof}

\begin{lemma}\label{lemma:stalk-open}
Let~$\varphi$ be a geometric formula.
Then $\varphi^\Box \Leftrightarrow \Box\varphi$.\end{lemma}
\begin{proof}The~``$\Leftarrow$'' direction is by Lemma~\ref{lemma:open-stalk}.
The~``$\Rightarrow$'' direction is an induction on the formula structure. By way of example, we verify
the case of~``$\bigvee$''. So assume~$\Box(\bigvee_i \varphi_i^\Box)$; we are
to show that~$\Box(\bigvee_i \varphi_i)$. Since this is a boxed statement, we
may in fact assume~$\bigvee_i \varphi_i^\Box$, so for some index~$j$, it holds
that~$\varphi_j^\Box$. By the induction hypothesis, it follows
that~$\Box\varphi_j$. By~$\varphi_j \Rightarrow \bigvee_i \varphi_i$ and the
monotonicity of~$\Box$, it follows that~$\Box(\bigvee_i \varphi_i)$.
\end{proof}

An analogous argument for infinite conjunctions is not valid:
Assume~$(\bigwedge_i \varphi_i)^\Box$. So for all~$j$,~$\varphi_j^\Box$ holds.
By the induction hypothesis,~$\Box\varphi_j$ holds for any~$j$. But from this
we may not deduce~$\Box\bigwedge_i \varphi_i$, since the axioms for a modal
operator only require commutativity with finite conjunctions. This failure also
has a geometric interpretation, for instance in the special case~$\Box =
\neg\neg$: Given dense open subsets~$U_i$ on which formulas~$\varphi_i$ hold,
we may not conclude that there exists a single dense open subset~$U$ on which
all the formulas~$\varphi_i$ hold.

\begin{rem}In the special case that~$\Box$ is the double negation modality,
Lemma~\ref{lemma:stalk-open} holds with slightly weaker hypotheses: Namely, implications may occur
in~$\varphi$, provided that for their antecedents~$\psi$ it holds that~$\psi
\Rightarrow \psi^\Box$. This is because for the double negation modality,
the formula~$\Box(\psi \Rightarrow \chi)$ is equivalent to~$\psi \Rightarrow
\Box\chi$. (In general, for an arbitrary modality, only the former implies the latter, but not vice versa.) The case
for~``$\Rightarrow$'' in the inductive proof then goes as follows:
Assume~$(\psi \Rightarrow \chi)^\Box$. Then~$\psi \Rightarrow \psi^\Box
\Rightarrow \chi^\Box \Rightarrow \Box\chi$, so~$\Box(\psi \Rightarrow \chi)$.
\end{rem}

\begin{lemma}\label{lemma:stalk-open-with-hypothesis}
Let~$\varphi, \varphi', \psi$ be formulas. Assume that:
\begin{enumerate}
\item The formula $\varphi'$ is geometric. (More generally, it suffices for~$(\varphi')^\Box$
to imply~$\Box\varphi'$.)
\item There is an intuitionistic proof that~$\varphi$
and~$\varphi'$ are equivalent under the (only) hypothesis~$\psi$.
\item Both~$\Box\psi$ and~$\psi^\Box$ hold.
\end{enumerate}
Then $\varphi^\Box \Rightarrow \Box\varphi$.
\end{lemma}
\begin{proof}
Assume~$\varphi^\Box$. Since~$\psi^\Box$, $(\varphi \wedge \psi)^\Box$. Because
the~$\Box$-translation is sound with respect to intuitionistic logic
(Lemma~\ref{lemma:box-translation-sound})
it follows that~$(\varphi')^\Box$. As~$\varphi'$ is geometric, it follows
that~$\Box\varphi'$. Since~$\Box\psi$ holds, it follows that~$\Box\varphi$.
\end{proof}

\begin{ex}\label{ex:module-zero-geometric}
Let~$M$ be an~$R$-module. The statement that~$M$ is zero is not
geometric: $\varphi \defequiv (\forall x\?M\_ x = 0)$. But if~$M$ is generated by some finite
family~$x_1,\ldots,x_n\?M$, then~$\varphi$ is equivalent to the
statement~$\varphi' \defequiv (x_1 = 0
\wedge \cdots \wedge x_n = 0)$ which is geometric; and there is an
intuitionistic proof of this equivalence. Since no implication signs occur
in~$\psi \defequiv \speak{$M$ is generated by~$x_1,\ldots,x_n$}$, Lemma~\ref{lemma:stalk-open-with-hypothesis} is
applicable and shows that~$\varphi^\Box$ implies~$\Box\varphi$.
This example will gain geometric meaning in
Lemma~\ref{lemma:module-zero-point-neighbourhood}.
\end{ex}

\begin{lemma}For the modality~$\Box$ defined by~$\Box\varphi \defequiv ((\varphi
      \Rightarrow \alpha) \Rightarrow \alpha)$, where~$\alpha$ is a fixed
proposition, the~$\Box$-translation of the law of excluded middle holds.
In particular, this applies to the double negation modality~$\Box = \neg\neg$, where~$\alpha =
\bot$.\end{lemma}
\begin{proof}We are to show that~$(\varphi \vee \neg\varphi)^\Box$, \ie that
\[ ((\varphi^\Box \vee (\varphi^\Box \Rightarrow \alpha)) \Longrightarrow
    \alpha) \Longrightarrow \alpha. \]
So assume that the antecedent holds. If~$\varphi^\Box$ holds, then in
particular~$\varphi^\Box \vee (\varphi^\Box \Rightarrow \alpha)$ and thus~$\alpha$
hold. Therefore it follows that~$(\varphi^\Box \Rightarrow \alpha)$. This
implies~$\varphi^\Box \vee (\varphi^\Box \Rightarrow \alpha)$ and
thus~$\alpha$.
\end{proof}

\subsection{\texorpdfstring{Truth at stalks \vs truth on neighborhoods}{Truth
at stalks vs. truth on neighborhoods}}\label{sect:spreading}
We now state the crucial property of the~$\Box$-translation. Recall
that~``$X_\Box$'' denotes the sublocale of~$X$ induced by~$\Box$
(Definition~\ref{defn:subspace-by-nucleus}).
\begin{thm}\label{thm:box-translation-semantically}
Let~$X$ be a topological space. Let~$\Box$ be a modal operator
in~$\Sh(X)$. Let~$\varphi$ be a formula over~$X$. Then
\[ \Sh(X) \models \varphi^\Box \quad\text{iff}\quad
  \Sh(X_\Box) \models \varphi, \]
where on the right hand side, all parameters occurring in~$\varphi$ were pulled
back to~$X_\Box$ along the inclusion~$X_\Box \hookrightarrow X$.
\end{thm}

We have not yet explicitly stated the Kripke--Joyal semantics for a sheaf topos
over a locale, which~$X_\Box$ is in general. The definition is exactly the same
as in the case for sheaf toposes over a topological space, only that any
mention of ``open sets'' has to be substituted by the more general ``opens''
and any mention of the union operator~``$\bigcup$'' has to be interpreted by
the supremum operator in the frame of opens of the locale. For~$X_\Box$, this
is~$\sup U_i = j(\bigcup_i U_i)$. Before giving a proof of Theorem~\ref{thm:box-translation-semantically}, we want
to discuss some of its consequences.

\begin{cor}\label{cor:spreading}
Let~$X$ be a topological space.
\begin{enumerate}
\item Let~$U \subseteq X$ be an open subset and let~$\Box\varphi \defequiv (U
\Rightarrow \varphi)$. Then
\[ \Sh(X) \models \varphi^\Box \quad\text{iff}\quad \Sh(U) \models \varphi. \]
\item Let~$A \subseteq X$ be a closed subset and let~$\Box\varphi \defequiv
(\varphi \vee A^c)$. Then
\[ \Sh(X) \models \varphi^\Box \quad\text{iff}\quad \Sh(A) \models \varphi. \]
\item Let~$\Box\varphi \defequiv \neg\neg\varphi$. Then
\[ \Sh(X) \models \varphi^\Box \quad\text{iff}\quad \Sh(X_{\neg\neg}) \models \varphi. \]
\item Let~$x \in X$ be a point and let~$\Box\varphi \defequiv ((\varphi
\Rightarrow \notat{x}) \Rightarrow \notat{x})$. Then
\[ \Sh(X) \models \varphi^\Box \quad\text{iff}\quad \text{$\varphi$ holds
at~$x$}. \]
\end{enumerate}
\end{cor}
\begin{proof}Combine Theorem~\ref{thm:box-translation-semantically} and
Table~\ref{table:nuclei}.\end{proof}

We want to discuss the last case of Corollary~\ref{cor:spreading} in more detail. Let~$x$ be a
point of a topological space~$X$ and let~$\varphi$ be a formula. Let~$\Box$ be
the modal operator given in the corollary. Then~$\varphi$ \emph{holds at~$x$}
if and only if, from the internal perspective of~$\Sh(X)$, the translated
formula~$\varphi^\Box$ holds; and~$\varphi$ \emph{holds on some open
neighborhood of~$x$} if and only if, from the internal perspective, the
formula~$\Box\varphi$ holds.

Thus the question whether the truth of~$\varphi$ at the point~$x$ spreads to
some open neighborhood can be formulated in the following way:
\begin{quote}
\emph{Does~$\varphi^\Box$ imply~$\Box\varphi$ in the internal language
of~$\Sh(X)$?}
\end{quote}
Phrased this way, technicalities like appropriately shrinking open
neighborhoods are blinded out. A purposefully trivial example to illustrate
this is the following. Let~$X$ be a scheme (or a ringed space). Let~$f,g \in
\Gamma(X,\O_X)$ be global functions. Suppose that the germs of~$f$ and~$g$ are
zero in some stalk~$\O_{X,x}$; we want to show that they are zero on a common
open neighborhood of~$x$.

\begin{proof}[Usual proof]Since the germ of~$f$ vanishes in~$\O_{X,x}$, there
is an open neighborhood~$U_1$ of~$x$ such that~$f|_{U_1} = 0$
in~$\Gamma(U_1,\O_X)$. Since furthermore the germ of~$g$ vanishes in the same stalk,
there exists an open neighborhood~$U_2$ of~$x$ such that~$g|_{U_2} = 0$. The
intersection of both neighborhoods is still an open neighborhood of~$x$; on
this neighborhood both~$f$ and~$g$ vanish.
\end{proof}

\begin{proof}[Proof in the internal language]We may suppose that~$(f = 0 \wedge
g = 0)^\Box$, that is $\Box(f=0) \wedge \Box(g=0)$, and have to prove
that~$\Box(f=0 \wedge g=0)$. (To this end, we could simply invoke the third
axiom on a modal operator, but we want to stay close to the given external
proof.) So by assumption, both~$\Box(f=0)$ and~$\Box(g=0)$ hold. Since our goal
is to prove a boxed statement, we may in fact assume that~$f = 0$ and~$g = 0$.
Thus~$f = 0 \wedge g = 0$.\end{proof}

By using the internal language with its modal operators, we can thus reduce
basic facts of scheme theory which deal with stalks and neighborhoods to facts
of algebra in a \emph{modal intuitionistic context}. As with using the internal
language in its basic form without modalities, this brings conceptual clarity
and reduced technical overhead. There are, however, two more distinctive
advantages. Firstly, many internal proofs do not require specific properties of
the modal operator and thus work with any modal operator. By interpreting such
a proof using different operators, one obtains an entire family of external
statements without any additional work (see
Lemma~\ref{lemma:module-zero-point-neighbourhood} for an example).

Secondly, the following corollary gives a general metatheorem which is
applicable to a wide range of cases. It allows to decide whether spreading will
occur (or is likely not to occur) simply by looking at the \emph{logical form}
of the statement in question.

\begin{cor}\label{cor:geometric-spreading}
Let~$X$ be a topological space. Let~$\varphi$ be a formula.
If~$\varphi$ is geometric, truth of~$\varphi$ at a point~$x \in X$ implies
truth of~$\varphi$ on some open neighborhood of~$x$, and vice versa.\end{cor}
\begin{proof}By the purely logical lemmas of Section~\ref{sect:box-translation}, it holds
that~$\varphi^\Box \Leftrightarrow \Box\varphi$.
\end{proof}

\begin{cor}
Let~$X$ be a topological space. Let~$\varphi$ be a formula.
If~$\varphi$ is geometric, the property ``$\varphi$ holds at a point~$x \in
X$'' is open.
\end{cor}
\begin{proof}This is just a reformulation of the previous corollary:
If~$\varphi$ holds at a point~$x \in X$, it holds on some open
neighborhood~$U$ of~$x$ as well. Going back to stalks, it follows
that~$\varphi$ holds at every point of~$U$.\end{proof}

\begin{ex}Let~$X$ be a scheme (or a ringed space). Since the condition for a
function~$f\?\O_X$ to be nilpotent is geometric (it is~$\bigvee_{n \geq 0} f^n
= 0$), nilpotency of~$f$ at a point is equivalent to nilpotency on some open
neighborhood.\end{ex}

Combined with Lemma~\ref{lemma:stalk-open-with-hypothesis}, this metatheorem is
quite useful. We will illustrate it with several examples in the next subsection.

An important special case of spreading from stalks to neighborhoods is the
case of spreading from the generic point (should it exist) to a dense open
subset. Whether this occurs can be phrased by
Lemma~\ref{lemma:negneg-generic-point} as follows:
\begin{quote}
\emph{Does~$\varphi^{\neg\neg}$ imply~$\neg\neg\varphi$ in the internal language
of~$\Sh(X)$?}
\end{quote}
This question is a question of ordinary (non-modal) intuitionistic algebra.

\begin{ex}We have seen in Remark~\ref{rem:epi-in-box-sheaves} that a
morphism~$f : A \to B$ in~$\Sh(X_\Box) \simeq \Sh_\Box(\Sh(X))$ is an
epimorphism if and only if \[ \Sh(X) \models \forall y\?B\_ \Box(\exists x\?X\_
f(x) = y). \] We can now understand a simple proof of this fact:
\begin{align*}
  &\ \text{$f$ is an epimorphism in~$\Sh_\Box(\Sh(X))$} \\
  \Longleftrightarrow&\
    \Sh_\Box(\Sh(X)) \models \speak{$f$ is surjective} \\
  \Longleftrightarrow&\
    \Sh(X) \models \left(\speak{$f$ is surjective}\right)^\Box \\
  \Longleftrightarrow&\
    \Sh(X) \models \forall y\?B\_ \Box(\exists x\?X\_ \Box(f(x) = y)) \\
  \Longleftrightarrow&\
    \Sh(X) \models \forall y\?B\_ \Box(\exists x\?X\_ f(x) = y).
\end{align*}
The ultimate equivalence is by Lemma~\ref{lemma:stalk-open}, applied to the
geometric subformula~``$\exists x\?X\_ f(x) = y$''.
\end{ex}

\begin{rem}Theorem~\ref{thm:box-translation-semantically} can also be motivated
by purely logical considerations. Namely, one can check that interpreting a
formula~$\varphi$ by $\Sh(X) \models \varphi^\Box$ gives rise to a model of
intuitionistic logic -- if~$\varphi$ intuitionistically implies~$\psi$,
then~$\Sh(X) \models \varphi^\Box$ implies~$\Sh(X) \models \psi^\Box$.
It is therefore a natural question whether there exists a
topos~$\E$ such that~$\E \models \varphi$ if and only if~$\Sh(X) \models
\varphi^\Box$. Theorem~\ref{thm:box-translation-semantically} gives an
affirmative answer to this question, explicitly stating that~$\E \defeq
\Sh(X_\Box)$ is such a topos.\end{rem}

\begin{proof}[Proof of Theorem~\ref{thm:box-translation-semantically}]
A fancy proof goes as follows. First, one shows intuitionistically that for a
modal operator~$\Box$ in~$\Set$, it holds that
\[ \Set \models \varphi^\Box \quad\Longleftrightarrow\quad
  \Sh_\Box(\Set) \models \varphi. \]
This can be verified by induction on the structure of
formulas~$\varphi$. Then one interprets this result in the sheaf topos~$\Sh(X)$:
\begin{align*}
  &\ \Sh(X) \models \varphi^\Box \\
  \Longleftrightarrow&\
  \Sh(X) \models \speak{$\Set \models \varphi^\Box$} &&\text{by idempotency}\\
  \Longleftrightarrow&\
  \Sh(X) \models \speak{$\Sh_\Box(\Set) \models \varphi$} &&\text{by the first step} \\
  \Longleftrightarrow&\
  \Sh_\Box(\Sh(X)) \models \varphi &&\text{by idempotency} \\
  \Longleftrightarrow&\
  \Sh(X_\Box) \models \varphi &&\text{since~$\Sh_\Box(\Sh(X)) \simeq
  \Sh(X_\Box)$}
\end{align*}
By \emph{idempotency}, we mean that internally employing the Kripke--Joyal
semantics to interpret doubly-internal statements is the same as using the
Kripke--Joyal semantics once. However, we do not want to discuss this here any further;
some details can be found in the original article on the stack
semantics~\cite[Lemma~7.20]{shulman:stack}, but the statement given there is not
general enough to justify the second use of idempotency above. For this, one
would have to extend the stack semantics to support internal statements about
locally internal categories like~$\Sh(X_\Box) \hookrightarrow \Sh(X)$ (which
then look like locally small categories from the internal point of view). This
is worthwhile for other reasons too, but shall not be pursued here.

Therefore, we give a more explicit proof. By induction, we are going to prove
that for any open subset~$U \subseteq X$ and any formula~$\varphi$ over~$U$, it
holds that
\[ U \models_X \varphi^\Box \quad\Longleftrightarrow\quad j(U) \models_{X_\Box}
\varphi, \]
where the internal statements are to be interpreted by the Kripke--Joyal
semantics of~$X$ and~$X_\Box$ respectively and~$j$ is the nucleus associated
to~$\Box$. We may assume that any sheaves occurring in~$\varphi$ as domains of
quantifications are in fact~$\Box$-sheaves; we justify this with a separate lemma
below.

The cases~$\varphi \equiv \top$,~$\varphi \equiv (\psi \wedge \chi)$,
and~$\varphi \equiv \bigwedge_i \psi_i$ are trivial. For~$\varphi \equiv \bot$,
the claim is that~$U \models_X \Box\bot$ if and only if~$j(U)
\models_{X_\Box} \bot$. The former means~$U \subseteq j(\emptyset)$ and the
latter means~$j(U) \leq \sup \emptyset = j(\emptyset)$, so the claim follows from
the first two axioms for a nucleus.

We omit the verification of the remaining cases.
\end{proof}

\begin{lemma}Let~$\Box$ be a modal operator. Let~$\varphi$ be a formula.
Let~$\psi \defequiv \varphi^\Box$ be the~$\Box$-translation of~$\varphi$.
Let~$\psi'$ be the formula obtained from~$\psi$ by substituting any occurring
domain of quantification by its~$\Box$-sheafification, as syntactically defined
in Definition~\ref{defn:plus-construction}. Then~$\psi$
and~$\psi'$ are intuitionistically equivalent.
\end{lemma}
\begin{proof}
For any formula~$\varphi$, we denote by~``$\varphi^\boxplus$'' the result of
first applying the~$\Box$-translation to~$\varphi$ and then substituting any
set~$F$ occurring in~$\varphi$ as a domain of quantification by the plus
construction~$F^+$. Recall that for any such~$F$ there is a canonical map~$F
\to F^+,\ x \mapsto [\{x\}]$. We are going to show by induction that for any
formula~$\varphi(x_1,\ldots,x_n)$ in which elements~$x_i\?F_i$ may occur as
terms, it holds that~$\varphi^\Box(x_1,\ldots,x_n)$ is equivalent
to~$\varphi^\boxplus([\{x_1\}],\ldots,[\{x_n\}])$. This suffices to prove the
lemma.

The cases for
\[ \top \quad \bot \quad \wedge \quad \bigwedge \quad \vee \quad \bigvee \quad \implies \]
are trivial. The cases for unbounded~``$\forall$'' and~``$\exists$'' are
trivial as well. The case for~``$=$'' is slightly more interesting; let~$\varphi(x,y)
\equiv (x = y)$. Then we are to show that~$\varphi^\Box(x,y) \equiv \Box(x=y)$
(equality in some set~$F$) is equivalent to~$\varphi^\boxplus([\{x\}],[\{y\}])
\equiv \Box([\{x\}] = [\{y\}])$ (equality in~$F^+$). This follows by the
definition of the plus construction. The case for~``$\in$'' is similar.

Let~$\varphi \equiv (\exists x\?F\_ \psi(x))$, where we have dropped further
variables occurring in~$\psi$ for simplicity. Then we are to show
that~$\varphi^\Box \equiv \Box(\exists x\?F\_ \psi^\Box(x))$ is equivalent
to~$\varphi^\boxplus \equiv \Box(\exists \bar x\?F^+\_ \psi^\boxplus(\bar x))$.
The ``only if'' direction is trivial (set~$\bar x \defeq [\{x\}]$). For the ``if''
direction, we may assume that there exists~$\bar x\?F^+$ such
that~$\psi^\boxplus(\bar x)$, since we want to prove a boxed statement. By
definition of the plus construction, it holds that~$\Box(\speak{$\bar x$ is a
singleton})$. So, again since we want to prove a boxed statement, we may assume
that~$\bar x$ is actually a singleton. Therefore there exists~$x\?F$ such
that~$\bar x = [\{x\}]$ and that~$\psi^\boxplus([\{x\}])$ holds. By the induction
hypothesis, it follows that~$\psi^\Box(x)$. From this the claim follows.

The case for~``$\forall$'' is similar.
\end{proof}

\begin{ex}\label{ex:translation-equivalence}Let~$X$ be a scheme. Let~$f$ be a
global function on~$X$. Let~$\varphi \defequiv \neg(\speak{$f$ \inv})$
and~$\varphi' \defequiv \speak{$f$ nilpotent}$. Then, by Proposition~\ref{prop:cond-zero}, we
have~$\Sh(X) \models (\varphi \Leftrightarrow \varphi')$. But in general, this
does not imply that~$\Sh(X) \models (\varphi^\Box \Leftrightarrow
(\varphi')^\Box)$. Consider for instance the modal operator given by~$\Box\alpha
\defequiv ((\alpha \Rightarrow \notat{x}) \Rightarrow \notat{x})$ associated to a
point~$x \in X$. Then~$\Sh(X) \models (\varphi^\Box \Leftrightarrow
(\varphi')^\Box)$ means that the equivalence~$\varphi \Leftrightarrow \varphi'$
holds at the point~$x$. This is false for~$X = \Spec \ZZ$,~$f = 2$, and~$x =
(2)$, since in the local ring~$\O_{X,x} = \ZZ_{(2)}$, the element $f$ is not invertible
while also not being nilpotent.
\end{ex}

\subsection{Internal proofs of common lemmas}

\begin{lemma}\label{lemma:module-zero-point-neighbourhood}
Let~$X$ be a scheme (or a ringed space). Let~$\F$ be an~$\O_X$-module
of finite type.
\begin{enumerate}
\item Let~$x \in X$ be a point. Then the stalk~$\F_x$ is zero if and
only if~$\F$ is zero on some open neighborhood of~$x$.
\item Let~$A \subseteq X$ be a closed subset. Then the restriction~$\F|_A$ (\ie
the pullback of~$\F$ to~$A$) is zero if and only if~$\F$ is zero on some open
subset of~$X$ containing~$A$.
\end{enumerate}
\end{lemma}
\begin{proof}\emph{Both} statements are simply internalizations of
Example~\ref{ex:module-zero-geometric}, using the modal operators~$\Box =
(\placeholder \vee A^c)$ and~$\Box = ((\placeholder \Rightarrow
\notat{x}) \Rightarrow \notat{x})$.
\end{proof}

\begin{rem}Lemma~\ref{lemma:module-zero-point-neighbourhood} fails if one drops the hypothesis
that~$\F$ is of finite type. Indeed, in this case one cannot reformulate the
condition that~$\F$ is zero in a geometric way.\end{rem}

In a remark after the proof of Proposition~\ref{prop:rank-function-internally},
we promised to present a simpler proof of it once we would have developed the theory for
doing so. We can now follow up on this promise.
\begin{lemma}\label{lemma:gen-family-n}
Let~$X$ be a scheme (or a ringed space). Let~$\F$ be an~$\O_X$-module
of finite type. Let~$x \in X$ be a point. Let~$n$ be a natural number. Then the
following statements are equivalent:
\begin{enumerate}
\item There exists a generating family for~$\F_x$ consisting of~$n$ elements.
\item There exists an open neighborhood~$U$ of~$x$ such that
\[ U \models \speak{there exists a generating family for~$\F$ consisting of~$n$
elements}. \]
\end{enumerate}
\end{lemma}
\begin{proof}Using the modal operator~$\Box$ defined by~$\Box\varphi \defequiv
((\varphi \Rightarrow \notat{x}) \Rightarrow \notat{x})$, we have to show that
the following statements in the internal language are equivalent:
\begin{enumerate}
\item $\speak{there exists a generating family
for~$\F$ consisting of~$n$ elements}^\Box$.
\item $\Box(\speak{there exists a generating family
for~$\F$ consisting of~$n$ elements})$.
\end{enumerate}
By Lemma~\ref{lemma:open-stalk}, the second statement implies the first, since
in a formal spelling of the statement in quotes,
\begin{equation}\label{eqn:finitely-generated}\tag{$\star$}
  \exists x_1,\ldots,x_n\?\F\_
  \forall x\?\F\_
  \exists a_1,\ldots,a_n\?\O_X\_
  x = \textstyle\sum_i a_i x_i,
\end{equation}
no implication signs occur. To show the converse direction,
we may assume that there is a generating family~$y_1,\ldots,y_m\?\F$ for~$\F$
(since~$\F$ is, externally speaking, of finite type). Then
the~$\Box$-translation of the statement that the~$y_i$ generate~$\F$ holds as
well (again by Lemma~\ref{lemma:open-stalk}). Since there is an intuitionistic
proof of
\begin{multline*}
  \speak{$y_1,\ldots,y_m$ generate~$\F$} \Longrightarrow \\
  \bigl(\speak{there exist $x_1,\ldots,x_n\?\F$ which generate~$\F$}
    \Longleftrightarrow \\
    \exists x_1,\ldots,x_n\?\F\_
    \exists A\?\O^{m \times n}\_ \speak{$\vec y = A \vec x$}\bigr),
\end{multline*}
Lemma~\ref{lemma:stalk-open-with-hypothesis} can substitute the non-geometric formula~\eqref{eqn:finitely-generated} by the geometric
formula
\[ \exists x_1,\ldots,x_n\?\F\_ \exists A\?\O^{m \times n}\_ \speak{$\vec
y = A \vec x$} \]
(Lemma~\ref{lemma:stalk-open-with-hypothesis}). Thus the claim follows.
\end{proof}

\begin{lemma}Let~$X$ be a scheme (or a ringed space). Let~$\alpha : \F \to \G$ be
a morphism of~$\O_X$-modules. Let~$\G$ be of finite type and assume
that~$\alpha_x : \F_x \to \G_x$ is surjective for some point~$x \in X$.
Then~$\alpha$ is an epimorphism on some open neighborhood of~$x$.\end{lemma}
\begin{proof}In the presence of generators~$y_1,\ldots,y_n\?\G$, the
non-geometric surjectivity condition ($\forall y\?\G\_ \exists x\?\F\_
\alpha(x) = y$) can be reformulated in a geometric way: $\bigwedge_{i=1}^n
\exists x\?\F\_ \alpha(x) = y_i$. Thus the claim follows by
Lemma~\ref{lemma:stalk-open-with-hypothesis}.\end{proof}


\begin{lemma}\label{lemma:pushforward-finite-type}
Let~$i : A \hookrightarrow X$ be a closed immersion of schemes (or
ringed spaces). Let~$\F$ be an~$\O_A$-module. Then~$i_*\F$ is of finite type if
and only if~$\F$ is of finite type.\end{lemma}
\begin{proof}
Let~$\Box$ be the modal operator defined by~$\Box\varphi \defequiv (\varphi \vee
A^c)$. From the internal perspective, we have a surjective ring homomorphism~$i^\sharp
: \O_X \to \O_A$, where we omit the forgetful functor~$i_*$ from~$\Box$-sheaves
to arbitrary sets in the notation, and an~$\O_A$-module~$\F$. Furthermore, we
may assume that~$\F$ is a~$\Box$-sheaf. We can regard~$\F$ as an~$\O_X$-module
by~$i^\sharp$.

Note that~$A^c \Rightarrow (\F = 0)$, by~$\Box$-separatedness of~$\F$.

We are to show that~$\F$ is a finitely generated~$\O_X$-module if and only if
the~$\Box$-translation of ``$\F$ is a finitely generated~$\O_A$-module'' holds.
In explicit terms, we have to show the equivalence of the following statements:
\begin{enumerate}
\item $\bigvee_{n \geq 0} \exists x_1,\ldots,x_n\?\F\_
  \forall x\?\F\_ \exists a_1,\ldots,a_n\?\O_X\_ x = \sum_i i^\sharp(a_i) x_i$.
\item $\Box(\bigvee_{n \geq 0} \Box(\exists x_1,\ldots,x_n\?\F\_
  \forall x\?\F\_ \Box(\exists b_1,\ldots,b_n\?\O_A\_ \Box(
    x = \sum_i b_i x_i))))$.
\end{enumerate}
It is clear that the first statement implies the second. For the converse
direction, we just have to repeatedly use the observation that~$\Box\varphi$
implies~$\varphi \vee (\F = 0)$ (once for each occurrence of~$\Box$). So in each
step, we either obtain the statement we want or may assume
that~$\F$ is the trivial module, in which case any subclaim trivially follows. By
surjectivity of~$i^\sharp$, we may write any~$b\?\O_A$ as~$b =
i^\sharp(a)$ for some~$a\?\O_X$.
\end{proof}

\begin{lemma}\label{lemma:fp-hom-geometric}
Let~$X$ be a scheme (or a ringed space). Let~$\F$ and~$\G$ be~$\O_X$-modules. Let~$x
\in X$. Then $\HOM_{\O_X}(\F,\G)_x \cong \Hom_{\O_{X,x}}(\F_x,\G_x)$ if~$\F$ is
of finite presentation around~$x$.\end{lemma}
\begin{proof}It suffices to give an intuitionistic proof of the following fact:
The construction~$\Hom_R(M,\placeholder)$ is geometric if~$M$ is a finitely
presented~$R$-module. So assume that~$M$ is the cokernel of a presentation
matrix~$(a_{ij}) \? R^{n \times m}$. Then we can describe the Hom with
any~$R$-module~$N$ as
\[ \Hom_R(M,N) \cong \Bigl\{ x \? N^n \ \Big|\ \bigwedge_{j=1}^m \sum_{i=1}^n a_{ij}
x_i = 0 \? N \Bigr\}, \]
and this construction is patently geometric, as a set comprehension with respect to
a geometric formula.
\end{proof}

\begin{lemma}Let~$X$ be a scheme (or a ringed space). Let~$\F$ be an~$\O_X$-module of finite
presentation. Let~$x \in X$. Then the stalk~$\F_x$ is a finite
free~$\O_{X,x}$-module if and only if~$\F$ is finite locally free on some open
neighborhood of~$x$.\end{lemma}
\begin{proof}The internal statement that~$\F$ is a finite free module is not geometric:
\[ \bigvee_{n \geq 0}
  \exists x_1,\ldots,x_n\?\F\_
  \forall x\?\F\_
  \exists! a_1,\ldots,a_n\?\O_X\_
  x = \textstyle\sum_i a_i x_i. \]
But it can equivalently be reformulated as
\[ \bigvee_{n \geq 0}
  \exists \alpha\?\HOM_{\O_X}(\F,\O_X^n)\_
  \exists \beta\?\HOM_{\O_X}(\O_X^n,\F)\_
  \alpha \circ \beta = \id \wedge \beta \circ \alpha = \id. \]
This reformulation is geometric, therefore it holds at~$x$ if and only if it
holds on some open neighborhood of~$x$. The claim follows since, by the
previous proposition, taking stalks commutes with
calculating~$\HOM_{\O_X}(\F,\placeholder)$ \resp~$\HOM_{\O_X}(\O_X^n,\placeholder)$;
thus the pulled back formula indeed expresses that~$\F_x$ is finite free as
an~$\O_{X,x}$-module.
\end{proof}

\begin{lemma}\label{lemma:torsion-module-generic-stalk}
Let~$X$ be an integral scheme with generic point~$\xi$. Let~$\F$
be a quasicoherent~$\O_X$-module. Then~$\F$ is a torsion module if and only if
its generic stalk~$\F_\xi$ vanishes.
\end{lemma}
\begin{proof}The generic stalk vanishes if and only if the internal
statement~``$(\F = 0)^{\neg\neg}$'' holds. Therefore it suffices to give an
intuitionistic proof of the following internal statement: The module~$\F$ is
torsion if and only if any element of~$\F$ is \notnot zero.

For the ``only if'' direction, let~$x\?\F$ be an arbitrary element. Since~$\F$
is a torsion module, there exists a regular element~$a\?\O_X$ such that~$ax =
0$. Since~$X$ is reduced, regularity is equivalent to not-not-invertibility.
Since we want to verify the~$\neg\neg$-stable statement~``$\neg\neg(x = 0)$'', we
may in fact assume that~$a$ is invertible. Then~$x = 0$ obviously follows.

For the ``if'' direction, let~$x\?\F$ be an arbitrary element; by assumption,~$x$
is \notnot zero. Since~$X$ is integral,
Lemma~\ref{lemma:dense-standard-reflection} is applicable. Therefore there
exists an element~$a\?\O_X$ such that~$a$ is \notnot invertible and such that
invertibility of~$a$ implies~$x = 0$. Since~$\F$ is quasicoherent, for some
natural number~$n$ it holds that~$a^n x = 0$ (Theorem~\ref{thm:qcoh-sheafchar}
below). Since~$a$ is \notnot invertible,
it is regular (see Lemma~\ref{lemma:regular-notnot-invertible} below for a short
and self-contained proof), and therefore~$a^n$ is regular. So~$x \in \F_\tors$.
\end{proof}

By simply using a different modal operator than~``\notnot'', we will -- without
any additional work -- obtain a more general form of this lemma, applicable to
non-integral schemes (see Lemma~\ref{lemma:torsion-module-generic-stalk-generalized}).

\section{Compactness and metaproperties}
\label{sect:compactness}

\subsection{Quasicompactness}

As stated in the introduction, quasicompactness of a space cannot be detected
by the internal language: There cannot exist a formula~$\varphi$ such that a
topological space is quasicompact if and only if~$\Sh(X) \models \varphi$,
since the latter is always a local property on~$X$ while quasicompactness is not.
However, quasicompactness can be characterized by a \emph{metaproperty} of the
internal language.

This result is best stated in a way which does not explicitly refer to a notion
of finiteness. So recall that quasicompactness of a topological space~$X$ can
be phrased in the following way: For any directed set~$I$ and any monotone
family~$(U_i)_{i \in I}$ of open subsets, if~$X = \bigcup_i U_i$ then~$X = U_i$
for some~$i \in I$. As usual, a \emph{directed set} is an inhabited partially
ordered set such that for any two elements there exists a common upper bound.
A family~$(U_i)_{i \in I}$ is \emph{monotone} if and only if~$i \preceq j$
implies~$U_i \subseteq U_j$.

\begin{prop}\label{prop:quasicompact-meta}
Let~$X$ be a topological space. Then~$X$ is quasicompact if and
only if the internal language of~$\Sh(X)$ has the following metaproperty:
For any directed set~$I$ and any monotone family~$(\varphi_i)_{i \in I}$ of
formulas over~$X$,
\[ \Sh(X) \models \bigvee_{i \in I} \varphi_i
  \quad\text{implies}\quad
  \text{for some~$i \in I$, $\Sh(X) \models \varphi_i$}. \]
The monotonicity condition means that~$\Sh(X) \models (\varphi_i \Rightarrow
\varphi_j)$ for any~$i \preceq j$ in~$I$.
\end{prop}

Stated more succinctly, a topological space~$X$ is quasicompact if and only
if~``$\Sh(X) \models$'' commutes with directed~``$\bigvee_{i \in I}$'''s.

\begin{proof}For the ``only if'' direction, let such a family of formulas be
given. Declare~$U_i$ to be the largest open subset of~$X$ where~$\varphi_i$
holds. Then by assumption, the sets~$U_i$ form a monotone family and cover~$X$. By
quasicompactness of~$X$, some single member~$U_i$ covers~$X$ as well, whereby the
corresponding formula~$\varphi_i$ holds on~$X$.

For the ``if'' direction, we observe that a monotone family~$(U_i)$ of open subsets
induces a monotone family of formulas by defining~$\varphi_i \defequiv U_i$,
employing the notational convention set out in
Section~\ref{sect:modalities-geometric-meaning}. This
correspondence is such that~$\Sh(X) \models \bigvee_i \varphi_i$ holds if and
only if~$X = \bigcup_i U_i$ and such that~$\Sh(X) \models \varphi_i$ if and
only if~$X = U_i$. With these observations the claim is obvious.
\end{proof}

\begin{ex}\label{ex:nilpotency-directed}
Let~$X$ be a quasicompact scheme (or quasicompact ringed space).
Let~$f \in \Gamma(X,\O_X)$ be a global function. Let the set of natural
numbers be endowed with the usual ordering. Then the family of formulas given by~$(f^n =
0)_{n \in \NN}$ is monotone. Thus, if it internally holds that~$f$ is
nilpotent, then~$f$ is nilpotent as an element of~$\Gamma(X,\O_X)$ as
well.\end{ex}

\begin{prop}\label{prop:locally-qc}
Let~$X$ be a topological space. Let~$K \subseteq X$ be an open
subset which is \emph{locally quasicompact} in the sense that there exists an open
covering~$X = \bigcup_j U_j$ such that each~$K \cap U_j$ is quasicompact. Then the
internal language of~$\Sh(X)$ has the following metaproperty: For any
directed set~$I$ and any monotone family~$(\varphi_i)_{i \in I}$ of formulas
over~$X$ it holds that
\[ \Sh(X) \models \bigl(K \Rightarrow \bigvee_i \varphi_i\bigr)
  \quad\text{implies}\quad
  \Sh(X) \models \bigvee_i (K \Rightarrow \varphi_i). \]
If additionally for any open subset~$V \subseteq X$ the set~$K \cap V$ is
locally quasicompact in~$V\!$, the following stronger and purely internal
statement holds:
\[ \Sh(X) \models \bigl(K \Rightarrow \bigvee_i \varphi_i\bigr)
  \Longrightarrow
  \bigvee_i (K \Rightarrow \varphi_i). \]
\end{prop}
\begin{proof}Assume that~$\Sh(X) \models (K \Rightarrow \bigvee_i \varphi_i)$.
This is equivalent to~$K \models \bigvee_i \varphi_i$. By the locality of the
internal language, it follows that~$K \cap U_j \models \bigvee_i \varphi_i$ for each~$j$.
Since~$K \cap U_j$ is quasicompact, it follows by Proposition~\ref{prop:quasicompact-meta} that
there exists an index~$i_j \in I$ such that~$K \cap U_j \models \varphi_{i_j}$.
This is equivalent to~$U_j \models (K \Rightarrow \varphi_{i_j})$. In
particular, it holds that~$U_j \models \bigvee_i (K \Rightarrow \varphi_i)$.
Since this is true for any~$j$, it follows that~$X \models \bigvee_i (K
\Rightarrow \varphi_i)$, again by the locality of the internal language.

The second statement is a corollary of the first one.
\end{proof}

\begin{ex}Any retrocompact subset of a scheme is locally quasicompact in the
sense of Proposition~\ref{prop:locally-qc}.\end{ex}

\begin{ex}\label{ex:df-locally-compact}
Let~$X$ be a scheme and~$f \in \Gamma(X,\O_X)$ be a global function.
Then the open set~$D(f) = \{ x \in X \,|\, \text{$f_x$ is invertible in~$\O_{X,x}$}
\}$ is locally quasicompact in the sense of Proposition~\ref{prop:locally-qc}, even in the
stronger sense: Let~$V \subseteq X$ be any open set. Consider a covering~$V = \bigcup_i
U_i$ by open affine subsets~$U_i = \Spec A_i$. Then~$D(f) \cap U_i \cong \Spec
A_i[f^{-1}]$ is quasicompact.\end{ex}

From this example it will trivially follow that the nilradical~$\sqrt{(0)}
\subseteq \O_X$ of a scheme and indeed the radical of any quasicoherent
sheaf of ideals is quasicoherent (Example~\ref{ex:radical-qcoh}). This example is also
pivotal for giving a simple description of the quasicoherator
(Proposition~\ref{prop:quasicoherator-arbitrary-algebra}), which in turn is
needed for an internal understanding of the relative
spectrum (Section~\ref{sect:relative-spectrum}).

\begin{rem}In applications, the open set~$K$ of Proposition~\ref{prop:locally-qc} is often given
as the largest open subset on which some formula~$\psi$ holds. (For instance,
in Example~\ref{ex:df-locally-compact},~$K$ was given by the formula~$\speak{$f$
is invertible in $\O_X$}$.)
Then the conclusion of the proposition is that \emph{assuming that~$\psi$ holds commutes
with directed disjunctions}.\end{rem}

\subsection{Locality}

A stronger condition on a topological space~$X$ than quasicompactness is
locality: A topological space is \emph{local} if and only if for any open
covering~$X = \bigcup_i U_i$ (not necessarily directed) a certain single subset~$U_i$
covers~$X$ as well. For instance, the spectrum of a ring~$A$ is local if and only
if~$A$ is a local ring. Locality has the following characterization as a metaproperty
of~$\Sh(X)$.

\begin{prop}\label{prop:local-meta}Let~$X$ be a topological space. Then~$X$ is local if and
only if the internal language of~$\Sh(X)$ has the following metaproperty:
For any set~$I$ and any family~$(\varphi_i)_{i \in I}$ of
formulas over~$X$, it holds that
\[ \Sh(X) \models \bigvee_{i \in I} \varphi_i
  \quad\text{implies}\quad
  \text{for some~$i \in I$, $\Sh(X) \models \varphi_i$}. \]
In this case, the internal language has additionally the following (weaker) metaproperty: For any
sheaf~$\F$ on~$X$ and any formula~$\varphi(s)$ containing a variable~$s\?\F$,
it holds that
\[ \Sh(X) \models \exists s\?\F\_ \varphi(s)
  \quad\text{implies}\quad
  \text{for some~$s \in \Gamma(X,\F)$, $\Sh(X) \models \varphi(s)$}. \]
\end{prop}
\begin{proof}The proof of the first part is very similar to the proof of
Proposition~\ref{prop:locally-qc}. For the ``only if'' direction of the second part, note
that the antecedent implies that there exist local sections~$s_i \in
\Gamma(U_i,\F)$ such that~$U_i \models \varphi(s_i)$ for some open covering~$X
= \bigcup_i U_i$. By locality of~$X$, one such~$U_i$ suffices to cover~$X$; so
the corresponding section~$s_i$ is actually a global section and verifies~$X
\models \varphi(s_i)$.
\end{proof}

\begin{rem}The second metaproperty stated in the proposition is indeed weaker
than the condition that~$X$ is local. For instance, let~$X$ be a space consisting
of two discrete points. Then~$\Sh(X)$ has the second metaproperty, but~$X$ is
not local.\end{rem}

\subsection{Irreducibility}

In intuitionistic logic, De Morgan's law~$\neg(\alpha \wedge \beta)
\Rightarrow \neg\alpha \vee \neg\beta$ is not generally justified; therefore we
can't use it when working internally to the topos of sheaves on a general scheme~$X$.
The following proposition demonstrates that if~$X$ is irreducible, the law
does hold.

\begin{prop}\label{prop:irreducibility-internally}
A topological space~$X$ is irreducible if and only if the internal
language of~$\Sh(X)$ has the following metaproperty: For any
formulas~$\varphi$ and~$\psi$
\[ \Sh(X) \models \neg(\varphi \wedge \psi)
  \quad\text{implies}\quad
  \Sh(X) \models \neg\varphi \text{ or }
  \Sh(X) \models \neg\psi, \]
and not $\Sh(X) \models \bot$.
Furthermore, in this case the following internal logical principle holds:
\[ \Sh(X) \models \forall \alpha,\beta \? \Omega\_
  \neg(\alpha \wedge \beta) \Rightarrow (\neg\alpha \vee \neg\beta). \]
\end{prop}
\begin{proof}The statement ``$\Sh(X) \models \neg(\varphi \wedge \psi)$'' means
that~$U \cap V = \emptyset$, where~$U$ and~$V$ are the largest open subsets on
which~$\varphi$ respectively~$\psi$ hold. The disjunction ``$\Sh(X) \models
\neg\varphi$ or $\Sh(X) \models \neg\psi$'' means that~$U = \emptyset$ or~$V =
\emptyset$. And ``$\Sh(X) \models \bot$'' is equivalent to~$X = \emptyset$.

Therefore, if~$X$ is irreducible, then the internal language has the claimed metaproperty. The converse
can be seen by instantiating~$\varphi$ and~$\psi$ with the formulas associated
to given open subsets having empty intersection. It then follows that one of
these formulas is false in the internal language; thus the associated subset is
empty.

The stated internal logical principle holds since nonempty open subsets of irreducible spaces are
irreducible.
\end{proof}

\subsection{Internal proofs of common lemmas}

\begin{lemma}Let~$X$ be an irreducible reduced scheme. Then all local
rings~$\O_{X,x}$ are integral domains.\end{lemma}
\begin{proof}It suffices to give a proof of the following statement: Let~$R$ be
a local ring such that elements which are not invertible are nilpotent. Furthermore
assume that~$R$ is reduced. Then~$R$ is an integral domain in the weak sense of
Definition~\ref{defn:integral-domain}.

This proof may, additionally to the rules of intuitionistic logic, use the
classical axiom stated in Proposition~\ref{prop:irreducibility-internally}.

So let arbitrary elements~$x,y \? R$ with~$xy = 0$ be given. Then it is not the
case that~$x$ and~$y$ are both invertible: If they were, their product~$xy$
would be invertible as well, contradicting~$1 \neq 0$. By the classicality
principle, it follows that~$x$ is not invertible or that~$y$ is not invertible.
Thus~$x$ or~$y$ is nilpotent and therefore zero.
\end{proof}

\begin{lemma}Let~$X$ be a scheme (or a ringed space). Let~$(\E_i)_i$ be a
directed system of~$\O_X$-modules such that~$\E \defeq \colim_i \E_i$ is of
finite type. If~$X$ is quasicompact, there is an index~$i$ such that~$\E_i$ is
of finite type and such that the coprojection~$\E_i \to \E$ is an epimorphism
of sheaves of modules.
\end{lemma}

\begin{proof}Since the usual proof of the statement ``if a directed colimit of
modules~$(M_i)_i$ is finitely generated, then so is one of the modules and its
coprojection into~$\colim_i M_i$ is surjective'' is intuitionistic, it can be
applied in the internal universe of~$\Sh(X)$. Hence we have
\[ \Sh(X) \models \bigvee_i \speak{$\E_i$ is finitely generated and~$\E_i \to \E$
is surjective}. \]
Therefore we can conclude by Proposition~\ref{prop:finite-type-and-co} and
by Proposition~\ref{prop:quasicompact-meta}.
\end{proof}

\begin{caveat}There's a lemma stating that on a quasicompact and quasiseparated
scheme, every quasicoherent sheaf of modules is a filtered colimit of finitely
presented sheaves of modules~\stacksproject{07V9}. There's also the
corresponding statement for modules, whose standard proof is intuitionistic:
Any module~$M$ is a filtered colimit of finitely presented modules (namely the
finitely presented modules mapping to~$M$).

However, the stated lemma does not immediately follow by applying the statement
for modules in the internal universe. This only yields that any sheaf of
modules is an \emph{internal} filtered colimit; those are more
general.\footnote{Any small category~$\I$ induces a small
category~$\ul{\I}$ internal to~$\Sh(X)$ in such a way that the category of
diagrams over~$\I$ coincides with the category of internal diagrams
over~$\ul{\I}$ and the corresponding notions of limit and colimit agree.
However, not every internal small category is of the form~$\ul{\I}$. Therefore
internal limits and colimits are more flexible than external ones.

For instance, it's not true that any sheaf of sets is a colimit of a suitable
system of copies of the terminal sheaf. In contrast, internally, any set is a
colimit of a suitable system of copies of the singleton set.}
\end{caveat}

\section{Quasicoherent sheaves of modules}
\label{sect:qcoh}

Recall that an~$\O_X$-module~$\F$ on a ringed space~$X$ is \emph{quasicoherent}
if and only if there exists a covering of~$X$ by open subsets~$U$ such that on
each such set~$U$, there exists an exact sequence
\[ (\O_X|_U)^J \longrightarrow (\O_X|_U)^I \longrightarrow \F|_U \longrightarrow 0 \]
of~$\O_X|_U$-modules, where~$I$ and~$J$ are arbitrary sets (which may depend
on~$U$).

If~$X$ is a scheme, quasicoherence can also be characterized in
terms of inclusions of distinguished open subsets of affines:
An~$\O_X$-module~$\F$ is quasicoherent if and only if for any open affine
subscheme~$U = \Spec A$ of~$X$ and any function~$f \in A$, the canonical map
\[ \Gamma(U,\F)[f^{-1}] \longrightarrow \Gamma(D(f),\F),\
  \tfrac{s}{f^n} \longmapsto f^{-n} s|_{D(f)} \]
is an isomorphism of~$A[f^{-1}]$-modules. Here~$D(f) \subseteq U$ denotes the
standard open subset~$\{ \ppp \in \Spec A \,|\, f \not\in \ppp \}$. Both
conditions can be internalized.

\begin{prop}Let~$X$ be a ringed space. Let~$\F$ be an~$\O_X$-module. Then~$\F$
is quasicoherent if and only if
\[ \Sh(X) \models \exists I,J\ \mathrm{lc}\_ \speak{there exists an
  exact sequence~$\O_X^J \to \O_X^I \to \F \to 0$}. \]
The ``\textnormal{lc}'' indicates that when interpreting this internal statement with the
Kripke--Joyal semantics,~$I$ and~$J$ should only be instantiated with
\emph{locally constant} sheaves.
\end{prop}
\begin{proof} We only sketch the proof.
The translation of the internal statement is that there exists a covering
of~$X$ by open subsets~$U$ such that for each such~$U$, there exist sets~$I,J$
and an exact sequence
\[ (\O_X|_U)^{\ul{J}} \longrightarrow (\O_X|_U)^{\ul{I}} \longrightarrow \F|_U
\longrightarrow 0 \]
where~$\ul{I}$ and~$\ul{J}$ are the constant sheaves associated to~$I$
respectively~$J$. The term~``$(\O_X|_U)^{\ul{I}}$'' refers to the internally
defined free~$\O_X$-module with basis the elements of~$\ul{I}$. By exploiting
that~$\ul{I}$ is a discrete set from the internal point of view (\ie any two
elements are either equal or not), one can show that this is the same
as~$(\O_X|_U)^I$; similarly for~$J$. With this observation, the statement
follows.
\end{proof}

\begin{rem}The restriction to locally constant sheaves is really necessary: The
internal statement~$\Sh(X) \models \exists I,J\_ \speak{there exists an
exact sequence~$\O_X^J \to \O_X^I \to \F \to 0$}$ is true for
\emph{any}~$\O_X$-module~$\F$. This is because the usual proof of the fact that
any module admits a resolution by (not necessarily finite) free modules is
intuitionistically valid and thus also valid in the internal
universe.\end{rem}

I don't know a useful internal characterization of
locally constant sheaves (but see Section~\ref{sect:lc-sheaves}). The
alternative internal condition given by the following
theorem does not need such a characterization.

\begin{thm}\label{thm:qcoh-sheafchar}
Let~$X$ be a scheme. Let~$\F$ be an~$\O_X$-module. Then~$\F$ is
quasicoherent if and only if, from the internal perspective, for any~$f\?\O_X$,
the localized module~$\F[f^{-1}]$ is a sheaf for the modal operator~$(\speak{$f$ \inv}
\Rightarrow \placeholder)$.
\end{thm}

In detail, the internal condition is that for any~$f\?\O_X$, it holds that
\[ \forall s\?\F[f^{-1}]\_
  (\speak{$f$ \inv} \Rightarrow s = 0) \Longrightarrow s = 0 \]
and for any subsingleton~$\S \subseteq \F[f^{-1}]$ it holds that
\[ (\speak{$f$ \inv} \Rightarrow \speak{$\S$ inhabited}) \Longrightarrow
  \exists s\?\F[f^{-1}]\_
  (\speak{$f$ \inv} \Rightarrow s \in \S). \]
Unlike with the internalizations of finite type, finite presentation and
coherence, this condition is \emph{not} a standard condition of commutative
algebra. In fact, in classical logic, this condition is always satisfied --
for trivial logical reasons if~$f$ is invertible, and because~$\F[f^{-1}]$ is
the zero module if~$f$ is not invertible (since~$f$ is nilpotent then, by
Proposition~\ref{prop:neginvnilpotent}).

That this condition in not known in commutative algebra is to be expected:
Quasicoherence is a condition on sheaves of modules, ensuring
that they are locally isomorphic to sheaves of the form~$M^\sim$,
where~$M$ is a plain module. But in commutative algebra, one \emph{only} studies plain
modules (and not sheaves of modules). The quasicoherence condition is imported
into the realm of commutative algebra only by the internal language.

We give the proof of Theorem~\ref{thm:qcoh-sheafchar} below, after first discussing some examples
and consequences. The proof will explain the origin of this condition.
The localized module~$\F[f^{-1}]$ appearing in the theorem
is externally a certain sheaf. If~$f \in \Gamma(U,\O_X)$, then it is the
sheafification of the presheaf on~$U$ given by~$V \mapsto
\Gamma(V,\F)[f^{-1}]$.

\begin{ex}The zero~$\O_X$-module is quasicoherent, since (it and) all
localizations of it are singleton sets from the internal perspective and
thus~$\Box$-sheaves for any modal operator~$\Box$
by Example~\ref{ex:special-sets-sheaves}.\end{ex}

\begin{cor}\label{cor:submodule-qcoh}
Let~$X$ be a scheme. Let~$\F$ be a quasicoherent~$\O_X$-module.
Let~$\G \subseteq \F$ be a submodule. Then~$\G$ is quasicoherent if and only
if
\[ \Sh(X) \models \forall f\?\O_X\_
  \forall s\?\F\_
  (\speak{$f$ \inv} \Rightarrow s \in \G) \Longrightarrow
  \bigvee_{n \geq 0} f^n s \in \G. \]
\end{cor}
\begin{proof}We can give a purely internal proof. Let~$f\?\O_X$.
Since subpresheaves of separated sheaves are separated, the module~$\G[f^{-1}]$
is in any case separated with respect to the modal operator~$\Box$
with~$\Box\varphi \defequiv (\speak{$f$ \inv} \Rightarrow \varphi)$.

Now suppose that~$\G$ is quasicoherent. Let~$f\?\O_X$. Let $s\?\F$ and assume that
if~$f$ were invertible,~$s$ would be an element of~$\G$. Define the
subsingleton~$S \defeq \{ t\?\G[f^{-1}] \,|\, \speak{$f$ \inv} \wedge t=s/1 \}$.
Then~$S$ would be inhabited by~$s/1$ if~$f$ were invertible. Since~$\G[f^{-1}]$
is a~$\Box$-sheaf, it follows that there exists an element~$u/f^n$ of~$\G[f^{-1}]$
such that, if~$f$ were invertible, it would be the case that~$u/f^n = s/1 \in
\G[f^{-1}] \subseteq \F[f^{-1}]$.
Since~$\F[f^{-1}]$ is~$\Box$-separated, it follows that it actually holds that~$u/f^n
= s/1 \in \F[f^{-1}]$. Therefore there exists~$m\?\NN$ such that $f^m f^n s =
f^m u \in \F$. Thus~$f^{m+n} s$ is an element of~$\G$.

For the converse direction, assume that~$\G$ fulfills the stated condition.
Let$f\?\O_X$. Let~$S \subseteq \G[f^{-1}]$ be a subsingleton which would be
inhabited if~$f$ were invertible. By regarding~$S$ as a subset of~$\F[f^{-1}]$,
it follows that there exists an element~$u/f^n \in \F[f^{-1}]$ such that,
if~$f$ were invertible, $u/f^n$ would be an element of~$S$. In particular,~$u$
would be an element of~$\G$. By assumption
it follows that there exists~$m\?\NN$ such that~$f^m u \in G$. Thus~$(f^m u) /
(f^m f^n)$ is an element of~$\G[f^{-1}]$ such that, if~$f$ were invertible, it
would be an element of~$S$.
\end{proof}

\begin{ex}\label{ex:annihilator-qcoh}
Let~$X$ be a scheme and~$s$ be a global section of~$\O_X$. Then the
annihilator of~$s$, \ie the sheaf of ideals internally defined by the
formula
\[ I \defeq \Ann_{\O_X}(s) = \{ t\?\O_X \,|\, st = 0 \} \subseteq \O_X \]
is quasicoherent. To prove this in the internal language it suffices to
verify the condition of Corollary~\ref{cor:submodule-qcoh}.
So let~$f\?\O_X$ and~$t\?\O_X$ be arbitrary and assume~$\speak{$f$ \inv} \Rightarrow t \in I$,
\ie assume that if~$f$ were invertible, then~$st$ would be zero. By
Proposition~\ref{prop:cond-zero} it follows that~$f^n st = 0$ for
some~$n\?\NN$, \ie that~$f^n t \in I$.
\end{ex}

\begin{ex}\label{ex:radical-qcoh} Let~$X$ be a scheme and~$\I \subseteq \O_X$
be a quasicoherent sheaf of ideals.  Then the radical of~$\I$, internally definable
as \[ \sqrt{\I} \defeq \Bigl\{ s\?\O_X \,\Big|\, \bigvee_{n \geq 0} s^n \in \I \Bigr\}, \] is
quasicoherent as well: Let~$f\?\O_X$ and~$s\?\O_X$ be arbitrary and
assume~$\speak{$f$ \inv} \Rightarrow s \in \sqrt{\I}$, \ie assume that if~$f$
were invertible, some power~$s^n$ would be an element of~$\I$. Since
\emph{assuming that~$f$ is invertible commutes with directed disjunctions}
(Example~\ref{ex:df-locally-compact}), it follows that for some natural
number~$n$, it holds that~$\speak{$f$ \inv} \Rightarrow s^n \in \I$. By
quasicoherence of~$\I$, we may deduce that for some natural number~$m$, it
holds that~$f^m s^n \in \I$. Thus~$fs \in \sqrt{\I}$.\end{ex}


\begin{prop}Let~$X$ be a scheme of dimension~$\leq 0$. Then any~$\O_X$-module
is quasicoherent.\end{prop}
\begin{proof}By Corollary~\ref{cor:scheme-dimension-zero}, any
element~$f\?\O_X$ is invertible or nilpotent. Therefore the quasicoherence
condition of Theorem~\ref{thm:qcoh-sheafchar} is trivially satisfied for any~$\O_X$-module.
\end{proof}

\begin{rem}\label{rem:qcoh-in-constructive-mathematics}
In general intuitionistic mathematics -- not inside the internal universe of a
scheme -- the notion of quasicoherence as given by the internal condition of
Theorem~\ref{thm:qcoh-sheafchar}
does not seem to be very interesting: For many important rings, there are few
quasicoherent modules in this sense. For instance, let~$M$ be a module over a
ring~$R$ in which every element is invertible or not invertible. (The
ring~$\ZZ$ is such a ring.) Then~$M$ is quasicoherent if and only if for any~$f
\? R$ which is not invertible, the localized module~$M[f^{-1}]$ is the zero
module, \ie any element of~$M$ is annihilated by some power~$f^n$. As a
concrete example, any~$\ZZ$-submodule of~$\ZZ$ which contains a nonzero element
fails to be quasicoherent.
\end{rem}

Incidentally, the internal condition of Theorem~\ref{thm:qcoh-sheafchar}
provides a way to distinguish the topos of sheaves over any nontrivial
topological, smooth or complex manifold from the little Zariski topos of any
scheme. This is because the topos of sheaves over a smooth manifold satisfies
the condition (referring to the sheaf of smooth functions) if and only if the
manifold is empty, basically because for no number~$n \geq 0$ the singularity
of the function~$x^n e^{1/x}$ can be removed.

\begin{proof}[Proof of Theorem~\ref{thm:qcoh-sheafchar}]
By the well-known characterization of quasicoherence in terms of inclusions of
distinguished open subsets, an~$\O_X$-module~$\F$ is quasicoherent if and only
if for any affine open subset~$U \subseteq X$ and any function~$f \in
\Gamma(U,\O_U)$, the canonical map
\begin{equation}\label{eqn:restr-map}
  \Gamma(U,\F)[f^{-1}] \lra \Gamma(D(f),\F), \ s/f^n \longmapsto
  f^{-n} s|_{D(f)}
\end{equation}
is bijective. We will see that this map is injective for all such~$U$ and~$f$
if and only if from the internal perspective, for any~$f\?\O_X$, the set~$\F[f^{-1}]$ is a
separated presheaf with respect to the modal operator~$(\speak{$f$ \inv}
\Rightarrow \placeholder)$; and we will see that in this
case, the map is additionally surjective for all such~$U$ and~$f$ if the full
sheaf condition is fulfilled.

Since the sheaf~$\F[f^{-1}]$ does not appear in the stated characterization, we
will first reformulate the separatedness and the sheaf condition in terms
of~$\F$ instead of~$\F[f^{-1}]$. To this end, we observe that the separatedness
condition is equivalent to
\begin{equation}\label{eqn:separated}
  \forall f\?\O_X\_ \forall s\?\F\_
  (\speak{$f$ \inv} \Rightarrow s = 0 \? \F) \Longrightarrow
  \bigvee_{n \geq 0} f^n s = 0 \? \F.
\end{equation}
The equivalence can easily be proven in the internal language. The sheaf
condition is equivalent to the conjunction of the separatedness condition and
\begin{multline}\label{eqn:sheaf}
  \forall f\?\O_X\_ \forall \K \subseteq \F\_
  (\speak{$f$ \inv} \Rightarrow \speak{$\K$ is a singleton})
  \Longrightarrow \\
  \bigvee_{n \geq 0} \exists s\?\F\_
  \speak{$f$ \inv} \Rightarrow f^{-n} s \in \K.
\end{multline}
In one direction, a set~$\S \subseteq \F[f^{-1}]$ is given; construct~$\K \defeq \{
s\?\F \,|\, s/1 \in \S \} \subseteq \F$. In the other direction, a set~$\K
\subseteq \F$ is given; construct~$S \defeq \{ s\?\F[f^{-1}] \,|\, \exists
s'\?\F\_ s' \in \K \wedge s = s'/1 \} \subseteq \F[f^{-1}]$. The remaining
details can easily be filled in.

We now interpret the internal statement~\eqref{eqn:separated} with the
Kripke--Joyal semantics. Using the simplification rules, the external meaning
is that for any affine open subset~$U \subseteq X$ and any function~$f \in
\Gamma(U,\O_U)$ the following condition is satisfied: For any section~$s \in
\Gamma(U,\F)$ it should hold that
\[ U \models (\speak{$f$ \inv} \Rightarrow s = 0) \quad\text{implies}\quad
  U \models \bigvee_{n \geq 0} f^n s = 0. \]
The antecedent is equivalent to saying that~$s$ is zero in~$\Gamma(D(f),\F)$.
The consequent is (by quasicompactness of~$U$, see
Example~\ref{ex:nilpotency-directed}) equivalent to saying that for some~$n \geq 0$, the
section~$f^n s$ is zero in~$\Gamma(U,\F)$, \ie that~$s$ is zero
in~$\Gamma(U,\F)[f^{-1}]$. So this condition is precisely the injectivity of
the canonical map~\eqref{eqn:restr-map}.

The external meaning of statement~\eqref{eqn:sheaf} is that for any affine open
subset~$U \subseteq X$ and any function~$f \in \Gamma(U,\O_U)$ the following
condition is satisfied: For any subsheaf~$\K \subseteq \F|_U$ it should hold
that
\begin{multline*}
  U \models (\speak{$f$ \inv} \Rightarrow \speak{$\K$ is a singleton})
  \quad\text{implies} \\
  U \models \bigvee_{n \geq 0} \exists s\?\F\_
  \speak{$f$ \inv} \Rightarrow f^{-n} s \in \K.
\end{multline*}
Given the injectivity of the canonical map~\eqref{eqn:restr-map} (for any
affine open subset, not only~$U$), this condition is equivalent to its
surjectivity: To see that surjectivity is sufficient, let a subsheaf~$\K
\subseteq \F|_U$ verifying the antecedent be given. Since~$\K|_{D(f)}$ is a
singleton sheaf, we can consider its unique section~$u \in \Gamma(D(f),\K)
\subseteq \Gamma(D(f),\F)$. By surjectivity, there exists a preimage, \ie a
fraction~$s/f^n \in \Gamma(U,\F)[f^{-1}]$ such that~$u = f^{-n} s|_{D(f)}$
in~$\Gamma(D(f),\F)$. Thus~$U \models f^{-n}s \in \K$ holds and the consequent
is verified.

To see that surjectivity is necessary, let a section~$u \in \Gamma(D(f),\F)$ be
given. Define a subsheaf~$\K \subseteq \F|_U$ by setting~$\Gamma(V,\K) \defeq \{
u|_V \,|\, V \subseteq D(f) \}$. Then~$\K$ verifies the antecedent. Thus the
consequent holds: There exists an open covering~$U = \bigcup_i U_i$ such that
for each~$i$, there exists a natural number~$n_i$ and a section~$s_i \in
\Gamma(U_i,\F)$ such that~$f^{-n_i} s_i = u$ on~$U_i \cap D(f)$. Without loss of
generality, we may assume that the~$U_i$ are distinguished open subsets~$D(g_i)
\subseteq U$; that they are finite in number; and that the natural
numbers~$n_i$ agree with each other and thus equal some number~$n$. Since~$s_i
= s_j$ in~$\Gamma(U_i \cap U_j \cap D(f), \F)$, injectivity of the canonical
map~\eqref{eqn:restr-map} (on the affine set~$U_i \cap U_j = D(g_i g_j)$)
implies that~$s_i = s_j$ in~$\Gamma(U_i \cap U_j, \F)[f^{-1}]$. Thus for
any indices~$i,j$ there exists a natural number~$m_{ij}$ such that~$f^{m_{ij}} s_i =
f^{m_{ij}} s_j$ in~$\Gamma(U_i \cap U_j, \F)$. We may assume that the
numbers~$m_{ij}$ equal some common number~$m$; thus the local sections~$f^m s_i$
glue to a section~$s \in \Gamma(U,\F)$. The sought preimage of~$u$ is the
fraction~$s/f^{n+m}$, since~$f^{-(n+m)} s|_{D(f)}$ equals~$u$
in~$\Gamma(D(f),\F)$ (as this is true on the covering~$D(f) = \bigcup_i (D(f)
\cap U_i)$).
\end{proof}

\subsection{The quasicoherator for radical ideals}

For applications in Section~\ref{sect:relative-spectrum} about interpreting the
relative spectrum as an internal spectrum, we want to specialize to radical
sheaves of ideals. In particular, we want to describe the \emph{quasicoherator} --
the left adjoint to the inclusion of the quasicoherent radical ideals in the
poset of all radical ideals -- in simple terms.

\begin{caveat}The quasicoherator we refer to does \emph{not} coincide with the
quasicoherator
of~$\O_X$-modules~\stacksproject{077P},~\cite{thomason-trobaugh}, which is the
\emph{right} adjoint to the inclusion of category of
quasicoherent~$\O_X$-modules in the category of all~$\O_X$-modules. We discuss
this in more detail in Example~\ref{ex:quasicoherator-of-ideals} below.
\end{caveat}

\begin{prop}\label{prop:quasicoherator-structure-sheaf}
Let~$X$ be a scheme. Let~$\I \subseteq \O_X$ be a radical ideal.
\begin{enumerate}
\item The ideal~$\I$ is quasicoherent if and only if
\[ \Sh(X) \models \forall s\?\O_X\_ (\speak{$s$ \inv} \Rightarrow s\in\I)
\Rightarrow s\in\I. \]
\item The reflection of~$\I$ in the poset of quasicoherent radical ideals is
the sheaf~$\overline{\I}$ given by the internal expression
\[ \overline{\I} \defeq \{ s\?\O_X \,|\, \speak{$s$ \inv} \Rightarrow s\in\I
\}. \]
\end{enumerate}
\end{prop}
\begin{proof}Both claims can be verified by purely internal reasoning. The
first claim is a straightforward calculation using the characterization given in
Corollary~\ref{cor:submodule-qcoh}. We discuss the second one in more detail.

Firstly, it's obvious that~$\overline{\I}$ contains~$\I$ and
that~$\overline{\I}$ is a radical ideal. To verify that~$\overline{\I}$ is
quasicoherent, let~$s\?\O_X$ be given such that, if~$s$ were invertible,
then~$s$ would be an element of~$\overline{\I}$. Symbolically, we have
\[ \speak{$s$ \inv} \Longrightarrow (\speak{$s$ \inv} \Rightarrow s\in\I), \]
which of course implies
\[ \speak{$s$ \inv} \Longrightarrow s\in\I. \]
This is precisely the condition for~$s$ to be an element of~$\overline{I}$.

To verify that the construction~$\I \mapsto \overline{\I}$ is really left
adjoint to the inclusion, let a quasicoherent radical ideal~$\J$ be given such
that~$\I \subseteq \J$. We have to show that~$\overline{\I} \subseteq \J$. This
is straightforward.
\end{proof}

\begin{ex}\label{ex:quasicoherator-of-ideals}
Let~$X \defeq \AA^1_k = \Spec k[T]$ be the affine line over a field~$k$.
Let~$j : U \defeq \AA^1_k \setminus \{ 0 \} \hookrightarrow X$ be the open inclusion
of the punctured line. Then~$\I \defeq j_! \O_U \hookrightarrow
\O_X$ is the standard example of a radical sheaf of ideals which is not
quasicoherent. The quasicoherator of modules maps~$\I$
to~$(\Gamma(X,\I))^\sim$, so to the zero module. In contrast, the reflection
of~$\I$ in the poset of quasicoherent radical ideals is~$(T)$.
\end{ex}

Generally, the reflection~$\overline{\I}$ of a radical ideal~$\I$ is the unique
radical ideal such that~$\overline{\I}$ is quasicoherent and such that~$D(\I) =
D(\overline{\I})$. Explicitly, it is the subsheaf of~$\O_X$ given by
\[ U \longmapsto \{ f \in \O_X(U) \,|\, 1 \in \I(D(f)) \}. \]

For arbitrary~$\O_X$-algebras~$\A$, the description of the quasicoherator for
radical ideals of~$\A$ is more involved than the description given in
Proposition~\ref{prop:quasicoherator-structure-sheaf}(b), but still
sufficiently explicit for the applications in
Section~\ref{sect:relative-spectrum}.

\begin{prop}\label{prop:quasicoherator-arbitrary-algebra}
Let~$X$ be a scheme. Let~$\A$ be a quasicoherent~$\O_X$-algebra.
Then the reflection of a radical ideal~$\I \subseteq \A$ in the poset of
quasicoherent radical ideals of~$\A$ is given by the internal expression
\[ \overline{\I} \defeq \bigcup_{n \geq 0} \I_n, \]
where~$(\I_n)$ is the family of radical ideals defined recursively by
\begin{align*}
  \I_0 &\defeq \I, \\
  \I_{n+1} &\defeq \textnormal{the radical ideal generated by
  $\{ fs \,|\, f\?\O_X, s\?\A, (\speak{$f$ \inv} \Rightarrow s \in \I_n) \}$}.
\end{align*}
\end{prop}
\begin{proof}We argue internally. The set~$\overline{\I}$ contains~$\I$ and is
a radical ideal, as an ascending union of radical ideals. To verify
that~$\overline{\I}$ is quasicoherent, let~$f\?\O_X$ and~$s\?\A$ be given such
that, if~$f$ were invertible, then~$s$ would be an element of~$\overline{\I}$.
This means that we have
\[ \speak{$f$ \inv} \Longrightarrow \bigvee_{n \geq 0} s \in \I_n. \]
Since assuming that~$f$ is invertible commutes with directed disjunctions
(Example~\ref{ex:df-locally-compact}), there is a natural number~$n$ such that
\[ \speak{$f$ \inv} \Longrightarrow s \in \I_n. \]
Therefore~$fs \in \I_{n+1} \subseteq \overline{\I}$.

Finally, to verify that the construction~$\I \mapsto \overline{\I}$ is indeed
left adjoint to the inclusion of the quasicoherent radical ideals in all
radical ideals, let a quasicoherent radical ideal~$\J$ be given such that~$\I
\subseteq \J$. By induction we can show that~$\I_n \subseteq \J$ for all
natural numbers~$n$. Therefore~$\overline{\I} \subseteq \J$.
\end{proof}

\begin{rem}\label{rem:reflector-single-element}
If the goal was to close a given radical ideal under the condition
\[ \forall s\?\A\_ (\speak{$f$ \inv} \Rightarrow s \in \I) \Longrightarrow fs \in \I, \]
where~$f \? \O_X$ is a fixed element, no infinite iteration would be necessary.
The closure would in this case simply be given by
\[ \overline{\I}^f \defeq \text{the radical ideal generated by the set~$\{ fs
\,|\, s \? \A, (\speak{$f$ \inv} \Rightarrow s \in \I) \}$}. \]\end{rem}

There is also a purely formal description of the reflector, given by
\[ \I \longmapsto \bigcap \{ \J \subseteq \A \,|\,
  \text{$\J$ is a quasicoherent radical ideal such that~$\I \subseteq \J$} \}. \]
Verifying that this construction has the universal property of
the reflector is straightforward. However, it is not sufficiently concrete for
calculations. In particular, we don't see a way to prove the following
corollary without the explicit description given by
Proposition~\ref{prop:quasicoherator-arbitrary-algebra}.

\begin{cor}\label{cor:quasicoherator-meet}
Let~$X$ be a scheme. Let~$\A$ be a quasicoherent~$\O_X$-algebra. Let~$\I$
and~$\J$ be radical ideals of~$\A$. Then~$\overline{\I \cap \J} =
\overline{\I} \cap \overline{\J}$.
\end{cor}
\begin{proof}The claim is not purely formal. As a left adjoint, the reflector
preserves arbitrary suprema (as a map from the poset of all radical ideals into the poset
of all quasicoherent radical ideals); but the claim is that it preserves (finite) intersections.

Since the reflector is monotone, it is clear that~$\overline{\I
\cap \J} \subseteq \overline{\I} \cap \overline{\J}$.

To verify the converse
direction, we show by induction that~$\I_n \cap \J_m \subseteq \overline{\I
\cap \J}$ for all natural numbers~$n$ and~$m$. The base case is trivial,
since~$\I_0 \cap \J_0 = \I \cap \J$. For the induction step let~$x \in \I_{n+1}
\cap \J_m$. Then~$x^\ell = \sum_i f_i s_i$ for some natural number~$\ell$ and
elements~$f_i \? \O_X$, $s_i \? \A$ such that~$\speak{$f_i$ \inv} \Rightarrow
s_i \in \I_n$. In particular we have~$\speak{$f_i$ \inv} \Rightarrow s_i x \in
\I_n \cap \J_m$, so by the induction hypothesis~$\speak{$f_i$ \inv} \Rightarrow
s_i x \in \overline{\I \cap \J}$. This implies~$f_i s_i x \in \overline{\I \cap
\J}$, since~$\overline{\I \cap \J}$ is quasicoherent. Therefore~$x^{\ell+1} \in
\overline{\I \cap \J}$ and thus~$x \in \overline{\I \cap \J}$.
\end{proof}

\begin{rem}\label{rem:quasicoherator-knaster-tarski}
If in the situation of
Proposition~\ref{prop:quasicoherator-arbitrary-algebra} the algebra~$\A$ is not
quasicoherent, the construction~$\I \mapsto \overline{\I}$ is still left
adjoint to the inclusion of the radical sheaves of ideals which satisfy the (then
somewhat unmotivated) internal condition given in
Corollary~\ref{cor:submodule-qcoh} in the poset of all radical sheaves of ideals.
Also Corollary~\ref{cor:quasicoherator-meet} remains valid.
This is even the case if~$X$ is an arbitrary ringed space; in this case,
the proofs of Proposition~\ref{prop:quasicoherator-arbitrary-algebra} and
Corollary~\ref{cor:submodule-qcoh} have to be modified, since then we may not
suppose that assuming that an element of~$\O_X$ is invertible commutes with
directed disjunctions.

Instead, the reflector~$\I \mapsto \overline{\I}$ has to
be characterized by
\[ \overline{\I} \defeq \text{least fixed point of~$P$ above~$\I$}, \]
where~$P$ is the monotone operator on the set of radical ideals which takes a
radical ideal~$\I$ to the radical ideal generated by~$\{ fs \,|\, f\?\O_X,
s\?\A, (\speak{$f$ \inv} \Rightarrow s \in \I) \}$. The existence of these
fixed points is guaranteed by the Knaster--Tarski theorem, which is
intuitionistically valid in the version we need~\cite{bauer:lumsdaine:bourbaki-witt}.

The following proof scheme is useful for verifying properties of the least
fixed point. Let~$\varphi(\J)$ be a statement on radical ideals~$\J$ such that
$\varphi(\sup_i \J_i) \Leftrightarrow \bigvee_i \varphi(\J_i)$ for every
family~$(\J_i)_i$ of radical ideals. If
\[ \varphi(P(\J)) \Longrightarrow \varphi(\J) \]
for all radical ideals~$\J$ containing~$\I$, then~$\varphi(\overline{\I})
\Rightarrow \varphi(\I)$. This proof scheme is a special case of the following
more general scheme, which is also sometimes needed for reasoning about the
least fixed point.

Let~$L$ be a complete partial order. Let~$\alpha$ be a map from the set of
radical ideals to~$L$ such that $\alpha(\sup_i \J_i) = \sup_i \alpha(\J_i)$ for
every family~$(\J_i)_i$ of radical ideals. If
\[ \alpha(P(\J)) \preceq \alpha(\J) \]
for all radical ideals~$\J$ containing~$\I$, then~$\alpha(\overline{\I})
\preceq \alpha(\I)$.
\end{rem}

\begin{rem}The reflector can also be given by the formula
\[ \overline{\I} = \bigcap_\J\,
  \Bigl(\J : \bigcap_{f\?\O_X} (\J : \overline{\I}^f)\Bigr), \]
where~$\overline{\I}^f$ is as in Remark~\ref{rem:reflector-single-element} and
the first intersection is indexed by all radical ideals~$\J \subseteq \A$.
This identity follows by the description of~$\overline{\I}$ as a least fixed
point and the explicit formula for the least fixed point from the proof of its
existence~\cite{bauer:lumsdaine:bourbaki-witt}. It also follows from the
observation that the operation~$\J \mapsto \overline{\J}$ is the nucleus
associated to the intersection of the sublocales given by the nuclei~$\J
\mapsto \overline{\J}^f$, which in turn is evident from the description of the
relative spectrum as a classifying locale given in
Proposition~\ref{prop:local-spectrum-classify}.
\end{rem}

\subsection{Characterizing locally constant sheaves}
\label{sect:lc-sheaves}

We don't think that there is a characterization of locally constant
sheaves in the internal language of an arbitrary topos of sheaves, other than
the following trivial one: A sheaf~$\E$ on a topological space~$X$ (or a locale,
or a site) is locally constant if and only if
\[ \Sh(X) \models \bigvee_M \speak{$\E \cong \ul{M}$}, \]
where the disjunction is over \emph{all sets} and~$\ul{M}$ is the constant
sheaf associated to the set~$M$. Strictly speaking, because of the class-sized
disjunction, this statement is not even well-formed; however one can still make
sense of its Kripke--Joyal translation.

In the special case that~$X$ is a scheme, however, there might be an internal
characterization. We failed to disprove the following speculation:

\begin{speculation}\label{speculation:lc-quasicoherent}
Let~$X$ be a scheme. Let~$\E$ be a sheaf of sets on~$X$.
Then~$\E$ is locally constant if and only if~$\O_X\langle\E\rangle$, the
free~$\O_X$-module on~$\E$ (constructed internally), is quasicoherent.
\end{speculation}

The free module occurring in this speculation is the sheafification of the
presheaf
\[ U \longmapsto \Gamma(U, \O_X)\langle\Gamma(U, \E)\rangle \]
and can also be described as~$f_! f^{-1} \O_X$, where~$f : \operatorname{\acute
Et}(X)
\to X$ is the projection of the étalé space associated to~$\E$
(and~$\operatorname{\acute Et}(X)$ is equipped with a scheme structure by exploiting
that~$\operatorname{\acute Et}(X)$ is locally homeomorphic to~$X$).

The ``only if'' direction of Speculation~\ref{speculation:lc-quasicoherent}
certainly holds; in fact, if~$\E$ is locally constant,
then~$\O_X\langle\E\rangle$ is even locally free. There are the following
indications that the converse might hold.

If~$X$ happens to be local as a topological space, then the converse
holds: Exploiting that in this case~$\Gamma(X, \O_X\langle\E\rangle) \cong
\Gamma(X,\O_X)\langle\Gamma(X,\E)\rangle$ one can show that the canonical
morphism~$\ul{\Gamma(X,\E)} \to \mathcal{E}$ is an isomorphism. Returning to
the general situation, we see that the pullback of~$\E$ to any of the
$\Spec(\O_{X,x})$ is constant if~$\O_X\langle\E\rangle$ is quasicoherent.
Thus~$\E$ is ``constant on all infinitesimal neighborhoods''.

If~$\O_X\langle\E\rangle$ is not only quasicoherent, but
even locally free (locally isomorphic to a module of the form~$\O_X^{\oplus
M}$), then locally we have~$\O_{X,x}\langle\E_x\rangle \cong \O_{X,x}\langle
M\rangle$, so~$\E_x \cong M$, so at least the stalks are locally constant.
Similarly, if~$\O_X\langle\E\rangle$ is of finite presentation, then~$\E$ is
locally constant (with finite stalks).

Finally, let~$j : V \hookrightarrow X$ be the inclusion of an open subset.
Let~$\E$ be~$j_!(1)$, the extension of the terminal sheaf on~$V$ by the
empty set. This sheaf is locally constant iff~$V$ is a clopen subset. Now
furthermore assume that~$X$ is integral. In this case one can check
that~$\O_X\langle\E\rangle = j_!(\O_V)$ (extension by zero) is quasicoherent
iff~$V$ is a clopen subset. Thus the converse holds in this case.


\section{Rational functions and Cartier divisors}
\label{sect:rational-functions}

\subsection{The sheaf of rational functions} Recall that the sheaf~$\K_X$ of rational
functions on a scheme~$X$ (or a ringed space) can be defined as the sheafification
of the presheaf
\[ \text{$U \subseteq X$ open} \quad\longmapsto\quad \Gamma(U,\O_X)[\Gamma(U,\S)^{-1}], \]
where~$\Gamma(U,\S)$ is the multiplicative set of those sections of~$\O_X$ on~$U$
which are regular in each stalk~$\O_{X,x}$, $x \in U$. Recall also that there are
some wrong definitions in the literature~\cite{kleiman:misconceptions}.

Using the internal language, we can give a simpler definition of~$\K_X$.
Recall that we can associate to any ring~$R$ its total quotient ring, \ie
its localization at the multiplicative subset of regular elements. Since from
the internal perspective~$\O_X$ is an ordinary ring, we can associate to it its
total quotient ring $\O_X[\S^{-1}]$,
where~$\S$ is internally defined by the formula
\[ \S \defeq \{ s\?\O_X \,|\, \speak{$s$ is regular} \} \subseteq \O_X. \]
Externally, this ring is the sheaf~$\K_X$.
\begin{prop}\label{prop:kx-internally}
Let~$X$ be a scheme (or a ringed space). The sheaf of rings defined
in the internal language by localizing~$\O_X$ at its set of regular elements is
(canonically isomorphic to) the sheaf~$\K_X$ of rational functions.
\end{prop}
\begin{proof}Internally, the ring~$\O_X[\S^{-1}]$ has the following
universal property: For any ring~$R$ and any homomorphism~$\O_X \to R$ which
maps the elements of~$\S$ to units, there exists exactly one
homomorphism~$\O_X[\S^{-1}] \to R$ which renders the evident diagram commutative.
\[ \xymatrix{
  \O_X \ar[rr] \ar[dr] && R \\
  & \O_X[\S^{-1}] \ar@{-->}[ru]
} \]
The translation using the Kripke--Joyal semantics gives the following universal
property: For any open subset~$U \subseteq X$, any sheaf of rings~$\R$ on~$U$ and any
homomorphism~$\O_X|_U \to \R$ which maps all elements of~$\Gamma(V,\S)$ for open subsets $V
\subseteq U$ to units, there exists exactly one homomorphism~$\O_X[\S^{-1}]|_U \to
\R$ which renders the evident diagram commutative.
It is well-known that the sheaf~$\K_X$ as usually defined has
this universal property as well.
\end{proof}

\begin{prop}\label{prop:stalks-kx}
Let~$X$ be a scheme (or a ringed space). Then the stalks of~$\K_X$
are given by
\[ \K_{X,x} = \O_{X,x}[\S_x^{-1}]. \]
The elements of~$\S_x$ are exactly the germs of those local sections which are
regular not only in~$\O_{X,x}$, but in all rings~$\O_{X,y}$ where~$y$
ranges over some open neighborhood of~$x$ (depending on the section).\end{prop}
\begin{proof}
Since localization is a geometric construction, the first statement is made entirely
trivial by our framework. The second statement follows since
\[ \Gamma(U,\S) = \{ s\in\Gamma(U,\O_X) \,|\, U \models \speak{$s$ is regular}
\} \]
and since regularity is a geometric implication, so that
$U \models \speak{$s$ is regular}$ if and only if the germ~$s_y$ is regular
in~$\O_{X,y}$ for all~$y \in U$.
\end{proof}

\begin{rem}Speaking internally, the multiplicative set~$\S$ is saturated.
Therefore an element~$s/t \? \K_X$ is invertible in~$\K_X$ if and only if the
numerator~$s$ belongs to~$\S$, that is if~$s$ is an regular element of~$\O_X$.\end{rem}

%

\subsection{Regularity of local functions}
It is well-known that on a locally Noetherian scheme, regularity spreads from
stalks to neighborhoods, that is that a section of~$\O_X$ is regular
in~$\O_{X,x}$ if and only if it is regular on some open neighborhood of~$x$.
This fact has a simple proof in the internal language.
\begin{prop}\label{prop:regularity-spreading}
Let~$X$ be a locally Noetherian scheme. Let~$s \in \Gamma(U,\O_X)$
be a local function on~$X$. Let~$x \in U$. Then the following statements are
equivalent:
\begin{enumerate}
\item The section~$s$ is regular in~$\O_{X,x}$.
\item The section~$s$ is regular in all local rings~$\O_{X,y}$ where~$y$ ranges
over some open neighborhood of~$x$.
\end{enumerate}
\end{prop}
\begin{proof}
Let~$\Box$ be the modal operator defined by~$\Box\varphi \defequiv ((\varphi
\Rightarrow {!x}) \Rightarrow {!x})$. By Corollary~\ref{cor:spreading}, we are
to show that the following statements of the internal language are equivalent:
\begin{enumerate}
\item $(\speak{$s$ is regular})^\Box$, \ie
$\forall t\?\O_X\_ st = 0 \Rightarrow \Box(t = 0)$.
\item $\Box(\speak{$s$ is regular})$, \ie
$\Box(\forall t\?\O_X\_ st = 0 \Rightarrow t = 0)$.
\end{enumerate}
It is clear that the second statement implies the first -- in fact, this is true
without any assumptions on~$X$: Let~$t\?\O_X$ be such that~$st = 0$. Since we want to
prove the boxed statement~$\Box(t=0)$, we may assume that~$s$ is regular and
prove~$t = 0$. This is immediate. (This direction also follows simply by
examining the logical form and applying Lemma~\ref{lemma:open-stalk}.)

For the converse direction, consider the annihilator of~$s$, \ie the ideal
\[ I \defeq \Ann_{\O_X}(s) = \{ t\?\O_X \,|\, st = 0 \} \subseteq \O_X. \]
This ideal satisfies the quasicoherence condition (this was Example~\ref{ex:annihilator-qcoh}),
thus~$I$ is a quasicoherent submodule of a finitely generated module. Since~$X$ is
locally Noetherian, it follows that~$I$ is finitely generated as well, say by~$x_1,\ldots,x_n \? I$. By
assumption, each generator~$x_i \? I$ fulfills~$\Box(x_i = 0)$. Since we want
to prove a boxed statement, we may in fact assume~$x_i = 0$. Thus~$I = (0)$ and
the assertion that~$s$ is regular follows.
\end{proof}

The proof critically depends on the ideal~$I$ being finitely
generated, since a modal operator need only commute with finite
conjunctions. Intuitively, each time we use the modus ponens rule~$(\Box\varphi \wedge
(\varphi \Rightarrow \psi)) \Rightarrow \Box\psi$, we restrict to a smaller open
neighborhood of~$x$. Since infinite intersections of open sets need not be
open, we cannot expect an infinitary modus ponens rule to hold.

\begin{cor}Let~$X$ be a locally Noetherian scheme. Then the stalks~$\K_{X,x}$
of the sheaf of rational functions are given by the total quotient rings of the
local rings~$\O_{X,x}$.\end{cor}
\begin{proof}Follows by combining Proposition~\ref{prop:stalks-kx} and
Proposition~\ref{prop:regularity-spreading}.\end{proof}

\subsection{Normality}\label{sect:normality}
Recall that a ring~$R$ is \emph{normal} if and only if
it is integrally closed in its total quotient ring. Recall also that a
scheme~$X$ (or a ringed space) is \emph{normal} if and only if all
rings~$\O_{X,x}$ are normal.

\begin{prop}\label{prop:normal-int-ext}A locally Noetherian scheme is normal if and only if the
ring~$\O_X$ is normal from the internal perspective.\end{prop}
\begin{proof}The condition of normality can be put into a form which is almost
a geometric implication:
\begin{multline*}
  \forall s,t\?\O_X\_
  \bigl(\speak{$t$ regular} \wedge {} \\
  \exists a_0,\ldots,a_{n-1}\?\O_X\_
  s^n + a_{n-1} t s^{n-1} + \cdots + a_1 t^{n-1} s + a_0 t^n = 0\bigr)
  \Longrightarrow \\
  \exists u\?\O_X\_ s = ut.
\end{multline*}
The only non-geometric subpart is the condition on~$t$ to be regular. However,
by Proposition~\ref{prop:regularity-spreading}, for the purposes of comparing
its truth at points \vs on neighborhoods, it behaves just like a geometric
formula. Therefore the claim follows.
\end{proof}

\subsection{Geometric interpretation of rational functions} Recall that on
integral schemes, rational functions (\ie sections of~$\K_X$) are the same
thing as regular functions defined on dense open subsets. This amounts to
saying that~\emph{$\K_X$ is the~$\neg\neg$-sheafification of~$\O_X$}
(see Proposition~\ref{prop:negneg-sheaves}). We want to rederive this result,
as far as possible in the internal language, and generalize it to arbitrary
(not necessarily locally Noetherian) schemes.

\begin{lemma}\label{lemma:regular-notnot-invertible}Let~$X$ be a reduced scheme. Then:
\begin{enumerate}
\item $\O_X$ is~$\neg\neg$-separated.
\item Internally, an element~$s\?\O_X$ is regular
if and only if it is \notnot invertible.
\end{enumerate}
\end{lemma}
\begin{proof}Recall from Corollary~\ref{cor:field-reduced} that
\begin{equation}\label{eqn:field-condition}\tag{${\!\Join\!}$}
  \Sh(X) \models \forall s\?\O_X\_ \neg(\speak{$s$ invertible}) \Leftrightarrow
  s=0.
\end{equation}
From this we can deduce that~$\O_X$ is~$\neg\neg$-separated:
Assume~$\neg\neg(s=0)$ for~$s\?\O_X$. If~$s$ were invertible, we would
have~$\neg\neg(1=0)$ and thus~$\bot$. Therefore~$s$ is not invertible and thus
zero.

For the ``only if'' direction of the second statement,
note that a regular element is not zero (if it were, then the true statement~$0
\cdot 0 = 0 \cdot 1$ would imply the false statement~$0 = 1$) and thus \notnot
invertible (by the contrapositive of equivalence~\eqref{eqn:field-condition}). For the ``if''
direction, let~$st = 0$ in~$\O_X$. Since~$s$ is \notnot invertible, it follows
that~$t$ is \notnot zero. Since~$\O_X$ is~$\neg\neg$-separated, this implies
that~$t$ really is zero.
\end{proof}

For the following, we need two technical conditions. Say that an affine
scheme~$\Spec A$ has property~$(\star)$ if and only if:
\begin{quote}
Every open dense subset~$U \subseteq \Spec A$ contains a
\emph{standard open} dense subset.
\end{quote}
Say that~$\Spec A$ has property~$(\star\star)$ if and only if:
\begin{quote}
Every open scheme-theoretically dense subset~$U \subseteq \Spec A$ contains a
\emph{standard open} scheme-theoretically dense subset.
\end{quote}
The first condition is satisfied if~$A$ is an irreducible ring (\ie if~$\Spec A$
is irreducible) or more generally if~$A$ contains only finitely many minimal
prime ideals. Both conditions are satisfied if~$A$ is integral or if~$A$ is
Noetherian; for convenience, we give a proof in the
latter case.

\begin{prop}Let~$A$ be a Noetherian ring. Then~$\Spec A$ has properties~$(\star)$
and~$(\star\star)$.
\end{prop}
\begin{proof}Recall that, under the Noetherian hypothesis, an open subset of~$\Spec A$ is dense if and only if it
contains all minimal prime ideals (this fact holds more generally if there are only finitely many minimal prime ideals) and that it
is scheme-theoretically dense if and only if it contains all associated prime
ideals. There are only a finite number of these prime ideals. Therefore the
claim is reduced to the following statement:

Let~$\ppp_1,\ldots,\ppp_n$ be a
finite number of points of an open subset~$U \subseteq \Spec A$. Then there
exists a standard open subset~$D(f) \subseteq U$ which also contains these
points.

The proof of this statement is a direct application of the prime
avoidance lemma.
\end{proof}

\begin{prop}
\label{prop:kx-is-negneg-sheafification}
Let~$X$ be a reduced scheme. Assume that every open affine subscheme has
property~$(\star)$. (For instance, this condition is satisfied if~$X$ is
integral, the set of irreducible components is locally finite, or if~$X$ is
locally Noetherian.) Then~$\K_X$ is the~$\neg\neg$-sheafification
of~$\O_X$.\end{prop}
\begin{proof}
We first show that~$\K_X$ is~$\neg\neg$-separated,
so assume~$\neg\neg(a/s = 0)$ for~$a/s \? \K_X$. Since~$\K_X$ is obtained
from~$\O_X$ by localizing at regular elements, the fraction~$a/s$ vanishes
in~$\K_X$ if and only if~$a = 0$ in~$\O_X$. Thus it follows that~$\neg\neg(a =
0)$ in~$\O_X$ and therefore~$a = 0$ in~$\O_X$; in particular, $a/s = 0$ in~$\K_X$.

We defer the proof that~$\K_X$ is a~$\neg\neg$-sheaf to the end and first
verify the universal property of~$\neg\neg$-sheafification.
So let~$G$ be a~$\neg\neg$-sheaf and let~$\alpha : \O_X \to G$ be a map. We
define an extension~$\bar\alpha : \K_X \to G$ in the following way:
Let~$f \? \K_X$. Define the subsingleton~$S \defeq \{ x \? G \,|\, \exists
b\?\O_X\_ f = b/1 \wedge x = \alpha(b) \} \subseteq G$. Since~$f$ can be
written in the form~$a/s$ with~$s$ \notnot invertible, it follows that~$S$
is \notnot inhabited. Since~$G$ is a~$\neg\neg$-sheaf, there exists a
unique~$x\?G$ such that~$\neg\neg(x \in S)$. We declare~$\bar\alpha(f)$ to be
this~$x$. It is straightforward to check that the composition~$\O_X \to \K_X
\to G$ equals~$\alpha$ and that~$\bar\alpha$ is unique with this property.

Up to this point, the proof did not need that~$X$ is a scheme -- it was enough
for~$X$ to be a ringed space such that equivalence~\eqref{eqn:field-condition} holds and
such that~$\neg(0 = 1)$ in~$\O_X$. Only now, in showing that~$\K_X$ is
a~$\neg\neg$-sheaf, the scheme condition enters. To this end, we first
reformulate the sheaf condition in a way such that it only refers to~$\O_X$,
not~$\K_X$: The quotient ring~$\K_X$ is a~$\neg\neg$-sheaf if and only if
\begin{multline*}
  \Sh(X) \models \forall T \subseteq \O_X\_
  \speak{$T$ is a subsingleton} \wedge \neg\neg(\speak{$T$ is inhabited})
  \Longrightarrow \\
  \exists a,b\?\O_X\_ \speak{$b$ is regular} \wedge \neg\neg(b^{-1} a \in T).
\end{multline*}
This is done just as in the proof of Theorem~\ref{thm:qcoh-sheafchar}.
The expression~``$b^{-1}$'' refers to the inverse of~$b$ which indeed exists in a doubly
negated context, since~$b$ is assumed regular. More explicitly, we should write
\[ \neg\neg(\exists c\?\O_X\_ bc = 1 \wedge ca \in T)
  \quad\text{instead of}\quad
  \neg\neg(b^{-1} a \in T). \]
To verify the Kripke--Joyal interpretation of the rewritten sheaf condition, let
an affine open subset~$U = \Spec A \subseteq X$ having property~$(\star)$ and a subsheaf~$T
\hookrightarrow \O_X|_U$ be given such that~$T$ is internally a subsingleton
and \notnot inhabited. We may glue the unique germs in the inhabited
stalks of~$T$ to obtain a section~$s \in \Gamma(V,\O_X)$ where~$V \subseteq U$
is a dense open subset. Since~$U$ has property~$(\star)$, we may assume that~$V
= D(f)$ is a standard open subset. Because~$V$ is dense and~$A$ is reduced, the
function~$f$ is a regular element of~$A$.
Since~$\Gamma(V,\O_X) = A[f^{-1}]$, we can write~$s = a/f^n$ with~$a \in A$
and~$n \geq 0$.

By Lemma~\ref{lemma:regular-affine}, the function~$b \defeq f^n$ is also regular as an
element of~$\O_U$ from the internal point of view. The function~$b$ is invertible
on~$V$, since~$V = D(f) = D(b)$. It follows that on the dense
open subset~$V \subseteq U$, the sections~$s$ and~$b^{-1} a$ agree.
This observation concludes the proof.
\end{proof}

\begin{cor}Let~$X$ be a reduced scheme such that any open affine subscheme has
property~$(\star)$. Then~$\K_X$ is the result of
pulling back~$\O_X$ to the sublocale~$X_{\neg\neg}$ and then pushing forward
again. If~$X$ is irreducible with generic point~$\xi$, then~$\K_X$ is the
constant sheaf associated to the set~$\O_{X,\xi}$.\end{cor}
\begin{proof}Recall from Section~\ref{sect:internal-sheaves} that pulling back
to~$X_{\neg\neg}$ is equivalent to sheafifying with respect to the double
negation modality; and that pushing forward is equivalent to forgetting the
sheaf property. Therefore the first statement holds.

For the second statement, recall from Lemma~\ref{lemma:negneg-generic-point} that the
sublocale~$X_{\neg\neg}$ is given by the subspace~$\{\xi\}$; that the
sheafification functor~$\Sh(X) \to \Sh(\{\xi\}) \simeq \Set$ is given by
calculating the stalk at~$\xi$; and that the inclusion functor~$\Set \simeq
\Sh(\{\xi\}) \hookrightarrow \Sh(X)$ is given by the constant sheaf
construction.
\end{proof}

If~$X$ is a general scheme (not necessarily reduced),
we can describe~$\K_X$ in a similar way as a sheafification
of~$\O_X$; specifically, it is the sheafification with respect to the modal
operator defined by
\[ \sdense\varphi \defequiv \speak{$\O_X$ is~$(\varphi \Rightarrow
\placeholder)$-separated} \]
in the internal language of~$\Sh(X)$, \ie
\[ \sdense\varphi \defequiv (\forall s\?\O_X\_ (\varphi \Rightarrow s = 0)
\Rightarrow s = 0). \]
This modal operator has an explicit scheme-theoretic description.

\begin{lemma}\label{lemma:scheme-theoretical-denseness}
Let~$U$ be an open subset of a scheme~$X$. Then~$\Sh(X) \models
\sdense U$ if and only if~$U$ is scheme-theoretically dense in~$X$.
\end{lemma}
\begin{proof}We have the following chain of equivalences.
\begin{align*}
  &\ X \models \sdense U \\
  \Longleftrightarrow&\
    \speak{$\O_X$ is~$(U \Rightarrow \placeholder)$-separated} \\
  \Longleftrightarrow&\
    X \models \speak{$\O_X \to \O_X^{+}$ is injective} \\
  &\qquad\qquad\text{(where the plus construction is wrt.\@ the modality~$(U \Rightarrow \placeholder)$)} \\
  \Longleftrightarrow&\
    X \models \speak{$\O_X \to \O_X^{++}$ is injective} \\
  &\qquad\qquad\text{(by the factorization~$\O_X \to \O_X^{+} \to \O_X^{++}$)} \\
  \Longleftrightarrow&\
    \text{the canonical morphism~$\O_X \to j_* \O_U$
    (with $j : U \hookrightarrow X$) is injective} \\
  \Longleftrightarrow&\
    \text{$U$ is scheme-theoretically dense in~$X$.} \qedhere
\end{align*}
\end{proof}

Using the internal language of a scheme, talking about scheme-theoretically
dense open subsets is therefore just as easy as talking about ordinary
topologically dense open subsets; the difference simply amounts to using the
modal operator~$\sdense$ instead of~``\notnot''.

\begin{prop}\label{prop:kx-is-box-sheafification}
Let~$X$ be a ringed space. Then:
\begin{enumerate}
\item The operator~$\sdense$ fulfills the axioms for a modal operator.
\item $\O_X$ is~$\sdense$-separated.
\item $\K_X$ is~$\sdense$-separated.
\item Internally, it holds that~$\sdense(\speak{$f$ \inv})$ implies that~$f$ is
regular for any~$f\?\O_X$.
\end{enumerate}
Suppose furthermore that~$X$ is a scheme. Then:
\begin{enumerate}
\addtocounter{enumi}{4}
\item The converse in~(4) holds.
\item If every open affine subscheme of~$X$ has
property~$(\star\star)$, then~$\K_X$ is the~$\sdense$-sheafification of~$\O_X$.
\end{enumerate}
\end{prop}
\begin{proof}The first four properties are entirely formal; we thus skip over
some details. For the first property, we verify the second axiom on a modal
operator. So we assume~$\sdense\sdense\varphi$ and have to show~$\sdense\varphi$. To
this end, let~$s\?\O_X$ be arbitrary such that~$\varphi \Rightarrow (s=0)$; we
have to prove that~$s = 0$. If~$\O_X$ were separated with respect to the modal
operator~$(\varphi \Rightarrow \placeholder)$, it would follow that~$s = 0$. So
unconditionally it holds that~$\sdense\varphi \Rightarrow (s=0)$. Since by
assumption~$\O_X$ is~$(\sdense\varphi \Rightarrow \placeholder)$-separated, the claim follows.

For the second property, let~$s\?\O_X$ be arbitrary such that~$\sdense(s = 0)$.
Obviously it holds that~$(s = 0) \Rightarrow (s = 0)$. Thus, since~$\O_X$ is
separated with respect to~$((s = 0) \Rightarrow \placeholder)$, it follows
that~$s = 0$. The proof of the third property is similar.

For the fourth property, assume~$\sdense(\speak{$f$ \inv})$ and let~$h\?\O_X$ be
arbitrary such that~$fh = 0$. Then, trivially, it holds that~$\speak{$f$ \inv}
\Rightarrow h = 0$. Since~$\O_X$ is separated with respect to~$(\speak{$f$
\inv} \Rightarrow \placeholder)$, it follows that~$h = 0$.

We now suppose that~$X$ is a scheme. To verify the fifth property, let a
regular element~$f\?\O_X$ be given. We have to show that~$\O_X$ is separated
with respect to the modality~$(\speak{$f$ \inv} \Rightarrow \placeholder)$. So
assume that~$\speak{$f$ \inv} \Rightarrow (s = 0)$ for some~$s\?\O_X$. By
Proposition~\ref{prop:cond-zero} it follows that~$f^n s = 0$ for some natural
number~$n$. Since~$f$ is regular, we may conclude that~$s = 0$.

The verification of the universal property of~$\K_X$ is done analogously as in
the case that~$X$ is reduced: For the proof of
Proposition~\ref{prop:kx-is-negneg-sheafification}, it was critical that
regular elements of~$\O_X$ are \notnot invertible. We now need (and have) that
regular elements of~$\O_X$ are~$\sdense(\speak{invertible})$.

Thus it only remains to verify that~$\K_X$ is a~$\sdense$-sheaf. We may again imitate
the proof of Proposition~\ref{prop:kx-is-negneg-sheafification}; using the same
notation, we may now suppose that~$V$ is a standard open subset such that~$U \models \sdense
V$ (previously, we supposed that~$U \models \neg\neg V$). The proof that the
denominator~$b$ is regular (as seen from the internal perspective, as an
element of~$\O_U$) now goes as follows: We have~$V \subseteq
D(b)$. Therefore~$U \models \sdense V$ implies~$U \models \sdense(\speak{$b$ \inv})$. By
the fourth property, it follows that~$U \models \speak{$b$ is regular}$.
\end{proof}

\begin{rem}The modal operator~$\sdense$ is the largest (weakest) operator such
that~$\O_X$ is~$\sdense$-separated, \ie if~$\sdenseother$ is any modal operator
such that~$\O_X$ is~$\sdenseother$-separated, then~$\sdenseother\varphi
\Rightarrow \sdense\varphi$ for any proposition~$\varphi$.\end{rem}

In the special case that~$X$ is a reduced scheme,
Proposition~\ref{prop:kx-is-box-sheafification} recovers
the result of Proposition~\ref{prop:kx-is-negneg-sheafification}:

\begin{prop}\label{prop:relation-sdense-notnot}
Let~$X$ be a scheme. Then~$\Sh(X) \models \forall\varphi\?\Omega\_
\sdense\varphi \Rightarrow \neg\neg\varphi$.
The converse holds if~$X$ is reduced, so that in
this case the modal operator~$\sdense$ coincides with the double negation modality.\end{prop}
\begin{proof}We argue internally. Let~$\varphi$ be an arbitrary truth value and
assume that~$\sdense\varphi$. The negation~$\neg\varphi$ (which is defined as~$\varphi \Rightarrow \bot$) is
equivalent to~$\varphi \Rightarrow (1 =
0)$. Since by assumption~$\O_X$ is separated with respect to the~$(\varphi
\Rightarrow \placeholder)$-modality, this in turn is equivalent to~$1 = 0 \?
\O_X$, \ie to~$\bot$. Thus~$\neg\neg\varphi$.

For the converse direction, let~$\varphi \Rightarrow (s = 0)$ for some~$s\?\O_X$;
we have to show that in fact~$s = 0$. Since by assumption~$\neg\neg\varphi$, it
follows that~$s$ is \notnot zero. Since~$X$ is reduced,~$\O_X$
is~$\neg\neg$-separated, so this implies that~$s$ is really zero.
\end{proof}

As a corollary, we can reprove the following basic lemma about
scheme-theoretical denseness.
\begin{lemma}Let~$U$ be an open subset of a scheme~$X$. If~$U$ is
scheme-theoretically dense, then~$U$ is also dense in the plain topological
sense. The converse holds if~$X$ is reduced.\end{lemma}
\begin{proof}The set~$U$ is scheme-theoretically dense if and only if~$\Sh(X)
\models \sdense U$ and is dense if and only if~$\Sh(X) \models \neg\neg U$.
Therefore the claim follows from Proposition~\ref{prop:relation-sdense-notnot}.
\end{proof}

\begin{prop}\label{prop:kx-ass}
Let~$X$ be a scheme such that any open affine subscheme has
property~$(\star\star)$. Then~$\K_X$ is the result of
pulling back~$\O_X$ to the sublocale~$X_\sdense$ associated to the modal
operator~$\sdense$ and then pushing forward again. If~$X$ is locally Noetherian,
this sublocale is the subspace of associated points in~$X$.
\end{prop}

In formulas, the proposition states that the canonical map
\[ \K_X \longrightarrow i_* i^{-1} \O_X \]
is an isomorphism, where~$i : X_\sdense \hookrightarrow X$ is the inclusion of
the sublocale~$X_\sdense$. This result requires a cover with
property~$(\star\star)$, but no Noetherian hypothesis.

\begin{proof}The first statement follows trivially by the results of
Section~\ref{sect:internal-sheaves} and the fact that~$\K_X$ is
the~$\sdense$-sheafification of~$\O_X$.

For the second statement, we need to verify that the nucleus~$j_{\Ass(\O_X)}$
associated to the subspace of associated points coincides with the
nucleus~$j_\sdense$ associated to the modal operator~$\sdense$. Recall from
Subsection~\ref{sect:subspace-to-modal-operator} that the latter is given by
\begin{align*}
  j_\sdense(U) &= \text{largest open subset of~$X$ on which~$\sdense U$ holds} \\
  &= \bigcup\ \{ V \subseteq X \ |\
  \text{$V$ open},\ V \models \sdense U \}
\intertext{and note that the former is given by}
  j_{\Ass(\O_X)}(U) &= \bigcup\ \{ V \subseteq X \ |\
  \text{$V$ open},\ V \cap \Ass(\O_X) \subseteq U \}.
\end{align*}
This is a general fact of locale theory, not depending on particular properties
of~$\Ass(\O_X)$. To verify this, one needs to check that~$j_{\Ass(\O_X)}$ is indeed a
nucleus and that the canonical map
\[ \{ U \in \Open(X) \,|\, j_{\Ass(\O_X)}(U) = U \} \longrightarrow \Open(\Ass(\O_X)),\ U \longmapsto \Ass(\O_X) \cap U \]
is an isomorphism of frames with inverse given by~$\Ass(\O_X) \cap U \mapsto
j_{\Ass(\O_X)}(U)$.

The equivalence thus follows from a standard result on the set of associated
points on locally Noetherian schemes:
\begin{align*}
  &\ V \cap \Ass(\O_X) \subseteq U \\
  \Longleftrightarrow&\
    \Ass(\O_V) \subseteq U \\
  \Longleftrightarrow&\
    \text{$U \cap V$ is scheme-theoretically dense in~$U$} \\
  &\qquad\qquad\text{(this step requires the Noetherian assumption)} \\
  \Longleftrightarrow&\
    V \models \sdense U. \qedhere
\end{align*}
\end{proof}

\begin{lemma}
Let~$X$ be a scheme such that any open affine subscheme has
property~$(\star\star)$. Let~$j : U \hookrightarrow X$ be the inclusion of an
open subset containing the sublocale~$X_\sdense$. (If~$X$ is locally
Noetherian, this is equivalent to requiring that~$U$ contains~$\Ass(\O_X)$.)
Then the canonical morphism~$\K_X \to j_* \K_U$ is an isomorphism.
\end{lemma}
\begin{proof}Write~$i : X_\sdense \hookrightarrow X$ and~$i' : X_\sdense
\hookrightarrow U$ for the inclusions. By Proposition~\ref{prop:kx-ass}, the
sheaf~$\K_X$ is given by~$i_* i^{-1} \O_X$. Similarly, the sheaf~$j_* \K_U$ is
given by~$j_* i'_* i'^{-1} j^{-1} \O_X$. The claim follows since~$j \circ i' =
i$.
\end{proof}

\begin{lemma}\label{lemma:dense-standard-reflection-generalized}
Let~$X$ be a scheme such that any open affine subscheme has
property~$(\star)$ respectively~$(\star\star)$. Then
\[ \Sh(X) \models \forall\varphi\?\Omega\_ \neg\neg\varphi \Longrightarrow \exists f\?\O_X\_
  \neg\neg(\speak{$f$ \inv}) \wedge (\speak{$f$ \inv} \Rightarrow \varphi) \]
respectively
\[ \Sh(X) \models \forall\varphi\?\Omega\_ \sdense\varphi \Longrightarrow \exists f\?\O_X\_
  \sdense(\speak{$f$ \inv}) \wedge (\speak{$f$ \inv} \Rightarrow \varphi). \]
\end{lemma}
\begin{proof}The proof of Lemma~\ref{lemma:dense-standard-reflection} carries
over, \emph{mutatis mutandis}.
\end{proof}

\begin{prop}\label{prop:boolean-dim0-continued}
Let~$X$ be a scheme of dimension~$\leq 0$ such that the set of irreducible components is locally finite or such that~$X$ is locally Noetherian. Then
the internal language of~$\Sh(X)$ is Boolean. (The converse holds as well and
was already stated as Corollary~\ref{cor:boolean-dim0}.)
\end{prop}
\begin{proof}
It suffices to verify the principle of double negation elimination, since the
law of excluded middle is equivalent to it.\footnote{This is a standard fact of
intuitionistic logic. Assume that the principle of double negation elimination
holds. We want to verify the law of excluded middle, so let an arbitrary
formula~$\varphi$ be given. Even intuitionistically it holds
that~$\neg\neg(\varphi \vee \neg\varphi)$. By double negation elimination it
follows that~$\varphi \vee \neg\varphi$.}
So let~$\varphi$ be an arbitrary formula and assume~$\neg\neg\varphi$. By the
previous lemma there exists an element~$f\?\O_X$ such that~$f$ is \notnot
invertible and such that~$(\speak{$f$ \inv} \Rightarrow \varphi)$. Since~$\dim
X \leq 0$, this element is invertible or nilpotent
(Corollary~\ref{cor:scheme-dimension-zero}).  In the first case, we are done.
In the second case, some power~$f^n$ is zero and therefore in particular
\notnot zero. Since~$f$ is \notnot invertible, this implies that \notnot~$1 =
0$. On the other hand~$1 \neq 0$, so we obtain a contradiction; from this
contradiction~$\varphi$ trivially follows.
\end{proof}

\begin{lemma}\label{lemma:torsion-module-generic-stalk-generalized}
Let~$X$ be a locally Noetherian scheme. Let~$\F$ be a
quasicoherent~$\O_X$-module. Then~$\F$ is a torsion module if and only if the
restriction of~$\F$ to~$\Ass(\O_X)$ vanishes.
\end{lemma}
\begin{proof}By Proposition~\ref{prop:kx-ass} and
Lemma~\ref{lemma:dense-standard-reflection-generalized} it suffices to repeat
the proof of Lemma~\ref{lemma:torsion-module-generic-stalk} with
``\notnot'' substituted by~``$\sdense$''.
\end{proof}

\subsection{Cartier divisors} Let~$X$ be a scheme (or a ringed space). Recall
that a \emph{Cartier divisor} on~$X$ is a global section of the sheaf of
groups~$\K_X^\times / \O_X^\times$. This sheaf can be constructed internally, with the
same notation: It is the quotient of the group of invertible elements of the
ring~$\K_X$ by the subgroup of invertible elements of the ring~$\O_X$. So an
arbitrary section of~$\K_X^\times/\O_X^\times$ is internally of the form~$[s/t]$
with~$s,t\?\O_X$ being regular elements; this is a simpler description than the
usual external one as a family~$(f_i)_i$ of functions~$f_i \in
\Gamma(U_i,\K_X^\times)$ such that~$f_i^{-1}|_{U_i \cap U_j} \cdot f_j|_{U_i \cap
U_j} \in \Gamma(U_i \cap U_j, \O_X^\times)$ for all~$i,j$.

We can sketch the basic theory of Cartier divisors completely from the internal
perspective. In accordance with common practice, we write the group
operation of~$\K_X^\times/\O_X^\times$ (which is induced by multiplication of elements
in~$\K_X^\times$) additively.

\begin{defn}\label{defn:effective-cartier-divisor}
A Cartier divisor is \emph{\effective} if and only if, from the
internal perspective, it can be written in the form~$[s/1]$ with~$s\?\O_X$
being a regular element.\end{defn}

Thus a Cartier divisor~$[s/t]$ is \effective if and only if~$s$ is
an~$\O_X$-multiple of~$t$.

\begin{defn}A Cartier divisor~$D$ is \emph{principal} if and only if there
exists a global section~$f \in \Gamma(X,\K_X^\times)$ such that internally,~$D = [f]$.
Two Cartier divisors are \emph{linearly equivalent} if and only if their
difference is a principal divisor.
\end{defn}

Decidedly, principality is a global notion: For \emph{any} divisor~$D$
there exists an open covering~$X = \bigcup_i U_i$ and local sections~$f_i \in
\Gamma(U_i, \K_X^\times)$ such that~$D|_{U_i} = [f_i]$.

\begin{defn}\label{defn:line-bundle-of-divisor}
The \emph{line bundle associated to a Cartier divisor}~$D$
is the~$\O_X$-submodule
\[ \O_X(D) \defeq \{ g \? \K_X \,|\, g D \in \O_X \} = D^{-1} \O_X \subseteq \K_X
\]
of~$\K_X$. Here we are abusing language for~``$gD \in \O_X$'' to mean that~$gf
\in \O_X$ if~$D = [f]$ with~$f\?\K_X$; and for~``$D^{-1} \O_X$'' to
mean~$f^{-1}\O_X$. This condition respectively submodule does not depend on the
representative~$f$, since~$f$ is well-defined up to multiplication by an element
of~$\O_X^\times$.\end{defn}

The submodule~$\O_X(D)$ is indeed locally free of rank~$1$, since
internally~$f^{-1}$ gives a one-element basis. The divisor~$D$ is \effective if
and only if~$\O_X(-D)$ is a subset of~$\O_X$ from the internal perspective
(this comparison makes sense, since~$\O_X(-D)$ and~$\O_X$ are both canonically
embedded in~$\K_X$). In
this case, we can define the \emph{support} of~$D$ to be the closed subscheme
of~$X$ associated to the sheaf of ideals~$\O_X(-D) \subseteq \O_X$.

The line bundle~$\O_X(D)$ can also be written in the familiar form
\[ \O_X(D) = \{ g \? \K_X \,|\, \operatorname{div}(g) + D \geq 0 \}, \]
if we define~``$\operatorname{div}(g)$'' as the equivalence class~$[g] :
\K_X/\O_X^\times$, interpret the left-hand side of the inequality as an element
of~$\K_X/\O_X^\times$, and declare that~$[s/t] \geq 0$ if and only if~$s$ is
an~$\O_X$-multiple of~$t$.

On the other hand, a description like
\[ \text{``$\O_X(D) = \{ 0 \} \cup \{ g \? \K_X^\times \,|\,
  \operatorname{div}(g) + D \geq 0 \}$''} \]
is not possible, since the case distinction necessary for a verification of the
inclusion~``$\subseteq$'' is not possible intuitionistically.


\begin{defn}The \emph{Cartier divisor associated to a free~$\O_X$-submodule~$\L \subseteq
\K_X$ of rank~1} is~$D \defeq [f^{-1}]$, where~$f\?\K_X$ is the unique element of
some one-element basis of~$\L$.\end{defn}

The basis element~$f\?\K_X$ does indeed lie in~$\K_X^\times$: Write~$f
= s/t$ with~$s,t \? \O_X$. It suffices to show that~$s$ is a regular element
of~$\O_X$. So let~$h\?\O_X$ be such that~$sh = 0$ in~$\O_X$. Then in
particular~$hf = 0$ in~$\K_X$. By linear independence, it follows that~$h = 0$
in~$\K_X$ and thus~$h = 0$ in~$\O_X$.

Furthermore, the associated divisor does not depend on the choice of~$f$,
since~$f$ is well-defined up to multiplication by an element of~$\O_X^\times$: If~$f
\O_X = g \O_X \subseteq \K_X$, then there exist elements~$u,v\?\O_X$ such that~$fu = g$
and~$gv = f$ in~$\K_X$. It follows that~$uv = fuvf^{-1} = gvf^{-1} = ff^{-1} =
1$ in~$\K_X$ and thus in~$\O_X$, by injectivity of the localization morphism~$\O_X \to
\K_X$. Therefore~$u$ and~$v$ are elements of~$\O_X^\times$.

\begin{lemma}Let~$D$ and~$D'$ be divisors on~$X$. Then~$\O_X(D) \otimes_{\O_X}
\O_X(D') \cong \O_X(D + D')$.\end{lemma}
\begin{proof}The wanted morphism of sheaves~$\O_X(D) \otimes \O_X(D') \to
\O_X(D + D')$ is given by multiplication. That this is well-defined and an
isomorphism can be checked from the internal point of view, where the claims
are obvious.\end{proof}

\begin{prop}The association~$D \mapsto \O_X(D)$ defines a one-to-one
correspondence between Cartier divisors on~$X$ and rank-one submodules
of~$\K_X$. This correspondence descends to a one-to-one correspondence between
Cartier divisors up to linear equivalence and rank-one submodules of~$\K_X$ up
to isomorphism (as abstract~$\O_X$-modules, ignoring their embedding
into~$\K_X$).\end{prop}
\begin{proof}The first statement is obvious from the definitions. For the
second statement, it suffices to show that~$\O_X(D)$ is isomorphic to~$\O_X$ if
and only if~$D$ is principal. An isomorphism~$\O_X \to \O_X(D)$ gives a
global section~$f \in \Gamma(X,\K_X^\times)$ (by considering the image of the unit element)
such that internally,~$D = [f^{-1}]$; this shows that~$D$ is principal. The
converse is similar.
\end{proof}


For the following definition, recall that we can localize an~$\O_X$-module~$\L$
at the set~$\S \subseteq \O_X$ of regular elements to obtain
a~$\K_X$-module~$\L[S^{-1}]$.

\begin{defn}Let~$f\?\L[\S^{-1}]$ be a rational section of a line bundle~$\L$
on~$X$. Assume that~``$f$ is nontrivial'', that is multiplication by~$f$ is an
injective map~$\O_X \to \L[\S^{-1}]$. Then the \emph{associated divisor} of~$f$
is~$\operatorname{div}(f) \defeq [\psi(s)/t]$, where~$f = s/t$ with~$s\?\L$ and~$t\?\O_X$
and~$\psi : \L \to \O_X$ is an isomorphism.\end{defn}

One can check that~$\psi(s)$ is a regular element of~$\O_X$; this statement is
equivalent to the multiplication map~$\O_X \to \L[\S^{-1}]$ being
injective. Furthermore one can check that~$[\psi(s)/t]$ does not depend on the
choice of~$s$,~$t$, and~$\psi$.

\begin{prop}Let~$f\?\L[\S^{-1}]$ be a nontrivial rational section of a line
bundle~$\L$ on~$X$. Then multiplication by~$f$ induces an
isomorphism~$\O_X(\operatorname{div}(f)) \to \L$.\end{prop}
\begin{proof}The isomorphism should map a rational function~$g$ to~$g f$. This
is a priori an element of~$\L[\S^{-1}]$; we have to check that it can be
regarded as an element of~$\L$. Just as in the definition
of~$\operatorname{div}(f)$, write~$f = s/t$ and fix an isomorphism~$\psi : \L
\to \O_X$. Write~$g = (t/\psi(s)) \cdot h$ for some function~$h\?\O_X$. Then~$g
f = sh/\psi(s) = h\psi^{-1}(1)$, since~$s = \psi^{-1}(\psi(s)) = \psi(s)
\cdot \psi^{-1}(1)$. The element~$h\psi^{-1}(1)$ can indeed be considered as an
element of~$\L$.

Injectivity of the map~$\O_X(\operatorname{div}(f)) \to \L$ is by
nontriviality of~$f$. For surjectivity, we observe that~$(t/\psi(s)) \cdot \psi(v)$ is a
preimage to~$v\?\L$, since~$(t/\psi(s)) \cdot \psi(v) \cdot f = \psi(v) \psi(s)
\psi^{-1}(1) / \psi(s) = v$.
\end{proof}

\begin{prop}Let~$\L$ be a line bundle on~$X$. Assume that~$\L$ can be embedded
into~$\K_X$. Then~$\L$ possesses a nontrivial rational section.
\end{prop}
\begin{proof}Let~$i : \L \to \K_X$ be the given injection. Let~$(v)$ be an
one-element basis for~$\L$. Write~$i(v) = s/t$. Then~$s$ is regular,
since~$hs = 0$ implies~$i(hv) = 0$ and thus~$h = 0$, for any~$h\?\O_X$.
Therefore~$f \defeq tv/s$ is a well-defined element of~$\L[\S^{-1}]$.
Furthermore it is nontrivial in the desired sense: If~$h \cdot (tv/s) = 0$,
then~$htv = 0$, thus~$ht = 0$ and~$h = 0$.

It remains to check that~$f$ is independent of the choice of~$v$ and of the
representation~$i(v) = s/t$; else we defined only local sections which might not
glue to a single nontrivial rational section (externally speaking). This verification is
trivial.
\end{proof}

\begin{prop}Let~$D$ be an \effective divisor on~$X$. Then the complement of its
support is scheme-theoretically dense.\end{prop}
\begin{proof}The complement of the support of~$D$, that is the open
subset~$D(\O_X(-D))$ (where we consider~$\O_X(-D)$ as an ideal of~$\O_X$), is
the truth value of the statement~``$1 \in \O_X(-D)$''. By
Lemma~\ref{lemma:scheme-theoretical-denseness}, we therefore have to verify
that~$\O_X$ is separated with respect to the modal operator~$(1 \in \O_X(-D)
\Rightarrow \placeholder)$.

Let~$s \? \O_X$ be given such that~$1 \in \O_X(-D) \Rightarrow s = 0$; we have
to show that~$s = 0$. Writing~$D = [f/1]$ where~$f \? \O_X$ is a regular
element, this condition is equivalent to~$\speak{$f$ \inv} \Rightarrow s = 0$.
By Proposition~\ref{prop:cond-zero} it follows that~$f^n s = 0$ for some~$n
\geq 0$. Since~$f$ is regular, we may cancel~$f^n$ in this equation.
\end{proof}

\begin{prop}Assume that~$X$ is an integral scheme. Then any line bundle on~$X$
is (uncanonically) a submodule of~$\K_X$.\end{prop}
\begin{proof}Let~$\xi$ be the generic point of~$X$ and let~$\Box \defeq \neg\neg$
denote the modal operator such that internal sheafification with respect
to~$\Box$ is the same as pulling back to~$\{\xi\}$ and then pushing forward
to~$X$ again (see Section~\ref{sect:negneg-sheaves}). Let~$\L$ be a line bundle on~$X$. Since~$\L_\xi \cong
\O_{X,\xi}$ (uncanonically), there is some injection~$\L_\xi \to \K_{X,\xi}$;
this corresponds internally to an injection~$\L^{++} \to \K_X^{++}$.
Since~$\K_X$ is already a~$\Box$-sheaf
(Proposition~\ref{prop:kx-is-negneg-sheafification}) and~$\L$ is~$\Box$-separated
(being isomorphic to~$\O_X$), we have the global injection
\[ \L \lhra \L^{++} \lhra \K_X^{++} \stackrel{({\cong})^{-1}}{\longrightarrow} \K_X. \qedhere \]
\end{proof}

\section{Subschemes}

\subsection{Sheaves on open and closed subspaces} It is well-known that sheaves
defined on open or closed subspaces of a topological space~$X$ can be related
with certain sheaves on~$X$, by using appropriate extension functors. We can
define these functors and show their basic properties in the internal
language. Recall from Section~\ref{sect:modalities-geometric-meaning} that we
have defined a formula~``$U$'' for any open subset~$U \subseteq X$ such that
$V \models U$ if and only if $V \subseteq U$.

\begin{lemma}\label{lemma:extension-by-empty-set}
Let~$X$ be a topological space. Let~$j : U \hookrightarrow X$ be the inclusion
of an open subspace. Then there is a canonical functor~$j_! : \Sh(U) \to
\Sh(X)$ called \emph{extension by the empty set} with the following properties:
\begin{enumerate}
\item The functor~$j_!$ is left adjoint to the restriction functor~$j^{-1} : \Sh(X) \to
\Sh(U)$.
\item The composition~$j^{-1} \circ j_! : \Sh(U) \to \Sh(U)$ is (canonically
isomorphic to) the identity.
\item The essential image of~$j_!$ consists of exactly those sheaves on~$X$
whose stalks are empty at all points of~$U^c$. For those sheaves~$\F$ it holds
that~$j_!j^{-1}\F \cong \F$ (canonically).
\end{enumerate}
\end{lemma}
\begin{proof}Internally, for a set~$\F$, we can define~$j_!(\F)$ simply to be the
set comprehension
\[ j_!(\F) \defeq \{ x\?\F \,|\, U \}. \]
Externally, the sections of the thus defined sheaf on an open subset~$V
\subseteq X$ are given by~$\{ x \in \Gamma(V,\F) \,|\, V \subseteq U \}$,
\ie all of~$\Gamma(V,\F)$ if~$V \subseteq U$ and the empty set otherwise.
With this short internal description, all of the stated properties can be
easily verified in the internal language.

For instance, recall that internally the functor~$j^{-1}$ is given by
sheafifying with respect to the modal operator~$\Box \defequiv (U \Rightarrow
\placeholder)$. Thus, to show the second statement, we have to give a
bijection~$(j_!(\F))^{++} \to \F$ for any~$\Box$-sheaf~$\F$. (This map has to
be given explicitly, to not only show a weaker statement about a local
isomorphism -- see Section~\ref{sect:internal-constructions}). To this end, we can use the composition
\[ (j_!(\F))^{++} \lhra \F^{++} \stackrel{({\cong})^{-1}}{\lra} \F, \]
where the first map is injective since sheafifying is exact. It is also
surjective, since the~$\Box$-translation of the statement~$\speak{$j_!(\F) \to
\F$ is surjective}$ holds: For any element~$x\?\F$, it holds
that~$\Box(\speak{$x$ possesses a preimage})$.

For the third property, we observe that a sheaf~$\F$ on~$X$ fulfills the stated
condition on stalks if and only if, from the internal perspective, it holds
that~$U \Rightarrow \speak{$\F$ is inhabited}$. We omit further details.
\end{proof}

\begin{lemma}\label{lemma:extension-by-zero}
Let~$X$ be a ringed space. Let~$j : U \hookrightarrow X$ be the inclusion
of an open subspace. Then there is a canonical functor~$j_! : \Mod_U(\O_U) \to
\Mod_X(\O_X)$ called \emph{extension by zero} with the following properties:
\begin{enumerate}
\item The functor~$j_!$ is left adjoint to the restriction functor~$j^{-1} :
\Mod_X(\O_X) \to \Mod_U(\O_U)$.
\item The composition~$j^{-1} \circ j_! : \Mod_U(\O_U) \to \Mod_U(\O_U)$ is (canonically
isomorphic to) the identity.
\item The essential image of~$j_!$ consists of exactly those~$\O_X$-modules
whose stalks are zero at all points of~$U^c$. For those sheaves~$\F$ it holds
that~$j_!j^{-1}\F \cong \F$ (canonically).
\end{enumerate}
\end{lemma}
\begin{proof}Internally, a sheaf of modules on~$\O_U$ is simply a module
on~$\O_X^{++}$ which is a~$\Box$-sheaf, where~$\Box \defequiv (U \Rightarrow
\placeholder)$. The suitable internal definition for the extension by zero of
such a module~$\F$ is
\[ j_!(\F) \defeq \{ x\?\F \,|\, (x = 0) \vee U \}. \]
With this description, all necessary verifications are easy. Note that
an~$\O_X$-module~$\F$ fulfills the stated condition on stalks if and only if
internally, it holds that~$\forall x\?\F\_ ((x = 0) \vee U)$.
\end{proof}

\begin{lemma}\label{lemma:essim-closed-immersion}
Let~$X$ be a topological space. Let~$i : A \hookrightarrow X$ be the inclusion
of a closed subspace. The essential image of the
inclusion~$i_* : \Sh(A) \to \Sh(X)$ consists of exactly those sheaves whose support
is a subset of~$A$. For those sheaves~$\F$ it holds that~$i_* i^{-1} \F \cong \F$
(canonically).\end{lemma}
\begin{proof}Recall that the modal operator associated to~$A$ is~$\Box\varphi
\defequiv (\varphi \vee A^c)$, and that by Section~\ref{sect:internal-sheaves} the
essential image of~$i_*$ consists of exactly those sheaves which
are~$\Box$-sheaves from the internal perspective. Let~$\F$ be a sheaf on~$X$.
Then it holds that
\[ \supp\F \subseteq A \quad\Longleftrightarrow\quad
  A^c \subseteq X \setminus \supp\F \quad\Longleftrightarrow\quad
  A^c \subseteq \Int(X \setminus \supp\F). \]
Since the interior of the complement of~$\supp\F$ can be characterized as the
largest open subset of~$X$ on which the internal statement~``$\F$ is a
singleton'' holds (Remark~\ref{rem:support-sheaf-of-sets}), the condition on
the support is fulfilled if and only if
\[ \Sh(X) \models (A^c \Rightarrow \speak{$\F$ is a singleton}). \]
We thus have to show that this internal condition is equivalent to~$\F$ being
a~$\Box$-sheaf. For the ``if'' direction, assume~$A^c$. Then the empty subset~$S
\subseteq \F$ trivially verifies the condition that~$\Box(\speak{$S$ is a
singleton})$. There thus exists an element~$x\?\F$ (such that~$\Box(x \in S)$).
If we're given a further element~$y\?\F$, it trivially holds that~$\Box(x =
y)$. By~$\Box$-separatedness, it thus follows that~$x = y$. Thus~$\F$ is the
singleton~$\{x\}$. The proof of the ``only if'' direction is similar.

The second statement claims that internally, sheafifying a~$\Box$-sheaf with
respect to the modal operator~$\Box$ and then forgetting that the result is
a~$\Box$-sheaf amounts to doing nothing. This is obvious.
\end{proof}

\subsection{Closed subschemes} Let~$X$ be a ringed space. Recall
that a sheaf of ideals~$\I \subseteq \O_X$ defines a closed subset~$V(\I) = \{ x
\in X \,|\, \I_x \neq (1) \subseteq \O_{X,x} \}$, a sheaf of
rings~$\O_X/\I$, and a ringed space~$(V(\I), \O_{V(\I)})$ where~$\O_{V(\I)}$ is
the pullback of~$\O_X/\I$ to~$V(\I)$. In the internal universe, we can
reify~$V(\I)$ by giving a modal operator~$\Box$ such that externally, the
subspace~$X_\Box$ coincides with~$V(\I)$.

\begin{prop}\label{prop:basics-closed-subspace}
Let~$X$ be a ringed space. Let~$\I \subseteq \O_X$ be a sheaf of ideals.
Then:
\begin{enumerate}
\item The subspace of~$X$ associated to the modal operator~$\Box$ defined
by~$\Box\varphi \defequiv (\varphi \vee (1 \in \I))$ is~$V(\I)$.
\item The support of~$\O_X/\I$ is exactly~$V(\I)$.
\item The canonical morphism~$i : V(\I) \to X$ is a closed immersion
of ringed spaces.
\end{enumerate}\end{prop}
\begin{proof}For any open subset~$U \subseteq X$, it holds that~$U \models 1
\in \I$ if and only if~$U \subseteq D(\I) = X \setminus V(\I)$. Thus~$D(\I)$
can be characterized as the largest open subset on which~``$1 \in \I$'' holds.
According to Table~\ref{table:nuclei} on page~\pageref{table:nuclei}, the
stated modal operator thus defines the subspace~$D(\I)^c$, \ie~$V(\I)$.

For the second statement, we observe that since~$\O_X/\I$ is a sheaf of rings, its
support is closed. Therefore the largest open subset of~$X$ where the internal
statement~``$\O_X/\I = 0$'' holds is the complement of the support
(Proposition~\ref{prop:characterization-support}). Since~$D(\I)$ is the largest
open subset where the internal statement~``$\I = (1)$'' holds, it suffices to
show that internally,~$\O_X/\I = 0$ if and only if~$\I = (1)$. This is obvious.

The topological part of the third statement is clear. For the ring-theoretic
part, we have to show that the canonical ring homomorphism~$\O_X \to i_*
\O_{V(\I)}$, that is the canonical projection~$\O_X \to \O_X/\I$, is an
epimorphism of sheaves. This is obvious.
\end{proof}

By Lemma~\ref{lemma:essim-closed-immersion}, the sheaf~$\O_X/\I$ is
thus a~$\Box$-sheaf from the internal perspective.

\begin{prop}Let~$X$ be a locally ringed space. Let~$\I \subseteq \O_X$ be a
sheaf of ideals. Then the ringed space~$(V(\I), \O_{V(\I)})$ is locally
ringed as well.\end{prop}
\begin{proof}We have to show that
\[ \Sh(V(\I)) \models \speak{$\O_{V(\I)}$ is a local ring}. \]
By Theorem~\ref{thm:box-translation-semantically}, this is equivalent to
\[ \Sh(X) \models (\speak{$\O_X/\I$ is a local ring})^\Box, \]
where~$\Box$ is the modal operator given by~$\Box\varphi \defequiv (\varphi \vee
(1 \in \I))$. We therefore have to give an intuitionistic proof of the fact
\[ \forall x,y\?\O_X/\I\_ \speak{$x+y$ \inv} \Longrightarrow
  \Box(\speak{$x$ \inv} \vee \speak{$y$ \inv}). \]
So let~$x = [s], y = [t] \? \O_X/\I$ such that~$x + y$ is invertible
in~$\O_X/\I$. This means that there exists~$u\?\O_X$ and~$v\?\I$ such that~$us
+ ut + v = 1$ in~$\O_X$. Since~$\O_X$ is a local ring, it follows
that~$us$,~$ut$, or~$v$ is invertible. In the first two cases, it follows
that~$x$ respectively~$y$ are invertible in~$\O_X/\I$. In the third case, it
follows that~$1 \in \I$ and thus any boxed statement is trivially true.
\end{proof}

If~$X$ is a scheme and~$\I \subseteq \O_X$ is a sheaf of ideals, it is well-known
that the locally ringed space~$V(\I)$ is a scheme if and only if~$\I$ is
quasicoherent. We cannot give an internal proof of this fact since we lack an
internal characterization of being a scheme.

\begin{lemma}\label{lemma:closed-subspace-reduced}
Let~$X$ be a scheme (or a ringed space). Let~$\I \subseteq \O_X$ be
a sheaf of ideals. The ringed space~$V(\I)$ is reduced if and only if, from the
internal perspective of~$\Sh(X)$, the ideal~$\I$ is a radical ideal.\end{lemma}
\begin{proof}The following chain of equivalences holds:
\begin{align*}
  &\ \Sh(V(\I)) \models \speak{$\O_{V(\I)}$ is a reduced ring} \\
  \Longleftrightarrow&\
    \Sh(V(\I)) \models \bigwedge_{n \geq 0} \forall s\?\O_{V(\I)}\_
      s^n = 0 \Longrightarrow s = 0 \\
  \Longleftrightarrow&\
    \Sh(X) \models \bigl(\bigwedge_{n \geq 0} \forall s\?\O_X/\I\_ s^n = 0
    \Rightarrow s = 0\bigr)^\Box \\
  \Longleftrightarrow&\
    \Sh(X) \models \bigwedge_{n \geq 0} \forall s\?\O_X/\I\_ s^n = 0 \Rightarrow \Box(s = 0) \\
  \Longleftrightarrow&\
    \Sh(X) \models \bigwedge_{n \geq 0} \forall s\?\O_X\_ s^n \in \I
    \Rightarrow \Box(s \in \I) \\
  \Longleftrightarrow&\
    \Sh(X) \models \bigwedge_{n \geq 0} \forall s\?\O_X\_ s^n \in \I
    \Rightarrow s \in \I \\
  \Longleftrightarrow&\
    \Sh(X) \models \speak{$\I$ is a radical ideal}
\end{align*}
In the second-to-last step, we used that~$\Box(s \in \I) \equiv ((s \in \I) \vee
(1 \in \I))$ implies~$s \in \I$. This is trivial in both cases of the
disjunction.
\end{proof}

\begin{lemma}\label{lemma:reduced-subspace}
Let~$X$ be a scheme (or a ringed space).
\begin{enumerate}
\item There exists a reduced closed sub-ringed space~$X_\mathrm{red}
\hookrightarrow X$ having the same underlying topological space as~$X$ with
the following universal property: Any morphism~$Y \to X$
of (ringed or locally ringed) spaces such that~$Y$ is reduced factors uniquely
over the closed immersion~$X_\mathrm{red} \hookrightarrow X$.
\item Let~$A \subseteq X$ be a closed subset. Then there exists a structure of
a reduced closed ringed subspace on~$A$ with a similar universal
property.
\end{enumerate}
\end{lemma}
\begin{proof}For the first statement, let~$\N \subseteq \O_X$ be the nilradical
of~$\O_X$. This can internally be simply defined by~$\N \defeq \sqrt{(0)} = \{
s\?\O_X \,|\, \bigvee_{n \geq 0} s^n = 0 \}$. Define~$X_\mathrm{red}$ as the closed
subspace associated to this sheaf of ideals. This ringed space is reduced by the
previous lemma. If~$X$ is a scheme, then quasicoherence of~$\N$ (which is
necessary and sufficient for~$X_\mathrm{red}$ to be a scheme) can be shown
internally (Example~\ref{ex:radical-qcoh}).
The proof of the universal property can also be done in the
internal language, by using the basic fact of locale theory that the
category of locales over~$X$ is equivalent to internal locales in~$\Sh(X)$; but
we do not want to discuss this further.

For the second statement, internally define the ideal~$\I \defeq \{ s\?\O_X
\,|\, \speak{$s$ \inv} \Rightarrow A^c \} \subseteq \O_X$.
Then~$1 \in \I$ if and only if~$A^c$, thus by
Proposition~\ref{prop:basics-closed-subspace} the closed ringed subspace defined
by~$\I$ has~$A$ as underlying topological space. It is reduced since~$\I$ is a
radical ideal.\end{proof}

\begin{rem}By Proposition~\ref{prop:quasicoherator-structure-sheaf}, the
ideal~$\I$ defined in the proof of Lemma~\ref{lemma:reduced-subspace} is
internally quasicoherent. Therefore the closed ringed subspace defined by~$\I$
is a scheme if~$X$ is.\end{rem}

\begin{lemma}Let~$X$ be a scheme of dimension~$\leq n$. Let~$V(\I)
\hookrightarrow X$ be a closed subscheme which is locally cut out by a regular
equation. Then~$\dim V(\I) \leq n-1$.\end{lemma}
\begin{proof}By Proposition~\ref{prop:dimension-scheme-ox}, it suffices to give
an intuitionistic proof of the following fact of dimension theory: Let~$A$ be
an arbitrary ring of dimension~$\leq n$. Let~$I = (s) \subseteq A$ be an ideal
which is generated by a regular element~$s\?A$. Then the~$\Box$-translation
of~``$A/I$ is of dimension~$\leq n-1$'' holds. In fact, we can show that~$A/I$
really is of dimension~$\leq n-1$; since no implication signs occur in a formal rendering of ``being of dimension~$\leq n-1$'',
Lemma~\ref{lemma:open-stalk} is applicable and implies
that this a stronger statement.

For this, let a sequence~$([a_0],\ldots,[a_{n-1}])$ of elements in~$A/I$ be
given. We can lift and extend this sequence to the
sequence~$(a_0,\ldots,a_{n-1},s)$ of elements of~$A$. Since~$\dim A \leq n$,
there exists a complementary sequence~$(b_0,\ldots,b_{n-1},b_n)$.
Since~$s$ is regular, the inclusion~$\sqrt{(s b_n)} \subseteq \sqrt{(0)}$
given by the definition of complementarity implies that~$b_n$ is nilpotent.
Thus we have that~$\sqrt{(a_{n-1}b_{n-1})} \subseteq \sqrt{(s,b_n)} =
\sqrt{(s)}$ in~$A$, which translates to~$\sqrt{([a_{n-1}] [b_{n-1}])} \subseteq
\sqrt{(0)}$ in~$A/I$.  Therefore~$([b_0],\ldots,[b_{n-1}])$ is a complementary
sequence to~$([a_0],\ldots,[a_{n-1}])$ in~$A/I$.
\end{proof}

\begin{lemma}\label{lemma:dim-closed-subscheme}
Let~$X$ be a scheme. Let~$\I$ be a sheaf of~$\O_X$-modules. Then:
\[ \dim V(\I) \leq n \quad\Longleftrightarrow\quad
  \Sh(X) \models \speak{$\O_X/\I$ is of Krull dimension~$\leq n$}. \]
\end{lemma}
\begin{proof}By Proposition~\ref{prop:dimension-scheme-ox}, the condition~$\dim
V(\I) \leq n$ is equivalent to
\[ \Sh(V(\I)) \models \speak{$\O_{V(\I)}$ is of Krull dimension~$\leq n$}. \]
By Theorem~\ref{thm:box-translation-semantically} this is equivalent to
\[ \Sh(X) \models (\speak{$\O_X/\I$ is of Krull dimension~$\leq n$})^\Box, \]
where~$\Box$ is the modal operator given by~$\Box\varphi \defeq (\varphi \vee
(1\in\I))$. The claimed equivalence then follows by
Lemma~\ref{lemma:open-stalk} (for~``$\Leftarrow$'') and by direct inspection
similar to the proof of Lemma~\ref{lemma:pushforward-finite-type}
(for~``$\Rightarrow$'').
\end{proof}

\section{Transfer principles}
\label{sect:transfer-principles}

Let~$M$ be an~$A$-module. A natural question is how properties of~$M$
relate to properties of the induced quasicoherent sheaf~$M^\sim$
on~$\Spec A$. For instance it is well-known that
\begin{itemize}
\item $M$ is finitely generated iff~$M^\sim$ is of finite type,
\item $M$ is flat over~$A$ iff~$M^\sim$ is flat over~$\O_{\Spec A}$, and
\item $M$ is torsion iff~$M^\sim$ is a torsion sheaf.
\end{itemize}
Using the internal language of the little Zariski topos of~$\Spec A$, we can
give a simple, conceptual, and uniform explanation of these equivalences.
Namely, from the internal point of view, the module~$M^\sim$ is obtained from
the constant sheaf~$\ul{M}$ by localizing at the \emph{generic filter}, a
particular multiplicative subset to be introduced below, and the set~$M$ and
the sheaf~$\ul{M}$ share the same properties (by
Lemma~\ref{lemma:properties-of-constant-sheaves} below).

This makes it obvious that, for instance, properties which are stable under
localization pass from~$M$ to~$M^\sim$.\footnote{More precisely, we should
write that properties for which there is an intuitionistic proof that they are
stable under localization pass from~$M$ to~$M^\sim$.}

\subsection{Internal properties of constant sheaves}

\begin{lemma}\label{lemma:properties-of-constant-sheaves}Let~$\varphi$ be a
formula in which arbitrary sets and elements may occur as parameters. Let~$X$
be a topological space and let~$U \subseteq X$ be an open subset. Then
\[ U \models \varphi \quad\text{iff}\quad (\text{$U$ inhabited} \Rightarrow
\varphi). \]
\end{lemma}
We are abusing notation on the left hand side: The parameters
of~$\varphi$, which are sets and elements, must be read as the induced constant
sheaves and constant functions (sections of those sheaves).
Unbounded quantifiers have to be read as ranging only over locally constant
sheaves, not all sheaves.
\begin{proof}By induction on the structure of~$\varphi$. By way of example, we
give the argument in the case that~$\varphi \equiv (a = b)$, where~$a$ and~$b$ are
elements of some set~$M$. Then~$U \models \varphi$ means by definition that the
constant functions~$U \to M$ with value~$a$ respectively~$b$ coincide. This is
equivalent to saying that~$a$ and~$b$ coincide if~$U$ is inhabited.
\end{proof}

The lemma in particular implies that constant sheaves enjoy several
classical properties from the internal point of view (if they are present in the metatheory), even though the internal
language only supports intuitionistic reasoning in general. For instance, for a
constant sheaf~$\ul{M}$ it holds that
\[ \Sh(X) \models \forall x,y\?\ul{M}\_ \neg\neg(x = y) \Rightarrow x = y \]
and even
\[ \Sh(X) \models \forall x,y\?\ul{M}\_ x = y \vee x \neq y. \]

\begin{rem}Lemma~\ref{lemma:properties-of-constant-sheaves} is also valid for
locales instead of topological spaces. If one works in an intuitionistic
metatheory, one has to add the additional requirement that the locale is
\emph{overt}; classically, every locale is overt, and intuitionistically,
at least locales arising from topological spaces are overt. We'll revisit this
subtle point in Section~\ref{sect:constructive-scheme-theory}, where we sketch
how scheme theory can be developed in an intuitionistic context, and more
specifically in Lemma~\ref{lemma:properties-of-constant-sheaves-over-locales}.
\end{rem}

\subsection{The generic filter}
\label{sect:generic-filter}

Let~$A$ be a ring.

\begin{defn}\label{defn:filter}
A \emph{filter} of~$A$ is a subset~$F \subseteq A$ such that
\begin{itemize}
\item $0 \not\in F$,
\item $x + y \in F \Longrightarrow (x \in F) \vee (y \in F)$, and
\item $1 \in F$,
\item $xy \in F \Longleftrightarrow (x \in F) \wedge (y \in F)$
\end{itemize}
for all~$x,y \? A$.
\end{defn}

In classical logic, the complement of a prime ideal is a filter and
furthermore every filter is of such a form. In constructive mathematics however,
it is useful to axiomatize complements of prime ideals directly, avoiding
negations. Intuitionistically, since De Morgan's law~$\neg(\alpha \wedge \beta)
\Rightarrow \neg\alpha \vee \neg\beta$ is not available, one can neither show
that the complement of a prime ideal is a filter nor that the complement of a
filter is a prime ideal.

A filter is in particular a multiplicative subset. Inverting the
elements of a filter results in a local ring, while intuitionistically
the localization of a ring at a prime ideal cannot in general be verified to be local.

\begin{defn}The \emph{generic filter}~$\F$ is the subsheaf of~$\ul{A}$
on~$\Spec A$ given by
\[ \Gamma(U, \F) \defeq \{ f : U \to A \,|\,
  \text{$f(\ppp) \not\in \ppp$ for all $\ppp \in U$} \}. \]
\end{defn}

\begin{prop}\label{prop:basics-univ-filter}\ \begin{enumerate}
\item Let~$f \in A$ and~$x \in A$. Then~$D(f) \models x \in \F$ if and only
if~$f \in \sqrt{(x)}$.
\item The stalk~$\F_\ppp$ at a point~$\ppp \in \Spec A$
is in canonical bijection with~$A \setminus \ppp$.
\item From the internal point of view of~$\Sh(\Spec A)$, the generic
filter is indeed a filter of~$\ul{A}$.
\end{enumerate}
\end{prop}
\begin{proof}By definition~$D(f) \models x \in \F$ means that~$x \not\in \ppp$
for all prime ideals~$\ppp$ with~$f \not\in \ppp$. This is well-known to be
equivalent to~$f \in \sqrt{(x)}$.

For the claim about stalks, we observe that the canonical map~$\F_\ppp \to A
\setminus \ppp$ sending a germ~$[f]$ to~$f(\ppp)$ is invertible with inverse
being the map which sends an element~$x \not\in \ppp$ to the germ of the
constant function with value~$x$ (defined on~$D(x)$).

Regarding the third statement we only verify the axiom regarding sums, the
other verifications being easier. Interpreting this axiom with the Kripke--Joyal
semantics and restricting without loss of generality to open subsets where
given locally constant functions are constant, let elements~$x,y \in A$ be
given such that~$D(f) \models x+y \in \F$. By the first statement~$f \in
\sqrt{(x+y)}$. Therefore~$D(f) \subseteq D(x) \cup D(y)$, and on~$D(x)$ it
holds that~$x \in \F$ and on~$D(y)$ it holds that~$y \in \F$.
\end{proof}

The significance of the generic filter is given by the following proposition.
\begin{prop}\label{prop:tilde-construction-internally}
From the internal point of view of~$\Sh(\Spec A)$,
\begin{enumerate}
\item the structure
sheaf~$\O_{\Spec A}$ is the localization of the constant sheaf~$\ul{A}$ at the
generic filter:~$\O_{\Spec A} = \ul{A}[\F^{-1}]$, and
\item the quasicoherent sheaf of modules~$M^\sim$ associated to
an~$A$-module~$M$ is the localization of the constant sheaf~$\ul{M}$ at the
generic filter.
\end{enumerate}\end{prop}
\begin{proof}Ignoring the ring respectively module structure, the second
statement is more general; therefore we prove this one. We didn't discuss the
case of quotients in Section~\ref{sect:internal-constructions}. However it
should be perspicuous that the interpretation of~$\ul{M}[\F^{-1}]$ is defined
as the colimit of~$\E \twoheadrightarrow \ul{M} \times \F$, taken in the
category of sheaves on~$\Spec A$, where~$\E$ is the subsheaf of~$\F \times
(\ul{M} \times \F) \times (\ul{M} \times \F)$ given by~$\E(U) \defeq \{
(s,(x,t),(y,u)) \,|\, sux = sty \}$.

This colimit can be obtained as the sheafification of the similarly defined
presheaf colimit~$\E' \twoheadrightarrow \ul{M}_\mathrm{pre} \times \F$,
where~$\ul{M}_\mathrm{pre}$ is the constant \emph{presheaf} associated to~$M$.
On an open subset~$U$ this presheaf colimit is simply the
localization~$\Gamma(U, \ul{M}_\mathrm{pre})[\Gamma(U, \F)^{-1}] = M[\Gamma(U,
\F)^{-1}]$. In the special case that~$U = D(f)$ is a standard open subset,
Proposition~\ref{prop:basics-univ-filter}(a) shows that this module is
canonically isomorphic to~$M[f^{-1}]$. The quasicoherent sheaf~$M^\sim$ of
modules admits the same description.
\end{proof}

Recognizing~$\O_{\Spec A}$ as a localization of~$\ul{A}$ fits nicely into the
following abstract algebraic motivation for schemes: Does the ring~$A$ admit a
\emph{universal localization}, \ie a homomorphism~$A \to A'$ into a local ring
such that every homomorphism~$A \to B$ into a local ring factors via a local
map over~$A \to A'$? Intuitively speaking, can we localize a ring at all prime
ideals at once, or equivalently at all filters at once? The answer is \emph{no}
in general,\footnote{Assume that the universal localization~$A'$ of a ring~$A$ exists as an ordinary ring in~$\Set$. Then any
two prime ideals~$\ppp$ and~$\qqq$ of~$A$ are equal: Let~$s \not\in \ppp$.
Since~$s$ is invertible in the local ring~$A_\ppp$ and the map~$A' \to A_\ppp$
induced by~$A \to A_\ppp$ is local, it is also invertible in~$A'$. Therefore the image of~$s$ in~$A_\qqq$ is
invertible as well. Thus~$s \not\in \qqq$.} but always \emph{yes} if we are
willing to change the topos in which we look for a solution: The universal
localization of~$A$ is given by the ring~$\O_{\Spec A}$ in the topos~$\Sh(\Spec
A)$; this ring is constructed by localizing~$\ul{A}$ at the generic
filter, a filter which exists in~$\Sh(\Spec A)$ but not in~$\Set$.

We expand on this point of view in Section~\ref{sect:relative-spectrum} on the
relative spectrum.

For transferring properties of~$M^\sim$ to~$M$, the following metatheorem is
crucial.
\begin{prop}\label{prop:metaproperty-of-the-generic-filter}
Let~$\I$ be an ideal in~$\ul{A}$ such that, for all inhabited
open subsets~$U \subseteq \Spec A$ and elements~$x \in A$, the set~$\Gamma(X,\I)$
contains the constant function with value~$x$ if~$\Gamma(U,\I)$ does. Then
\[ D(f) \models \speak{$\I \cap \F$ is inhabited}
  \quad\text{implies}\quad
  \text{for some~$n \geq 0$, $D(f) \models f^n \in \I$.} \]
\end{prop}

Lemma~\ref{lemma:properties-of-constant-sheaves} gives a simple and purely
syntactical criterion for the hypothesis on~$\I$: It suffices for~$\I$ to be
internally defined by an expression of the form~$\{ a\?\ul{A} \,|\, \varphi(a)
\}$, where~$\varphi$ is a formula which refers only to constant sheaves.

The metatheorem reflects the following well-known fact of classical
ring theory: If an ideal meets every filter (that is, the complement of every
prime ideal), it is the unit ideal. In this particular formulation the statement can't be
proven intuitionistically; the occurrence of ``every filter'' has to be
replaced by ``generic filter''. Intuitively, the generic filter is a
reification of the abstract idea of an ``arbitrary filter'', a filter about
which nothing is known except that it satisfies the filter axioms.

\begin{proof}
Let~$D(f) \models \speak{$\I \cap \F$ is inhabited}$. Then there
exists an open cover~$D(f) = \bigcup_i D(f_i)$ and elements~$x_i \in A$ such
that~$D(f_i) \models x_i \in \F$ and~$D(f_i) \models x_i \in \I$. By
Proposition~\ref{prop:basics-univ-filter} we have that~$f_i \in \sqrt{(x_i)}$
and therefore~$D(f_i) \models f_i^{m_i} \in \I$ for some~$m_i \geq 0$. We may
assume that all the~$D(f_i)$ are inhabited and that the exponents~$m_i$ are all
equal to some number~$m$. The assumption on~$\I$ implies~$D(f) \models f_i^m
\in \I$ for all~$i$. By a standard argument we can write~$f^n = \sum_i a_i
f_i^m$ for some coefficients~$a_i$; thus~$D(f) \models f^n \in \I$.
\end{proof}

\begin{rem}The stronger statement
\[ D(f) \models (\speak{$\I \cap \F$ is inhabited} \Rightarrow \bigvee_{n \geq
0} (f^n \in \I)) \]
does not hold in general. Indeed, consider the example~$f \defeq 1$ and~$\I \defeq
\brak{(g)} \defeq \brak{\{ a\?\ul{A} \,|\, \exists b\?\ul{A}\_ a = bg \}}$,
where~$g$ is a fixed element of~$A$ which is not nilpotent and not invertible.
Since~$D(g) \models g \in \I \cap \F$, the stronger statement would imply~$D(g)
\models 1 \in \I$. By Lemma~\ref{lemma:properties-of-constant-sheaves}, this is
equivalent to~$g$ being invertible in~$A$.
\end{rem}

\begin{rem}Recall from Proposition~\ref{prop:kx-internally} that the
sheaf~$\K_{\Spec A}$ of rational functions can internally be obtained by
localizing~$\O_{\Spec A}$ at the set of regular elements. Since~$\O_{\Spec A}$
is itself a localization, the sheaf~$\K_{\Spec A}$ is therefore obtained by a
two-step process. It can also be obtained in a single step by
localizing~$\ul{A}$ at~$\T$, where~$\T$ is the subsheaf of~$\ul{A}$ defined
by
\[ \Gamma(U,\T) = \{ f : U \to A \,|\, \text{$f(\ppp)$ is regular in~$A_\ppp$
for all~$\ppp \in U$} \}. \]
This subsheaf is characterized by the property that, for all~$f \in A$ and~$x
\in A$,~$D(f) \models x \in \T$ if and only if~$x$ is regular in~$A[f^{-1}]$.
\end{rem}

\subsection{Internal proofs of common lemmas}
\label{sect:common-lemmas-transfer-principles}

\begin{lemma}Let~$A$ be a ring. Then~$A$ is reduced if and only if the
scheme~$\Spec A$ is reduced.\end{lemma}
\begin{proof}By Proposition~\ref{prop:reduced-ring} the scheme~$\Spec A$ is
reduced if and only if~$\O_{\Spec A}$ is a reduced ring
from the internal point of view of~$\Sh(\Spec A)$.

For the ``only if'' direction assume that~$A$ is reduced. Then~$\ul{A}$ is
reduced as well, by Lemma~\ref{lemma:properties-of-constant-sheaves}. Since
localizations of reduced rings are reduced (and this fact has an intuitionistic
proof), in particular~$\O_{\Spec A} = \ul{A}[\F^{-1}]$ is reduced.

For the ``if'' direction let~$x \in A$ be an element such that~$x^n = 0$.
Since~$\O_{\Spec A} = \ul{A}[\F^{-1}]$ is reduced from the internal point of
view, the element~$x$ is zero in that ring, that is
\[ \Sh(\Spec A) \models \exists s\?\F\_ sx = 0. \]
Therefore the ideal internally defined by
\[ \I \defeq \{ a\?\ul{A} \,|\, ax = 0 \} \]
meets the generic filter. By
Proposition~\ref{prop:metaproperty-of-the-generic-filter} it follows
that~$\Sh(\Spec A) \models 1 \in \I$. By
Lemma~\ref{lemma:properties-of-constant-sheaves} this is equivalent to~$1 \cdot
x = 0$ as elements of~$A$.
\end{proof}

The ``if'' direction also admits a shorter proof, by simply
considering the Kripke--Joyal interpretation of~$\Sh(\Spec A) \models
\speak{$\O_{\Spec A}$ is reduced}$ and using the canonical isomorphism~$\Gamma(\Spec A, \O_{\Spec A})
\cong A$. We included the given proof to give a simple example of the mixed
internal/external reasoning with the generic filter. In a similar way we could
reprove Lemma~\ref{lemma:regular-affine}, the statement that a ring
element~$f \in A$ is regular in~$A$ if and only if, from the internal point of
view, it is regular in~$\O_{\Spec A}$.

\begin{lemma}\label{lemma:finite-type-using-universal-filter}
Let~$M$ be an~$A$-module. Then~$M^\sim$ is of finite type if and
only if~$M$ is finitely generated.\end{lemma}
\begin{proof}First assume that~$M$ is finitely generated over~$A$.
Then~$\ul{M}$ is finitely generated over~$\ul{A}$, by
Lemma~\ref{lemma:properties-of-constant-sheaves}. Since localizations of
finitely generated modules are finitely generated (over the localized ring),
the module~$M^\sim = \ul{M}[\F^{-1}]$ is finitely generated from the internal
point of view. By Proposition~\ref{prop:finite-type-and-co} this means
that~$M^\sim$ is of finite type from the external point of view.

For the ``only if`` direction, we assume that~$M^\sim$ is finitely generated
over~$\O_{\Spec A}$ from the internal point of view and have to verify that~$M$
is finitely generated over~$A$. So it holds that
\[ \Sh(\Spec A) \models \bigvee_{n \geq 0}
  \exists x_1,\ldots,x_n\?\ul{M}[\F^{-1}]\_
  \speak{the $x_i$ span~$\ul{M}[\F^{-1}]$ over~$\ul{A}[\F^{-1}]$}. \]
Since multiplying a generating family by an unit results again in a generating
family, we have in fact that
\[ \Sh(\Spec A) \models \bigvee_{n \geq 0}
  \exists x_1,\ldots,x_n\?\ul{M}\_
  \speak{the $x_i/1$ span~$\ul{M}[\F^{-1}]$ over~$\ul{A}[\F^{-1}]$} \]
or equivalently
\[ \Sh(\Spec A) \models \bigvee_{n \geq 0, x_1,\ldots,x_n \in M}
  \speak{the $x_i/1$ span~$\ul{M}[\F^{-1}]$ over~$\ul{A}[\F^{-1}]$}. \]
Since this is a directed disjunction and~$\Spec A$ is quasicompact,
Proposition~\ref{prop:quasicompact-meta} is applicable and shows that there
exists a natural number~$n \geq 0$ and elements~$x_1,\ldots,x_n \in M$ such
that
\[ \Sh(\Spec A) \models \speak{the $x_i/1$ span~$\ul{M}[\F^{-1}]$
over~$\ul{A}[\F^{-1}]$}. \]
We claim that these~$x_i$ also span~$M$ as an~$A$-module. So let~$x \in M$ be
arbitrary. By elementary linear algebra we can deduce that
\[ \Sh(\Spec A) \models \exists s\in\F\_ \exists a_1,\ldots,a_n\?\ul{A}\_
  sx = \sum_i a_i x_i. \]
Therefore the ideal internally defined by
\[ \I \defeq \{ s\?\ul{A} \,|\, \exists a_1,\ldots,a_n\?\ul{A}\_
  sx = \textstyle\sum_i a_i x_i \} \]
meets the generic filter.
Proposition~\ref{prop:metaproperty-of-the-generic-filter} shows that~$\Sh(\Spec A)
\models 1 \in \I$, that is~$x$ is an element of the~$A$-span of the~$x_i$.
\end{proof}

\begin{rem}If~$M^\sim$ can be generated by~$\leq n$ elements over~$\O_{\Spec
A}$ from the internal point of view, it needn't be the case that~$M$ can be
generated by~$\leq n$ elements over~$A$. It is instructive to see where the
appropriately modified version of the above proof fails: In this case we still
have
\[ \Sh(\Spec A) \models \bigvee_{x_1,\ldots,x_n \in M}
  \speak{the $x_i/1$ span~$\ul{M}[\F^{-1}]$ over~$\ul{A}[\F^{-1}]$}, \]
but this disjunction is no longer directed.
\end{rem}

\begin{lemma}Let~$X$ be a scheme. Then kernels and cokernels of morphisms
between quasicoherent~$\O_X$-modules are quasicoherent.\end{lemma}

\begin{proof}We may assume that~$X = \Spec A$ is affine. A morphism between
quasicoherent~$\O_X$-modules is of the form~$\ul{\varphi}[\F^{-1}] :
\ul{M}[\F^{-1}] \to \ul{N}[\F^{-1}]$, where~$\varphi : M \to N$ is a linear map
between~$A$-modules. Since taking constant sheaves and localization are exact,
we have the chain of isomorphisms
\[ (\ul{\Kernel(\varphi)})[\F^{-1}] \cong
  (\Kernel(\ul{\varphi}))[\F^{-1}] \cong
  \Kernel(\ul{\varphi}[\F^{-1}]), \]
and similarly for the cokernel.
\end{proof}

\subsection{An application to constructive mathematics}
\label{sect:eliminating-prime-ideals}

The generic filter has a practical application in constructive mathematics.
Recall that intuitionistically prime and maximal ideals don't work very well,
since one often needs the axiom of choice or related set-theoretical principles
in dealing with them. This is unfortunate, since calculations with prime and maximal ideals are
often quite useful. For example:
\begin{itemize}
\item To verify that a ring element is nilpotent, it suffices to verify that it
is an element of every prime ideal. For instance, this is calculationally simpler
when proving that the coefficients of a nilpotent polynomial are
themselves nilpotent.
\item To verify that there is a relation of the form~$1 = p_1f_1 + \cdots +
p_mf_m$ among polynomials~$f_1,\ldots,f_m \in K[X_1,\ldots,X_n]$ where~$K$ is
an algebraically closed field, it suffices to show that the~$f_i$ don't have a
common zero.
\end{itemize}

One could of course blithely switch to classical logic in this case. However this
might not be desirable, as a constructive proof would contain more information:
For instance, if we have classically proven that an element~$x$ is an element
of every prime ideal, then we know that some power~$x^n$ is zero. But from such
a proof we can't directly read off any upper bound on~$n$. From a constructive
proof of nilpotency, we could.

There is a way to combine some of the powerful tools of classical ring theory
with the advantages that constructive reasoning provides. Namely we can devise
a language in which we can usefully talk about prime ideals, but which
substitutes all non-constructive arguments by constructive arguments ``behind
the scenes''. The key idea is to substitute the phrase ``for all prime ideals''
(or equivalently ``for all filters'') by ``for the generic filter''.

This was already explored by Coquand, Coste, Lombardi, Roy, and
others under the theme of \emph{dynamical methods in
algebra}~\cite{clr:dynamicalmethod,cl:logical}. Here we show how one can use
the generic filter, as reified by a sheaf in the little Zariski topos, to
achieve similar effects.

\begin{prop}Let~$M$ and~$N$ be~$A$-modules. Let~$\alpha : M \to N$ be a linear
map. The interpretations of the statements in the second column of
Table~\ref{table:generic-filter-statements} in the internal language
of~$\Sh(\Spec A)$ are intuitionistically equivalent to the statements given in
the third column.\end{prop}
\begin{proof}To demonstrate the technique we verify the first and the last claim.
To make the following proofs constructive we have to define~$\Spec A$,
its sheaf topos, and the generic filter in a constructive fashion, not using prime ideals. This can
be done, by constructing~$\Spec A$ as a locale instead of a topological space.
We expand on this in Section~\ref{sect:spectrum-as-a-locale} and in
Section~\ref{sect:constructive-scheme-theory}.

The interpretation of~$\Sh(\Spec A) \models x \not\in \F$ by the Kripke--Joyal
semantics is that~$D(f) \models x \in \F$ implies~$D(f) = \emptyset$ for
all~$f \in A$. By Proposition~\ref{prop:basics-univ-filter}(a) this is
equivalent to
\[ \forall f \in A\_ f \in \sqrt{(x)} \Rightarrow f \in \sqrt{(0)}, \]
that is the statement that~$x$ is nilpotent in~$A$.

Assume that~$\alpha : M \to N$ is surjective. By
Lemma~\ref{lemma:properties-of-constant-sheaves} the induced map~$\ul{M}
\to \ul{N}$ is surjective from the internal point of view. Since localization
preserves surjectivity, also the map~$\ul{M}[\F^{-1}] \to \ul{N}[\F^{-1}]$ is
surjective.

Conversely, assume that~$\ul{M}[\F^{-1}] \to \ul{N}[\F^{-1}]$ is surjective
from the internal point of view. To verify that~$\alpha : M \to N$ is
surjective, let~$y \in N$. The assumption implies that the ideal internally
defined by
\[ \I \defeq \{ s\?\ul{A} \,|\, \exists x\?\ul{A}\_ sy = \ul{\alpha}(x) \} \]
meets the generic filter. By
Proposition~\ref{prop:metaproperty-of-the-generic-filter} this implies
that~$\Sh(\Spec A) \models 1 \in \I$, that is there exists an element~$x \in A$
such that~$\alpha(x) = y$.
\end{proof}

\begin{table}
  \centering
  \renewcommand{\arraystretch}{1.3}
  \small
  \begin{tabular}{lll}
    \toprule
    Statement & constructive substitution & meaning \\\midrule
    $x \in \ppp$ for all~$\ppp$. &
    $x \not\in \F$. &
    $x$ is nilpotent. \\
    $x \in \ppp$ for all~$\ppp$ such that~$y \in \ppp$. &
    $x \in \F \Rightarrow y \in \F$. &
    $x \in \sqrt{(y)}$. \\
    $x$ is regular in all stalks~$A_\ppp$. &
    $x$ is regular in~$\ul{A}[\F^{-1}]$. &
    $x$ is regular in~$A$. \\
    The stalks~$A_\ppp$ are reduced. &
    $\ul{A}[\F^{-1}]$ is reduced. &
    $A$ is reduced. \\
    The stalks~$M_\ppp$ vanish. &
    $\ul{M}[\F^{-1}] = 0$. &
    $M = 0$. \\
    The stalks~$M_\ppp$ are fin.\@ gen.\@ over~$A_\ppp$. &
    $\ul{M}[\F^{-1}]$ is fin.\@ gen.\@ over
    $\ul{A}[\F^{-1}]$. &
    $M$ is fin.\@ gen.\@ over~$A$. \\
    The stalks~$M_\ppp$ are flat over~$A_\ppp$. &
    $\ul{M}[\F^{-1}]$ is flat over~$\ul{A}[\F^{-1}]$. &
    $M$ is flat over~$A$. \\
    The maps~$M_\ppp \to N_\ppp$ are injective. &
    $\ul{M}[\F^{-1}] \to \ul{N}[\F^{-1}]$ is injective. &
    $M \to N$ is injective. \\
    The maps~$M_\ppp \to N_\ppp$ are surjective. &
    $\ul{M}[\F^{-1}] \to \ul{N}[\F^{-1}]$ is surjective. &
    $M \to N$ is surjective. \\
    \bottomrule
  \end{tabular}
  \vspace{0.5em}

  \caption{\label{table:generic-filter-statements}Substituting the use of prime
  ideals by the generic filter.}
\end{table}

\begin{rem}As is apparent from Table~\ref{table:generic-filter-statements},
there is a slight mismatch between the external ``for any prime ideal'' and
the internal ``for the generic filter''. It's not true that a module is
finitely generated if and only if all its stalks are finitely generated (a
counterexample is the~$\ZZ$-module~$\bigoplus_p \ZZ/(p)$). But it is true that
an~$A$-module~$M$ is finitely generated if and only if, internally
to~$\Sh(\Spec(A))$, the generic stalk~$\ul{M}[\F^{-1}]$ is finitely generated.

Intuitively, verifying a statement about the generic stalk doesn't only mean
that it holds for all (ordinary) stalks; it means that it holds for the
ordinary stalks in a \emph{uniform manner}. This extra bit of rigidity is what
allows to draw slightly stronger conclusions.

The other entries in Table~\ref{table:generic-filter-statements} don't show
this slight difference in semantics.
\end{rem}

The sheaf-theoretical approach using the generic filter is different from the
dynamical methods in the following aspect. We have to reword classical
arguments using (the generic) filter instead of (the generic) prime ideal.
Depending on the situation this might be a nuisance. One might be tempted to
employ the complement of the generic filter, but this is only an ideal, not a
prime ideal from the internal point of view.\footnote{One can check that the
complement of~$\F$ in~$\ul{A}$ is the subsheaf~$\P$ defined by~$\Gamma(U, \P)
\defeq \{ f : U \to A \,|\, \text{$f(\ppp) \in \ppp$ for all~$\ppp \in U$} \}$
and that~$D(f) \models x \in \P$ if and only if~$fx$ is nilpotent. This can be
used to show that the statement~$\Sh(\Spec A) \models \forall x,y\?\ul{A}\_ xy
\in \P \Rightarrow x \in \P \vee y \in \P$ is false in general. A concrete
counterexample is given by~$A = \ZZ[U,V]/(UV)$. Then~$\Sh(\Spec A) \models [U]
\cdot [V] \in \P$, but $\Sh(\Spec A) \not\models [U] \in \P \vee [V] \in \P$.}

\subsection{An internal proof of Grothendieck's generic freeness lemma}\label{sect:generic-freeness}
The goal of this subsection is to give a simple proof of Grothendieck's generic
freeness lemma in the following general form.

\begin{thm}\label{thm:generic-freeness}
Let~$A$ be a reduced ring. Let~$B$ be an~$A$-algebra of finite type.
Let~$M$ be a finitely generated~$B$-module. Then there is a dense open
subset~$U \subseteq \Spec(A)$ such that over~$U$,
\begin{enumerate}
\item $B^\sim$ is finitely presented as an~$\O_{\Spec(A)}$-algebra,
\item $M^\sim$ is of finite presentation over~$B^\sim$, and
\item $M^\sim$ is (not necessarily finite) locally free as an~$\O_{\Spec(A)}$-module.
\end{enumerate}
\end{thm}

The usual proofs of Grothendieck's generic freeness lemma proceed using a series
of reduction steps which are arguably not very memorable or
straightforward, see for instance~\stacksproject{051Q}
or~\cite{staats:generic-freeness}. In particular, there doesn't seem to be a
published proof which tackles the Noetherian and non-Noetherian cases in one go.
Employing the internal language, Grothendieck's generic freeness lemma can be
proved in a simple, conceptual, and constructive way without any reduction steps.

This section was prompted by a MathOverflow thread~\cite{mo:kernel} and
greatly benefited from discussions with Brandenburg.

\begin{figure}[h]
  \centering
  \begin{tikzpicture}[scale=1.4]
    \draw[step=1cm,gray,very thin] (0,0) grid (8,8);

    \draw
      (0,8) -- (0,0) -- (8,0);

    \fill[pattern=north east lines,pattern color=gray,opacity=0.5]
      (2,8) -- (2,5) -- (5,5) -- (5,3) -- (6,3) -- (6,2) -- (8,2) -- (8,8);
    \draw
      (2,8) -- (2,5) -- (5,5) -- (5,3) -- (6,3) -- (6,2) -- (8,2);

    \fill[fill=blue!40!white,opacity=0.5]
      (0,8) -- (0,0) -- (4,0) -- (4,3) -- (3,3) -- (3,5) -- (2,5) -- (2,8);

    \fill[fill=red!40!white,opacity=0.5]
      (3,3) -- (4,3) -- (4,4) -- (3,4);

    \begin{scope}[yshift=4cm, xshift=4cm]
      \matrix[matrix of nodes,nodes={inner sep=0pt,text width=1.4cm,align=center,minimum height=1.4cm}]{
        $x^0y^7v_1$ & $x^1y^7v_1$ & $x^2y^7v_1$ & $x^3y^7v_1$ & $x^4y^7v_1$ & $x^5y^7v_1$ & $x^6y^7v_1$ & $x^7y^7v_1$ & \\
        $x^0y^6v_1$ & $x^1y^6v_1$ & $x^2y^6v_1$ & $x^3y^6v_1$ & $x^4y^6v_1$ & $x^5y^6v_1$ & $x^6y^6v_1$ & $x^7y^6v_1$ & \\
        $x^0y^5v_1$ & $x^1y^5v_1$ & $x^2y^5v_1$ & $x^3y^5v_1$ & $x^4y^5v_1$ & $x^5y^5v_1$ & $x^6y^5v_1$ & $x^7y^5v_1$ & \\
        $x^0y^4v_1$ & $x^1y^4v_1$ & $x^2y^4v_1$ & $x^3y^4v_1$ & $x^4y^4v_1$ & $x^5y^4v_1$ & $x^6y^4v_1$ & $x^7y^4v_1$ & \\
        $x^0y^3v_1$ & $x^1y^3v_1$ & $x^2y^3v_1$ & $x^3y^3v_1$ & $x^4y^3v_1$ & $x^5y^3v_1$ & $x^6y^3v_1$ & $x^7y^3v_1$ & \\
        $x^0y^2v_1$ & $x^1y^2v_1$ & $x^2y^2v_1$ & $x^3y^2v_1$ & $x^4y^2v_1$ & $x^5y^2v_1$ & $x^6y^2v_1$ & $x^7y^2v_1$ & \\
        $x^0y^1v_1$ & $x^1y^1v_1$ & $x^2y^1v_1$ & $x^3y^1v_1$ & $x^4y^1v_1$ & $x^5y^1v_1$ & $x^6y^1v_1$ & $x^7y^1v_1$ & \\
        $x^0y^0v_1$ & $x^1y^0v_1$ & $x^2y^0v_1$ & $x^3y^0v_1$ & $x^4y^0v_1$ & $x^5y^0v_1$ & $x^6y^0v_1$ & $x^7y^0v_1$ & \\
      };
    \end{scope}
  \end{tikzpicture}
  \caption{\label{fig:single-step}A single step in the iterative process used
  in the proof of Theorem~\ref{thm:generic-freeness}, in the special case~$n = 2, m = 1$. The hatched cells
  indicate vectors which have already been removed from the generating family.
  The vector in the red cell was found to be expressible as a linear
  combination of vectors with smaller index (blue cells). It is therefore
  about to be removed, along with the vectors in all cells to the top and to
  the right of the red cell.}
\end{figure}

\begin{figure}[h]
  \centering
  \begin{tikzpicture}[scale=0.4]
    \begin{scope}
      \draw[step=1cm,gray,very thin] (0,0) grid (8,8);

      \draw
        (0,8) -- (0,0) -- (8,0);

      \fill[fill=blue!40!white,opacity=0.5]
        (0,8) -- (0,0) -- (6,0) -- (6,5) -- (5,5) -- (5,8);

      \fill[fill=red!40!white,opacity=0.5]
        (5,5) -- (6,5) -- (6,6) -- (5,6);

      \node[shape=circle,draw,inner sep=2pt] at (1,1) {1};
    \end{scope}

    \begin{scope}[xshift=9cm]
      \draw[step=1cm,gray,very thin] (0,0) grid (8,8);

      \draw
        (0,8) -- (0,0) -- (8,0);

      \fill[pattern=north east lines,pattern color=gray,opacity=0.5]
        (5,8) -- (5,5) -- (8,5) -- (8,8);
      \draw
        (5,8) -- (5,5) -- (8,5);

      \fill[fill=blue!40!white,opacity=0.5]
        (0,8) -- (0,0) -- (3,0) -- (3,5) -- (2,5) -- (2,8);

      \fill[fill=red!40!white,opacity=0.5]
        (2,5) -- (3,5) -- (3,6) -- (2,6);

      \node[shape=circle,draw,inner sep=2pt] at (1,1) {2};
    \end{scope}

    \begin{scope}[xshift=18cm]
      \draw[step=1cm,gray,very thin] (0,0) grid (8,8);

      \draw
        (0,8) -- (0,0) -- (8,0);

      \fill[pattern=north east lines,pattern color=gray,opacity=0.5]
        (2,8) -- (2,5) -- (8,5) -- (8,8);
      \draw
        (2,8) -- (2,5) -- (8,5);

      \fill[fill=blue!40!white,opacity=0.5]
        (0,8) -- (0,0) -- (7,0) -- (7,2) -- (6,2) -- (6,5) -- (2,5) -- (2,8);

      \fill[fill=red!40!white,opacity=0.5]
        (6,2) -- (7,2) -- (7,3) -- (6,3);

      \node[shape=circle,draw,inner sep=2pt] at (1,1) {3};
    \end{scope}

    \begin{scope}[yshift=-9cm]
      \draw[step=1cm,gray,very thin] (0,0) grid (8,8);

      \draw
        (0,8) -- (0,0) -- (8,0);

      \fill[pattern=north east lines,pattern color=gray,opacity=0.5]
        (2,8) -- (2,5) -- (6,5) -- (6,2) -- (8,2) -- (8,8);
      \draw
        (2,8) -- (2,5) -- (6,5) -- (6,2) -- (8,2);

      \fill[fill=blue!40!white,opacity=0.5]
        (0,8) -- (0,0) -- (6,0) -- (6,3) -- (5,3) -- (5,5) -- (2,5) -- (2,8);

      \fill[fill=red!40!white,opacity=0.5]
        (5,3) -- (6,3) -- (6,4) -- (5,4);

      \node[shape=circle,draw,inner sep=2pt] at (1,1) {4};
    \end{scope}

    \begin{scope}[yshift=-9cm,xshift=9cm]
      \draw[step=1cm,gray,very thin] (0,0) grid (8,8);

      \draw
        (2,8) -- (2,5) -- (5,5) -- (5,3) -- (6,3) -- (6,2) -- (8,2);

      \draw
        (0,8) -- (0,0) -- (8,0);

      \fill[pattern=north east lines,pattern color=gray,opacity=0.5]
        (2,8) -- (2,5) -- (5,5) -- (5,3) -- (6,3) -- (6,2) -- (8,2) -- (8,8);

      \fill[fill=blue!40!white,opacity=0.5]
        (0,8) -- (0,0) -- (4,0) -- (4,3) -- (3,3) -- (3,5) -- (2,5) -- (2,8);

      \fill[fill=red!40!white,opacity=0.5]
        (3,3) -- (4,3) -- (4,4) -- (3,4);

      \node[shape=circle,draw,inner sep=2pt] at (1,1) {5};
    \end{scope}

    \begin{scope}[yshift=-9cm,xshift=18cm]
      \draw[step=1cm,gray,very thin] (0,0) grid (8,8);

      \draw
        (0,8) -- (0,0) -- (8,0);

      \fill[pattern=north east lines,pattern color=gray,opacity=0.5]
        (2,8) -- (2,5) -- (3,5) -- (3,3) -- (6,3) -- (6,2) -- (8,2) -- (8,8);
      \draw
        (2,8) -- (2,5) -- (3,5) -- (3,3) -- (6,3) -- (6,2) -- (8,2);

      \node[shape=circle,draw,inner sep=2pt] at (1,1) {6};
    \end{scope}
  \end{tikzpicture}
  \caption{\label{fig:iteration}The iterative process used
  in the proof of Theorem~\ref{thm:generic-freeness}, in the special case~$n = 2, m = 1$.
  The process terminates after reducing the generating family
  a finite number of times.}
\end{figure}
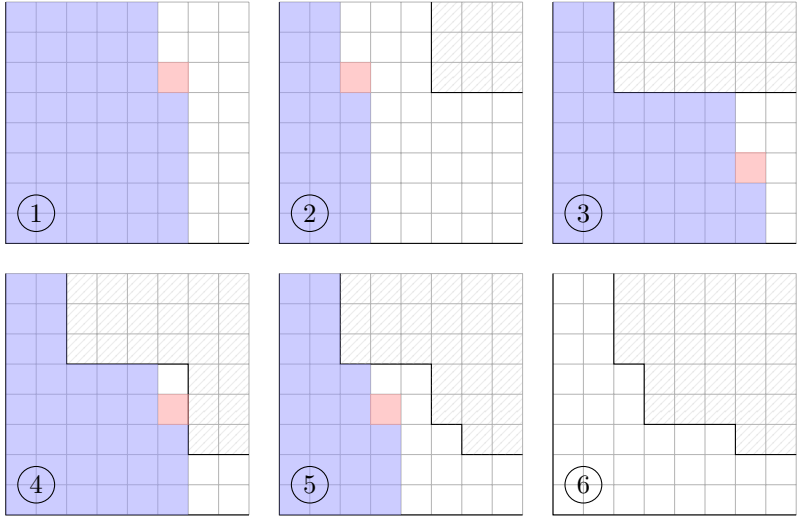

\begin{proof}[Proof of Theorem~\ref{thm:generic-freeness}]
Since ``dense open'' translates to ``not not'' in the internal language
(Proposition~\ref{prop:modops-kripke}), it suffices to prove that, from the
internal point of view of~$\Sh(\Spec(A))$, it's \notnot the case that
\begin{enumerate}
\item $B^\sim$ is of finite presentation over~$\O_{\Spec(A)}$,
\item $M^\sim$ is finitely presented as a~$B^\sim$-module, and
\item $M^\sim$ is (not necessarily finite) free over~$\O_{\Spec(A)}$.
\end{enumerate}

Since~$B^\sim$ is finitely generated as an~$\O_{\Spec(A)}$-algebra, it is
isomorphic to an algebra of the form~$\O_{\Spec(A)}[X_1,\ldots,X_n]/\aaa$ for
some number~$n \geq 0$ and some ideal~$\aaa$. By
Proposition~\ref{prop:ox-anonymously-noetherian} and Theorem~\ref{thm:hilbert}, the
ring~$\O_{\Spec(A)}[X_1,\ldots,X_n]$ is anonymously Noetherian. Therefore~$\aaa$ is
\notnot finitely generated, showing that~$B^\sim$ is \notnot of finite
presentation over~$\O_{\Spec(A)}$.

Similarly, the module~$M^\sim$ is of the form~$(B^\sim)^m/U$ for some number~$m
\geq 0$ and some submodule~$U$. Since~$(B^\sim)^m$ is anonymously Noetherian as a
direct sum of anonymously Noetherian modules, the submodule~$U$ is \notnot finitely
generated. Thus~$M^\sim$ is \notnot a finitely presented~$B^\sim$-module.

The basic idea to show that~$M^\sim$ is \notnot free over~$\O_{\Spec(A)}$ is as
follows. Since~$\O_{\Spec(A)}$ is a field in the sense that noninvertible
elements are zero, minimal generating families are already linearly independent;
we observed this in the proof of Lemma~\ref{lemma:rank-functor-locally-constant}.
By the finiteness hypotheses, the module~$M^\sim$ admits a countable
generating family. It's \notnot the case that either one of these vectors can be
expressed as a linear combination of the others, or not. In the second case
we're done; in the first case, we remove the redundant vector and continue in
the same fashion.

However, if we shrink the given generating family in this naive fashion, the
process may not terminate in finitely many steps. In a classical context,
Zorn's lemma could be used to iterate the process transfinitely and eventually
obtain a minimal generating family, but Zorn's lemma is not available in the
internal universe of the little Zariski topos. We therefore have to pick the
vectors we'll remove in a more systematic fashion.

Let~$(x_1,\ldots,x_n)$ be a generating family for~$B^\sim$ as
an~$\O_{\Spec(A)}$-algebra and let $(v_1,\ldots,v_m)$ be a generating
family for~$M^\sim$ as a~$B^\sim$-module. We endow the set
\[ I \defeq \{ (j, i_1,\ldots,i_n) \,|\,
  j \in \{ 1,\ldots,m \},
  i_1,\ldots,i_n \in \{ 0,1,\ldots \} \} \]
with the lexicographic order. We choose the family~$(x_1^{i_1} \cdots x_n^{i_n}
v_j)_{j,i_1,\ldots,i_n}$ as the starting point of the shrinking process. In
each step, we use that it's \notnot the case that
\begin{itemize}
\item either one of the vectors of the generating family can be expressed as a
linear combination of vectors in the family with a smaller index,
\item or not.
\end{itemize}

In the second case, the generating family is linearly independent: For any
linear combination summing to zero, we can show that all coefficients are zero,
beginning with the coefficient which is paired with the vector of greatest
index.

Figure~\ref{fig:single-step} illustrates our action in the first case. We remove the redundant vector~$x_1^{i_1} \cdots x_n^{i_n}
v_j$ \emph{and also any vector with greater powers of the~$x_1,\ldots,x_n$}
from the generating family. The resulting family will still be a generating
family, since the linear combination witnessing the redundancy of~$x_1^{i_1}
\cdots x_n^{i_n} v_j$ successively gives rise to linear combinations witnessing
the redundancy of the vectors~$x_1^{\geq i_1} \cdots x_n^{\geq i_n} v_j$;
we maintain the invariant that any member of the starting generating family can
be expressed as a linear combination of vectors of the current generating
family with smaller or equal index.

As indicated in Figure~\ref{fig:iteration}, this process terminates after
finitely many steps. This fact is related to the fact that the ordinal~$\omega^n$ is
well-founded; the formal statement ensuring termination is known as
\emph{Dickson's Lemma} (see, for instance,~\cite[Thm.~1.1]{veldman:kruskal}).
\end{proof}

Since the given internal proof was (necessarily) intuitionistically valid, the
internal language machinery is intuitionistically valid, and the construction
of the spectrum can be set up in an intuitionistically sensible way
(Section~\ref{sect:relative-spectrum}), an intuitionistic external proof not
employing the topos machinery can be extracted from the given argument.
The resulting proof will verify Grothendieck's generic freeness lemma in the
following form.

\begin{thm}\label{thm:generic-freeness-constructively}
Let~$A$ be a reduced ring. Let~$B$ be an~$A$-algebra of finite type.
Let~$M$ be a finitely generated~$B$-module. Assume that the only element~$f \in
A$ such that
\begin{enumerate}
\item $B[f^{-1}]$ is of finite presentation over~$A[f^{-1}]$,
\item $M[f^{-1}]$ is finitely presented as a~$B[f^{-1}]$-module, and
\item $M[f^{-1}]$ is free over~$A[f^{-1}]$
\end{enumerate}
is~$f = 0$. Then~$A = 0$.
\end{thm}

In classical logic, this form implies Grothendieck's generic freeness lemma in its more abstract
formulation by a routine argument: Let~$U \subseteq \Spec(A)$ be the union over
all standard open subsets~$D(f)$ such that the statements~(1),~(2), and~(3) in
Theorem~\ref{thm:generic-freeness-constructively} hold. The
statements~(1),~(2), and~(3) of Theorem~\ref{thm:generic-freeness} hold on this
open subset, therefore it remains to show that~$U$ is dense.

So let a nonempty open subset~$V$ of~$\Spec(A)$ be given. This
contains a standard open subset~$D(g) \subseteq V$ such that~$g$ is not
nilpotent. Therefore the localized ring~$A[g^{-1}]$ is not zero. Thus the
conclusion of Theorem~\ref{thm:generic-freeness-constructively} is not satisfied.
Since we assume classical logic, there is a nonzero element~$f \in A[g^{-1}]$
such that statements~(1),~(2), and~(3) in
Theorem~\ref{thm:generic-freeness-constructively} hold for~$A[g^{-1}][f^{-1}]$,
$B[g^{-1}][f^{-1}]$, and~$M[g^{-1}][f^{-1}]$.
Writing~$f = h/g^n$, we see that~$U \cap V$ contains the nonzero open
subset~$D(gh)$.

We refrain from giving the resulting explicit proof of
Theorem~\ref{thm:generic-freeness-constructively} here, but will report on it
in the future~\cite{blechschmidt:generic-freeness}. A part of the proof was
included by Brandenburg in a paper of his~\cite{brandenburg:schur}.

\begin{rem}There is no hope that there is an intuitionistic proof of
Gro\-then\-dieck's generic freeness lemma in the form of
Theorem~\ref{thm:generic-freeness} even if the spectrum is constructed in an
intuitionistically sensible way, since there is the following Brouwerian
counterexample. Let~$\varphi$ be an arbitrary statement. Then
the~$\ZZ$-module~$M \defeq \ZZ/\aaa$, where~$\aaa \defeq \{ x \in \ZZ \,|\, (x
= 0) \vee \varphi \}$ as in Footnote~\ref{fn:z-principal-ideal-domain} on
page~\pageref{fn:z-principal-ideal-domain}, is finitely generated. By
assumption, there exists a nonzero element~$f \in \ZZ$ such that~$M[f^{-1}]$ is
a finite free module over~$A[f^{-1}]$ of some rank~$n$. If~$n = 0$, then~$f^m
\in \aaa$ for some~$m \geq 0$, therefore~$\varphi$ holds. If~$n \geq 1$,
then~$\neg\varphi$ holds, since~$\varphi$ would imply~$\aaa = \ZZ$ and
therefore~$M[f^{-1}] = 0$. Since~$n = 0 \vee n \geq 1$, it follows
that~$\varphi \vee \neg\varphi$.
\end{rem}

\subsection{\texorpdfstring{A note on~$\QQ$-algebras which are finitely
generated over~$\ZZ$}{A note on~ℚ-algebras which are finitely generated
over~ℤ}}

In this section, we want to show how the internal language of~$\Sh(\Spec(\ZZ))$
can be used to give a proof of the following fact.

\begin{prop}\label{prop:fingen-algebra-q}
Let~$A$ be a finitely presented~$\ZZ$-algebra such that any positive natural
number is invertible in~$A$. Then~$1 = 0$ in~$A$.
\end{prop}

A slick classical proof runs as follows: Assume that~$1 \neq 0$ in~$A$. Then
there exists a maximal ideal~$\mmm \subseteq A$. The preimage of~$\mmm$
in~$\ZZ$ is maximal since~$\ZZ$ is a Jacobson ring~\stacksproject{00GB} and
therefore of the form~$(p)$ for a prime number~$p$. Thus~$p \in \mmm$.
Since~$p$ is also invertible in~$A$, it follows that~$1 = 0$ in~$A$.

We intend the following proof as an example of how one can extend, in some
cases, the applicability of theorems about fields to rings using the internal
language. If one sets up the spectrum in an intuitionistically sensible way, as
described in Section~\ref{sect:relative-spectrum}, the proof avoids the axiom
of choice.

\begin{proof}[Proof of Proposition~\ref{prop:fingen-algebra-q}]
Noether's normalization lemma is intuitionistically valid in the following
form: ``Let~$K$ be a ring such that~$1 \neq 0$ and such that~$\neg(\speak{$s$
\inv}) \Rightarrow s = 0$ for all~$s \? K$. Let~$\aaa \subseteq
K[X_1,\ldots,X_n]$ be an ideal. Then it's \notnot the case that either~$\aaa =
(1)$ or that there exists a number~$r \geq 0$ and a finite injective
homomorphism~$K[Y_1,\ldots,Y_r] \to K[X_1,\ldots,X_n]/\aaa$ of~$K$-algebras.''

In the internal universe of~$\Sh(\Spec(\ZZ))$, the structure
sheaf~$\O_{\Spec(\ZZ)}$ satisfies the assumption on~$K$ by
Corollary~\ref{cor:field-reduced}. We can therefore apply the Noether
normalization lemma to the~$\O_{\Spec(\ZZ)}$-algebra~$A^\sim$.

Writing~$A =
\ZZ[X_1,\ldots,X_n]/(f_1,\ldots,f_m)$, we thus obtain that, internally,
it's \notnot the case that~$1 \in (f_1,\ldots,f_m)$ or that there
is a finite injective
morphism~$\O_{\Spec(\ZZ)}[Y_1,\ldots,Y_r] \to
\O_{\Spec(\ZZ)}[X_1,\ldots,X_n]/(f_1,\ldots,f_m) = A^\sim$ for some number~$r
\geq 0$. The latter case can't occur: Assume that it does. Since any positive natural number is invertible in~$A^\sim$ and finite
injective homomorphisms of rings reflect invertibility, any positive natural
number is also invertible in~$\O_{\Spec(\ZZ)}[Y_1,\ldots,Y_r]$ and therefore
in~$\O_{\Spec(\ZZ)}$. This, however, is false.

We therefore have~$\Sh(\Spec(\ZZ)) \models \neg\neg(1 \in (f_1,\ldots,f_m))$.
Thus~$D(h) \models 1 \in (f_1,\ldots,f_m)$ for some dense open subset~$D(h)
\subseteq \Spec(\ZZ)$. This implies that~$h^l \in (f_1,\ldots,f_m)$ for some~$l
\geq 0$. Therefore~$h^l = 0$ in~$A$; on the other hand,~$h^l$ is invertible
in~$A$. Thus~$1 = 0$ in~$A$.
\end{proof}

\section{Relative spectrum}
\label{sect:relative-spectrum}

Recall that if~$\A$ is a quasicoherent~$\O_X$-algebra on a scheme~$X$, one can
construct the \emph{relative spectrum}~$\RelSpec_X{\A}$ by appropriately
gluing the spectra~$\Spec \Gamma(U,\A)$ where~$U$ ranges over the affine opens
of~$X$. This relative spectrum comes equipped with a canonical
morphism~$\RelSpec_X{\A} \to X$.

From the internal point of view of~$\Sh(X)$, the sheaf~$\A$ looks just like a
plain algebra, to which therefore the usual (absolute) spectrum construction
can be applied. One could hope that this construction yields the relative
spectrum.

In this section, we discuss generalities on how to make sense of this internal
construction; we show that this proposed construction is too naive and doesn't
yield the relative spectrum; we give a refined internal construction which does
yield the relative spectrum, discuss its relation to the naive construction,
and phrase it in topos-theoretic terms; and we deduce, as an application, a
description of limits in the category of locally ringed spaces. We also cover
the relative Proj construction.

In much of the following, it's not actually necessary that~$X$ is a scheme
and~$\A$ is a quasicoherent algebra. If~$X$ is not a scheme or~$\A$ is not
quasicoherent, then~$\RelSpec_X(\A)$ might fail to be a scheme and can of
course not be constructed by gluing usual spectra, but it still exists as
a more general kind of space and still verifies a meaningful universal
property. We give details on this generalization below.

\subsection{Internal locales}\label{sect:internal-locales}
Let~$X$ be a topological space (or a locale). A
fundamental fact in the theory of locales is that there is a canonical
equivalence between the category of \emph{locales over~$X$} -- that is
locales~$Y$ equipped with a morphism~$Y \to X$ -- and \emph{internal locales
in~$\Sh(X)$}~\cite[p.~49]{johnstone:point}. An internal locale in a topos~$\E$ is given by an object~$L$ of~$\E$
(the internal frame of opens of the locale) together with a binary
relation~$(\preceq) \hookrightarrow L \times L$ such that the axioms for a
locale hold from the internal point of view. (For our purposes, we do not need a
precise wording of these axioms.)

The equivalence is described as follows: A locale~$f : Y \to X$ over~$X$
induces an internal locale~$I(Y)$ with object of opens given by~$\Open(I(Y)) \defeq
f_* \Omega_{\Sh(Y)} \in \Sh(X)$, where~$f_*$ is the pushforward functor
and~$\Omega_{\Sh(Y)}$ is the object of truth values in the topos of sheaves
on~$Y$. Conversely, an internal locale given by an internal frame~$\L \in \Sh(X)$ induces an (external)
locale~$E(\L)$ with frame of opens given by~$\Open(E(\L)) \defeq \Gamma(X,\L)$.
This comes equipped with a canonical morphism~$Y \to X$ of locales which we do
not need to describe explicitly~\cite[Section~C1.6]{johnstone:elephant}.

As a special case, the internalization of the trivial locale~$\id : X \to X$
over~$X$ has as frame of opens the object~$\id_* \Omega_{\Sh(X)} =
\Omega_{\Sh(X)} = \P(1)$. This is precisely the frame of opens of the
one-point space. Thus~$I(X) \cong \pt$. This illustrates the intuition
behind working internally in~$\Sh(X)$: From the perspective
of~$\Sh(X)$, the space~$X$ looks like the one-point space (even if in fact it
is not).

One can associate to an internal locale~$T$ in a topos~$\E$ a topos of internal
sheaves on it:~$\Sh_\E(T)$. The correspondence is made in such a way that the topos of
sheaves on a locale~$Y$ over~$X$ is equivalent to the topos of sheaves on the
internal locale~$I(Y)$: $\Sh(Y) \simeq \Sh_{\Sh(X)}(I(Y))$.

There is no similarly nice correspondence between topological spaces
over~$X$ and internal topological spaces
in~$\Sh(X)$~\cite[Corollary~C1.6.7]{johnstone:elephant}. This is one of the
reasons why locales are better suited for working internally and for switching
between internal and external perspectives.

For verification of properties of such sheaves, the \emph{idempotency} of the
internal language is useful: If~$\varphi$ is a formula over~$Y$, then
\[ \Sh(Y) \models \varphi \qquad\text{if and only if}\qquad
  \Sh(X) \models \speak{$\Sh(I(Y)) \models \varphi$}. \]
Here we're abusing notation in two ways. Firstly, the formula~$\varphi$ has to
be appropriately interpreted in the expression~``$\Sh(I(Y)) \models \varphi$''.
Secondly, the expression~``$\Sh(I(Y))$'' doesn't actually refer to the
category~$\Sh_{\Sh(X)}(I(Y))$, but to the \emph{locally internal} category induced by
the canonical geometric morphism~$\Sh_{\Sh(X)}(I(Y)) \to \Sh(X)$. We give some
details on this point in Section~\ref{sect:relation-big-little}. However, in
the situations encountered in this section, the meaning will always be
reasonably clear.

\subsection{The spectrum of a ring as a locale}
\label{sect:spectrum-as-a-locale}
Recall that the spectrum
of a ring~$A$ is usually constructed as the set
\[ \Spec A \defeq \{ \ppp \subseteq A \,|\,
  \text{$\ppp$ is a prime ideal} \} \]
endowed with a certain topology and a sheaf of rings~$\O_{\Spec A}$. From an
intuitionistic (and thus internal) point of view, this construction does not
work well: Prime ideals are intuitionistically much more elusive than
classically, where one can appeal to Zorn's lemma to obtain maximal (and thus
prime) ideals. More to the point, one cannot show that this construction of
the spectrum as a topological space verifies the expected universal property,
namely
\[ \Hom_\LRS(X, \Spec A) \cong \Hom_\Ring(A, \Gamma(X, \O_X)) \]
for all locally ringed spaces~$X$ (or some variant of this property involving
more general kinds of spaces).

On the other hand, the frame of opens of~$\Spec A$ admits a simple
description not requiring the notion of prime ideals:
\[ \Open(\Spec A) \cong \{ \aaa \subseteq A \,|\,
  \text{$\aaa$ is a radical ideal} \}. \]
An open subset~$U \subseteq \Spec A$ corresponds to the radical ideal~$\{ h \in
A \,|\, D(h) \subseteq U \}$ (so in particular, the open subset~$D(f)$
corresponds to the radical ideal~$\sqrt{(f)}$); conversely, a radical ideal~$\aaa$
corresponds to the open subset~$\bigcup_{h \in \aaa} D(h)$.

Thus, in an intuitionistic context, we will construct the spectrum of a ring~$A$
as a locale, not as a topological space, and adopt the following definition.

\begin{defn}\label{defn:spectrum-as-a-locale}
The \emph{spectrum}~$\Spec(A)$ of a ring~$A$ is the locale whose frame of
opens is the frame of radical ideals of~$A$. We endow it with the structure
sheaf~$\O_{\Spec(A)} \defeq \ul{A}[\F^{-1}]$, where~$\F$ is the generic filter
as described in Section~\ref{sect:generic-filter}.\end{defn}

This construction has the expected
universal property, namely that it is adjoint to the global functions functor:
\[ \Hom_{\LRL}(X, \Spec A) \cong \Hom_{\Ring}(A, \Gamma(X,\O_X)). \]
Here, ``$\LRL$'' refers to the category of \emph{locally ringed locales}, \ie
locales~$X$ equipped with a sheaf of rings~$\O_X$ such that from the internal point of
view of~$\Sh(X)$, the ring~$\O_X$ is a local ring. A morphism~$Y \to X$ of
locally ringed locales consists of a locale morphism~$f : Y \to X$ and a
morphism~$f^\sharp : f^{-1} \O_X \to \O_Y$ of sheaves of rings on~$Y$ such that, from the
internal point of view of~$\Sh(Y)$, the ring homomorphism~$f^\sharp$ is a local
homomorphism. The notion of a locally ringed locale is thus a straightforward
generalization of that of a locally ringed space.

Schemes are usually regarded as locally ringed spaces, not
as locally ringed locales. However, in a classical
context where the axiom of choice is available, schemes are \emph{sober}
topological spaces~\stacksproject{01IS}. For sober topological spaces, the passage from the space to
its induced locale (forgetting the set of points and only keeping the frame of
open subsets) doesn't lose information: The category of sober topological
spaces with arbitrary continuous maps embeds into the category of locales as a
full subcategory. Therefore the category of schemes can just as well be viewed
as a full subcategory of the category of locally ringed locales.

Describing morphisms between locally ringed locales is just as simple as
describing morphisms between locally ringed spaces; using the viewpoint of
classifying locales, one may even pretend that it suffices to give a map of
points. We expand on this in Section~\ref{sect:constructive-scheme-theory}.

The importance of a locale-theoretic approach to spectra of rings, especially in
relative situations, has also been stressed by Lurie~\cite[p.~37]{lurie:dag5}.

\begin{rem}In contrast to the prime spectrum, the spectrum of maximal ideals
can't in general be realized as a locale. This is because the maximal spectrum
is in general not sober; its soberification is the subspace of the prime
spectrum consisting of those prime ideals which are intersections of maximal
ideals. (For Jacobson rings, every prime ideal is of this form.)
\end{rem}

{\tocless

\subsection*{Points of the locale-theoretic spectrum}
Constructing the spectrum as a locale instead of a topological space
sidesteps any issues with prime ideals, since points are not a defining
ingredient of a locale. However, points are still meaningful as a \emph{derived
concept}: A point of locale~$X$ is a morphism~$\pt \to X$, where~$\pt$ is the
terminal locale, the locale corresponding to the one-point topological space
with frame of opens~$\P(1) = \Omega$ (where~$1 = \{ \star \}$ is a one-element set).
Therefore it's still an interesting question what the points of the
locale~$\Spec(A)$ look like.

\begin{prop}\label{prop:points-spectrum}
Let~$A$ be a ring. Then the points of the locale~$\Spec(A)$ are in
canonical one-to-one correspondence with the filters of~$A$
(as in Definition~\ref{defn:filter}), even intuitionistically.\end{prop}

\begin{proof}The points of a locale~$X$ are in canonical one-to-one
correspondence with the \emph{completely prime filters} of~$\Open(X)$,
subsets~$K \subseteq \Open(X)$ which are upward-closed, downward-directed, and
have the property that, whenever a supremum of a set~$M \subseteq \Open(X)$ is
contained in~$K$, then so is some element of~$M$.

Such a completely prime filter~$K \subseteq \Open(\Spec(A))$ corresponds to the
ring-theoretic filter
\[ F \defeq \{ s \? A \,|\, \sqrt{(s)} \in K \} \subseteq A, \]
and a ring-theoretic filter~$F \subseteq A$ corresponds to the completely prime
filter
\[ K \defeq \{ \aaa \? \Open(\Spec(A)) \,|\,
  \text{$\aaa \cap F$ is inhabited} \}. \]
We omit the required routine verifications.
\end{proof}

In classical logic, where complementation yields a one-to-one correspondence
between filters and prime ideals, the points of~$\Spec(A)$ are therefore in
canonical bijection with the prime ideals of~$A$, just as one would expect.

Observing that intuitionistically the points of the locale~$\Spec(A)$ are
filters, not prime ideals, one might wonder: Is the locale-theoretic approach
really necessary? Wouldn't it suffice to define~$\Spec(A)$ as the topological
space of filters of~$A$? Indeed, for some time this was
believed~\cite[Section~3]{lawvere:icm-address}; however, this hope turned out
to be too naive: Joyal gave an explicit example of a nontrivial ring in a
certain topos without any filters~\cite[pp.~200f.]{tierney:spectrum}, thus
showing that the construction can't have the expected universal property and
that therefore a true pointfree approach as provided by lattice theory/locale theory~\cite{cls:spectral-schemes}, topos
theory, or formal topology~\cite{schuster:formal-zariski} is necessary to construct the spectrum in an
intuitionistic context.\footnote{When following
reference~\cite{tierney:spectrum}, note that Tierney calls ``primes'' what we
call ``filters''. Joyal's example was none other than the ring~$\afflz$
in the functor category~$[\Ring_\mathrm{fp}, \Set]$. The big Zariski topos
of~$\Spec(\ZZ)$, when defined using the parsimonious sites, is a subtopos of
that topos; in it, the ring~$\affl$ does have filters, for instance the filter
of units. These two facts are not contradictory, since not having any filters
is not a geometric implication and is therefore not guaranteed to be preserved by inverse image
parts of geometric morphisms.}

\subsection*{The spectrum as a classifying locale}

The fact that the points of~$\Spec(A)$ are in canonical one-to-one
correspondence with the filters of~$A$ is a shadow of a more general fact.
Namely, for any locale~$X$ (and in fact any topos), maps~$X \to \Spec(A)$ are
in canonical one-to-one correspondence with the internal filters of~$A$
in~$\Sh(X)$, that is subsheaves of the constant sheaf~$\ul{A}$ satisfying the
filter axioms from the point of view of the internal language of~$\Sh(X)$:
The locale~$\Spec(A)$ is the \emph{classifying locale of the theory of filters
of~$A$}.

The fact about the points of~$\Spec(A)$ can be recovered from this
observation as follows. A point of~$\Spec(A)$ is a map~$\pt \to \Spec(A)$ and
therefore corresponds to a subsheaf of the constant sheaf~$\ul{A}$
satisfying the filter axioms from the point of view of~$\Sh(\pt)$. Since~$\Sh(\pt)
\simeq \Set$, such a subsheaf amounts to a subset of~$A$ satisfying the filter
axioms.

The notion of classifying locales provides a pleasant way of approaching the
problem of constructing a space of models of a \emph{propositional geometric
theory} (in the case of the spectrum the theory of filters), simultaneously
streamlining the usual topological approach and generalizing it to work in an
intuitionistic context: Instead of first constructing the \emph{set} of models
(filters of~$A$) and then manually endowing this set with a suitable topology
(the Zariski topology), one can simply consider the \emph{locale} of models,
that is the classifying locale of the theory. Its sets of points coincides with
the set of models of the topological approach, but the locale is not determined
by its sets of points, facilitating a better behavior in contexts where the
points might be elusive.

Put more concisely, the topological space of filters doesn't work well in an
intuitionistic context, but the locale of filters does. This phenomenon is an
instance of Vickers's motto ``if you define points by a geometric theory, then
the topology is implicit''. A lucid expository account of the theory of
classifying locales can be found in a survey article by
him~\cite{vickers:locales-toposes}.

\begin{rem}\label{rem:theory-of-filters}
For comparison with a refined geometric theory discussed below, we describe the
geometric theory of filters of~$A$ here explicitly. It has one atomic
proposition~``$s \in F$'' for each element~$s \? A$, and its axioms are given by the
following axiom schemes:
\begin{enumerate}
\item $\top \vdash 1 \in F$
\item $st \in F \dashv\vdash s \in F \wedge t \in F$ (two axioms for each $s,t\?A$)
\item $0 \in F \vdash \bot$
\item $s+t \in F \vdash s \in F \vee t \in F$ (one axiom for each $s,t\?A$)
\end{enumerate}
\end{rem}

\subsection*{A trivial case}
For later use, we study the question when the spectrum is the one-point space.
The answer is well-known classically, but since we want to use this result in
an internal context, we have to give an intuitionistic proof.
\begin{lemma}\label{lemma:spectrum-one-point}
Let~$A$ be a ring. Its spectrum is a one-point
space (as a locale) if and only if~$1 \neq 0$ in~$A$ and any element of~$A$ is nilpotent or
invertible.\end{lemma}
\begin{proof}The locale~$\Spec A$ is a one-point space if and only if the
unique continuous map~$\Spec(A) \to \pt$ of locales is an isomorphism. This is the case if and only
if the canonical frame homomorphism
\[ \begin{array}{@{}rcl@{}}
  \Omega = \P(1) &\longrightarrow& \Open(\Spec A) \\[0.1em]
  \varphi &\longmapsto& \aaa_\varphi \defeq \sup\{\sqrt{(1)} \,|\, \varphi \} =
  \{ x \? A \,|\, \speak{$x$ nilpotent} \vee \varphi \}
\end{array} \]
is surjective and reflects the ordering (and is therefore automatically
injective). If~$1 = 0$ in~$A$, this homomorphism is not injective, since~$\bot$
and~$\top$ get both mapped to~$\sqrt{(0)}$. For the rest of the proof, we'll
therefore assume that~$1 \neq 0$ in~$A$.

Under this assumption, the homomorphism reflects the ordering: If~$\aaa_\varphi \subseteq \aaa_\psi$,
then~$(1 \in \aaa_\varphi) \Rightarrow (1 \in \aaa_\psi)$. Since the unit
of~$A$ is not nilpotent, this amounts to~$\varphi \Rightarrow \psi$.

The homomorphism is surjective if and only if for any radical ideal~$\aaa \subseteq A$,
it holds that~$\aaa = \{ x \? A \,|\, \speak{$x$ nilpotent} \vee \varphi \}$
for some proposition~$\varphi$. By considering the condition~``$1 \in \aaa$'',
it follows that this proposition~$\varphi$ must be equivalent to the
proposition~``$1 \in \aaa$'' (if it is at all possible to write~$\aaa$ in such
a way).

So the map is surjective if and only if for any radical ideal~$\aaa \subseteq
A$ and any element~$x$ of~$A$ it holds that
\[ x \in \aaa \quad\Longleftrightarrow\quad
  \speak{$x$ nilpotent} \vee (1 \in \aaa). \]
The ``if'' direction always holds. If any element of~$A$ is nilpotent or
invertible, the ``only if'' direction holds as well (for any~$\aaa$ and any~$x$).
Conversely, if the ``only if'' direction holds, then any element of~$A$ is
nilpotent or invertible. This follows by
considering the radical ideal~$\sqrt{(f)}$ for an element~$f \? A$.
\end{proof}

\begin{rem}The structure sheaf~$\O_X$ of a scheme fulfills almost, but not
quite, the condition given in Lemma~\ref{lemma:spectrum-one-point}: By
Proposition~\ref{prop:neginvnilpotent}, it has the property that
any element which is not invertible is nilpotent. In classical logic, this
statement is equivalent to the statement that every element is nilpotent or invertible.
However, intuitionistically the former is a weaker statement than the latter.
This observation entails that the internally constructed spectrum does
\emph{not} coincide with the relative spectrum, and that instead a refined
approach is necessary. Section~\ref{sect:rel-spec-as-ordinary-spec} is devoted
to studying this difference.
\end{rem}

}

\subsection{Digression: Further topologies on the set of prime ideals}
\label{sect:flat-constructible-topologies}

The Zariski topology is not the only interesting topology on the set of prime
ideals. For instance, the constructible topology and the flat topology studied
by Tarizadeh~\cite{tarizadeh:flat} (also called ``co-Zariski topology'') too
have their uses.

The universal properties given in the following two propositions should be
compared with the following way of phrasing the universal property of the
ordinary locale-theoretic spectrum. The usual phrasing employs the
categories~$\RL$ and~$\LRL$ of (locally) ringed locales, therefore emphasizing
the spatial character. But the dual categories~$\RL^\op$ and~$\LRL^\op$ can be
used just as well; since the morphisms in~$\RL^\op$ and~$\LRL^\op$ go in the
direction of the ring-theoretic parts, they can be thought of as the category
of \emph{all} rings respectively \emph{all} local rings, where~``all'' refers
to the fact that these categories don't only include the (local) rings
in~$\Set$, but the (local) rings in arbitrary localic sheaf toposes.

Formulated using~$\RL^\op$ and~$\LRL^\op$, and adopting the notation to
suppress mention of the involved spaces (instead of the involved sheaves
of rings), the universal property of~$\Spec(A)$ reads as follows: For any
local ring~$\O_Y$ over any locale~$Y$,
\[ \Hom_{\LRL^\op}(\O_{\Spec(A)}, \O_Y) \cong
  \Hom_{\RL^\op}(A, \O_Y). \]
The morphism~$A \to \O_{\Spec(A)}$ in~$\RL^\op$ is therefore the
\emph{universal localization} of~$A$.

\begin{prop}Let~$A$ be a ring. The locale given by the space of prime ideals
of~$A$ with the flat topology is the classifying locale of prime ideals of~$A$.
Equipped with~$\ul{A}/\P$ as structure sheaf, where~$\P$ is the generic prime
ideal, it is the universal way of mapping~$A$ to an integral domain in the weak
sense (as defined in Section~\ref{sect:integrality}; with morphisms of
weak integral domains taken to be injective ring homomorphisms).
\end{prop}

\begin{proof}See~\cite[Proposition~4.5]{johnstone:rings-fields-and-spectra}.\end{proof}

\begin{prop}Let~$A$ be a ring. The locale given by the space of prime ideals
of~$A$ with the constructible topology is the classifying locale of detachable
prime ideals (or equivalently detachable filters) of~$A$. Equipped
with~$\ul{A}/\P$ as structure sheaf, where~$\P$ is the generic prime ideal, it
is the universal way of mapping~$A$ to an integral domain in the strong sense.
Equipped with~$\ul{A}[\F^{-1}]$, where~$\F$ is the generic filter, it is the
universal way of mapping~$A$ to a local ring in which invertibility is
decidable.\end{prop}

\begin{proof}This is mostly covered
in~\cite[p.~253]{johnstone:rings-fields-and-spectra}.\end{proof}

In constructive mathematics, a subset~$U \subseteq A$ is \emph{detachable} if
and only if for every element~$a\?A$, either $a \in U$ or~$a \not\in U$. While
intuitionistically the complement of a filter might fail to be a prime ideal
and the complement of a prime ideal might fail to be a filter, the complement
of a detachable filter is a detachable prime ideal, and vice versa.

\subsection{The relative spectrum as an ordinary spectrum\except{toc}{ from the internal
point of view}}\label{sect:rel-spec-as-ordinary-spec}

Let~$X$ be a scheme and~$\A$ be a quasicoherent~$\O_X$-algebra.
Since~$\A$ looks like a plain algebra from the internal perspective
of~$\Sh(X)$, we can consider its internally defined spectrum. This is a locale
internal to~$\Sh(X)$; we might hope that its externalization is precisely the
relative spectrum of~$\A$ (considered as a locale):
\[ E(\Spec \A) \stackrel{?}{\cong} \RelSpec_X{\A}. \]
However, this turns out to be too naive. The locale~$E(\Spec(\A))$ is equipped
with a map to~$X$, being an externalization of a locale internal to~$\Sh(X)$,
and it is equipped with a sheaf of rings (because we can transport the
internally defined structure sheaf along the
equivalence~$\Sh_{\Sh(X)}(\Spec(\A)) \simeq \Sh(E(\Spec(A)))$. Furthermore,
this sheaf of rings is local, since we know
\[ \Sh(X) \models \speak{$\Sh(\Spec(\A)) \models \speak{%
  $\O_{\Spec(\A)}$ is a local ring}$} \]
which by idempotency of the internal language is equivalent to
\[ \Sh(E(\Spec(\A))) \models \speak{$\O_{\Spec(\A)}$ is a local ring}. \]

However, the map~$E(\Spec(\A)) \to X$ is only part of a morphism of ringed
locales, not of locally ringed locales (even though domain and codomain happen
to be locally ringed): Internally, the morphism~$(\Spec(\A),\O_{\Spec(\A)}) \to
(\pt,\O_X)$ of ringed locales, which is defined using the~$\O_X$-algebra
structure of~$\A$, is not a morphism of locally ringed locales (even though
domain and codomain happen
to be locally ringed).

In contrast, the true relative spectrum~$\RelSpec_X(\A)$ is equipped with a
morphism of locally ringed locales to~$X$.

It's illuminating to compare the different universal properties
of~$E(\Spec(\A))$ and~$\RelSpec_X(\A)$. There is a canonical
morphism~$E(\Spec(\A)) \to E(\Spec(\O_X))$ of locally ringed locales (the
externalization of the canonical morphism~$\Spec(\A) \to \Spec(\O_X)$ given by
the~$\O_X$-algebra structure of~$\A$), but in general, the
locales $E(\Spec(\O_X))$ and~$X$ are not isomorphic.

As we justify below, the externalization of the internally
defined spectrum has the universal property
\[
  \Hom_{\LRL/E(\Spec\O_X)}(Y, E(\Spec\A)) \cong
    \Hom_{\O_X}(\A, \mu_* \O_Y)
\]
for all locally ringed locales~$Y$ over~$E(\Spec\O_X)$. Here,~$\mu$ is the
structure morphism~$Y \to \Spec\O_X$, $E(\Spec\O_X)$ is the locally ringed
locale associated to the internally defined spectrum of~$\O_X$,
and~$\LRL_{\Sh(X)}$ is the category of locally ringed locales internal
to~$\Sh(X)$. In contrast, the relative spectrum has the different universal property
\[
  \Hom_{\LRL/X}(Y, \RelSpec_X{\A}) \cong
    \Hom_{\O_X}(\A, \mu_* \O_Y)
\]
for all locally ringed locales~$Y$ over~$X$.\footnote{If~$X$ is a scheme
and~$\A$ is quasicoherent, this universal property is well-known, even though
it's usually only stated for schemes~$Y$ over~$X$ instead of general locally
ringed locales over~$X$. In any case, we take this universal property as the
definition of what the relative spectrum should be.} The crucial
difference is that in general, the internally defined locally ringed
locale~$\Spec\O_X$ does \emph{not} coincide with the internal locally ringed
locale~$(\pt,\O_X)$ (which is simply~$(X,\O_X)$ from the external point of
view). More succinctly, the functor~$E \circ \Spec$ is an adjoint to the
pushforward-of-sheaf-of-functions functor~$\LRL/E(\Spec\O_X) \to \Alg(\O_X)^\op$, while the
relative spectrum functor is an adjoint to the analogous functor~$\LRL/X
\to \Alg(\O_X)^\op$.

The universal property of~$E(\Spec(\A))$ can be determined as follows.
From the internal point of view of~$\Sh(X)$, the locally ringed
locale~$E(\Spec(\A))$ looks like the ordinary locale-theoretic
spectrum~$\Spec(\A)$ and therefore has the universal property
\[ \Hom_{\LRL}(Y, \Spec(\A)) \cong
  \Hom_{\Ring}(\A, \Gamma(Y, \O_Y)) \]
for any locally ringed locale~$Y$.\footnote{Externally, this implies that for any
locally ringed locale over the underlying locale of~$X$ (that is, for any
locale~$Y$ equipped with a morphism~$\mu : Y \to X$ and a local sheaf of
rings), we have
\[ \Hom_{\mathrm{LR}(\mathrm{L}/X)}(Y, E(\Spec(\A))) \cong
  \Hom_{\Ring_{\Sh(X)}}(\A, \mu_*\O_Y). \]}
If we restrict the right-hand side to the
set of~$\O_X$-algebra homomorphisms, the left-hand side restricts to the set of
morphisms~$Y \to \Spec(\A)$ of locally ringed locales over the locally ringed
locale~$\Spec(\O_X)$. So we have
\[ \Hom_{\LRL/\Spec(\O_X)}(Y, \Spec(\A)) \cong \Hom_{\Alg(\O_X)}(\A,
\Gamma(Y,\O_Y)). \]

This discussion took place in the internal universe of~$\Sh(X)$. Externally,
the displayed universal property implies that for any locally ringed
locale~$\mu : Y \to X$ over~$E(\Spec(\O_X))$,
\[ \Hom_{\LRL/E(\Spec(\O_X))}(Y, E(\Spec(\A))) \cong
  \Hom_{\Alg(\O_X)}(\A, \mu_*\O_Y), \]
as claimed above.

\begin{defn}\label{defn:local-spectrum}
Let~$R$ be a ring. Let~$A$ be an~$R$-algebra. The \emph{local
spectrum} of~$A$ over~$R$ is the locale~$\Spec(A|R)$ with frame of opens
given by
\begin{multline*}
  \quad\qquad\Open(\Spec(A|R)) \defeq
    \{ \aaa \subseteq A \,|\,
      \text{$\aaa$ is a radical ideal such that} \\
  \text{$\forall f\?R\_ \forall s\?A\_
    (\speak{$f$ \inv} \Rightarrow s \in \aaa) \Rightarrow fs \in \aaa$} \}.\qquad\quad
\end{multline*}
\end{defn}

We'll equip the local spectrum with the structure of a locally ringed locale
below. It is this refined construction which correctly internalizes the
relative spectrum:

\begin{thm}\label{thm:local-spectrum-yields-relative-spectrum}
Let~$X$ be a scheme (or a locally ringed locale). Let~$\A$ be
an~$\O_X$-algebra. Then the externalization~$E(\Spec(\A|\O_X))$ coincides
with~$\RelSpec_X(\A)$ as locally ringed locales over~$X$.\end{thm}

Before giving the proof, we want to clarify some details of the
construction.

Firstly, the base ring~$R$ directly enters the construction. This is in
contrast to the usual spectrum: If~$A$ is an~$R$-algebra, the construction
of~$\Spec(A)$ does not depend on the~$R$-algebra structure of~$A$. The algebra
structure only enters in the construction of a morphism~$\Spec(A) \to
\Spec(R)$.

Secondly, in the case that~$X$ is a scheme and~$\A$ is a
quasicoherent~$\O_X$-algebra, we can compare the externalization
of~$\Spec(\A|\O_X)$ with the result of the construction of~$\RelSpec_X(\A)$
by gluing spectra:

\begin{prop}Let~$X$ be a scheme. Let~$\A$ be a quasicoherent~$\O_X$-algebra.
Then~$E(\Spec(\A|\O_X))$ coincides with~$\RelSpec_X(\A)$ as locales over~$X$.
\end{prop}

\begin{proof}The condition
\[ \forall f\?\O_X\_ \forall s\?\A\_
    (\speak{$f$ \inv} \Rightarrow s \in \aaa) \Longrightarrow fs \in \aaa \]
appearing in Definition~\ref{defn:local-spectrum} is precisely the internal
quasicoherence condition of Corollary~\ref{cor:submodule-qcoh} (slightly
simplified in view that~$\aaa$ is a radical ideal). The sections of the
sheaf~$\brak{\Open(\Spec(\A|\O_X))}$ on an open subset~$U \subseteq X$ are
therefore precisely the quasicoherent sheaves of radical ideals~$\aaa
\hookrightarrow \A|_U$. Let~$\pi : \RelSpec_X(\A) \to X$ be the canonical
morphism. If~$U$ is affine, then
\[ \pi^{-1}U \cong \RelSpec_X(\A) \times_X U \cong \RelSpec_U(\A|_U) \cong
  \Spec(\Gamma(U,\A)) \]
is affine as well and
\begin{align*}
  \Gamma(U, \Open(I(\RelSpec_X(\A)))) &=
  \Gamma(U, \pi_* \Omega_{\RelSpec_X(\A)}) =
  \Omega_{\RelSpec_X(\A)}(\pi^{-1}U) \\
  &\cong \text{set of open subsets of~$\pi^{-1}U$} \\
  &\cong \text{set of open subsets of~$\Spec(\Gamma(U,\A))$} \\
  &\cong \text{set of radical ideals of~$\Gamma(U,\A)$} \\
  &\cong \text{set of quasicoherent sheaves of radical ideals of~$\A|_U$} \\
  &\cong \Gamma(U, \brak{\Open(\Spec(\A|\O_X))}).
\end{align*}
Therefore~$I(\RelSpec_X(\A))$ is canonically isomorphic to~$\Spec(\A|\O_X)$ as
locales internal to~$\Sh(X)$. Expressed externally: The relative
spectrum~$\RelSpec_X(\A)$ coincides with the externalization
of~$\Spec(\A|\O_X)$ as locales over~$X$, as claimed.
\end{proof}

Thirdly, the partial order~$\Open(\Spec(A|R))$ is indeed a frame. A quick
way to verify this is to recognize that it is related to the frame of opens
of~$\Spec(A)$ by the formula
\[ \Open(\Spec(A|R)) =
  \{ \aaa \? \Open(\Spec(A)) \,|\, \aaa = \overline{\aaa} \}, \]
where~$(\aaa \mapsto \overline{\aaa})$ is the quasicoherator described in
Remark~\ref{rem:quasicoherator-knaster-tarski}. Since the quasicoherator
satisfies the axioms for a nucleus, this formula exhibits~$\Spec(A|R)$ as a
sublocale of~$\Spec(A)$. In particular, suprema are computed
in~$\Open(\Spec(A|R))$ by applying the quasicoherator to the suprema computed
in~$\Open(\Spec(A))$. We denote the inclusion~$\Spec(A|R) \hookrightarrow
\Spec(A)$ by~``$i$''.

Lastly, it's interesting to know what the points of~$\Spec(A|R)$ are, even though these
don't determine~$\Spec(A)$.

\begin{defn}\label{defn:over-the-filter-of-units}
Let~$R$ be a ring. Let~$\varphi : R \to A$ be an algebra. A
filter~$F \subseteq A$ \emph{lies over the filter of units} if and only if
$\varphi^{-1}F \subseteq R^\times$, that is if
\[ \varphi(r) \in F \Longrightarrow \text{$r$ is invertible in $R$} \]
for all~$r \? R$. (The reverse inclusion~``$\varphi^{-1}F \supseteq R^\times$''
holds automatically.)\end{defn}

This definition will mostly be used in situations where the ring~$R$ is local,
in which case the subset~$R^\times$ is actually a filter and the phrase
``filter of units'' is therefore justified.

It's illuminating to consider Definition~\ref{defn:over-the-filter-of-units} in
a classical context, even though the use case we have in mind is to apply it in
the internal language of the little Zariski topos of a base scheme.
Classically, a filter~$F$ lies over the filter of units if and only
if~$\varphi^{-1}\ppp \supseteq R \setminus R^\times$, where~$\ppp \defeq
F^c = A \setminus F$ is the prime ideal associated to~$F$. If~$R$ is local, the
set~$R \setminus R^\times$ is the unique maximal ideal~$\mmm$ of~$R$. Thus~$F$
lies over the filter of units if and only if~$\ppp$ lies over the maximal
ideal.

\begin{prop}Let~$R$ be a ring. Let~$\varphi : R \to A$ be an~$R$-algebra. Then
the points of~$\Spec(A|R)$ are intuitionistically in canonical one-to-one
correspondence with those filters of~$A$ which lie over the filter of units.
\end{prop}

\begin{proof}The correspondence outlined in
Proposition~\ref{prop:points-spectrum} can be adapted to the situation at hand.
A completely prime filter~$K \subseteq \Open(\Spec(A|R))$ corresponds to the
ring-theoretic filter
\[ F \defeq \{ s\?A \,|\, \overline{\sqrt{(s)}} \in K \} \]
and a ring-theoretic filter~$F$ corresponds to the completely prime filter
\[ K \defeq \{ \aaa\?\Open(\Spec(A|R)) \,|\, \text{$\aaa \cap F$ is inhabited} \}. \]
It's instructive to perform some of the necessary verifications, to see how
the quasicoherator is used, even though
Proposition~\ref{prop:local-spectrum-classify} will subsume this
correspondence.

The filter~$F$ corresponding to~$K$ has the displayed property for the
following reason. Let~$\varphi(r) \in F$. We want to verify that~$r$ is
invertible in~$R$. Under the assumption that~$r$ is invertible in~$R$,
it's trivial that~$1$ is an element of
\begin{align*}
  \aaa &\defeq \sup \{ \sqrt{(1)} \,|\, \text{$r$ is invertible in $R$} \} \\
  &\phantom{\vcentcolon}= \{ s\?A \,|\, \text{$s$ is nilpotent or $r$ is invertible in $R$} \}
  \in \Open(\Spec(A)).
\end{align*}
Therefore, without any assumption on~$r$, we have that~$r \cdot 1 = \varphi(r)$ is an
element of~$\overline{\aaa}$ and therefore~$\overline{\sqrt{(\varphi(r))}} \subseteq
\overline{\aaa}$. Since~$K$ is upward-closed, it follows that~$\overline{\aaa}
\in K$. Since~$\overline{\aaa}$ is the supremum of the set~$\{ \sqrt{(1)} \,|\,
\text{$r$ is invertible} \}$ in~$\Open(\Spec(A|R))$ and~$K$ is completely prime, it
follows that this set is inhabited. Thus~$r$ is invertible in~$R$.

The set~$K$ corresponding to a ring-theoretic filter~$F$ is completely prime
for the following reason. Let~$\sup_i \aaa_i = \overline{\sqrt{\sum_i \aaa_i}}
\in K$. Then~$\overline{\sqrt{\sum_i \aaa_i}} \cap F$ is inhabited. By the
special assumption on~$F$, the intersection~$\sqrt{\sum_i \aaa_i} \cap F$ is inhabited
as well: In the case that~$X$ is a scheme, this follows easily using the
description of the quasicoherator given in
Proposition~\ref{prop:quasicoherator-arbitrary-algebra}. In the general case,
we use the proof scheme outlined in
Remark~\ref{rem:quasicoherator-knaster-tarski} -- using the notation of that
remark, if~$P(\bbb) \cap F$ is inhabited, then~$\bbb \cap F$ is as well.

A short calculation using the filter axioms then shows that there
exists an index~$i$ such that~$\aaa_i \cap F$ is inhabited.
\end{proof}

\begin{prop}\label{prop:local-spectrum-classify}
Let~$R$ be a ring. Let~$\varphi : R \to A$ be an algebra. Then~$\Spec(A|R)$ is
the classifying locale of the theory of filters of~$A$ which lie over the
filter of units, that is of the geometric theory with atomic propositions~``$s
\in F$'' for~$s\?A$ and axioms given by the following axiom schemes:
\begin{enumerate}
\item $\top \vdash 1 \in F$
\item $st \in F \dashv\vdash s \in F \wedge t \in F$ (two axioms for each $s,t\?A$)
\item $0 \in F \vdash \bot$
\item $s+t \in F \vdash s \in F \vee t \in F$ (one axiom for each $s,t\?A$)
\item $\varphi(r) \in F \vdash \bigvee \{ \top \,|\, \textnormal{$r$ is invertible in
$R$} \}$ (one axiom for each $r\?R$)
\end{enumerate}
\end{prop}

\begin{proof}The frame of the classifying locale of the given theory~$T$ is the
free frame on generators~``$s \in F$'' for~$s\?A$ subject to the relations
given by the axioms of the theory. More explicitly, it's the Lindenbaum
algebra~$L(T)$ of the theory, so its elements are the formulas of the theory up
to provable equivalence and the ordering is defined by~$[\varphi] \preceq
[\psi] \vcentcolon\Leftrightarrow (\varphi \vdash \psi)$. We want to verify that this frame is
isomorphic to~$\Open(\Spec(A|R))$.

We define a frame homomorphism~$L(T) \to \Open(\Spec(A|R))$ by sending the
generators~$[s \in F]$ to the radical ideal~$\overline{\sqrt{(s)}}$. This
respects the relations and therefore gives a well-defined map. The map is
surjective, since a preimage to~$\aaa \? \Open(\Spec(A|R))$
is~$[\bigvee_{s\in\aaa} (s \in F)]$. To verify that it is an isomorphism of
frames, we therefore only have to verify that it reflects the ordering.

By the axiom schemes~(1) and~(2), any formula of~$T$ is provably equivalent to
a formula of the form~$\bigvee_i (s_i \in F)$. It therefore suffices to verify
that, for any families~$(s_i)_i$ and~$(t_j)_j$ such that
$\overline{\sqrt{(s_i)_i}} \subseteq \overline{\sqrt{(t_j)_j}}$, the sequent
$\bigvee_i (s_i \in F) \vdash \bigvee_j (t_j \in F)$ is derivable. We'll show
more generally: If~$\aaa$ and~$\bbb$ are radical ideals such
that~$\overline{\aaa} \subseteq \overline{\bbb}$,
then~$\bigvee_{s\in\aaa}(s\in F) \vdash \bigvee_{t\in\bbb}(t\in F)$. This
follows from the following chain of deductions:
\[ \bigvee_{s\in\aaa}(s\in F) \vdash
  \bigvee_{s\in\overline{\aaa}}(s\in F) \vdash
  \bigvee_{s\in\overline{\bbb}}(s\in F) \vdash
  \bigvee_{s\in\bbb}(s\in F). \]
All but the final step are trivial. The final step is an application of the
general proof scheme outlined in
Remark~\ref{rem:quasicoherator-knaster-tarski}. In the notation of that remark,
we set~$\alpha(\J) \defeq [\bigvee_{s\in\J}(s\in F)]$ and exploit that, if~$s
\in P(\J)$, then~$s \in F \vdash \bigvee_{t\in\J} (t \in F)$. This is
because~$s$ can be written as~$s^n = \sum_j a_j f_j u_j$ such that, for each~$j$,
if~$f_j$ is invertible in~$R$ then~$u_j \in \J$, and we have the
following chain of deductions.
\begin{align*}
  s \in F &\vdash s^n \in F \\
  &\vdash \bigvee_j (t_j f_j u_j \in F) \\
  &\vdash \bigvee_j (\varphi(f_j) \in F \wedge u_j \in F) \\
  &\vdash \bigvee_j \bigl(\bigvee\{\top\,|\,\text{$f_j$ invertible in $R$}\} \wedge u_j \in F\bigr) \\
  &\vdash \bigvee_j \bigvee\{(u_j\in F) \,|\, \text{$f_j$ invertible in $R$}\} \\
  &\vdash \bigvee_{t \in \J} (t \in F). \qedhere
\end{align*}
\end{proof}

\begin{lemma}\label{lemma:universal-property-local-spectrum}
Let~$R$ be a local ring. Let~$\varphi : R \to A$ be an~$R$-algebra.
Then, intuitionistically, the locale~$\Spec(A|R)$ carries a canonical structure
as a locally ringed locale over~$(\pt,R)$ and has the following universal
property: For any locally ringed locale~$(Y,\O_Y)$ over~$(\pt,R)$,
\[ \Hom_{\LRL/(\pt,R)}(Y, \Spec(A|R)) \cong \Hom_{\Alg(R)}(A, \Gamma(Y,\O_Y)). \]
\end{lemma}

\begin{proof}Since~$\Spec(A|R)$ is a sublocale of~$\Spec(A)$, we can
equip~$\Spec(A|R)$ with the restriction of~$\O_{\Spec(A)}$ to~$\Spec(A|R)$ as
the structure sheaf:
\[ \O_{\Spec(A|R)} \defeq
  i^{-1}\O_{\Spec(A)} =
  i^{-1}(\ul{A}[\F^{-1}]) \cong
  (i^{-1}\ul{A})[(i^{-1}\F)^{-1}] \cong
  \ul{A}[(i^{-1}\F)^{-1}]. \]
The generic filter~$\F$ was described in Section~\ref{sect:generic-filter}.
The penultimate isomorphism is because localizing is a geometric construction.
Since locality of a ring is a geometric implication, this structure sheaf is
indeed a local sheaf of rings. Thus~$\Spec(A|R)$ is a locally ringed locale.

Next, we have to describe a morphism~$(\Spec(A|R), \O_{\Spec(A|R)}) \to
(\pt,R)$. Locale-theoretically, this morphism is given by the unique map~$! :
\Spec(A|R) \to \pt$. The ring-theoretic part is given by the composition
\[ !^{-1}R = \ul{R} \longrightarrow
  \ul{A} \longrightarrow
  \ul{A}[(i^{-1}\F)^{-1}] =
  \O_{\Spec(A|R)}. \]
This homomorphism of rings which happen to be local is indeed a local
homomorphism, that is, it reflects invertibility. More precisely,
\[ \Spec(A|R) \models
  \forall f\?\ul{R}\_
  \speak{$\ul{\varphi}(f)$ is \inv\@ in~$\O_{\Spec(A|R)}$} \Rightarrow
  \speak{$f$ is \inv\@ in~$\ul{R}$}. \]
Denoting the modal operator associated to the sublocale inclusion~$\Spec(A|R)
\hookrightarrow \Spec(A)$ by~``$\Box$'', this statement is equivalent to
\[ \Spec(A) \models (\forall f\?\ul{R}\_
  \ul{\varphi}(f) \in \F \Rightarrow \speak{$f$ is \inv\@ in~$\ul{R}$})^\Box \]
by Theorem~\ref{thm:box-translation-semantically} and
Lemma~\ref{lemma:box-translation-sound}. To verify this, let~$s\?A$ and~$f\?R$
be given such that~$\sqrt{(s)} \models \varphi(f) \in \F$, that is,~$s \in
\sqrt{(\varphi(f))}$. We are to show that~$\sqrt{(s)} \models \Box(\speak{$f$
is invertible in~$\ul{R}$})$.

The largest open in~$\Spec(A)$ on which~$\speak{$f$ is invertible in~$\ul{R}$}$
holds is
\begin{align*}
  \aaa &\defeq \sup \{ \sqrt{(1)} \,|\, \text{$f$ is invertible in~$R$} \} \\
  &\phantom{\vcentcolon}=
  \{ t\?A \,|\, \text{$t$ is nilpotent or $f$ is invertible in $R$} \} \in
  \Open(\Spec(A)),
\end{align*}
by Lemma~\ref{lemma:properties-of-constant-sheaves}. Under the assumption
that~$f$ is invertible in~$R$, trivially~$1 \in \aaa$. Therefore, without any
assumptions on~$f$, we have that~$\varphi(f) \in \overline{\aaa}$.
Thus~$\sqrt{(\varphi(f))} \subseteq \overline{\aaa}$ and
therefore~$\sqrt{(\varphi(f))} \models \Box(\speak{$f$ is invertible
in~$\ul{R}$})$. Since~$\sqrt{(s)} \subseteq \sqrt{(\varphi(f))}$, the
monotonicity of the internal language implies~$\sqrt{(s)} \models
\Box(\speak{$f$ is invertible
in~$\ul{R}$})$.

Finally, we verify the universal property. Let~$Y$ be a locally ringed locale
over~$(\pt,R)$ and let a morphism~$A \to \Gamma(Y,\O_Y)$ of~$R$-algebras be
given. We like this data to uniquely induce a morphism~$Y \to \Spec(A|R)$ of
locally ringed locales over~$(\pt,R)$.

To obtain a locale-theoretic map~$f : Y \to \Spec(A|R)$, by
Proposition~\ref{prop:local-spectrum-classify} we need to specify a
filter of~$\ul{A}$ in~$\Sh(Y)$ which lies over the filter of units.
The given morphism~$A \to \Gamma(Y,\O_Y)$ induces a morphism~$\alpha : \ul{A}
\to \O_Y$ in~$\Sh(Y)$. Since~$\O_Y$ is a local ring, the subsheaf~$\O_Y^\times$
is a filter. Its preimage~$F \defeq \alpha^{-1}\O_Y^\times$ is the sought filter of~$\ul{A}$.
It lies over the filter of units because the composition~$\ul{R} \to \ul{A} \to
\O_Y$ is local. By the general theory, the pullback of the generic filter
in~$\Sh(\Spec(A|R))$ to~$\Sh(Y)$ along~$f$ is~$F$.

The ring-theoretic part of the sought morphism~$Y \to \Spec(A|R)$ of locally ringed
locales over~$(\pt,R)$ is the canonical homomorphism
\[ f^{-1}\O_{\Spec(A|R)} =
  f^{-1}(\ul{A}[(i^{-1}\F)^{-1}]) =
  \ul{A}[F^{-1}] \longrightarrow \O_Y \]
of local rings.

This finishes the description of the construction. We omit further
verifications that the construction works as claimed.
\end{proof}

\begin{rem}The modal operator~$\Box$ associated to the inclusion~$\Spec(A|R)
\hookrightarrow \Spec(A)$ can be defined in the internal language
of~$\Sh(\Spec(A))$. Namely, it's the smallest operator such that
the~$\Box$-translated statement
\[ (\speak{the morphism $\ul{R} \to \O_{\Spec(A)}$ is local})^\Box \]
holds. It is thus the smallest operator such that for
any~$f\?\ul{R}$ with~$\ul{\varphi}(f) \in \F$, $\Box(\speak{$f$ is invertible
in $\ul{R}$})$. The sublocale~$\Spec(A|R)$ is therefore the largest
sublocale of~$\Spec(A)$ on which the morphism~$\ul{R} \to \O_{\Spec(A)}$ is
local.
\end{rem}

\begin{rem}In classical logic, the sublocale~$\Spec(A|R)$ is closed
in~$\Spec(A)$, coinciding with~$V(\mmm_R A)$ (see
Remark~\ref{rem:local-spectrum-in-classical-logic}). But we don't think that
this property can be verified intuitionistically.\footnote{The nucleus
corresponding to the sublocale~$\Spec(A|R)$ is the quasicoherator~$\aaa \mapsto
\overline{\aaa}$. Assume that the sublocale is closed. Then there is a radial
ideal~$\bbb$ such that~$\overline{\aaa} = \sqrt{\aaa + \bbb}$ for all radical
ideals~$\aaa \subseteq A$. This radical ideal~$\bbb$ is uniquely determined
by~$\bbb = \overline{\sqrt{(0)}}$. Thus we obtain the simple
description~$\overline{\aaa} = \sqrt{\aaa + \overline{\sqrt{(0)}}}$ for the
quasicoherator. We don't believe that this is right in general.} For technical
reasons, it would be nice to know that~$\Spec(A|R)$ is an essential sublocale,
since the pullback functor admits a simpler description for essential
sublocales. It is an intersection of open, hence essential, sublocales, but the
intersection of essential sublocales needn't be
essential~\cite{kelly-lawvere:essential-localizations}. We are grateful
to Guilherme Frederico Lima de Carvalho e Silva for valuable discussions and
references on this matter.
\end{rem}

\begin{proof}[Proof of Theorem~\ref{thm:local-spectrum-yields-relative-spectrum}]
Follows immediately by interpreting the intuitionistic proof of
Lemma~\ref{lemma:universal-property-local-spectrum} in the internal language
of~$\Sh(X)$, applied to~$R \defeq \O_X$ and~$A \defeq \A$.
Then~``$(\pt,\O_X)$'' actually refers to the locally ringed locale~$(X,\O_X)$
and ``$\Gamma(Y,\O_Y)$'' refers to~$\mu_*\O_Y$, where~$\mu : (Y,\O_Y) \to
(X,\O_X)$ is a locally ringed locale over~$(X,\O_X)$.
\end{proof}

Theorem~\ref{thm:local-spectrum-yields-relative-spectrum} settles the question
how the little Zariski topos of~$\RelSpec_X(\A)$ looks like from the
internal point of view of~$\Sh(X)$. A related question is how the big Zariski
topos looks like. We give the answer in
Theorem~\ref{thm:big-zariski-topos-of-relative-spectrum}.

A basic fact about the ordinary spectrum is that the ring of global sections
of~$\O_{\Spec(A)}$ is canonically isomorphic to~$A$. This is not true for the
local spectrum: A trivial example is given by any algebra over a
ring~$R$ which is strongly not local over~$R$ in the sense that there are some nonunits
of~$R$ which generate the unit ideal in~$A$.
In this case the theory which~$\Spec(A|R)$ classifies is
inconsistent. Thus~$\Spec(A|R)$ is the empty locale and the ring of global
sections of~$\O_{\Spec(A|R)}$ is the zero ring.\footnote{This failure is not
entirely unexpected, since Coste's general result on sheaf
representations~\cite[Theorem~5.1.1]{coste:sheaf-representation}, which would
immediately guarantee that the global sections of~$\O_{\Spec(A|R)}$ are in
canonical one-to-one correspondence with the elements of~$A$, is not
applicable: The theory which~$\Spec(A|R)$ classifies is not a coherent
theory. The set-indexed disjunction appearing in axiom scheme~(5) of the
description given in Proposition~\ref{prop:local-spectrum-classify} can't be
rewritten as a finite disjunction.}

\begin{prop}\label{prop:local-spectrum-global-functions}
Let~$R$ be a local ring. Let~$A$ be an~$R$-algebra. The canonical
homomorphism~$A \to \Gamma(\Spec(A|R), \O_{\Spec(A|R)})$ is an isomorphism in
any of the following situations:
\begin{enumerate}
\item The algebra~$A$ satisfies the quasicoherence condition given in
Theorem~\ref{thm:qcoh-sheafchar}.
\item The algebra~$A$ is local and the structure morphism~$R \to A$ is local.
\end{enumerate}
Furthermore, the sheaf~$\O_{\Spec(A)}$ is a sheaf
for the modal operator associated to the sublocale~$\Spec(A|R) \hookrightarrow
\Spec(A)$ if and only if~$A$ satisfies the quasicoherence condition.
\end{prop}

\begin{proof}We only verify the claim in the second situation. In this
case~$1 \in \overline{\aaa}$ implies~$1 \in \aaa$ for any radical ideal of~$A$,
as can be checked using the proof scheme given in
Remark~\ref{rem:quasicoherator-knaster-tarski}.
Hence~$\Spec(A|R)$ is a local locale, meaning that for any covering~$\sqrt{(1)}
= \bigvee_i \aaa_i = \overline{\sqrt{\sum_i \aaa_i}}$ of the top element
of~$\Open(\Spec(A|R))$, there is an index~$i$ such that~$\aaa_i = \sqrt{(1)}$.

The locale~$\Spec(A|R)$ thus has an initial (locale-theoretic) point. This
focal point can be explicitly described: it is the filter~$A^\times$ (which
lies over the filter of units because~$R \to A$ is local). As is generally the
case for local locales, taking global sections is the same as taking the stalk
at the focal point. Therefore we can conclude by the following string of
isomorphisms.
\[ \Gamma(\Spec(A|R), \O_{\Spec(A|R)}) \cong
  \O_{\Spec(A|R), A^\times} \cong
  A[(A^\times)^{-1}] \cong
  A. \qedhere \]
\end{proof}

\begin{cor}\label{cor:pushforward-relative-spectrum}
Let~$X$ be a scheme. Let~$\A$ be an~$\O_X$-algebra.
Let~$f : \RelSpec_X(\A) \to X$ be the canonical projection morphism.
The canonical morphism~$\A \to f_*\O_{\RelSpec_X(\A)}$ of~$\O_X$-algebras is an
isomorphism in any of the following situations:
\begin{enumerate}
\item The algebra~$\A$ is quasicoherent.
\item From the point of view of~$\Sh(X)$, the algebra~$\A$ is local and the
homomorphism~$\O_X \to \A$ is local. (This means that, for every point~$x \in
X$, the stalk~$\A_x$ is local and the homomorphism~$\O_{X,x} \to \A_x$ is
local.)
\end{enumerate}
\end{cor}

\begin{proof}This is just the interpretation of
Proposition~\ref{prop:local-spectrum-global-functions} internal to~$\Sh(X)$.
\end{proof}

\begin{rem}Naively one might think that the canonical morphism~$\A \to
f_*\O_{\RelSpec_X(\A)}$ of Corollary~\ref{cor:pushforward-relative-spectrum} is
the canonical morphism from~$\A$ to the quasicoherization of~$\A$. This is not
the case. Firstly, the canonical morphism obtained by quasicoherization goes in
the other direction. Secondly, as stated in
Corollary~\ref{cor:pushforward-relative-spectrum}, the canonical morphism can
be an isomorphism even if~$\A$ is not quasicoherent.
\end{rem}

\subsection{Comparing the different spectrum constructions}

For rings and algebras, there are at least the following spectrum
constructions.

\begin{itemize}
\item The ordinary spectrum of a ring, possibly realized as a locale instead of
a topological space in order to work in an intuitionistic setting: $\Ring^\op \to
\LRS$ or $\Ring^\op \to \LRL$
\item The local spectrum of an algebra: $\Alg(R)^\op \to \LRL/(\pt,R)$
\item (A special case of) Gillam's spectrum of a sheaf of algebras~\cite{gillam:localization}:
$\Alg(\O_X)^\op \to \LRS/(X,\O_X)$
\item Hakim's spectrum of a ringed topos~\cite{hakim:relative-schemes},
yielding a locally ringed topos: $\RT \to \LRT$.
\item Cole's general framework for spectrum constructions~\cite{cole:spectra}
(also reported on at~\cite[Theorem~6.58]{johnstone:topos-theory})
\end{itemize}

These are related as follows.

As described in Section~\ref{sect:rel-spec-as-ordinary-spec}, the ordinary
spectrum construction cannot only be applied to rings, but also to sheaves of
rings and indeed ring objects internal to arbitrary elementary toposes equipped
with a natural numbers object, by employing the internal language. Applied to a
ring~$\O$ internal to such a topos~$\E$, it yields a locally ringed locale
internal to~$\E$, or equivalently a locally ringed localic topos internal
to~$\E$.  Externally, this corresponds to a locally ringed topos which is
equipped with a localic geometric morphism to~$\E$.

The ordinary spectrum construction can therefore be used to turn a ringed
topos~$(\E,\O)$ (with a natural numbers object) into a locally ringed topos
(which will be equipped with a morphism of ringed toposes to~$(\E,\O)$, but
which will, even if~$\O$ happens to be a local ring, not be equipped with a
morphism of locally ringed toposes to~$(\E,\O)$).

By comparing the universal properties one sees that this kind of internal
application of the ordinary spectrum construction coincides with the result of
Hakim's spectrum construction. In fact, it can be interpreted as a simultaneous
simplification and generalization of Hakim's construction: It's simpler, since
it's just the familiar spectrum construction and no explicit site calculations
are required; and it's more general, since Hakim's construction only applies to
ringed Grothendieck toposes whereas the internally-performed construction of
the ordinary spectrum applies to ringed elementary toposes with natural numbers
object.

Gillam's spectrum coincides with internally performing the construction of the
local spectrum, with the caveat that Gillam's construction starts with and
yields a locally ringed space, whereas ours starts with and yields a locally
ringed locale.\footnote{More generally, the local spectrum construction can be
applied to any algebra over a local ring~$\O$ internal to an elementary
topos~$\E$ with a natural numbers objects and yields a locally ringed topos
equipped with a morphism of locally ringed toposes to~$(\E,\O)$.} More
precisely:

For a locale~$Y$, let~$Y_P$ be the topological space of points of~$Y$, and for
a topological space~$T$, let~$T_L$ be the induced locale.  Let~$(X,\O_X)$ be a
sober locally ringed topological space. Let~$\A$ be an~$\O_X$-algebra. Then we
have a morphism~$E(\Spec(\A|\O_X)) \to X_L$ of locally ringed locales.
Since~$X \cong (X_L)_P$, there is an induced morphism~$E(\Spec(\A|\O_X))_P \to X$ of
locally ringed spaces. The adjunction~$(\placeholder)_L \dashv
(\placeholder)_P$ relating locales and topological spaces then yields, for any
locally ringed space~$\mu : Y \to X$ over~$X$,
\begin{align*}
  \Hom_{\LRS/X}(Y, E(\Spec(\A|\O_X))_P) &\cong
  \Hom_{\LRL/X_L}(Y_L, E(\Spec(\A|\O_X))) \\
  &\cong \Hom_{\Alg(\O_X)}(\A, \mu_*\O_Y).
\end{align*}
This is precisely the universal property which Gillam's spectrum enjoys.

Cole's framework for spectrum constructions is sufficiently general to
encompass both the ordinary spectrum and the local spectrum, and by extension
Hakim's spectrum and Gillam's spectrum. As is well-known, the ordinary spectrum
can be obtained from Cole's framework by applying it to the geometric
theory~$\mathbb{S}$ of rings, its quotient theory~$\mathbb{T}$ of local rings,
and the admissible class~$\mathbb{A}$ of local homomorphisms (notation as
in~\cite[Theorem~6.58]{johnstone:topos-theory}). The local spectrum can be
obtained by applying it to the geometric theory~$\mathbb{S}$
of~$\O_X$-algebras, its quotient theory~$\mathbb{T}$ of local~$\O_X$-algebras
which are local over~$\O_X$, and the admissible class of local homomorphisms.
For this to make sense, one has to interpret Cole's framework in the internal
language of~$\Sh(X)$, since there are no external geometric theories of
(local)~$\O_X$-algebras.

In general, the local spectrum doesn't coincide with the usual spectrum and
Gillam's spectrum doesn't coincide with Hakim's spectrum. (This was partly
observed before, for instance in~\cite[Remark~A.1.2]{chai:compactification}.) However, if the base
space is a scheme of dimension~$\leq 0$, they do coincide.

\begin{prop}\label{prop:local-spectrum-full-spectrum}
Let~$X$ be a scheme. Then~$E(\Spec(\O_X)) \cong X$ as locales
over~$X$ if and only if~$\dim X \leq 0$.\end{prop}
\begin{proof}The externalization of~$\Spec\O_X$ coincides with~$X$ if and only
if from the internal point of view, the locale~$\Spec\O_X$ coincides with the
one-point locale. By interpreting Lemma~\ref{lemma:spectrum-one-point} in the
internal language of~$\Sh(X)$, it follows that this is the case if and only if
\[ \Sh(X) \models \forall f\?\O_X\_ \speak{$f$ nilpotent} \vee \speak{$f$
invertible}. \]
(Internally, it always holds that~$\neg(1 = 0)$ in~$\O_X$, even if~$X$ happens
to be the empty scheme. Therefore the lemma is indeed applicable.) By
Corollary~\ref{cor:scheme-dimension-zero}, this condition is equivalent to the
dimension of~$X$ being less than or equal to zero (\ie to~$X$ being empty or
having dimension exactly zero).
\end{proof}

\begin{cor}Let~$X$ be a scheme. Then the relative spectrum of~$\O_X$-algebras
can be computed by the internal spectrum (instead of the internal local
spectrum) if and only if~$\dim X \leq 0$.\end{cor}
\begin{proof}The externalization of the internal spectrum of
arbitrary~$\O_X$-algebras~$\A$ coincides with the relative spectrum if and
only if it coincides in the special case~$\A = \O_X$. This is apparent by the
universal properties of both constructions. Thus the claim follows from
Proposition~\ref{prop:local-spectrum-full-spectrum}.
\end{proof}

Which construction is more fundamental, the ordinary spectrum of a ring or the
local spectrum of an algebra? The ordinary spectrum~$\Spec(A)$ can be expressed
as the local spectrum~$\Spec(A^\sim|\O_{\Spec(\ZZ)})$, where~$A^\sim$ is the
induced quasicoherent algebra on~$\Sh(\Spec(\ZZ))$. This fact is well-known
in the alternate form~``$\RelSpec_{\Spec(\ZZ)}(A^\sim) \cong \Spec(A)$''.

Fast and loose reasoning as follows could lead one to believe that it's
similarly possible to express the local spectrum as an ordinary spectrum.
Let~$R$ be a local ring. Let~$\varphi : R \to A$ be an algebra. The points
of~$\Spec(A|R)$ are those filters~$F \subseteq A$ such that~$\varphi^{-1}F =
R^\times$. Illicitly assuming classical logic, the points of~$\Spec(A|R)$ are in
canonical one-to-one correspondence with those prime ideals~$\ppp \subseteq A$
such that~$\varphi^{-1}\ppp = \mmm_R$. The points of~$\Spec(A|R)$ are therefore
in canonical one-to-one correspondence with the points of~$\Spec(A \otimes_R
k)$, where~$k = R/\mmm_R$ is the residue field of~$R$. Therefore~$\Spec(A|R)$
and~$\Spec(A \otimes_R k)$ might coincide.

However, we have the following negative result.\footnote{Intuitionistically,
it's still true that the prime ideals of a quotient ring~$A/\ppp$ are in
one-to-one correspondence with those prime ideals of~$A$ which contain~$\ppp$.
However, the analogous statement ``filters of~$A/F$ correspond to those filters
of~$A$ which are contained in~$F$'' can't be shown intuitionistically, if~$A/F$
is defined as~$A/F^c$. However, informally speaking, this failure is not the
fault of the statement, but of the definition of~$A/F$. The definition
raises red flags from an intuitionistic point of view, since not~$F$, but only
its complement~$F^c$ enters the construction.

The statement can be salvaged by defining~``$A/F$'' to mean the set~$A$
equipped with a new \emph{apartness relation} defined by $a \apart b
\vcentcolon\Leftrightarrow a-b \in F$. (A basic example for a
ring-with-apartness-relation is the field of real numbers equipped with~$x
\apart y \vcentcolon\Leftrightarrow \exists q \in \QQ\_ |x-y| \geq q > 0$.) A filter~$G$ of this
ring-with-apartness-relation~$A$ is by definition a subset~$G \subseteq A$ which
verifies the filter axioms and which is \emph{open with respect to the
apartness relation} in that for any elements~$a,b \? A$, the implication~$a \in
G \Rightarrow ((b \in G) \vee (a \apart b))$ holds.

This construction provides one of several motivations for developing the theory
of rings using apartness relations and anti-ideals; one can even define the
spectrum of a ring-with-apartness-relation. However, we'll not pursue these
ideas further here.}

\begin{prop}In general, the local spectrum of an algebra can't be expressed as
an ordinary spectrum.\end{prop}

\begin{proof}It is well-known that the ordinary spectrum is always quasicompact. The local spectrum,
however, can fail to be quasicompact. A quick way to see this is to notice
that, if that was the case, the locale-theoretic part of the projection
morphism~$\RelSpec_X(\A) \to X$ would always be a proper map of
locales~\cite{vermeulen:locales}.

There's also a more direct way of seeing this, which in fact proves a slightly
stronger statement. Let~$X$ be a scheme. Let~$f\in\Gamma(X,\O_X)$.
From the internal point of view of~$\Sh(X)$, the local
spectrum~$\Spec(\O_X[f^{-1}]|\O_X) \hookrightarrow \Spec(\O_X|\O_X) \cong \pt$
is the open sublocale of~$\pt$ corresponding to the truth value of~``$f$~is
invertible''.  Explicitly, the frame of opens of~$\Spec(\O_X[f^{-1}]|\O_X)$ is
isomorphic to $\{ \psi \? \Omega \,|\, \psi \Rightarrow \text{$f$ is
invertible} \}$.

The ordinary spectrum always has the Frobenius reciprocity property, being
quasicompact. In contrast, the locale~$\Spec(\O_X[f^{-1}]|\O_X)$ has this
property if and only if~$f$ is nilpotent or invertible.
\end{proof}

\begin{rem}\label{rem:local-spectrum-in-classical-logic}
Even in classical logic, where the local spectrum~$\Spec(A|R)$
and the ordinary spectrum~$\Spec(A \otimes_R R/\mmm_R)$ coincide as locales,
they do not coincide as locally ringed locales. The structure sheaf
of~$\Spec(A|R)$, regarded as a sheaf on~$\Spec(A \otimes_R R/\mmm_R)$,
is~$i^{-1}\O_{\Spec(A)}$, where~$i : \Spec(A \otimes_R R/\mmm_R)
\hookrightarrow \Spec(A)$ is the closed immersion corresponding to the
inclusion~$\Spec(A|R) \hookrightarrow \Spec(A)$. It's in general
not~$\O_{\Spec(A \otimes_R R/\mmm_R)}$.
\end{rem}

Finally, we want to restate the universal properties of the ordinary spectrum
and the local spectrum in ring-theoretic language, employing the dual
categories~$\RL^\op$ and~$\LRL^\op$, as in
Section~\ref{sect:flat-constructible-topologies}.

Let~$A$ be a ring. The morphism~$A \to \O_{\Spec(A)}$ in~$\RL^\op$ (the
ring-theoretic part of the canonical morphism~$(\Spec(A), \O_{\Spec(A)}) \to
(\Set, A)$) is the \emph{universal localization} of~$A$: The
ring~$\O_{\Spec(A)}$ is local, and for any morphism~$A \to \B$ into a local
ring~$\B$ (over any locale), there is a unique local morphism~$\O_{\Spec(A)}
\to \B$ rendering the diagram
\[ \xymatrix{
  A \ar[rd] \ar[rrr] &&& {\substack{\text{local}\\\text{\normalsize$\B$}\\\phantom{\text{local}}}} \\
  & {\substack{\text{\normalsize$\O_{\Spec(A)}$}\\\text{local}}} \ar@{-->}_[@!29]{\text{local}}[rru]
} \]
commutative. In contrast, the universal property of the local spectrum is as
follows. Let~$R$ be a ring. Let~$A$ be an~$R$-algebra. The morphism~$A \to
\O_{\Spec(A|R)}$ is the universal way of turning~$A$ into a local ring
\emph{which is local over~$R$}: The ring~$\O_{\Spec(A|R)}$ is local, the
composition~$R \to A \to \O_{\Spec(A|R)}$ is local, and for any morphism~$A
\to \B$ into a local ring (over any locale) such that the composition~$R \to A
\to \B$ is local, there is a unique local morphism~$\O_{\Spec(A|R)} \to \B$
such that the diagram
\[ \xymatrix{
  R \ar[r]\ar@/^2pc/[rrrr]^{\text{local}}\ar@/_/[rrd]_[@!-26]{\text{local}} &
    A \ar[rd] \ar[rrr] &&&
    {\substack{\text{local}\\\text{\normalsize$\B$}\\\phantom{\text{local}}}} \\
  && {\substack{\text{\normalsize$\O_{\Spec(A|R)}$}\\\text{local}}} \ar@{-->}[rru]_[@!28]{\text{local}}
} \]
commutes.

\begin{rem}It's possible to state the universal property of the structure sheaf
of the big Zariski topos of a ring~$A$, more precisely of the canonical
morphism~$(\Zar(A),\affla) \to (\Set,A)$ of ringed toposes, in a similar manner,
employing the dual categories~$\RT^\op$ and~$\LRT^\op$ of the categories of
(locally) ringed toposes. However, unlike the universal property of the
spectrum, this universal property looks slightly odd from an algebraic point of
view: For any morphism~$A \to \B$ into a local ring (over any topos~$\E$), there is
a unique bijective homomorphism~$\affla \to \B$ rendering the diagram
\[ \xymatrix{
  A \ar[rd] \ar[rrr] &&& {\substack{\text{local}\\\text{\normalsize$\B$}\\\phantom{\text{local}}}} \\
  & {\substack{\text{\normalsize$\affla$}\\\text{local}}} \ar@{-->}_[@!31]{\text{bijective}}[rru]
} \]
commutative. By ``bijective'' we mean that the ring-theoretic part~$f^\sharp :
f^{-1}\affla \to \B$ of the morphism~$f:(\E,\B) \to (\Zar(A),\affla)$ is
bijective as seen from the internal point of view of~$\E$.
\end{rem}

\subsection{The spectrum of the generic ring}

Let~$\Set[\Ring]$ be the classifying topos of the theory of
rings; explicitly, it's the topos of presheaves on~$\Ring_\fp^\op$, the dual of
the category of finitely presented rings. This topos contains the \emph{generic
ring}~$U$ (explicitly the presheaf~$R \mapsto R$): any ring in any topos is the
pullback of~$U$ along a suitable geometric morphism.

Let~$\Set[\LocRing]$ be the classifying topos of the theory of local rings.
Explicitly, it's the big Zariski topos~$\Zar(\Spec(\ZZ))$ (built using one of
the \emph{parsimonious sites}, as described in
Section~\ref{sect:proper-choice-of-site}). This topos contains the
\emph{generic local ring}~$U'$: any local ring in any topos is the pullback
of~$U'$ along a suitable geometric morphism.

Let~$A$ be a ring. By the universal property of~$\Set[\Ring]$, there is a
geometric morphism~$g : \Set \to \Set[\Ring]$ such that~$g^{-1}U \cong A$.
Since~$U'$ is in particular a ring, again by the universal property
of~$\Set[\Ring]$, there is a geometric morphism~$f : \Set[\LocRing] \to
\Set[\Ring]$ such that~$f^{-1}U \cong U'$.
By the universal property of~$\Set[\LocRing]$, the topos of
sheaves over the spectrum of~$A$ admits a geometric morphism~$g'$
to~$\Set[\LocRing]$ such that~$(g')^{-1}U' \cong \O_{\Spec(A)}$.

The resulting solid diagram
\[ \xymatrix{
  \E \ar@/^1pc/@{.>}[rrd]^{\tilde g} \ar@/_1pc/@{.>}[ddr]_{\tilde f} \ar@{.>}[rd]^h \\
  & \Sh(\Spec(A)) \ar[r]^{g'} \ar[d]_{f'} & \Set[\LocRing] \ar[d]^f \\
  & \Set \ar[r]_g \ar@{}[ur]^(.3){}="a"^(.7){}="b" \ar@{=>}_\eta "a";"b" & \Set[\Ring]
} \]
commutes up to a non-invertible natural transformation~$\eta$; under the
equivalence
\begin{multline*}
  \qquad\text{category of geometric morphisms~$\Sh(\Spec(A)) \to \Set[\Ring]$} \simeq \\
  \text{category of ring objects in~$\Sh(\Spec(A))$}\qquad
\end{multline*}
this transformation corresponds to the non-invertible localization homomorphism~$\ul{A} \to
\ul{A}[\F^{-1}] = \O_{\Spec(A)}$. It is folklore that this square is a lax
pullback square in the 2-category of Grothendieck toposes (for instance, this
is reported on at~\cite{arndt:lax-pullback}); however, this is not true.

Given a topos~$\E$ together with geometric morphisms~$\tilde f : \E \to \Set$
and~$\tilde g : \E \to \Set[\LocRing]$ and a natural transformation~$\tilde
\eta : \tilde f^{-1} \circ g^{-1} \Rightarrow \tilde g^{-1} \circ f^{-1}$
(these data correspond to a local ring~$\O_\E$ in~$\E$ together with a ring
homomorphism~$\varphi : \ul{A} \to \O_\E$), there is a canonical geometric morphism~$h :
\E \to \Sh(\Spec(A))$ (determined by requiring that~$h^{-1}\F \cong \F_0 \defeq
\varphi^{-1}[\O_\E^\times]$), and this morphism renders the lower left triangle
commutative up to a natural isomorphism, but it renders the upper right
triangle commutative only up to a non-invertible natural transformation
(corresponding to the non-invertible ring homomorphism~$\ul{A}[\F_0^{-1}] \to
\O_\E$).

The observation that the square is not a lax pullback is joint with Peter Arndt
and Matthias Hutzler. The observation raises two questions: What is the lax
pullback (which exists by general theory), if it's not~$\Sh(\Spec(A))$?
And how can~$\Sh(\Spec(A))$ be described as a pullback? The following two
propositions answer these questions. The geometric morphism~$\Set \to
\Set[\Ring]$ which they implicitly refer to is the morphism~$g$ mentioned
above.

\begin{prop}Let~$A$ be a ring. The lax pullback~$(\Set \Rightarrow_{\Set[\Ring]}
\Set[\LocRing])$ is the big Zariski topos of~$\Spec(A)$ (built using one of the
parsimonious sites, as described in Section~\ref{sect:proper-choice-of-site}).
\end{prop}

\begin{proof}The claim can be checked by hand, but it's more instructive to
employ the general theory of classifying toposes. In the situation
\[ \xymatrix{
  (\Set[T] \Rightarrow_{\Set[T_0]} \Set[T']) \ar[r] \ar[d] & \Set[T'] \ar[d]^f \\
  \Set[T] \ar[r]_g \ar@{}[ur]^(.3){}="a"^(.7){}="b" \ar@{=>}_\eta "a";"b" & \Set[T_0],
} \]
where~$T_0$, $T$, and~$T'$ are arbitrary geometric theories, the lax pullback
classifies the geometric theory whose models consist of a model~$M$ of~$T$, a
model~$N$ of~$T'$, and a homomorphism~$G(M) \to F(N)$ of~$T_0$-models. The
constructions~$G$ and~$F$ are given by the geometric morphisms~$g$ and~$f$ as
follows:

Any object of~$\Set[T]$ can be obtained by geometric constructions from~$U_T$,
the universal model of~$T$ in~$\Set[T]$. In particular, the pullback~$g^{-1}
U_{T_0}$, which is a model of~$T_0$, can be obtained by geometric constructions
from~$U_T$. Therefore the geometric morphism~$g$ displays a way to turn the
generic model of~$T$ into a model of~$T_0$ using only geometric constructions.
The same constructions can be applied to any model~$M$ of~$T$, yielding a
model~$G(M)$ of~$T_0$.

In the concrete situation at hand, the theory~$T$ is the empty theory
(admitting in any topos a unique model~$M$), the theory~$T'$ is the theory of local
rings, and~$T_0$ is the theory of rings. The~$T_0$-model~$G(M)$ is the
ring~$A$. The~$T_0$-model~$F(N)$ of a local ring~$N$ is the underlying ring
of~$N$.

Therefore the lax pullback~$(\Set \Rightarrow_{\Set[\Ring]} \Set[\LocRing])$
classifies ring homomorphisms~$A \to R$ where~$R$ is a local ring, that is,
local~$A$-algebras. It's well-known that~$\Zar(\Spec(A))$ classifies these as
well.
\end{proof}

\begin{rem}The pseudo pullback of the geometric morphism~$\Set[\LocRing] \to
\Set[\Ring]$ along~$\Set \to \Set[\Ring]$ is not very interesting: It is the
largest subtopos of~$\Set$ where the given ring~$A$ is local. Assuming
classical logic, this subtopos is either the trivial topos (if~$A$ is not
local) or~$\Set$ (if~$A$ is local).\end{rem}

\begin{prop}\label{prop:spectrum-as-pullback}
Let~$A$ be a ring. The pullback of the spectrum of the generic ring
along~$\Set \to \Set[\Ring]$ is the spectrum of~$A$.
\end{prop}

\begin{proof}There are two related ways of making the statement precise.
Firstly, the spectrum of the generic ring~$U$ can be interpreted as a (locally
ringed) locale internal to~$\Set[\Ring]$. Locales can be pulled back along
geometric morphisms (even though the pullback of a frame along a geometric
morphism typically fails to be a frame)~\cite{vickers:case-study}. In this
way~$\Spec(U)$ pulls back to a locale internal to~$\Set$, that is an ordinary
external locale. The claim is that this locale is canonically isomorphic
to~$\Spec(A)$.

A second way to interpret the statement of the proposition is to regard the
spectrum of the generic ring as a localic geometric morphism with
codomain~$\Set[\Ring]$. The claim is then that the diagram
\[ \xymatrix{
  \Sh(\Spec(A)) \ar[r] \ar[d] & \Sh_{\Set[\Ring]}(\Spec(U)) \ar[d] \\
  \Set \ar[r]_g & \Set[\Ring]
} \]
is a pullback diagram in the 2-category of toposes.

Using the language of classifying locales and classifying toposes, both claims
are easy to establish. The pulled-back locale (or topos) classifies the
pulled-back geometric theory~\cite[Corollary~5.4]{vickers:case-study}. Since
the description of the theory which~$\Spec(U)$ classifies -- the theory of
filters of~$U$ -- is itself geometric, the pulled-back theory is the theory of
filters of~$g^{-1}U \cong A$.\footnote{In the notation
of~\cite[Section~5]{vickers:case-study}, the theory of filters of~$U$ is
represented by a GRD system with~$G = U$ and $R = 1 \amalg U^2 \amalg U^2 \amalg 1
\amalg U^2$ (one summand for each axiom scheme).}
\end{proof}

\begin{prop}\label{prop:local-spectrum-generic}\begin{enumerate}
\item Let~$A$ be an~$R$-algebra. The local spectrum~$\Spec(A|R)$ is the
pullback of~$\Spec(U''|R)$ along the geometric
morphism~$\Set \to \E$ given by~$A$,
where~$U''$ is the \emph{generic~$R$-algebra}
contained in the classifying topos~$\E$ of~$R$-algebras.
\item Let~$X$ be a scheme (or a locally ringed locale). Let~$\A$ be
an~$\O_X$-algebra. The relative spectrum~$\RelSpec_X(\A)$ is the pullback
of~$\Spec(U''|\O_X)$ along the geometric
morphism~$\Sh(X) \to \E$ given by~$\A$,
where~$U''$ is the generic~$\O_X$-algebra contained in
the classifying $\Sh(X)$-topos~$\E$ of~$\O_X$-algebras.
\end{enumerate}
\end{prop}

\begin{proof}Straightforward modification of the proof of
Proposition~\ref{prop:spectrum-as-pullback}.
\end{proof}

\begin{rem}The big Zariski topos~$\Zar(\Spec(A))$ can be obtained as the
pullback of the big Zariski topos of the generic ring~$U$, if both toposes
are understood to be defined using the parsimonious sites as described in
Section~\ref{sect:proper-choice-of-site}.
\end{rem}

\subsection{Limits in the category of locally ringed locales}

The category of ringed locales has small limits, by the naive construction. For
instance, the fiber product~$X \times_Z Y$ of ringed locales is given by the
fiber product of the underlying locales and the
structure sheaf~$\pi_X^{-1}\O_X \otimes_{\pi_Z^{-1}\O_Z} \pi_Y^{-1}\O_Y$. More
generally, the limit of a small diagram of ringed locales is given by the limit~$L$
of the underlying locales and the colimit of the pulled-back structure sheaves
(computed in the category of sheaves of rings on~$L$).

However, when applied to a diagram of locally ringed locales, the ringed locale
which this simple construction yields is in general not locally ringed. This
can be nicely understood from the internal point of view: Let~$R$ be a local
ring. Let~$R \to A$ and~$R \to B$ be local~$R$-algebras which are furthermore
local over~$R$. Then the tensor product~$A \otimes_R B$ is in general not a
local ring. Indeed, this fails even in the easiest case, where all rings
involved are fields: The rings~$\RR$ and~$\CC$ are local, and the
inclusion~$\RR \to \CC$ is local, but~$\CC \otimes_\RR \CC \cong \CC
\otimes_\RR \RR[X]/(X^2+1) \cong \CC[X]/(X^2+1) \cong \CC \times
\CC$ is not.

The following proposition explains that the true limit in the category of
locally ringed locales is obtained by \emph{relocalizing} the limit in the
category of ringed locales.

\begin{prop}\label{prop:lrl-complete}
The category of locally ringed locales has all small limits.
\end{prop}

\begin{proof}For notational simplicity, we describe how products in the
category of locally ringed locales can be constructed. The general case is
entirely analogous.

Let~$X$ and~$Y$ be locally ringed locales. Their product~$P$ as ringed locales has
two defects: Firstly, it's not locally ringed. Secondly, the ring-theoretic
parts of the projection morphisms~$\pi_X : P \to X$ and~$\pi_Y : P \to Y$
aren't local, that is, don't reflect invertibility.

The first issue could be solved by constructing, internally to~$\Sh(P)$, the
ordinary spectrum of~$\O_P$. From the external point of view, this would yield
a locally ringed locale equipped with morphisms of ringed, but not of locally
ringed, locales to~$X$ and~$Y$.

To solve both issues, we need to employ a refined spectrum construction,
similar to the modification required by the internal account of the relative
spectrum: Internally to~$\Sh(P)$, we construct the classifying locale of the
theory of those filters of~$\O_P$ which simultaneously lie over the filter of
units of~$\pi_X^{-1}\O_X$ and which lie over the filter of units
of~$\pi_Y^{-1}\O_Y$. This locale is a sublocale of~$\Spec(\O_P)$, the largest
such that the morphisms to~$(\pt,\pi_X^{-1}\O_X)$ and to~$(\pt,\pi_Y^{-1}\O_Y)$ are
morphisms of locally ringed locales.

The externalization of the internal locally ringed locale obtained in this way
is the sought product of~$X$ and~$Y$ in the category of locally ringed locales.
\end{proof}

\begin{rem}The category of locally ringed locales embeds as a (non-full)
coreflective subcategory into the category of ringed locales; the coreflector
maps a ringed locale~$(X,\O_X)$ to the externalization of~$\Spec(\O_X)$
(constructed internally to~$\Sh(X)$). However, as is familiar in situations where the embedding
is not full~\cite{adamek:rosicky:reflective}, it's in general not the case that limits in~$\LRL$ are computed by
applying the coreflector to the limit computed in~$\RL$. Employing the language
of the proof of Proposition~\ref{prop:lrl-complete}, applying the coreflector only
solves the first issue, but not the second.
\end{rem}

It's instructive to determine the points of limits in~$\LRL$, even though a
locale is of course not determined by its points. For instance, the construction in
Proposition~\ref{prop:lrl-complete} shows that the points of the product~$X
\times Y$ of locally ringed locales in~$\LRL$ are in canonical one-to-one
correspondence with tuples~$(x,y,F)$, where~$x$ is a point of~$X$,~$y$ is a
point of~$y$, and~$F$ is a filter of~$\O_{X,x} \otimes_\ZZ \O_{Y,y}$ which lies
over the filter of units of~$\O_{X,x}$ and of~$\O_{Y,y}$. In classical logic,
those tuples are in canonical one-to-one correspondence with
tuples~$(x,y,\ppp)$, where~$x$ and~$y$ are as before and~$\ppp$ is a prime
ideal of~$k(x) \otimes_\ZZ k(y)$.

Similarly, points of the fiber product~$X \times_Z Y$ are in canonical
one-to-one correspondence with tuples~$(x,y,F)$, where~$x$ is a point of~$X$
and $y$ is a point of~$y$ such that both map to the same point~$z$ of~$Z$,
and~$F$ is a filter of~$\O_{X,x} \otimes_{\O_{Z,z}} \O_{Y,y}$ lying over the
filter of units of~$\O_{X,x}$ and of~$\O_{Y,y}$ (and therefore automatically
of~$\O_{Z,z}$). In classical logic, those tuples are in canonical one-to-one
correspondence with tuples~$(x,y,\ppp)$, where~$x$ and~$y$ are as before
and~$\ppp$ is a prime ideal of~$k(x) \otimes_{k(z)} k(y)$.

\begin{rem}By the adjunction~$(\placeholder)_L \dashv (\placeholder)_P$
relating locales and topological spaces, limits of locally ringed spaces which
happen to be sober can be computed by regarding them as locally ringed
locales by~$(\placeholder)_L$, computing their limit in~$\LRL$, and taking
the associated topological space of the limit by~$(\placeholder)_P$.

Small diagrams of arbitrary locally ringed spaces admit limits as well.
Indeed, the proof of Proposition~\ref{prop:lrl-complete} was adapted from
Gillam's proof of this fact~\cite[Corollary~5]{gillam:localization}.\end{rem}

\subsection{Relative Proj construction} Similar issues as with the relative
spectrum arise with the Proj construction: The standard definition of the Proj
construction as a topological space of homogeneous prime ideals gives rise to a
space which can't intuitionistically be shown to satisfy the expected
universal property. The construction has to be reimagined as a locale
instead of a topological space. A certain sublocale of this locale then yields
the relative Proj construction when interpreted in the internal language of
the little Zariski topos of a base scheme (or a locally ringed locale).

\begin{defn}The \emph{Proj construction} of an~$\NN$-graded ring~$S$ is the
locale with frame of opens given by
\begin{multline*}
  \qquad\Open(\Proj(S)) \defeq
    \{ \aaa \subseteq S \,|\,
      \text{$\aaa$ is a homogeneous radical ideal such that} \\
  \forall x\?S\_
    x S_+ \subseteq \aaa \Rightarrow x \in \aaa \},\qquad
\end{multline*}
where~$S_+ = \bigoplus_{i > 0} S_i$ is the irrelevant ideal.
\end{defn}

A quick way to see that the partial order~$\Open(\Proj(S))$ is a frame is to
recognize that it's the frame of opens of a sublocale of~$\Spec(S)$.
The associated nucleus~$j : \Open(\Spec(S)) \to \Open(\Spec(S))$ is given by
\[ j(\aaa) \defeq
  (\sqrt{\aaa^h} : S_+), \]
where~$\aaa^h$ is the homogenization of~$\aaa$, the ideal of~$S$ generated by
all homogeneous components of the elements of~$\aaa$. Since~$\aaa \subseteq
\aaa^h \subseteq \sqrt{\aaa^h} \subseteq j(\aaa)$, a radical ideal~$\aaa$
is an element of~$\Open(\Proj(S))$ if and only if~$\aaa = j(\aaa)$.

One way to derive this definition is to start, within a classical context, with
the general expression for the nucleus associated to the subspace of~$\Spec(S)$
consisting of those prime ideals which are homogeneous and don't
contain~$S_+$, and then rewrite this expression to not refer to prime ideals.

\begin{defn}A filter~$F \subseteq S$ in an~$\NN$-graded ring~$S$ is
\emph{homogeneous} if and only if, for any element~$a \? S$, the filter~$F$
contains~$a$ if it contains at least one of the homogeneous components of~$a$.
It \emph{meets the irrelevant ideal} if and only if~$F \cap S_+$ is
inhabited.\end{defn}

In classical logic, a subset is a homogeneous filter meeting the irrelevant
ideal if and only if its complement is a homogeneous prime ideal not containing
the irrelevant ideal. Intuitionistically, neither direction can be shown.

\begin{prop}\label{prop:proj-classifying-locale}
Let~$S$ be an~$\NN$-graded ring. Then~$\Proj(S)$ is the classifying
locale of any of the following geometric theories.
\begin{enumerate}
\item The theory of homogeneous filters of~$S$ meeting the irrelevant ideal,
that is the theory of Remark~\ref{rem:theory-of-filters} supplemented by the
following two axiom schemes:
\begin{itemize}
\item $\bigvee_i (a_i \in F) \vdash a \in F$ (one axiom for each
decomposition~$a = \sum_i a_i$ of an element of~$S$ into homogeneous components)
\item $\top \vdash \bigvee_{a \in S_+} (a \in F)$ (one axiom)
\end{itemize}
\item The theory given by one atomic proposition~``$a \in F_i$'' for each
homogeneous element~$a$ of degree~$i$ in~$S$ and axioms given by the following
axiom schemes:
\begin{itemize}
\item $\top \vdash 1 \in F_0$ (one axiom)
\item $st \in F_{i+j} \dashv\vdash s \in F_i \wedge t \in F_j$ (two axioms for
each $i, j \geq 0$, $s \in S_i$, $t \in S_j$)
\item $0 \in F_i \vdash \bot$ (one axiom for each~$i \geq 0$)
\item $s+t \in F_i \vdash s \in F_i \vee t \in F_i$ (one axiom for each $i \geq
0$, $s,t \in A_i$)
\item $\top \vdash \bigvee_{i \geq 1} \bigvee_{a \in S_i} (a \in F_i)$ (one
axiom)
\end{itemize}
\item The same theory as in~(2), but with atomic propositions only for
homogeneous elements of degree~$\geq 1$ and without the first axiom~``$\top
\vdash 1 \in F_0$''.
\end{enumerate}
\end{prop}

\begin{proof}That~$\Proj(S)$ coincides with the classifying locale of the
theory given in~(1), can be verified by a direct calculation. By the general
theory, the nucleus associated to the quotient theory given in~(1) maps a
radical ideal~$\aaa \? \Open(\Spec(S))$ to the least fixed point above
of~$\aaa$ of the map
\[ \bbb \longmapsto
  \bbb \vee
  \bigvee_{a \? S} \Bigl(\sqrt{(a_i)_i} \cap \bigl(\sqrt{(a)} \rightarrow \bbb\bigr)\Bigr) \vee
  \Bigl(\sqrt{(a)_{a \in S_+}} \to \bbb\Bigr), \]
where~$(\ccc \to \bbb) = (\bbb : \ccc)$ is the Heyting implication
and~``$\vee$'' is the join in~$\Open(\Spec(S))$. We omit the intermediate steps
of the calculation.

The theories given in~(1) and in~(2) are bi-interpretable. The interpretation
of the atomic propositions~``$a \in F_i$'' of theory~(2) using the signature of
theory~(1) is~``$a \in F$''. Verifying the axioms is straightforward.
Conversely, the interpretation of~``$a \in F$'' in the signature of theory~(2)
is~``$\bigvee_i (a_i \in F_i)$'', where~$a = \sum_i a_i$ is the decomposition
into homogeneous components. For verifying the axioms, one needs the lemma that
\[ \bigvee_i (s_i \in F_i)  \wedge  \bigvee_j (t_j \in F_j)
    \ \dashv\vdash\ \bigvee_n \Bigl(\sum_{i+j=n} s_i t_j \in F_n\Bigr) \]
is derivable in theory~(2), for any decompositions~$s = \sum_i s_i$ and~$t =
\sum_j t_j$ of elements of~$S$ into homogeneous components. In the
guise~``$\sqrt{(s_i)_i} \cap \sqrt{(t_j)_j} = \sqrt{(\sum_{i+j=n} s_i t_j)_n}$''
this is a familiar fact on the content of
polynomials~\cite[Proposition~1]{banaschewski:vermeulen:radical-content}.

Also theories~(2) and~(3) are bi-interpretable. The interpretation of~``$a \in
F_0$'' in the signature of theory~(3) is~``$\bigvee_{i \geq 1} \bigvee_{h \in
S_i} (ha \in F_i)$''.
\end{proof}

\begin{cor}Let~$S$ be an~$\NN$-graded ring. The points of~$\Proj(S)$ are in
canonical one-to-one correspondence with the homogeneous filters
of~$S$ meeting the irrelevant ideal.\end{cor}

\begin{proof}Points of~$\Proj(S)$ are given by models of the theory of
homogeneous filters of~$S$ meeting the irrelevant ideal in~$\Set$.
\end{proof}

\begin{rem}The same presentation as in
Proposition~\ref{prop:proj-classifying-locale}(3) has been used to
construct~$\Proj(S)$ not as a locale, but as a distributive
lattice~\cite{cls:projective-spectrum}.\end{rem}

\begin{defn}Let~$S$ be an~$\NN$-graded ring. The \emph{generic homogeneous
filter meeting the irrelevant ideal} is the subsheaf~$\F \hookrightarrow
\ul{S}$ over~$\Proj(S)$ generated by the sections~$a$ over~$D_+(a) \defeq j(\sqrt{(a)})$.
\end{defn}

Equivalently, the generic homogeneous filter meeting the irrelevant ideal is
the pullback of the generic filter in~$\Sh(\Spec(S))$ to~$\Sh(\Proj(S))$.

\begin{defn}Let~$S$ be an~$\NN$-graded ring. The structure sheaf of~$\Proj(S)$ is the homogeneous
localization~$\ul{S}[\F^{-1}]_0$ of the ring~$\ul{S}$ at the generic
homogeneous filter meeting the irrelevant ideal, that is the degree-zero part
of~$\ul{S}[\F^{-1}]$. The \emph{tilde construction} of a graded~$S$-module~$M$
is~$M^\sim \defeq \ul{M}[\F^{-1}]_0$.
\end{defn}

The locally ringed locale~$\Proj(S)$ and the tilde construction defined in this
way enjoy their familiar properties. For instance, we have the following lemma.

\begin{lemma}Let~$S$ be an~$\NN$-graded ring.
\begin{enumerate}
\item Let~$f \? S$ be homogeneous of degree~$d \geq 1$. Then~$D_+(h) \cong
\Spec(S[f^{-1}]_0)$.
\item Assume that~$S$ is generated as an~$S_0$-algebra by~$S_1$. Let~$M$
and~$N$ be graded~$S$-modules. Then~$M^\sim \otimes_{\O_{\Proj(S)}} N^\sim
\cong (M \otimes_S N)^\sim$.
\item Under the same assumption as in~(2), the twisting sheaves~$\O(m) \defeq
(S(m))^\sim$ are finite locally free of rank~$1$.
\end{enumerate}
\end{lemma}

\begin{proof}For the first statement, it suffices to verify that the theories
of homogeneous filters of~$S$ meeting the irrelevant ideal and containing~$h$
and of filters of~$S[f^{-1}]_0$ are bi-interpretable. It's slightly more
convenient to use the presentation given by
Proposition~\ref{prop:proj-classifying-locale}(2) for the former theory.

The interpretation of~``$q \in F$'' for~$q \? S[f^{-1}]_0$ in the signature of
the theory given by Proposition~\ref{prop:proj-classifying-locale}(2) is
\[ \bigvee\{ (x \in F_{di}) \,|\, \text{$q = x/f^i$ for some $x \? S$, $i \geq
0$} \}. \]
Conversely, the interpretation of~``$a \in F_i$'' in the signature of the
theory of filters of~$S[f^{-1}]_0$ is~``$x^d / h^i \in F$''.

The second statement follows from the calculation
\begin{align*}
  M^\sim \otimes_{\O_{\Proj(S)}} N^\sim &=
  \ul{M}[\F^{-1}]_0 \otimes_{\ul{S}[\F^{-1}]_0} \ul{N}[\F^{-1}]_0 \\
  &\cong (\ul{M} \otimes_{\ul{S}} \ul{N})[\F^{-1}]_0
  \cong (\ul{M \otimes_S N})[\F^{-1}]_0
  = (M \otimes_S N)^\sim.
\end{align*}
The first isomorphism maps~$x/s \otimes y/t$ to~$(x \otimes y)/(st)$. By the
assumption that~$S$ is generated as an~$S_0$-algebra by~$S_1$, the generic
filter contains a homogeneous element~$h$ of degree~$1$ from the internal point
of view of~$\Sh(\Proj(S))$. Therefore the
map has an inverse sending~$(a \otimes b) / u$, where~$a$ and~$b$ are
homogeneous of degrees~$i$ and~$j$, to~$(h^j a)/u \otimes b/h^j$. The second
isomorphism is because the tensor product is a geometric construction and
therefore commutes with constructing the constant sheaf.

For the proof of the third statement, we show that~$(S(m))^\sim$ is a finite
free module of rank~$1$ from the internal point of view. We again use that the
generic filter contains a homogeneous element~$h \? \ul{S}$ of degree~$1$ from
the internal point of view. Such an element allows to define an
isomorphism~$\O_{\Proj(S)} = \ul{S}[\F^{-1}]_0 \to \ul{S(m)}[\F^{-1}]_0 = \O(m)$ by
mapping~$x/s$ to~$(h^m x)/s$ if~$m \geq 0$ and to~$x/(h^{-m} s)$ otherwise.
\end{proof}

\begin{defn}Let~$R$ be a ring. Let~$S$ be an~$\NN$-graded~$R$-algebra. The
\emph{local Proj construction} of~$S$ over~$R$ is the sublocale~$\Proj(S|R)$
of~$\Proj(S)$ with frame of opens given by
\[
  \Open(\Proj(S|R)) \defeq
    \{ \aaa \? \Open(\Proj(S)) \,|\,
      \forall f\?R\_ \forall s\?S\_
        (\speak{$f$ \inv} \Rightarrow s \in \aaa) \Rightarrow fs \in \aaa \}
\]
and with the pullback of~$\O_{\Proj(S)}$ as the structure sheaf.
\end{defn}

\begin{prop}\label{prop:local-proj-classifying-locale}
Let~$R$ be a ring. Let~$S$ be an~$\NN$-graded~$R$-algebra.
Then the local Proj construction~$\Proj(S|R)$ is the classifying locale of the theory of homogeneous
filters of~$S$ meeting the irrelevant ideal and lying over the filter of
units.\end{prop}

\begin{proof}Direct calculation similar to the proof of
Proposition~\ref{prop:proj-classifying-locale}.\end{proof}

Since pullback and localization commute, the structure sheaf of~$\Proj(S|R)$
can also be described as~$\ul{S}[\F^{-1}]_0$, where by abuse of notation we
mean by~``$\F$'' the pullback of the generic filter on~$\Proj(S)$
to~$\Proj(S|R)$. This filter has the special property
\[ \Sh(\Proj(S|R)) \models
  \forall r\?\ul{R}\_
    r \in \F \Rightarrow \speak{$r$ \inv\@ in $\ul{R}$}. \]

\begin{thm}\label{thm:local-proj-yields-relative-proj}
Let~$X$ be a scheme (or a locally ringed locale). Let~$\S$ be
an~$\NN$-graded $\O_X$-algebra. Then the externalization~$E(\Proj(\S|\O_X))$ coincides
with the relative Proj construction~$\RelProj_X(\S)$ as locally ringed locales
over~$X$.\end{thm}

\begin{proof}For simplicity, we assume that~$\S$ is generated as
an~$\S_0$-algebra by~$\S_1$. In this case, the expected universal property
of the relative Proj construction is that it's a locally ringed locale over~$X$ such
that, for all locally ringed locales~$\mu : Y \to X$ over~$X$, the
set~$\Hom_{\LRL/X}(Y, \RelProj_X(\S))$ is canonically isomorphic (by pullback
of the standard such datum on~$\RelProj_X(\S)$) to the set of pairs~$(\L,
\psi)$ such that
\begin{itemize}
\item $\L$ is a line bundle on~$Y$ and
\item $\psi : \mu^*\S \to \bigoplus_{n\geq0} \L^{\otimes n}$ is a graded
morphism of~$\O_Y$-algebras such that the degree-$1$ part of~$\psi$ is a
surjective morphism~$\mu^*\S_1 \to \L$
\end{itemize}
modulo equivalence. For instance, it is known that this property is satisfied
in the case that~$X$ is a scheme and~$\S$ is
quasicoherent~\stacksproject{01O4}.

We verify that~$E(\Proj(\S|\O_X))$ enjoys the same property, even if~$X$ is not
a scheme or~$\S$ is not quasicoherent. For the rest of the proof, we switch to
the internal universe of~$\Sh(X)$.

The local Proj construction is a locally ringed locale over~$(\pt, \O_X)$ by
the unique morphism~$! : \Proj(\S|\O_X) \to \pt$ of locales and by the canonical
morphism~$!^\sharp : \ul{\O_X} \to \ul{\S}_0 \to \ul{S}[\F^{-1}]_0 = \O_{\Proj(\S|\O_X)}$
of local rings.

As the standard datum on~$\Proj(\S|\O_X)$, we choose the line bundle~$\O(1)$
(pulled back to~$\Proj(\S|\O_X)$) together with the canonical morphism~$!^* \S
\to \oplus_{n \geq 0} \O(1)^{\otimes n}$.

Let~$Y$ be a locally ringed locale over~$(\pt, \O_X)$. Let a pair~$(\L,\psi)$
be given. In the internal language of~$\Sh(Y)$, we define a filter by the formula
\[ \F' \defeq \{ s \? \ul{S} \,|\,
  \speak{there exists~$i$ such that~$(\psi(s_i \otimes 1))$ is a basis
  of~$\L^{\otimes i}$} \} \subseteq \ul{S}, \]
where~$s_i$ refers to the homogeneous component of~$s$ of degree~$i$.
Since~$\L^{\otimes i}$ is finite free of rank~$1$, a one-element family
in~$\L^{\otimes i}$ is a basis if and only if it's a generating family. This
observation can be repeatedly used to verify that~$\F'$ is
homogeneous, meets the irrelevant ideal, and lies over the filter of
units. Since~$\Proj(\S|\O_X)$ is the classifying locale of such filters
(Proposition~\ref{prop:local-proj-classifying-locale}), we obtain a
morphism~$f : Y \to \Proj(\S|\O_X)$ of locales which is unique with the property
that~$f^{-1}\F = \F'$.

To obtain a morphism~$Y \to \Proj(\S|\O_X)$ of locally ringed locales, it
remains to define a morphism~$f^\sharp : f^{-1}\O_{\Proj(\S|\O_X)} =
\ul{S}[\F'^{-1}]_0 \to \O_Y$. A canonical choice is
\[ x/s \mapsto \speak{
  the coefficient of~$\psi(x \otimes 1)$
  with respect to the basis~$(\psi(s \otimes 1))$}. \]
We omit further verifications.
\end{proof}

\subsection{A constructive account of scheme theory}
\label{sect:constructive-scheme-theory}

Scheme theory as classically set up heavily relies on prime ideals and
therefore only works well in a classical context, where the law of excluded
middle and (at least some forms of) the axiom of choice are available. However,
the actual mathematical ideas often do not fundamentally require classical logic; we
don't begin the proof that the kernel of a morphism between quasicoherent
sheaves of modules is quasicoherent by supposing that it's not. Instead,
classical logic is only needed because the usual foundations of scheme theory
involving locally ringed spaces require it.

In this section, we sketch how scheme theory can be developed in an
intuitionistic metatheory; there are several reasons why it's desirable to have
such an account. Firstly, as is familiar from constructive treatments of other
subjects, the constraint to set up all definitions in an intuitionistically
sensible way is a useful guiding principle which can increase the perceived
elegance of the theory and result in more direct proofs.

It would be interesting to know which advanced results in algebraic geometry
\emph{actually} require classical logic (or at least classicality hypotheses on
the ground ring); to study this question, one has to use a foundation
which doesn't itself require classical logic just for organizational
purposes.\footnote{For instance, some results in linear algebra can
intuitionistically only be shown for \emph{discrete} fields -- fields such that
any element is zero or not zero. Such hypotheses are computationally meaningful
and will entail similar hypotheses for some results in algebraic geometry.}
McLarty and other researchers study a similar question: Which axioms of set
theory are actually needed for algebraic geometry, in particular for proving
Fermat's Last Theorem?~\cite{mclarty:what-does-it-take}

Secondly, one might be interested in concrete computations and might therefore
leverage the fact that one can mechanically extract algorithms from
constructive proofs. For instance, an intuitionistic proof that some cohomology
is finite dimensional yields an algorithm for computing the dimension and even
a basis.

Finally, one might want to apply scheme theory in the intuitionistic internal
universe of the little Zariski topos of a base scheme, in order to generalize
results of absolute scheme theory to relative scheme theory with little effort
and no duplication of proofs. The starting point for such a transfer is that
locally ringed locales over a locally ringed locale~$X$ look like locally
ringed locales over the point from the internal point of view of~$\Sh(X)$, as
discussed in Section~\ref{sect:internal-locales}.

The internal language of the big Zariski topos, presented in
Part~\ref{part:big-zariski}, is too a vehicle for relative scheme theory;
however, its language looks quite different from what one is accustomed to.

In this section, we only sketch how the basics of constructive scheme theory
could look like. Some parts are certainly folklore among constructive
mathematicians, but to the best of our knowledge no coherent summary appeared
in print before.

There is a vast literature on algorithmic computations in algebraic geometry
(to exemplarily cite just two references, Eisenbud's textbook on
syzygies~\cite{eisenbud:syzygies} and the GAP project~\cite{gap} are
well-known). However, these results are often still set in a classical context,
relying on classical logic for termination or correctness proofs. They
therefore don't contain an intuitionistic development of scheme theory.

{\tocless

\subsection*{Constructive algebra} Any constructive development of scheme
theory needs to rest on a constructive development of commutative algebra.
Such an account is readily
available~\cite{mines-richman-ruitenburg:constructive-algebra,lombardi:quitte:constructive-algebra}.

\subsection*{Local models} As discussed in
Section~\ref{sect:spectrum-as-a-locale}, defining the spectrum of a ring as a
topological space isn't sensible from a constructive point of view. A working
alternative is defining the spectrum as a locally ringed locale, employing the
frame of radical ideals. By considering sheaves over it, this yields a locally
ringed topos; this topos can also be presented by a more parsimonious site,
namely the site whose objects are the elements of the ring and whose coverings
are those finite families~$(g_i \to f)_i$ such that~$\sqrt{(f)} =
\sqrt{(g_i)_i}$.
This construction is due to
Joyal~\cite{joyal:spectrum,espanol:spectrum,tierney:spectrum} and was further explored
by several researchers~\cite{cls:spectral-schemes,cls:projective-spectrum}.
It's also possible to employ the framework of formal
topology~\cite{schuster:formal-zariski}.

The universal property of the localic spectrum ensures that morphisms~$\Spec(B)
\to \Spec(A)$ of locally ringed locales (or locally ringed toposes) are in
canonical one-to-one correspondence with ring homomorphisms~$A \to B$, as it
should be. There are two ways for explicitly constructing a morphism between
spectra. One is to specify a morphism of frames going in the other direction.
For instance, given a ring homomorphism~$\varphi : A \to B$, one can map a
radical ideal~$\aaa \subseteq A$ to~$\sqrt{\aaa B} \subseteq B$; this yields a
morphism~$\Open(\Spec(A)) \to \Open(\Spec(B))$.

Using the device of classifying locales, there is also another way which more
closely mimics the classical approach of taking preimages of prime ideals. To
give a morphism~$\Spec(B) \to \Spec(A)$ of locales amounts to give a model of
the theory of filters of~$A$ in~$\Sh(\Spec(B))$. The sheaf topos
over~$\Spec(B)$ contains the generic filter~$\F$ of~$B$; given a ring
homomorphism~$\varphi : A \to B$, this filter can be turned into a filter
of~$A$ by taking the preimage~$\varphi^{-1}[\F]$.

Classically, the induced map~$\Spec(B) \to \Spec(A)$ would be described
by~$\ppp \mapsto \varphi^{-1}[\ppp]$, where~$\ppp$ ranges over all prime ideals
of~$B$, and after defining the map in this way, one would have to verify its
continuity; constructively, we can describe it as~$\F \mapsto \varphi^{-1}[\F]$,
where~$\F$ is just a single special filter, and get continuity for free.
For more on this way of pretending that morphisms between locales are just maps
between points, we highly recommend an expository survey by Vickers on this
topic~\cite{vickers:continuity}.

As we have seen in Section~\ref{sect:transfer-principles}, for deriving
transfer principles it's useful to be able to quickly gauge properties of
constant sheaves over~$\Spec(A)$. For topological spaces,
Lemma~\ref{lemma:properties-of-constant-sheaves} could be used to this effect.
For locales in an intuitionistic metatheory, the lemma has to be modified
slightly.

\begin{defn}\begin{enumerate}
\item A locale~$X$ is \emph{overt} if and only if the unique morphism~$X
\to \pt$ of locales is open.
\item A \emph{positivity predicate} on a frame~$P$ is a predicate on the set of
elements of~$P$, written~``$U > 0$'' for~$U \? P$, such that for any element~$U
\? P$ and any subset~$M \subseteq P$,
\begin{itemize}
\item if~$U > 0$ and~$U
\preceq \bigvee M$, then there exists an element~$V \in
M$ such that~$V > 0$, and
\item if~$U > 0 \Longrightarrow U \preceq \bigvee M$, then~$U \preceq \bigvee M$.
\end{itemize}
\end{enumerate}
\end{defn}

\begin{ex}The frame of open subsets of a topological space~$X$ has a positivity
predicate, given by declaring~$U > 0$ if and only if~$U$ is inhabited.\end{ex}

\begin{ex}Assuming classical logic, any frame admits the positivity predicate
given by declaring~$U > 0$ if and only if~$U \neq \bot$.\end{ex}

\begin{prop}A locale~$X$ is overt if and only if its frame of opens admits a
\emph{positivity predicate}.\end{prop}

\begin{proof}Instructive unraveling of the definitions.
\end{proof}

\begin{lemma}\label{lemma:properties-of-constant-sheaves-over-locales}
Let~$\varphi$ be a first-order formula in which arbitrary sets and elements may
occur as parameters. Let~$X$ be a locale and let~$U$ be an open of~$X$.
Consider the following statements:
\begin{enumerate}
\item $U \models \varphi$ (with the same abuse of notation as in
Lemma~\ref{lemma:properties-of-constant-sheaves}).
\item $U \preceq \bigvee \{ \top \,|\, \varphi \}$.
\item (If~$X$ has a positivity predicate.) $U > 0 \Longrightarrow \varphi$.
\item $\varphi$.
\end{enumerate}
Then:
\begin{itemize}
\item ``(4)~$\Rightarrow$~(2)''.
\item If~$\varphi$ is a geometric formula, then~``(1)~$\Leftrightarrow$~(2)''.
\item If all subformulas of~$\varphi$ appearing as antecedents of implications
satisfy ``(1)~$\Rightarrow$~(2)'' (for instance, because they are geometric
formulas or because~$\varphi$ doesn't contain any~``$\Rightarrow$'' signs),
then~``(2)~$\Rightarrow$~(1)''.
\item If~$X$ is overt, then~``(1)~$\Leftrightarrow$~(2)~$\Leftrightarrow$~(3)''.
\end{itemize}
\end{lemma}

\begin{proof}The implication~``(4)~$\Rightarrow$~(2)'' is trivial. The other
claims can be checked by induction on the structure of~$\varphi$.
\end{proof}

\begin{rem}It is no coincidence that the conditions in
Lemma~\ref{lemma:properties-of-constant-sheaves-over-locales} are reminiscent
of the conditions in Lemma~\ref{lemma:open-stalk}. In fact, associated to any
locale~$X$ is a modal operator~$\Box_X$ on~$\Set$, which implicitly appeared in
Lemma~\ref{lemma:properties-of-constant-sheaves-over-locales}. It maps a truth
value~$\varphi$ to the truth value~$\brak{X \preceq
\bigvee\{\top\,|\,\varphi\}}$. The associated sublocale~$\pt_{\Box_X}$ of the
one-point locale is the image of the unique locale morphism~$X \to \pt$.\end{rem}

\begin{prop}\label{prop:spectrum-overt}
The localic spectrum of a ring~$A$ is overt if and only if any
element of~$A$ is nilpotent or not nilpotent.\end{prop}

\begin{proof}For the ``if'' direction, we can define a positivity predicate
by declaring~$\aaa > 0$ if and only if~$\aaa$ contains an element which is not nilpotent.

For the ``only if'' direction, let~$f \? A$ be an arbitrary element. Then
\begin{align*}
  \sqrt{(f)} &\subseteq \bigvee \{ \sqrt{(1)} \,|\, \sqrt{(f)} > 0 \} \\
  &\subseteq \bigvee \{ \sqrt{(1)} \,|\, \text{$f$ is not nilpotent} \} \\
  &= \{ s \? A \,|\, \text{$s$ is nilpotent or $f$ is not nilpotent} \}.
\end{align*}
Considering that~$f \in \sqrt{(f)}$, it follows that~$f$ is nilpotent or~$f$ is
not nilpotent.
\end{proof}

\begin{rem}I don't know when the local spectrum~$\Spec(A|R)$ is overt. This
question is related to openness of morphisms between schemes as follows.
Let~$X$ be a scheme (in a classical context). Let~$\A$ be a
quasicoherent~$\O_X$-algebra. Then the relative spectrum~$\RelSpec_X(\A)$
exists as a topological space, and is given by the externalization of the local
spectrum~$\Spec(\A|\O_X)$. If the canonical morphism~$\RelSpec_X(\A) \to X$ is
open, then the induced morphism of locales is open as well (the converse
doesn't hold in
general~\cite[Proposition~IX.7.5]{moerdijk-maclane:sheaves-logic}). This is the
case if and only if~$\Spec(\A|\O_X)$ is an overt locale from the internal point
of view of~$\Sh(X)$.
\end{rem}

Since Proposition~\ref{prop:spectrum-overt} shows that the spectrum of a ring
is in general not overt,
Lemma~\ref{lemma:properties-of-constant-sheaves-over-locales} is not applicable
to the spectrum in its full power. However, there is a substitute which is often
sufficient: For a ring element~$f \? A$, it holds that
\[ \sqrt{(f)} \subseteq \bigvee \{ \top \,|\, \varphi \}
  \qquad\text{if and only if}\qquad
  \speak{$f$ is nilpotent} \vee \varphi. \]
The case that~$f$ is nilpotent often trivializes the situation, allowing to
extend Lemma~\ref{lemma:properties-of-constant-sheaves-over-locales}, at least
morally. For instance, it still holds that an~$A$-module~$M$ is finitely
generated if and only if~$\ul{M}$ is finitely
generated as an~$\ul{A}$-module from the internal point of view
of~$\Spec(A)$. (This then implies that~$M$ is finitely generated if and only
if~$M^\sim$ is of finite type, as in the proof of
Lemma~\ref{lemma:finite-type-using-universal-filter}.) The ``only if''
direction is straightforward. For the ``if'' direction, we may assume that
we're given a covering~$\sqrt{(1)} = \bigvee_i \sqrt{(f_i)}$ such that, for
each~$i$, there are elements~$x_{i1},\ldots,x_{i,n_i} \? M$ satisfying
\[ \sqrt{(f_i)} \models \forall x\?\ul{M}\_
  \exists a_1,\ldots,a_{n_i}\?\ul{A}\_
  x = \textstyle\sum_j a_j x_{ij}. \]
Without loss of generality, we may assume that the covering is finite. We can
then verify that the joint system~$(x_{ij})_{ij}$ generates~$M$. Let~$x \? M$. For
each index~$i$, there exists a finite covering~$\sqrt{(f_i)} = \bigvee_k
\sqrt{(g_{ik})}$ such that, for each index~$k$, there exist
elements~$a_1,\ldots,a_{n_i} \? A$ such that
\[ \text{$g_{ik}$ is nilpotent} \qquad\text{or}\qquad
  x = \textstyle\sum_j a_j x_{ij}. \]
If the second case occurs for at least one pair~$(i,k)$ of indices, we are
done. Else all the~$g_{ik}$ are nilpotent. This implies that all the~$f_i$ are
nilpotent, which in turn implies that the unit of~$A$ is nilpotent. Thus~$A$ is
the zero ring. In this case~$x = 0$; thus we are done as well.

\subsection*{Gluing} The following definition is intuitionistically sensible:

\begin{defn}\label{defn:scheme}
An \emph{affine scheme} is a locally ringed locale which is isomorphic to the
spectrum of a ring. A \emph{scheme} is a locally ringed locale which is locally
(on an open cover) isomorphic to the spectrum of a ring.
\end{defn}

It's crucial that we're able to verify the affine communication
lemma~\cite[Lemma~5.3.2]{vakil:foag} in this setting; this is the lemma which
ensures that for many properties, there is no difference between mandating that
they hold for the members of some open affine cover and that they hold on any
affine open. Its validity rests solely on the following technical statement.

\begin{prop}Let~$(X,\O_X)$ be a locally ringed locale. Let~$U$ and~$V$ be opens of~$X$
such that~$(U,\O_X|_U)$ and~$(V,\O_X|_V)$ are affine. Then the meet~$U \wedge
V$ admits a covering by opens which are simultaneously standard opens
of~$(U,\O_X|_U)$ and of~$(V,\O_X|_V)$.
\end{prop}

\begin{proof}Since~$U$ is affine,
\begin{align*}
  U \wedge V &= \bigvee \{ W \preceq U \wedge V \,|\,
  \text{$W \hookrightarrow U$ is a standard open} \}. \\
\intertext{For any such open~$W$,}
  W &= \bigvee \{ W' \preceq W \,|\,
  \text{$W' \hookrightarrow W \hookrightarrow V$ is a standard open} \}
\end{align*}
since~$V$ is affine. We show that any such open~$W'$ is also standard open
in~$U$; this suffices to establish the claim.

Since~$W$ is standard open in~$U$, there is a function~$f \? \Gamma(U,\O_X)$
such that~$W = D(f)$, where
\[ D(f) = \bigvee \{ A \preceq X \,|\, A \models \speak{$f$ \inv} \}. \]
Since~$W'$ is standard open in~$V$, there is a function~$g \? \Gamma(V,\O_X)$
such that~$W' = D(g)$. The restriction~$g|_W$ can be regarded
as an element of~$\Gamma(U,\O_X)[f^{-1}]$; as such, it is of the form~$h/f^n$.
Then~$W' = D(f) \wedge D(h) = D(fh)$. The
open~$W'$ therefore coincides, as an open of~$U$, with~$\sqrt{(fh)}$ and is thus
standard open in~$U$.
\end{proof}

\subsection*{Properties of sheaves} Schemes in the sense of
Definition~\ref{defn:scheme} can't intuitionistically be shown to have enough
points~\cite{tierney:spectrum}. Classically, they can; this ensures that
classically there is no difference between the category of schemes as usually
defined and the category of schemes in the sense of
Definition~\ref{defn:scheme}.

As a consequence, properties of morphisms of sheaves can't be checked on
stalks. For instance, for a morphism~$\alpha : \G \to \H$ of sheaves on a
locale~$X$ to be an epimorphism it's not enough that~$\alpha_x : \G_x \to \H_x$
is surjective for all locale-theoretic points of~$X$. Instead, for every local
section~$s \? \H(U)$ there has to be a covering~$U = \bigvee_i U_i$ such that,
for each~$i$, there is a preimage of~$s|_{U_i}$.

Many of the results in Section~\ref{sect:sheaves-of-rings} and
Section~\ref{sect:sheaves-of-modules} thus have to be made into definitions.
For instance, a sheaf of modules should be declared \emph{flat} if and only if
it is flat as an ordinary module from the internal point of view.

The results in Section~\ref{sect:eliminating-prime-ideals} can be used to keep
up the appearance that testing on stalks suffices. For instance, let~$\alpha :
M^\sim \to N^\sim$ be a morphism of quasicoherent sheaves on~$\Spec(A)$. The
points of~$\Spec(A)$ are the filters of~$A$; but as remarked it doesn't suffice
to test the stalks~$\alpha_F : M[F^{-1}] \to N[F^{-1}]$. However, it does
suffice to test, internally to~$\Sh(\Spec(A))$, the map~$\ul{M}[\F^{-1}] \to
\ul{N}[\F^{-1}]$ (which is just~$\alpha$), where~$\F$ is the generic filter.

\subsection*{Big toposes} Just as classically, the big Zariski topos of a
scheme~$S$ can be defined as the topos of sheaves over the parsimonious
sites~$(\Aff/S)_\lfp$ or~$(\Sch/S)_\lfp$ (details about the possible choices
for the site are in Section~\ref{sect:proper-choice-of-site}). The proof that
this topos classifies local rings over~$S$ is intuitionistically valid.

We don't know whether all of the common subtoposes of the big Zariski topos
corresponding to finer topologies like the étale or fppf topology have all the
properties which are classically expected of them. In any case, if it's
classically known that the subtopos of the big Zariski topos corresponding to
a finer topology classifies a certain explicitly presented geometry
theory, one could adopt such a result as a definition in an intuitionistic
context. For instance, the big étale topos of a scheme~$S$ can be defined as
the classifying topos of separably closed local rings over~$S$ and the big fppf
topos can be defined as the classifying topos of fppf-local rings over~$S$
(Section~\ref{sect:beyond-zariski}).

\subsection*{Cohomology} We don't know how a general constructive framework for
cohomology might look like (besides Čech cohomology, which has its well-known
shortcomings) and can only remark that Grothendieck's approach using injective
resolutions can't work, since it's consistent with Zermelo--Fraenkel set theory
that no nontrivial injective abelian groups exist~\cite{blass:inj-proj-axc},
and that the account of Kempf~\cite{kempf:cohomology} looks promising, since he
employs flabby resolutions instead of injective ones.

However, Barakat and Lange-Hegermann pioneered constructive approaches to
cohomology of certain base schemes, which are not only mathematically elegant
but also work very well in practice (much more efficiently than Čech methods).
We refer to their articles for details~\cite{barakat-lh:homalg,barakat-lh:ext}.

}

\section{Higher direct images and other derived functors}

\subsection{Flabby sheaves}

Recall that a sheaf~$\F$ of sets on a topological space (or a locale)~$X$ is
\emph{flabby} if and only if, for any open subset~$U \subseteq X$ the
restriction map~$\F(X) \to \F(U)$ is surjective.

Flabbiness is a local property, even though it doesn't seem like that
at first sight: If the restrictions~$\F|_{U_i}$ of~$\F$ to the members of an
open covering~$X = \bigcup_i U_i$ are flabby, then the verification that~$\F$ is
flabby can't proceed as follows. ``Let~$s \in \F(U)$ be an arbitrary section.
Since each~$\F|_{U_i}$ is flabby, the section~$s|_{U \cap U_i}$ extends to a
section on~$U_i$.'' The reason is that the individual extensions obtained in
this way might not glue.

A correct proof employs Zorn's lemma in a typical way, considering a maximal
extension and then verifying that the subset this maximal extension is defined
on is all of~$X$.

Since flabbiness is a local property, it's not unreasonable to expect that
flabbiness can be characterized in the internal language. The following
proposition shows that this is indeed the case.

\begin{prop}\label{prop:internal-char-flabbiness}
Let~$\F$ be a sheaf of sets on a topological space~$X$ (or a locale).
Then the following statements are equivalent:
\begin{enumerate}
\item $\F$ is flabby.

\item ``Any section of~$\F$ can be locally extended'':
For any open~$U \subseteq X$ and any section~$s \in \F(U)$ there is
an open covering~$X = \bigcup_i V_i$ such that, for each~$i$, there is an
extension of~$s$ to~$U \cup V_i$ (that is, a section~$s' \in \F(U \cup V_i)$
such that~$s'|_U = s$).

(If~$X$ is a space instead of a locale, this can be equivalently formulated as
follows: For any open subset~$U \subseteq X$, any section~$s \in \F(U)$, and any
point~$x \in X$, there is an open neighborhood~$V$ of~$x$ and an extension
of~$s$~to~$U \cup V$.)

\item From the point of view of the internal language of~$\Sh(X)$, for any
subsingleton~$K \subseteq \F$ there exists an element~$s \? \F$ such that~$s \in
K$ if~$K$ is inhabited. More precisely,
\begin{multline*}
  \qquad\qquad\Sh(X) \models
  \forall K \subseteq \F\_
  (\forall s,s'\?K\_ s = s') \Longrightarrow \\
  \exists s\?\F\_ (\text{$K$ is inhabited} \Rightarrow s \in K).\qquad\qquad
\end{multline*}

\item The canonical map~$\F \to \P_{\leq 1}(\F), s \mapsto \{s\}$ is
final from the internal point of view, that is
\[ \Sh(X) \models
  \forall K \? \P_{\leq 1}(\F)\_
  \exists s \? \F\_
  K \subseteq \{s\}, \]
where~$\P_{\leq 1}(\F)$ is the object of subsingletons of~$\F$.
\end{enumerate}
\end{prop}

\begin{proof}
The implication~``(1)~$\Rightarrow$~(2)'' is trivial. The converse direction uses a
typical argument with Zorn's lemma, considering a maximal extension. The
equivalence~``(2)~$\Leftrightarrow$~(3)'' is routine, using the Kripke--Joyal
semantics to interpret the internal statement. Condition~(4) is a straightforward
reformulation of Condition~(3).
\end{proof}

Condition~(2) of the proposition is, unlike the standard definition of flabbiness,
manifestly local. Also its equivalence with Condition~(3) and Condition~(4) is
intuitionistically valid; therefore one might consider to adopt Condition~(2) as
the definition of flabbiness.

The object~$\P_{\leq 1}(\F)$ of subsingletons of $\F$ can be
interpreted as the object of \emph{partially-defined
elements} of $\F$. In this view, the empty subset is the maximally undefined
element and a singleton is a maximally defined element. In classical logic,
there are no further examples of partially-defined elements, but
intuitionistically, there might; and indeed, in the model of intuitionistic
logic provided by~$\Sh(X)$, there are many more. An explicit description of the
sheaf~$\P_{\leq 1}(\F)$ is given in Remark~\ref{rem:godement-construction}.

The proposition shows that a sheaf~$\F$ is flabby if and only if any
partially-defined element of~$\F$ can be refined to an honest element of~$\F$.

\subsection{Injective sheaves}

Recall that an object~$I$ of a category~$\C$ is \emph{injective} if and only if,
for any monomorphism~$X \hookrightarrow Y$ in~$\C$ and any morphism~$X \to I$, there
is an extension such that the diagram
\[ \xymatrix{
  X \ar@{^{(}->}[r]\ar[d] & Y \ar@{-->}[dl] \\
  I
} \]
commutes. Equivalently, an object~$I$ is injective if and only if the Hom
functor~$\Hom_\C(\placeholder, I) : \C^\op \to \Set$ maps monomorphisms in~$\C$ to
surjective maps. This general definition is often specialized to one of these cases:
to the category of modules over a ring, to the category of set-valued sheaves
on a topological space, and to the category of sheaves of~$\O_X$-modules on a
ringed space~$(X,\O_X)$.

The definition is seldom applied in the category of sets, since in a classical
context it's easy to show that a set is injective if and only if it's
inhabited, thereby completely settling the question which objects are
injective in a trivial manner.

The question is more interesting in an intuitionistic setting, since
intuitionistically one cannot prove that inhabited sets are
injective~\cite{aczel-et-al:injective}; but one can still verify that any set embeds
into an injective set: The powerset~$\P(X)$ and even the smaller
set~$\P_{\leq1}(X)$ of subsingletons of a given set~$X$ are injective. This
fact is well-known in the constructive mathematics community, but for
convenience we spell out the proof as Lemma~\ref{lemma:enough-flabby}.

For a cartesian or monoidal closed category~$\C$, there is also the notion of an
\emph{internally injective} object. This is an object~$I$ such that the
internal Hom functor~$[\placeholder, I] : \C^\op \to \C$ maps monomorphisms
in~$\C$ to epimorphisms. In the special case that~$\C$ is a elementary topos
with a natural numbers object, such as the topos of set-valued sheaves on a
space, this condition can be rephrased in several ways. The following
proposition lists five of these conditions. The equivalence of the first four is
due to Harting~\cite{harting}.

\begin{prop}\label{prop:notions-of-internal-injectivity}
Let~$\E$ be an elementary topos. Then
the following statements about an object~$I \in \E$ are equivalent.
\begin{enumerate}
\item $I$ is internally injective.
\item The functor~$[\placeholder, I] : \E^\op \to \E$ maps monomorphisms in $\E$
to morphisms for which any global element of the target locally (after change of
base along an epimorphism) possesses a preimage.
\item For any morphism $p : A \to 1$ in $\E$, the object $p^*I$ has property~(1)
as an object of $\E/A$.
\item For any morphism $p : A \to 1$ in $\E$, the object $p^*I$ has property~(2)
as an object of $\E/A$.
\item From the point of view of the internal language of~$\E$, the object~$I$
is injective.\footnote{In Section~\ref{sect:internal-language}, we have only
introduced the internal language for sheaf toposes. The general definition is
in~\cite[Section~7]{shulman:stack}.}
\end{enumerate}
\end{prop}

\begin{proof}
The implications ``(1)~$\Rightarrow$~(2)'', ``(3)~$\Rightarrow$~(4)'',
``(3)~$\Rightarrow$~(1)'', and ``(4)~$\Rightarrow$~(2)'' are trivial.

The equivalence ``(3)~$\Leftrightarrow$~(5)'' follows directly from the
interpretation rules of the stack semantics.

The implication ``(2)~$\Rightarrow$~(4)'' employs the
extra left adjoint $p_! : \E/A \to \E$ of $p^* : \E
\to \E/A$~(which maps an object~$(X \to A)$ to~$X$), as in the usual proof that
injective sheaves remain injective when
restricted to smaller open subsets: We have that $p_* \circ [\placeholder, p^*I]_{\E/A}
\cong [\placeholder, I]_\E \circ p_!$, the functor $p_!$ preserves monomorphisms, and one
can check that $p_*$ reflects the property that global elements locally possess
preimages. Details are in~\cite[Thm.~1.1]{harting}.\footnote{Harting formulates
the statement for abelian group objects, and has to assume that~$\E$ contains a
natural numbers object to ensure the existence of an abelian version of~$p_!$.}

The implication ``(4)~$\Rightarrow$~(3)'' follows by performing an extra change of
base, since any non-global element becomes a global element after a suitable
change of base.
\end{proof}

Somewhat surprisingly, and in stark contrast with the situation for internally
projective objects (which are defined dually), internal injectivity coincides
with external injectivity for sheaf toposes over spaces.

\begin{thm}\label{thm:char-injectivity}
Let~$X$ be a topological space (or a locale). An object~$\I \in \Sh(X)$ is
injective if and only if it is internally injective.
\end{thm}

\begin{proof}For the ``only if'' direction, let~$\I$ be an injective sheaf of
sets. Then~$\I$ satisfies Condition~(2) in
Proposition~\ref{prop:notions-of-internal-injectivity}, even without having to
pass to covers.

For the ``if'' direction, let~$\I$ be an internally injective object. Let~$m :
\E \hookrightarrow \F$ be a monomorphism in~$\Sh(X)$ and let~$k : \E \to \I$ be an arbitrary
morphism. We want to show that there exists an extension $\F \to \I$ of~$k$
along~$m$. To this end, we consider the sheaf defined by the internal expression
\[ \G \defeq \brak{\{ k' \? [\F,\I] \,|\, k' \circ m = k \}}. \]
Global sections of~$\G$ are extensions of the kind we're looking for.
Therefore it suffices to show that~$\G$ is flabby. We do this by verifying
Condition~(3) of Proposition~\ref{prop:internal-char-flabbiness} in the internal
language of~$\Sh(X)$.

Let~$K \subseteq \G$ be a subsingleton. We consider the injectivity diagram
\[ \xymatrix{
  m[\E] \cup \F' \ar@{^{(}->}[r]\ar[d] & \F \ar@{-->}[ld] \\
  \I
} \]
where~$\F'$ is the set~$\{ s \? \F \,|\, \text{$K$ is inhabited} \}$ and the solid
vertical arrow is defined in the following way: It should map an element~$s \in
\F'$ to~$k'(s)$, where~$k'$ is any element of~$K$; and it should map an
element~$m(u) \in m[\E]$ to~$k(u)$. These prescriptions determine a well-defined
map.

Since~$\I$ is injective from the internal point of view we're taking up here,
there exists a dotted map rendering the diagram commutative. This map is an
element of~$\G$. Furthermore, if~$K$ is inhabited, then this map is an element of~$K$.
\end{proof}

\begin{thm}\label{thm:char-injectivity-modules}
Let~$(X,\O_X)$ be a ringed topological space (or a ringed locale).
An~$\O_X$-module~$\I$ is injective if and only if it is internally injective.
\end{thm}

\begin{proof}Proposition~\ref{prop:notions-of-internal-injectivity} can be
adapted from sheaves to sets to sheaves of modules, with the same proof.
The extra left adjoint~$p_! : \Mod_{\Sh(X)/A}(\O_X \times A) \to
\Mod_{\Sh(X)}(\O_X)$ required by the proof maps a module~$M \to A$ to the
internal direct sum~$\bigoplus_{a \? A} M(a)$.

The proof of Theorem~\ref{thm:char-injectivity} can be adopted as well.
It suffices to change~``$[\F, \I]$'' to~``$[\F,\I]_{\Mod(\O_X)}$'' (denoting
the Hom sheaf of~$\O_X$-linear morphisms~$\F \to \I$)
and~``$m[\E] \cup \F'$'' to~``$m[\E] + \F''$'', where~$\F''
\defeq \{ s \? \F \,|\, \text{$s = 0$ or $K$ is inhabited} \}$.
\end{proof}

\begin{rem}The proof of Theorem~\ref{thm:char-injectivity} crucially rests on
Proposition~\ref{prop:internal-char-flabbiness} and therefore on Zorn's lemma,
to ensure that the sheaf~$\G$ defined in the proof, which has the property that
any of its sections can be locally extended, admits a global section. The proof
is therefore not intuitionistically valid.

On a related note, we don't think that the statement of
Theorem~\ref{thm:char-injectivity} can be generalized to arbitrary
(Grothendieck) toposes. The proof gradually refines the trivial
generalized element of~$\G$ (defined on the empty stage) to a global element.
Such a procedure is not really meaningful for sheaf toposes over sites for
which not any object is a subobject of the terminal object.
\end{rem}

\subsection{Internal proofs of common lemmas}

\begin{lemma}A sheaf of sets or a sheaf of modules is injective if and only if it
is locally injective.\end{lemma}

\begin{proof}By Theorem~\ref{thm:char-injectivity} respectively
Theorem~\ref{thm:char-injectivity-modules}, injectivity can be characterized in the
internal language. Any such property is local.\end{proof}

\begin{lemma}Let~$X$ be a topological space (or a locale).
\begin{enumerate}
\item Let~$\I$ be an injective sheaf of sets over~$X$. Let~$\F$ be an arbitrary
sheaf of sets. Then~$\HOM(\F,\I)$ is flabby.
\item Let~$\I$ be an injective sheaf of modules over some sheaf~$\O_X$ of rings
over~$X$. Let~$\F$ be an arbitrary sheaf of modules.
Then~$\HOM_{\O_X}(\F,\I)$ is flabby.
\end{enumerate}
\end{lemma}

\begin{proof}
We first cover the case of sheaves of sets. By Theorem~\ref{thm:char-injectivity}
and Proposition~\ref{prop:internal-char-flabbiness}, it suffices to give an
intuitionistic proof of the following statement: If~$I$ is an injective set
and~$F$ is an arbitrary set, then partially defined elements of the set~$[F,I]$
of all maps~$F \to I$ can be refined to honest elements.

Thus let a subsingleton~$K \subseteq [F,I]$ be given. We consider the
injectivity diagram
\[ \xymatrix{
  F' \ar[r]\ar[d] & F \ar@{-->}[ld] \\
  I
} \]
where~$F'$ is the subset~$\{ s \? F \,|\, \text{$K$ is inhabited} \} \subseteq F$ and the
solid vertical map sends~$s \in F'$ to~$f(s)$, where~$f$ is an arbitrary element
of~$K$. This association is well-defined. Since~$I$ is injective, a dotted lift
as indicated exists. If~$K$ is inhabited, this lift is an element of~$K$.

The same kind of argument applies to the case of sheaves of modules, relying on
Theorem~\ref{thm:char-injectivity-modules} and defining~$F'$ as the
submodule~$\{ s \? F \,|\, \text{$s = 0$ or $K$ is inhabited} \}$.
\end{proof}

\begin{cor}Injective sheaves of sets and injective sheaves of modules are
flabby.\end{cor}

\begin{proof}Follows from the previous lemma by considering the special cases~$\F
\defeq 1$ respectively~$\F \defeq \O_X$.\end{proof}

\begin{lemma}\label{lemma:enough-flabby}
Let~$X$ be a topological space (or a locale). Any sheaf of sets
over~$X$ can be embedded into an injective (therefore flabby) sheaf of sets.
\end{lemma}

\begin{proof}By Proposition~\ref{prop:internal-char-flabbiness}, it suffices to
give an intuitionistic proof of the following statement: Any set~$F$ can be
embedded into an injective set.

As already indicated, there are at least two simple ways how~$F$ can be embedded
into an injective set: by embedding~$F$ in its powerset~$\P(F)$ or by
embedding~$F$ in~$\P_{\leq1}(F)$, the set of subsingletons of~$F$. For
conciseness, we only verify that~$\P_{\leq1}(F)$ is injective.

So let~$m : A \hookrightarrow B$ be an injective map and let~$k : A \to
\P_{\leq1}(F)$ be an arbitrary map. Then we can extend~$k$ to a map~$k' : B \to
\P_{\leq1}(F)$ by defining for~$y \? B$
\begin{align*}
  k'(y) &\defeq \bigcup k[m^{-1}[\{y\}]] \\
  &\phantom{\vcentcolon}= \{ s\?F \,|\, \text{$s \in k(x)$ for some~$x \in A$ such that~$m(x) = y$} \}.
  \qedhere
\end{align*}
\end{proof}

\begin{rem}\label{rem:godement-construction}
The \emph{Godement construction} provides a well-known way of embedding an
inhabited sheaf of sets~$\F$ into an injective sheaf, namely
the sheaf of not necessarily continuous sections of the étale space of~$\F$:
\[ U \subseteq X \quad\longmapsto\quad
  \prod_{x \in U} \F_x. \]
The sheaf~$\P_{\leq1}(\F)$ does not coincide with this construction.
Instead by Definition~\ref{defn:interpretation-internal-constructions}, it is the sheaf with
\[ U \subseteq X \quad\longmapsto\quad
  \{ \langle V, s \rangle \,|\,
    \text{$V \subseteq U$ open, $s \in \F(V)$} \}. \]
It's not possible to describe the Godement construction in the internal language
of~$\Sh(X)$, since the Godement construction depends on the underlying set of~$X$,
but the sheaf topos of~$X$ doesn't remember this set. For instance, if~$X$
is an inhabited indiscrete topological space, then~$\Sh(X)$ is equivalent
to~$\Set$.
\end{rem}

\begin{rem}It's not known to me whether it's possible to intuitionistically
prove that any module can be embedded into a module which satisfies the internal
flabbiness criterion of Proposition~\ref{prop:internal-char-flabbiness}. This
would give an internal proof of the well-known fact that any sheaf of modules can be embedded into a
flabby sheaf of modules. The naive candidates don't work: The
set~$\P_{\leq1}(F)$ doesn't admit a canonical module structure (though it does admit the structure of a commutative monoid), and the free
module over that set is not flabby in general.

Since by the Godement construction the statement that any sheaf of modules can
be embedded into a flabby sheaf of modules is true in many models of
intuitionistic logic, the sheaf toposes over topological
spaces,\footnote{However it should be noted that the Godement construction
doesn't work in an intuitionistic metatheory, unless the underlying set of
points of the given topological space is postulated to have decidable
equality.}
it's not entirely unreasonable to believe
that such an intuitionistic proof is possible.\footnote{There is a metatheorem
guaranteeing that a statement is intuitionistically provable if and only if it
holds in the sheaf topos over any topological
space~\cite[Theorem~B]{awodey-butz:completeness}. However, this metatheorem
requires the considered statements to be of a certain form, which in particular
forbids them from mentioning the object of truth values. The internal
statements given in Proposition~\ref{prop:internal-char-flabbiness} depend on
this object in a crucial way.}

On the other hand, it's certainly not possible to intuitionistically prove that
any module can be embedded into an injective module, since it's consistent with
Zermelo--Fraenkel set theory that no nontrivial injective abelian groups
exist~\cite{blass:inj-proj-axc}.
\end{rem}

\begin{lemma}Let~$X$ be a ringed space (or a ringed locale). Let~$0 \to \E'
\xra{\alpha} \E \xra{\beta} \E'' \to 0$ be a short exact sequence
of~$\O_X$-modules. If~$\E'$ is flabby, then the induced sequence
\[ 0 \longrightarrow \Gamma(X,\E')
  \longrightarrow \Gamma(X,\E)
  \longrightarrow \Gamma(X,\E'') \longrightarrow 0 \]
is exact.
\end{lemma}

\begin{proof}Since taking global sections is left exact (being a right adjoint
functor), it suffices to verify that the map~$\Gamma(X,\E) \to \Gamma(X,\E'')$
is surjective. We'll do this by showing, in the internal language of~$\Sh(X)$,
that the sheaf of preimages of a given global section~$s \in \Gamma(X,\E'')$ is
flabby and therefore has a global section.

In the internal language, this sheaf has the description~$F \defeq \{ u \? \E
\,|\, \beta(u) = s \}$. To verify the internal condition of
Proposition~\ref{prop:internal-char-flabbiness}, let a subsingleton~$K \subseteq
F$ be given. Since~$\beta$ is surjective, there is a preimage~$u_0 \in F$.
The translated set~$K - u_0 \subseteq \E$ is still a subsingleton, and its
preimage under~$\alpha$ is as well. By the assumption on~$\E'$, there is an
element~$v \? \E'$ such that~$v \in \alpha^{-1}[K - u_0]$ if~$\alpha^{-1}[K - u_0]$ is
inhabited. We'll now verify that~$u_0 + \alpha(v) \in K$ if~$K$ is inhabited.

So assume that~$K$ is inhabited. Then~$K - u_0$ is as well. Since the image of
its unique element under~$\beta$ is zero and the given sequence is exact, the
set~$\alpha^{-1}[K - u_0]$ is inhabited as well. Therefore~$v \in \alpha^{-1}[K
- u_0]$. Thus~$u_0 + \alpha(v) \in K$.
\end{proof}

\begin{lemma}Let~$X$ be a ringed space (or a ringed locale). Let~$0 \to \E'
\xra{\alpha} \E \xra{\beta} \E'' \to 0$ be a short exact sequence
of~$\O_X$-modules. If~$\E'$ and~$\E''$ are flabby, then~$\E$ is flabby as well.
\end{lemma}

\begin{proof}We verify the condition of
Proposition~\ref{prop:internal-char-flabbiness} in the internal language
of~$\Sh(X)$.

Let~$K \subseteq \E$ be a subsingleton. Then its image~$\beta[K] \subseteq \E''$
is a subsingleton as well. Since partial elements of~$\E''$ can be refined to
honest elements, there is an element~$s \? \E''$ such that~$\beta[K] \subseteq
\{ s \}$.

Since~$\beta$ is surjective, there is an element~$t_0 \? \E$ such
that~$\beta(t_0) = s$.

The preimage~$\alpha^{-1}[K - t_0] \subseteq \E'$ is a subsingleton. This
partial element can be refined to an honest element, so there exists an
element~$u \? \E'$ such that~$\alpha^{-1}[K - t_0] \subseteq \{u\}$.

The partial element~$K$ can thereby refined to the honest element~$t \defeq t_0 + \alpha(u)$.
\end{proof}

\subsection{Tor and sheaf Ext}
\label{sect:sheaf-ext-and-tor}

The following lemma expresses a prototype result for constructing sheaves in
the internal language. We'll use it to internally define derived functors.

\begin{lemma}\label{lemma:descent-for-modules}
Let~$X$ be a ringed topological space (or a ringed locale).
Let~$\varphi(\E)$ be a property of sheaves of~$\O_X$-modules, formulated in the internal
language of~$\Sh(X)$. Let~$\psi(f)$ be a property of morphisms
of sheaves of~$\O_X$-modules, formulated in the internal
language of~$\Sh(X)$ and stable under composition. Assume that
\begin{enumerate}
\item $\Sh(X) \models \exists \textnormal{$\E$ $\O_X$-module}\_ \varphi(\E)$
and
\item $\Sh(X) \models \forall \textnormal{$\E,\E'$ $\O_X$-modules}\_ \varphi(\E) \wedge \varphi(\E')
  \Longrightarrow
  \exists! \textnormal{$f : \E \to \E'$ linear}\_ \psi(f)$.
\end{enumerate}
Then there exists a sheaf~$\E$ of~$\O_X$-modules such that~$\Sh(X) \models \varphi(\E)$,
and any two such sheaves are isomorphic via a unique isomorphism which
satisfies~$\psi$ from the internal point of view.
\end{lemma}

\begin{proof}This is a reformulation of the well-known fact that we have
descent for sheaves of~$\O_X$-modules. By the first assumption, there is an
open covering~$X = \bigcup_i U_i$ such that for each~$i$, there is
an~$\O_X|_{U_i}$-module~$\E_i$ with~$U_i \models \varphi(\E_i)$. By the second
assumption and by Proposition~\ref{prop:simplification}, for each pair~$(i,j)$
of indices there is a unique morphism~$f_{ij} : \E_i|_{U_i \cap U_j} \to
\E_j|_{U_i \cap U_j}$ such that~$U_i \cap U_j \models \psi(f_{ij})$. Since the
property~$\psi$ is stable under composition, these morphisms are isomorphisms
which satisfy the cocycle condition. Thus the~$\O_X$-modules~$\E_i$ glue to a
global~$\O_X$-module~$\E$, which satisfies property~$\varphi$ because it does
so locally.

The uniqueness claim is immediate by Proposition~\ref{prop:simplification} and
by the assumption that property~$\psi$ is stable under composition.
\end{proof}

Lemma~\ref{lemma:descent-for-modules} can be generalized in two
ways: from sheaves of modules to other kinds of algebraic structures, for
instance complexes of sheaves of modules; and from sheaf toposes over locales
to more general Grothendieck toposes, by the descent theorem for Grothendieck
toposes. We will use the former, but not the latter generalization.

\begin{lemma}\label{lemma:internal-resolutions}
Let~$X$ be a ringed space (or a ringed locale). From the internal point of view
of~$\Sh(X)$, any~$\O_X$-module admits a resolution by injective~$\O_X$-modules,
and any two such are related by a morphism of complexes which is unique up to
homotopy with the property that it induces the identity on the resolved module.
\end{lemma}

\begin{proof}There can't be an intuitionistic proof of this fact, since it's
consistent with Zermelo--Fraenkel set theory
that no nontrivial injective abelian groups exist~\cite{blass:inj-proj-axc}.
But working in a classical metatheory, it's well-known that, for any open
subset~$U \subseteq X$, the category of sheaves of~$\O_X|_U$-modules has enough
injectives. Since externally injective sheaves of modules look like injective
modules from the internal point of view, by (the easy part of)
Theorem~\ref{thm:char-injectivity-modules}, the internal
statement~``any~$\O_X$-module can be embedded into an injective~$\O_X$-module''
holds.

Under the assumption of the existence of enough injectives, the usual proof
that any object admits a resolution by injective objects is intuitionistically
valid. We can therefore interpret this proof in the internal language
of~$\Sh(X)$ and conclude -- were it not for a subtle fine point regarding the
failure of the axiom of countable choice in~$\Sh(X)$.

A resolution is an infinite complex of modules. The assumption of the existence
of enough injectives allows us to extend any finite partially-constructed
injective resolution to a longer one; but collecting all of the resulting
injective objects into a complete resolution requires some form of choice.

There are two ways to counter this problem. If one wants to prove the lemma
exactly as stated, one has to construct the injective resolution externally
(and appeal to the axiom of choice in the metatheory) instead of internally
constructing it step by step. But for many purposes, there's also an
alternative: Often, a full injective resolution isn't actually needed. For
instance, for evaluating an~$n$-th derived functor on an object, it suffices
to have a finite partial injective resolution. If one adopts this stance, then
it's enough to adopt the statement ``any~$\O_X$-module can be embedded into an
injective~$\O_X$-module'' as an axiom in the internal language.
\end{proof}

\begin{rem}It's known that the axiom of choice suffices to construct injective
resolutions of abelian groups,\footnote{More precisely, it is intuitionistically
provable that any abelian group admits a resolution by divisible groups, even a
canonical such. Some form of choice is needed to verify that such a resolution
is actually a resolution by injective abelian groups.}
and also that the axiom of choice implies the law of excluded middle in the
presence of the other axioms of set theory (Diaconescu's theorem).
Lemma~\ref{lemma:internal-resolutions} shows that the axiom ``any abelian group
can be embedded into an injective abelian group'' does not imply the law of
excluded middle, since (assuming the axiom of choice in the metatheory) this
statement is true in the internal language of the sheaf topos over any
topological space and such toposes typically do not satisfy the law of excluded
middle.
\end{rem}

We use Lemma~\ref{lemma:internal-resolutions} as follows to construct the sheaf
Ext in the internal language. Let~$\E$ and~$\F$ be~$\O_X$-modules over a ringed
space (or a ringed locale). Internally, we define~$\EXT^n(\E,\F) \defeq
H^n([\E,\I^\bullet]_{\Mod(\O_X)})$ where~$0 \to \F \to \I^\bullet$ is an
injective resolution and~$[\E,\I^k]_{\Mod(\O_X)}$ is the set of~$\O_X$-linear
maps~$\E \to \I^k$. The module constructed in this way depends on the chosen
injective resolution, but for any two such resolutions, there is a unique
isomorphism in cohomology which is induced by a morphism of resolutions.

Externally, this definition gives rise to a well-defined sheaf on~$X$, by
arguing similarly as in the proof of Lemma~\ref{lemma:descent-for-modules}: We
obtain~$\O_X$-modules on an open cover; on intersections, we find coherent
isomorphisms by the uniqueness statement; therefore we can glue.
The~$\O_X$-module constructed in this way coincides with the sheaf Ext as
usually conceived.

Along the same lines, we can construct Tor sheaves internally. Let~$\E$
and~$\F$ be~$\O_X$-modules. Assume that~$\O_X$ and~$\F$ are coherent. Internally, we
define~$\TOR_n(\E,\F) \defeq H_n(\E \otimes_{\O_X} \P_\bullet)$,
where~$\P_\bullet \to \F \to 0$ is a projective resolution.
Such a resolution exists; in fact, we can resolve~$\F$ by
finite free modules (which are projective even without the axiom of choice) by
the coherence assumptions on~$\F$ and~$\O_X$. As
with~$\EXT^n$, the module constructed in this way is unique up to a unique
isomorphism induced by a morphism of resolutions.

\subsection{Higher direct images}\label{sect:higher-direct-images}

Higher direct images are thought of as a relative ``fiberwise'' version of
cohomology. One way to make this precise is to show that higher direct images
can be made sense of internally to the topos of sheaves over the base, where
they then look like ordinary cohomology.

Let~$f : Y \to X$ be a morphism of ringed locales. As discussed in
Section~\ref{sect:internal-locales}, there is a locale~$I(Y)$ internal
to~$\Sh(X)$ mirroring~$Y$; from the point of view of~$\Sh(X)$, the given
morphism~$f$ looks like the unique morphism~$I(Y) \to \pt$.

Let~$\E$ be a sheaf of~$\O_Y$-modules on~$Y$. This sheaf corresponds to a sheaf
on~$I(Y)$. Internally to~$\Sh(X)$, we can take an injective
resolution~$\J^\bullet$ of this sheaf and define~$H^n(I(Y), \E) \defeq
H^n(\Gamma(I(Y), \J^\bullet))$. Just like with sheaf Ext and Tor presented in
Section~\ref{sect:sheaf-ext-and-tor}, this internal description gives rise to
a sheaf of~$\O_X$-modules on~$X$; the sheaf constructed in this way coincides
with the higher direct image~$R^n f_*(\E)$ as usually defined.

The conception of higher direct images as internal cohomology entails that
basic statements about cohomology yield corresponding statements about higher
direct images. For instance:

\begin{itemize}
\item That~$R^{\geq1} \id_*(\E)$ vanishes is a reflection of the fact
that~$H^{\geq1}(\pt, \E)$ vanishes.
\item Čech methods to compute cohomology entail Čech methods to compute higher
direct images.
\item The computation of the cohomology of projective~$n$-space over a ring
immediately yields the higher direct images of the canonical morphism~$\PP^n_S
\to S$ for any base scheme~$S$.
\end{itemize}

Also the failure of higher direct images to commute with arbitrary base change
gains a logical interpretation. Let
\[ \xymatrix{
  Y' \ar[r]^{g'}\ar[d]_{f'} & Y \ar[d]^f \\
  X' \ar[r]_g & X
} \]
be a pullback diagram of locally ringed locales. Let~$\E$ be a sheaf of modules on~$Y$.
If taking cohomology was a geometric construction, then it would be preserved
by pullback along arbitrary geometric morphisms. Since higher direct images are
just cohomology from the internal point of view, we would therefore have a
canonical isomorphism
\[ g^*(R^nf_*(\E)) \cong R^nf'_*((g')^*\E). \]
However, taking cohomology is not a geometric construction, and indeed in
general there is no such isomorphism. It's an open question whether the
well-known cases where there is such an isomorphism can be treated by a purely
or mostly logical framework.

\chapter{The big Zariski topos}\label{part:big-zariski}

The preceding part demonstrated that working in the internal universe of
the little Zariski topos of a scheme~$S$, the topos of sheaves on~$S$, is
useful for simplifying local work on~$S$. The basic tenet was that sheaves of
modules look just like plain modules and that theorems of intuitionistic
algebra yield theorems about sheaves.

But the little Zariski topos is not particularly suited for dealing with
\emph{schemes} over~$S$. For this, we need a related topos. For the scope of
this introduction only, we blithely employ the following slightly problematic
definition which we'll correct in Section~\ref{sect:proper-choice-of-site}.
We'll keep the base scheme~$S$ fixed throughout this part.

For some material in this part, we assume basic familiarity with classifying
toposes.

\section{Basics}

\begin{defn}[provisional]The \emph{big Zariski topos}~$\Zar(S)$ of a scheme~$S$ is the
topos of sheaves on the Grothendieck site~$\Sch/S$ of schemes over~$S$.
\end{defn}

Explicitly, an object of~$\Zar(S)$ is a functor~$F : (\Sch/S)^\op \to \Set$
satisfying the gluing condition with respect to ordinary Zariski coverings:
If~$X = \bigcup_i U_i$ is a cover of an~$S$-scheme~$X$ by open subsets, the
canonical diagram
\[ F(X) \longrightarrow \prod_i F(U_i) \xbigtoto{} \prod_{j,k} F(U_j \cap U_k)
\]
should be an equalizer diagram.

{\tocless

\subsection*{Internal language} Just like the topos of sheaves on a topological
space or on a locale admits an internal language, so does the big Zariski topos.
The necessary modifications of the Kripke--Joyal semantics
(Definition~\ref{defn:kripke-joyal}) are straightforward. Instead of
defining recursively the meaning of~``$U \models \varphi$'' for open subsets~$U
\subseteq S$, we define the meaning of~``$T \models \varphi$''
for~$S$-schemes~$T$ and slightly rewrite the rules for implication and universal quantification.
Instead of
\begin{align*}
  U \models \varphi \Rightarrow \psi \quad&\Ll\quad
  \text{for all open~$V \subseteq U$:} \\
  &\qquad\qquad\qquad\qquad\text{$V \models \varphi$ implies $V \models \psi$} \\
  U \models \forall s\?\F\_ \varphi(s) \quad&\Ll\quad
  \text{for all sections~$s \in \Gamma(V, \F)$ on open $V \subseteq U$:} \\
  &\qquad\qquad\qquad\qquad V \models \varphi(s) \\
\intertext{they have to read as follows.}
  T \models \varphi \Rightarrow \psi \quad&\Ll\quad
  \text{for all morphisms~$T' \to T$ in~$\Sch/S$:} \\
  &\qquad\qquad\qquad\qquad\text{$T' \models \varphi$ implies $T' \models \psi$} \\
  T \models \forall s\?F\_ \varphi(s) \quad&\Ll\quad
  \text{for all morphisms~$T' \to T$ in~$\Sch/S$ and} \\
  &\phantom{{}\Ll\quad\text{for }} \text{all sections~$s \in \Gamma(T', F)$:} \\
  &\qquad\qquad\qquad\qquad T' \models \varphi(s)
\end{align*}

The analogues of Proposition~\ref{prop:locality-of-the-internal-language} and
Proposition~\ref{prop:soundness-of-the-internal-language} are true for the internal
language of the big Zariski topos:

\begin{prop}\label{prop:basic-properties-language-big}
Let~$T$ be an~$S$-scheme and~$\varphi$ be a formula over~$T$.
\begin{enumerate}
\item If~$T \models \varphi$ and if there is an intuitionistic proof
that~$\varphi$ implies a further formula~$\psi$, then~$T \models \psi$.
\item Let~$T' \to T$ be a morphism of~$S$-schemes. If~$T \models \varphi$,
then~$T' \models \varphi$.
\item If~$T = \bigcup_i T_i$ is an open covering and if~$T_i \models \varphi$
for all~$i$, then~$T \models \varphi$.
\end{enumerate}
\end{prop}

\begin{proof}The proofs of
Proposition~\ref{prop:locality-of-the-internal-language} and
Proposition~\ref{prop:soundness-of-the-internal-language} carry over.
\end{proof}

When working with the internal language of the little Zariski topos, we often
used the fact that if a formula holds on some open subset~$U$, then it also
holds on all open subsets contained in~$U$.
Proposition~\ref{prop:basic-properties-language-big}(2) states a stronger version
of this: All properties which can be expressed using the internal language of
the big Zariski topos are automatically \emph{stable under base change}.

\subsection*{Important objects in the big Zariski topos}\label{page:important-objects} It's convenient to
introduce notation for objects which often appear when working with the big
Zariski topos.

Let~$X$ be an~$S$-scheme. Its functor of points, which maps
an~$S$-scheme~$T$ to~$\Hom_S(T,X)$, is an object of~$\Zar(S)$. We denote it
by~``$\ul{X}$''.

From the internal point of view of~$\Zar(S)$, such a functor~$\ul{X}$ looks
like a single set. It can be pictured as the ``set of points of~$X$'',
where~``point'' doesn't mean ``point of the underlying topological space
of~$X$'', but rather~``$T$-point of~$X$'', where~$T$ varies over
all~$S$-schemes. The internal language of the big Zariski topos hides any
explicit mentions of the stage~$T$; it is therefore a device for reifying the
multitude of points of~$X$, defined on varying stages, as a single entity.

Particularly important is~$\affl$, the functor of points of the affine line
over~$S$. The object~$\ul{S}$ is the terminal object in~$\Zar(S)$. This fits
into the philosophy: From the point of view of the big Zariski topos, the base
scheme should simply look like a point. The functor of points of~$S \amalg S$
looks like a two-element set from the internal point of view.

Let~$\F$ be a sheaf of sets on~$S$. For reasons explained in
Section~\ref{sect:relation-big-little}, we denote by~``$\pi^{-1}\F$''
the induced sheaf on~$\Sch/S$ mapping an~$S$-scheme
$(f : T \to S)$ to~$\Gamma(T, f^{-1}\F)$.

Let~$\F$ be a sheaf of~$\O_S$-modules. We denote by~``$\F^\Zar$''
the induced sheaf on~$\Sch/S$ mapping an~$S$-scheme
$(f : T \to S)$ to~$\Gamma(T, f^*\F)$.\label{page:induced-sheaf-on-zar}

\subsection*{A first example illustrating the Kripke--Joyal translation rules}
Since all the sets~$\affl(T) \cong \Gamma(T,\O_T)$ carry ring structures and do so in a compatible way, the
object~$\affl$ can be endowed with a canonical structure as a ring object
in~$\Zar(S)$. For a particular~$S$-scheme~$T$, the ring~$\affl(T)$ will almost
never be a field, but the system of these rings, conceptualized as a single
entity from the internal point of view, does satisfy a field axiom. In the
case~$S = \Spec\ZZ$, this was first observed by
Kock~\cite{kock:univ-proj-geometry}.

\begin{prop}\label{prop:affl-field-informal}
The ring~$\affl$ is a field from the internal point of view
of~$\Zar(S)$, in the sense that
\[ \Zar(S) \models
  \forall f\?\affl\_
    \neg(f = 0) \Rightarrow \speak{$f$ \inv}. \]
\end{prop}

\begin{proof}According to the Kripke--Joyal semantics of~$\Zar(S)$, we
have to show that for any~$S$-scheme~$T$ and any function~$f \in
\Gamma(T,\O_T)$ the statement~$T \models \neg(f = 0)$ implies~$T \models
\speak{$f$ \inv}$. The antecedent states that, for any~$T$-scheme~$T'$, if the
pullback of~$f$ to~$T'$ vanishes, then~$T'$ is the empty scheme.

As with the analogous statement about the little Zariski topos
(Lemma~\ref{lemma:internal-invertibility}), the consequent means that~$f$ is
invertible in~$\Gamma(T,\O_T)$.

The claim follows by considering the particular~$T$-scheme~$T' \defeq V(f)$.
Since~$f$ vanishes on~$V(f)$, this subscheme is empty and therefore its
complement~$D(f)$ is all of~$T$.
\end{proof}

The field property can be interpreted as follows. A function~$f$ not being the
zero function does not imply that it's invertible. But if~$f$ is
\emph{universally nonzero} in that the only scheme such that pullback of~$f$ to
that scheme vanishes is the empty scheme, then~$f$ is indeed invertible.

We'll revisit the field property in Section~\ref{sect:special-properties-affl};
it turns out that it has a deeper reason than the manual proof given here showed.

}

\section{On the proper choice of a big Zariski site}
\label{sect:proper-choice-of-site}

Unlike with the construction of the little Zariski topos, set-theoretical
issues of size arise when constructing the big Zariski topos. These can be
solved in several different manners, yielding toposes which are not equivalent,
and actually differ in some important aspects, but otherwise enjoy very similar
properties.

{\tocless

\subsection*{Naive approach}
Some authors construct the big Zariski topos of~$S$ as the topos of
sheaves over the site~$\Sch/S$ of all schemes over~$S$. This option is quite
attractive since the Yoneda embedding~$\Sch/S \to \Sh(\Sch/S)$, which sends
an~$S$-scheme to its functor of points, is fully faithful, therefore the
internal language of~$\Sh(\Sch/S)$ can distinguish arbitrary schemes.

However, because~$\Sch/S$ is not essentially small, forming the sheaf topos is
not possible in plain Zermelo--Fraenkel set theory.

Since it's still possible to meaningfully speak of individual
functors~$(\Sch/S)^\op \to S$, we can attach a Kripke--Joyal semantics
to~$\Sh(\Sch/S)$, as long as we keep in mind that~$\Sh(\Sch/S)$ might
not contain a subobject classifier and might not be cartesian closed. From the
internal point of view, powersets and function sets might therefore not exist.

\subsection*{Using Grothendieck universes}
We could also assume the existence of a Gro\-then\-dieck universe~$\U$
containing~$S$ and construct~$\Zar(S)$ as the topos of sheaves over the
small site~$\Sch_\U/S$, the category of~$S$-schemes contained in~$U$.

By the \emph{comparison lemma} (see, for instance,~\cite[Theorem~3.7]{caramello:preliminaries}), we could also construct~$\Zar(S)$ as
the topos of sheaves over~$\Aff_\U/S$, the category of~$S$-schemes in~$U$ which
are affine (as absolute schemes), and obtain an equivalent topos.

In this case, the Yoneda functor~$\Sch/S \to \Zar(S)$ might not be faithful,
but the restricted Yoneda functor~$\Sch_\U/S \to \Zar(S)$ will.

\subsection*{Approach of the Stacks Project}
The Stacks Project proposes a more nuanced approach, namely expanding a given
set~$M_0$ of schemes containing~$S$ to a superset~$M$ which is closed (up to
isomorphism) under several constructions~\cite[Tag~000H]{stacks-project}: fiber
products, countable coproducts, domains of open and closed immersions and of
morphisms of finite type, spectra of local rings~$\O_{X,x}$, spectra of residue
fields, and others.

The Stacks Project then defines~$\Zar(S)$ as~$\Sh(\Sch_M/S)$, where~$\Sch_M/S$
is the small category of~$S$-schemes in~$M$, or equivalently
as~$\Sh(\Aff_M/S)$. This approach has the advantage that one doesn't have to
assume the existence of a Grothendieck universe; the \emph{partial
universe}~$M$ can be constructed entirely within ZFC set theory using
transfinite recursion.

\subsection*{Employing parsimonious sites}
From a topos-theoretical point of view, it's natural to settle for an even
more parsimonious site: the site~$(\Sch/S)_\lfp$ consisting of the~$S$-schemes
which are locally of finite presentation over~$S$, or equivalently the
essentially small site~$(\Aff/S)_\lfp$ of the~$S$-schemes which are locally of
finite presentation over~$S$ and affine (as absolute schemes).\footnote{It's
not reasonable to restrict to the even smaller site consisting of the finitely
presented~$S$-schemes, since open immersions can fail to be finitely presented.
We want the site used to construct~$\Zar(S)$ to be closed under domains of open immersions,
for instance to facilitate a comparison with the little
Zariski topos~$\Sh(S)$, whose site does contain all open subsets of~$S$.
Furthermore, since a finitely presented~$S$-scheme might not admit an open
covering by finitely presented~$S$-schemes which are affine (as absolute
schemes), the toposes~$\Sh((\Sch/S)_\fp)$ and~$\Sh((\Aff/S)_\fp)$ can differ.}

In the special case that~$S = \Spec(A)$ is affine, this site is the
dual of the category of finitely presented~$A$-algebras; in this case the
topos-theoretic points of the resulting topos are precisely the local~$A$-algebras,
and moreover, the resulting topos is the classifying topos of the theory of
local~$A$-algebras, such that for any Grothendieck topos~$\E$, geometric
morphisms~$\E \to \Sh((\Aff/S)_\lfp)$ correspond to local~$A$-algebras internal
to~$\E$. A textbook reference for these facts
is~\cite[Section~VIII.6]{moerdijk-maclane:sheaves-logic}.

In contrast, the toposes arising when using the larger sites have categories of
points which contain further objects in addition to all local~$A$-algebras; and
no simple description of the theory they classify is known.\footnote{The
category of points of a presheaf topos~$[\C^\op,\Set]$ coincides
with~$\Ind(\C^\op)$, the ind-completion of~$\C^\op$. This general fact explains that
in the case that~$\C$ is the category of finitely presented~$A$-algebras,
the category of points coincides with the category of~$A$-algebras,
since~$\Ind(\Alg(A)_\fp) \simeq \Alg(A)$. For the larger sites, understanding
the structure of their points is therefore tantamount to understanding the
structure of the ind-completion of their dual category (and understanding which
points of the presheaf topos are also points of the sheaf topos).}

A further advantage of these parsimonious sites is that they don't require arbitrary
choices of a starting set~$M_0$ or of a way of expanding~$M_0$ to a sufficiently
ample set~$M$ of schemes.

However, the parsimonious sites also have a serious disadvantage, namely that
with them, the Yoneda functor is only fully faithful when restricted
to~$(\Sch/S)_\lfp$. For instance, in the case~$S = \Spec(\ZZ)$, the
schemes~$\Spec(\QQ)$ and the empty scheme have isomorphic functors of points by Proposition~\ref{prop:fingen-algebra-q}.
Therefore~$\Spec(\QQ)$ looks like the empty set from the internal point of
view.

In the following, we do not commit to a single one of these options for
resolving the set-theoretical size issues, but rather keep any of them
in mind. This approach will sometimes necessitate phrases such as ``for
any~$S$-scheme~$T$ contained in the site used to define~$\Zar(S)$'', which might
seem awkward to a topos-theorist when taken out of context, since
the site used to construct a Grothendieck topos is not at all uniquely
determined by the resulting topos.

We will indicate the few places where the choice of site makes a difference.
When the definition of the Kripke--Joyal semantics for~$\Zar(S)$ refers
to~$S$-schemes, it actually refers only to the~$S$-schemes contained in the
site. Similarly, one has to restrict oneself to such schemes in the statement
of Proposition~\ref{prop:basic-properties-language-big}.
Proposition~\ref{prop:affl-field-informal} holds for any choice of site.

It's possible to define the big Zariski topos of a scheme without recourse to
classical scheme theory; we discuss this in
Section~\ref{sect:big-zariski-without-classical-scheme-theory}.

\begin{rem}Some authors define the big Zariski topos of~$S$ as the topos of
sheaves over the category of affine~$S$-schemes (that is, $S$-schemes~$f : X
\to S$ where~$f$ is affine) contained in some universe. In case that the
diagonal morphism~$S \to S \times S$ is affine, the resulting topos is
equivalent to what we regard as the big Zariski topos of~$S$, when defined using the
same universe. This is because in this case morphisms of the form~$\Spec(A) \to
S$ are affine and affine morphisms with codomain~$S$ can be refined by such
morphisms.
\end{rem}

}

\section{Relation between the big and little Zariski toposes}
\label{sect:relation-big-little}

The big Zariski topos~$\Zar(S)$ is a topos over the little Zariski
topos~$\Sh(S)$ in that there is a canonical geometric morphism
\[ \pi : \Zar(S) \longrightarrow \Sh(S) \]
with direct and inverse image parts given by
\[ \pi_*E = E|_{\Sh(S)} \qquad\text{and}\qquad
  \pi^{-1}\F = ((T \xra{f} S) \mapsto \Gamma(T, f^{-1}\F)). \]
Since~$\pi^{-1}$ is fully faithful, this geometric morphism is connected; and
furthermore, it is a local geometric morphism (a further right adjoint~$\pi^!$
which is fully faithful exists).

By general results on local geometric morphisms, the adjoint pair~$(\pi_*
\dashv \pi^!)$ is a geometric morphism which is right inverse to~$\pi$ and
which exhibits~$\Sh(S)$ as a subtopos of~$\Zar(S)$, similarly to how~$\Set$ is
a subtopos of a sheaf topos over a local topological space. In this context,
it's customary to introduce notation for the idempotent monad~$\sharp$ and the
idempotent comonad~$\flat$ arising from the adjoint triple~$\pi^{-1} \dashv
\pi_* \dashv \pi^!$:
\[ \sharp E = \pi^!(E|_{\Sh(S)}) \qquad\text{and}\qquad
  \flat E = \pi^{-1}(E|_{\Sh(S)}). \]

In the case that~$S = \Spec(A)$ is an affine scheme and we employ
one of the parsimonious sites to construct~$\Zar(S)$, it's well-known
that~$\Sh(S)$ classifies local localizations of~$A$ and that~$\Zar(S)$ classifies
arbitrary local~$A$-algebras. On points, the morphism~$\pi$ sends a
local~$A$-algebra~$\varphi : A \to R$ to the local localization~$A \to
A[(\varphi^{-1}[R^\times])^{-1}]$, and its right inverse sends a local
localization~$A \to A[F^{-1}]$ to itself.

\subsection{Recovering the big Zariski topos from the little Zariski topos}

What does~$\Zar(S)$ classify in the case that~$S$ is an arbitrary scheme?
We don't know a nontautologous answer to this question, but we can answer a
related one: What does~$\Zar(S)$ classify as seen from the internal point of
view of~$\Sh(S)$?

To make sense of this question, we employ a slight extension of Shulman's
stack semantics which allows to refer to locally internal
categories~\cite{penon:locally-internal-categories} over a base topos~$\E$ from
the internal language. Using this extension, a locally internal category
over~$\E$ looks like a locally small category from the internal point of view
of~$\E$. In particular, a geometric morphism~$f : \F \to \E$ gives rise to a
locally internal category (which over an object~$A \in \E$ is given by
the~$\E/A$-enriched category~$\F/f^{-1}A$) which will look like an ordinary
topos from the internal point of view of~$\E$.

For instance, the trivial~$\E$-topos~$\E$ will look like~$\Set$ and the
slice topos~$\E/X$ will look like~$\Set/X$ from the internal point of
view of~$\E$.

\begin{thm}\label{thm:zar-classifies}
In the situation that the site used to construct~$\Zar(S)$ is one of the
parsimonious sites, the big Zariski topos~$\Zar(S)$ is, from the internal point
of view of~$\Sh(S)$, the classifying topos of
the theory of local~$\O_S$-algebras which are local over~$\O_S$.
\end{thm}

For an arbitrary topos~$\F$ over~$\Set$, the concept of an~``$\O_S$-algebra
in~$\F$'' doesn't make any sense -- in contrast to the concept of
an~$A$-algebra in~$\F$, which can either be defined as a ring
homomorphism~$\ul{A} \to R$ in~$\F$ (where~$\ul{A}$ is the pullback of~$A \in
\Set$ to~$\F$) or as a ring object which is equipped with an~$A$-indexed family of
endomorphisms satisfying suitable axioms. However, for a~$\Sh(S)$-topos~$f : \F \to
\Sh(S)$, the concept of an~$\O_S$-algebra in~$\F$ is meaningful: It's a ring
homomorphism~$f^{-1}\O_S \to R$ in~$\F$.

Similarly, there is no absolute ``geometric theory of~$\O_S$-algebras''.
However, there is a geometric theory of~$\O_S$-algebras \emph{internal
to~$\Sh(S)$}. Theorem~\ref{thm:zar-classifies} should be viewed in this light.
A detailed account of internal geometric theories and internal classifying toposes
is given in~\cite[Chapter~II]{henry:classifying-topoi}.

The proviso ``local over~$\O_S$'' is as in the discussion of the relative
spectrum from the internal point of view
(Section~\ref{sect:relative-spectrum}).

\begin{proof}[Proof of Theorem~\ref{thm:zar-classifies}]
We have to verify that, from the point of view of~$\Sh(S)$, the topos~$\Zar(S)$
contains a canonical local and local-over-$\O_S$ $\O_S$-algebra and that for
any Grothendieck topos~$\F$, pulling back this canonical algebra yields an
equivalence between the category of geometric morphisms~$\F \to \Zar(S)$ and
the category of local and local-over-$\O_S$ $\O_S$-algebras in~$\F$.

The canonical local and local-over-$\O_S$ $\O_S$-algebra in~$\Zar(S)$ is
the algebra~$\flat\affl \to \affl$. Indeed, the ring~$\affl$ is local and
the homomorphism~$\flat\affl \to \affl$ is local, since its restriction to
any sheaf topos~$\Sh(X)$, where~$f : X \to S$ is an~$S$-scheme contained in the site used
to define~$\Zar(S)$, is local: It's the morphism $f^\sharp : f^{-1}\O_S \to \O_X$.

We now want to verify the universal property, which expressed internally
to~$\Sh(S)$ reads as
\[ \Hom(\E, \Zar(S)) \simeq
  \text{category of local and local-over-$\O_S$ $\O_S$-algebras in~$\E$}. \]
Externally, this means that for any open subset~$U \subseteq S$ and any
topos~$\E$ over~$\Sh(S)/\ul{U}$,
\begin{multline*}
  \qquad\Hom_{\Sh(S)/\ul{U}}(\E, \Zar(S)/\pi^{-1}\ul{U}) \simeq \\
  \text{category of local and local-over-$\pi^{-1}\O_S$ $\pi^{-1}\O_S$-algebras
  in~$\E$}.\qquad
\end{multline*}
We will verify this equivalence in the case that~$S = \Spec(A)$ is affine and
that~$S = U$. This suffices to establish the theorem, since~$\Sh(S)/\ul{U} \simeq
\Sh(U)$,~$\Zar(S)/\pi^{-1}\ul{U} \simeq \Zar(U)$, and since the internal
language is local.

So let~$f : \E \to \Sh(\Spec(A))$ be a~$\Sh(\Spec(A))$-topos. By the universal
property of~$\Zar(\Spec(A))$ as the classifying topos of local~$A$-algebras, a geometric
morphism~$g : \E \to \Zar(\Spec(A))$ is uniquely determined by a
local~$A$-algebra~$\varphi : \ul{A} \to \B$ in~$\E$. By the universal property
of~$\Sh(\Spec(A))$ as the classifying topos of local localizations of~$A$, the
composition~$\pi \circ g : \E \to \Sh(\Spec(A))$ is uniquely determined by
the local localization~$\ul{A} \to g^{-1} \pi^{-1} \O_{\Spec(A)} = g^{-1}(\flat
\affl)$ in~$\E$.

In the composition
\[ \ul{A} \lra \flat \affl \lra \affl, \]
the first morphism is a local localization and the second morphism is local.
Since these properties can be formulated as geometric
implications,\footnote{A ring homomorphism~$\alpha : R \to T$ is a localization
(that is, isomorphic to the canonical localization morphism~$R \to R[S^{-1}]$
for some multiplicative subset~$S$) if and only if the canonical comparison
morphism~$R[(\alpha^{-1}T^\times)^{-1}] \to T$ is bijective. This is the case
if and only if
\begin{multline*}
  \quad\qquad\forall y\?T\_
  \exists x\?R\_ \exists s\?R\_
  \speak{$\alpha(s)$ \inv} \wedge \alpha(s) y = x
  \qquad\text{and} \\
  \forall x\?R\_ \alpha(x) = 0 \Rightarrow
  \exists s\?R\_ \speak{$\alpha(s)$ \inv} \wedge sx = 0.\qquad\quad
\end{multline*}}
they are preserved by the functor~$g^{-1}$. Since
furthermore such a factorization is unique, the localization~$\ul{A} \to
g^{-1}(\flat \affl)$ which determines~$\pi \circ g$ coincides with the
localization~$\ul{A}[(\varphi^{-1}[\B^\times])^{-1}]$. Referring directly
to the involved filters, the filter~$f^{-1}\F$ which determines~$\pi \circ g$
(where~$\F$ is the generic filter of~$\ul{A}$ in~$\Sh(\Spec(A))$) coincides
with the filter~$\varphi^{-1}[\B^\times]$. This explains the first
equivalence in the chain
\begin{align*}
  &\mathrel{\phantom{\simeq}} \Hom_{\Sh(\Spec(A))}(\E, \Zar(\Spec(A))) \\
  &\simeq
    \text{category of local algebras~$\varphi : \ul{A} \to \B$ in~$\E$
  such that $\varphi^{-1}\B^\times = f^{-1}\F$} \\
  &\simeq
    \text{category of local algebras~$\psi : f^{-1}\O_{\Spec(A)} \to \B$
    in~$\E$ such that~$\psi$ is local}.
\end{align*}
The second equivalence maps an algebra~$\varphi$
to~$\ul{A}[(\varphi^{-1}\B^\times)^{-1}] \to \B$; conversely, an algebra~$\psi$
is mapped to the composition~$\ul{A} \to f^{-1}\O_{\Spec(A)} \xra{\psi} \B$.
\end{proof}

Similarly to how Theorem~\ref{thm:zar-classifies} shows how the big Zariski
topos of~$S$ looks like from the point of view of~$\Sh(S)$, it's possible to
give an internal description of what the big Zariski topos of an arbitrary
relative spectrum over~$S$ looks like. We state and verify such a description
in Theorem~\ref{thm:big-zariski-topos-of-relative-spectrum}.

It is well-known that the points of~$\Zar(\Spec(R))$, when constructed using
one of the parsimonious sites, are in canonical bijection with the
local~$R$-algebras; for instance, this follows from the description
of~$\Zar(\Spec(R))$ as the classifying topos of the theory of
local~$R$-algebras. For the case of a general base scheme, we introduce the
following definition.

\begin{defn}\label{defn:ring-over-s}A \emph{ring over~$S$} is a ring~$A$ together with a
morphism~$\Spec(A) \to S$ of locally ringed locales. A \emph{morphism of rings
over~$S$} is a ring homomorphism which is compatible with the structure
morphisms to~$S$.\end{defn}
%

In the special case that~$S = \Spec(B)$ is an affine scheme, a ring over~$S$ in
the sense of Definition~\ref{defn:ring-over-s} is the same as an~$B$-algebra.

\begin{prop}\label{prop:points-of-big-zariski}
In the situation that one of the parsimonious sites is used to
define~$\Zar(S)$, the category of points of~$\Zar(S)$ is canonically equivalent
to the full subcategory of the rings over~$S$ whose underlying ring is
local.\end{prop}

\begin{proof}By Theorem~\ref{thm:zar-classifies}, a point of~$\Zar(S)$ is given
by a point of~$\Sh(S)$, that is by a point~$s$ of~$S$,\footnote{The
topos-theoretic points of the topos of sheaves over a topological space~$T$ are
in canonical bijection with the locale-theoretic points of~$T$, that is with
locale morphisms~$1 \to T$. If~$T$ is sober, such points are in canonical
bijection with the elements of the underlying set of~$T$. In a classical
metatheory, schemes are sober~\stacksproject{01IS}. If one wants the proof to
work intuitionistically, the base scheme~$S$ has to be defined in a
intuitionistically sensible way, for instance as a locally ringed locale.
Correspondingly, the point~$s$ of~$S$ has to be interpreted in the
locale-theoretic sense.} together with a local~$\O_{S,s}$-algebra~$A$ which is
local over~$\O_{S,s}$. These data define a ring over~$S$, namely the ring~$A$
together with the composite~$\Spec(A) \to \Spec(\O_{S,s}) \to S$.
Since the structure morphism~$\O_{S,s} \to A$ is local, this composite maps the
focal point of~$\Spec(A)$ to the given point~$s \in S$.

Conversely, let a local ring~$A$ together with a morphism~$f : \Spec(A) \to S$ of
locally ringed locales be given. Let~$x \in \Spec(A)$ be the focal point
of~$\Spec(A)$. We set~$s \defeq f(x)$; then~$A$ is an~$\O_{S,s}$-algebra
by~$(f^\sharp)_x$. It is local over~$\O_{S,s}$ since~$f^\sharp$ is a local
homomorphism.

These constructions are mutually inverse since the morphisms $\Spec(\O_{S,s})
\to S$ are monomorphisms in the category of locally ringed locales.
\end{proof}

\begin{rem}\label{rem:zar-classifies-absolute}
In the situation that one of the parsimonious sites is used to
define the big Zariski topos of~$S$, it classifies the theory of local rings
over~$S$. This is a restatement of Theorem~\ref{thm:zar-classifies}.
Explicitly, the theory of local rings over~$S$ is given by:
\begin{enumerate}
\item A theory which~$\Sh(S)$ classifies.
\item Structure and axioms for a ring~$R$.
\item Structure and axioms which guarantee that the interpretation of~$R$ in
any cocomplete topos coincides with the pullback of~$\O_S$.
\item Structure and axioms for a local ring~$A$ and a local homomorphism~$R \to A$.
\end{enumerate}
The fourth item can be substituted by:
\begin{enumerate}
\item[(4')] Structure and axioms for a local ring~$A$ and a morphism~$\Spec(A)
\to (\pt,R)$ of locally ringed locales.
\end{enumerate}
This is because such a morphism is given by a local homomorphism~$\ul{R}
\to \O_{\Spec(A)}$ of sheaves of rings which in turn is given by a local ring
homomorphism~$R \to \Gamma(\Spec(A), \O_{\Spec(A)}) = A$. (Taking global
sections of a local homomorphism of sheaves of rings yields a homomorphism of
rings which will typically fail to be local. However, here taking global
sections coincides with calculating the stalk at the focal point of~$\Spec(A)$,
and pullback preserves locality of ring homomorphisms.)
\end{rem}

Corollary~\ref{cor:pp1-classifies} gives a description of the theory which the
big Zariski topos of~$\PP^1_\ZZ$ classifies, building upon
Remark~\ref{rem:zar-classifies-absolute}.

\subsection{Recovering the little Zariski topos from the big Zariski topos}

Theorem~\ref{thm:zar-classifies} shows that~$\Zar(S)$ can be reconstructed
from~$\Sh(S)$ (and its structure sheaf~$\O_S$). Similarly, it's possible to
reconstruct~$\Sh(S)$ from~$\Zar(S)$ (and the canonical morphism~$\flat \affl
\to \affl$).

\begin{thm}\label{thm:reconstruct-little-topos}
In the situation that the site used to construct~$\Zar(S)$ is one of the
parsimonious sites, the little Zariski topos~$\Sh(S)$ is the largest subtopos
of~$\Zar(S)$ where the canonical morphism~$\flat \affl \to \affl$ is an
isomorphism.
\end{thm}

In other words, the little Zariski topos is the largest subtopos~$\E
\hookrightarrow \Zar(S)$ such that $\Zar(S) \models (\speak{$\flat\affl \to
\affl$ is bijective})^\Box$ (where~$\Box$ is the modal operator corresponding
to the subtopos), that is that the pullback of the canonical
morphism~$\flat\affl \to \affl$ to~$\E$ is an isomorphism.

In the case that~$S
= \Spec(A)$ is affine, we also have the ring~$\ul{A}$ in~$\Zar(S)$ available.
In this case the condition is equivalent to
\[ \Zar(S) \models \speak{$\ul{A} \to \affl$ is a localization}^\Box, \]
since in the composition~$\ul{A} \to \flat\affl \to \affl$ the first morphism
is a localization.

\begin{proof}[Proof of Theorem~\ref{thm:reconstruct-little-topos}]
The little Zariski topos is a subtopos of the big Zariski topos via the right
inverse~$s$ of~$\pi : \Zar(S) \to \Sh(S)$, the geometric morphism~$(\pi_* \dashv
\pi^!)$. The pullback of~$\flat\affl \to \affl$ to~$\Sh(S)$ is therefore the
morphism~$(\flat\affl)|_{\Sh(S)} \to \affl|_{\Sh(S)}$, that is~$\O_S \to \O_S$,
which is an isomorphism.

Let~$f : \E \hookrightarrow \Zar(S)$ be any subtopos such that the pullback
of~$\flat\affl \to \affl$ to~$\E$ is an isomorphism. We want to verify that~$f$
factors over the inclusion~$s : \Sh(S) \hookrightarrow \Zar(S)$.
\[ \xymatrix{
  \E \ar@{^{(}->}[rr]^f \ar@{-->}[rd] && \Zar(S) \\
  & \Sh(S) \ar@{^{(}->}[ru]_s
} \]
A candidate for a morphism~$\E \to \Sh(S)$ witnessing this factorization is the
composite~$\pi \circ f$. It remains to show that~$s \circ (\pi \circ f) = f$.
Both~$s \circ (\pi \circ f)$ and~$f$ are morphisms of~$\Sh(S)$-toposes,
where~$\E$ is regarded as a~$\Sh(S)$-topos by the composition~$\pi \circ f$.
By the universal property of the big Zariski topos given in
Theorem~\ref{thm:zar-classifies}, they are therefore uniquely determined by
the~$\O_S$-algebra they classify.

The morphism~$s \circ (\pi \circ f)$ classifies
the~$\O_S$-algebra~$f^{-1}\pi^{-1}s^{-1}\affl = f^{-1}(\flat\affl)$. The
morphism~$f$ classifies the~$\O_S$-algebra~$f^{-1}\affl$.
Since~$f^{-1}(\flat\affl) \to f^{-1}\affl$ is an isomorphism, these algebras
coincide.
\end{proof}

\subsection{Change of base}
\label{sect:change-of-base}

Let~$f : X \to S$ be a morphism of schemes. In any of the situations that
\begin{enumerate}
\item the parsimonious sites are used to construct the big Zariski
toposes and~$f$ is locally of finite presentation, or
\item the same (Grothendieck or partial) universe is used for constructing both
Zariski toposes and both~$X$ and~$S$ are contained in the universe,
\end{enumerate}
the morphism~$f$ induces an essential geometric morphism~$\Zar(X) \to \Zar(S)$
which we again denote by~``$f$''. Explicitly, the big Zariski toposes are
related by the adjoint triple~$f_! \dashv f^{-1} \dashv f_*$ with
\begin{align*}
  f_* : \Zar(X) \lra \Zar(S),\ F &\longmapsto ((T \xra{g} S) \mapsto F(T \times_S X)), \\
  f^{-1} : \Zar(S) \lra \Zar(X),\ E &\longmapsto ((T \xra{g} X) \mapsto E(T \xra{g} X \xra{f} S)), \\
  f_! : \Zar(X) \lra \Zar(S),\ F &\longmapsto ((T \xra{g} S) \mapsto \coprod_{h : T \to X} F(T \xra{h} X)).
\end{align*}
In situation~(2), the well-definedness of these functors is trivial. In
situation~(1), the well-definedness rests on the lemma that an~$S$-morphism~$h : T \to X$ is locally of
finite presentation if~$T$ and~$X$ are locally of finite presentation over~$S$~\stacksproject{02FV}.

The objects of~$\Zar(S)$ listed on page~\pageref{page:important-objects} pull back
as expected:
\begin{itemize}
\item Let~$Y$ be an~$S$-scheme. Then~$f^{-1} \ul{Y} = \ul{Y \times_S X}$, by
the universal property of the fiber product.
\item In particular, $f^{-1} \affl = \afflx$, since~$\AA^1_S \times_S X =
\AA^1_X$.
\item Let~$\F$ be a sheaf of sets on~$S$. Then~$f^{-1} \pi_S^{-1} \F =
\pi_X^{-1} f^{-1} \F$.
\item Let~$\F$ be a sheaf of~$\O_S$-modules. Then~$f^{-1} \F^\Zar =
(f^* \F)^\Zar$.
\end{itemize}

The functors~$f_! \dashv f^{-1}$ induce an equivalence
\[ \Zar(X) \simeq \Zar(S)/\ul{X}, \]
explicitly described by
\[ \begin{array}{r@{}c@{}l}
  F &{}\longmapsto{}& (f_!F \to f_!1), \\
  ((T \xra{g} X) \mapsto \{ s \in (f^{-1}E)(T) \,|\, \alpha(s) = g \}) &{}\longmapsfrom{}& (E
  \xra{\alpha} \ul{X}).
\end{array} \]
From the internal point of view of~$\Zar(S)$, the big Zariski topos of~$X$ is
therefore simply~$\Set/\ul{X}$, the category of~$\ul{X}$-indexed families of
sets or equivalently the category of sheaves on~$\ul{X}$ considered as a
\emph{discrete} locale. This fits nicely with the philosophy that~``$S$-schemes are
plain unstructured sets from the internal point of view
of~$\Zar(S)$''.\footnote{The requirement on~$f$ mentioned at the beginning of
this subsection is really necessary. For instance, let~$f$ be the unique
morphism~$\Spec(\QQ) \to \Spec(\ZZ)$. This morphism is not locally of finite
presentation. By Proposition~\ref{prop:fingen-algebra-q}, if the parsimonious
sites are used, the functor of
points~$\ul{\QQ} \in \Zar(\Spec(\ZZ))$ coincides with the functor of points of
the empty scheme and is thus the initial
object in~$\Zar(S)$. Therefore~$\Zar(\Spec(\ZZ))/\ul{\QQ}$ is the trivial
topos. In contrast,~$\Zar(\Spec(\QQ))$ is not.}

In contrast, for the little Zariski toposes, there is no similarly simple
description of the little Zariski topos of~$X$ as a slice of the little Zariski
topos of~$S$. From the internal point of view of~$\Sh(S)$, the topos~$\Sh(X)$
looks like the topos of sheaves over a locale which is not discrete, and the
topos~$\Zar(X)$ doesn't even look like a topos of sheaves over an arbitrary
locale (discrete or not).

The internal language of a slice topos~$\E/I$ admits a
simple description from the point of view of~$\E$. Namely, for
any formula~$\varphi$ over~$\E/I$,
\[ \E/I \models \varphi \qquad\text{iff}\qquad
  \E \models \forall i\?I\_ \varphi(i). \]
For the right hand side to make sense, it has to be interpreted in the
following way. Any object~$(p : M \to I)$ of~$\E/I$ which appears in~$\varphi$,
for instance as a domain of quantification, has to be substituted by the
internal expression~``$p^{-1}[\{i\}]$'' denoting the fiber of~$p$
over~$i\?I$.\footnote{This substitution is less ad~hoc as it might at first
appear. The internal language of a topos~$\E$ is \emph{dependently typed}, meaning
that the types one can quantify over may depend on previously introduced
values. Types in the empty context, depending on no values, correspond to
objects of~$\E$. Types in the context of a variable~$i\?I$ correspond to
objects~$(p : M \to I)$ of~$\E/I$. For instance, in this case one can form
formulas of the form~``$\forall i\?I\_ \forall m\?M(i)\_ \psi(i,m)$''. If in
the translation process using the Kripke--Joyal semantics a formal variable~$i$
was substituted by a generalized element~$i_0 : A \to I$, the
expression~``$M(i_0)$'' has to be interpreted as the pullback~$i_0^* M$.}
For example, if~$(p : M \to I)$ is such an object of~$\E/I$,
\begin{align*}
  &\phantom{\text{iff }} \E/I \models \speak{$M$ is inhabited} \\
  &\text{iff } \phantom{/I}\E \models \forall i\?I\_ \speak{the fiber of $p$ over $i$ is inhabited} \\
  &\text{iff } \phantom{/I}\E \models \forall i\?I\_ \exists m\?M\_ p(m) = i.
\end{align*}

Thanks to this description of the internal language of a slice topos, the
equivalence~$\Zar(X) \simeq \Zar(S)/\ul{X}$ is useful for lifting internal
characterizations concerning properties of~$S$-schemes to properties of
morphisms of~$S$-schemes. For instance, we will see in
Proposition~\ref{prop:char-surjective-morphisms} that the structure morphism of
an~$S$-scheme~$f : Y \to S$ is surjective if and only if
$\Zar(S) \models \neg\neg(\speak{$\ul{Y}$ is inhabited})$.
This automatically implies (Corollary~\ref{cor:char-surjective-morphisms-relative})
that a morphism~$p : Y \to X$ of~$S$-schemes is surjective if and only if
\[ \Zar(S) \models \forall x\?\ul{X}\_ \neg\neg(\speak{the fiber of~$\ul{p}$
over~$x$ is inhabited}). \]

Many properties of morphisms in algebraic geometry, and all properties which
can be characterized using the internal language of the big Zariski topos, are
stable under base change. For those kinds of properties~$P$, if a morphism~$Y
\to X$ is~$P$, then for any point~$x \in X$ the base change~$Y_x \to
\Spec(k(x))$ along~$\Spec(k(x)) \to X$ is~$P$ as well. The converse is usually
false, but the motto~``a morphism is~$P$ if all its fibers are~$P$ in a
continuous fashion'' is still useful for intuition. The equivalence~$\Zar(X)
\simeq \Zar(S)/\ul{X}$ makes this motto precise: For any morphism~$p : Y \to X$
of~$S$-schemes and any formula~$\varphi(M)$ of~$\Zar(S)$ containing a free
variable~$M$,
\[ \Zar(X) \models \varphi(\ul{Y}) \qquad\text{iff}\qquad
  \Zar(S) \models \forall x\?\ul{X}\_ \varphi(\ul{p}^{-1}[\{x\}]), \]
that is $Y$ has property~$\varphi$ when regarded as an~$X$-scheme if and only
if all the fibers of~$Y \to X$ have property~$\varphi$ when regarded
as~$S$-schemes.

The family of geometric morphisms
\[ \Sh(X) \to \Zar(X) \to \Zar(S) \]
with inverse image~$E \mapsto E|_{\Sh(X)}$, where~$X$ ranges over
all~$S$-schemes contained in the site used to define~$\Zar(S)$, is jointly
surjective. That is, the restriction functors~$\Zar(S) \to \Sh(X)$ are jointly
conservative. This fact, together with the fact that these functors commute
with geometric constructions, is frequently useful to relate truth in~$\Zar(S)$
with truth in all of the~$\Sh(X)$.


\begin{lemma}\label{lemma:zar-tensor-product-commutes}
\begin{enumerate}
\item The functor~$\pi^{-1}$ from~$\O_S$-modules to~$\pi^{-1}\O_S$-modules commutes
with tensor product.
\item The functor~$\pi_*$ from~$\affl$-modules to~$\O_S$-modules commutes with
tensor product.
\item The~$\affl$-module~$\E^\Zar$ associated to an~$\O_S$-module~$\E$ (defined
on page~\pageref{page:induced-sheaf-on-zar}) is canonically isomorphic
to~$\pi^{-1}\E \otimes_{\pi^{-1}\O_S} \affl$.
\item The functor~$\E \mapsto \E^\Zar$ from~$\O_S$-modules to~$\affl$-modules
commutes with tensor product.
\end{enumerate}
\end{lemma}

\begin{proof}Forming the tensor product is a geometric construction. It is
therefore preserved by~$\pi^{-1}$ and by~$\pi_*$, since both functors are
inverse-image-parts of geometric morphisms. Claim~(3) follows from the
calculation
\begin{align*}
  (\pi^{-1}\E \otimes_{\pi^{-1}\O_S} \affl)|_{\Sh(X)}
  &\cong (\pi^{-1}\E)|_{\Sh(X)} \otimes_{(\pi^{-1}\O_S)|_{\Sh(X)}} \affl|_{\Sh(X)} \\
  &\cong (f^{-1}\E) \otimes_{f^{-1}\O_S} \O_X \\
  &\cong f^*\E
\end{align*}
for~$S$-schemes~$(f : X \to S)$ contained in the site used to define~$\Zar(S)$.
Finally, claim~(4) follows by using a standard property for tensor products
(whose proof is intuitionistic and therefore valid in~$\Zar(S)$):
\begin{align*}
  (\E \otimes_{\O_S} \F)^\Zar
  &\cong \pi^{-1}(\E \otimes_\S \F) \otimes_{\pi^{-1}\O_S} \affl \\
  &\cong (\pi^{-1}\E \otimes_{\pi^{-1}\O_S} \pi^{-1}\F) \otimes_{\pi^{-1}\O_S} \affl \\
  &\cong (\pi^{-1}\E \otimes_{\pi^{-1}\O_S} \affl) \otimes_\affl (\pi^{-1}\F \otimes_{\pi^{-1}\O_S} \affl) \\
  &\cong \E^\Zar \otimes_\affl \F^\Zar. \qedhere
\end{align*}
\end{proof}

\begin{ex}Constructing, internally in~$\Zar(S)$, the
module~$\Omega^1_{\affl|\flat\affl}$ of Kähler differentials of the ring
morphism~$\flat\affl \to \affl$ yields the ``universal cotangent sheaf''
\[ X \longmapsto \Gamma(X, \Omega^1_{X|S}). \]
This is because constructing the module of Kähler differentials is a geometric
construction and the restriction functors~$\Zar(S) \to \Zar(X)$ commute with
geometric constructions, so
\[ \brak{\Omega^1_{\affl|\flat\affl}}|_{\Sh(X)} \cong
  \Omega^1_{(\affl|_{\Sh(X)}) | ((\flat\affl)|_{\Sh(X)})} \cong
  \Omega^1_{\O_X | f^{-1}\O_S} \cong \Omega^1_{X|S} \]
for any~$S$-scheme~$(f : X \to S)$ contained in the site used to
define~$\Zar(S)$.

This universal cotangent sheaf doesn't enjoy any finiteness properties from the
internal point of view. For instance, it isn't finitely generated, because else
all the individual cotangent sheaves~$\Omega^1_{X|S}$ would be of finite type
(even with a uniform bound on the number of generators, if~$S$ is quasicompact) by
Proposition~\ref{prop:locally-free-big-zariski}.

We feel that the significance of the canonical morphism~$\flat\affl \to
\affl$ hasn't been adequately explored yet.
\end{ex}

\begin{rem}Some care is needed when dealing with the modalities~$\flat$
and~$\sharp$, since they are not compatible with change of base.
If~$f : X \to S$ is a morphism of schemes, then in general~$f^{-1}(\flat E)
\not\cong \flat(f^{-1}E)$, since
\begin{align*}
  f^{-1}(\flat E) &= ((T \xra{g} X) \mapsto \Gamma(T,
  g^{-1}f^{-1}(E|_{\Sh(S)}))), \quad\text{but}\\
  \flat(f^{-1} E) &= ((T \xra{g} X) \mapsto \Gamma(T, g^{-1}(E|_{\Sh(X)}))).
\end{align*}
A special case in which the canonical morphism~$f^{-1}(\flat E) \to
\flat(f^{-1}E)$ is an isomorphism is when~$f$ is an open immersion.

A consequence of the fact that~$\flat$ and~$\sharp$ aren't compatible with
change of base is that these modalities can't be defined in the internal
language of~$\Zar(S)$, since any construction which can be described in the
internal language is automatically compatible with change of base. However, the
modalities can still be used and their general properties can even be
elementarily axiomatized~\cite{awodey:birkedal:local-maps}.
\end{rem}

\subsection{The big Zariski topos of a relative spectrum}

\begin{thm}\label{thm:big-zariski-topos-of-relative-spectrum}
Let~$\A$ be a quasicoherent~$\O_S$-algebra. In the situation that the
parsimonious sites are used for constructing big Zariski toposes, the big
Zariski topos of~$\RelSpec_S(\A)$ is, from the internal point of view
of~$\Sh(S)$, the classifying topos of the theory of
local~$\A$-algebras which are local over~$\O_S$.
\end{thm}

\begin{proof}The proof is similar to the proof of
Theorem~\ref{thm:zar-classifies}. Let~$X = \RelSpec_S(\A)$ and~$f : X \to S$ be the
canonical morphism. The big Zariski topos of~$\RelSpec_S(\A)$ is
a~$\Sh(S)$-topos by the composition~$\Zar(\RelSpec_S(\A)) \to \Zar(S) \to
\Sh(S)$. The pullback of~$\O_S$ along this geometric morphism
is~$f^{-1}(\flat\affl)$. A canonical~$\O_S$-algebra in~$\Zar(\RelSpec_S(\A))$
is therefore
\[ f^{-1}(\flat\affl) \lra \flat\afflx \lra \afflx. \]
This algebra is indeed local and local over~$f^{-1}(\flat\affl)$.

For verifying the universal property, it suffices to restrict to the case
that~$S = \Spec(R)$ is affine, as in the proof of
Theorem~\ref{thm:zar-classifies}, and consider a geometric morphism~$f : \E \to
\Sh(S)$. In this case~$\A = A^\sim$
and~$X = \RelSpec_S(\A) = \Spec(A)$. Let~$\alpha : R \to A$ be the structure
morphism of~$A$. We then have the chain of equivalences
\begin{align*}
  &\mathrel{\phantom{\simeq}} \Hom_{\Sh(S)}(\E, \Zar(X)) \\
  &\simeq \text{cat.\@ of local algebras~$\varphi : \ul{A} \to \B$
  in~$\Zar(X)$ such that~$\ul{\alpha}^{-1} \varphi^{-1} \B^\times = f^{-1}\F$} \\
  &\simeq \text{cat.\@ of local algebras~$\psi : f^{-1}\A \to \B$ such that
  $f^{-1}\O_S \to f^{-1}\A \to \B$ is local}.
\end{align*}
The first equivalence maps a geometric morphism~$g : \E \to \Zar(X)$ to~$\ul{A}
\to g^{-1}\afflx$. The second equivalence acts as follows. Given an
algebra~$\varphi : \ul{A} \to \B$ such that~$\ul{\alpha}^{-1} \varphi^{-1}
\B^\times = f^{-1}\F$, we can factor~$\ul{R} \to \ul{A} \to \B$ uniquely as a
localization~$\ul{R} \to C$ followed by a local homomorphism~$C \to \B$. By the
condition on filters, the localization~$C$ is isomorphic to~$f^{-1}\O_S$. From
the description~$\A = \ul{A}[\F^{-1}]$ it is apparent that~$\ul{A} \to \B$
factors over~$\ul{A} \to f^{-1}\A$. In this way, we obtain
morphisms~$f^{-1}\O_S \to f^{-1}\A \to \B$.
\end{proof}

The only reason why we have supposed that~$\A$ is quasicoherent in the
statement of Theorem~\ref{thm:big-zariski-topos-of-relative-spectrum} is
because else~$\RelSpec_S(\A)$ might fail to be a scheme, whereby the notion
``big Zariski topos of~$\RelSpec_S(\A)$'' is not defined.

In fact, we propose the following definition:\label{page:big-zariski-of-lrl} If~$(X,\O_X)$ is an arbitrary
locally ringed locale (or even a locally ringed topos), then the big Zariski
topos of~$X$ should be the classifying~$\Sh(X)$-topos of the theory (internal
to~$\Sh(X)$) of local~$\O_X$-algebras which are local over~$\O_X$.
The following proposition shows that this definition is consistent with
Theorem~\ref{thm:zar-classifies} and
with~Theorem~\ref{thm:big-zariski-topos-of-relative-spectrum}.

\begin{prop}Let~$\A$ be an~$\O_S$-algebra. The following constructions,
performed internally to~$\Sh(S)$, yield canonically equivalent toposes:
\begin{enumerate}
\item Constructing first the local spectrum~$X \defeq \Spec(\A|\O_S)$ and then,
internally to~$\Sh_{\Sh(S)}(X)$, the classifying topos of the theory
of~$\O_X$-algebras which are local over~$\O_X$.
\item Constructing the classifying topos of the theory of~$\A$-algebras which
are local over~$\O_S$.
\end{enumerate}
If furthermore~$\A$ is finitely presented as an~$\O_S$-algebra from the
internal point of view of~$\Sh(S)$, then the following construction yields the
same result as well:
\begin{enumerate}
\addtocounter{enumi}{2}
\item Constructing first the big Zariski topos of~$S$ as the classifying topos
of local~$\O_S$-algebras which are local over~$\O_S$ and then constructing,
internally to that topos, the slice topos over~$[\A^\Zar,
\affl]_{\Alg(\affl)}$.
\end{enumerate}
\end{prop}

\begin{proof}
If~$S$ is indeed a scheme, as is supposed throughout Part~\ref{part:big-zariski}, and~$\A$ is
quasicoherent, then all three constructions yield the big
Zariski topos of~$\RelSpec_S(\A)$ (defined using one of the parsimonious
sites). For the first construction, this is by
Theorem~\ref{thm:local-spectrum-yields-relative-spectrum} and
Theorem~\ref{thm:zar-classifies}; for the second construction, this is by
Theorem~\ref{thm:big-zariski-topos-of-relative-spectrum}; and for the third
construction, this is by Theorem~\ref{thm:zar-classifies},
Proposition~\ref{prop:relative-spectrum-big-zariski}, and the description of
the slice topos in Section~\ref{sect:change-of-base}. However, the claim also
holds if~$\A$ is not quasicoherent or if~$S$ is an arbitrary locally ringed
locale, and it's instructive to see the proof in this more general situation.

We work in the internal universe of~$\Sh(S)$. Let~$\E$ be an arbitrary
(Gro\-then\-dieck) topos. Then~$\E$-valued points of the three toposes are given
by:
\begin{enumerate}
\item a filter~$F \subseteq \A$ lying over the filter of units of~$\O_S$
together with a local~$\A_F$-algebra~$R$ which is local over~$\A_F$
\item a local~$\A$-algebra which is local over~$\O_S$
\item a local~$\O_S$-algebra~$R$ which is local over~$\O_S$ together with an
element of the stalk of~$[\A^\Zar, \affl]_{\Alg(\affl)}$ at the point
corresponding to~$R$
\end{enumerate}
In the case that~$\A$ is finitely presented, the stalk appearing in
description~(3) is canonically isomorphic to the set of~$R$-algebra
homomorphisms~$\A \otimes_{\O_S} R \to R$, as discussed in
Lemma~\ref{lemma:fp-hom-geometric}.

With these descriptions, the equivalence is immediate. For instance, a
datum~$(F \subseteq \A, \A_F \to R)$ as in description~(1) gives rise to the
datum~$(\O_S \to \A_F \to R)$ as in description~(2). Conversely, the structure
morphism of a datum as in description~(2) can be factored as a localization
followed by a local homomorphism to yield a datum as in~(1).
\end{proof}

\subsection{Constructing the big Zariski topos without recourse to classical
scheme theory}\label{sect:big-zariski-without-classical-scheme-theory}

Given a scheme~$S$, is it possible to construct the big Zariski topos of~$S$
without recourse to classical scheme theory? Without employing a site which
refers to schemes as classically conceived?

Taken literally, this question is ill-posed, since the datum~$S$ is given as a
classical scheme. A better question is: Is it possible to setup the basics of
the theory of schemes using only big Zariski toposes, preferably even in an
intuitionistic fashion?

This is indeed possible, and we wish to sketch how this can be done. Given a
base ring~$A$, the big Zariski topos of~$\Spec(A)$ can be defined as the topos
of sheaves over the parsimonious site~$\Alg(A)_\fp^\op$ consisting of (formal
duals) of finitely presented~$A$-algebras. We can then declare an~$A$-scheme to
be an object of~$\Zar(\Spec(A))$ having certain properties, for instance being
a \emph{finitely presented synthetic scheme}
(Definition~\ref{defn:synthetic-scheme}).

The big Zariski topos of such an~$A$-scheme~$X$ can then simply be defined as
the slice topos~$\Zar(\Spec(A))/X$, in accordance with the equivalence noted in
Section~\ref{sect:change-of-base}. This slice topos can serve as the base over
which further schemes and their big Zariski toposes can be constructed.

Inaccessible to this approach to scheme theory are schemes which are not
locally of finite presentation over the base ring. If one wants to account for
such schemes, one has to substitute the parsimonious site for a larger one;
however, some problems remain, as indicated in
Section~\ref{sect:qc-qs-morphisms}.

\section{The double negation modality}
\label{sect:double-negation-modality-big-zariski}

\begin{prop}\label{prop:notnot-in-big-zariski-topos}
Let~$\varphi$ be a formula over~$S$. Consider the following
statements:
\begin{enumerate}
\item $\Zar(S) \models \neg\neg\varphi$.
\item For all points~$s \in S$, there is a field extension~$K \fieldext k(s)$ such
that~$\Spec(K) \to \Spec(k(s)) \to S$ is contained in the site used to
define~$\Zar(S)$ and such that~$\Spec(K) \models \varphi$.
\item For all closed points~$s \in S$, there is a finite field extension~$K \fieldext
k(s)$ such that~$\Spec(K) \models \varphi$.
\end{enumerate}
Then:
\begin{itemize}
\item Condition~(2) implies condition~(1). The converse holds if the site used
to define~$\Zar(S)$ is closed under taking spectra of residue fields (this
is satisfied for all sites listed in Section~\ref{sect:proper-choice-of-site}
except for the parsimonious sites).
\item If one of the parsimonious sites is used to define~$\Zar(S)$ and~$S$ is
locally Noetherian, condition~(1) implies condition~(3). The converse holds if
additionally~$S$ is locally of finite type over a field.
\end{itemize}
\end{prop}

\begin{proof}We begin with showing that condition~(2) implies condition~(1). By
the Kripke--Joyal translation, we need to verify that
\[ \forall (X \to S)\_
  \Bigl(\forall (T \to X)\_ (T \models \varphi) \Rightarrow T = \emptyset\Bigr)
  \Longrightarrow X = \emptyset, \]
where the universal quantifiers range over all schemes contained in the site
used to define~$\Zar(S)$. So let such an~$S$-scheme~$f : X \to S$ be given. We
show that the fiber over any point~$s \in S$ is empty. By assumption, there is
a field extension~$K \fieldext k(s)$ such that~$\Spec(K) \to \Spec(k(s)) \to S$
is contained in the site used to define~$\Zar(S)$ and such that~$\Spec(K)
\models \varphi$. The base change~$T$ of the fiber~$X_s$ to~$\Spec(K)$ as
indicated in the diagram
\[ \xymatrix{
  T \ar[r]\ar[d] & X_s \ar[r]\ar[d] & X \ar[d] \\
  \Spec(K) \ar[r] & \Spec(k(s)) \ar[r] & S
} \]
is contained in the site used to define~$\Zar(S)$ as well, therefore
saying~``$T \models \varphi$'' is meaningful. And indeed~$T \models \varphi$,
since~$\Spec(K) \models \varphi$. Therefore~$T = \emptyset$. Thus~$X_s =
\emptyset$ as claimed.

For the direction~``(1)~$\Rightarrow$~(2)'', let a point~$s \in S$ be given.
Since we assume that the site used to define~$\Zar(S)$ contains
the~$S$-scheme~$X \defeq \Spec(k(s))$ and since~$X \neq \emptyset$, the
assumption implies that there exists a nonempty~$X$-scheme~$T$ such that~$T
\models \varphi$.  Since~$T$ is nonempty, there exists a point~$t \in T$. By
the morphism~$\Spec(k(t)) \to T \to X$, the field~$K \defeq k(t)$ is an
extension of~$k(s)$, and since~$\Spec(K) \to T \to X \to S$ is contained in the
site, we have~$\Spec(K) \models \varphi$.

The proof that condition~(1) implies condition~(3) in the case that one of the
parsimonious sites is used to define~$\Zar(S)$ and that~$S$ is locally
Noetherian is similar. For a closed point~$s \in S$, the residue field~$k(s)$
can be computed as~$A/\mmm$, where~$A$ is the ring of functions of an open
affine neighborhood of~$s$ and~$\mmm$ is a maximal ideal in~$A$. Since~$A$ is
Noetherian, the ideal~$\mmm$ is finitely generated and therefore~$A/\mmm$ is
finitely presented as an~$A$-algebra. Thus the canonical morphism~$\Spec(k(s))
\to \Spec(A) \to S$ is locally of finite presentation and thereby contained in
the parsimonious site. The hypothesis is therefore applicable to~$X \defeq
\Spec(k(s))$ and yields a nonempty~$X$-scheme~$T$ which is locally of finite
presentation over~$X$ such that~$T \models \varphi$.

Since the structure morphism~$T \to X$ is locally of finite presentation, the
scheme~$T$ inherits the property to be locally Noetherian from~$X$. Let~$U
\subseteq T$ be a nonempty open affine subset and let~$t \in U$ be a point
which is closed in~$U$. With the same reasoning as above, the canonical
morphism~$\Spec(k(t)) \to U \to T$ is therefore contained in the parsimonious
site. Thus~$\Spec(k(t)) \models \varphi$. The field~$K \defeq k(t)$ is finitely
presented as a~$k(s)$-algebra. By Noether normalization, it is also of finite
dimension as a~$k(s)$-vector space.

Finally, we verify that condition~(3) implies condition~(1) if one of the
parsimonious sites is used to define~$\Zar(S)$ and if~$S$ is locally
of finite type over a field (and therefore in particular Noetherian). We adopt
the notation of the proof of~``(2)~$\Rightarrow$~(1)''. The argument there
shows that all fibers of~$f$ over closed points are empty. If~$X$ is not empty,
it contains a closed point~$x$ (since~$X$ is locally of finite type over a field,
any point which is closed in an open affine neighborhood will do). Since~$X$
is locally of finite type over a field, the point~$f(x)$ is closed in~$S$.
Therefore~$x$ is contained in the fiber over a closed point; a contradiction.
\end{proof}

\begin{rem}The proof of Proposition~\ref{prop:notnot-in-big-zariski-topos} uses
classical logic in a substantial way, since repeatedly the lemma that a scheme
is trivial if it doesn't have any points was used. Even if scheme theory is
set up in an intuitionistically sensible way (for instance by defining a scheme to be
a locally ringed locale which is locally isomorphic to the locale-theoretic
spectra of rings as discussed in Section~\ref{sect:spectrum-as-a-locale}), one
should therefore not expect the proposition to admit an intuitionistic proof
without additional hypotheses.
\end{rem}

\begin{lemma}\label{lemma:image-coincides}
Let~$f : X \to S$ and~$g : Y \to S$ be~$S$-schemes which are locally contained
in the site used to define~$\Zar(S)$. In the case that the site is one of the
parsimonious sites, further assume that~$f$ and~$g$ are quasicompact and quasiseparated.
The following statements are equivalent:
\begin{enumerate}
\item The image of~$f$ coincides with the image of~$g$ topologically.
\item $\Zar(S) \models \neg\neg(\speak{$\ul{X}$ inhabited}) \Leftrightarrow
  \neg\neg(\speak{$\ul{Y}$ inhabited})$.
\end{enumerate}
\end{lemma}

\begin{proof}By Proposition~\ref{prop:char-surjective-morphisms}, which we'll
prove below, statement~(2) is equivalent to:
\begin{indentblock}
For any~$S$-scheme~$h : T \to S$ contained in the site used to
define~$\Zar(S)$, the morphism~$X \times_S T \to T$ is surjective if and only
if~$Y \times_S T \to T$ is.
\end{indentblock}
We verify that this statement implies statement~(1). Let~$s \in \im(f)$. Then the canonical
morphism~$X_s \to \Spec(k(s))$ is surjective. Therefore there exists
an~$S$-scheme~$h : T \to S$ which is contained in the site used to
define~$\Zar(S)$ such that~$X \times_S T \to T$ is surjective and such that~$s
\in \im(h)$: If~$\Zar(S)$ is defined using a Grothendieck or partial universe,
this claim is trivial, since we can take~$T \defeq \Spec(k(s))$. If~$\Zar(S)$
is defined using one of the parsimonious sites, we employ the technique of
relative approximation.\footnote{More specifically, we may assume that~$S$ is
affine. Then the lemma on relative approximation~\stacksproject{09MV} can be
applied to write~$\Spec(k(s))$ as a directed limit of an inverse system of
finitely presented~$S$-schemes~$T_i$ with affine transition maps. Let~$U
\subseteq X$ be an open affine subset containing a preimage of~$s$. The
property that~$U \hookrightarrow X \to \Spec(k(s))$ is surjective descends to
one of the morphisms~$U \times_S T_i \to T_i$~\stacksproject{07RR}. In
particular, the morphism~$X \times_S T_i \to T_i$ is surjective. We can
therefore take~$T \defeq T_i$. The image of~$T_i \to S$ contains~$s$
since~$\Spec(k(s)) \to S$ factors over~$T_i \to S$.}

The assumption yields that the induced morphism~$Y \times_S T \to T$ is
surjective. Since~$s \in \im(h)$, also~$s \in \im(g)$.

The proof of the converse containment relation is analogous.

The direction~``(1)~$\Rightarrow$~(2)'' is immediate, since
\[ \im(X \times_S T \to T) = h^{-1} \im(f) = h^{-1} \im(g) = \im(Y \times_S T \to T).
  \qedhere \]
\end{proof}

\begin{rem}Let~$f : X \to S$ be contained in the site used to define~$\Zar(S)$.
In the case that the site is one of the parsimonious sites, further assume
that~$f$ is quasicompact and quasiseparated.
The expression~``$\neg\neg(\speak{$\ul{X}$ is inhabited})$'' of the internal
language of~$\Zar(S)$ denotes the subfunctor of the terminal
functor~$\ul{S} = 1 \in \Zar(S)$ given by
\[ (h : T \to S) \longmapsto \{ \star \,|\, \im(h) \subseteq \im(f) \}. \]
If~$f$ is an open immersion, then this functor coincides with the functor of
points of~$X$, since the set-theoretic image of a morphism of schemes is
contained in an open subset~$U \subseteq S$ if and only if it factors over the
open immersion~$U \hookrightarrow S$.

If~$f$ is a closed immersion, this functor is the functor of points of the
formal completion of~$S$ along~$X$. More generally, for an
arbitrary~$S$-scheme~$X$ and a closed subscheme~$Z \hookrightarrow X$ (such
that both~$X$ and~$S$ are locally contained in the site used to
define~$\Zar(S)$), the internal expression~``$\{ x \? \ul{X} \,|\, \neg\neg(x \in
\ul{Z}) \}$'' denotes the functor of points of the formal completion of~$X$
along~$Z$. For instance, the expression
\[ \{ f \? \affl \,|\, \neg\neg(f = 0) \} =
  \{ f \? \affl \,|\, \speak{$f$ is nilpotent} \} \]
denotes the formal neighborhood of the origin in the affine line~$\AA^1_S$.
(The equivalence~$\neg\neg(f = 0) \Leftrightarrow \speak{$f$ is nilpotent}$ is
by Proposition~\ref{prop:a1-nilp}.)
\end{rem}

\section{Sheaves of rings, algebras, and modules}

\begin{prop}\label{prop:locally-free-big-zariski}
Let~$\E$ be an~$\O_S$-module. Properties of~$\E$ and of the
induced~$\affl$-module~$\E^\Zar$ are related as follows:
\begin{itemize}
\item $\E$ is finite locally free if and only if~$\E^\Zar$ is finite free as
an~$\affl$-module from the internal point of view of~$\Zar(S)$.
\item $\E$ is of finite type if and only if~$\E^\Zar$ is finitely generated as
an~$\affl$-module from the internal point of view of~$\Zar(S)$.
\item $\E$ is of finite presentation if and only if~$\E^\Zar$ is finitely
presented as an~$\affl$-module from the internal point of view of~$\Zar(S)$.
\end{itemize}
\end{prop}

\begin{proof}The ``if'' directions follow just as in
Proposition~\ref{prop:finite-type-and-co}. The proofs of the ``only if''
directions can too proceed by hand, further exploiting only that generators and
relations are stable under base change. But there's also a more conceptual
proof: From the point of view of~$\Sh(S)$, the~$\affl$-module~$\E^\Zar$ admits
the description~$\E^\Zar = \ul{\E} \otimes_{\ul{\O_S}} \affl$, where the
underline denotes the constant sheaf construction (which externally is
interpreted by the functor~$\pi^{-1}$ described in
Section~\ref{sect:relation-big-little}). By
Lemma~\ref{lemma:properties-of-constant-sheaves-over-locales} (generalized from
locales to toposes), the~$\ul{\O_S}$-module~$\ul{\E}$ inherits any property
from~$\E$ which is of a certain logical form, and the properties under
discussion are stable under tensoring.
\end{proof}

An analogous result holds for sheaves of algebras and their finiteness
conditions.

Coherence is missing from the list in
Proposition~\ref{prop:locally-free-big-zariski} since coherence is not stable
under pullback (for instance, the pullback of~$\O_{\Spec(\ZZ)}$ is not stable
under any morphism~$X \to \Spec(\ZZ)$ for which~$\O_X$ is not coherent) and can
therefore not be characterized by any formula in the internal language of a
big Zariski topos whose site is sufficiently encompassing.

However, coherence is stable under pullback along locally finitely presented
morphisms. Therefore we do have the following result.

\begin{prop}\label{prop:coherent-big-zariski}
Let~$\E$ be an~$\O_S$-module. If the induced~$\affl$-module~$\E^\Zar$ is
coherent from the internal point of view of~$\Zar(S)$, then~$\E$ is coherent.
The converse holds if~$\Zar(S)$ is defined using one of the parsimonious sites.
\end{prop}

\begin{prop}\label{prop:affl-modules-fully-faithful}
The inclusion functor
\[ \Mod_{\Sh(S)}(\O_S) \longrightarrow \Mod_{\Zar(S)}(\affl),\ \E \longmapsto
\E^\Zar \]
is fully faithful.
\end{prop}

\begin{proof}This follows from the following string of isomorphisms, employing
the adjunction~$\pi^{-1} \dashv \pi_*$ and the fact that~$\pi_*$ commutes with
tensor product by Lemma~\ref{lemma:zar-tensor-product-commutes}.
\begin{align*}
  \Hom_{\Mod(\affl)}(\E^\Zar, \F^\Zar) &=
  \Hom_{\Mod(\affl)}(\pi^{-1}\E \otimes_{\pi^{-1}\O_S} \affl, \pi^{-1}\F \otimes_{\pi^{-1}\O_S} \affl) \\
  &\cong
  \Hom_{\Mod(\pi^{-1}\O_S)}(\pi^{-1}\E, \pi^{-1}\F \otimes_{\pi^{-1}\O_S} \affl) \\
  &\cong
  \Hom_{\Mod(\O_S)}(\E, \pi_*(\pi^{-1}\F \otimes_{\pi^{-1}\O_S} \affl)) \\
  &\cong
  \Hom_{\Mod(\O_S)}(\E, \pi_*\pi^{-1}\F \otimes_{\pi_*\pi^{-1}\O_S} \pi_*\affl)) \\
  &\cong
  \Hom_{\Mod(\O_S)}(\E, \F \otimes_{\O_S} \O_S) \\
  &\cong
  \Hom_{\Mod(\O_S)}(\E, \F). \qedhere
\end{align*}
\end{proof}

\begin{caveat}The inclusion functor from the category of~$\O_S$-modules to the
category of~$\affl$-modules is right exact, and preserves even arbitrary
colimits, but is not left exact.
\end{caveat}

\subsection{Relative spectrum}

\begin{defn}In the context of a specified local ring~$R$, the \emph{synthetic spectrum} of an~$R$-algebra~$A$ is
\[ \Spec(A) \defeq [A, R]_{\mathrm{Alg}(R)}, \]
the set of~$R$-algebra homomorphisms from~$A$ to~$R$.\end{defn}

\begin{ex}The synthetic spectrum of~$R$ is the one-element set.
More generally, the synthetic spectrum of the algebra~$R[X_1,\ldots,X_n]/(f_1,\ldots,f_m)$
is the solution set~$\{ x \? R^n \,|\, f_1(x) = \cdots = f_n(x) = 0 \}$.
\end{ex}

\begin{ex}The synthetic spectrum of~$R/(f)$ is~$\brak{f = 0}$, the truth value
of the formula~``$f = 0$'', the subsingleton set~$\{ \star \,|\, f = 0 \}$.
If classical logic is available, then this set contains~$\star$ or is empty,
depending on whether~$f$ is zero or not. Similarly, the synthetic spectrum
of~$R[f^{-1}]$ is~$\brak{\speak{$f$ \inv}}$.\end{ex}

\begin{prop}\label{prop:relative-spectrum-big-zariski}
Let~$\A_0$ be an~$\O_S$-algebra (not necessarily quasicoherent).
Then the synthetic spectrum of the~$\affl$-algebra~$(\A_0)^\Zar$, as constructed
in the internal language of~$\Zar(S)$, is the functor of points of~$\RelSpec_S \A_0$.
\end{prop}

\begin{proof}The Hom set occurring in the definition of the synthetic spectrum is
interpreted by the internal Hom when using the internal language. For
any~$S$-scheme~$f : T \to S$ contained in the site used to define~$\Zar(S)$, we
have the following chain of isomorphisms.
\begin{align*}
  \brak{\Spec(A)}(T) &= [(\A_0)^\Zar, \affl]_{\mathrm{Alg}(\affl)}(T) \\
  &\cong
  \Hom_{\Zar(S)}(\ul{T}, [(\A_0)^\Zar, \affl]_{\mathrm{Alg}(\affl)}) \\
  &\cong
  \Hom_{\Zar(S)}(\ul{T} \times (\A_0)^\Zar, \affl)_{\ldots} \\
  &\cong
  \Hom_{\Zar(S)/\ul{T}}(\ul{T} \times (\A_0)^\Zar, \ul{T} \times \affl)_{\ldots} \\
  &\cong
  \Hom_{\mathrm{Alg}_{\Zar(T)}(\afflt)}((f^*\A_0)^\Zar, \afflt) \\
  &=
  \Hom_{\mathrm{Alg}_{\Zar(T)}(\afflt)}(\pi^{-1}(f^*\A_0) \otimes_{\pi^{-1}\O_T}
  \afflt, \afflt) \\
  &\cong
  \Hom_{\mathrm{Alg}_{\Zar(T)}(\pi^{-1}\O_T)}(\pi^{-1}(f^*\A_0), \afflt) \\
  &\cong
  \Hom_{\mathrm{Alg}_{\Sh(T)}(\O_T)}(f^*\A_0, \O_T) \\
  &\cong
  \Hom_{\mathrm{Alg}_{\Sh(S)}(\O_S)}(\A_0, f_*\O_T) \\
  &\cong
  \Hom_S(T, \RelSpec_S \A_0).
\end{align*}
The omitted subscripts~``$\ldots$'' should denote that we're only taking the
subset of the Hom set where, for each fixed first argument, the morphisms are
morphisms of~$\affl$-algebras.
\end{proof}

If~$X \in \Zar(S)$ is an arbitrary object, there is a canonical morphism~$X \to
\Spec([X,\affl])$. In the internal language of~$\Zar(X)$ it looks like the
``inclusion into the double dual'':
\[ x \longmapsto \placeholder(x),
  \quad\text{where $\placeholder(x) : [X,\affl] \to \affl,\ \varphi \mapsto
  \varphi(x)$.} \]
The following proposition shows that bijectivity of this map is related to~$X$
being the functor of points of an affine~$S$-scheme (an~$S$-scheme whose
structure morphism to~$S$ is affine).

\begin{prop}\label{prop:char-affine-zar}
Let~$X \in \Zar(S)$ be a Zariski sheaf. Consider the following statements:
\begin{enumerate}
\item The sheaf~$X$ is isomorphic to the functor of points of an
affine~$S$-scheme.
\item In the internal language of~$\Zar(S)$,
the~$\affl$-algebra~$[X,\affl]$ is synthetically quasicoherent
(Definition~\ref{defn:synth-qcoh}) and the canonical map~$X \to
\Spec([X,\affl])$ is bijective.
\end{enumerate}
If one of the parsimonious sites is used to define~$\Zar(S)$ or if, internally,
the algebra~$[X,\affl]$ is finitely presented, then~``(2)~$\Rightarrow$~(1)''.
The converse holds in any of the following
situations:
\begin{itemize}
\item The affine~$S$-scheme which~$X$ represents is of finite presentation
over~$S$.
\item The site used to define~$\Zar(S)$ is defined using a partial universe
and the affine~$S$-scheme which~$X$ represents is of finite type over~$S$.
\item The affine~$S$-scheme which~$X$ represents is contained in the site used
to define~$\Zar(S)$. (This situation subsumes the previous ones.)
\end{itemize}\end{prop}

\begin{proof}The direction~``(2)~$\Rightarrow$~(1)'' is straightforward, since the assumption
expresses~$X$ as the functor of points of the relative spectrum of a
quasicoherent~$\O_S$-algebra. Theorem~\ref{thm:qcoh-big-char}, which is needed
to obtain the sought quasicoherent~$\O_S$-algebra, is applicable.

For the converse direction, we abuse notation and denote the given
affine~$S$-scheme whose functor of points is~$X$ by~``$f : X \to S$''.
Then~$f_*\O_X$ is quasicoherent and the canonical morphism~$X \to
\RelSpec_S{f_*\O_X}$ is an isomorphism. In any of the listed situations, the
internal Hom~$[X,\affl]$ is canonically isomorphic to~$(f_*\O_X)^\Zar$, since
for any object~$T \xra{g} S$ of the site used to define~$\Zar(S)$ we have that
\begin{align*}
  [X,\affl](T) &\cong
  \Hom_{\Zar(S)}(\underline{T}, [X,\affl]) \cong
  \Hom_{\Zar(S)}(\underline{T} \times X,\affl) \\
  &\cong \Hom_{\Zar(S)}(\underline{T} \times \underline{X},\affl) \cong
  \Hom_{\Zar(S)}(\underline{T \times_S X},\affl) \\
  &\cong \affl(T \times_S X) \cong
  (g^* f_* \O_X)(T) =
  (f_* \O_X)^\Zar(T).
\end{align*}
Therefore~$[X,\affl]$ is quasicoherent. The map induced by the isomorphism~$X \to
\RelSpec_S{f_*\O_X}$ on the level of functors of points is precisely the
canonical map~$X \to \Spec([X,\affl])$ as defined in the internal language;
therefore this map is bijective from the internal point of view.
\end{proof}

\begin{rem}\label{rem:local-representability}
The condition in Proposition~\ref{prop:char-affine-zar} that~$\Spec(A)$ is
representable by an object of the site used to define~$\Zar(S)$ is slightly
unnatural from a topos-theoretic point of view, since the conclusion of the
Scholium depends only on the topos over the site and not the site itself.
In fact, the condition can be weakened and made more natural at the
same time: It suffices to require that~$\Spec(A)$ is \emph{locally}
representable by an object of the site.
\end{rem}

The following corollary answers a question by
Madore~\cite[entry~2002-07-07:044]{madore:diary}.

\begin{cor}A morphism~$f : X \to S$ of schemes is finite if, from
the internal point of view of~$\Zar(S)$, the canonical map~$\ul{X} \to
\Spec([\ul{X},\affl])$ is bijective and the~$\affl$-algebra~$[\ul{X}, \affl]$
is synthetically quasicoherent and finitely generated. The converse holds
if~$X$ is locally contained in the site used to define~$\Zar(S)$.
\end{cor}

\begin{proof}Immediate using Proposition~\ref{prop:char-affine-zar} and
Proposition~\ref{prop:locally-free-big-zariski}.
\end{proof}

\begin{rem}Let~$\A_0$ be an~$\O_S$-algebra. Then one can form, internally
to~$\Zar(S)$, two locales related to~$\A_0$: the discrete locale on the synthetic
spectrum of~$(\A_0)^\Zar$, and the local spectrum of~$(\A_0)^\Zar$ over~$\affl$ as described
in~Definition~\ref{sect:rel-spec-as-ordinary-spec}. These locales don't coincide. In fact,
the pullback of a discrete locale is discrete, whereas the pullback of
the local spectrum to any of the little Zariski toposes~$\Sh(X)$, where~$f : X
\to S$ is an~$S$-scheme contained in the site used to define~$\Zar(S)$, is the
relative spectrum~$\RelSpec_X(f^*\A_0)$, which is typically not discrete as
an~$X$-locale. (This is because the local spectrum construction is geometric,
by Proposition~\ref{prop:local-spectrum-generic}.)

There is, however, a comparison morphism from the discrete locale on the
synthetic spectrum to the local spectrum. On points, it sends
an~$\affl$-algebra homomorphism~$\varphi : (\A_0)^\Zar \to \affl$ to the
filter~$\varphi^{-1}[(\affl)^\times]$.

One can also form, internally to~$\Zar(S)$, the classifying topos
of~$(\A_0)^\Zar$-algebras which are local over~$\affl$. This topos doesn't
coincide with the (toposes of sheaves over) the mentioned two locales, either.
The pullback of that classifying topos to any of the~$\Sh(X)$ is the big
Zariski topos of~$\RelSpec_X(f^*\A_0)$ (built using one of the parsimonious
sites).
\end{rem}

\subsection{Relative Proj construction}

\begin{defn}In the context of a specified local ring~$R$, the \emph{synthetic
Proj} of an~$\NN$-graded~$R$-algebra~$A$ which is generated as an~$A_0$-algebra
by~$A_1$ is the set
\[ \Proj(A) \defeq
  (\text{set of all surj.\@ graded $R$-algebra homomorphisms~$A \to R[T]$})/R^\times. \]
\end{defn}

\begin{ex}\label{ex:proj-polynomial-ring}
The synthetic Proj of~$R[X_0,\ldots,X_n]$ is canonically isomorphic
to the set of points~$[x_0 \hg \cdots \hg x_n]$ with at
least one invertible coordinate.
\end{ex}

\begin{prop}\label{prop:relative-proj-big-zariski}
Let~$\A$ be an~$\NN$-graded~$\O_S$-algebra (not necessarily quasicoherent).
Assume that~$\A$ is generated as an~$\A_0$-algebra by~$\A_1$.
Then the synthetic Proj of the~$\affl$-algebra~$\A^\Zar$, as constructed
in the internal language of~$\Zar(S)$, is the functor of points of~$\RelProj_S \A$.
\end{prop}

\begin{proof}We omit the somewhat tedious verification.
\end{proof}

The following corollary was prompted by a question on
MathOverflow~\cite{mo:pp1}. We are grateful to Yuhao Huang for the impulse.

\begin{cor}\label{cor:pp1-classifies}
The big Zariski topos of the projective line~$\PP^1_\ZZ$, when constructed using a parsimonious site, classifies the theory of ``a local
ring together with a point~$[a \hg b]$'' (that is a pair~$(a,b)$ of ring elements,
where at least one coordinate is invertible, up to multiplication by units).
Explicitly, this theory is given by:
\begin{enumerate}
\item A sort~$A$ together with function symbols, constants, and axioms expressing
that~$A$ is a local ring.
\item A sort~$P$ (to be thought of as the set of $[a \hg b]$ with $a,b\?A$ where at
least one coordinate is invertible) together with a
relation~$\langle\cdot,\cdot,\cdot\rangle$ on $A \times A \times P$ and the
following axioms:
\begin{itemize}
\item $\speak{$a$ \inv} \vee \speak{$b$ \inv}
\dashv\vdash_{a,b{:}A} \exists p\?P\_ \langle a,b,p \rangle$
\item $\langle a,b,p \rangle \wedge \langle a,b,p' \rangle
\vdash_{a,b{:}A,\,p,p'{:}P} p = p'$
\item $\top \vdash_{p{:}P} \exists a,b\?A\_ \langle a,b,p \rangle$
\item $\langle a,b,p \rangle \wedge \langle a',b',p \rangle
\dashv\vdash_{a,a',b,b'{:}A,\, p{:}P} \exists s\?A\_ \speak{$s$ \inv} \wedge a' = s a \wedge b' = s b$
\end{itemize}
\item A constant of sort~$P$.
\end{enumerate}
\end{cor}

\begin{proof}The big Zariski topos of~$\PP^1_\ZZ$ is a topos over the big
Zariski topos of~$\Spec(\ZZ)$; from the point of view of~$\Zar(\Spec(\ZZ))$, it
is the classifying topos of a point~$[a \hg b]$ where~$a,b\?\afflz$,
since~$\Zar(\PP^1_\ZZ) \simeq \Zar(\Spec(\ZZ))/\ul{\PP^1_\ZZ}$
as discussed in Section~\ref{sect:change-of-base} and~$\ul{\PP^1_\ZZ} \cong \{ [a \hg b] \,|\,
a,b\?\afflz \}$ by Proposition~\ref{prop:relative-proj-big-zariski} and
Example~\ref{ex:proj-polynomial-ring}. The big Zariski topos of~$\Spec(\ZZ)$
classifies local rings. Therefore the claim follows by considering the
combined geometric theory.

An alternative proof builds upon Remark~\ref{rem:zar-classifies-absolute} and
the description of the theory which the little Zariski topos of~$\PP^1_\ZZ$
classifies (Proposition~\ref{prop:proj-classifying-locale}). Combining these,
we see that~$\Zar(\PP^1_\ZZ)$ classifies the theory of a homogeneous filter~$F$
of~$\ZZ[X,Y]$ meeting the irrelevant ideal together with a local
homomorphism~$\alpha : \ZZ[X,Y][F^{-1}]_0 \to A$ into a local ring~$A$. Such data
gives rise to a point~$[\alpha(X/u) : \alpha(Y/u)]$, where~$u$ is a
homogeneous element of degree~1 contained in~$F$; and conversely any
point~$[a \hg b]$ gives rise to a filter
\[ F \defeq \{ f \in \ZZ[X,Y] \,|\, \text{$f_n(a,b)$ is invertible in~$A$
for some~$n \geq 0$} \}, \]
where~$f_n$ is the homogeneous component of degree~$n$ of~$f$, and a local
homomorphism~$\alpha : \ZZ[X,Y][F^{-1}]_0 \to A$ mapping~$f/g$
to~$f(a,b)/g(a,b)$.
\end{proof}

\subsection{Quasicoherence} The goal of this section is to give an internal
characterization of quasicoherence. We'll build several notions of synthetic
algebraic geometry on quasicoherence; it is therefore central to the theory.

\begin{lemma}\label{lemma:qcoh-local}
Let~$E$ be an~$\affl$-module. Let~$S = \bigcup_i U_i$ be an open covering such
that the restrictions~$E|_{\Zar(U_i)}$ are quasicoherent, that is of the
form~$(\E_i)^\Zar$ for quasicoherent~$\O_{U_i}$-modules~$\E_i$. Then~$E$ is
quasicoherent, that is of the form~$(\E_0)^\Zar$ for a
quasicoherent~$\O_S$-module~$\E_0$.
\end{lemma}


\begin{proof}The given modules~$\E_i$ glue to a
quasicoherent~$\O_S$-module~$\E_0$, and the sheaf condition ensures that~$E$ is
isomorphic to~$(\E_0)^\Zar$. Details are given
in~\stacksproject{03DN}.
\end{proof}

\begin{defn}\label{defn:synth-qcoh}
An~$R$-module~$E$ is \emph{synthetically quasicoherent} if and only if,
for any finitely presented~$R$-algebra~$A$, the canonical~$R$-module
homomorphism
\[ E \otimes_R A \longrightarrow [\Spec(A), E] = [[A, R]_{\Alg(R)}, E] \]
which maps a pure tensor~$x \otimes f$ to the function~$(\varphi \mapsto \varphi(f) x)$ is
bijective. Here and in the following, the set~$[\Spec(A), E]$ is the set of all
maps~$\Spec(A) \to E$, and~$[A,R]_{\Alg(R)}$ is the set of all~$R$-algebra
homomorphisms~$A \to R$.\end{defn}

This definition has the following interpretation. The codomain of the displayed
canonical map is the set of all~$E$-valued functions on~$\Spec(A)$. Elements
of~$E \otimes_R A$ induce such functions; these induced functions can
reasonably be called ``algebraic''. In a synthetic context, there should be no
other~$E$-valued functions as these algebraic ones, and different algebraic
expressions should yield different functions. This is precisely what the
postulated bijectivity expresses.

The notion of synthetic quasicoherence is only meaningful in an intuitionistic
context. For instance, even~$R$ itself can't be synthetically quasicoherent in
the presence of the law of excluded middle, since it forces the canonical evaluation
morphism $R[T] \to [R, R]$ (obtained by setting~$A \defeq R[T]$ in the
definition of synthetic quasi\-co\-he\-rence) to never be bijective: If~$R$ is
finite, then the evaluation morphism isn't injective, since~$\prod_{x \in R} (T - x)$
is mapped to the same function as the zero polynomial is. If~$R$ is
infinite, then~$R[T]$ has cardinality~$|R|$ while the set of functions~$R \to
R$ has strictly greater cardinality.

\begin{thm}\label{thm:qcoh-big-char}
Let~$E \in \Zar(S)$ be an~$\affl$-module.
If~$E$ is quasicoherent, that is of
the form~$(\E_0)^\Zar$ for some quasicoherent~$\O_S$-module~$\E_0$,
then~$E$ is synthetically quasicoherent from the internal point of view of~$\Zar(S)$.
The converse holds in any of the following situations:
\begin{enumerate}
\item The site used to construct~$\Zar(S)$ is one of the parsimonious sites.
\item The functor~$E$ maps directed limits of inverse systems of~$S$-schemes with
affine transition morphisms to colimits in~$\Set$.
\item From the internal point of view of~$\Zar(S)$, the module~$E$ is even
finitely presented.
\end{enumerate}
\end{thm}

\begin{proof}Let~$E = (\E_0)^\Zar$ for a quasicoherent~$\O_S$-module~$\E_0$.
To verify that~$E$ is synthetically quasicoherent, we have to
verify a condition for~$\affl$-algebras~$A$ in any
slice~$\Zar(S)/\ul{T}$. If such an algebra is finitely presented from
the internal point of view, then there is a covering~$T = \bigcup_i T_i$ such
that each of the restrictions of the algebra to the schemes~$T_i$ is of the
form~$(\A_0)^\Zar$ for some finitely presented~$\O_{T_i}$-algebra~$\A_0$.
Without loss of generality, we will just assume that~$A$ itself is of the
form~$(\A_0)^\Zar$ for a finitely presented~$\O_S$-algebra~$\A_0$.

By Proposition~\ref{prop:relative-spectrum-big-zariski}, the interpretation~$\brak{\Spec(A)}$ of the internal spectrum is the functor of
points of~$\RelSpec_S \A_0$. For any~$S$-scheme~$f : T \to S$ contained in the site
used to define~$\Zar(S)$, we consider the fiber product
\[ \xymatrix{
\RelSpec_T(f^*\A_0) \ar[r]^{f'} \ar[d]_{p'} & \RelSpec_S\A_0 \ar[d]^p \\
T \ar[r]_f & S.
} \]
Since~$\RelSpec_T(f^*\A_0) \to S$ is contained in the site (for any of our
admissible sites), we may conclude using the following chain of isomorphisms:
\begin{align*}
[\Spec(A), E](T) &\cong
\Hom_{\Zar(S)}(\ul{T}, [\Spec(A), E])
\cong \Hom_{\Zar(S)}(\ul{T} \times \Spec(A), E) \\
&\cong \Hom_{\Zar(S)}(\underline{T \times_S \RelSpec_S\A_0}, E)
\cong E(\RelSpec_T(f^*\A_0)) \\
&\cong \Gamma(\RelSpec_T(f^*\A_0), (p')^* f^* \E_0)
\cong \Gamma(T, (p')_* (p')^* f^* \E_0) \\
&\cong \Gamma(T, f^*\E_0 \otimes_{\O_T} f^*\A_0)
\cong (\E_0 \otimes_{\O_S} \A_0)^\Zar)(T) \\
&\cong ((\E_0)^\Zar \otimes_\affl (\A_0)^\Zar)(T)
\cong (E \otimes_\affl A)(T).
\end{align*}
The antepenultimate isomorphism is because pullback of modules in~$\Sh(S)$ to
modules in~$\Sh(T)$ commutes with tensor product. The penultimate isomorphism
is because pullback of a sheaf in~$\Sh(S)$ to a sheaf in~$\Zar(S)$ commutes
with tensor product (Lemma~\ref{lemma:zar-tensor-product-commutes}).

For the converse direction, we first verify that the restriction~$E|_{\Sh(T)}$
to the little Zariski topos of any~$S$-scheme~$T$ contained in the site used
to define~$\Zar(S)$ is a quasicoherent~$\O_T$-module. For this, we employ the
quasicoherence criterion of Theorem~\ref{thm:qcoh-sheafchar}: For any open
affine subset~$T' \subseteq T$ and any function~$h \in \Gamma(T', \O_T)$ we
verify that the canonical morphism
\begin{equation}\label{eqn:want-iso}\tag{$\dagger$}
E|_{\Sh(T)}[h^{-1}] \longrightarrow j_*(E|_{\Sh(D(h))})
\end{equation}
is an isomorphism, where~$j : D(h) \hookrightarrow T'$ denotes the inclusion.
This follows from the assumption of synthetic quasicoherence by considering
the~$\affl$-algebra~$A \defeq \affl[h^{-1}]$ (in the slice~$\Zar(S)/\ul{T'}$):
This expresses that the canonical morphism
\begin{equation}\label{eqn:have-iso}\tag{$\ddagger$}
E \otimes_\affl \affl[h^{-1}] \longrightarrow [\Spec(A), E]
\end{equation}
is an isomorphism (of~$\affl$-modules in~$\Zar(S)/\ul{T'}$). Restricting the
domain to~$\Sh(T')$ yields the sheaf~$E|_{\Sh(T')} \otimes_{\O_{T'}}
\O_{T'}[h^{-1}]$, since restricting commutes with the geometric constructions
``forming the tensor product'' and ``localizing away from~$h$''.
Since~$\Spec(A)$ is the functor of points of~$D(h)$, restricting
the codomain to~$\Sh(T')$ yields the sheaf~$j_*(E|_{\Sh(D(h))})$.
The canonical morphism~\eqref{eqn:want-iso} which we want to recognize as an
isomorphism is therefore the restriction of the canonical
morphism~\eqref{eqn:have-iso} which we know to be an isomorphism.

A natural candidate for a quasicoherent~$\O_S$-module~$\E_0$ with~$E \cong
(\E_0)^\Zar$ is~$\E_0 \defeq E|_{\Sh(S)}$. We'll show that this candidate indeed fits.
Let~$f : T \to S$ be any~$S$-scheme contained in the site used to
define~$\Zar(S)$. We assume for the time being that~$f$ is of finite
presentation and affine, so~$T \cong \RelSpec_S \A_0$ for some finitely
presented~$\O_S$-algebra~$\A_0$. We want to verify that the canonical morphism
\begin{equation}\label{eqn:want-iso2}\tag{§}
f^*(E|_{\Sh(S)}) \longrightarrow E|_{\Sh(T)}
\end{equation}
is an isomorphism. Since the functor~$f_*$ from quasicoherent~$\O_T$-modules to
quasicoherent~$\O_S$-modules reflects isomorphisms (the morphism~$f$ being affine)
and the domain and codomain of morphism~\eqref{eqn:want-iso2} are quasicoherent, it suffices to
verify that its image under~$f_*$ is an isomorphism. This image is the
canonical morphism
\[ E|_{\Sh(S)} \otimes_{\O_S} \A_0 \longrightarrow f_*(E|_{\Sh(T)}). \]
The assumption of synthetic quasicoherence, applied to the finitely presented~$\affl$-algebra~$A
\defeq (\A_0)^\Zar$, shows that this morphism is an isomorphism.

In situation~(1), the only step left to do is to generalize the argument in the
previous paragraph to morphisms~$f : T \to S$ which are locally of finite
presentation. This works out because there are open covers of~$S$ and~$T$ such
that the appropriate restrictions of~$f$ are of finite presentation and affine.
The assumption of synthetic quasicoherence then needs to be applied
to~$\affl$-algebras in suitable slices of~$\Zar(S)$, showing that the canonical
morphism~\eqref{eqn:want-iso2} is locally an isomorphism and therefore globally
as well.

In situation~(2), we may by Lemma~\ref{lemma:qcoh-local} assume without loss of generality that~$S$ is affine.
We then employ the technique of approximating general~$S$-schemes
by~$S$-schemes of finite presentation. Specifically, let~$f : T \to S$ be an
arbitrary~$S$-scheme contained in the site used to define~$\Zar(S)$. Without
loss of generality, we may assume that~$T$ is an affine scheme. Thus~$T$ is
quasicompact and quasiseparated, and~$S$ is quasiseparated since it is affine. We may therefore
apply the lemma of relative approximation~\stacksproject{09MV} to deduce
that~$T$ is a directed limit of an inverse system of~$S$-schemes~$f_i : T_i \to
S$ of finite presentation with affine transition maps. These~$S$-schemes
are contained in the site used to define~$\Zar(S)$. Furthermore, they inherit
quasicompactness and quasiseparatedness from~$S$. Therefore we can apply a
comparison result on the categories of quasicoherent
modules~\stacksproject{01Z0}:
\[
E(T) = E(\lim_i T_i) \cong \colim_i E(T_i) \cong
\colim_i \Gamma(T_i, f_i^* \E_0)
\cong \Gamma(T, f^* \E_0).
\]

In situation~(3), we may assume by Lemma~\ref{lemma:qcoh-local} that~$E$ is
the cokernel of a morphism~$\alpha : (\affl)^m \to (\affl)^n$ of~$\affl$-modules.
This morphism induces a morphism~$\alpha|_{\Sh(S)} : \O_S^m \to \O_S^n$
of~$\O_S$-modules. One can then check that~$E$ is canonically isomorphic
to~$(\cok(\alpha|_{\Sh(S)}))^\Zar$, by using that the restriction
functors~$\Mod_{\Zar(S)}(\affl) \to \Mod_{\Sh(X)}(\O_X)$ are jointly
conservative and right exact.
\end{proof}

\begin{scholium}\label{scholium}
Let~$E \in \Zar(S)$ be a quasicoherent~$\affl$-module.
Let~$A \in \Zar(S)$ be a quasicoherent~$\affl$-algebra such that~$\brak{\Spec(A)} \in
\Zar(S)$ is representable by an object of the site used to define~$\Zar(S)$.
Then the canonical morphism
\[ E \otimes_\affl A \longrightarrow [\Spec(A), E] \]
is an isomorphism.
\end{scholium}

\begin{proof}The second paragraph of the proof of
Theorem~\ref{thm:qcoh-big-char} applies.\end{proof}

\begin{rem}As noted in Remark~\ref{rem:local-representability} in a slightly
different context, the condition in Scholium~\ref{scholium} that~$\Spec(A)$ is
representable by an object of the site is unnatural from a topos-theoretic
point of view and should be weakened to require only local representability.

However, the condition can't be dropped completely. For instance, if we employ the
parsimonious sites and consider~$S = \Spec \ZZ$,~$E = \affl$, and $A =
(\K_S)^\Zar$ (where~$\K_S$ is the sheaf of rational functions on~$S$, which in
this case is the constant sheaf~$\ul{\QQ}$), then~$\brak{\Spec(A)}$ is the functor of
points of the~$\ZZ$-scheme~$\Spec(\QQ)$. By Proposition~\ref{prop:fingen-algebra-q}, this functor coincides with the
functor of points of the empty~$\ZZ$-scheme on the parsimonious sites;
therefore~$\Spec(A) = \emptyset$ from the internal point of view. Thus the
codomain of the canonical morphism is the zero algebra, but the domain is not.
\end{rem}

The internal quasicoherence condition for the little Zariski topos
(Theorem~\ref{thm:qcoh-sheafchar}) is related to the notion of synthetic
quasicoherence as follows. Recall that an~$\O_S$-module~$\E_0$ is quasicoherent
if and only if, from the internal point of view of~$\Sh(S)$, the localized
module~$\E_0[f^{-1}]$ is a sheaf with respect to the modal operator~$(\speak{$f$
\inv} \Rightarrow \placeholder)$ for any function~$f\?\O_S$. The sublocale
associated with this modal operator is the open sublocale~$j :
\Spec(\O_S[f^{-1}]|\O_S) \hookrightarrow \pt$. The condition can therefore also
be put in the form 
\[ \Sh(S) \models
  \forall f\?\O_S\_
  \speak{$\E_0[f^{-1}] \lra j_* j^{-1} (\E_0[f^{-1}])$ is bijective}. \]
One can verify that the functor~$j_* \circ
j^{-1}$ is canonically isomorphic to the functor~$[\brak{\speak{$f$ \inv}},
\placeholder]$, hence to the functor~$[\Spec(\O_S[f^{-1}]), \placeholder]$.
Since~$\E_0[f^{-1}] \cong \E_0 \otimes_{\O_S} \O_S[f^{-1}]$, the condition can
therefore also be put in the form
\[ \Sh(S) \models
  \forall f\?\O_S\_
  \E_0 \otimes_{\O_S} \O_S[f^{-1}] \lra [\Spec(\O_S[f^{-1}]), \E_0]. \]
The synthetic quasicoherence condition therefore implies the condition which
characterizes quasicoherence in the little Zariski topos as a special case.

\begin{lemma}\label{lemma:qcoh-ideals}
Let~$J \hookrightarrow \affl$ be an ideal such that~$\affl/J$ is
of the form~$(\E_0)^\Zar$ for an~$\O_S$-module~$\E_0$. Let~$\I$ be the kernel
of the epimorphism~$\O_S \twoheadrightarrow \E_0$ induced by the quotient
morphism~$\affl \to \affl/J$. Then, for any~$S$-scheme~$(f : X \to S)$
contained in the site used to define~$\Zar(S)$, there is a canonical
isomorphism
\[ \im(f^*\I \to \O_X) \cong J|_{\Sh(X)}. \]
\end{lemma}

\begin{proof}The short exact sequence~$0 \to J \to \affl \to \affl/J \to 0$
of~$\affl$-modules in~$\Zar(S)$ remains exact when restricted to~$\Sh(X)$, since
restricting to~$\Sh(X)$ is taking the inverse image along a geometric morphism.
Hence the sequence~$0 \to J|_{\Sh(X)} \to \O_X \to f^*\E_0 \to 0$ is exact.
On the other hand, the short exact sequence~$0 \to \I \to \O_S \to \E_0 \to 0$
yields the short exact sequence~$0 \to \im(f^*\I \to \O_X) \to \O_X \to f^*\E_0
\to 0$.
\end{proof}

\begin{rem}\label{rem:radical-not-qcoh}
The quotient~$\affl/\sqrt{(0)}$ in~$\Zar(S)$ is an example
for a sheaf of modules which is not quasicoherent even though all of its
restrictions to the little Zariski toposes~$\Sh(X)$ for morphisms~$f : X \to S$
are:
Since taking the quotient and taking the radical of an ideal are geometric
constructions, we have~$(\affl/\sqrt{(0)})|_{\Sh(X)} \cong \O_X/\sqrt{(0)}$.
These sheaves of modules are quasicoherent (Example~\ref{ex:radical-qcoh}).
However, in general,~$f^*(\O_S/\sqrt{(0)}) \not\cong \O_X/\sqrt{(0)}$.
A specific counterexample is~$S = \Spec(k)$ and~$X = \Spec(k[T]/(T^2))$.
In this case~$f^*(\O_S/\sqrt{(0)}) \cong f^*(\O_S) \cong \O_X \not\cong
\O_X/\sqrt{(0)}$.
\end{rem}

\begin{caveat}The kernel of a morphism of quasicoherent~$\affl$-modules,
computed in the category of all~$\affl$-modules, is in general not
quasicoherent. This fact is evident from Lemma~\ref{lemma:qcoh-ideals}. In the
notation of that lemma, the~$\affl$-module~$J$ is in general not quasicoherent,
since if it was, the restriction~$J|_{\Sh(X)}$ would be canonically isomorphic
to~$f^*\I$.

The category of quasicoherent~$\affl$-modules possesses kernels, since it is
equivalent to the category of quasicoherent~$\O_S$-modules by
Proposition~\ref{prop:affl-modules-fully-faithful}, but the inclusion from
quasicoherent~$\affl$-modules to arbitrary~$\affl$-modules does not preserve
them.
\end{caveat}

\begin{lemma}\label{lemma:tensor-product-qcoh}
Let~$\Zar(S)$ be defined using one of the parsimonious sites. Then tensor
products of synthetically quasicoherent modules and cokernels of morphisms
between synthetically quasicoherent modules are synthetically quasicoherent
from the internal point of view of~$\Zar(S)$.
\end{lemma}

\begin{proof}By Theorem~\ref{thm:qcoh-big-char} and the description of the
slice toposes of~$\Zar(S)$ given in Section~\ref{sect:change-of-base}, the
first claim reduces to the fact that the tensor product of
quasicoherent~$\affl$-modules is quasicoherent. Indeed, the tensor product
of~$(\E_0)^\Zar$ with~$(\F_0)^\Zar$ is~$(\E_0 \otimes_{\O_S} \F_0)^\Zar$ by
Lemma~\ref{lemma:zar-tensor-product-commutes} and the tensor product of
quasicoherent~$\O_S$-modules is quasicoherent.

The second claim reduces to the statement that the cokernel of a morphism~$\alpha :
(\E_0)^\Zar \to (\F_0)^\Zar$, calculated in the category of~$\affl$-modules,
coincides with~$(\cok(\alpha|_{\Sh(S)}))^\sim$ and that the cokernel of a
morphism of quasicoherent~$\O_S$-modules is quasicoherent.
\end{proof}

It's somewhat embarassing that we didn't give an internal proof of
Lemma~\ref{lemma:tensor-product-qcoh}. Also we don't know whether the result
holds when using one of the larger sites (the given proof doesn't generalize to
this situation since Theorem~\ref{thm:qcoh-big-char} can't be applied). We
expand on this in Section~\ref{sect:outlook}.

\begin{lemma}Let~$A$ and~$B$ be~$\affl$-algebras. The canonical
map \[ \Hom_{\Alg(\affl)}(A,B) \lra
\Hom_{\Zar(S)}(\brak{\Spec(B)},\brak{\Spec(A)}) \] is bijective in the following
situations:
\begin{enumerate}
\item The algebra~$B$ is finitely presented.
\item The algebra~$B$ is quasicoherent and the functor~$\brak{\Spec(B)} \in
\Zar(S)$ is locally representable by an object of the site used to
define~$\Zar(S)$. (This situation subsumes the previous one.)
\end{enumerate}
\end{lemma}

\begin{proof}By Scholium~\ref{scholium}, the canonical morphism~$B \to
[\Spec(B), \affl]$ is an isomorphism. Hence the claim follows by the following
entirely formal calculation:
\begin{align*}
  &\mathrel{\phantom{=}}
  \Hom(\brak{\Spec(B)}, \brak{\Spec(A)}) \\
  &=
  \Hom([B,\affl]_{\Alg(\affl)}, [A,\affl]_{\Alg(\affl)}) \\
  &\cong
  \Hom([B,\affl]_{\Alg(\affl)} \times A, \affl)_{\text{$\affl$-homomorphism in
  the second argument}} \\
  &\cong
  \Hom_{\Alg(\affl)}(A, [[B, \affl]_{\Alg(\affl)}, \affl]) \\
  &\cong
  \Hom_{\Alg(\affl)}(A, [\brak{\Spec(B)}, \affl]) \\
  &\cong
  \Hom_{\Alg(\affl)}(A, B). \qedhere
\end{align*}
\end{proof}


It's a basic fact that for an~$\O_S$-algebra~$\A$, the canonical map
\[ \Hom_{\Alg(\O_S)}(\A_0, f_*\O_X) \lra
  \Hom_{\LRL/S}(X, \RelSpec_S(\A_0)) \]
is bijective for any~$S$-scheme~$(f : X \to S)$. One should not expect that the
similar statement that for an~$\affl$-algebra~$A$, the canonical map
\[ \Hom_{\Alg(\affl)}(A, [X,\affl]) \lra \Hom_{\Zar(S)}(X, \brak{\Spec(A)}) \]
is bijective for arbitrary functors~$X \in \Zar(S)$ is bijective, holds. This
is because the objects of~$\Zar(S)$ are quite a bit more general than locally
ringed spaces (or locales) over~$S$.

\subsection{Properties of the affine line}
\label{sect:special-properties-affl}

The ring object~$\affl$ in the big Zariski topos enjoys several special
properties, some of which are unique in that they're only possible in an
intuitionistic context. We compile here a short list of such
properties. As was already mentioned, at least one of them, the field
property, was already noticed in the 1970s by Kock~\cite{kock:univ-proj-geometry}.

The statements and proofs in this subsection are formulated in the internal
language. The proofs only use the fact that~$\affl$ is a synthetically
quasicoherent local ring. This supports the meta-claim that synthetic
quasicoherence is a strong and meaningful condition.

\begin{prop}\label{prop:a1-field}
$\affl$ is a field in the sense that any element which is not zero is
invertible: $\forall x\?\affl\_ \neg(x = 0) \Rightarrow \speak{$x$ \inv}$. More generally,
for any number~$n \geq 0$,
\begin{multline*}
  \qquad\qquad \forall x_1,\ldots,x_n\?\affl\_
  \neg(x_1 = 0 \wedge \cdots \wedge x_n = 0) \Longrightarrow \\
  (\speak{$x_1$ \inv} \vee \cdots \vee \speak{$x_n$ \inv}). \qquad\qquad
\end{multline*}
\end{prop}

\begin{proof}Let $x\?\affl$ be such that~$\neg(x=0)$. We consider the quasicoherence
condition for the finitely presented~$\affl$-algebra~$A \defeq \affl/(x)$.
Since~$\Spec(A) \cong \brak{x=0} = \brak{\bot} = \emptyset$, the condition posits
that the canonical homomorphism
\[ \affl/(x) \longrightarrow [\emptyset, \affl] \]
is an isomorphism. Since its codomain is the zero algebra, so is~$\affl/(x)$.
Therefore~$1 \in (x)$, that is,~$x$ is invertible.

The more general statement follows in the same way, by using the quasicoherence
condition for~$A \defeq \affl/(x_1,\ldots,x_n)$. This yields~$1 \in
(x_1,\ldots,x_n)$. Since~$\affl$ is a local ring, one of the~$x_i$ is
invertible.\end{proof}

\begin{prop}\label{prop:a1-not-reduced}
$\affl$ is not a reduced ring:
\[ \neg \bigl(\forall x\?\affl\_ \bigl(\bigvee_{n \geq 0} x^n = 0\bigr) \Rightarrow x = 0\bigr). \]
\end{prop}

\begin{proof}Assume that~$\affl$ is reduced. Then the set~$\Delta \defeq \{
\varepsilon \in \affl \,|\, \varepsilon^2 = 0 \}$ is equal to~$\{ 0 \}$.
By the quasicoherence criterion applied to the finitely
presented~$\affl$-algebra~$A \defeq \affl[T]/(T^2)$, the canonical map
\[ \affl[T]/(T^2) \longrightarrow [\Spec(\affl[T]/(T^2)), \affl] \cong
  [\Delta, \affl] \]
is an isomorphism. It maps~$[T]$ to zero (the value of~$T$ at~$0 \in \Delta$).
Thus~$T \in (T^2)$ and therefore~$1 = 0$ in~$\affl$. This is a contradiction.
\end{proof}

In classical logic, Proposition~\ref{prop:a1-field} and
Proposition~\ref{prop:a1-not-reduced} would directly contradict each
other; only an intuitionistic context allows for fields which are not reduced.

That~$\affl$ is not reduced, irrespective of the reducedness of the base
scheme~$S$, should not come as a surprise: Reducedness is not stable under base
change, but all statements of the internal language of~$\Zar(S)$ are.
If~$\affl$ was reduced, then all~$S$-schemes (at least those contained in the
site used to construct~$\Zar(S)$) would be reduced as well. In contrast, the
structure sheaf~$\O_S$ is reduced from the point of view of the little Zariski
topos if and only if~$S$ is reduced (Proposition~\ref{prop:reduced-ring}).

\begin{prop}\label{prop:a1-nilp}
The following statements about an element~$x\?\affl$ are
equivalent:
\begin{enumerate}
\item $x$ is not invertible.
\item $x$ is nilpotent.
\item $x$ is \notnot zero.
\end{enumerate}
\end{prop}

\begin{proof}Let~$x \? \affl$ be not invertible. We consider the quasicoherence
condition for the finitely presented~$\affl$-algebra~$A \defeq \affl[x^{-1}]$.
Since~$\Spec(A) \cong \brak{\speak{$x$ \inv}} = \emptyset$, it follows
that~$\affl[x^{-1}] = 0$, similarly to the proof of
Proposition~\ref{prop:a1-field}. Thus~$x$
is nilpotent.

Let~$x \? \affl$ be a nilpotent element. Thus~$x^n = 0$ for some number~$n \geq
0$. If~$x$ was nonzero, then~$x$ and therefore~$x^n$ would be invertible, in
contradiction to~$0 \neq 1$ since~$\affl$ is a local ring.

Let~$x \? \affl$ be \notnot zero. Then~$x$ is not invertible, since if~$x$ was
invertible, then~$x$ would be nonzero.
\end{proof}

Summarizing, the following facts about nilpotents hold in the internal
universe of the big Zariski topos. Firstly, it's not true that~$\affl$ is
reduced. But this doesn't mean that there actually exist nilpotent elements
which are not zero. In fact, any nilpotent is \notnot zero.

\begin{prop}Any function~$\affl \to \affl$ is given by a unique polynomial
in~$\affl[T]$.
\end{prop}

\begin{proof}Immediate by considering the quasicoherence condition for the finitely
presented~$\affl$-algebra~$A \defeq \affl[T]$ and noticing that~$\Spec(A) \cong
\affl$.\end{proof}

This statement too cannot be satisfied in classical logic: for infinite
fields the existence part fails and for finite fields the uniqueness part
fails.

In synthetic differential geometry, the \emph{axiom of microaffinity} is
central to the theory. It is fulfilled by the image of~$\RR^1$ in any
well-adapted model of synthetic differential geometry, and also by~$\affl \in
\Zar(S)$. This fact is well-known; we include the proof only to show that it
is a consequence of synthetic quasicoherence.

\begin{prop}\label{prop:affl-microaffine}
$\affl$ fulfills the axiom of microaffinity: Let~$\Delta = \{ \varepsilon \?
\affl \,|\, \varepsilon^2 = 0 \}$. Let~$f : \Delta \to \affl$ be an arbitrary
function. Then there are unique elements~$a, b \? \affl$ such
that~$f(\varepsilon) = a + b \varepsilon$ for all~$\varepsilon \? \Delta$.
\end{prop}

\begin{proof}Immediate from the definition of synthetic quasicoherence,
considering the finitely presented~$\affl$-algebra~$\affl[T]/(T^2)$.
\end{proof}

\begin{prop}\label{prop:affl-anonymously-algebraically-closed}
$\affl$ is \emph{anonymously algebraically closed}, in the following sense:
Any monic polynomial~$p \? \affl[T]$ of degree at least one does \notnot have
a zero.\end{prop}

\begin{proof}Let~$p \? \affl[T]$ be a monic polynomial of degree at least one. Assume
that~$p$ doesn't have a zero in~$\affl$. Then the spectrum of~$A \defeq
\affl[T]/(p)$ is empty. The quasicoherence condition for~$A$ therefore implies
that~$\affl[T]/(p)$ is zero. This means that~$p$ is invertible in~$\affl[T]$.
A basic lemma in commutative algebra (whose standard proof is constructive)
then implies that with the exception of the constant term in~$p$, all
coefficients are nilpotent. This contradicts the assumption that~$p$ is monic
of degree at least one.\end{proof}

\begin{prop}$\affl$ is infinite in the following sense: For any number~$n \geq 0$
and any given elements~$x_1,\ldots,x_n \? \affl$, there is \notnot an element~$y$
which is distinct from all of the~$x_i$.
\end{prop}

\begin{proof}The polynomial~$f(T) \defeq (T - x_1) \cdots (T - x_n) + 1$
does \notnot have a zero~$y\?\affl$, since~$\affl$ is anonymously algebraically
closed. This element cannot equal any of the elements~$x_i$, since~$f(x_i) = 1$ is not zero.
\end{proof}

\begin{prop}$\affl$ fulfills the following version of the Nullstellensatz:
Let $f_1,\ldots,f_m \in \affl[X_1,\ldots,X_n]$ be polynomials without a common
zero in~$(\affl)^n$. Then there are polynomials~$g_1,\ldots,g_m \in
\affl[X_1,\ldots,X_n]$ such that~$\sum_i g_i f_i = 1$.
\end{prop}

\begin{proof}We consider the quasicoherence condition for the finitely
presented~$\affl$-algebra~$A \defeq \affl[X_1,\ldots,X_n]/(f_1,\ldots,f_m)$.
Since~$\Spec(A) \cong \{ x \in (\affl)^n \,|\, f_1(x) = \ldots = f_m(x) = 0 \}
= \emptyset$, the condition implies that~$A$ is the zero algebra just as in the
verification of Proposition~\ref{prop:a1-field}.
\end{proof}

\begin{rem}The Krull dimension of the ring~$\O_S$ of the little Zariski topos
coincides with the dimension of~$S$
(Proposition~\ref{prop:dimension-scheme-ox}). The analogous statement
for~$\affl$ in the big Zariski topos is false. Unless~$S$ is the empty scheme,
the internal statement
\[ \Zar(S) \models \speak{$\affl$ is of Krull dimension~$\leq n$} \]
is false for any natural number~$n \geq 0$: Since the property of having Krull
dimension~$\leq n$ is a geometric implication, this statement would imply that
for any~$S$-scheme~$X$ (contained in the site used to define~$\Zar(S)$) the
ring~$\O_X$ in~$\Sh(X)$ is of Krull dimension~$\leq n$.

The ring~$\affl$ in~$\Zar(S)$ is therefore an example for a ring of infinite
Krull dimension which nevertheless fulfills a field condition. A ring in the
big Zariski topos which does reflect the dimension of~$S$ is~$\flat\affl$. The
scheme~$S$ is of dimension~$\leq n$ if and only if~$\flat\affl$ is of Krull
dimension~$\leq n$ from the internal point of view of~$\Zar(S)$.
\end{rem}

\section{Basic constructions of relative scheme theory}
\label{sect:basic-constructions}

With~$\affl$ at hand, we can perform many of the usual constructions of
(relative) scheme theory internally.

\subsection{Group schemes} The functors associated to the standard group schemes~$\GG_\text{a}$, $\GG_\text{m}$,
$\mathrm{GL}_n$, and~$\mu_n$ are given by the internal expressions
\begin{align*}
  \GG_\text{a} &\defeq \affl \text{ (as an additive group)}, \\
  \GG_\text{m} &\defeq \{ x\?\affl \,|\, \speak{$x$ \inv} \}, \\
  \mathrm{GL}_n &\defeq \{ M \? (\affl)^{n \times n} \,|\, \speak{$M$ \inv} \}, \\
  \mu_n &\defeq \{ x \? \affl \,|\, x^n = 1 \}.
\end{align*}

\subsection{Affine and projective space}
Affine~$n$-space over~$S$ is given by~$(\affl)^n$, \ie internally the set
of~$n$-tuples of elements of~$\affl$. The functor of points of
projective~$n$-space over~$X$, with all its nontrivial topological and
ring-theoretical structure, is described by the naive expression
\[ \PP^n \defeq \{ (x_0,\ldots,x_n) \? (\affl)^{n+1} \,|\,
  x_0 \neq 0 \vee \cdots \vee x_n \neq 0 \}/{\sim}, \]
where the equivalence relation is the usual rescaling relation from the
internal point of view. This example was suggested by Zhen~Lin Low (private
communication).

More generally, for an~$S$-scheme~$X$, affine and projective~$n$-space
over~$X$ are given by~$\ul{X} \times (\affl)^n$ and~$\ul{X} \times \PP^n$,
respectively.

\subsection{Tangent bundle}

For an~$S$-scheme~$X$, the internal Hom~$[\Delta,\ul{X}] \in \Zar(S)$ describes the
tangent bundle of~$X$, \ie the~$S$-scheme~$\RelSpec_X{\operatorname{Sym}(\Omega^1_{X/S})} \to X \to S$, as can be seen by
chasing the definitions~\cite[Lemma~5.12.1]{brandenburg:tensor-foundations}.
Intuitively, a map~$f : \Delta \to \ul{X}$ from the internal point of view is
given by slightly more data than merely the point~$f(0)$; one also has to
specify first-order information.

This description of the (not necessarily
locally trivial) tangent bundle fits nicely with the intuition of tangent
vectors as infinitesimal curves, and in fact is precisely the definition of the
tangent bundle in synthetic differential geometry~\cite[Def.~7.1]{kock:sdg}.

\subsection{The tilde construction}\label{sect:relative-tilde-construction}
Let~$A$ be an~$\affl$-algebra.
If the functor~$\brak{\Spec(A)}$ is representable by a scheme contained in the
site used to define~$\Zar(S)$, then Section~\ref{sect:change-of-base} described
the big Zariski topos of that scheme as the slice topos~$\Zar(S)/\Spec(A)$.
From the point of view of~$\Zar(S)$, this topos looks like~``$\Set/\Spec(A)$'',
the topos of~$\Spec(A)$-indexed families of sets.

In this picture, the three functors~$f_! \dashv f^{-1} \dashv f_*$
relating~$\Zar(S)$ and~$\Zar(S)/\Spec(A)$ look as follows. We set~$X
\defeq \Spec(A)$.
\begin{enumerate}
\item The functor~$f_!$ maps a family~$(M_x)_{x{:}X}$ to its \emph{dependent
sum}~$\coprod_{x{:}X} M_x$, the ``total space'' of the family.
\item The functor~$f^{-1}$ maps a set~$N$ to the constant family~$(N)_{x{:}X}$.
\item The functor~$f_*$ maps a family~$(M_x)_{x\?X}$ to its \emph{dependent
product}, the set~$\prod_{x{:}X} M_x$ of sections of the
projection~$(\coprod_{x{:}X} M_x) \to X$.
\end{enumerate}

The affine line over~$X$ is the constant family~$\afflx = (\affl)_{x{:}X}$.

We can canonically associate an~$\afflx$-module~$E^\sim$ to a
given~$A$-module~$E$. This module lives in the topos~$\Set/X$, so
is an~$X$-indexed family of~$\affl$-modules. Explicitly, it is the family
\[ E^\sim \defeq (E \otimes_A \affl)_{x{:}X}, \]
where for the tensor product in the component~$x\?X$ the ring~$\affl$ is
regarded as an~$A$-algebra by the~$\affl$-homomorphism~$x : A \to \affl$.
To render the dependence of the tensor product on~$x$ explicit, we use the
familiar notation~$E(x) \defeq E \otimes_A \affl$ to denote the fiber of~$E$
over~$x$.

Conversely, given an~$\afflx$-module~$F = (F_x)_{x{:}X}$, there is a canonically
associated~$A$-module given by calculating the dependent product. In analogy
with the classical situation, we write
\[ \Gamma(X, F) \defeq \prod_{x{:}X} F_x \]
for this module.

A morphism~$F \to F'$ of~$\afflx$-modules is an~$X$-indexed family~$(F_x \to
F'_x)_{x{:}X}$ of morphism of~$\affl$-modules.

\begin{prop}The functor~$(\placeholder)^\sim : \Mod(A) \to \Mod(\afflx)$ is
left adjoint to the functor~$\Gamma(X, \placeholder) : \Mod(\afflx) \to
\Mod(A)$.
\end{prop}

\begin{proof}For an~$A$-module~$E$ and an~$\afflx$-module~$F$, the required
bijection~$\Hom_\afflx(E^\sim,F) \cong \Hom_A(E,\Gamma(X,\F))$ is given by
\[ (\alpha_x : E(x) \to F_x)_{x{:}X} \longmapsto
  (v \mapsto (\alpha_x([v \otimes 1]))_{x{:}X}) \]
with inverse given by
\[ (v \otimes f \mapsto f \cdot \beta(v)_x) \longmapsfrom
  \beta. \qedhere \]
\end{proof}

It's fruitful to study the unit of the adjunction~$(\placeholder)^\sim \dashv
\Gamma(X,\placeholder)$ in more detail.

\begin{defn}An~$A$-module~$E$ has \emph{cohesive fibers} if and only if the
canonical linear map~$E \to \prod_{x{:}X} E(x)$ is a bijection.
\end{defn}

Ordinary modules in classical logic almost never have cohesive fibers.
Intuitively, the law of excluded middle allows to define ``discontinuous''
elements of~$\prod_{x{:}X} E(x)$. These can't be induced from elements~$v\?E$,
whose values~$[v \otimes 1] \? E(x)$ always ``vary continuously'' in~$x\?X$.
In particular, the trivial~$A$-module~$A$ won't have cohesive fibers in a
classical context, since~$A$ might not be a product of copies of the base ring.

However, in the intuitionistic context of the internal language of the big
Zariski topos, the law of excluded middle is not available to define
discontinuous families of values. And indeed, we have the following
proposition.

\begin{prop}\label{prop:cohesive-fibers}
Assume that there is a quasicoherent~$\O_S$-algebra~$\A_0$ such
that~$A = \A_0^\Zar$. Further assume that~$\brak{\Spec(A)}$ is locally
representable by an object of the site used to define~$\Zar(S)$.
Let~$\E_0$ be a quasicoherent~$\A_0$-module. Then~$(\E_0)^\Zar$ has cohesive
fibers from the internal point of view of~$\Zar(S)$.
\end{prop}

\begin{proof}This just reflects the fact that the canonical morphism~$\E_0 \to
f_*((\E_0)^\sim)$, where~$f : \RelSpec_S(\A_0) \to S$ is the structure morphism
of the relative spectrum and~$(\E_0)^\sim$ is the result of applying the
(ordinary) relative tilde construction to~$\E_0$, is an
isomorphism~\stacksproject{01SB}.
\end{proof}

\begin{cor}\label{cor:higher-typed-quasicoherence}
Let~$\Zar(S)$ be defined using one of the parsimonious sites. Assume that~$A$
is a finitely presented~$\affl$-algebra. Then the following statements hold
from the internal point of view of~$\Zar(S)$.
\begin{enumerate}
\item The algebra~$A$ satisfies the following ``higher-typed'' version of the
synthetic quasicoherence condition: Any~$A$-module which is synthetically
quasicoherent as an~$\affl$-module (for instance because it is finitely
presented) has cohesive fibers.
\item For any~$\afflx$-module~$F = (F_x)_{x{:}X}$ such that~$F_x$ is a synthetically
quasicoherent~$\affl$-module for all~$x\?X$, the canonical
map~$\Gamma(X,F)^\sim \to F$ is an isomorphism.
\end{enumerate}
\end{cor}

\begin{proof}The first claim is immediate from
Proposition~\ref{prop:cohesive-fibers} and Theorem~\ref{thm:qcoh-big-char}. The
second claim is an internal rendition of~\stacksproject{01SB}.
\end{proof}

\begin{rem}Since the internal language of~$\Set/\Spec(A)$ is just ordinary
language applied to all fibers, the following properties of
an~$\afflx$-module~$F = (F_x)_{x{:}X}$ are equivalent:
\begin{enumerate}
\item For any~$x\?X$, the~$\affl$-module~$F_x$ is synthetically quasicoherent.
\item From the point of view of~$\Set/\Spec(A)$, the~$\afflx$-module~$F$ is
synthetically quasicoherent.
\end{enumerate}
\end{rem}

We believe that the condition that~$A$ is synthetically quasicoherent
is not strong enough to allow an internal proof of
Corollary~\ref{cor:higher-typed-quasicoherence}. The conclusion of
Corollary~\ref{cor:higher-typed-quasicoherence} seems quite natural from a
synthetic point of view; we therefore propose to adopt it as an axiom for
synthetic algebraic geometry.

\subsection{Open immersions}

A basic concept in the functor-of-points approach to algebraic geometry is the
concept of an \emph{open subfunctor}. It is used to delimit schemes from more
general kinds of spaces: A functor is deemed to be a scheme if and only if it
admits a covering by open subfunctors which are representable.

The following definition is phrased in such a way as to apply to any of the
several ways to define the big Zariski topos~$\Zar(S)$. In particular, it
applies to the definition using the site consisting of affine schemes which are
locally of finite presentation over~$S$. If~$S$ is affine, the definition only
refers to affine schemes and open subschemes of affine schemes and is therefore
suitable if one wants to found the theory of schemes using the functorial
approach.

\begin{defn}[{\cite[Définition~I.1.3.6 on page~10]{demazure:gabriel},
\cite[Tag~01JI]{stacks-project}}]
A subfunctor~$U \hookrightarrow X$ in~$\Zar(S)$ is an \emph{open subfunctor} if
and only if for any object~$(T \to S)$ of the site used to define~$\Zar(S)$
and any~$x \in X(T)$ there exists an open subscheme~$T_0 \subseteq T$
such that for any object~$(T' \xra{f} T \to S)$ of the site used to
define~$\Zar(S)$ the map~$T' \to T$ factors over~$T_0$ if and only if~$X(f)(x)
\in U(T')$.
\end{defn}

The open subschemes~$T_0 \subseteq T$ appearing in this definition are uniquely
determined by their universal property. The relation of open subfunctors to
open immersions is as follows.

\begin{prop}\label{prop:char-open-immersion}
Let~$X$ be an~$S$-scheme.
\begin{enumerate}
\item Let~$U \subseteq X$ be an open subscheme. Then the subfunctor~$\ul{U}
\hookrightarrow \ul{X}$ is open.
\item If~$\ul{X}$ is locally representable by an object of the site used to
define~$\Zar(S)$, any open subfunctor~$U \hookrightarrow \ul{X}$ is isomorphic
to the open subfunctor associated to an open subscheme of~$X$.
\end{enumerate}
\end{prop}

\begin{proof}For the first claim, let~$(T \to S)$ be an object of the site used
to define~$\Zar(S)$ and let~$x \in \ul{X}(T)$. The open
subscheme~$T_0 \subseteq T$ required by the definition of an open subfunctor can
then be chosen as~$x^{-1}[U]$.

For the second claim, assuming for notational simplicity that~$\ul{X}$ is
directly representable without having to pass to a cover, the desired open
subscheme of~$X$ can be obtained as the witnessing subscheme~``$T_0$'' as it
appears in the definition of an open subfunctor in the special case~$(T \to S)
\defeq (X \to S)$.
\end{proof}

From the point of view of the internal language of~$\Zar(S)$, a subfunctor~$U
\hookrightarrow X$ looks like the inclusion of a subset. The natural question
how one can characterize those inclusions which externally correspond to open
subfunctors is answered as follows.

\begin{defn}In the context of a specified local ring~$R$, as for instance~$\affl$
of the big Zariski topos of a scheme, a truth value~$\varphi$ is
\emph{open} if and only if there exists an ideal~$J \subseteq R$ such
that~$R/J$ is synthetically quasicoherent
(Definition~\ref{defn:synth-qcoh}) and such that~$\varphi$ holds if and only if~$1 \in
J$. (Section~\ref{sect:basics-on-truth-values} contains generalities on truth
values.)\end{defn}

\begin{ex}Let~$f \? R$. Then~``$f$ is invertible'' is an open truth value with
witnessing ideal~$J = (f)$. The quotient~$R/J$ is indeed synthetically
quasicoherent, since it is finitely presented. More generally,
let~$f_1,\ldots,f_n \? R$. Then~``one of the~$f_i$ is invertible''
is an open truth value with witnessing ideal~$J =
(f_1,\ldots,f_n)$.

In case that~$R$ fulfills the same field condition as~$\affl$ does,
one can write this truth value also as~``$f_1 \neq 0 \vee \cdots \vee f_n \neq
0$''.
\end{ex}

\begin{defn}In the context of a specified local ring, a map~$U \to X$ of sets is a
\emph{synthetic open immersion} if and only if it is injective and for any~$x\?X$
the truth value of~``the fiber of~$x$ is inhabited'' is open.\end{defn}

\begin{ex}The inclusion~$R^\times \hookrightarrow R$ of the
invertible elements is a synthetic open immersion, since for~$x\?R$ the
truth value of~``the fiber of~$x$ is inhabited'' equals the truth value of~``$x$
is invertible''.\end{ex}

\begin{ex}Let~$X$ be a set. Let~$f : X \to R$ be a function.
The inclusion~$\{ x\?X \,|\, \text{$f(x)$ is invertible} \} \hookrightarrow X$ is a
synthetic open immersion.\end{ex}

\begin{prop}\label{prop:characterization-open-subfunctor}
Let~$X \in \Zar(S)$ be a Zariski sheaf. If a subfunctor~$U
\hookrightarrow X$ is open, then the map~$U \hookrightarrow X$
is a synthetic open immersion from the internal point of view of~$\Zar(S)$,
that is
\begin{multline*}
  \qquad\Zar(S) \models
  \forall x\?X\_
  \exists J \subseteq \affl\_ \\
  \speak{$J$ is an ideal} \wedge
  \speak{$\affl/J$ is synth.\@ quasicoherent} \wedge
  (x \in U \Leftrightarrow 1 \in J). \end{multline*}
The converse holds if one of the parsimonious sites is used to
define~$\Zar(S)$ or if the ideals~$J$ are required to be finitely generated
from the internal point of view.
\end{prop}

We postpone the proof of this proposition in order to give a bit of context
first.

Firstly, the displayed condition is only meaningful in an intuitionistic context as
provided by the big Zariski topos. In classical logic, the condition would be
trivially satisfied for any subfunctor~$U \hookrightarrow X$: Classically, we
have~$(x \in U) \vee (x \not\in U)$. If~$x \in U$, we can pick~$J = (1)$, and
if~$x \not\in U$, we can pick~$J = (0)$ (whereby the quotient~$\affl/J$ is
isomorphic to~$\affl$, thus finitely presented and therefore in particular
synthetically quasicoherent).\footnote{Strictly speaking, incompatibility with classical
logic surfaces even earlier: in our synthetic quasicoherence condition. The
map~$E \otimes_\affl A \to [\Spec(A), E]$ which the condition demands to be
bijective has hardly any chance to be surjective if the law of excluded middle
is available to define maps~$\Spec(A) \to E$ by case distinction.}

Proposition~\ref{prop:characterization-open-subfunctor} is often used in the
following form, which is weaker because it
only gives one direction, but which is applicable for any of our choices for
the site used to define~$\Zar(S)$.

\begin{cor}\label{cor:sufficient-criterion-open-subfunctor}
Let~$X \in \Zar(S)$ be a Zariski sheaf. Let~$U \hookrightarrow X$ be
a subfunctor. If
\[
  \Zar(S) \models
  \forall x\?X\_
  \bigvee_{n \geq 0} \exists f_1,\ldots,f_n\?\affl\_
  (x \in U \Leftrightarrow \bigvee_i \speak{$f_i$ \inv}),
\]
then the subfunctor is open.
\end{cor}

\begin{proof}
We show that the assumption implies the displayed condition of
Proposition~\ref{prop:characterization-open-subfunctor} in the internal
language. Given elements~$f_1,\ldots,f_n$ as in the assumption, we
construct the ideal~$J \defeq (f_1,\ldots,f_n) \subseteq \affl$. The
quotient~$\affl/J$ is finitely presented, hence synthetically quasicoherent, and
the statement that~$1 \in J$ is equivalent to one of
the~$f_i$ being invertible by locality of~$\affl$.
\end{proof}

The internal condition appearing in
Corollary~\ref{cor:sufficient-criterion-open-subfunctor} reflects basic
intuition about openness in algebraic geometry: Intuitively, a subset is open
if it is given by inequations, so that to decide whether a point belongs to the
subset one has to check that at least one of some numbers is not zero.

Of course, in classical scheme theory, one would put some condition on these
numbers in order not to trivialize the notion. For instance, one would require
that they depend continuously on the point in some sense or, more
specifically, that these numbers are given by evaluating certain locally
defined regular functions at the point.

On first sight, such a condition seems to be lacking in
Corollary~\ref{cor:sufficient-criterion-open-subfunctor}. However, it's
implicitly built into the language, since by the Kripke--Joyal semantics the
external meaning of~``$\exists f\?\affl$'' is that there exist, locally on an
open cover, suitable elements of~$\affl(T)$, that is regular functions on~$T$.

It's useful to give a name for the kind of the subfunctors appearing in
Corollary~\ref{cor:sufficient-criterion-open-subfunctor}.

\begin{defn}\begin{itemize}
\item A subfunctor~$U \hookrightarrow X$ in~$\Zar(S)$ is a \emph{quasicompact open subfunctor} if
and only if for any object~$(T \to S)$ of the site used to define~$\Zar(S)$
and any~$x \in X(T)$ there exists an open subscheme~$T_0 \subseteq T$
such that the open immersion~$T_0 \hookrightarrow T$ is quasicompact and such
that for any object~$(T' \xra{f} T \to S)$ of the site used to define~$\Zar(S)$
the map~$T' \to T$ factors over~$T_0$ if and only if~$X(f)(x) \in U(T')$.
\item In the context of a specified local ring, a truth value~$\varphi$ is
\emph{quasicompact open} if and only if there exists a finitely generated
ideal~$J \subseteq R$ such that~$\varphi$ holds if and only if~$1 \in
J$.
\item In the context of a specified local ring, a map~$U \to X$ of sets is a
\emph{synthetic quasicompact open immersion} if and only if it is injective and for any~$x\?X$
the truth value of~``the fiber of~$x$ is inhabited'' is quasicompact open.
\end{itemize}
\end{defn}

\begin{proof}[Proof of Proposition~\ref{prop:characterization-open-subfunctor}]
We begin with the ``only if'' direction. Let~$T$ be an~$S$-scheme contained in
the site used to define~$\Zar(S)$ and let~$x \in X(T)$. By assumption there is
an open subscheme~$T_0 \subseteq T$ such that, for any object~$(T' \xra{f} T
\to S)$ of the site the morphism~$f$ factors over~$T_0$ if and only if~$X(f)(x)
\in U(T')$.

There is a quasicoherent sheaf~$\I$ of ideals on~$T$ such that~$T_0 = D(\I)$.
We set~$J \defeq \im(\I^\Zar \to \afflt)$. Then the quotient module~$\afflt/J$
coincides with~$(\O_T/\I)^\Zar$ and is therefore quasicoherent, hence
synthetically quasicoherent from the internal point of view of~$\Zar(T)$ by
Theorem~\ref{thm:qcoh-big-char}.

To verify that~$T \models (x \in U \Leftrightarrow 1 \in J)$, let an
object~$(T' \xra{f} T \to S)$ of the site be given. Then we have the
chain of equivalences
\begin{align*}
  &\mathrel{\phantom{\Longleftrightarrow}}
  X(f)(x) \in U(T') \\
  &\Longleftrightarrow
  \text{$T' \to T$ factors over~$T_0$} \\
  &\Longleftrightarrow
  T' \models 1 \in \im(f^*\I \to \O_{T'}) \\
  &\Longleftrightarrow
  T' \models 1 \in J,
\end{align*}
where the last step follows by Lemma~\ref{lemma:qcoh-ideals}.

For the converse direction, let~$T$ be an~$S$-scheme contained in the site and
let~$x \in X(T)$. By assumption there is an open covering~$T = \bigcup_i V_i$
and ideals~$J_i \in \Zar(V_i)$ such that~$\afflvi/J_i$ is synthetically
quasicoherent and such that~$V_i \models (x \in U \Leftrightarrow 1 \in J_i)$.

By Theorem~\ref{thm:qcoh-big-char}, there are quasicoherent modules~$\E_i$
over~$V_i$ such that~$\afflvi/J_i \cong (\E_i)^\Zar$. We set~$\I_i \defeq
\Kernel(\O_{V_i} \to \E_i)$. One can check that~$D(\I_i) \cap V_j =
D(\I_j) \cap V_i$ for all~$i$ and~$j$. We set~$T_0 \defeq \bigcup_i D(\I_i)$.

To verify that this construction satisfies what is expected of it, let~$(T'
\xra{f} T \to S)$ be an object of the site. We then have the following chain of
equivalences:
\begin{align*}
  &\mathrel{\phantom{\Longleftrightarrow}}
  \text{$T' \to T$ factors over~$T_0$} \\
  &\Longleftrightarrow
  \text{for all~$i$, $f_i : f^{-1}V_i \to V_i$ factors over~$D(\I_i)$} \\
  &\Longleftrightarrow
  \text{for all~$i$, $f^{-1}V_i \models 1 \in \im(f_i^*\I_i \to \O_{f^{-1}V_i})$} \\
  &\stackrel{\star}{\Longleftrightarrow}
  \text{for all~$i$, $f^{-1}V_i \models 1 \in J_i$} \\
  &\Longleftrightarrow
  \text{for all~$i$, $f^{-1}V_i \models x \in U$} \\
  &\Longleftrightarrow
  \text{$T' \models x \in U$} \\
  &\Longleftrightarrow
  \text{$X(f)(x) \in U(T')$}
\end{align*}
For this calculation, it's not important that the ideals~$\I_i$ are
quasicoherent. The assumption that the quotient modules~$\afflvi/J_i$ are
quasicoherent is only needed to ensure that Lemma~\ref{lemma:qcoh-ideals}
can be applied in the marked step.
\end{proof}

\begin{prop}\label{cor:sufficient-criterion-qc-open-subfunctor}
Let~$X \in \Zar(S)$ be a Zariski sheaf. A subfunctor~$U \hookrightarrow X$ is a
quasicompact open subfunctor if and only if~$U \hookrightarrow X$ is a
synthetic quasicompact open immersion from the internal point of view
of~$\Zar(S)$.
\end{prop}

\begin{proof}The proof of
Proposition~\ref{prop:characterization-open-subfunctor} can be adapted.
\end{proof}

\begin{lemma}\label{lemma:open-truth-values-stable}
Finite conjunctions and finite disjunctions of quasicompact open
truth values are quasicompact open. If~$\Zar(S)$ is defined using one of the
parsimonious sites, then furthermore finite conjunctions and finite
disjunctions of open truth values are open from the internal point of view
of~$\Zar(S)$.
\end{lemma}

\begin{proof}Let~$\varphi$ and~$\psi$ be quasicompact open truth values with
witnessing finitely generated ideals~$I$ and~$J$. Because the ring is
local,~$\varphi \wedge \psi$ is equivalent to~$1 \in I \cdot J$, and~$\varphi \vee
\psi$ is equivalent to~$1 \in I + J$. The ideals~$I \cdot J$ and~$I + J$ are
finitely generated. Therefore~$\varphi \wedge \psi$ and~$\varphi \vee \psi$ are
quasicompact open truth values.

For the case of arbitrary open truth values, we need to verify that~$\affl/(I
\cdot J)$ and~$\affl/(I + J)$ are synthetically quasicoherent if~$\affl/I$
and~$\affl/J$ are. This second claim from Lemma~\ref{lemma:tensor-product-qcoh},
since~$\affl/(I + J) \cong \affl/I \otimes_\affl \affl/J$. The first claim
follows from a similar calculation.
\end{proof}

The notion of open truth values is not unique to our account of synthetic
algebraic geometry. Rather, it's a concept in the established and more general
framework of synthetic topology~\cite{escardo,lesnik} which aims to do topology in a
synthetic fashion: Any set should have an intrinsic topology and any map should
be automatically continuous with respect to this intrinsic topology.

This automatic continuity reflects as stability of open subfunctors under
pullbacks:

\begin{lemma}Let~$f : X \to Y$ be a morphism in$~\Zar(S)$. Let~$U
\hookrightarrow Y$ be an open subfunctor. Then its pullback along~$f$,
denoted~``$f^{-1}U \hookrightarrow X$'', is too an open subfunctor.\end{lemma}

\begin{proof}From the internal point of view of~$\Zar(S)$, the
subfunctor~$f^{-1}U \hookrightarrow X$ looks like the inclusion of the
preimage~$f^{-1}[U] \subseteq X$.

So, to verify the claim, let internally an element~$x\?X$ be given. We are to
show that the truth value of~``$x \in f^{-1}[U]$'' is open. This truth value
equals the truth value of~``$f(x) \in U$'' which is open by assumption, and is
therefore open.
\end{proof}

\begin{rem}In the internal language of toposes used to carry out synthetic
differential geometry, there is the concept of an \emph{Penon-open}
subset~\cite[Chapitre~III]{penon}: A subset~$U \subseteq X$ is Penon-open if
and only if
\[ \forall x \in U\_ \forall y\?X\_
  (x \neq y) \vee (y \in U). \]
This notion is not useful in synthetic algebraic geometry, since it is much too
weak: Any subset of the one-element set~$1$ is Penon-open. However, not every
subfunctor of the terminal functor in~$\Zar(S)$ is an open subfunctor.
\end{rem}

In many flavors of synthetic topology, open truth values~$\varphi$
are~$\neg\neg$-stable in that~$\neg\neg\varphi$ implies~$\varphi$. With a
small caveat, this is true for open truth values in the big Zariski topos as
well.

\begin{prop}\label{prop:open-truth-values-stable}
Let~$U \hookrightarrow 1$ be a subfunctor in~$\Zar(S)$ such
that~$\Zar(S) \models \neg\neg U$. Then in any of the following
situations it follows that~$\Zar(S) \models U$:
\begin{enumerate}
\item $U$ is a quasicompact open truth value.
\item $U$ is an arbitrary open truth value and the site used to
define~$\Zar(S)$ is closed under domains of closed immersions. (This is for
instance satisfied for the sites built using a Grothendieck or a partial
universe. It is satisfied for the parsimonious sites if~$S$ is locally
Noetherian.)
\end{enumerate}
\end{prop}

\begin{proof}We give two proofs, an internal one and an external one, since
they employ different ideas.

\emph{Internal proof.} Since~$U$ is an open truth value, there exists an ideal~$J
\subseteq \affl$ such that~$\affl/J$ is synthetically quasicoherent and such that~$U$ holds
if and only if~$1 \in J$. By assumption, the element~1 is \notnot an element
of~$J$; we want to verify that it's actually an element of~$J$.

By Scholium~\ref{scholium}, the canonical homomorphism
\[ \affl/J \longrightarrow [\Spec(\affl/J), \affl] \]
is bijective; the assumptions of that scholium are satisfied in either of
the two situations. The set~$\Spec(\affl/J)$ is isomorphic to~$\brak{J = (0)}$.
Since~$\neg\neg(1 \in J)$, we also have~$\neg(J = (0))$.
Therefore~$\Spec(\affl/J)$ is empty and the codomain of the displayed
isomorphism is the zero algebra. Thus~$\affl/J$ is trivial as well, showing~$1
\in J$.

\emph{External proof.} Since~$U \hookrightarrow 1$ is an open subfunctor, there
is an open subscheme~$S_0 \subseteq S$ such that a morphism~$f : T \to S$
factors over~$S_0$ if and only if~$U(T)$ is inhabited. In both situations it's
possible to endow~$X \defeq S \subseteq S_0$ with the structure of a closed
subscheme such that~$X$ is contained in the site used to define~$\Zar(S)$.
By the universal property of~$S_0$, we have~$X \models \neg U$. Since~$\Zar(S)
\models \neg\neg U$, it follows that~$X$ is empty. Therefore~$S_0 = S$ and~$U$
is globally inhabited.
\end{proof}

\begin{cor}\label{cor:open-subfunctors-d1}
Let~$\gamma : \Delta \to X$ be a morphism in~$\Zar(S)$.
Let~$U \hookrightarrow X$ be an open subfunctor such that~$\Zar(S) \models
\gamma(0) \in U$. Then, in any of the situations in
Proposition~\ref{prop:open-truth-values-stable}, the morphism~$\gamma$ factors
over~$U$.
\end{cor}

\begin{proof}We give an internal proof. Let~$\varepsilon \in \Delta$.
Then~$\neg\neg(\varepsilon = 0)$. Therefore~$\neg\neg(\gamma(\varepsilon) \in
U)$. Since being an element of~$U$ is~$\neg\neg$-stable, it follows
that~$\gamma(\varepsilon) \in U$.
\end{proof}

\begin{rem}Subobjects~$U \hookrightarrow X$ for which any morphism~$\gamma :
\Delta \to X$ with~$\gamma(0) \in U$ factors over~$U$ are called~``$D_1$-open''
in the literature on synthetic differential
geometry~\cite[p.~60]{reyes:wraith:note-tangent-bundles}.
Corollary~\ref{cor:open-subfunctors-d1} shows that open subfunctors
are~$D_1$-open.
\end{rem}

In ordinary scheme theory, an inclusion of a standard open subset~$D(f)
\hookrightarrow X$ is isomorphic to the structure morphism of the relative
spectrum~$\RelSpec_X \O_X[f^{-1}]$. Inclusions of more general open subsets
can typically not be described using the relative spectrum construction, the
standard example being the inclusion~$\AA^2_k \setminus \{ 0 \} \hookrightarrow
\AA^2_k$ whose domain is not affine.

An interesting feature of the internal universe of the big Zariski topos is
that it's flexible enough to express \emph{any} open subset as a spectrum.
The contradiction is only apparent since the algebra used for constructing
such a spectrum is not in general quasicoherent.

\begin{prop}Let~$U \hookrightarrow 1$ be an open truth value. In any of the
situations of Proposition~\ref{prop:open-truth-values-stable}, there is a
(not necessarily quasicoherent)~$\affl$-algebra~$A$ such that the inclusion is
isomorphic to the morphism~$\Spec(A) \to 1$.\end{prop}

\begin{proof}The open truth value~$U$ is given by an ideal~$J \subseteq \affl$
such that~$\affl/J$ is synthetically quasicoherent and such that~$U$ holds if
and only if~$1 \in J$. We set~$A \defeq \affl[M^{-1}]$, where~$M$ is the
multiplicatively closed subset
\[ M \defeq \{ f \? \affl \,|\, 1 \in J \Rightarrow \speak{$f$ \inv} \} \subseteq \affl. \]
The spectrum of~$A$ is inhabited if and only if~$M \subseteq (\affl)^\times$,
in which case the unique element of~$\Spec(A)$ is the inverse of the
localization morphism~$\affl \to \affl[M^{-1}]$. Thus~$\Spec(A)$ is isomorphic
to~$\brak{M \subseteq (\affl)^\times}$. Therefore we have to verify that~$U$
holds if and only if~$M \subseteq (\affl)^\times$.

The ``only if'' direction is trivial.

For the ``if'' direction, we exploit the~$\neg\neg$-stability of~$U$. If~$\neg
U$, then~$\neg(1 \in J)$, so~$M = \affl$, and since~$M \subseteq (\affl)^\times$ by
assumption, it follows that zero is invertible. This is a contradiction.
Thus~$\neg\neg U$.
\end{proof}

\begin{rem}The radical~$\sqrt{J}$ of the ideal~$J$ appearing in
Proposition~\ref{prop:characterization-open-subfunctor} is unique:
It is equal to the radical ideal
\[ K \defeq \{ f\?\affl \,|\, \speak{$f$ \inv} \Rightarrow (x \in U) \}
\subseteq \affl. \]
It's obvious that~$J \subseteq K$ and therefore~$\sqrt{J} \subseteq K$.
For the converse direction, let~$f \in K$ be given. Since~$\affl/J$ is
synthetically quasicoherent, the canonical map
\[ (\affl/J)[f^{-1}] \longrightarrow [ \Spec(\affl[f^{-1}]), \affl/J ] \]
is bijective. Since~$\Spec(\affl[f^{-1}]) \cong \brak{\speak{$f$ \inv}}$, the image
of~$1$ is zero: If~$\Spec(\affl[f^{-1}])$ is inhabited, the element~$f$ is invertible and
therefore~$x$ is an element of~$U$. This implies that~$1 \in J$.
Thereby~$\affl/J = 0$. By injectivity of the canonical map, the
algebra~$(\affl/J)[f^{-1}]$ is zero. Therefore~$f^n \in J$ for some natural
number~$n$.\end{rem}

\begin{rem}In view of the previous remark, one might hope to be able to simplify the
condition in Proposition~\ref{prop:characterization-open-subfunctor} as
follows: ``For any~$x\?X$, the quotient~$\affl/K$ modulo the ideal~$K = \{
f\?\affl \,|\, \speak{$f$ \inv} \Rightarrow (x \in U) \}$ is synthetically quasicoherent.''
However, this doesn't work out.
This statement implies the condition in the proposition, but the converse
direction does not hold, since~$\affl/K \cong \affl/\sqrt{J}$ might fail to be
synthetically quasicoherent. For instance that's the case if~$U = \emptyset$; then~$K =
\sqrt{(0)}$ by Proposition~\ref{prop:a1-nilp}. The
quotient~$\affl/\sqrt{(0)}$ is not synthetically quasicoherent by
Remark~\ref{rem:radical-not-qcoh}.\end{rem}

\begin{rem}\label{rem:open-geometric-morphism}
There is the notion of an open geometric morphism of toposes. For the big
Zariski toposes, this notion is not related to open morphisms or open
immersions between schemes: If~$X \to S$ is any morphism of schemes, the
induced geometric morphism~$\Zar(X) \to \Zar(S)$ is isomorphic to the canonical
morphism~$\Zar(S)/\ul{X} \to \Zar(S)$, as detailed in
Section~\ref{sect:change-of-base}. Geometric morphisms of the form~$\E/A \to
\E$ are always open.\end{rem}

\subsection{Closed immersions}

\begin{defn}In the context of a specified local ring~$R$, as for instance~$\affl$
of the big Zariski topos of a scheme, a truth value~$\varphi$ is
\emph{closed} if and only if there exists an ideal~$J \subseteq R$ such
that~$R/J$ is synthetically quasicoherent
(Definition~\ref{defn:synth-qcoh}) and such that~$\varphi$ holds if and only if~$J =
(0)$.\end{defn}

In other words, a truth value~$Z \subseteq 1$ is closed if and only if~$Z$ is
isomorphic to the spectrum of a synthetically quasicoherent quotient algebra
of~$R$.

\begin{ex}Let~$f \? R$. Then~``$f = 0$'' is a closed truth value with
witnessing ideal~$J = (f)$. More generally, if~$f_1,\ldots,f_n \? R$, the
truth value~``$f_1 = \cdots = f_n = 0$'' is closed with witnessing
ideal~$(f_1,\ldots,f_n)$.\end{ex}

\begin{defn}In the context of a specified local ring, a map~$Z \to X$ of sets is a
\emph{synthetic closed immersion} if and only if it is injective and for
any~$x\?X$ the truth value of~``the fiber of~$x$ is inhabited'' is
closed.\end{defn}

\begin{ex}The inclusion~$\{0\} \hookrightarrow R$ is a synthetic closed
immersion. More generally, for any functions~$f_1,\ldots,f_m : R^n \to
R$, the inclusion of the set of their common zeros in~$R^n$ is a
synthetic closed immersion.\end{ex}

\begin{ex}Let~$X$ be a set. Let~$f : X \to R$ be a function.
The inclusion~$\{ x\?X \,|\, f(x) = 0 \} \hookrightarrow X$ is a
synthetic closed immersion.\end{ex}

\begin{prop}\label{prop:char-closed-immersion}
Let~$X$ be an~$S$-scheme.
\begin{enumerate}
\item Let~$Z \hookrightarrow X$ be a closed subscheme. Then the subfunctor~$\ul{Z}
\hookrightarrow \ul{X}$ is a synthetic closed immersion from the internal point
of view of~$\Zar(S)$.
\item Assume that~$\ul{X}$ is locally representable by an object of the site used to
define~$\Zar(S)$. Let~$Z \hookrightarrow \ul{X}$ be a synthetic closed
immersion. If the site is one of the parsimonious sites or if the witnessing
ideals for the immersion are finitely generated from the internal point of
view, then~$Z \hookrightarrow \ul{X}$ is
isomorphic to the subfunctor associated to a closed subscheme of~$X$.
\end{enumerate}
\end{prop}

\begin{proof}To verify the first claim, let a
quasicoherent~$\O_X$-algebra~$\J_0$ be given such that the closed subscheme~$Z
\hookrightarrow X$ is the vanishing scheme of~$\J_0$. Following the translation with
the Kripke--Joyal semantics, let~$f : T \to S$ be an object of the site used to
define~$\Zar(S)$ and let~$x \in \ul{X}(T)$. We define~$J \defeq (f^*
\J_0)^\Zar \in \Zar(S)/\ul{T}$. Then $T \models \speak{$\affl/J$ is
synthetically quasicoherent}$ and $T \models (x \in \ul{Z} \Leftrightarrow J =
(0))$, therefore~``$x \in \ul{Z}$'' is a closed truth value.

For the converse direction, we may assume that~$X = S$ since~$\Zar(S)/\ul{X}
\simeq \Zar(X)$, as discussed in Section~\ref{sect:change-of-base}. We then
observe that the problem is local on~$S$, since we can glue matching schemes
defined over the members of an open covering of~$S$~\stacksproject{01JJ}. We
may therefore assume that we are given an~$\affl$-algebra~$J$ such
that~$\affl/J$ is synthetically quasicoherent from the internal point of view
and such that~$Z = \brak{J = (0)}$. The assumptions ensure that
Theorem~\ref{thm:qcoh-big-char} is applicable. Therefore there is a
quasicoherent~$\O_S$-module~$\E$ such that~$\affl/J \cong \E^\Zar$. We set~$\J
\defeq \Kernel(\O_S \twoheadrightarrow \E)$. Lemma~\ref{lemma:qcoh-ideals} then
implies that~$Z$ is the functor of points of the~$S$-scheme~$V(\J)$.
\end{proof}

\begin{prop}\label{prop:closed-truth-value-ideal}
Let~$\Zar(S)$ be defined using a parsimonious site and assume
that~$S$ is locally Noetherian. Then the witnessing ideals of closed truth
values are uniquely determined by the truth value.
\end{prop}

\begin{proof}We argue internally. Let~$\varphi$ be a closed truth value.
Let~$J \subseteq \affl$ be an ideal such that~$\affl/J$ is synthetically
quasicoherent and such that~$\varphi$ holds if and only if~$J = (0)$.

We set~$K \defeq \{ f \? \affl \,|\, \varphi \Rightarrow f = 0 \}$. Then,
trivially,~$J \subseteq K$. For the converse containment relation, let~$f \in
K$. The canonical homomorphism~$\affl/J \to [\brak{\varphi}, \affl]$ is
bijective by synthetic quasicoherence of~$\affl/J$ and by
Scholium~\ref{scholium} (it is here that we need that~$S$ is locally Noetherian
-- else we can't ensure that the~$S$-scheme~$\RelSpec_S((\affl/J)|_{\Sh(S)})$
is locally finitely presented). The image of~$[f]$ is zero, hence~$f$ is an
element of~$J$.
\end{proof}

\begin{rem}Without any assumptions on~$S$ or on the site used to
define~$\Zar(S)$, Proposition~\ref{prop:closed-truth-value-ideal} holds for
those closed truth values which admit a finitely generated witnessing ideal.
For those closed truth values, any given finitely generated witnessing ideals
are equal.
\end{rem}

\begin{lemma}If~$\Zar(S)$ is defined using one of the parsimonious sites, then
finite conjunctions of closed truth values are closed from the internal point
of view of~$\Zar(S)$.
\end{lemma}

\begin{proof}We argue internally. Let~$\varphi$ and~$\psi$ be closed truth
values with witnessing ideals~$I$ and~$J$. Then~$\varphi \wedge \psi$ is
equivalent to~$I + J = (0)$ and the module~$\affl/(I + J)$ is synthetically
quasicoherent (as in the proof of Lemma~\ref{lemma:open-truth-values-stable}).
Hence~$\varphi \wedge \psi$ is a closed truth value.
\end{proof}

It's in general not the case that finite disjunctions of closed truth
values are closed, and it's instructive to see why. The external interpretation
of this failure is the following: Let~$A \hookrightarrow S$ and~$B
\hookrightarrow S$ be closed subschemes. Then the functor of points of the
union~$A \cup B$ does \emph{not} coincide with the union of the
subfunctors~$\ul{A} \hookrightarrow \ul{S}$ and~$\ul{B} \hookrightarrow
\ul{S}$. An explicit description of the former functor is
\[ (X \xra{f} S) \longmapsto
  \{ \star \,|\, \text{$f$ factors over~$A \cup B$} \} \]
and of the latter is
\[ (X \xra{f} S) \longmapsto
  \{ \star \,|\, \text{locally on the target, $f$ factors over~$A$ or over~$B$} \}. \]
The inclusion~$A \cup B \hookrightarrow S$ trivially factors over~$A \cup B$,
but in general there isn't an open covering~$S = \bigcup_i U_i$ such that for
each~$i$, the restriction~$(A \cup B) \cap U_i \hookrightarrow U_i$ factors
over~$A \cap U_i$ or over~$B \cap U_i$. For instance, this isn't the case
if~$S$ is the affine plane over a field and~$A$ and~$B$ are the two axes.

This phenomenon doesn't happen for open subschemes, which explains why finite
disjunctions of open truth values are open.

\begin{rem}\label{rem:closed-geometric-morphism}
There is the notion of a closed geometric morphism of toposes. For
an arbitrary topos~$\E$ and an object~$A \in \E$, the canonical geometric
morphism~$\E/A \to E$ is closed if and only if
\[ \forall U \subseteq A\_
  \forall \varphi \? \Omega\_ \quad
  A \subseteq (U \cup \{ x \in A \,|\, \varphi \}) \quad\Longrightarrow\quad
  (A \subseteq U) \vee \varphi \]
from the internal point of view of~$\E$. If~$X \to S$ is a closed morphism of
schemes, then the induced geometric morphism~$\Sh(X) \to \Sh(S)$ between the
little Zariski toposes is closed in this sense.

However, the induced geometric morphism~$\Zar(X) \to \Zar(S)$ is typically not
closed. For instance, if~$X \to S$ is the embedding of a closed subset~$V(f)$
with~$f \in \Gamma(S,\O_S)$, then the morphism~$\Zar(X) \to \Zar(S)$ is
isomorphic to~$\Zar(S)/\ul{V(f)} \to \Zar(S)$, as discussed in
Section~\ref{sect:change-of-base}. In the special case~$U \defeq \emptyset$ and
$\varphi \defeq \brak{f = 0}$, the displayed closedness condition simplifies
to~$\speak{$f$ \inv} \vee (f = 0)$. This is typically not true in the internal
language of~$\Zar(S)$. A specific counterexample is given in
Example~\ref{ex:translation-equivalence}.
\end{rem}

\subsection{Surjective morphisms}

\begin{prop}\label{prop:char-surjective-morphisms}
Let~$f : X \to S$ be an arbitrary~$S$-scheme. Consider the following
statements:
\begin{enumerate}
\item The morphism~$f$ is surjective.
\item From the internal point of view of~$\Zar(S)$, it's not the case that~$\ul{X}$ is
empty, that is
\[ \Zar(S) \models \neg\neg(\speak{$\ul{X}$ is inhabited}). \]
\end{enumerate}
If~$X$ is locally contained in the site used to define~$\Zar(S)$ (for instance,
if~$X$ is contained in the universe used to define~$\Zar(S)$ or if one of the
parsimonious sites is used and~$X$ is locally of finite presentation over~$S$),
then~(1) implies~(2). The converse holds if the site is closed under
taking spectra of residue fields or if one of the parsimonious sites is used
and~$f$ is quasicompact and quasiseparated.
\end{prop}

\begin{proof}The translation of the internal statement using the Kripke--Joyal
semantics is:
\begin{indentblock}For any~$S$-scheme~$T$ of the site used to define~$\Zar(S)$,
if~\mbox{$\ul{X \times_S T} = \ul{\emptyset}$} (as functors of points
of~$T$-schemes), then~$T = \emptyset$.
\end{indentblock}
In the case that the site used to define~$\Zar(S)$ is closed under taking
spectra of residue fields, this implies that~$f$ is surjective as follows.
Let~$s \in S$ be an arbitrary point. The~$S$-scheme~$T \defeq \Spec(k(s))$ is
not empty. Therefore the fiber~$X_s = X \times_S T$ of~$f$ over~$s$ is not empty.

If one of the parsimonious sites is used to define~$\Zar(S)$, we can't apply
the assumption to the~$S$-scheme~$T = \Spec(k(s))$ since it might not be
locally of finite presentation over~$S$. We therefore argue as follows. Without
loss of generality, we may assume that~$S$ is affine. Writing~$k(s)$ as the
canonical filtered colimit of all finitely presented~$\Gamma(S,\O_S)$-algebras
mapping to~$k(s)$ (and then rewriting this filtered colimit as a directed
colimit~\cite[Theorem~1.5]{adamek:rosicky:presentable}), we see
that~$\Spec(k(s))$ is the directed limit of an inverse system of finitely
presented affine~$S$-schemes~$T_i$ with affine transition maps.
In particular, the structure morphisms~$T_i \to S$ are quasicompact and
quasiseparated. By the assumption that the morphism~$X \to S$ is quasicompact and
quasiseparated as well, the schemes~$X \times_S T_i$ are quasicompact and
quasiseparated (as absolute schemes). Therefore, if~$X_s = X \times_S T
= \lim_i (X \times_S T_i)$ is empty, then so is~$X \times_S T_i$ for
some~$i$~\stacksproject{01ZC}. Thus~$T_i = \emptyset$ and
hence~$\Spec(k(s)) = \emptyset$; this is a contradiction.

For the converse direction, let an~$S$-scheme~$T$ contained in the site used to
define~$\Zar(S)$ be given such that~$\ul{X \times_S T} = \ul{\emptyset}$ as
functors of points of~$T$-schemes. The assumption that~$X$ is locally contained
in the site used to define~$\Zar(S)$ implies that~$X \times_S T = \emptyset$ as
schemes. Since the base change~$X \times_S T \to T$ of~$f$ is surjective, this
implies that~$T$ is empty.
\end{proof}

\begin{cor}\label{cor:char-surjective-morphisms-relative}
Let~$p : X \to Y$ be a morphism of~$S$-schemes. Assume that~$Y$ is locally
contained in the site used to define~$\Zar(S)$.
Further assume that the site used to define the big Zariski toposes are closed
under taking spectra of residue fields or that the parsimonious sites are used
and that~$p$ is quasicompact and quasiseparated.
Consider the following statements:
\begin{enumerate}
\item The morphism~$p$ is surjective.
\item From the internal point of view of~$\Zar(S)$ all fibers of~$\ul{p}$ are
nonempty, that is
\[ \Zar(S) \models \forall y\?\ul{Y}\_
  \neg\neg \exists x\?\ul{X}\_ \ul{p}(x) = y. \]
\end{enumerate}
If~$X$ is locally contained in the
site used to define~$\Zar(Y)$, then~(1) implies~(2). The converse holds if the
sites used to define the big Zariski toposes are closed under taking spectra of
residue fields or that the parsimonious sites are used and that~$p$ is
quasicompact and quasiseparated.
\end{cor}

\begin{proof}Immediate using Proposition~\ref{prop:char-surjective-morphisms}
and the equivalence~$\Zar(Y) \simeq \Zar(S)/\ul{Y}$,
as explained in Section~\ref{sect:change-of-base}.
\end{proof}

\begin{rem}Combining Proposition~\ref{prop:notnot-in-big-zariski-topos} and
Proposition~\ref{prop:char-surjective-morphisms} yields a proof of the fact
that a quasicompact, quasiseparated, and locally finitely presented morphism~$X \to S$, where~$S$ is locally of
finite type over a field, is surjective if it is surjective on closed points.
\end{rem}

\begin{rem}In the case that the parsimonious sites are used, the assumption in
Proposition~\ref{prop:char-surjective-morphisms} that the morphism~$f$ is
quasicompact can't be dropped. For instance, let~$k$ be an
algebraically closed field. Then the canonical morphism
\[ X \defeq \coprod_{a \in k} \Spec(k[X]/(X-a)) \lra \Spec(k[X]) =\vcentcolon S \]
is surjective on closed points. By Proposition~\ref{prop:notnot-in-big-zariski-topos},
it's not the case that~$\ul{X}$ is empty from the internal point of view
of~$\Zar(S)$. However, the morphism is not surjective.
\end{rem}

\subsection{Universally injective morphisms}

\begin{prop}\label{prop:char-univ-injective-morphisms}
Let~$f : X \to S$ be an~$S$-scheme which is locally contained in the site
used to define~$\Zar(S)$. In the case that the parsimonious sites are used to
define~$\Zar(S)$, further assume that~$f$ is quasicompact and quasiseparated.
Then the following statements are equivalent:
\begin{enumerate}
\item The morphism~$f$ is universally injective.
\item The diagonal morphism~$X \to X \times_S X$ is surjective.
\item From the internal point of view of~$\Zar(S)$, any given elements of~$\ul{X}$
are \notnot equal, that is
\[ \Zar(S) \models \forall x\?\ul{X}\_ \forall x'\?\ul{X}\_ \neg\neg(x = x'). \]
\end{enumerate}
\end{prop}

\begin{proof}The equivalence~``(1)~$\Leftrightarrow$~(2)'' is
well-known~\stacksproject{01S4}. The equivalence~``(2)~$\Leftrightarrow$~(3)''
is by Corollary~\ref{cor:char-surjective-morphisms-relative} and the fact that,
internally, there is \notnot a preimage for any element of~$\ul{X} \times
\ul{X}$ under the diagonal map~$\ul{X} \to \ul{X} \times \ul{X}$ if and only if
any given elements of~$X$ are \notnot equal.
\end{proof}

\begin{cor}Let~$p : X \to Y$ be a morphisms of~$S$-schemes which are locally
contained in the site used to define~$\Zar(S)$. In the case that the
parsimonious sites are used to define~$\Zar(S)$, further assume that~$f$ is
quasicompact and quasiseparated. Then the following statements
are equivalent:
\begin{enumerate}
\item The morphism~$p$ is universally injective.
\item From the internal point of view of~$\Zar(S)$, any given elements of any
fiber of~$p$ are \notnot equal.
\end{enumerate}
\end{cor}

\begin{proof}Immediate using Proposition~\ref{prop:char-univ-injective-morphisms}
and the equivalence~$\Zar(Y) \simeq \Zar(S)/\ul{Y}$,
as explained in Section~\ref{sect:change-of-base}.
\end{proof}

\subsection{Universally closed morphisms}

\begin{defn}In the context of a specified local ring, a set~$X$ is
\emph{synthetically closed} if and only if, for any synthetic closed
immersion~$Z \hookrightarrow X$, there is a closed truth value~$\varphi$ such
that~$Z$ is \notnot inhabited if and only if~$\neg\neg\varphi$.
\end{defn}

\begin{ex}\label{ex:singleton-synth-closed}
Any singleton set is synthetically closed.\end{ex}

\begin{ex}\label{ex:r-not-synth-closed}
The specified local ring~$R$ is typically not synthetically closed.
For let~$f \? R$ be an element. Then the inclusion~$Z \defeq \{ g \? R \,|\,
fg - 1 = 0 \} \hookrightarrow R$ is a synthetic closed immersion. The set~$Z$ is \notnot inhabited if and
only if~$f$ is \notnot invertible if and only if~$f$ is invertible (by
Proposition~\ref{prop:a1-field}); typically, there is no closed truth
value~$\varphi$ such that~$\neg\neg\varphi$ is equivalent to the open truth
value~``$f$ is invertible''.\end{ex}

\begin{prop}\label{prop:char-closed-image}
Assume that~$S$ is locally Noetherian. Let~$f : X \to S$ be a
finitely presented morphism. In the situation that one the parsimonious
sites is used to define~$\Zar(S)$, the following statements are equivalent:
\begin{enumerate}
\item The morphism~$f$ has closed image.
\item The morphism~$f$ has universally closed image, that is for
any~$S$-scheme~$T$ the image of the induced morphism~$X \times_S T \to T$ is
closed.
\item $\Zar(S) \models \exists \varphi\?\Omega\_
  \speak{$\varphi$ is a closed truth value} \wedge
  (\neg\neg(\speak{$\ul{X}$ inhabited}) \Leftrightarrow \neg\neg\varphi)$.
\end{enumerate}
\end{prop}

\begin{proof}The direction~``(2)~$\Rightarrow$~(1)'' is trivial, and the
direction~``(1)~$\Rightarrow$~(2)'' is immediate, since the image of~$X
\times_S T \to T$ is the preimage of the image of~$f$ under the structure
morphism~$T \to S$.

For the~``(1)~$\Rightarrow$~(3)'' direction, we may pick the subfunctor
of~$\ul{S}$ induced by the closed immersion~$\im(f) \hookrightarrow S$ as the
sought truth value. This truth value is closed by
Proposition~\ref{prop:char-closed-immersion} and its double negation is
equivalent to~$\speak{$\ul{X}$ is inhabited}$ by
Lemma~\ref{lemma:image-coincides}.

For the converse direction, we see that, after passing to an open covering
of~$S$ which we won't reflect in the notation, there is a closed subfunctor~$Z
\hookrightarrow 1$ such that~$\Zar(S) \models \neg\neg(\speak{$\ul{X}$
inhabited}) \Leftrightarrow \neg\neg(\speak{$Z$ is inhabited})$. By
Proposition~\ref{prop:char-closed-immersion}, this subfunctor is the functor of
points of a closed subscheme of~$S$. Since~$S$ is locally Noetherian, this
subscheme is locally of finite presentation over~$S$. Therefore
Lemma~\ref{lemma:image-coincides} is applicable and yields that the image of~$X
\to S$ coincides with the found closed subscheme. This concludes the proof.
\end{proof}


\begin{cor}\label{cor:univ-closed}
Assume that~$S$ is locally Noetherian. Let~$f : X \to S$ be a
finitely presented morphism. In the situation that one the parsimonious
sites is used to define~$\Zar(S)$, the following statements are equivalent:
\begin{enumerate}
\item The morphism~$f$ is universally closed.
\item $\Zar(S) \models \speak{$\ul{X}$ is synthetically closed}$.
\end{enumerate}
\end{cor}

\begin{proof}Immediate using Proposition~\ref{prop:char-closed-immersion} and
Proposition~\ref{prop:char-closed-image}.
\end{proof}

In view of Corollary~\ref{cor:univ-closed},
Example~\ref{ex:singleton-synth-closed} and Example~\ref{ex:r-not-synth-closed}
have a geometric interpretation. The first example reflects the fact that the
identity morphism~$S \to S$ is universally closed. The second example reflects
the fact that the projection morphism~$\affl \to S$ is typically not
universally closed.

\subsection{Quasicompact and quasiseparated morphisms}
\label{sect:qc-qs-morphisms}

\begin{defn}\label{defn:synthetic-affine-scheme}
In the context of a specified local ring~$R$,
a \emph{synthetic affine scheme} is a set which is isomorphic (as a set)
to the synthetic spectrum of a synthetically
quasicoherent~$R$-algebra.
\end{defn}

\begin{defn}\label{defn:synthetic-scheme}
In the context of a specified local ring~$R$:
\begin{enumerate}
\item A \emph{quasicompact synthetic scheme} is a set~$X$ which admits a finite
open covering~$X = \bigcup_{i=1}^n U_i$ by synthetic affine schemes~$U_i$.
\item A \emph{locally finitely presented quasicompact synthetic scheme} is a
set~$X$ which admits a finite open covering~$X = \bigcup_{i=1}^n U_i$ such that
the sets~$U_i$ are isomorphic to spectra of finitely presented~$R$-algebras.
\item A \emph{finitely presented synthetic scheme} is a set~$X$ which admits a
finite open covering~$X = \bigcup_{i=1}^n U_i$ such that the sets~$U_i$ are
isomorphic to spectra of finitely presented~$R$-algebras and such that the
intersections~$U_i \cap U_j$ can be covered by finitely many open subsets which
are isomorphic to spectra of finitely presented~$R$-algebras.
\end{enumerate}
\end{defn}

\begin{prop}\label{prop:synthetic-schemes}
Let~$X \in \Zar(S)$ be a Zariski sheaf. Consider the following
statements:
\begin{enumerate}
\item $X$ is the functor of points of a quasicompact~$S$-scheme.
\item $X$ is the functor of points of a locally finitely presented quasicompact~$S$-scheme.
\item $X$ is the functor of points of a finitely presented~$S$-scheme (locally
finitely presented, quasicompact, and quasiseparated over~$S$).
\item[(1')] $\Zar(S) \models \speak{$X$ is a quasicompact synthetic scheme}$.
\item[(2')] $\Zar(S) \models \speak{$X$ is a locally finitely presented quasicompact synthetic scheme}$.
\item[(3')] $\Zar(S) \models \speak{$X$ is a finitely presented synthetic scheme}$.
\end{enumerate}
Then:
\begin{itemize}
\item (1)~$\Rightarrow$~(1'), (2)~$\Rightarrow$~(2'), and~(3)~$\Rightarrow$~(3'). 
\item If the parsimonious sites are used to define~$\Zar(S)$, then all three
converses hold.
\item If the given synthetic open immersions are even quasicompact open
immersions, then~(2')~$\Rightarrow$~(2) and~(3')~$\Rightarrow$~(3).
\end{itemize}
\end{prop}

\begin{proof}For proving~(1)~$\Rightarrow$~(1'), (2)~$\Rightarrow$~(2'),
and~(3)~$\Rightarrow$~(3'), we may assume that~$S$ is affine, since the
internal language is local. Let~$X_0$ be an~$S$-scheme representing the functor~$X$. Since
the structure morphism~$X_0 \to S$ is quasicompact and~$S$ is quasicompact,
there exist finitely many open affine subschemes~$U_i \subseteq X_0$ which
cover~$X_0$. By Proposition~\ref{prop:char-open-immersion}, the
subfunctors~$\ul{U_i} \hookrightarrow X$ are synthetic open immersions from the
internal point of view of~$\Zar(S)$. The internal union~$\bigcup_i \ul{U_i}
\hookrightarrow X$ is the functor
\[ T/S \longmapsto \{ f : T \to X_0 \,|\,
  \text{locally, the morphism $f$ factors over one of the opens~$U_i$} \} \]
and therefore coincides with~$X$.

Since each scheme~$U_i$ can be realized as a relative spectrum of a
quasicoherent~$\O_S$-algebra, both~$U_i$ and~$S$ being affine, the sets~$U_i$
are synthetic affine schemes from the internal point of view.
Hence~(1)~$\Rightarrow$~(1'). If~(2) holds, then the~$U_i$ are spectra of
finitely presented~$\O_S$-algebras, so~(2') holds. If~(3) holds, then the
intersections~$U_i \cap U_j$ can be covered by finitely many open subschemes
which are spectra of finitely presented~$\O_S$-algebras, so~(3') holds.

For the converse directions, we first note that the problem is local on~$S$,
since we can glue matching schemes defined over the members of an open covering
of~$S$~\stacksproject{01JJ}. We may therefore assume that we are given
subfunctors~$U_1,\ldots,U_n \hookrightarrow X$ such that, from the internal
point of view of~$\Zar(S)$, the subsets~$U_i \subseteq X$ are synthetic affine
schemes and the inclusions~$U_i \hookrightarrow X$ are synthetic open
immersions. The extra assumptions ensure that
Proposition~\ref{prop:characterization-open-subfunctor} and
Proposition~\ref{prop:char-affine-zar} are applicable. These imply that the
functors~$U_i$ are representable by affine~$S$-schemes. These can be glued to
yield an~$S$-scheme which represents the functor~$X$~\stacksproject{01JJ}.

In case that internally further finiteness conditions are satisfied, the
resulting scheme~$X$ will satisfy the corresponding external finiteness
conditions.
\end{proof}

One can reasonably wonder why we didn't include the following notion
in Definition~\ref{defn:synthetic-scheme}: A \emph{synthetic scheme} is a
set~$X$ which admits an arbitrary open covering by synthetic affine schemes.
The reason is that, with this definition, any subset~$X$ of the singleton
set~$1 = \{\star\}$ is a synthetic scheme, since it admits the open affine
covering~$X = \bigcup \{ 1 \,|\, \star \in X \}$. But not any subfunctor of the
terminal functor is representable by a scheme.

This phenomenon is well-known in synthetic topology; one has to put some
restrictions on the kind of allowed open coverings. Being finite is the
simplest such condition.

\begin{prop}\label{prop:big-char-qcqs}
Let~$\Zar(S)$ be defined using the parsimonious sites. Let~$f : X \to
S$ be a locally finitely presented morphism. Then:
\begin{enumerate}
\item The morphism~$f$ is quasicompact if and only if~$\ul{X}$ is a
quasicompact synthetic scheme from the internal point of view of~$\Zar(S)$.
\item The morphism~$f$ is quasiseparated if and only if, internally, for any
elements~$x,y \? \ul{X}$ the set~$\brak{x = y}$ is a quasicompact synthetic
scheme.
\item The morphism~$f$ is separated if and only if, internally, for any
elements~$x,y \? \ul{X}$ the truth value~$\brak{x = y}$ is closed. If this is
the case, the witnessing ideal can be chosen to be finitely generated.
\end{enumerate}
\end{prop}

\begin{proof}The first statement is by Proposition~\ref{prop:synthetic-schemes}
and by the fact that the scheme representing a representable functor is unique.

The second statement follows from the first by applying it to the diagonal
morphism~$\Delta : X \to X \times_S X$. More precisely, we have the following
chain of equivalences:
\begin{align*}
  &\mathrel{\phantom{\Longleftrightarrow}}
  \Zar(S) \models \forall x,y\?\ul{X}\_
    \speak{$\brak{x = y}$ is a quasicompact synthetic scheme} \\
  &\Longleftrightarrow
  \Zar(X \times_S X) \models
    \speak{$\ul{X}$ is a quasicompact synthetic scheme} \\
  &\Longleftrightarrow
  \text{$\Delta : X \to X \times_S X$ is quasicompact}
\end{align*}
The first step is by the discussion in Section~\ref{sect:change-of-base} and
the second by applying the first statement; the diagonal morphism is locally of
finite presentation~\stacksproject{0818}, as required.

The third statement follows in a similar way, but employing
Proposition~\ref{prop:char-closed-immersion} instead of the first statement.
\end{proof}

\begin{lemma}Internally in~$\Zar(S)$, a set~$X$ is separated if and only if all of
its points are closed. More formally, the following two statements are
equivalent:
\begin{enumerate}
\item For all~$x\?X$, the inclusion~$\{x\} \hookrightarrow X$ is a synthetic
closed immersion.
\item The diagonal morphism~$X \hookrightarrow X \times X$ is a synthetic
closed immersion.
\end{enumerate}
\end{lemma}

\begin{proof}We argue internally. Let~$x\?X$ be fixed. The inclusion~$\{x\}
\hookrightarrow X$ is a synthetic closed immersion if and only if, for
all~$y\?X$, the truth value~``the fiber over~$y$ is inhabited'' is closed. This
truth value is equal to the truth value of~``$x = y$''.

With these observations, the equivalence of the two statements is immediate.
\end{proof}

\subsection{Proper morphisms}
\label{sect:proper-morphisms}

\begin{prop}Let~$\Zar(S)$ be defined using the parsimonious sites. Assume
that~$S$ is locally Noetherian. Let~$X \in \Zar(S)$ be a Zariski sheaf. Then
the following statements are equivalent:
\begin{enumerate}
\item The sheaf~$X$ is the functor of points of a proper~$S$-scheme.
\item $\Zar(S) \models \speak{$X$ is a finitely presented synthetic scheme and
synthetically closed}$.
\end{enumerate}
\end{prop}

\begin{proof}Follows immediately from Proposition~\ref{prop:synthetic-schemes},
Proposition~\ref{prop:big-char-qcqs}, and Corollary~\ref{cor:univ-closed}.
\end{proof}

\section{Case studies}

\subsection{Punctured plane}

\begin{defn}The \emph{synthetic punctured plane} is the set~$P \defeq (\affl)^2
\setminus \{ 0 \}$.\end{defn}

\begin{prop}The synthetic punctured plane, as constructed by the internal
language of~$\Zar(S)$, is the functor of points of the ordinary punctured plane
over~$S$, that is the open subscheme~$D(X) \cup D(Y) \hookrightarrow
\AA^2_S$.\end{prop}

\begin{proof}The set~$P$ can be written as~$\{ (x,y) \? (\affl)^2 \,|\,
\text{$x$ \inv} \vee \text{$y$ \inv} \}$.
\end{proof}

\begin{prop}The evaluation morphism~$\affl[X,Y] \to [P, \affl]$ is bijective.\end{prop}


\begin{proof}The synthetic punctured plane can be expressed as the pushout
\[ P \cong D(X) \amalg_{D(X) \cap D(Y)} D(Y). \]
Therefore we have the chain of isomorphisms
\begin{align*}
  [P,\affl] &\cong
  [D(X) \amalg_{D(X) \cap D(Y)} D(Y), \affl] \\
  &\cong [D(X), \affl] \times_{[D(X) \cap D(Y), \affl]} [D(Y), \affl] \\
  &\cong \affl[X,X^{-1}] \times_{\affl[XY, (XY)^{-1}]} \affl[Y,Y^{-1}] \\
  &\cong \affl[X,Y].
\end{align*}
The penultimate isomorphism exploits the synthetic quasicoherence of~$\affl$,
which ensures that the canonical map
\[ \affl[X,X^{-1}] \longrightarrow [\Spec(\affl[X,X^{-1}]), \affl] \cong
  [D(X), \affl] \]
is bijective. The ultimate isomorphism rests on the purely algebraic argument
that elements of~$\affl[X,X^{-1}]$ and~$\affl[Y,Y^{-1}]$ which agree as
elements of~$\affl[(XY),(XY)^{-1}]$ are both given by an element
of~$\affl[X,Y]$ and in fact by the same element.
\end{proof}

\begin{cor}The punctured plane is not affine.
\end{cor}

\begin{proof}The canonical map~$P \to \Spec([P, \affl])$ is isomorphic to the
strict inclusion~$P \hookrightarrow (\affl)^2$ and therefore not bijective.
\end{proof}

\subsection{Cohomology of projective space}

Let~$\V$ be a finite locally free~$\O_S$-module. In this section, we give
internal descriptions of the projectivization of~$\V$ and of Serre's twisting
bundles on~$\PP(V)$, and show how their cohomology groups can be computed
internally. Let~$V \defeq \V^\Zar$.

\begin{defn}The \emph{synthetic projective space} is the set
\[ \PP(V) \defeq \{ \ell \subseteq V \,|\, \speak{$\ell$ is a one-dimensional
subspace} \}. \]
\end{defn}

The twisting bundles are modules over the affine line of~$\PP(V)$. As such they
are~$\PP(V)$-indexed families of~$\affl$-modules, as in
Section~\ref{sect:relative-tilde-construction}. We set~$X \defeq \PP(V)$.

\begin{defn}\label{defn:twisting-bundle}
For~$m \in \ZZ$, the \emph{twisting bundle}~$\O(m)$ is
the~$\afflx$-module~$(\ell^{\otimes(-m)})_{\ell{:}X}$. In particular, the
\emph{tautological bundle}~$\O(-1)$ is~$(\ell)_{\ell{:}X}$ and its dual
is~$\O(1) = (\ell^\vee)_{\ell{:}X}$.
\end{defn}

We stress that Definition~\ref{defn:twisting-bundle} formalizes the common
geometric intuition about the twisting bundles. In algebraic geometry, one
often pictures a bundle as the collection of its fibers, well knowing that this
picture is formally incomplete and that a bundle is not really determined by
the collection of its fibers. In our account of synthetic algebraic geometry,
the statement that a bundle is given by its collection of fibers is \emph{just
true}, and in fact is so by definition.

\begin{prop}\label{prop:twisting-bundles}
\begin{enumerate}
\item The canonical map~$\affl \to \Gamma(X, \O(0)) = [X, \affl]$ is an isomorphism.
\item The canonical map~$V^\vee \to \Gamma(X, \O(1)),\,\vartheta \mapsto
(\vartheta|_\ell)_{\ell{:}X}$ is an isomorphism.
\item The module~$\Gamma(X, \O(-1))$ is zero.
\end{enumerate}
\end{prop}

\begin{proof}For notational convenience, we assume that~$V$ is isomorphic
to~$(\affl)^2$. In this case~$X = \{ [x \hg y] \,|\, x \neq 0 \vee y \neq 0
\}$. Let~$X = U_0 \cup U_1$ be the standard open covering with
\[ U_0 = \{ [x \hg y] \,|\, x \neq 0 \} = \{ [1 : y] \,|\, y \? \affl \} \cong
\affl \]
and similarly with~$U_1$.

For the first claim, let~$f : X \to \affl$ be an arbitrary function.
Since~$\affl$ is synthetically quasicoherent, the canonical morphism~$\affl[T]
\lra [U_0, \affl]$ is bijective. Hence there is a unique polynomial~$p_0 \?
\affl[T]$ such that~$f([x \hg y]) = p_0(y/x)$ for all~$[x \hg y] \in U_0$.
Analogously, there is a unique polynomial~$p_1 \? \affl[T]$ such that~$f([x \hg
y]) = p_1(x/y)$ for all~$[x \hg y] \in U_1$.

Since~$\affl$ is synthetically quasicoherent, the canonical morphism
\[ \affl[T, T^{-1}] \lra [U_0 \cap U_1, \affl] \]
is bijective. Because~$p_0(y/x) = p_1(x/y)$ as functions on~$U_0 \cap U_1$, the
polynomials~$p_0(T)$ and~$p_1(1/T)$ agree as elements of~$\affl[T, T^{-1}]$.
Hence~$p_0$ and~$p_1$, and thus~$f$, are constant.

This shows that the canonical morphism~$\affl \to [X,\affl]$ is surjective.
This morphism is trivially injective since~$X$ is inhabited.

The other claims are verified similarly.
\end{proof}

In the same spirit as in the proof of Proposition~\ref{prop:twisting-bundles},
we can compute the cohomology of the twisting bundles, if we define it as
the cohomology of the Čech complex associated to the standard open covering.
We want to indicate how the general idea of the usual proof can be carried out
in the synthetic account.

For notational convenience, we assume that~$V$ is isomorphic to~$(\affl)^2$.
The Čech complex for computing~$H^\bullet(X, \O(-2))$ is
\[ 0 \lra \Gamma(U_0,\O(-2))\times\Gamma(U_1,\O(-2)) \lra
  \Gamma(U_0 \cap U_1,\O(-2)) \lra 0. \]
The differential maps~$(f,g)$ to~$g|_{U_0 \cap U_1} - f|_{U_0 \cap U_1}$.

The verification that~$H^0(X, \O(-2))$ is zero proceeds similarly to the proof of
Proposition~\ref{prop:twisting-bundles}. To verify that~$H^1(X, \O(-2)) \cong
\affl$, we explicitly describe the isomorphism. Given a function~$h \?
\Gamma(U_0 \cap U_1,\O(-2))$, there is for any number~$y \? (\affl)^\times$ a
unique value~$h_0(y)$ such that
\[ h([1 \hg y]) = h_0(y) \cdot \left(\begin{pmatrix}1\\y\end{pmatrix}
\otimes \begin{pmatrix}1\\y\end{pmatrix}\right) \?
\operatorname{span}\left(\begin{pmatrix}1\\y\end{pmatrix}\right)^{\otimes 2}. \]
By the principle of unique choice, the mapping~$y \mapsto h_0(y)$
defines a well-defined function on~$(\affl)^\times$ as the notation suggests.
Because~$\affl$ is synthetically quasicoherent, the canonical map
\[ \affl[T, T^{-1}] \lra [U_0 \cap U_1, \affl] \]
is bijective. Hence there is a unique Laurent polynomial~$p_0 \in
\affl[T,T^{-1}]$ such that~$h_0(y) = p_0(y)$ for all~$y \in
(\affl)^\times$. The sought isomorphism maps~$h$ to the coefficient of~$T^{-1}$
in~$p_0$.


\subsection{Grassmannian}

Let~$\V$ be a finite locally free~$\O_S$-module.
We want to illustrate the synthetic approach by verifying the basic fact that the Grassmannian~$\Gr(\V,r)$ of
rank-$r$ locally free quotients of~$\V$, defined as a certain functor of points,
is representable by a locally finitely presented~$S$-scheme using the internal
language of~$\Zar(S)$.

\begin{defn}The \emph{Grassmannian}~$\Gr(\V,r)$ is the functor which associates
to an~$S$-scheme~$f : T \to S$ the set
\[ \Gr(\V,r)(T) \defeq \{
  \text{$U \subseteq f^*\V$ sub-$\O_T$-module} \,|\,
  \text{$(f^*\V)/U$ is locally free of rank~$r$} \}. \]
\end{defn}

\begin{defn}The \emph{synthetic Grassmannian} of rank-$r$ quotients of a
module~$V$ is the set
\[ \Gr(V,r) \defeq \{ \text{$U \subseteq V$ submodule} \,|\,
  \text{$V/U$ is free of rank~$r$} \}. \]
\end{defn}

We could just as well define the synthetic Grassmannian somewhat more directly
as the set of free rank~$r$-quotients (up to isomorphism). This set is canonically
isomorphic to the Grassmannian as we chose to define it, by mapping a
quotient~$\pi : V \twoheadrightarrow Q$ to the kernel of~$\pi$.

\begin{prop}The synthetic Grassmannian of~$\V$, as constructed by the internal
language of~$\Zar(S)$ where~$\V$ looks like an ordinary free module, coincides
with the functorially defined Grassmannian.\end{prop}

\begin{proof}Immediate from
Definition~\ref{defn:interpretation-internal-constructions} and
Proposition~\ref{prop:locally-free-big-zariski}.
\end{proof}

Having established that the internally constructed synthetic Grassmannian
actually describes the external Grassmannian which we're interested in, we can
switch to a fully internal perspective. We'll reflect this switch notationally
by referring to the~$\affl$-module~$V \defeq \V^\Zar$ instead of~$\V$.

We define for any free submodule~$W \subseteq V$ of rank~$r$ the subset
\[ G_W \defeq \{ U \in \Gr(V,r) \,|\, \text{$W \to V \to V/U$ is bijective} \}. \]
This sets admits a more concrete description, since it is in canonical bijection
to the set
\[ G_W' \defeq \{ \pi : V \to W \,|\, \pi \circ \iota = \id \} \]
of all splittings of the inclusion~$\iota : W \hookrightarrow V$: An element~$U
\in G_W$ corresponds to the splitting~$V \twoheadrightarrow V/U
\xrightarrow{({\cong})^{-1}} W$. Conversely, a splitting~$\pi$ corresponds to~$U
\defeq \ker(\pi) \in G_W$.

\begin{prop}The union of the subsets~$G_W$ is~$\Gr(V,r)$.\end{prop}

\begin{proof}Let~$U \in \Gr(V,r)$. Then there exists a
basis~$([v_1],\ldots,[v_r])$ of~$V/U$. The family~$(v_1,\ldots,v_r)$ is
linearly independent in~$V$, therefore the submodule~$W \defeq
\operatorname{span}(v_1,\ldots,v_r) \subseteq V$ is free of rank~$r$. The
canonical linear map~$W \hookrightarrow V \twoheadrightarrow V/U$ maps
the basis~$(v_i)_i$ to the basis~$([v_i])_i$ and is therefore bijective. Thus~$U
\in G_W$.\end{proof}

\begin{prop}The sets~$G_W$ are (quasicompact-)open subsets of~$\Gr(V,r)$.
\end{prop}

\begin{proof}Let~$U \in \Gr(V,r)$. Then~$U \in G_W$ if and only if
the canonical linear map~$W \hookrightarrow V \twoheadrightarrow V/U$
is bijective. Since~$W$ and~$V/U$ are both free modules of rank~$r$, this map is
given by an~$(r \times r)$-matrix~$M$ over~$\affl$; therefore it's bijective if
and only if the determinant of~$M$ is invertible.

Thus we've found a number which is invertible if and only if~$U \in G_W$. By
Corollary~\ref{cor:sufficient-criterion-open-subfunctor}, the truth value of~``$U
\in G_W$'' is open.\end{proof}

\begin{prop}The sets~$G_W$ are synthetic affine schemes. Moreover, the algebras which the~$G_W$ are
spectra of are finitely presented.\end{prop}

\begin{proof}The set of all linear maps~$V \to W$ is the spectrum of
the~$\affl$-algebra~$A \defeq \Sym(\Hom_\affl(V,W)^\vee)$, since
\begin{align*}
  \Spec(A) &=
  \Hom_{\mathrm{Alg}(\affl)}(\Sym(\Hom_{\mathrm{Mod}(\affl)}(V,W)^\vee), \affl) \\
  &\cong \Hom_{\mathrm{Mod}(\affl)}(\Hom_{\mathrm{Mod}(\affl)}(V,W)^\vee, \affl) \\
  &=\Hom_{\mathrm{Mod}(\affl)}(V,W)^{\vee\vee} \\
  &\cong \Hom_{\mathrm{Mod}(\affl)}(V,W).
\end{align*}
In the last step the assumption that not only~$W$, but also~$V$ is a free module
of finite rank enters. (This is the first time in this development that we need
this assumption.)

The set~$G_W'$ is a closed subset of this spectrum, namely the locus where the
generic linear map~$V \to W$ is a splitting of the inclusion~$\iota : W
\hookrightarrow V$. If we choose bases of~$V$ and~$W$,
whereby~$\Sym(\Hom_\affl(V,W)^\vee)$ is isomorphic
to~$\affl[M_{11},\ldots,M_{rn}]$, we can be more explicit: The set~$G_W'$ is
isomorphic to
\[ \Spec(k[M_{11},\ldots,M_{rn}]/(MN-I)), \]
where~$I$ is the~$(r \times r)$ identity matrix, $M$ is the generic matrix~$M =
(M_{ij})_{ij}$, and~$N$ is the matrix of~$\iota$ with respect to the chosen
bases. The notation~``$(MN-I)$'' denotes the ideal generated by the entries
of~$MN-I$.
\end{proof}

\begin{cor}\label{cor:grassmannian-synthetic-scheme}
The Grassmannian~$\Gr(V,r)$ is a locally finitely presented
quasicompact synthetic scheme.\end{cor}

\begin{proof}We need to verify that~$\Gr(V,r)$ admits a finite covering by
spectra of finitely presented~$\affl$-algebras. We already know that~$\Gr(V,r)$
can be covered by the open subsets~$G_W$ and that these sets are spectra of
finitely presented algebras. Therefore it remains to prove that finitely many of
these subsets suffice to cover~$\Gr(V,r)$.

In fact, if we choose an isomorphism~$V \cong (\affl)^n$, we see
that~$\binom{n}{r}$ of these subsets suffice: namely those where~$W$ is one of
the standard submodules of~$(\affl)^n$ (generated by the standard basis
vectors). For if~$U \in \Gr((\affl)^n,r)$, the surjection~$V \to V/U$ maps
the basis of at least one of these standard submodules to a basis and is
therefore bijective. This is because from any surjective~$(r \times n)$-matrix
over a local ring one can select~$r$ columns which form an linearly independent
family.
\end{proof}

\begin{prop}Let~$U \in \Gr(V,r)$. Then the tangent space at~$U$ is given by
$T_U \Gr(V,r) \cong \Hom_{\affl}(U, V/U)$.\end{prop}

\begin{proof}For notational simplicity, we verify the claim in the case~$r = 1$, in
which case~$\Gr(V,r)$ is the projectivization of~$V$. Let~$\gamma : \Delta \to
\PP(V)$ be a tangent vector with base point~$\ell \defeq \gamma(0) \? \PP(V)$.
By Corollary~\ref{cor:open-subfunctors-d1} and
Corollary~\ref{cor:grassmannian-synthetic-scheme}, there's a
lift~$\overline{\gamma} : \Delta \to V$ such that~$\gamma(\varepsilon) =
\operatorname{span}(\overline{\gamma}(\varepsilon))$ for all~$\varepsilon \?
\Delta$.
We define a linear map~$\alpha : \ell \to V/\ell$ by setting
\[ x \longmapsto \alpha(x) = [x/\overline{\gamma}(0) \cdot
\overline{\gamma}'(0)]. \]
The expression ``$x/\overline{\gamma}(0)$'' should be read as follows: The
vector~$x$, being an element of~$\ell$, is some multiple~$\lambda$
of~$\overline{\gamma}(0)$. The expression ``$x/\overline{\gamma}(0)$'' denotes
this unique number~$\lambda$. It can be checked that the vector~$\alpha(x) \?
V/\ell$ does not depend on the choice of the lifting~$\overline{\gamma}$. The
element~$\alpha$ is therefore a well-defined element
of~$\Hom_\affl(\ell,V/\ell)$.

Conversely, let an element~$\alpha \? \Hom_\affl(\ell,V/\ell)$ be given. We
choose vectors~$x_0 \? V$ and~$z \? V$ such that~$\ell =
\operatorname{span}(x_0)$ and~$\alpha(x_0) = [z]$ and define~$\gamma : \Delta
\to \PP(V)$ by setting
\[ \varepsilon \longmapsto \gamma(\varepsilon) = \operatorname{span}(x_0 + \varepsilon z). \]
The definition of~$\gamma(\varepsilon)$ is invariant under scaling of~$x_0$ and
also under changing~$z$ to some other vector~$z + \lambda x_0$ in its
equivalence class:
\begin{align*}
  \operatorname{span}(x_0 + \varepsilon (z + \lambda x_0))
  &= \operatorname{span}((1 + \varepsilon \lambda) x_0 + \varepsilon z) \\
  &= \operatorname{span}(x_0 + \varepsilon / (1 + \varepsilon \lambda) z) \\
  &= \operatorname{span}(x_0 + \varepsilon z),
\end{align*}
since~$\varepsilon / (1 + \varepsilon \lambda) = \varepsilon (1 - \varepsilon
\lambda) = \varepsilon$. Therefore~$\gamma$ is a well-defined element
of~$T_\ell\PP(V)$ which only depends on~$\alpha$ and not on
the arbitrary choices of~$x_0$ and~$z$.

One can check that the two described constructions are mutually inverse.
\end{proof}

\begin{rem}The arguments given in this section are intended to be
applied internally to the big Zariski topos of a base scheme. However, one can
also apply them internally to well-adapted models for synthetic differential
geometry. In this way, one almost obtains that the differential-geometric
Grassmannian can be represented as a manifold; only a verification of
smoothness is missing.
\end{rem}

\section{Beyond the Zariski topology}
\label{sect:beyond-zariski}

The Zariski topology is of course not the only interesting topology
on~$\Sch/S$. For any finer topology~$\tau$, such as the Nisnevich, étale, or fppf
topology (a valuable hyperlinked chart of the various topologies is located
at~\cite{belmans:topologies}), the big~$\tau$-topos of~$S$, that is the topos
of sheaves on~$\Sch/S$ with respect to~$\tau$, is a subtopos
of the big Zariski topos. Therefore there is a modal operator~$\Box_\tau$
in~$\Zar(S)$ reflecting the topology~$\tau$. Explicitly, for an~$S$-scheme~$T$
and a formula~$\varphi$ over~$T$, the meaning of
\[ T \models \Box_\tau \varphi \]
is that there exists a~$\tau$-covering~$(T_i \to T)_i$ of~$T$ such that~$T_i
\models \varphi$ for all~$i$ (where parameters appearing in~$\varphi$ have to
be pulled back along~$T_i \to T$). Succinctly, the formula ``$\Box_\tau
\varphi$'' means that~$\varphi$ holds~$\tau$-locally. Generalizing
Theorem~\ref{thm:box-translation-semantically} from sheaves on locales to
sheaves on arbitrary Grothendieck sites we also have
\[ \Zar(S) \models \varphi^{\Box_\tau} \qquad\text{iff}\qquad
  \Sh((\Sch/S)_\tau) \models \varphi. \]

\subsection{The étale topology}

A basic illustration of these modal operators is provided by the Kummer sequence, that is the short sequence
\[ 1 \lra \mu_n \lra \GG_\text{m} \stackrel{(\underline{\ })^n}{\lra} \GG_\text{m} \lra 1 \]
of multiplicatively-written commutative group objects in~$\Zar(S)$. With the
internal description of~$\mu_n$ and~$\GG_\text{m}$, there is a purely internal
and straightforward proof that this sequence is exact at the first two terms.
But except for trivial cases, the~$n$-th power map~$\GG_\text{m} \to
\GG_\text{m}$ will fail to be an epimorphism;
internally speaking, the statement
\[ \forall f\?(\affl)^\times\_ \phantom{\Box_\text{ét}(}\exists
g\?(\affl)^\times\_ f = g^n\phantom{)} \]
is false in general. However, if~$n$ is invertible in~$\Gamma(S,\O_S)$, the
internal statement
\[ \forall f\?(\affl)^\times\_ \Box_\text{ét}(\exists g\?(\affl)^\times\_ f = g^n) \]
\emph{is} true. In fact, the more general statement
\begin{multline*}
  \forall p\?\affl[X]\_ \speak{$p$ is monic, of positive degree, and separable}
  \Longrightarrow \\
  \Box_\text{ét}(\exists x\?\affl\_ p(x) = 0 \wedge \speak{$p'(x)$ \inv})
\end{multline*}
is true from the internal point of view, where a polynomial~$p$ is called
\emph{separable} if and only if there exists a Bézout representation~$ap + bp'
= 1$. After simplifying, the interpretation of that statement with the
Kripke--Joyal semantics is that for any~$S$-scheme~$T$ and any monic separable
polynomial~$p \in \Gamma(T,\O_T)[X]$ of positive degree there exists an étale
covering~$(T_i \to T)_i$ of~$T$ such that the pullbacks of~$p$ to each of
the~$T_i$ possess a simple zero. The required covering is given
by the canonical surjective étale map~$\RelSpec_T{\O_T[X]/(p)} \to T$.

The following theorem shows that the modal operator~$\Box_\text{ét}$
corresponding to the étale topology admits a purely internal characterization
in~$\Zar(S)$, which furthermore resonates well with the intuition about the
étale topology.

\begin{thm}Let~$S$ be a scheme. Employ one of the parsimonious sites to
define~$\Zar(S)$. The modal operator~$\Box_\text{ét}$
in~$\Zar(S)$ corresponding to the étale topology is the smallest
operator~$\Box$ such that the~$\Box$-translation of the statement~``$\affl$ is
separably closed'' holds.\end{thm}

Here, a ring~$A$ is \emph{separably closed} if and only if
\begin{multline*}
  \forall p\?A[X]\_ \speak{$p$ is monic, of positive degree, and unramifiable}
  \Longrightarrow \\
  \exists x\?A\_ p(x) = 0 \wedge \speak{$p'(x)$ \inv}.
\end{multline*}
We call a polynomial~$p$ over a ring~$A$ \emph{unramifiable} if and
only if it admits at least one simple root in every algebraically closed field
over~$A$. Since quantifying over algebraically closed fields raises red flags
from an intuitionistic point of view, just as quantifying over maximal ideals
does, this condition has to be formulated in a sensible way. One possibility is
to use the \emph{hyperdiscriminants} of~$p$, \ie the elementary symmetric
polynomials in the values of~$p'$ at the roots of~$p$, resulting in a simple
existential statement involving only the coefficients of~$p$; in particular,
the condition for a polynomial to be unramifiable is a geometric formula.
See~\cite[p.~751]{wraith:generic-galois-theory} for the precise formulation.

In more detail, the claim is that firstly~$\Box_\text{ét}$ is a modal operator
such that the displayed formula holds and that secondly, if~$\Box$ is any modal
operator verifying the formula, internally it holds that~$\Box_\text{ét}\varphi
\Rightarrow \Box\varphi$ for any truth value~$\varphi\?\Omega$.

\begin{proof}For the proof we require some familiarity with the concept of
classifying toposes. We are grateful to Felix Geißler for contributing a key step of
the argument.

To verify the first statement, we observe that the displayed formula is a geometric
implication and that the big étale topos~$\Et(S)$ has \emph{enough points}.
Therefore it suffices to show that for any~$S$-scheme~$T$ and any geometric
point~$\bar t$ of~$T$, the stalk~$\O_{T,\bar t}$ is separably closed. It is
well-known that this is true.

For the second statement we may assume without loss of generality that~$S =
\Spec A$ is affine. It is well-known that, for any cocomplete topos~$\E$,
geometric morphisms~$\E \to \Zar(\Spec A)$ are in canonical one-to-one correspondence
with local algebras over~$\ul{A}$ in~$\E$ (where~$\ul{A}$ denotes the pullback
of~$A$ along the unique geometric morphism~$\E \to \Set$) and that geometric
morphisms~$\E \to \Et(\Spec A)$ are in canonical one-to-one correspondence
with algebras over~$\ul{A}$ which are local and separably closed from the
internal point of view of~$\E$;
see~\cite[Section~VIII.6]{moerdijk-maclane:sheaves-logic}
and~\cite{anel:factorization-systems}.

Therefore a geometric morphism~$\E \to \Zar(\Spec A)$ factors over the
geometric embedding~$\Et(\Spec A) \hookrightarrow \Zar(\Spec A)$ if and only if
the pullback of~$\affla$ along~$\E \to \Zar(\Spec A)$ is separably closed.

Let~$\Box$ be a modal operator in~$\Zar(\Spec A)$ such that the~$\Box$-translation
of~``$\affla$ is separably closed'' holds. Then the pullback of~$\affla$
along~$\Zar(\Spec A)_\Box \hookrightarrow \Zar(\Spec A)$ is separably closed
and therefore this geometric embedding factors over~$\Et(\Spec A)
\hookrightarrow \Zar(\Spec A)$. This shows that any~$\Box$-sheaf is also
a~$\Box_\text{ét}$-sheaf.

The claim that~$\Box_\text{ét}\varphi \Rightarrow \Box\varphi$ for any truth
value~$\varphi\?\Omega$ then follows by combining the following two basic
observations of the theory of modal operators, valid for any modal
operator~$\Box$:
\begin{enumerate}
\item $\Box\varphi \Longleftrightarrow
  \forall \psi\?\Omega\_ ((\Box\psi\Rightarrow\psi) \wedge
  (\varphi\Rightarrow\psi)) \Rightarrow \psi$.
\item $(\Box\psi \Rightarrow \psi) \Longleftrightarrow
  \speak{$\{x \in 1 \,|\, \psi\}$ is a~$\Box$-sheaf}$. \qedhere
\end{enumerate}
\end{proof}

\subsection{The fppf topology}

The big Zariski topos of a scheme~$S$ is the classifying topos of local rings
over~$S$, when employing one of the parsimonious sites. The big fppf topos
of~$S$ is a particular subtopos of the big Zariski topos; it therefore
classifies a particular quotient theory of the theory of local rings of~$S$,
obtained by adding certain further axioms~\cite{caramello:lattices}. What are
these axioms?

The analogous question for the big étale topos has been answered by Wraith,
building upon work by Hakim~\cite{hakim:relative-schemes}: The big étale topos
classifies separably closed local rings~\cite{wraith:generic-galois-theory}.
Since the big fppf topos is a subtopos of the big étale topos (any étale
covering being in particular an fppf covering), the sought axioms need to at
least imply the axioms for separably closed rings.

Wraith conjectured that the big fppf topos classifies algebraically closed
local rings (local rings for which any monic polynomial of positive degree has
a zero, also called absolutely integrally closed local rings). We neither confirm nor
refute his conjecture, but we are able to give an alternative explicit
description: The big fppf topos of a scheme classifies fppf-local rings
over~$S$.

These kinds of rings where studied by Schröer and independently by Gabber and
Kelly~\cite{schroer:points-fppf,gabber:kelly:points}; we'll review the notion
below. Every fppf-local ring is algebraically closed, therefore the theory of
fppf-local rings encompasses the theory of algebraically closed local rings. It
is an open question whether these two theories coincide.

As discussed in Section~\ref{sect:proper-choice-of-site}, one can only expect
the big fppf topos to classify a simple theory if one employs parsimonious
sites.

\begin{defn}The \emph{big fppf topos} of a scheme~$S$ is the topos of sheaves
over the category~$(\Sch/S)_\lfp$ of locally finitely presented~$S$-schemes, where a
family~$(f_i : X_i \to X)_i$ of morphisms is deemed a covering if and only if
the morphisms~$f_i$ are flat, locally of finite presentation, and jointly surjective.\end{defn}

The condition that the morphisms~$f_i$ are locally of finite presentation is
automatically satisfied in our setup, since we require that source and target
are locally of finite presentation~\stacksproject{02FV}. One can equivalently
define the big fppf topos as the topos of sheaves over the
category~$(\Aff/S)_\lfp$ of locally finitely presented~$S$-schemes which are
affine as absolute schemes.

\begin{defn}Let~$A$ be a ring.
\begin{enumerate}
\item An \emph{fppf-algebra over~$A$} is an~$A$-algebra~$B$ such that the
structure morphism~$A \to B$ is faithfully flat and of finite presentation.
\item A \emph{basic fppf-algebra over~$A$} is an~$A$-algebra which
is finite free of positive rank as an~$A$-module.
\end{enumerate}
\end{defn}

Since algebras which are free as modules are also finitely presented as
algebras, a basic fppf-algebra is also an fppf-algebra and in fact an integral
fppf-algebra. Conversely, an integral fppf-algebra over a local ring is a basic
fppf-algebra, since integral algebras which are finitely presented as algebras
are also finitely presented as modules~\stacksproject{0564}, finitely presented
flat modules are projective~\stacksproject{058R}, and finitely generated
projective modules over local rings are finite free~\stacksproject{00NX}. This
equivalence even holds intuitionistically.

\begin{defn}\label{defn:fppf-local-ring}
An \emph{fppf-local ring} is a local ring~$A$ such that any finite system of
polynomial equations over~$A$ which has a solution in some basic fppf-algebra
over~$A$ admits a solution in~$A$.
\end{defn}

An fppf-local ring is algebraically closed, since
the~$A$-algebra~$A[X]/(f)$ is a basic fppf-algebra whenever~$f$ is a monic polynomial
of positive degree.

We refer to basic fppf-algebras instead of arbitrary integral fppf-algebras in
Definition~\ref{defn:fppf-local-ring} in order to ensure that the condition for
a ring to be fppf-local is a geometric implication; we'll expand on this below.
The standard definition of fppf-locality refers to arbitrary (not necessarily
integral) fppf-algebras~\cite[Definition~4.1]{schroer:points-fppf}. The
following proposition establishes the equivalence of our definition with the
standard one.

\begin{prop}\label{prop:char-fppf-local-ring}
Let~$A$ be a local ring. The following statements are equivalent.
\begin{enumerate}
\item[(1)] The ring~$A$ is fppf-local.
\item[(1')] Any finite system of polynomial equations over~$A$ which has a
solution in some fppf-algebra over~$A$ admits a solution in~$A$.
\item[(2)] The structure morphism~$A \to B$ of any basic fppf-algebra
has a retraction.
\item[(2')] The structure morphism~$A \to B$ of any
fppf-algebra has a retraction.
\item[(3)] The functor
\[ \Sch \lra \Set,\,X \longmapsto \Hom(\Spec(A), X) \]
maps fppf coverings to jointly surjective families. That is, the canonical
map~$\coprod_i \Hom(\Spec(A), X_i) \to \Hom(\Spec(A), X)$ is surjective for any
fppf covering~$(X_i \to X)_i$.
\end{enumerate}
Furthermore, for any scheme~$S$ and any morphism~$\Spec(A) \to S$, the
following statement is equivalent to the others:
\begin{enumerate}
\item[(4)] The functor
$\Sch/S \to \Set, X \mapsto \Hom_S(\Spec(A), X)$
maps fppf coverings to jointly surjective families.
\end{enumerate}
\end{prop}

\begin{proof}The directions~(1')~$\Rightarrow$~(1),~(2')~$\Rightarrow$~(2),
and~(3)~$\Leftrightarrow$~(4) are trivial.

For verifying~(1)~$\Rightarrow$~(2), let a basic fppf-algebra~$B$ over~$A$ be
given. Writing~$B \cong A[X_1,\ldots,X_n]/(f_1,\ldots,f_m)$, we see that the
polynomial system~``$f_1 = 0, \ldots, f_m = 0$'' has the tautologous
solution~$([X_1],\ldots,[X_n])$ in~$B$. Since~$A$ is fppf-local, it therefore
has a solution in~$A$. Such a solution gives rise to an~$A$-algebra
homomorphism~$A[X_1,\ldots,X_n]/(f_1,\ldots,f_m) \to A$, so to a retraction of
the structure morphism~$A \to B$.

The proof of the converse direction is similar: A solution of a polynomial
system of equations in a basic fppf-algebra~$B$ can be transported along a
retraction to yield a solution in~$A$.

The directions~(1')~$\Rightarrow$~(2') and~(2')~$\Rightarrow$~(1') can be
verified in exactly the same way.

We now verify~(2)~$\Rightarrow$~(3). Let~$(f_i : X_i \to X)_i$ be an fppf
covering and let a morphism~$g : \Spec(A) \to X$ be given. We want to show
that~$g$ factors over one of the morphisms~$f_i$. Since the fppf topology is
generated by Zariski coverings and singleton coverings~$(T \to W)$ where~$T \to
W$ is a surjective finite locally free morphism~\stacksproject{05WN}, we may
assume that the given covering~$(X_i \to X)_i$ consists entirely of open
immersions or is such a singleton covering.

In the first case, the morphism~$g$ lifts to one of the open subschemes~$X_i$
since the preimages~$g^{-1}X_i$ cover~$\Spec(A)$ and~$\Spec(A)$ is a local
topological space.

In the second case, we may assume that~$W$ is affine and that~$f$ is not only
finite locally free, but finite free. The left morphism in the pullback diagram
\[ \xymatrix{
  \Spec(A) \times_W T \ar[r]\ar[d] & T \ar[d] \\
  \Spec(A) \ar[r] & W
} \]
admits a section since it is the induced morphism on spectra of a basic
fppf-algebra. The sought lift is then the composite~$\Spec(A) \to \Spec(A)
\times_W T \to T$.

Finally, we verify~(3)~$\Rightarrow$~(2'). Let an fppf-algebra~$B$ over~$A$ be
given. The singleton family~$(\Spec(B) \to \Spec(A))$ is an fppf covering.
Therefore the identity morphism on~$\Spec(A)$ lifts to~$\Spec(B)$. This lift
yields the desired retraction.
\end{proof}

\begin{rem}A scheme~$T$ such that the functor~$\Hom(T, \placeholder) :
\Sch \to \Set$ maps Zariski coverings to jointly surjective
families is already the spectrum of a local ring. Thus
Proposition~\ref{prop:char-fppf-local-ring} implies that a scheme~$T$ such
that~$\Hom(T, \placeholder)$ maps fppf coverings to jointly surjective families
is the spectrum of an fppf-local ring.\end{rem}

\begin{rem}The proof of Proposition~\ref{prop:char-fppf-local-ring} is
intuitionistically valid, with the possible exception of the
part~(2)~$\Rightarrow$~(3). This part of the proof crucially rests upon
the description of the fppf topology given in the Stacks
Project~\stacksproject{05WN}, which is proved in the usual language
involving prime ideals and is therefore not obviously intuitionistically valid.
\end{rem}

\begin{lemma}\label{lemma:basic-fppf-algebra-geometric}
\begin{enumerate}
\item The condition for a finitely presented algebra to be a basic fppf-algebra can
be expressed as a geometric formula.
\item The condition for a ring to be fppf-local can be expressed as a
countable conjunction of geometric implications.
\end{enumerate}
\end{lemma}

\begin{proof}A finitely presented~$A$-algebra~$B \defeq
A[X_1,\ldots,X_n]/(f_1,\ldots,f_m)$ is a basic fppf-algebra if and only if
there exists a number~$r \geq 1$, polynomials~$g_1,\ldots,g_r \?
A[X_1,\ldots,X_n]$, vectors~$v_1,\ldots,v_n \? A^r$, a vector~$u \? A^r$, and
vectors~$w_{ij} \? A^r$ for~$i,j = 1,\ldots,r$ such that
\begin{itemize}
\item the multiplication defined on~$A^r$ by~$(e_i, e_j) \mapsto w_{ij}$ is
associative, commutative, and has~$u$ as neutral element,
\item the map~$B \to A^r$ given by~$[X_k] \mapsto v_k$ is well-defined, that
is~$f_l(v_1,\ldots,v_n) = 0$ for~$l = 1,\ldots,m$, and
\item the map~$B \to A^r$ and the map~$A^r \to B$ sending~$e_i$ to~$[g_i]$ are
inverse to each other.
\end{itemize}
Each of these conditions can be expressed by geometric formulas involving the
components of the data~$g_i, v_k, u, w_{ij}$.

A ring~$A$ is fppf-local if and only if it is local, and for any numbers~$n \geq 0,
m \geq 0$ and polynomials~$f_1,\ldots,f_m \in A[X_1,\ldots,X_n]$ the
implication
\begin{multline*}
  \quad \speak{$A[X_1,\ldots,X_n]/(f_1,\ldots,f_m)$ is a basic fppf-algebra} \Longrightarrow \\
    \exists x_1,\ldots,x_n \? A\_ \bigwedge_{l=1}^m f_l(x_1,\ldots,x_n) = 0 \quad
\end{multline*}
holds. Since the antecedent can be expressed as a geometric formula, this
formula is a geometric implication.
\end{proof}

\begin{thm}\label{thm:big-fppf-largest-subtopos}
The big fppf topos of a scheme~$S$ is the largest subtopos
of~$\Zar(S)$ where~$\affl$ is fppf-local.
\end{thm}

\begin{proof}Let~$\Box_\fppf$ be the modal operator associated to the fppf
topology; the big fppf topos of~$S$ is the subtopos~$\Zar(S)_{\Box_\fppf}
\hookrightarrow \Zar(S)$. We verify that
\begin{multline*}
  \quad\Zar(S) \models \bigwedge_{n \geq 0} \bigwedge_{m \geq 0}
  \forall f_1,\ldots,f_m \? \affl[X_1,\ldots,X_n]\_ \\
    \speak{$\affl[X_1,\ldots,X_n]/(f_1,\ldots,f_m)$ is a basic fppf-algebra}
    \Longrightarrow \\
      \Box_\fppf(\exists x_1,\ldots,x_n \? \affl\_
      \bigwedge_{j=1}^m f_j(x_1,\ldots,x_n) = 0)\quad
\end{multline*}
and we show that if~$\Box$ is any modal operator with this property, then
\[ \Zar(S) \models
  \forall \varphi \? \Omega\_
    \Box_\fppf\varphi \Rightarrow \Box\varphi. \]
For the first part we may
assume, by Lemma~\ref{lemma:basic-fppf-algebra-geometric}, that~$S = \Spec(A)$ is affine and that we're given
polynomials~$f_1,\ldots,f_m \in A[X_1,\ldots,X_n]$ such that~$B \defeq
A[X_1,\ldots,X_n]/(f_1,\ldots,f_m)$ is a basic fppf-algebra. Then, trivially,
\[ \Spec(B) \models \exists x_1,\ldots,x_n\?\affl\_ \bigwedge_{j=1}^m f_j(x_1,\ldots,x_n) = 0. \]
Since~$(\Spec(B) \to \Spec(A))$ is an fppf covering, we have
\[ \Zar(S) \models \Box_\fppf(\exists x_1,\ldots,x_n \? \affl\_
  \bigwedge_{j=1}^m f_j(x_1,\ldots,x_n) = 0) \]
as claimed.

For the second part, let an fppf covering~$(X_i \to X)_i$ be given such
that~$X_i \models \varphi$ for all~$i$. We want to show that~$X \models
\Box\varphi$. Since the fppf topology is generated by Zariski coverings and
singleton coverings~$(T \to W)$ where~$T \to W$ is a surjective finite locally
free morphism~\stacksproject{05WN} and the internal language of~$\Zar(S)$
is~(Zariski-)local, we may assume that the given covering is such a singleton
covering. Moreover, we may assume that~$W$ is affine and that~$T$ is of the
form~$\Gamma(W,\O_W)[X_1,\ldots,X_n]/(f_1,\ldots,f_m)$. Then
\begin{multline*}
  \quad W \models \speak{$\affl[X_1,\ldots,X_n]/(f_1,\ldots,f_m)$ is a basic
  fppf-algebra}\ \text{and} \\
  W \models \forall x_1,\ldots,x_n \? \affl\_
    \Bigl(\bigwedge_{j=1}^m f_j(x_1,\ldots,x_n) = 0\Bigr) \Rightarrow
      \varphi.\quad
\end{multline*}
The latter is because for any~$\Gamma(W,\O_W)$-algebra~$R$ such that there
are elements $x_1,\ldots,x_n \in R$ with~$f_j(x_1,\ldots,x_n) = 0$ for~$j =
1,\ldots,m$, the structure morphism~$\Spec(R) \to W$ factors over~$T \to W$.
The assumption on the modal operator~$\Box$ implies~$W \models \Box\varphi$.
\end{proof}

The next statement uses the concept of a local ring over~$S$. This was defined
on page~\pageref{defn:ring-over-s}.

\begin{cor}\label{cor:big-fpps-topos-classifies}
The big fppf topos of a scheme~$S$ is the classifying topos of the
theory of fppf-local rings over~$S$.\end{cor}

\begin{proof}The ring object~$\affl$ of the big fppf topos of~$S$ is an
fppf-local ring by Theorem~\ref{thm:big-fppf-largest-subtopos}. Equipped with the morphism~$\flat\affl
\to \affl$, it is thus an fppf-local ring over~$S$.

Let~$\E$ be an arbitrary cocomplete topos containing an fppf-local
ring~$A$ over~$S$. As detailed in
Remark~\ref{rem:zar-classifies-absolute}, this comprises an fppf-local ring~$A$, a model of the theory which~$\Sh(S)$
classifies (yielding a geometric morphism~$f : \E \to \Sh(S)$), and a local
homomorphism~$\alpha : f^{-1}\O_S \to A$. By Theorem~\ref{thm:zar-classifies},
this homomorphism gives rise to a unique geometric morphism~$g : \E \to \Zar(S)$
over~$\Sh(S)$ such that~$g^{-1}\affl \cong A$ and such that the induced
morphism~$g^{-1}(\flat\affl) \to g^{-1}\affl$ coincides with~$\alpha$.

The geometric morphism~$g$ factors over the inclusion of the big fppf topos
by Theorem~\ref{thm:big-fppf-largest-subtopos}, yielding a geometric morphism
from~$\E$ to the big fppf topos such that the pullback of~$\affl$ coincides
with~$A$ as rings over~$S$.

Uniqueness of the geometric morphism constructed in this way follows already
from the universal property of the big Zariski topos.
\end{proof}

\begin{cor}\label{cor:points-fppf-topos}
The points of the big fppf topos of a scheme~$S$ are in canonical
one-to-one correspondence with the fppf-local rings over~$S$, that is
fppf-local rings~$A$ equipped with a morphism~$\Spec(A) \to S$.\end{cor}

\begin{proof}By Proposition~\ref{prop:points-of-big-zariski}, the points
of~$\Zar(S)$ are in canonical one-to-one correspondence with the local rings
over~$S$. Such a point is contained in the big fppf topos (that is, the
associated geometric morphism~$\Set \to \Zar(S)$ factors over the inclusion of
the big fppf topos) if and only if the underlying local ring is fppf-local.
\end{proof}

We wish to record some algebraic and logical facts about fppf-local rings,
hoping that they entice the reader to tackle the question whether any
algebraically closed ring is already fppf-local.

\begin{itemize}
\item For any prime ideal~$\ppp$ of an algebraically closed ring, it holds
that~$\ppp^2 = \ppp$, since any element of~$\ppp$ possesses a square root.
\item If an algebraically closed ring~$A$ is Noetherian, it is already a field,
since by Krull's intersection theorem~$(0) = \bigcap_{n \geq 0} \mmm_A^n =
\mmm_A$.
\item Algebraically closed fields~$K$ are fppf-local: Let~$K \to B$ be a basic
fppf-algebra. Since~$B \neq 0$, there exists a maximal ideal~$\nnn \subseteq
B$. The quotient ring~$B/\nnn$ is an algebraic extension of~$K$. Since~$K$ is
algebraically closed, the identity morphism~$K \to K$ can be extended to a
morphism~$B/\nnn \to K$ of~$K$-algebras, yielding the desired retraction.
This statement can also be verified intuitionistically, see
Proposition~\ref{prop:alg-closed-field-is-fppf-local} below.
\item Since the condition that a ring is fppf-local is a conjunction of
geometric implications, a proof that any algebraically closed ring is
fppf-local using classical logic and the axiom of choice implies that
(nonconstructively) there also exists an intuitionistic proof of this statement.
\item The proof of Lemma~\ref{lemma:basic-fppf-algebra-geometric} yields a way
to phrase fppf-locality of a local ring in the language of linear algebra:
A local ring~$A$ is fppf-local if and only if for any number~$n \geq 1$
and any pairwise commuting~$(n \times n)$-matrices~$W_1,\ldots,W_n$ over~$A$ such that
the first matrix is the identity matrix and such that the~$i$-th row of~$W_j$
coincides with the~$j$-th row of~$W_i$, there is a common eigenvector of
the matrices~$W_1,\ldots,W_n$ whose first component is the unit
of~$A$.\label{page:fppf-linear-algebra}

An equivalent condition is the following: A local ring~$A$ is fppf-local if and
only if for any matrices~$W_1,\ldots,W_n$ as above there are
elements~$\lambda_1,\ldots,\lambda_n \in A$ such that~$\lambda_1 = 1$ and such
that~$W_i \cdot (\lambda_1,\ldots,\lambda_n)^T = \lambda_i \cdot
(\lambda_1,\ldots,\lambda_n)^T$.

The condition is nontrivial only for~$n \geq 3$.

\item In order to show that any algebraically closed local ring is fppf-local,
it would suffice to show that commuting matrices over algebraically closed local
rings admit a common nontrivial eigenvector. However, we suspect that this
stronger statement is false, since it would imply that any matrix over an
algebraically closed local ring admits a nontrivial eigenvector.

Owing to the logical form of this statement, if it admits a classical proof, then there is also a
constructive proof; however, there is probably no such constructive proof,
since it would imply that, internally to any topos, matrices over
algebraically closed local rings admit a
nontrivial eigenvector. Since the complex numbers (constructed using Cauchy
reals) form an algebraically closed local
ring~\cite[Theorem~3.13]{ruitenburg:roots}, this would imply that complex
matrices admit nontrivial eigenvectors. Since eigenvectors are in general
uncomputable~\cite[Proposition~12]{ziegler:brattka:spectrum}, there is probably
a suitable realizability topos in which this statement fails.
\end{itemize}

\subsection{The fpqc topology}

\begin{defn}The \emph{big fpqc topos} of a scheme~$S$ is the topos of sheaves
over the category~$(\Sch/S)_\lfp$ of locally finitely presented~$S$-schemes, where a
family~$(f_i : X_i \to X)_i$ of morphisms is deemed a covering if and only if
the morphisms~$f_i$ are flat, jointly surjective, and each affine open subset
of~$X$ is the union of the images of finitely many affine open subsets of some of
the~$X_i$.
\end{defn}

Every fppf covering is an fpqc covering. Conversely, since we employ the
parsimonious sites, every fpqc covering is an fppf covering. Therefore in our
setup there is no difference between the big fppf topos and the big fpqc topos.

\subsection{The ph topology}

The goal of this section is to give an explicit description of the theory which
the big ph~topos of a scheme classifies, conditional on a conjecture we
didn't prove.

\begin{defn}\begin{enumerate}
\item A \emph{standard ph~covering} of an~$S$-scheme~$X$ is a family of
the form~$(U_i \hookrightarrow T \to X)_i$, where~$T \to X$ is a proper
surjective morphism and~$T = \bigcup_i U_i$ is an open affine covering.
\item A family~$(f_i : X_i \to X)_i$ of morphisms between locally finitely
presented~$S$-schemes if a \emph{covering for the ph~topology} if and
only if for any affine open subset~$U \subseteq X$ the family~$(X_i \times_X U
\to U)_i$ can be refined by a standard ph~covering.
\item The \emph{big ph~topos}~$\Zar(S)_\ph$ of a scheme~$S$ is the topos of
sheaves over the category~$(\Sch/S)_\lfp$ of locally finitely
presented~$S$-schemes equipped with the ph~coverings.
\end{enumerate}
\end{defn}

The ph~topology is not subcanonical, therefore we mean by ``the affine line in
the big ph~topos'' the sheafification of~$\affl$. We write~``$\Box_\ph$'' for
the modal operator of~$\Zar(S)$ corresponding to the ph~topology.

\begin{defn}A \emph{valuation ring} is an integral domain~$R$ in the weak sense
(see Definition~\ref{defn:integral-domain}) such that, for any
elements~$a,b\?R$, $a \mid b$ or $b \mid a$.\end{defn}

With this definition, valuation rings are local rings.

\begin{prop}From the point of view of~$\Zar(S)_\ph$, the affine line is an
algebraically closed valuation ring and a field in the sense that nonzero
elements are invertible.
\end{prop}

\begin{proof}Since~$\Zar(S)_\ph$ is a subtopos of~$\Zar(S)_\fppf$, the affine
line is algebraically closed (and even fppf-local) from the internal point of
view of~$\Zar(S)$. The proof that nonzero elements are invertible proceeds just
as in Proposition~\ref{prop:affl-field-informal}, exploiting that an~$S$-scheme
for which the empty family is a covering is empty.

To show that~$\affl$ is an integral domain in the weak sense, we employ the
ph~covering~$(V(a) \to X, V(b) \to X)$ for~$S$-schemes~$X$ and functions~$a,b
\in \Gamma(X,\O_X)$ such that~$ab = 0$.

To show that for any given functions~$a,b \in \Gamma(X,\O_X)$ one divides the
other, we employ the ph~covering~$(\RelProj_X \O_X[U,V]/(bU-aV) \to X)$. This
covering ensures that, from the internal point of view, the set~$\{
[u\mathbin{:}v] \,|\, bu-av=0 \}$ is inhabited. If~$u$ is invertible, then~$a
\mid b$; if~$v$ is invertible, then~$b \mid a$.
\end{proof}

\begin{prop}\label{prop:affine-line-ph}
The affine line of~$\Zar(S)$ has the following closure property: Any finite
system of homogeneous polynomial equations which the projective Nullstellensatz
predicts to have a nontrivial solution~$\Box_\ph$-has a solution. Formally,
\begin{multline*}
  \Zar(S) \models \bigwedge_{n,m\geq0}
    \forall f_1,\ldots,f_m\?\affl[X_0,\ldots,X_n]\_ \\
      \Bigl(\speak{the $f_j$ are homogeneous} \wedge
      \neg\bigl((X_0,\ldots,X_n) \subseteq \sqrt{(f_1,\ldots,f_m)}\bigr)\Bigr)
      \Longrightarrow \\
        \Box_\ph\Bigl(\exists x_0,\ldots,x_n\?\affl\_
          \bigvee_{i=0}^n \speak{$x_i$ \inv} \wedge
          \bigwedge_{j=1}^m f_j(x_0,\ldots,x_n) = 0\Bigr).
\end{multline*}
\end{prop}

\begin{proof}Let homogeneous polynomials~$f_1,\ldots,f_m \in
\Gamma(T,\O_T)[X_0,\ldots,X_n]$ be given such that
\[ T \models \neg\bigl((X_0,\ldots,X_n) \subseteq
\sqrt{(f_1,\ldots,f_m)}\bigr). \]
Then the projection morphism~$T' \defeq \RelProj_T
\O_T[X_0,\ldots,X_n]/(f_1,\ldots,f_m) \to T$ is surjective. Since
\[ T' \models \exists x_0,\ldots,x_n\?\affl\_
  \bigvee_{i=0}^n \speak{$x_i$ \inv} \wedge
    \bigwedge_{j=1}^m f_j(x_0,\ldots,x_n) = 0 \]
and since~$(T' \to T)$ is a ph~covering, it follows that
\[ T \models \Box_\ph\Bigl(\exists x_0,\ldots,x_n\?\affl\_
  \bigvee_{i=0}^n \speak{$x_i$ \inv} \wedge
    \bigwedge_{j=1}^m f_j(x_0,\ldots,x_n) = 0\Bigr). \qedhere \]
\end{proof}

To proceed further, we need to assume the following statement.

\begin{conjecture}\label{conj:nullstellensatz-geometric}
For any natural numbers~$n$ and~$m$ and any finite set of variables organized
to yield the coefficients of~$m$ homogeneous polynomials~$f_1,\ldots,f_m$
in~$n+1$ variables, there is a geometric formula~$\alpha(f_1,\ldots,f_m)$ which
is equivalent to the formula
\[ \neg\Bigl(\bigwedge_{i=0}^n X_i \in \sqrt{(f_1,\ldots,f_m)}\Bigr) \]
in the intuitionistic first-order theory of algebraically closed valuation
rings which are fields in the sense that an element is invertible if and only
if its nonzero.\end{conjecture}

We believe that this conjecture can be verified by using classical resultant
theory. For instance, the conjecture holds in the easy case of two linear
homogeneous polynomials in two variables, where~$\alpha(f,g)$ can be chosen
as~$\operatorname{Res}(f,g) = 0$.

\begin{prop}\label{prop:ph-topos-classifies}
Let~$S$ be a locally Noetherian scheme. Assuming the existence of
geometric formulas~$\alpha(f,g)$ as in
Conjecture~\ref{conj:nullstellensatz-geometric}, the big ph~topos of~$S$
classifies algebraically closed valuation rings with the extra property that
\[ \alpha(f_1,\ldots,f_m) \Longrightarrow \exists x_0,\ldots,x_n\_
    \bigvee_{i=0}^n \speak{$x_i$ \inv} \wedge
      \bigwedge_{j=1}^m f_j(x_0,\ldots,x_n) = 0.
\]
\end{prop}

\begin{proof}We show that~$\Zar(S)_\ph$ is the largest subtopos of~$\Zar(S)$
where~$\affl$ has the properties mentioned in the statement of the proposition.
The claim then follows in the same way as
Corollary~\ref{cor:big-fpps-topos-classifies} follows from
Theorem~\ref{thm:big-fppf-largest-subtopos}.

Applying Proposition~\ref{prop:affine-line-ph}, one can see that the affine
line has these properties internally in~$\Zar(S)_\ph$.

Conversely, let a modal operator~$\Box$ be given such that~$\affl$ has these
properties in~$\Zar(S)_\Box$. We want to verify that
$\Zar(S) \models
  \forall \varphi \? \Omega\_
    \Box_\ph\varphi \Rightarrow \Box\varphi$.
Let a ph~covering~$(X_i \to X)_i$ of an~$S$-scheme~$X$ be given such that~$X_i
\models \varphi$ for all~$i$. We are to show that~$X \models \Box\varphi$.

Without loss of generality, we may assume that~$X$ is affine and that the
covering is a standard ph~covering. Furthermore, we may assume that the
covering is a singleton covering~$(T \to X)$ where~$T \to X$ is a proper
surjective morphism. By Chow's lemma, we can further assume that~$T \to X$ is
the canonical projection of a closed subscheme of some projective
space~$\PP^n_X$ to~$X$. Since~$X$ is locally Noetherian, the defining sheaf of
ideals is of finite type. Thus we may assume~$T = \RelProj_X
\O_X[U_0,\ldots,U_n]/(f_1,\ldots,f_m)$. Since~$T \models \varphi$, we have
\[ X \models \speak{the system ``$f_1 = \cdots = f_m = 0$'' has a nontrivial
solution}
\Rightarrow \varphi. \]
Since~$T \to X$ is surjective, we have
\[ X \models
    \neg\Bigl(\bigwedge_{i=0}^n U_i \in \sqrt{(f_1,\ldots,f_m)}\Bigr) \]
and therefore~$X \models \Box(\alpha(f_1,\ldots,f_m))$ (we might not have~$X
\models \alpha(f_1,\ldots,f_m)$, since the proof of equivalence may assume
that~$\affl$ is an algebraically closed valuation ring, which~$\affl \in
\Zar(S)$ is only in the pathological case~$S = \emptyset$). Combining these, we
see that~$X \models \Box\varphi$.
\end{proof}

A corollary of Proposition~\ref{prop:ph-topos-classifies} (assuming
Conjecture~\ref{conj:nullstellensatz-geometric}) is that the points of the big
ph~topos are those algebraically closed valuation rings over~$S$ which satisfy
the condition on solvability of systems of homogeneous polynomial equations.
Goodwillie and Lichtenbaum have already determined the points of the big
ph~topos to be the valuation rings with algebraically closed field of fractions
over~$S$~\cite[Proposition~2.2]{goodwillie-lichtenbaum:cohomological-bound},
without any extra condition on solvability of systems of equations.

This mismatch can be explained as follows. In classical logic one can show,
using the valuative criterion for properness, that the rings studied by
Goodwillie and Lichtenbaum automatically satisfy the condition on solvability
of systems of equations. However, we don't believe that the proof can be made
intuitionistic. Assuming this, it's no surprise that the theory which the big
ph~topos classifies contains further axioms.

\subsection{The surjective topology}

\begin{defn}A family~$(f_i : X_i \to X)_i$ of morphisms between locally finitely
presented~$S$-schemes is a \emph{covering for the surjective topology} if and
only if the morphisms~$f_i$ are jointly surjective and each affine open subset
of~$X$ is the union of the images of finitely many affine open subsets of some of
the~$X_i$.
\end{defn}

Equivalently, a family~$(X_i \to X)_i$ is a covering for the surjective
topology if and only if any affine open subset of~$X$ is the image of a
quasicompact open subset under the induced morphism~$\coprod_i X_i \to X$.

\begin{defn}An \emph{algebraically closed geometric field} is a ring such
that~$1 \neq 0$, any element is zero or invertible, and that any monic
polynomial of positive degree has a zero.
\end{defn}

In contrast with other field conditions in intuitionistic mathematics, the
condition for a ring to be an algebraically closed geometric field is the
(countable conjunction of) geometric implications.

\begin{prop}\label{prop:alg-closed-field-is-fppf-local}
Intuitionistically, an algebraically closed geometric field is
fppf-local.
\end{prop}

\begin{proof}Over geometric fields, the kernel of any matrix admits a finite
basis. Moreover, the kernel of any matrix of determinant zero contains a vector
which has at least one invertible component. Therefore the usual proof that
commuting matrices admit a common eigenvector (with at least one invertible
component) applies. This fact can be used to show that the linear algebra
problem stated on page~\pageref{page:fppf-linear-algebra} which characterizes
fppf-locality is solvable.
\end{proof}

\begin{thm}\label{thm:surjective-topology-classifies}
The topos of sheaves over~$(\Sch/S)_\lfp$ for the surjective topology is the
largest subtopos of~$\Zar(S)$ where~$\affl$ is an algebraically closed
geometric field.
\end{thm}

\begin{proof}Let~$\Box_\surj$ be the modal operator associated to the surjective
topology. We verify that
\begin{multline*}
  \quad\Zar(S) \models
    \bigl(\forall s\?\affl\_
      \Box_\surj(s = 0 \vee \speak{$s$ \inv})\bigr) \wedge \\
    \bigwedge_{n \geq 0} \forall a_0,\ldots,a_{n-1}\?\affl\_
      \Box_\surj(\exists x\?\affl\_ x^n + a_{n-1}x^{n-1} + \cdots + a_1x + a_0 = 0)
\end{multline*}
and we show that if~$\Box$ is any modal operator with this property, then
\[ \Zar(S) \models
  \forall \varphi \? \Omega\_
    \Box_\surj\varphi \Rightarrow \Box\varphi. \]
For the first part it suffices to prove the following two statements:
If~$s \in A$ is an element of a ring, then there is a covering~$(X_i \to
\Spec(A))_i$ for the surjective topology such that, for each~$i$,~$X_i \models
s = 0$ or~$X_i \models \speak{$s$ \inv}$. If~$p \in A[X]$ is a monic polynomial
of positive degree over a ring, then there is a covering~$(X_i \to
\Spec(A))_i$ for the surjective topology such that, for each~$i$,~$X_i \models
\exists x\?\affl\_ p(x) = 0$.

For the first claim, we may use the covering~$(D(s) \to \Spec(A), V(s) \to
\Spec(A))$. For the second claim, we may use the singleton
covering~$(\Spec(A[X]/(p)) \to \Spec(A))$.

For the second part, let a covering~$(X_i \to X)_i$ of the surjective topology
be given such that~$X_i \models \varphi$ for each~$i$. We want to show that~$X
\models \Box\varphi$. Without loss of generality, we may assume that~$X$ is
affine and that the given covering is a singleton covering~$(Y \to X)$
where~$Y$ is affine and therefore~$Y \to X$ is of finite presentation.

By the lemma on the existence of a flattening
stratification~\stacksproject{0ASY}, there exist finitely many locally closed
subschemes~$E_j = D(f_j) \cap V(g_{j1},\ldots,g_{j,m_j}) \subseteq X$ such that
the pullback~$Y_j \defeq Y \times_X E_j \to E_j$ is flat (and, being surjective, therefore
faithfully flat). Since~$Y_j \models \varphi$ and since~$(Y_j \to E_j)$ is an
fppf covering, we have~$E_j \models \Box_\fppf\varphi$. By
Proposition~\ref{prop:alg-closed-field-is-fppf-local} and
Theorem~\ref{thm:big-fppf-largest-subtopos}, we also have~$E_j \models
\Box\varphi$.

It's easily checked that~$E \models (\speak{$f_j$ \inv} \wedge g_{j1} = 0
\wedge \cdots \wedge g_{j,m_j} = 0) \Rightarrow \varphi$, for each~$j$. To
conclude that~$X \models \Box\varphi$, it therefore suffices to verify that
\[ X \models \Box\Bigl(\bigvee_{j=1}^n (\speak{$f_j$ \inv} \wedge
  g_{j1} = 0 \wedge \cdots \wedge g_{j,m_j} = 0)\Bigr). \]
This claim follows from distributivity of disjunction over conjunction\footnote{Let
statements~$\psi_{jk}$ where~$j = 1,\ldots,n$, $k = 1,\ldots,r_j$ be given.
We picture this situation as a ragged matrix of statements with~$\psi_{jk}$ located
at column~$k$ of row~$j$. Assume that, for any selection of one statement from
each row, at least one of the selected statements holds. By distributivity of
disjunction over conjunction, there is a row all of whose statements hold.}
and the elementary reformulation of the statement that~$X = \bigcup_{j=1}^n
E_j$: For any finite subset~$J \subseteq \{1,\ldots,n\}$ and any indices~$k_j
\in \{ 1, \ldots, m_j \}$ for~$j \in J$,
\[ \prod_{j \in J} g_{j, k_j} \in \sqrt{(f_j)_{j \in \{1,\ldots,n\} \setminus J}}. \qedhere \]
\end{proof}

\begin{cor}The topos of sheaves over~$(\Sch/S)_\lfp$ for the surjective
topology is the classifying topos of algebraically closed geometric fields
over~$S$. The points of that topos are the algebraically closed geometric
fields over~$S$.\end{cor}

\begin{proof}Follows from Theorem~\ref{thm:surjective-topology-classifies} in
the same way as Corollary~\ref{cor:big-fpps-topos-classifies} and
Corollary~\ref{cor:points-fppf-topos} follow from
Theorem~\ref{thm:big-fppf-largest-subtopos}.
\end{proof}

\subsection{The double negation topology}

As in Section~\ref{sect:subspace-to-modal-operator}, let~$\Zar(S)_{\neg\neg}$
be the smallest dense subtopos of~$\Zar(S)$ (defined using the parsimonious
sites). It is the Boolean topos of sheaves over~$(\Sch/S)_\lfp$ for the double
negation topology. The following proposition describes the double
negation topology in explicit terms.

\begin{prop}The subtopos~$\Zar(S)_{\neg\neg}$ can be presented as the topos of
sheaves over the site~$(\Sch/S)_\lfp$ whose covering families~$(X_i \to X)_i$
are precisely the families such that
\[ \label{eqn:almost-surjectivity}\tag{$\star$}
  \textnormal{$T \times_X X_i = \emptyset$ for all $i$} \quad\text{implies}\quad
  T = \emptyset \]
for all locally finitely presented~$X$-schemes~$T$.

In case that~$S$ is locally of finite type over a field,
condition~\eqref{eqn:almost-surjectivity} is satisfied if and only if for every
closed point~$x \in X$ there is a finite field extension~$K \fieldext k(x)$
such that~$x$ has a~$K$-valued preimage in one of the~$X_i$.
\end{prop}

\begin{proof}For the first claim, we note that by a general
fact~$\Zar(S)_{\neg\neg}$ can be presented as the topos of sheaves
over~$(\Sch/S)_\lfp$ whose covering families~$(f_i : X_i \to X)_i$ are
precisely the families such that~$(\ul{X_i}^{++} \to \ul{X}^{++})_i$ is a
jointly epimorphic family in~$\Zar(S)_{\neg\neg}$, where~$(\placeholder)^{++} :
\Zar(S) \to \Zar(S)_{\neg\neg}$ is the sheafification functor.

This is the case if and only if the morphisms~$\ul{X_i}^{++} \to \ul{X}^{++}$
are jointly surjective from the internal point of view of~$\Zar(S)_{\neg\neg}$.
By the generalization of Theorem~\ref{thm:box-translation-semantically} from
locales to toposes, this is the case if and only if
\[ \Zar(S) \models \forall x\?\ul{X}\_
  \neg\neg \bigl(\bigvee_i \exists u\?\ul{X_i}\_
    \ul{f_i}(u) = x\bigr). \]
As detailed in Section~\ref{sect:change-of-base}, this is the case if and only if
\[ \Zar(X) \models \neg\neg\bigl(\bigvee_i \speak{$\ul{X_i}$ is
inhabited}\bigr). \]
Similarly to the (easy part of)
Proposition~\ref{prop:char-surjective-morphisms}, this in turn is equivalent to
condition~\eqref{eqn:almost-surjectivity}.

The second claim follows from
Proposition~\ref{prop:notnot-in-big-zariski-topos}.
\end{proof}

\begin{prop}Let~$p$ be a point of~$\Zar(S)$. By
Proposition~\ref{prop:points-of-big-zariski}, the point~$p$ corresponds to a
local ring~$A$ over~$S$. If~$p$ is even a point of~$\Zar(S)_{\neg\neg}$,
then~$A$ is an algebraically closed geometric field.\end{prop}

\begin{proof}By Proposition~\ref{prop:a1-field} and
Proposition~\ref{prop:affl-anonymously-algebraically-closed},
the~$\neg\neg$-translation of~``$\affl$ is an algebraically closed geometric
field'' holds in~$\Zar(S)$. Therefore~$\affl$ is an algebraically closed
geometric field from the internal point of view of~$\Zar(S)_{\neg\neg}$. Since
being an algebraically closed geometric field is a (conjunction of) geometric
implications, this property is preserved under pullback along geometric
morphisms. This suffices to establish the claim.
\end{proof}

Since~$\affl$ is an algebraically closed geometric field from the internal
point of view of~$\Zar(S)_{\neg\neg}$,
Theorem~\ref{thm:surjective-topology-classifies} implies
that~$\Zar(S)_{\neg\neg}$ is a subtopos of the
subtopos~$\Zar(S)_\surj$.\footnote{This can also be seen more directly by
observing that~$\Zar(S)_\surj$ is a dense subtopos of~$\Zar(S)$.
Because~$\Zar(S)_{\neg\neg}$ is the smallest dense subtopos of~$\Zar(S)$, this
observation implies that~$\Zar(S)_{\neg\neg}$ is a subtopos of~$\Zar(S)_\surj$.
Denseness of~$\Zar(S)_\surj$ can be checked by verifying that
the internal statement~$\Zar(S) \models \neg \Box_\surj\bot$ holds; this
amounts to verifying that a scheme which admits a cover in the surjective
topology by empty schemes is empty.}

These toposes don't coincide, however, for~$\Zar(S)_{\neg\neg}$ is Boolean
while~$\Zar(S)_\surj$ is not. One way to see this is to appeal to the syntactic
characterization of when the classifying topos of a coherent theory is
Boolean~\cite[Theorem~D3.4.6]{johnstone:elephant}. In the coherent theory of an
algebraically closed geometric field (which is the theory which~$\Zar(S)_\surj$
classifies in the special case~$S = \Spec(\ZZ)$), the formulas~``$p \cdot 1 = 0$''
for prime numbers~$p$ are pairwise non-equivalent. By the cited characterization, the
existence of an infinite family of pairwise non-equivalent formulas is
sufficient to ensure that the classifying topos is not Boolean.

This observation settles a question by
Madore~\cite[entry~2002-03-16:036]{madore:diary}.

We don't know much more about~$\Zar(S)_{\neg\neg}$. Far from knowing which
theory~$\Zar(S)_{\neg\neg}$ classifies, we don't even have a full description
of its points. In the case that~$S$ is~$\Spec(\ZZ)$, determining this theory
would solve the purely algebro-logical problem of describing the
Booleanization of the theory of a local ring.

\section{Outlook}
\label{sect:outlook}

We believe that the approach of using the internal language in algebraic geometry
is already in its current form as developed in Part~\ref{part:little-zariski} and
Part~\ref{part:big-zariski} useful to working algebraic geometers
and interesting for topos theorists and logicians. However, there are many
intriguing avenues for further research.

{\tocless

\subsection*{Further axioms for synthetic algebraic geometry} We showed in
Section~\ref{sect:special-properties-affl} and in
Section~\ref{sect:basic-constructions} that from the single axiom
\begin{quote}
``The ring~$\affl$ is local and synthetically quasicoherent.''
\end{quote}
a surprising number of properties can be deduced and that the core of a synthetic
account for algebraic geometry, comprising notions such as affine schemes,
open and closed immersions, surjective, universally closed and proper
morphisms, and quasicompact synthetic schemes, can be built around it.

However, we still don't feel that we have a full understanding of the notion of
quasicoherence.

Firstly, there are several properties of synthetically quasicoherent modules
which we couldn't account for internally. For instance, the tensor product of
synthetically quasicoherent modules is synthetically quasicoherent, but our
only proof of this fact is external (Lemma~\ref{lemma:tensor-product-qcoh}).
It doesn't seem prudent to just add statements like these as axioms on a
case-by-case basis. We rather need to find further suitable general axioms from
which these statements can be deduced.

Secondly, we know of no properties of~$\affl$, formulated purely within the
first-order language of rings, which wouldn't follow from the axiom that~$\affl$
is local and synthetically quasicoherent. Have we just not looked hard enough,
or does this observation have a deeper reason? Determining the first-order
properties of the universal model of a geometric theory is a well-known
problem; maybe having a look at this special case can shed light on this
problem.

While speculating, we can just as well remark that the question can be
generalized considerably:

\begin{speculation}Let~$\TT$ be an equational theory. Let~$\TT'$ be a geometric
theory obtained from~$\TT$ by adding further axioms. Then every first-order
sequent which is true for the universal model~$U$ of~$\TT'$ can be deduced from
the following quasicoherence condition: ``For any finitely
presented~$\TT$-model~$A$, the canonical map
\[ A \lra [ [A,U]_{\Mod(\TT)}, U ] \]
is bijective.''
\end{speculation}

To make this speculation into a rigorous conjecture, we would have to specify
what kind of deductions from the quasicoherence condition (which is a
higher-order statement) are allowed.

%

\subsection*{Characterizing Zariski toposes among arbitrary toposes}
When is an arbitrary topos the big Zariski topos of a scheme? How can this property
be detected in the internal language, or in otherwise intrinsic terms? An
answer to this question would also tell us which axioms should be stipulated
for synthetic algebraic geometry.

Conversely, we proposed on page~\pageref{page:big-zariski-of-lrl} a definition
of what could be called the big Zariski topos of an arbitrary locally ringed
locale~$S$ (or even a locally ringed topos): the classifying~$\Sh(S)$-topos of
the~$\Sh(S)$-theory of local~$\O_S$-algebras which are local over~$\O_S$.
Internal to any such topos, one can try to conduct synthetic algebraic
geometry. Which properties of~$S$, shared by all schemes, are needed for which
results in synthetic algebraic geometry?

\subsection*{Proper morphisms} In Section~\ref{sect:proper-morphisms}, we
characterized proper morphisms in the
internal language of the big Zariski topos by mimicking the classical
definition as a separated universally
closed morphism (fulfilling a finiteness condition). Although we didn't tell so
in Part~\ref{part:little-zariski}, a characterization along similar lines is
possible in the internal language of the little Zariski topos.


However, the synthetic approach should facilitate a characterization which is
closer to the intuitive way of thinking about proper morphisms: that any
one-parameter family has a unique limit. This intuition can be formally
expressed by the valuative criterion for properness, which states that a
morphism~$X \to S$ of schemes is proper if and only if it is of finite type and
if for every valuation ring~$A$ with field of fractions~$K$, any solid diagram
\[ \xymatrix{
  \Spec(K) \ar[r]\ar[d] & X \ar[d] \\
  \Spec(A) \ar[r]\ar@{-->}[ru] & S
} \]
can be filled by a unique dashed morphism. It is evident from this formulation
of the valuative criterion that it depends only on the functor of points
of~$X$; therefore it should be perfectly suited to internalization in the big
Zariski topos of~$S$.

However, as of yet, we failed to do so. One difficulty is that we don't know
a description of the sheaf of rational functions in the internal language of
the big Zariski topos. We showed in Section~\ref{sect:rational-functions}
that~$\K_S$ is, from the point of view of the little Zariski topos, just the
localization of~$\O_S$ at the filter of regular elements. This description
can't be carried over to the big Zariski topos: Since regular elements
of~$\affl$ are already invertible by Proposition~\ref{prop:a1-field}, the
localization of~$\affl$ at the filter of regular elements is just~$\affl$.

Moreover, depending on the site used, the synthetic spectrum of~$\K_S^\Zar$ can
be empty from the internal point of view. A specific example is given by~$S
= \Spec(\ZZ)$ and the parsimonious sites. In this case, the synthetic spectrum
of~$\K_S^\Zar$ coincides with the functor of points of
the~$\ZZ$-scheme~$\Spec(\QQ)$, which is empty from the internal point of view
by Proposition~\ref{prop:fingen-algebra-q}.

A related problem is that we don't know how to properly describe punctured
germs of curves. For instance, a candidate for what could be considered a germ
of a curve is the synthetic spectrum of the formal power series
ring~$\affl\brak{T}$. This spectrum is canonically isomorphic to the subset of~$\affl$
consisting of those elements which are \notnot zero.\footnote{Given an
element~$a \? \affl$ which is \notnot zero, the evaluation map~$\affl\brak{T}
\to \affl,\ f \mapsto f(a)$ is a well-defined element of~$\Spec(\affl\brak{T})$
since~$a$ is nilpotent by Proposition~\ref{prop:a1-nilp}. Conversely,
let~$\varphi : \affl\brak{T} \to \affl$ be a homomorphism
of~$\affl$-algebras. If~$\varphi(T)$ is invertible, then so
is~$\varphi(\varphi(T) - T) = \varphi(T) - \varphi(T) = 0$, since~$\varphi(T) -
T$ is invertible in~$\affl\brak{T}$. Hence~$\varphi(T)$ is not invertible and
hence \notnot zero.} Removing the origin from this subset yields, however, the
empty set, and not some nontrivial punctured germ.

We strongly believe that these hurdles can be overcome and that there exists an
characterization of proper maps in the big Zariski topos which is close to
geometric intuition.

\subsection*{Cohomology, intersection theory, derived categories}
Real algebraic geometry starts where the basics of scheme theory end. Therefore
there should be an internal treatment of advanced tools like cohomology,
intersection theory, and bounded derived categories of coherent sheaves. Such
an account would be particularly interesting in the big Zariski topos, where
the largest simplifications can be expected.

\subsection*{Stacks and other kinds of generalized schemes} These notes only
deal with classical 1-categorical schemes. An internal treatment of higher
stacks and derived schemes is also desirable; it would probably rest upon a
version of homotopy type theory as the internal language of $\infty$-toposes.

If one wants to stay in the 1-categorical setting, then one could extend the
internal language to the various accounts of scheme theory over the field with
one element~\cite{pena-lorscheid:fun}.

\subsection*{Applications to constructive algebra} We demonstrated in
Section~\ref{sect:generic-freeness} using the example of Grothendieck's generic
freeness lemma how the internal language of the little Zariski topos can be
used to derive results in constructive algebra. Particularly useful is the
property
\[ \Sh(S) \models \forall f\?\O_S\_ \neg\neg(f = 0) \Longrightarrow f = 0, \]
valid for reduced schemes, which allows us to use classical logic to some
extent and still obtain intuitionistically valid results. The internal language
is therefore a tool for constructively pretending that a reduced ring is an
anonymously Noetherian
local ring with~$\neg\neg$-stable equality and the field
property~$\neg(\speak{$f$ \inv}) \Rightarrow f = 0$.

Any of the related toposes, such as the toposes of sheaves for related
topologies such as the constructible or the flat topology, or the various
subtoposes of the big Zariski topos, can be used for similar purposes.
Which statements of constructive algebra profit from the internal language of
theses toposes? Are there further useful toposes?

\subsection*{Constructive algebraic geometry} The book on homotopy type theory
states~\cite[Section~3.4]{hott}:
\begin{quote}
Thus, contrary to how it may appear on the surface, doing mathematics
``constructively'' does not usually involve giving up important theorems, but
rather finding the best way to state the definitions so as to make the important
theorems constructively provable. That is, we may freely use the [law of
excluded middle] when first investigating a subject, but once that subject is
better understood, we can hope to refine its definitions and proofs so as to
avoid that axiom.
\end{quote}
We believe that algebraic geometry has definitely reached the necessary
maturity alluded to in this quote and that it will be a rich and interesting
endeavor to try to give a completely constructive account of algebraic
geometry, including nontrivial results of current interest. We sketched why one
might be interested in such an account in
Section~\ref{sect:constructive-scheme-theory}.

The Stacks Project~\cite{stacks-project} already takes care to formulate
statements in proper generality, not needlessly requiring Noetherian hypotheses
or demanding that fields are algebraically closed. Having such a careful
treatment is a very valuable first step towards a fully constructive development.
Additionally, a constructive account of commutative algebra is readily
available~\cite{mines-richman-ruitenburg:constructive-algebra,lombardi:quitte:constructive-algebra};
a constructive account of algebraic geometry is therefore entirely within
reach.

}

\chapter*{Appendix}

\newcounter{saved-section-number}
\setcounter{saved-section-number}{\value{section}}

\begin{appendix}

\setcounter{section}{\value{saved-section-number}}

\section{Dictionary relating external notions and notions internal to the little Zariski topos}

{\small\renewcommand{\arraystretch}{1.3}
\begin{longtable}{@{}p{4.4cm}@{\qquad}p{6.7cm}@{\qquad}p{1.5cm}@{}}
  \toprule
  External & Internal & Reference \\ \midrule
  \textbf{Sheaves of sets} \\
  sheaf of sets & set \\
  $\alpha : \F \to \G$ monomorphism & $\alpha$ injective & Ex.\@~\ref{ex:injective-surjective} \\
  $\alpha : \F \to \G$ epimorphism & $\alpha$ surjective & Ex.\@~\ref{ex:injective-surjective} \\
  $\Int(X \setminus \supp\F)$ & truth value of ``$\F$ is a singleton'' & Rem.\@~\ref{rem:support-sheaf-of-sets} \\
  $f : X \to \NN$ upper semicont.\@ & element of~$\widehat\NN$ & Lemma~\ref{lemma:upper-semicontinuous-functions} \\
  $f : X \to \NN$ locally constant & element of~$\NN$ & Lemma~\ref{lemma:upper-semicontinuous-functions} \\\\

  \textbf{Sheaves of rings} \\
  sheaf of rings & ring & Prop.\@~\ref{prop:rings-internally} \\
  local sheaf of rings & local ring & Prop.\@~\ref{prop:local-ring} \\
  $X$ is reduced & $\O_X$ is reduced (and $\neg\text{invertible} \Rightarrow \text{zero}$) & Prop.\@~\ref{prop:reduced-ring} \\
  $\dim X \leq n$ & Krull dimension of~$\O_X$ is~$\leq n$ & Prop.\@~\ref{prop:dimension-scheme-ox} \\
  $X$ is integral at all points & $\O_X$ is an integral domain & Prop.\@~\ref{prop:internal-integrality} \\
  $X$ is locally Noetherian & $\O_X$ is processly Noetherian & Prop.\@~\ref{prop:internal-noetherianity} (only~``$\Rightarrow$'' holds) \\
  $X$ is normal & $\O_X$ is normal (assuming that~$X$ is locally Noetherian) & Prop.\@~\ref{prop:normal-int-ext} \\
  relative spectrum~$\RelSpec_X{\A}$ & local spectrum $\Spec(\A|\O_X)$ & Thm.\@~\ref{thm:local-spectrum-yields-relative-spectrum} \\
  relative Proj~$\RelProj_X{\A}$ & local Proj $\Proj(\A|\O_X)$ & Thm.\@~\ref{thm:local-proj-yields-relative-proj} \\\\

  \textbf{Sheaves of modules} \\
  sheaf of modules & module \\
  $\F$ is finite locally free & $\F$ is finite free & Prop.\@~\ref{prop:locally-free} \\
  $\F$ is of finite type & $\F$ is finitely generated & Prop.\@~\ref{prop:finite-type-and-co} \\
  $\F$ is of finite presentation & $\F$ is finitely presented & Prop.\@~\ref{prop:finite-type-and-co} \\
  $\F$ is coherent & $\F$ is coherent & Prop.\@~\ref{prop:finite-type-and-co} \\
  $\F$ is quasicoherent & $\F[f^{-1}]$ is a sheaf wrt.\@~$(\speak{$f$ \inv} \Rightarrow \placeholder)$ for~$f\?\O_X$ & Thm.\@~\ref{thm:qcoh-sheafchar} \\
  $\F$ is flat & $\F$ is flat & Prop.\@~\ref{prop:flatness} \\
  $\F$ is torsion & $\F$ is torsion & Prop.\@~\ref{prop:torsion-int-ext} \\
  $\F$ is flabby & partially-defined elements of~$\F$ can be refined & Prop.\@~\ref{prop:internal-char-flabbiness} \\
  $\F$ is injective & $\F$ is injective & Thm.\@~\ref{thm:char-injectivity-modules} \\
  rank function of~$\F$ & minimal number of generators for~$\F$ & Prop.\@~\ref{prop:rank-function-internally} \\
  $M^\sim$ & $\ul{M}[\F^{-1}]$ (localization at generic filter) & Prop.\@~\ref{prop:tilde-construction-internally} \\
  tensor product $\F \otimes \G$ & tensor product $\F \otimes \G$ & Prop.\@~\ref{prop:internal-tensor-product} \\
  dual~$\F^\vee = \HOM_{\O_X}(\F,\O_X)$ & dual $\F^\vee = \Hom_{\O_X}(\F,\O_X)$ \\
  $\Int(X \setminus \supp\F)$ & truth value of~``$\F = 0$'' & Prop.\@~\ref{prop:characterization-support} \\
  quasicoherator of~$\I$ & $\{ s\?\O_X \,|\, \speak{$s$ \inv}
  \Rightarrow s \in \I \}$ ($\I$ a radical ideal) & Prop.\@~\ref{prop:quasicoherator-structure-sheaf} \\
  $\Omega^1_{X|S}$ & $\Omega^1_{\O_X | f^{-1}\O_S}$ & Prop.\@~\ref{prop:kaehler} \\
  Tor, sheaf Ext & Tor, Ext & Sect.\@~\ref{sect:sheaf-ext-and-tor} \\
  higher direct images & cohomology & Sect.\@~\ref{sect:higher-direct-images} \\\\

  \multicolumn{3}{@{}l@{}}{\textbf{Subspaces} ($i : A \hookrightarrow X$ closed immersion, $j : U \hookrightarrow X$ open immersion)} \\
  sheaf supported on~$A$ & $\Box$-sheaf, where~$\Box = (\placeholder \vee A^c)$ & Lemma~\ref{lemma:essim-closed-immersion} \\
  sheaf of the form~$j_*(\F)$ & $\Box$-sheaf, where~$\Box = (U \Rightarrow
  \placeholder)$ & \\
  extension of~$\F$ by the empty set & $j_!(\F) = \{ x\?\F \,|\, U \}$ & Lemma~\ref{lemma:extension-by-empty-set} \\
  extension of~$\F$ by zero & $j_!(\F) = \{ x\?\F \,|\, (x = 0) \vee U \}$ & Lemma~\ref{lemma:extension-by-zero} \\
  sheaf with empty/zero stalks on~$U^c$ & sheaf of the form~$j_!(\F)$ \\
  sections of~$\F$ are equal if they agree on dense open & $\F$ is $\neg\neg$-separated & Prop.\@~\ref{prop:negneg-sheaves} \\
  sheaf of sections of~$\F$ defined on dense open subsets & $\F^+$ with respect to~$\Box = \neg\neg$ (assuming that~$\F$ is~$\neg\neg$-separated) & Prop.\@~\ref{prop:negneg-sheaves} \\
  $U$ is dense & $\neg\neg U$ & Prop.\@~\ref{prop:modops-kripke} \\
  $U$ is scheme-theoretically dense & $\sdense U$, \ie $\O_X$ is separated
  wrt.~$(U \Rightarrow \placeholder)$ & Lemma\@~\ref{lemma:scheme-theoretical-denseness} \\
  $V(\I)$ is reduced & $\I$ is a radical ideal & Lemma~\ref{lemma:closed-subspace-reduced} \\
  $\O_{X_\mathrm{red}}$ & $\O_X/\sqrt{(0)}$ & Lemma~\ref{lemma:reduced-subspace} \\\\

  \multicolumn{3}{@{}l@{}}{\textbf{Over-schemes and over-toposes}} \\
  relative spectrum~$\RelSpec_X(\A)$ & local spectrum $\Spec(\A|\O_X)$ & Thm.\@~\ref{thm:local-spectrum-yields-relative-spectrum} \\
  relative Proj~$\RelProj_X(\S)$ & local Proj construction $\Proj(\S|\O_X)$ & Thm.\@~\ref{thm:local-proj-yields-relative-proj} \\
  big Zariski topos of~$X$ & local Zariski topos~$\Zar(\O_X|\O_X)$ & Thm.\@~\ref{thm:big-zariski-topos-of-relative-spectrum} \\
  big Zariski topos of~$\RelSpec_X(\A)$ & local Zariski topos~$\Zar(\A|\O_X)$ & Thm.\@~\ref{thm:big-zariski-topos-of-relative-spectrum} \\\\

  \multicolumn{3}{@{}l@{}}{\textbf{Rational functions and Cartier divisors}} \\
  $\K_X$ & total quotient ring of~$\O_X$ & Prop.\@~\ref{prop:kx-internally} \\
  Cartier divisor & element of~$\K_X^\times/\O_X^\times$ \\
  \effective Cartier divisor & $[s/1]$ with~$s\?\O_X$ regular & Def.\@~\ref{defn:effective-cartier-divisor} \\
  line bundle~$\O_X(D)$ & $D^{-1} \O_X \subseteq \K_X$ & Def.\@~\ref{defn:line-bundle-of-divisor} \\\\

  \multicolumn{3}{@{}l@{}}{\textbf{Topological properties}} \\
  $X$ is quasicompact & ``$\Sh(X) \models$'' commutes with directed disjunctions & Prop.\@~\ref{prop:quasicompact-meta} \\
  $X$ is local & ``$\Sh(X) \models$'' commutes with arbitrary disjunctions & Prop.\@~\ref{prop:local-meta} \\
  $X$ is irreducible & if $\neg(\varphi \wedge \psi)$, then $\neg\varphi$ or~$\neg\psi$ & Prop.\@~\ref{prop:irreducibility-internally} \\
  \bottomrule
\end{longtable}}

\section{The inference rules of intuitionistic logic}
\label{appendix:inference-rules}


\begin{center}
  \textbf{Structural rules} \\
  \vspace{-0.5em}
  \phantom{a}\hfill
  \AxiomC{$\phantom{\seq{\vec x}}$}\UnaryInfC{$\varphi \seq{\vec x} \varphi$}\DisplayProof\hfill
  \AxiomC{$\varphi \seq{\vec x} \psi$}\UnaryInfC{$\varphi[\vec s/\vec x]
  \seq{\vec y} \psi[\vec s/\vec x]$}\DisplayProof\hfill
  \AxiomC{$\varphi \seq{\vec x} \psi$}\AxiomC{$\psi \seq{\vec x}
  \chi$}\BinaryInfC{$\varphi \seq{\vec x} \chi$}\DisplayProof
  \phantom{a}\hfill
  \vspace{2.0em}

  \textbf{Rules for nullary and binary conjunction} \\
  \vspace{-0.5em}
  \phantom{a}\hfill
  \AxiomC{$\phantom{\seq{\vec x}}$}\UnaryInfC{$\varphi \seq{\vec x} \top$}\DisplayProof\hfill
  \AxiomC{$\phantom{\seq{\vec x}}$}\UnaryInfC{$\varphi \wedge \psi \seq{\vec x} \varphi$}\DisplayProof\hfill
  \AxiomC{$\phantom{\seq{\vec x}}$}\UnaryInfC{$\varphi \wedge \psi \seq{\vec x} \psi$}\DisplayProof\hfill
  \AxiomC{$\varphi \seq{\vec x} \psi$}\AxiomC{$\varphi \seq{\vec x} \chi$}\BinaryInfC{$\varphi \seq{\vec x} \psi \wedge \chi$}\DisplayProof
  \phantom{a}\hfill
  \vspace{2em}

  \textbf{Rules for nullary and binary disjunction} \\
  \vspace{-0.5em}
  \phantom{a}\hfill
  \AxiomC{$\phantom{\seq{\vec x}}$}\UnaryInfC{$\bot \seq{\vec x} \varphi$}\DisplayProof\hfill
  \AxiomC{$\phantom{\seq{\vec x}}$}\UnaryInfC{$\varphi \seq{\vec x} \varphi \vee \psi$}\DisplayProof\hfill
  \AxiomC{$\phantom{\seq{\vec x}}$}\UnaryInfC{$\psi \seq{\vec x} \varphi \vee \psi$}\DisplayProof\hfill
  \AxiomC{$\varphi \seq{\vec x} \chi$}\AxiomC{$\psi \seq{\vec x} \chi$}\BinaryInfC{$\varphi \vee \psi \seq{\vec x} \chi$}\DisplayProof
  \phantom{a}\hfill
  \vspace{2em}

  \textbf{Rules for arbitrary set-indexed conjunction and disjunction} \\
  \vspace{-0.5em}
  \phantom{a}\hfill
  \AxiomC{$\phantom{\seq{\vec x}}$}\UnaryInfC{$\bigwedge_{i \in I} \varphi_i \seq{\vec x} \varphi_j$ for all~$j \in I$}\DisplayProof\hfill
  \AxiomC{$\varphi \seq{\vec x} \psi_j$ for all~$j \in I$}\UnaryInfC{$\varphi \seq{\vec x} \bigwedge_{i \in I} \psi_i$}\DisplayProof
  \phantom{a}\hfill
  \vspace{1em}

  \phantom{a}\hfill
  \AxiomC{$\phantom{\seq{\vec x}}$}\UnaryInfC{$\varphi_j \seq{\vec x} \bigvee_{i \in I} \varphi_i$ for all~$j \in I$}\DisplayProof\hfill
  \AxiomC{$\varphi_j \seq{\vec x} \psi$ for all~$j \in I$}\UnaryInfC{$\bigvee_{i \in I} \varphi_i \seq{\vec x} \psi$}\DisplayProof
  \phantom{a}\hfill
  \vspace{2em}

  \textbf{Double rule for implication} \\
  \vspace{-0.5em}
  \phantom{a}\hfill
  \Axiom$\varphi \wedge \psi\ \fCenter\seq{\vec x} \chi$
  \doubleLine
  \UnaryInf$\varphi\ \fCenter\seq{\vec x} \psi \Rightarrow \chi$
  \DisplayProof
  \phantom{a}\hfill
  \vspace{2em}

  \textbf{Double rules for bounded and unbounded quantification} \\
  \vspace{-0.5em}
  \phantom{a}\hfill
  \Axiom$\varphi\ \fCenter\seq{\vec x, y} \psi$
  \doubleLine
  \UnaryInf$\exists y\?Y\_\! \varphi\ \fCenter\seq{\vec x} \psi$
  \DisplayProof
  {\tiny ($y$ not occurring in~$\psi$)}
  \hfill
  \Axiom$\varphi\ \fCenter\seq{\vec x, y} \psi$
  \doubleLine
  \UnaryInf$\varphi\ \fCenter\seq{\vec x\phantom{, y}} \forall y\?Y\_\! \psi$
  \DisplayProof
  {\tiny ($y$ not occurring in~$\varphi$)}
  \hfill\phantom{a}
  \vspace{1em}

  \phantom{a}\hfill
  \Axiom$\varphi\ \fCenter\seq{\vec x, Y} \psi$
  \doubleLine
  \UnaryInf$\exists Y\_\! \varphi\ \fCenter\seq{\vec x} \psi$
  \DisplayProof
  {\tiny ($Y$ not occurring in~$\psi$)}
  \hfill
  \Axiom$\varphi\ \fCenter\seq{\vec x, Y} \psi$
  \doubleLine
  \UnaryInf$\varphi\ \fCenter\seq{\vec x\phantom{, Y}} \forall Y\_\! \psi$
  \DisplayProof
  {\tiny ($Y$ not occurring in~$\varphi$)}
  \hfill\phantom{a}
  \vspace{2em}

  \textbf{Rules for equality} \\
  \vspace{-0.5em}
  \phantom{a}\hfill
  \AxiomC{$\phantom{\seq{\vec x}}$}
  \UnaryInfC{$\top \seq{x} x = x$}
  \DisplayProof
  \hfill
  \AxiomC{$\phantom{\seq{\vec x}}$}
  \UnaryInfC{$(\vec x = \vec y) \wedge \varphi \seq{\vec z} \varphi[\vec y/\vec x]$}
  \DisplayProof
  \hfill\phantom{a} \\[0.5em]
  (``$\vec x = \vec y\,$'' is short for~``$x_1 = y_1 \wedge \cdots \wedge x_n =
  y_n$''.)
\end{center}

\end{appendix}

\nocite{*}
\printbibliography

\vfill
\centering
\rotatebox{90}{\tiny Carina Willbold}\!
\includegraphics[width=0.5\textwidth]{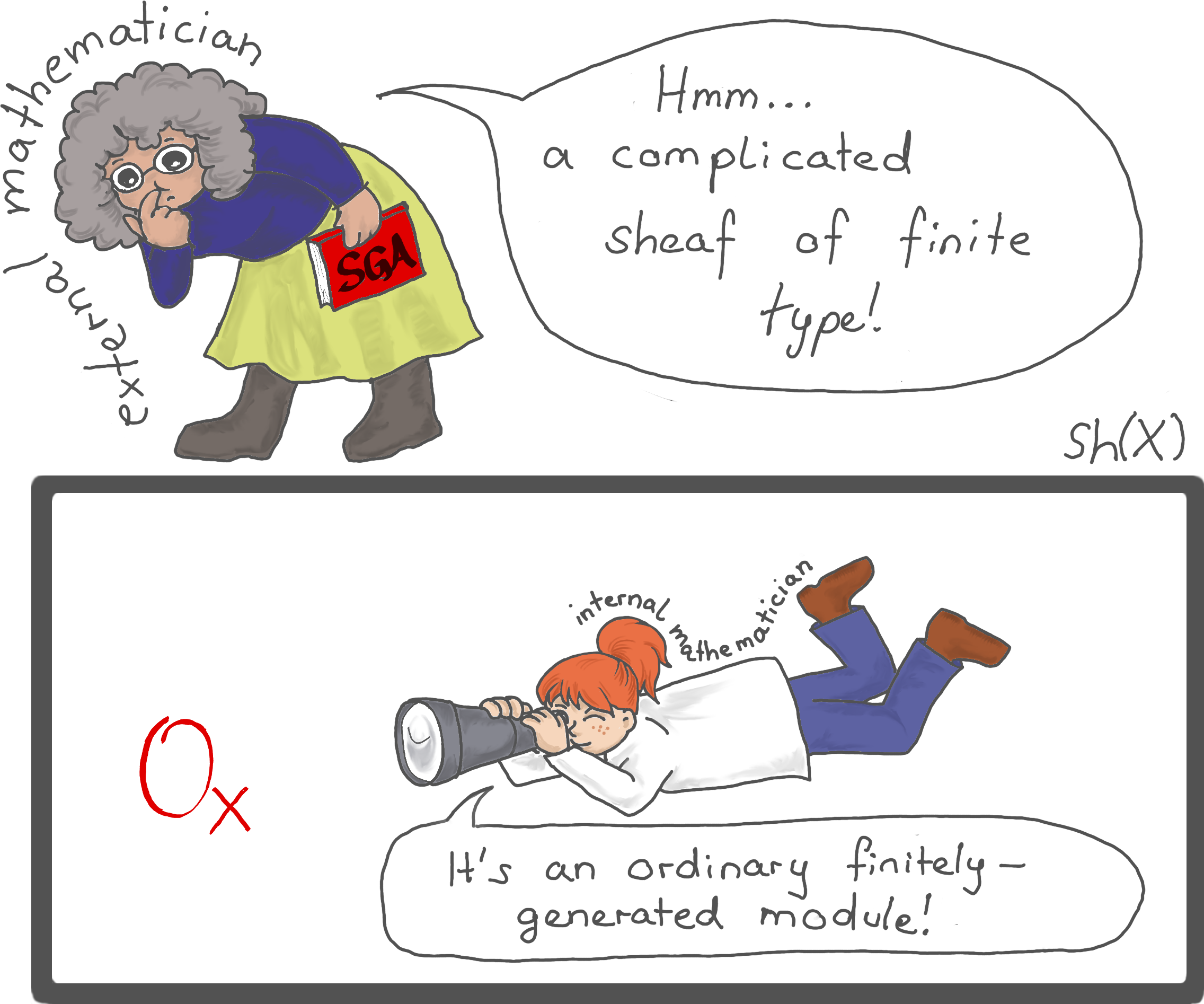}
\bigskip
\bigskip

\small
This text is licensed under the Creative Commons BY-SA 4.0 License.\par

\end{document}